\documentclass[11pt]{report}
\usepackage{a4wide}
\usepackage[all]{xy}
\usepackage{amsfonts}
\usepackage{amsthm}
\usepackage{amssymb}
\usepackage{amsmath}
\usepackage{appendix}
\usepackage[french,english]{babel}
\usepackage{caption}
\usepackage{comment}
\usepackage{dsfont}
\usepackage[T1]{fontenc}
\usepackage{geometry}
\usepackage{graphicx}

\usepackage[utf8]{inputenc}

\usepackage{mathrsfs}

\usepackage{stmaryrd}
\usepackage{ulem}

\usepackage[makeindex]{imakeidx}
\makeindex[columns=3]
\makeindex[title=Index of Definitions]
\makeindex[name=notation, title=Index of Symbols]
\makeindex[name=axiom, title=Index of Axioms]
\newcommand{\axiom}[1]{\index{#1}}
\usepackage{enumitem}

\newtheorem{Theorem}{Theorem}[section] 
\newtheorem{Proposition}[Theorem]{Proposition}
\newtheorem{Corollary}[Theorem]{Corollary}
\newtheorem{Definition}[Theorem]{Definition}

\newtheorem{Lemma}[Theorem]{Lemma}
\newtheorem{Example}[Theorem]{Example}
\newtheorem{Remark}[Theorem]{Remark}
\newtheorem{Definition/Proposition}[Theorem]{Definition/Proposition}

\newtheorem{Notation}[Theorem]{Notation}
\newtheorem{Assumption}[Theorem]{Assumption}

\theoremstyle{plain}

\makeatletter
\renewcommand\theequation%
{\thesection.\arabic{equation}}
\@addtoreset{equation}{section}
\makeatother

\makeatletter \@addtoreset{figure}{section}\makeatother

\usepackage[colorlinks=true,breaklinks=true,linkcolor=black]{hyperref}

\newcommand{\cA}{\mathcal{A}}
\newcommand{\fA}{\mathfrak{A}}

\newcommand{\sA}{\mathscr{A}}

\newcommand{\cB}{\mathcal{B}}
\newcommand{\fB}{\mathfrak{B}}

\newcommand{\sB}{\mathscr{B}}

\newcommand{\cC}{\mathcal{C}}

\newcommand{\sC}{\mathscr{C}}

\newcommand{\cD}{\mathcal{D}}

\newcommand{\cE}{\mathcal{E}}

\newcommand{\rE}{\mathrm{E}}

\newcommand{\cF}{\mathcal{F}}

\newcommand{\sF}{\mathscr{F}}

\newcommand{\fG}{\mathfrak{G}}

\newcommand{\cH}{\mathcal{H}}

\newcommand{\cI}{\mathcal{I}}

\newcommand{\cK}{\mathcal{K}}

\newcommand{\sL}{\mathscr{L}}

\newcommand{\cM}{\mathcal{M}}

\newcommand{\cN}{\mathcal{N}}
\newcommand{\fN}{\mathfrak{N}}

\newcommand{\sN}{\mathscr{N}}

\newcommand{\cO}{\mathcal{O}}

\newcommand{\cP}{\mathcal{P}}
\newcommand{\fP}{\mathfrak{P}}

\newcommand{\sP}{\mathscr{P}}

\newcommand{\fQ}{\mathfrak{Q}}

\newcommand{\cR}{\mathcal{R}}
\newcommand{\fR}{\mathfrak{R}}

\newcommand{\cS}{\mathcal{S}}
\newcommand{\fS}{\mathfrak{S}}

\newcommand{\sS}{\mathscr{S}}

\newcommand{\cT}{\mathcal{T}}
\newcommand{\fT}{\mathfrak{T}}

\newcommand{\cU}{\mathcal{U}}
\newcommand{\fU}{\mathfrak{U}}

\newcommand{\sU}{\mathscr{U}}

\newcommand{\cV}{\mathcal{V}}

\newcommand{\ff}{\mathfrak{f}}

\newcommand{\ffg}{\mathfrak{g}}

\newcommand{\fh}{\mathfrak{h}}

\newcommand{\fm}{\mathfrak{m}}

\newcommand{\fn}{\mathfrak{n}}

\newcommand{\fp}{\mathfrak{p}}

\newcommand{\fq}{\mathfrak{q}}

\newcommand{\fr}{\mathfrak{r}}

\newcommand{\fs}{\mathfrak{s}}

\newcommand{\bt}{\mathbf{t}}
\newcommand{\ft}{\mathfrak{t}}

\newcommand{\bv}{\mathbf{v}}

\newcommand{\bw}{\mathbf{w}}

\newcommand{\bx}{\mathbf{x}}

\newcommand{\FF}{\mathbb{F}}
\newcommand{\GG}{\mathbb{G}}

\newcommand{\QQ}{\mathbb{Q}}

\newcommand{\R}{\mathbb{R}}
\newcommand{\A}{\mathbb{A}}

\newcommand{\N}{{\mathbb{N}}}
\newcommand{\Z}{\mathbb{Z}}
\newcommand{\Q}{\mathbb{Q}}
\newcommand{\C}{\mathbb{C}}
\newcommand{\Ne}{{\Z}_{>0}}
\newcommand{\I}{\mathcal{I}}

\newcommand{\Id}{\mathrm{Id}}

\newcommand{\ines}{\mathrm{ines}}

\newcommand{\In}{\mathrm{Int}}
\newcommand{\conv}{\mathrm{conv}}
\newcommand{\cl}{\mathrm{cl}}
\newcommand{\fin}{\mathrm{fin}}

\newcommand{\htt}{\mathrm{ht}}
\newcommand{\Hom}{\mathrm{Hom}}
\newcommand{\qp}{\varpi}
\newcommand{\IN}{\I_{\sN}}
\newcommand{\IL}{\I_{\sL}}
\newcommand{\fho}{\mathring{\fh}}

\newcommand{\dw}{d^{W^v}}

\newcommand{\blue}[1]{{\color{blue}#1}}

\newcommand{\efface}[1]{}

\newcommand{\pr}{\mathrm{proj}}
\newcommand{\ve}{\mathrm{vert}}

\newcommand{\cur}{\curvearrowright}

\renewcommand{\emptyset}{\varnothing}

\renewcommand{\tilde}[1]{\widetilde{#1}}

\makeatletter
\def\Ddots{\mathinner{\mkern1mu\raise\p@
\vbox{\kern7\p@\hbox{.}}\mkern2mu
\raise4\p@\hbox{.}\mkern2mu\raise7\p@\hbox{.}\mkern1mu}}
\makeatother

\DeclareMathOperator{\supp}{Supp}

\DeclareMathOperator{\End}{End}
\DeclareMathOperator{\Int}{Int}

\DeclareMathOperator{\Aut}{Aut}

\DeclareMathOperator{\ad}{ad}

\newcommand{\Inv}{\mathrm{Inv}}

\newcommand{\kk}{\Bbbk}
\newcommand{\Stab}{\mathrm{Stab}}
\newcommand{\Fix}{\mathrm{Fix}}
\newcommand{\germ}{\mathrm{germ}}
\newcommand{\Ad}{\mathrm{Ad}}

\newcommand{\sch}{\mathrm{sch}}

\newenvironment{psmallmatrix}
  {\left(\begin{smallmatrix}}
  {\end{smallmatrix}\right)}

\newcommand{\alg}{\mathrm{alg}}

\newcommand{\ffgo}{\mathring{\ffg}}

\newcommand{\dd}{d}
\newcommand{\opp}{\ \mathrm{opp}\ }

\newcommand{\vw}{\mathrm{vw}}
\newcommand{\Isom}{\mathrm{Isom}}

\date{}

\title{Masures associated with split Kac-Moody groups over valued fields}
\author{Auguste \textsc{Hébert} \\Université de Lorraine, Institut Élie Cartan de Lorraine, F-54000 Nancy, France\\ UMR 7502,
auguste.hebert@univ-lorraine.fr}

\hypersetup{
pdfauthor={Hébert, Auguste},
pdftitle={Masures split Kac--Moody groups},
}
\begin{document}
\maketitle
\begin{abstract}
Masures are generalizations of Bruhat-Tits buildings adapted to the study of Kac--Moody groups over valued fields. They were introduced by Gaussent and Rousseau in 2007. Rousseau defined an axiomatic for these object and we simplified it. In this paper, which is mainly expository, we construct the masure associated with a split Kac--Moody group over a valued field, and we prove that it satisfies our axiomatic.
\end{abstract}
\tableofcontents
\section*{Introduction}
\subsection*{Kac--Moody groups, masures and applications}
\subsubsection*{Kac--Moody groups}
Kac--Moody groups are interesting generalizations of reductive group which integrate Kac--Moody (Lie) algebras, in some sense. There are several possible definitions of Kac--Moody groups, which mainly fall into two categories: the minimal ones and the completed ones. An example of a minimal Kac--Moody group is $\mathrm{SL}_n(\cK[t,t^{-1}])$, for $\cK$ a field. The maximal Kac--Moody group associated is $\mathrm{SL}_n(\cK(\!(t)\!))$. Here we are mainly interested in minimal Kac--Moody groups as defined by Tits in \cite{tits1987uniqueness}, but we also use Mathieu's completion as a tool to study them.

A split Kac--Moody group $\fG$ is a functor from the category of rings to the category of groups. It admits a split maximal torus $\fT$ and a  set of real roots $\Phi$. If $\cK$ is  a field, then $G:=\fG(\cK)$ is generated by $T$ and root subgroups $U_\alpha$, for $\alpha\in \Phi$, where for each $\alpha\in \Phi$, $U_\alpha$ is isomorphic to $(\cK,+)$. Note that unless $G$ is reductive, $\Phi$ is infinite and thus $G$ is infinite dimensional.

 For example if $G=\mathrm{SL}_n(\cK)$, one can choose $T$ as the subgroup of diagonal matrices, and one can choose  the root subgroups to be the subgroups $\fU_{i,j}(\cK)$ with ones on the diagonal, a coefficient in $\cK$ in position $(i,j)$, and $0$ elsewhere, with $(i,j)\in \llbracket 1,n\rrbracket$, $i\neq j$. If $G=\mathrm{SL}_n(\cK[t,t^{-1}])$, one can choose the torus to be the diagonal matrices with coefficients in $\cK^*$, and the root groups to be the $U_{i,j,k}=\fU_{i,j}(t^k \cK)$, with $i,j$ as above and  $k\in \Z$.
There are three types of indecomposable Kac--Moody groups: finite type, affine type and indefinite type. Apart from the finite type ones, which are reductive, the most studied and understood Kac--Moody groups are the untwisted affine groups. Up to central extensions and semi-direct product with $(\cK^\times)^n$ for some $n\in \N$, they are  of the form $\cK\mapsto \mathring{\fG}(\cK[t,t^{-1}])$ for $\mathring{\fG}$ a reductive group. Indefinite Kac--Moody groups are more mysterious. They are not linear in general (see \cite[9.3]{marquis2018introduction} for precise statements) and no ``explicit'' realization of them are known.

The study of Kac--Moody groups over non-Archimedean local fields began in the affine case with the works of Garland (\cite{garland1995cartan}). Braverman and Kazhdan then initiated the theory of Hecke algebras of affine Kac--Moody groups over local fields (\cite{braverman2011spherical}) and Gaussent and Rousseau generalized their works in the general Kac--Moody case, using the theory of masures, that they initiated in \cite{gaussent2008kac}.
\subsection*{Buildings and groups}
In order to study reductive and Kac--Moody groups, one can associate ``buildings'' to them. They are objects of geometric and combinatorial nature on which the groups under study act. One can then translate questions on the group in questions on the associated building. There are several types of buildings. Let $\fG$ be a reductive group (functor), that we assume to be split for simplicity. If $\cK$ is any field, then the Tits building of  $\fG(\cK)$ is a simplicial complex equipped with an action of $\fG(\cK)$. 

 If $\cF$ is a field equipped with a non-trivial real-valued valuation $\omega:\cF\rightarrow \R\cup \{+\infty\}$, one can then associate a Bruhat--Tits building $\I$ to $\fG(\cF)$ (this was done by Bruhat and Tits in \cite{bruhat1972groupes} and \cite{bruhat1984groupes}).  The Tits building of $\fG(\cF)$  can then be regarded as a boundary at infinity of $\I$. Let $\kk$ be the  residual field of $(\cF,\omega)$. Then the Tits building of $\fG(\kk)$ can be regarded as the tangent building at some vertex of $\I$. 
Let $\fB$ be a Borel subgroup of $\fG$. For example if $\fG=\mathrm{SL}_n$, with $n\in \Z_{\geq 2}$, one can take $\fB$ to be the subfunctor of $\fG$ consisting of the upper triangular matrices. The Iwahori subgroup $K_I$ of $\fG(\cF)$, is $\pr_{\kk}^{-1}(\fB(\kk))\subset \fG(\cO)$, if $\cO$ is the ring of integers of $\cF$ and $\pr_{\kk}:\fG(\cO)\rightarrow \fG(\kk)$ is the natural projection.  The Tits building of $\fG(\cK)$ is  a kind of geometric realization of $\fG(\cK)/\fB(\cK)$, and if $\omega$ is discrete, we can regard $\I$ as a geometric realization of $\fG(\cF)/K_I$.

Consider now a split Kac--Moody group $\fG$. If $\cK$ is a field, the analogue of the Tits building of $\fG(\cK)$ is a twin building. This is a pair $(\I_\cK^+,\I_{\cK}^-)$ of buildings, equipped with a codistance (a twinning) between $\I_{\cK}^+$ and $\I_{\cK}^-$. The fact that there are naturally two buildings instead of one comes from the fact that there are two conjugacy classes of Borel subgroups, contrary to the reductive case. If $\fB^+$ and $\fB^-$ are two opposite Borel subgroup, i.e there are Borel subgroups whose intersection is a maximal split torus, then one can take $\I_{\cK}^+=\fG(\cK)/\fB^+(\cK)$ and $\I_{\cK}^-=\fG(\cK)/\fB^-(\cK)$. 
Let $(\cF,\omega)$ be as above. The object corresponding to the Bruhat--Tits building  of $(\fG,\cF,\omega)$ is the masure $\I$. At infinity of the masure, we find the twin building of $\fG(\cF)$, and the tangent building at some vertex is the twin building of $\fG(\kk)$. 
\subsection*{Masures and applications}
Masures  enabled to prove many results in the theory of Kac--Moody groups over non-Archimedean local fields. To cite but a few, it enabled progresses in the following fields: \begin{itemize}
\item In the theory of MV-cycles for Kac--Moody groups, see \cite{gaussent2008kac},  \cite{bozec2023MV} and \cite{bouthier2025geometric}. Bouthier and Vasserot indirectly use masures in their work on the geometric Satake equivalence, since they use finiteness results proved using this object.

\item In the theory of Hecke algebras for Kac--Moody groups over valued fields. It enabled Bardy-Panse, Gaussent and Rousseau to define a spherical Hecke algebra (\cite{gaussent2014spherical}) and an Iwahori--Hecke algebra (\cite{bardy2016iwahori}) in general, extending the constructions of Braverman, Kazhdan and Patnaik in the untwisted affine case (see \cite{braverman2011spherical} and  \cite{braverman2016iwahori}). Hébert and Abdellatif also associated Hecke algebras to certain parahoric subgroups in \cite{abdellatif2019completed}. 
\item In Kazhdan--Lusztig theory for Kac--Moody groups, it enabled Muthiah, Hébert and Philippe to give a conjectural definition of Kazhdan--Lusztig $R$-polynomials (\cite{muthiah2019double} and \cite{hebert2024affine}), based on a   theory of twin masures (\cite{bardy2025twin}). 
\end{itemize}
\section*{Axioms of masures}
\subsection*{Bruhat--Tits buildings}
Let $(\cF,\omega)$ be a valued field, $\fG$ be a split reductive group or Kac--Moody group and $G=\fG(\cF)$.  We now describe how to construct the Bruhat--Tits building of $G$ or its masure and what are its properties.

First assume that $\fG$ is reductive. Let $\fT$ be a maximal split torus of $\fG$. Let $Y$ be the cocharacter lattice of $(\fG,\fT)$ and let $\A=Y\otimes \R$ be the \textbf{standard apartment}. Let $N$ be the normalizer of $T=\fT(\cF)$ in $G$. Then $N$ naturally acts on $\A$ by affine automorphisms.  Let $\Phi$ be the root system of $(\fG,\fT)$. For $\alpha\in \Phi$, we have an isomorphism $x_\alpha:\cF\rightarrow U_\alpha$.

 The Bruhat-Tits building $\I=\I(\fG,\cF,\omega)$ is defined as follows. For $x\in \A$, one sets $G_x=\langle N_x,U_{\alpha,x}\mid \alpha\in \Phi\rangle\subset G$, where $N_x$ is the fixator of $x$ in $N$ and $U_{\alpha,x}=x_{\alpha}(\{a\in \cF\mid \alpha(x)+\omega(a)\geq 0\})$ for $\alpha\in \Phi$. For example $G_0=\fG(\cO)$. Then $\I$ is defined as $(G\times \A)/\sim$, for some equivalence relation $\sim$ such that for all $x\in G$, $G_x$ is the fixator of $x$ in $\A$, where $\A$ is regarded as a subset of $\I$ via the embedding $\A \hookrightarrow \I$, $a\mapsto [a,1]$ (see Subsection~\ref{ss_d_mas}).

  For example, if $\fG=\mathrm{SL}_n$, with $n\in \N_{\geq 2}$, then  $N$ is the group of monomial matrices and $\Phi$ can be written $\Phi=\{\alpha_{i,j}\mid (i,j)\in \llbracket 1,n \rrbracket\}$. For $i,j\in \llbracket 1,n\rrbracket$ one can define $x_{\alpha_{i,j}}$ as follows. For $a\in \cF$,   $x_{\alpha_{i,j}}(a)=I_n+a E_{i,j}$, where $E_{i,j}$ is the matrix with a one in position $(i,j)$ and $0$ elsewhere. The subgroup $G_x$ is computed in Proposition~\ref{p_fix_pt}, when $n=2$.

\begin{figure}
\centering
\includegraphics[scale=0.3]{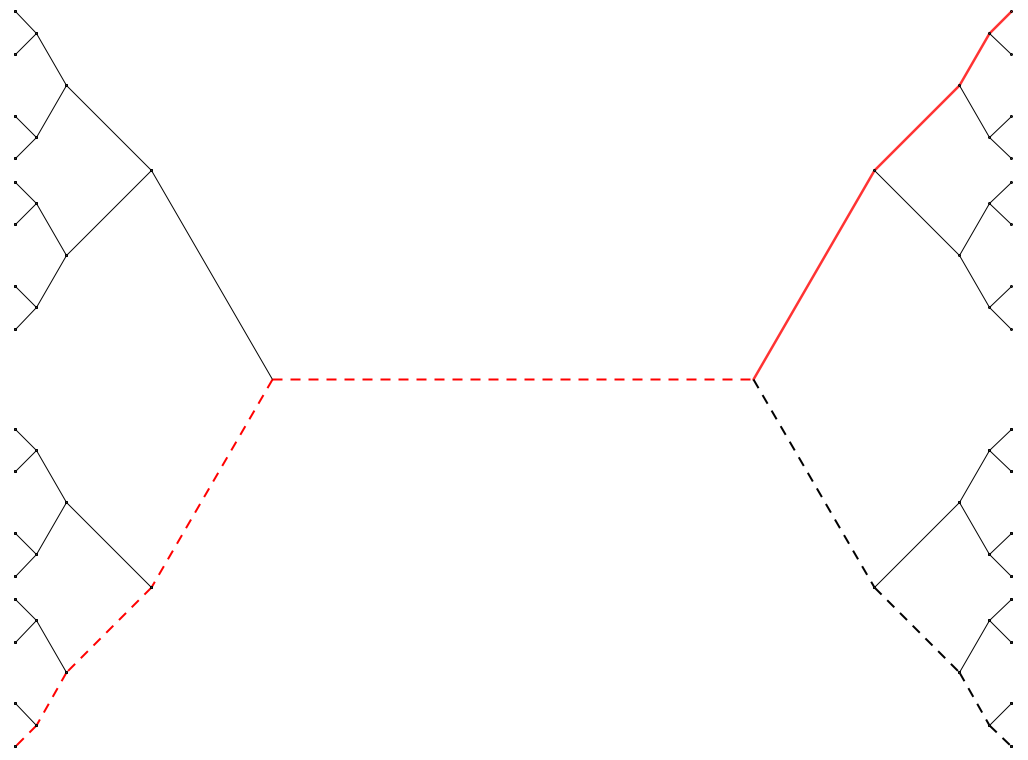}
\caption{The tree of $\mathrm{SL}_2(\cF)$ for $\cF=\Q_2$ or $\cF=\FF_2(\!(t)\!)$. The dotted ``line'' and the red ''line'' are two apartments}\label{f_tree_SL2}
\end{figure}

\begin{figure}[h]
\centering
\includegraphics[scale=0.3]{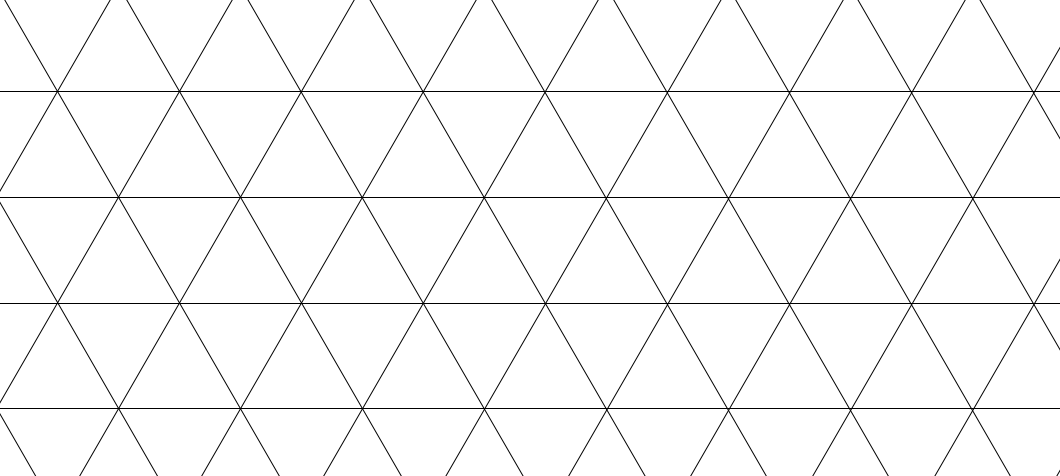}
\caption{Affine apartment of $\mathrm{SL}_3$.}\label{f_aprt_SL3}
\end{figure}

\begin{figure}[h]
\centering
\includegraphics[scale=0.3]{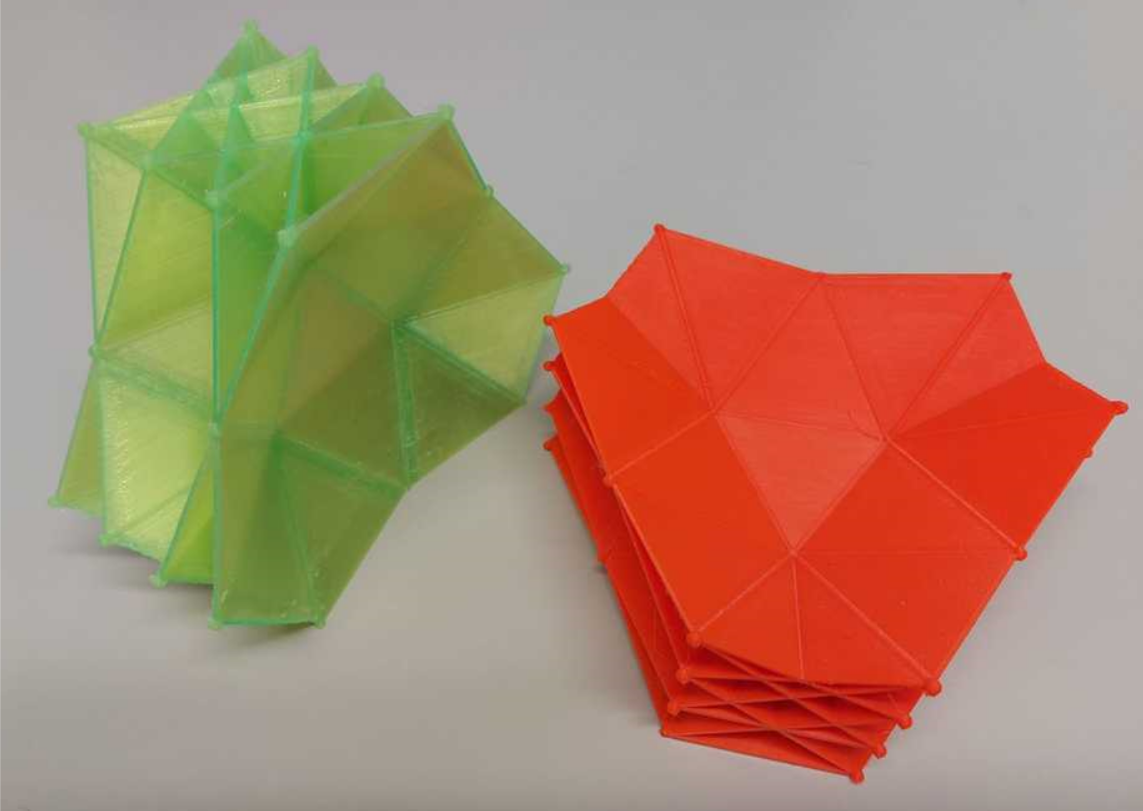}
\caption{Two 3D printed models of parts of the building of $\mathrm{SL}_3(\cF)$, for $\cF\in \{\QQ_2,\FF_2(\!(t)\!)\}$. }\label{f_building_SL3}
\end{figure}

The set $\I$ is then covered by its \textbf{apartments}, where an apartment is a set of the form $g.\A$, for some $g\in G$. Let $\Phi$ be the root system of $(\fG,\fT)$.  Then $\A$ is equipped with the hyperplane arrangement $\cH=\{\alpha^{-1}(\{k\})\mid k\in \omega(\cF^\times)\}$ whose elements are called \textbf{walls}. Then, using the action of $G$, we can equip every apartment with its wall arrangement. A \textbf{half-apartment} is then a half-space of an apartment delimited by a wall. When $\omega(\cF^\times)$ is discrete, $\cH$ naturally equips $\A$ with the structure of a CW-complex. Then $\I$ is an affine $\R$-building in the sense of \cite[§2.4.4]{rousseau2023euclidean}. In particular, it satisfies the following two properties: \begin{itemize}
\item[\textbf{($I_1$)}] for every two faces, there exists an apartment $A$ containing these two faces,
\item[\textbf{($I_2$)}] for every two apartments $A$ and $A'$, the set  $A\cap A'$ is a finite intersection of half-apartments of $A$ and  there exists $g\in G$ such that $g.A=A'$ and $g$ fixes $A\cap A'$. 
\end{itemize}
Note that property ($I_1$) is the building's theoretic translation of the Bruhat decomposition $G=G_F\cdot N \cdot G_{F'}$, for every two faces $F,F'$ of $\A$, where $G_F$ and $G_{F'}$ are respectively the fixators of $F$ and $F'$ in $G$ and $N$ is the normalizer of $\fT(\cF)$ in $G$.
For example, when $\fG=\mathrm{SL}_2$, $(\A,\cH)\simeq (\R,\Z)$ and $\I$ is a homogeneous tree of valency $|\kk|+1$, where $\kk$ is the residue field of $(\cF,\omega)$ (see Figure~\ref{f_tree_SL2}). When $\fG=\mathrm{SL}_3$, $\A\simeq \R^2$ and the walls of $\A$ form a tessellation of $\A$ with regular triangles (see Figure~\ref{f_aprt_SL3}). The building is then a union of apartments, which intersect along walls, see Figure~\ref{f_building_SL3}).

\subsection*{Masures}
Let $\fG$ be a maximal split reductive group and let $\cF$ be a field equipped with a non-trivial valuation $\omega$. Let $\fT$ be a maximal split torus, $G=\fG(\cF)$ and  $T=\fT(\cF)$. Let $\A=Y\otimes \R$, where $Y$ is the cocharacter lattice of $(\fG,\fT)$. In order to define the masure $\I$ of $G$, we proceed as in the reductive case. We define a group $G_x$ for each $x\in \R$, and we then define $\I$ as $G\times \A/\sim$, for some equivalence relation $\sim$ such that $G_x$ is the fixator of $x$ in $\A$. The definition of $\sim$ is the same as in the reductive case, the difficulty is to define the $G_x$, for $x\in \A$.

One defines an action of  $N$ on $\A$ by affine automorphisms. This enables to define the fixator $N_x$, for every $x\in \A$. This defines the $N$-component of $G_x$, similarly as in the reductive case. To define the whole  $G_x$ however, we use the positive and negative Mathieu's completions of $G$.

More precisely, let $\ffg$ be the  Kac--Moody algebra  over $\C$ associated with $G$. Let $\fh$ be a Cartan subalgebra of $\ffg$.  Then $\ffg$ admits a root space decomposition $\ffg=\bigoplus_{\alpha\in \Delta\cup \{0\}}\ffg_\alpha$, where $\ffg_\alpha=\{x\in \ffg\mid [h,x]=\alpha(x)x,\forall h\in \fh\}$, for $\alpha\in \Delta\cup \{0\}$. The set $\Delta$ is a subset of $\fh^*$ called the set of roots of $(\ffg,\fh)$.   There are two types of roots: the real ones (whose set is denoted $\Phi$) and the imaginary ones. To each real root $\alpha$ is associated a root subgroup $U_\alpha$, isomorphic to $(\cF,+)$ by some isomorphism $x_\alpha$. Then $G=\langle T,U_\alpha\mid \alpha\in \Phi\rangle$. Among $\Delta$, we can choose a subset $\Delta_+$ of positive roots, such that $\Delta=\Delta_+\sqcup \Delta_-$, where $\Delta_-=-\Delta_+$. Set $\Phi_{\pm}=\Delta_{\pm}\cap \Phi$ and $U^{\pm}:=\langle U_\alpha\mid \alpha\in \Phi_{\pm}\rangle\subset G$. Then $U^+$ and $U^-$ are the positive and negative standard unipotent subgroups of $G$. In the reductive case, $\Phi=\Delta$ and the natural multiplication map  $\prod_{\alpha\in \Phi_{\pm}}U_{\alpha}\rightarrow  U^{\pm}$ is a bijection, for any order on $\Phi_{\pm}$. In the Kac--Moody case however, this is no longer the case. One can define  completions $U^{ma+}$ and $U^{ma-}$ of $U^+$ and $U^-$ respectively. Then Mathieu's completion  $G^{ma\epsilon}$ satisfies  $G^{ma\epsilon}=\langle U^{ma\epsilon},U^{-\epsilon},T\rangle$, for both $\epsilon\in \{-,+\}$.
 An advantage of $U^{ma+}$ or $U^{ma-}$ is that we have ``coordinates'' to describe their elements. One can define a $\Z$-form $\ffg_{\Z}$ of $\ffg$, which also admits a root space decomposition $\ffg_{\Z}=\bigoplus_{\alpha\in \Delta\cup \{0\}}\ffg_{\alpha,\Z}$. Then using ``twisted exponentials'', Rousseau defines maps $X_\alpha:\ffg_{\alpha,\Z}\otimes \cF\rightarrow U^{ma+}$, for any $\alpha\in  \Delta_+$. When $\alpha\in \Phi_+$, one can take $X_\alpha=x_\alpha$. Then every element of $U^{ma+}$ decomposes uniquely as an infinite product $\prod_{\alpha\in \Delta_+}X_\alpha(\underline{u_\alpha})$, where $(\underline{u_\alpha})\in \prod_{\alpha\in \Delta_+} \ffg_{\alpha,\Z}\otimes \cF$, for any order on $\Delta_+$ satisfying certain conditions (and similarly for $\Delta_-$). If $x\in \A$ and $\epsilon\in \{-,+\}$, then one defines: \[U^{ma\epsilon}_x=\{\prod_{\alpha\in \Delta^\epsilon} X_\alpha(\underline{u_\alpha})\mid (\underline{u_\alpha})\in \prod_{\alpha\in \Delta^\epsilon} \ffg_{\alpha,\Z}\otimes \cF, \omega(\underline{u_\alpha})+\alpha(x)\geq 0\}\text{ and }U_{x,\infty}^\epsilon=U^{ma\epsilon}_x\cap G \]  Then $G_x$ is defined as $G_x=\langle U_{x,\infty}^+,U_{x,\infty}^-,N_x\rangle\subset G$.
  A crucial property of $G_x$ is that it decomposes as \[G_x= U_{x,\infty}^+\cdot U_{x,\infty}^-\cdot N_x= U_{x,\infty}^-\cdot U_{x,\infty}^+\cdot N_x\] (see Theorem~\ref{t_G_x}). Note that eventually, we obtain that $U_{x,\infty}^\epsilon=U_{x,\fin}^{\epsilon}:=\langle U_{\alpha,x}\mid \alpha\in \Phi^\epsilon\rangle$ (see Proposition~\ref{p_sphrcl_para}). We could thus define $G_x$ using the minimal group $G$ only. However, we use the masure to prove this equality, and we do not know how to prove it directly. Note that we change the notation compared to the litterature: $U_{x,\infty}^+$, $U_{x,\infty}^-$ and $U_{x,\fin}^{\pm}$ are usually denoted $U_x^{pm+}$, $U_{x}^{nm-}$  and $U_{x}^{\pm}$ respectively.

\begin{figure}[h]
\centering
\includegraphics[scale=0.3]{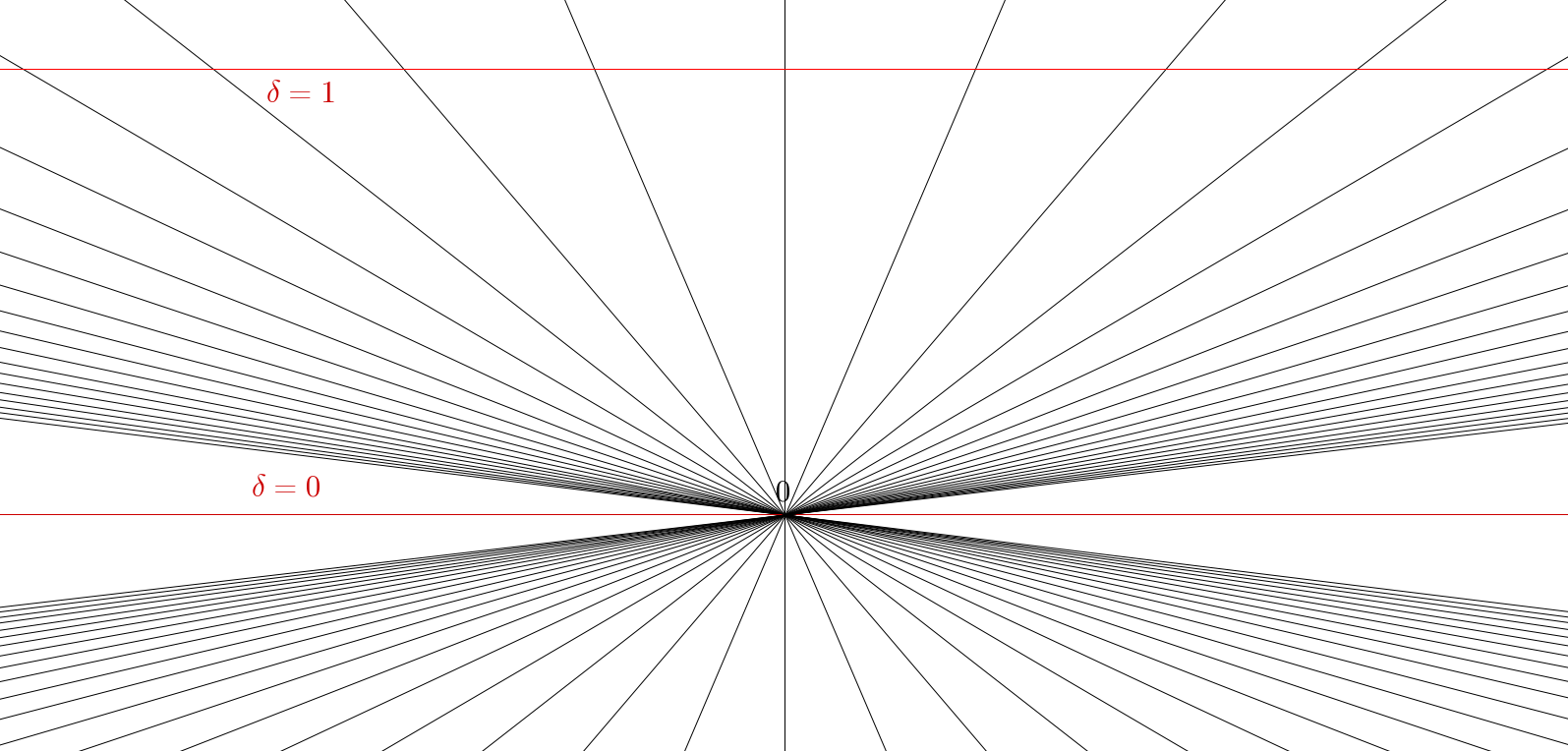}
\caption{Set of walls containing $0$ of the standard apartment of  affine $\mathrm{SL}_2$. Not all the walls are represented since there are infinitely many. If we draw a copy of $\Z$ in the line $\delta=1$, then the set of walls is the set of lines containing $0$ and an element of $\Z$. }\label{f_vect_aprt_SL2_aff}
\end{figure}

\begin{figure}[h]
\centering
\includegraphics[scale=0.5]{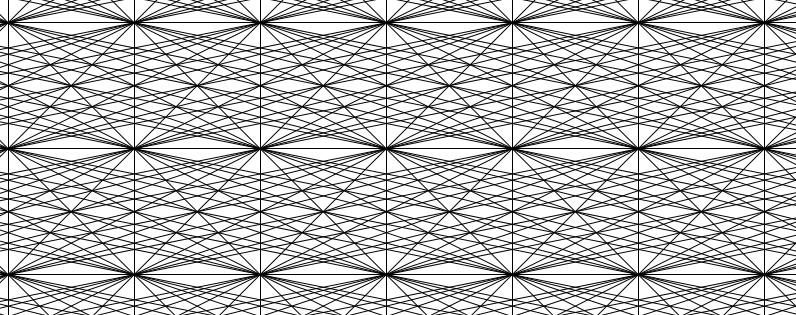}
\caption{Set of walls of the standard apartment of  affine $\mathrm{SL}_2$, when $\omega(\cF^\times)=\Z$. Not all the walls are represented since we would obtain a black picture. This hyperplane arrangement is obtained by drawing a copy of $\Z^2$ and the walls are then the lines containing two points of $\Z^2$.}\label{f_aff_aprt_SL2_aff}
\end{figure}

 Another difficulty in the Kac--Moody case is that (unless $\fG$ is reductive),  $\Phi$ is infinite and thus the hyperplane arrangement $\cH=\{\alpha^{-1}(\{k\})\mid k\in \omega(\cF^\times)\}$ can be dense, even if $\omega(\cF^\times)$ is discrete (see Figures~\ref{f_vect_aprt_SL2_aff} and \ref{f_aff_aprt_SL2_aff}). Therefore $\cH$ does not equip $\A$ with the structure of a CW-complex. A solution to this issue is to define the faces of $\A$ as filters on $\A$ (instead of subsets of $\A$).  Note that this solution is already used in the reductive case, when $\omega(\cF^\times)$ is not discrete.

 Then $\I$ is the union of its apartments, which are the $g.\A$ for $g\in G$. We prove that $\I$ satisfies the following two axioms (see Definition~\ref{d_w_mas} and Theorem~\ref{t_I_abst_mas}): \begin{itemize}
\item [(MA2)] Let $A$ and $B$ be two apartments. Then there exists finitely many half-apartments $D_1,\ldots,D_n$ of $A$ such that $A\cap B=\bigcap_{i=1}^n D_i$ and there exists $g\in G$ such that $g.A=B$ and $g$ fixes $A\cap B$. 
\item[(MA3)] Let $\fR_1$ be a splayed chimney germ and $\fR_2$ be a face or a chimney germ. Then there exists an apartment $A$ containing $\fR_1$ and $\fR_2$. 
\end{itemize}
Axiom~\ref{a_ma2} improves and replaces some of the axioms defined by Rousseau in \cite{rousseau2011masures}. Note that it involves only half-apartments defined by real roots (and not imaginary) although the construction of $\I$ uses Mathieu's completion of $G$, which involves imaginary roots. 
   Axiom~\ref{a_ma3} was introduced in \cite{rousseau2011masures} (in a slightly different version). It replaces the fact that two faces are always contained in a common apartment in a Bruhat-Tits building (axiom ($I_1$)). This fact is not true for masures since the Cartan decomposition does not hold on the entire $G$. Chimneys-germs are objects at the infinity of the masure including sector-germs. Axiom~\ref{a_ma3} is the masure theoretic translation of the decomposition \begin{equation}\label{e_MA3}
   G=G_{\fR_1}\cdot N\cdot G_{\fR_2},
   \end{equation} where $G_{\cV}$ denotes the fixator of $\cV$, if $\cV$ is a  filter   on $\A$.     For example, sector-germs at infinity are particular cases of chimney-germs at infinity and thus  the Bruhat and Birkhoff decompositions $G=U^\epsilon \cdot N \cdot U^{\epsilon'}$, for $\epsilon,\epsilon'\in \{-,+\}$ are particular cases of \eqref{e_MA3}. The Iwasawa decomposition $G=U^\epsilon \cdot N\cdot G_F$, where $F$ is a face of $\A$, is also an example of this decomposition.

\subsection*{Chronology of the theory of masures}
Masures were introduced by Gaussent and Rousseau in \cite{gaussent2008kac}. In this paper, they are associated with split Kac--Moody groups over discretely valued fields whose residual fields contains $\C$. Rousseau obtained a construction valid for any field $\cF$ equipped with an $\R$-valued valuation in \cite{rousseau2016groupes}.  Rousseau also defined masures associated with almost-split Kac--Moody groups in \cite{rousseau2017almost}, using works of Charignon (\cite{charignon2010immeubles}).

Rousseau gave an axiomatic definition of masures in \cite{rousseau2011masures} and proved that masures associated with  Kac--Moody groups satisfy its axioms. We simplified this definition in \cite{hebert2020new} (in the affine case) and in general in \cite{hebert2022new}: we proved that under a freeness assumption on the simple coroots, his axioms are equivalent to~\ref{a_ma2} and~\ref{a_ma3}.

Here we prove directly that $\I$ satisfies these axioms.  We start from scratch: we recall most of the necessary facts on Kac--Moody groups and we do not assume that the reader is familiar with the theory of Bruhat--Tits buildings or masures. We restrict to the case of split Kac--Moody groups, except in Chapter~\ref{C_w_mas}.

\section*{New results of this paper}
The aim of this paper is to  give a self-contained construction of the masure associated to a split Kac--Moody group and to give a direct proof that it satisfies~\ref{a_ma2} and~\ref{a_ma3}. While it is mainly expository, we obtain a few new results that we list here.

\subsection*{Filters have good fixators} For $\cV$ a filter on $\A$, we denote by $G_\cV$ its fixator in $G$. When they defined masures in \cite{gaussent2008kac}, Gaussent and Rousseau introduced the notion of ``good fixator'', which enabled them to prove many useful properties of masures. The filter $\cV$ is said to have a good fixator if $G_\cV$ acts transitively on the set of apartments containing $\cV$ and if $G_\cV$ decomposes as: \[G_\cV=(U^+\cap G_\cV)\cdot (U^-\cap G_\cV)\cdot (N\cap G_{\cV})=(U^-\cap G_\cV)\cdot (U^+\cap G_\cV)\cdot (N\cap G_{\cV}).\] We prove that actually, every filter on $\cV$ has a good fixator (see Theorem~\ref{t_fltrs_gd_fix}).

\subsection*{Hecke paths} Most of the applications of the theory of masures to Kac--Moody groups over local fields use Hecke paths. If $x\in \I$, the retraction $\rho_{+\infty}(x)$ of $x$ on $\A$ along $U^+$ is the unique element of $U^+.x\cap \A$. In \cite[Theorem 6.2]{gaussent2008kac}, the authors prove that if $x,y\in \I$ are preordered (a condition defined in \cite[5.8]{rousseau2016groupes}), then $\rho_{+\infty}([x,y])$ is a Hecke path, which means that it is a piecewise affine path satisfying precise folding conditions (see Definition~\ref{d_H_path}). We generalize this result by dropping the assumption on $x$ and $y$ (see Theorem~\ref{t_H_pth}). While for all the current applications in arithmetic, the segment is assumed to be preordered, this generalization is useful in our proof that $\I$ satisfies~\ref{a_ma2}.

\subsection*{On the axiom~\ref{a_ma2}}
In \cite{hebert2022new}, we proved that $\I$ satisfies~\ref{a_ma2}, provided that the sets of simple roots and the sets of simple coroots are free families. We drop the freeness assumption on the set of simple co-roots. To do so, we prove that $\I$ can be regarded as a quotient of a masure $\tilde{\I}$ whose set of coroots is free. We then deduce the result for $\I$ from the result for $\tilde{\I}$. 

\subsection*{Description of the Iwahori subgroup via the Borel subgroup}
Let $C_0^+$ be the fundamental alcove of $\A$. Let $I=G_{C_0^+}$ be the \textbf{Iwahori subgroup of $G$}.  Let $\fB$ be the Borel subgroup $\fT\cdot \fU^+$. Let $\cO$ be the ring of integers of $\cF$ and $\fm$ be its maximal ideal. Let $\kk$ be the residual field of $\cO/\fm$. Let $\pr_{\kk}:\fG(\cO)\rightarrow \fG(\kk)$ be the natural projection induced by the projection $\cO\rightarrow \kk$. We prove that $I=\pr_{\kk}^{-1}(\fB(\kk))$ (see Proposition~\ref{p_Iwahori}), which generalizes a well-known fact in the reductive case.

\subsection*{Description of the fixator of a spherical local face}
Let $C$ be an alcove of $\A$. In \cite[2.4]{bardy2025twin}, it is proved that $G_C=\langle N_C, U_{\alpha,C}\mid \alpha\in \Phi \rangle$, where $N_C=\fT(\cO)$ is the fixator of $C$ in $N$. We slightly generalize this result here. We prove that if $F$ is 
a spherical local face of $\A$ (i.e the fixator of its direction in $W^v$ is finite), then $G_F=\langle N_F,U_{\alpha,F}\mid \alpha\in \Phi\rangle$ (see Proposition~\ref{p_sphrcl_para}).

\section*{Organization of the paper}

This paper is organized as follows.

In Chapter~\ref{C_tree}, we construct and describe the tree of $\mathrm{SL}_2$, which is the ``simplest'' possible example of masure. This chapter is intended to give an intuition of what a masure is and is mainily independent of the others.

In Chapter~\ref{C_KM_alg}, we define and study Kac--Moody algebras. We then describe the vectorial apartments associated, which are the standard apartments of the masure equipped with the walls containing $0$.

In Chapter~\ref{C_splt_KM_grps}, we introduce minimal and completed Kac--Moody groups.

In Chapter~\ref{C_w_mas}, we introduce abstract axioms of masures. We then study masures from an axiomatic viewpoint and prove that a series of axioms implies~\ref{a_ma2}. Masures of the chapter are not associated with a group. In particular, our results also apply to masures associated with almost-split Kac--Moody groups.

In Chapter~\ref{C_d_mas}, we define and study the fixators $G_x$, for $x\in \A$, when $G=\fG(\cF)$ is a split Kac--Moody group over a valued field. We define the action of $N$ on $\A$ and we associate subgroups of $U^{ma+}$ to filters on $\A$.

In Chapter~\ref{C_prop}, we define the masure of $G$ and we prove that it satisfies our axiomatic. We also obtain descriptions of the fixator $G_\cV$, for $\cV$ a filter on $\A$. 

In Chapter~\ref{C_para_min_gp}, we describe the fixator $G_{\cV}$ of $\cV$, in certain particular cases, where $\cV$ is a subface of the fundamental alcove $C_0^+$.

\paragraph{Acknowledgments}

\begin{figure}[h!]
    \centering
   \def\svgwidth{0.6\linewidth}
    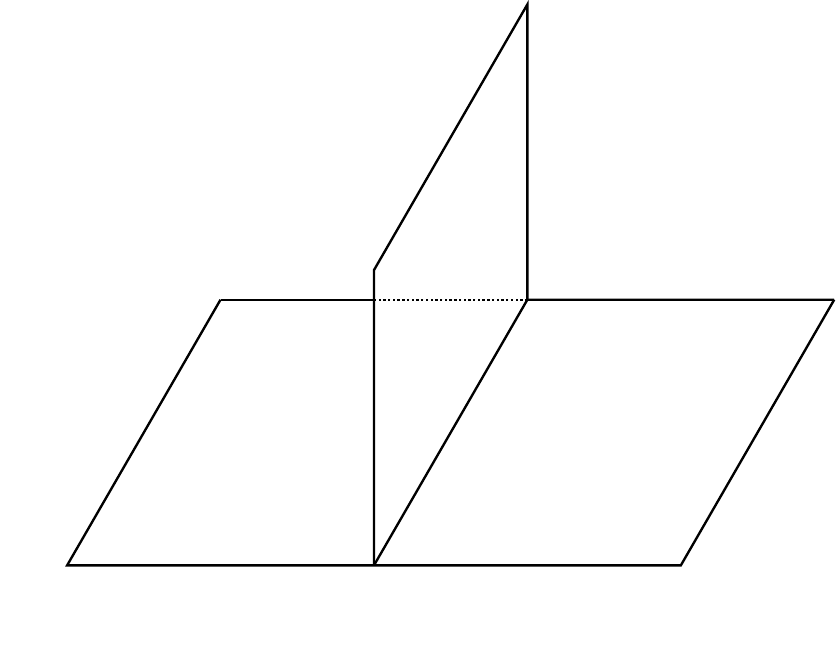
    \caption{This figure illustrates the alcoves associated with a sector-germ at a point. Here $A$ is an apartment, $\fQ$ is a sector-germ of $A$ and $\fQ'$ is a sector-germ of $\I$ adjacent to $\fQ$ but not contained in $A$. The half-apartments $D_1$ and $D_2$ are the half-apartments such that $D_i\cup \fQ'$ is contained in an apartment $A_i$ of $\I$, for both $i\in \{1,2\}$ and $M$ is the wall $D_1\cap D_2$. All the $\fQ_a$, $a\in \A$ are parallel, but the direction of $\fQ_a'$ depends on the position of $a$: if $a\in \mathring{D_1}$, then $\fQ_a'$ is an alcove of $A$ adjacent and different from $\fQ_a$, if $a\in \mathring{D_2}$, then $\fQ_a'=\fQ_a$ and if $a\in M$, then $\fQ'_a$ is not contained in $A$. In general if $\fQ''$ is a sector-germ of $\I$ and $(P_i,A_i,f_i)_{i\in \llbracket 1,k\rrbracket}$ is a decomposition of $A$ with respect to $\fQ''$ (which exists by Proposition~\ref{p_splt_apt}),  then for all $i\in \llbracket 1,k\rrbracket$ and $a_1,a_2\in \mathring{P_i}$, we have $\fQ''_{a_1},\fQ_{a_2}''\subset A$ and $\fQ''_{a_1}$ and $\fQ''_{a_2}$ are parallel. }
    \label{fig:l_adj_SG}
\end{figure}

I am very much indebted to  Guy Rousseau for the numerous discussions, clarifications, comments and answers to my questions, for his thorough reading of previous versions of the manuscript and for the many suggestions of corrections and  improvements. I am also indebted to  Stéphane Gaussent, for introducing me to the theory and for many exchanges on the subject. I also thank Nicole Bardy-Panse, Yvann Gaudillot--Estrada, Timothée Marquis, Dinakar Muthiah, Manish Patnaik and  Paul Philippe for stimulating conversations. I thank Raphael Appenzeller for Figures~\ref{f_building_SL3} and \ref{fig:l_adj_SG}. The figures of this paper were made using Geogebra and Inkscape.

\section*{General notations}

The set $\N$ is the set of non-negative integers.

If $a,b$ are in an affine space or in $\R\cup \{-\infty,+\infty\}$, we denote by $[a,b[$ (resp. $]a,b]$, $]a,b[$) the set $[a,b]\setminus\{b\}$ (resp. $[a,b]\setminus \{a\}$, $[a,b]\setminus \{a,b\}$). If $a,b\in \R$, we denote by $\llbracket a,b\rrbracket$ the set $[a,b]\cap \Z$. 

The notation ``$\ll 0$'' and ``$\gg 0$'' mean ``small enough'' and ``big enough''. For example if $\cP_k$ is a property, for $k\in \N$, then $\cP_k$ is true for $k\gg 0$ if there exists $N\in \N$ such that $\cP_k$ is true for all $k\geq N$. 

We denote by $\lfloor \cdot \rfloor:\R\rightarrow \Z$ and by $\lceil \cdot \rceil:\R\rightarrow \Z$ the floor function and the ceiling function respectively.

Let $A$ be a finite dimensional affine space. Let $C\subset A$ be non-empty. We denote by $\conv(C)$ the \textbf{convex hull} of $C$, i.e the smallest convex subset of $A$ containing $C$.  We denote the interior of $C$ by $\mathring{C}$ or by $\In(C)$, depending on the context. We denote by $\overline{C}$ the closure of $C$ in $A$ and by $\partial C$ the boundary of $C$, i.e $\overline{C}\setminus \mathring{C}$.   The \textbf{support of }$C$, denoted $\supp(C)$ is the smallest affine subspace of $A$ containing $C$. The \textbf{dimension} of $C$ is the dimension of its support.The \textbf{relative interior} of $C$ is its interior, regarded as a subset of its support. We denote it $\In_r(C)$. If $F,F'$ are subsets of $A$, we set that $F$ \textbf{dominates}\index{d@domination} $F'$ if $\overline{F}\supset \overline{F'}$. 

We usually denote the filters by $\cV, \cV',\ldots$ and elements of these filters by $\Omega,\Omega',\ldots$.

We usually denote a ring by $\cR$, a field by $\cK$ and a valued field by $(\cF,\omega)$. 

We usually denote functors with a gothic typographie and their value at $\cF$ or $\cK$ with a roman letter, for example $\fG$ (for a Kac--Moody functor), $\fU^\pm$ (for its positive and negative unipotent subgroups) etc. and $G$ (for its $\cF$ points), $U^{\pm}$ etc.

If $(E_j)_{j\in J}$ is a family of disjoint subsets of a set $E$, we often write $\bigsqcup_{j\in J}E_j$ instead of $\bigcup_{j\in J} E_j$. This implicitly means that the $E_j$ are disjoint.

\chapter{The tree of $\mathrm{SL}_2$}\label{C_tree}

In this chapter, we construct and study the tree of $\mathrm{SL}_2(\cF)$, when $\cF$ is a real-valued field. The aim of this chapter is to give an example of a Bruhat-Tits building in order to give an intuition of what Bruhat-Tits buildings look like. This chapter is mainly independent of the others.

 We use the ``norm approach'' introduced in \cite{goldman1963space} (when $\cF$ is a local field) and used in \cite{parreau2000immeubles} (see \cite{parreau2023affine} for an english translation) and \cite[6]{rousseau2023euclidean}. We then relate it to the ``lattice approach'' used in \cite{serre1977arbres} (see \cite{serre1980trees} for an english translation) (when $\cF$ is discretely valued) and \cite{chiswell2001introduction} (for a non-necessarily real-valued valuation). All these are very particular cases of the Bruhat-Tits buildings introduced in \cite{bruhat1972groupes} and \cite{bruhat1984groupes}. For simplicity, we choose to study only $\mathrm{SL}_2(\cF)$, but most of the results are easily adaptable to the case of $\mathrm{SL}_n(\cF)$, for $n\in \Z_{\geq 2}$. We use non necessarily discrete valuations: we choose this level of generality because in this set-up, faces are naturally filters, whereas when the valuation is discrete, they can be seen as sets. For masures, the faces are filters even if the valuation is discrete.

\section{Definition of the tree}

Let $\cF$\index[notation]{f@$\cF$} be a field equipped with a non-trivial valuation $\omega:\cF\rightarrow \R\cup \{+\infty\}$\index[notation]{o@$\omega$}. We set $\Lambda=\omega(\cF^*)$\index[notation]{l@$\Lambda$}.  Let $\cO=\{x\in \cF\mid \omega(x)\geq 0\}$\index[notation]{o@$\cO$} be the valuation ring of $(\cF,\omega)$. We now construct the Bruhat-Tits building of $(\mathrm{SL}_2(\cF),\omega)$ as a quotient of a space of norms on $\cF^2$. This follows Goldman--Iwahori's and Parreau's constructions (\cite{goldman1963space} and \cite{parreau2000immeubles}).  For $a\in \cF$, we set $|a|=e^{-\omega(a)}$\index[notation]{z@$\mid\cdot\mid $}. We set $G=\mathrm{SL}_2(\cF)$\index[notation]{g@$G$}.

\begin{Definition}
A \textbf{ultrametric norm} on $\cF^2$ is a map $\eta:\cF^2\rightarrow \R_{\geq 0}$ such that for all $u,v\in \cF^2$, $a\in \cF$, we have:  \begin{enumerate}
\item $\eta(v)=0$ if and only if $v=0$,

\item $\eta(av)=|a|\eta(v)$,

\item $\eta(u+v)\leq \max(\eta(u),\eta(v))$.
\end{enumerate}

Let $b=(b_1,b_2)$ be a basis of $\cF^2$ A norm $\eta$ is called \textbf{adapted to }$b$ if we have $\eta(\lambda b_1+\mu b_2)=\max  (\eta(\lambda b_1),\eta(\mu b_2))$, for all $\lambda,\mu\in \cF$.   A norm $\eta$ of $\cF^2$ is called \textbf{adaptable} if there exists a basis with respect to which $\eta$ is adapted. We denote by $\sN$\index[notation]{n@$\sN$} the set of adaptable ultrametric norms on $\cF^2$. 
\end{Definition}

Note that when $\cF$ is a local field,  every ultrametric norm on $\cF^2$ is adaptable, by \cite[Proposition 1.1]{goldman1963space}. More generally, if the valuation $\omega$ is discrete and if $\cF$ is  Henselian and complete, then every ultrametric norm on $\cF^2$ is adaptable, by \cite[Proposition 15.1.11]{kaletha2023bruhat}.

\begin{Definition}
Let $\eta,\eta'\in \sN$. We write $\eta\sim \eta'$  if there exists $a\in \R_{>0}$ such that $\eta'=a\eta$. Then $\sim$ is an equivalence relation. If $\eta\in \sN$, we denote by $[\eta]$ its equivalence class. We set $\I_{\sN}=\sN/\sim$.\index[notation]{i@$\I_{\sN}$} This is the \textbf{tree of }$(\mathrm{SL}_2(\cF),\omega)$. 

The group $\mathrm{GL}_2(\cF)$  acts on $\I_{\sN}$ via $g.\eta(v)=\eta(g^{-1}.v)$, for $g\in \mathrm{GL}_2(\R)$, $\eta\in \sN$ and $v\in \cF^2$. By restriction, this induces an action of $G$ on $\I_{\sN}$.
\end{Definition}

Let $b=(b_1,b_2)$ be a basis of $\cF^2$. For $x=(x_1,x_2)\in \R^2$, define $\eta_{b,x}\in \sN$ by \[\eta_{b,x}(v_1b_1+v_2b_2)=\max (e^{-x_1}|v_1|,e^{-x_2}|v_2|),\] for $(v_1,v_2)\in \cF$. We have $[\eta_{b,x}]=[\eta_{b,x+(y,y)}]$ for all $y\in \R$ and thus one sets $[\eta_{b,x+\R(1,1)}]=[\eta_{b,x}]$. Note that a norm $\eta$ is adapted to a basis $b$ if and only if $\eta=\eta_{b,x}$, for some $x\in \R^2$. Indeed, if $\eta$ is adapted to $b$, we have \begin{equation}\label{e_norm_adapted_basis}
\eta=\eta_{b,(-\ln(\eta(b_1)),-\ln(\eta(b_2))}, 
\end{equation} and the converse is clear.

\section{Apartments}

In this section, we  introduce the apartments of $\I_{\sN}$. These are subsets of $\I_{\sN}$, isometric to $\R$, which cover $\I_{\sN}$. 

Apartments are a crucial tool in the study of buildings and of the groups acting on them. For example we can translate certain decompositions of $G$ in terms of existence of apartments containing pairs of filters, see for instance Corollary~\ref{c_frndl_pr_tr} and Propositions~\ref{p_Bruhat} and \ref{p_Br_dec}.

\subsection{Apartments}
Let $b$ be a basis of $\cF^2$. 
 The   \textbf{apartment} associated with $b$ is the set $A_b=\{[\eta_{b,x}]\mid x\in \R^2/\R(1,1)\}$. It is easy to check that $x\mapsto \eta_{b,x}$ is a bijection between $\R^2/\R(1,1)$ and $A_b$. Therefore we can regard apartments as copies of $\R$.  The \textbf{standard apartment} $\A$\index[notation]{a@$\A$} is the apartment associated with the canonical basis $e=(e_1,e_2)$ of $\cF^2$. Explicitely, we  identify $\A$ and $\R$ via the map \begin{equation}\label{e_identification_A_R}
\R\mapsto \A, \ x\mapsto [\eta_{e,(0,x)}].
\end{equation}

An \textbf{apartment of $\I_{\sN}$} is a set of the form $A_b$, for some basis $b$ of $\cF^2$. If $\eta\in \sN$, then $[\eta]\in A_b$, for every basis adapted to $\eta$, by \eqref{e_norm_adapted_basis}. In particular, $\I_\sN$ is the union of all its apartments. 

If $a\in \cF^\times$, $x\in \R$ and $b=(b_1,b_2)$ is a basis of $\cF^2$, we have: \begin{equation}\label{e_relation_norms}
\eta_{ab,x}=e^{\omega(a)} \eta_{b,x}, \eta_{(b_1,ab_2),x}=\eta_{b,x-\omega(a)}\text{ and }\eta_{(b_2,b_1),x}=e^{-x}\eta_{(b_1,b_2),-x}.
\end{equation}

The apartments are subsets of $\IN$ which are isometric to a line: when $\omega$ is discrete, they are bi-infinite paths of edges without u-turns (see Figure~\ref{f_tree_SL2}). When $(\cF,\omega)$ is complete, one can show that the apartments are exactly the subsets of $\IN$ isometric to $\R$. 

\subsection{Action of the torus and of its normalizer on $\A$}\label{ss_Action_N_A}

Set $T_{\mathrm{GL_2}(\cF)}=\begin{psmallmatrix}\cF^\times & 0 \\ 0 & \cF^\times\end{psmallmatrix}\subset \mathrm{GL}_2(\cF)$\index[notation]{t@$T_{\mathrm{GL}_2(\cF)}$}.

Let $T$\index[notation]{t@$T$} be the standard torus of $G$, i.e $T=G\cap T_{\mathrm{GL_2}(\cF)}$. For $x\in \cF^\times$, we set $\bt_x=\begin{psmallmatrix} x & 0\\ 0 & x^{-1}\end{psmallmatrix}\in T.$

\begin{Lemma}\label{l_act_torus}
\begin{enumerate}

\item The subgroup $T_{\mathrm{GL}_2(\cF)}$ stabilizes $\A$ and acts by translation on $\A$. 

\item  Let $t=\begin{psmallmatrix} x & 0 \\ 0 & y\end{psmallmatrix}\in T_{\mathrm{GL}_2(\cF)}$ and $z\in \R$. Then $t.\eta_{e,z}=|x|^{-1}\eta_{e,z-\omega(x)+\omega(y)}$. In particular, $t.[\eta_{e,z}]=[\eta_{e,z-\omega(x)+\omega(y)}]$ and if $x\in \cF^\times$, then $\bt_x$ acts by translation by the vector $2\omega(x)$ on $\A$. 
\end{enumerate}
\end{Lemma}

\begin{proof}
We  prove (2), which implies (1). Let $v=(v_1,v_2)\in \cF^2$. Then \begin{align*}t.\eta_{e,z}(v)=\eta_{e,z}(t^{-1}.v)=\eta_{e,z}(x^{-1} v_1, y^{-1}v_2)&=\max(|x^{-1}v_1|,e^{-z}|y^{-1}|v_2)
\\&=|x^{-1}|\max(|v_1|,e^{-z}|xy^{-1}||v_2|)\\
&=|x^{-1}|\max(|v_1|,e^{-z}|xy^{-1}||v_2|)\\
&=|x^{-1}|\max(|v_1|,e^{-z-\omega(x)+\omega(y)}|v_2|)\\
&=|x^{-1}|\eta_{e,z+\omega(x)-\omega(y)}(v). \end{align*}

Therefore $t.\eta_{e,z}=|x^{-1}|\eta_{e,z+\omega(x)-\omega(y)}$ and the lemma follows.
\end{proof}

\begin{Lemma}\label{l_act_G_norms}
\begin{enumerate}
\item Let $b$ be a basis of $\cF^2$, $g\in \mathrm{GL}_2(\cF)$ and $x\in \R$. Then $g.\eta_{b,x}=\eta_{g.b,x}$.

\item The group $G$ acts transitively on the set of apartments of $\IN$. 

\end{enumerate}
\end{Lemma}

\begin{proof}
1) Let $\lambda,\mu\in \cF$. Then $g.\eta_{b,t}(\lambda g.b_1+\mu g.b_2)=\eta_{b,x}(\lambda b_1+\mu b_2)=\max(|\lambda|,|\mu|e^{-x})=\eta_{g.b,x}(\lambda g.b_1+\mu g.b_2)$, which proves 1). 

2) Let $b$ be a basis of $\cF^2$. As $\mathrm{GL}_2(\cF)$ acts transitively on the set of bases of $\cF^2$, there exists $g\in \mathrm{GL}$  such that $g.e=b$. Then $g.\A=A_b$. Let $t=\begin{psmallmatrix} 1/\det(g) & 0 \\ 0 & 1\end{psmallmatrix}$. Then $gt\in G$ and  by Lemma~\ref{l_act_torus}, $gt.\A=A_{b}$, which proves the lemma.

\end{proof}

Let $N=N_G(T)$\index[notation]{n@$N=N_G(T)$} be the normalizer of $T$ in $G$. We have $N=T\sqcup (\begin{psmallmatrix}0 & *\\ *&0\end{psmallmatrix}\cap G)$. Let $\tilde{s}=\begin{psmallmatrix}0 & 1\\ -1 & 0\end{psmallmatrix}$. Then $\tilde{s}$ stabilizes $\A$ and acts as $-\Id$ on $\A\simeq \R$. Therefore $N$ stabilizes $\A$.  We will see later (see Proposition~\ref{p_A2}) that $N$ is actually the stabilizer of $\A$ in $G$. 

\section{Vertices of $\I_\sN$ and lattice viewpoint}
We now give the lattice description of the tree of $G$. This appears in \cite{goldman1963space} and this is the viewpoint adopted in \cite{serre1977arbres} (where $\Lambda=\Z$) and \cite[4.3]{chiswell2001introduction} (in \cite{chiswell2001introduction}, $\omega$ does not necessarily take values in $\R$). We prove that the lattice tree can be naturally identified with the set of vertices of $\I_{\sN}$.

A \textbf{lattice} of $\cF^2$ is a subset of $\cF^2$ of the form $\cO b_1\oplus \cO b_2$, for some basis $(b_1,b_2)$ of $\cF^2$. We say that two lattices $L_1,L_2$ of $\cF^2$ are equivalent if there exists $\lambda\in \cF^\times$ such that $L_1=\lambda L_2$. If $L$ is a lattice of $\cF^2$, we denote by $[L]$ its equivalence class. The \textbf{lattice tree $\I_{\sL}$ of $G$ with respect to $\omega$} is the set of lattices of $\cF^2$ quotiented by the equivalence relation. The group $\mathrm{GL}_2(\cF)$ (and thus  $G$) acts on $\I_{\sL}$ via $g.[L]=[g.L]$ for $g\in \mathrm{GL}_2(\cF)$ and $L$ a lattice of $\cF^2$.

A \textbf{vertex} of $\I_\sN$ is a point of $\I_\sN$ of the form $[\eta_{b,0}]$, for some basis $b$ of $\cF^2$. We denote by $\ve(\I_\sN)$\index[notation]{v@$\ve(\I_\sN)$} the set of vertices of $\I_\sN$. We set $0_\I=[\eta_{e,0}]$, where $e$ is the canonical basis of $\cF^2$.

If $\eta$ is a norm on $\cF^2$ and $r\in \R$, we set $\overline{B_\eta}(r)=\{x\in \cF^2\mid \eta(x)\leq e^{-r}\}$ and $B_\eta(r)=\{x\in \cF^2\mid \eta(x)<e^{-r}\}$. For $r\in \R$, we set $\cF_{\geq r}=\{x\in \cF\mid \omega(x)\geq r\}$.

\begin{Lemma}\label{l_equivalence_lattices_norms}
Let $b,b'$ be two bases of $\cF^2$. Then:\begin{enumerate}
\item  $\cO b_1\oplus \cO b_2=\cO b_1'\oplus \cO b_2'$ if and only if $\eta_{b,0}=\eta_{b',0}$. 

\item $\eta_{b,0}\in \R_{>0} \eta_{b',0}$ if and only if $\eta_{b,0}\in e^{\Lambda} \eta_{b',0}$, where $\Lambda=\omega(\cF^\times)$.

\item for all $a\in \cF$, we have $\eta_{b,0}=e^{\omega(a)} \eta_{b',0}$ if and only if $\cO b_1\oplus \cO b_2=a (\cO b_1'\oplus \cO b_2')$. 
\end{enumerate}
\end{Lemma}

\begin{proof}
1) If $\eta:=\eta_{b,0}=\eta_{b',0}=:\eta'$, then we have $\overline{B_{\eta}}(0)=\overline{B_{\eta'}}(0)=\cO b_1\oplus \cO b_2=\cO b_1'\oplus \cO b_2'$· Conversely, let $b,b'$ be two basis of $\cF^2$ be such that $\cO b_1\oplus \cO b_2=\cO b_1'\oplus \cO b_2'$. Let $\eta=\eta_{b,0}$ and $\eta'=\eta_{b',0}$.  Then $b_1,b_2\in \cO b_1\oplus \cO b_2'$ and thus for every $r\in \R$, we have $\overline{B_{\eta}}(e^{-r})= \cF_{\geq r}b_1\oplus \cF_{\geq r} b_2\subset \cF_{\geq r} b_1'\oplus \cF_{\geq r}b_2'=\overline{B_{\eta'}}(e^{-r})$. By symmetry, we deduce $\overline{B_{\eta}}(e^{-r}) =\overline{B_{\eta'}}(e^{-r})$  and similarly, $B_{\eta}(e^{-r})=B_{\eta'}(e^{-r})$. Let now $x\in \cF^2$. Then \[\{r\in \R\mid x\in \overline{B_{\eta}}(r)\setminus B_{\eta}(r)\}=\{r\in \R\mid x\in \overline{B_{\eta'}}(r)\setminus B_{\eta'}(r)\}=\{-\ln(\eta(x))\}=\{-\ln(\eta'(x))\},\]  which proves that $\eta(x)=\eta'(x)$ and hence $\eta=\eta'$.

2) Assume that $\eta_{b,0}=t \eta_{b',0}$, for some $t\in \R_{>0}$. Write $b_1'=\lambda b_1+\mu b_2$, where $\lambda,\mu \in \cF$. Then $\eta_{b,0}(b_1')=\max (|\lambda|,|\mu|)=\exp(-\min(\omega(\lambda),\omega(\mu))=t\in e^{\Lambda}$.

3) By \eqref{e_relation_norms}, we have $\eta_{ab',0}=e^{\omega(a)}\eta_{b',0}$ and then (3) follows from (1). 
\end{proof}

\begin{Proposition}
The map $\psi:\I_{\sL}\rightarrow \ve(\I_{\sN})$ defined by $[\cO b_1 \oplus \cO b_2]\mapsto [\eta_{b,0}]$ is well-defined and is a $G$-equivariant bijection. 
\end{Proposition}

\begin{proof}
Let $b=(b_1,b_2),b'=(b_1',b_2')$ be two bases of $\cF^2$ such that $[\cO b_1\oplus \cO b_2]=[\cO b_1'\oplus \cO b_2']$. Then there exists $t\in \cF^\times $ such that $\cO tb_1'\oplus \cO tb_2'=\cO b_1\oplus \cO b_2$ and thus by Lemma~\ref{l_equivalence_lattices_norms}, we have $[\eta_{b,0}]=[\eta_{tb',0}]=[\eta_{b',0}]$. This proves that $\psi$ is well-defined. By definition, for every $x\in \ve(\I_{\sN})$, there exists a basis $b$ of $\cF^2$ such that $x=[\eta_{b,0}]$, and then $x=\psi(\cO b_1\oplus \cO b_2)$, which proves that $\psi$ is surjective.  

Let $L=\cO b_1\oplus \cO b_2$ and $L'=\cO b_1'\oplus \cO b_2'$ be two lattices of $\cF^2$. Assume that $\psi([L])=\psi([L'])$. Then we have $\eta_{b,0}=t\eta_{b',0}$, for some $t\in \R_{>0}$. By Lemma~\ref{l_equivalence_lattices_norms} (2), we can write $t=e^{-\omega(a)}$, with $a\in \cF^\times$. Then by Lemma~\ref{l_equivalence_lattices_norms} (3), $[\cO b_1\oplus \cO b_2]=[a(\cO b_1'\oplus \cO b_2')]=[L]=[L']$, which proves that $\psi$ is injective. The $G$-equivariance of $\psi$ is clear since $g.\eta_{b,0}=\eta_{g.b,0}$ for every basis $b$ of $\cF^2$ and $g\in G$. 
\end{proof}

Using $\psi$, we now regard $\I_{\sL}$ as a subset of $\IN$.

Let $b=(b_1,b_2)$ be a basis of $\cF^2$. The associated apartment $A^\sL_b$ is $\{[\cO b_1\oplus \cO xb_2]\mid x\in \cF^\times\}$. Let $\Lambda=\omega(\cF^\times)\subset \R$. For $\lambda\in \Lambda$, choose $x_\lambda\in \cF$ such that $\omega(x_\lambda)=\lambda$.  Then the map $\Lambda\rightarrow A^\sL_b$ defined by $\lambda\mapsto [\cO b_1\oplus \cO x_\lambda b_2]$, for $\lambda\in \Lambda$, is a  bijection.

\begin{Lemma}\label{l_description_A_L}
Let $\A_{\sL}=\A\cap \I_{\sL}$. Then we have $\A_{\sL}=\{[\eta_{e,y}]\mid y\in \Lambda\}\simeq \Lambda$. 
\end{Lemma}

\begin{proof}
Let $a\in \A_{\sL}$. Write $a=[\eta_{b,0}]$, where $b$ is a basis of $\cF^2$. Then by definition of $\A$, there exists $y\in \R$ such that $a=[\eta_{e,y}]$ and thus there exists $t\in \R_{>0}$ such that $\eta_{e,y}=t\eta_{b,0}$. Write $e_1=\lambda_1 b_1+\mu_1 b_2$ and $e_2=\lambda_2 b_1+\mu_2 b_2$, where $\lambda_1,\lambda_2,\mu_1,\mu_2\in \cF$. We have $\eta_{e,y}(e_1)=t\eta_{b,0}(e_1)=t\max(|\lambda_1|,|\mu_1|)$, which proves that $t\in e^{\Lambda}$. Then $\eta_{e,y}(e_2)=e^{-y}=t\max (|\lambda_2|,|\mu_2|)\in e^{\Lambda}e^{\Lambda}=e^{\Lambda}$ and thus $y\in \Lambda$. This proves that $\A_{\sL}\subset \{[\eta_{e,y}]\mid y\in \Lambda\}$. 

Let now $y\in \Lambda$. Write $y=e^{\omega(a)}$, with $a\in \cF$. Then by \eqref{e_relation_norms}, we have $\eta_{(e_1,a^{-1}e_2),0}=\eta_{e,y}$ and thus $[\eta_{e,y}]\in \A_\sL$. Therefore $\{[\eta_{e,y}]\mid y \in \Lambda\}\subset \A_{\sL}$, and the lemma follows. 
\end{proof}

\section{Fields extensions}

Let $\tilde{\cF}$ be a field extension of $\cF$ and $\tilde{\omega}$ be a valuation on $\tilde{\cF}$, with values in $\R$. We say that $(\tilde{\cF},\tilde{\omega})$ is an \textbf{extension} of $(\cF,\omega)$ if  $\tilde{\omega}|_{\cF}=\omega$. 

In this section, we prove that if $(\tilde{\cF},\tilde{\omega})$ is an extension of $(\cF,\omega)$, then the building of $(\mathrm{SL}_2,\cF,\omega)$ embeds in the one of $(\mathrm{SL}_2,\tilde{\cF},\tilde{\omega})$. This enables to consider every point of $\I_{\sN}$ as a vertex of some over-tree.

For $x\in \R$, we set $\cF_{\geq x}=\{a\in \cF\mid \omega(a)\geq x\}$ and if $(\tilde{\cF},\tilde{\omega})$ is an extension of $(\cF,\omega)$, we set  $\tilde{\cF}_{\geq x}=\{a\in \tilde{\cF}\mid \tilde{\omega}(a)\geq x\}$.

\begin{Lemma}\label{l_ext_nrms}
Let $(\tilde{\cF},\tilde{\omega})$ be an extension of $(\cF,\omega)$. Let $b=(b_1,b_2)$ and $b'=(b_1',b_2')$ be two bases of $\cF^2$, $x,x'\in \R$ and $t\in \R_{>0}$. Assume that $\eta_{b,x}^{(\cF,\omega)}=t\eta_{b',x'}^{(\cF,\omega)}$. Then $\eta_{b,x}^{(\tilde{\cF},\tilde{\omega})}=t\eta_{b',x'}^{(\tilde{\cF},\tilde{\omega})}$.
\end{Lemma}

\begin{proof}
Let $\eta=\eta_{b,x}^{(\cF,\omega)}$, $\eta'=\eta_{b',x'}^{(\cF,\omega)}$, $\tilde{\eta}=\eta_{b,x}^{(\tilde{\cF},\tilde{\omega})}$ and $\tilde{\eta'}=\eta_{b',x'}^{(\tilde{\cF},\tilde{\omega})}$.  

We first prove the existence of $\phi(b')\in \{(a_1b_1',a_2,b_2'),(a_1b_2',a_2b_1')\mid a_1,a_2\in \cF^\times \}$ such that \begin{equation}\label{e_rel}
t\eta'=\eta_{\phi(b'),x}^{(\cF,\omega)}. 
\end{equation}

 We have: \[\eta(\cF^2)\setminus\{0\}=\exp(\Lambda\cup (-x+\Lambda))=\exp(\ln(t))\exp(\Lambda\cup (-x'+\Lambda))\] and thus: \begin{equation}\label{e_Lambda_t}
\Lambda\cup(-x+\Lambda)=\ln(t)+\left(\Lambda\cup(-x'+\Lambda)\right).
\end{equation}

 Therefore: \[\ln(t) \in \Lambda\cup(-x+\Lambda).\]

If $\ln(t)\in \Lambda$, then \eqref{e_Lambda_t} becomes \begin{equation}\label{e_Lambda_t1}
\Lambda\cup (-x+\Lambda)=\Lambda\cup (-x'+\Lambda).
\end{equation}

If $x\in \Lambda$, this becomes $\Lambda\cup (-x'+\Lambda)=\Lambda$, thus $-x'\in \Lambda$ and hence $x-x'\in \Lambda$. If $x\notin \Lambda$, then $(-x+\Lambda)\cap \Lambda=\emptyset$, thus $-x+\Lambda\subset -x'+\Lambda$ and hence $x'-x\in \Lambda$.  Write $\ln(t)=\lambda$ and $x-x'=\omega(a)$, where $a\in \cF$. Set $\phi(b')=(\lambda b_1,\lambda ab_2)$.  Then by \eqref{e_relation_norms}, we have \eqref{e_rel}.

Assume now $\ln(t)\notin \Lambda$. Then $\ln(t)\in -x+\Lambda$ and \eqref{e_Lambda_t} becomes \begin{equation}\label{e_Lambda_t2}
\Lambda\cup (-x+\Lambda)=(-x+\Lambda)\cup (-(x+x')+\Lambda).
\end{equation}

If $-(x+x')\in -x+\Lambda$, then $-x'\in \Lambda$ and thus $-(x+x')+\Lambda=-x+\Lambda$. Therefore the right hand side of \eqref{e_Lambda_t2} becomes $-x+\Lambda$ and $\Lambda\subset (-x+\Lambda)$. This is impossible since $\ln(t)\in (x+\Lambda)\setminus \Lambda$ (which implies $(x+\Lambda)\cap \Lambda=\emptyset$). Therefore: \[-(x+x')\in \Lambda.\] Write $-x-x'=\omega(a)$ and $\ln(t)=-x+\lambda$, with $\lambda\in \Lambda$ and $a\in \cF$. Set $\phi(b')=(\lambda a b_2',\lambda b_1')$. Then by \eqref{e_relation_norms}, we have \eqref{e_rel} in this case.

We thus constructed $\phi(b')$ satisfying \eqref{e_rel} in any cases. 
Then by \eqref{e_relation_norms}, $t\eta_{b',x'}^{(\tilde{\cF},\tilde{\omega})}=\eta_{\phi(b'),x}^{(\tilde{\cF},\tilde{\omega})}$. Therefore up to replacing $(b',x')$ by $(\phi(b'),x)$, we can assume that $x=x'$ and $t=1$.

For $y\in \R$, denote by $B_y$, $B_y'$, $\tilde{B}_y$, $\tilde{B'_y}$ the  closed balls  of radius $e^{-y}$ and center $0$ with respect to $\eta$, $\eta'$, $\tilde{\eta}$ and $\tilde{\eta'}$ respectively. We have $B_y=\cF_{\geq y}b_1\oplus \cF_{\geq y-x}b_2$  As $\eta=\eta'$,  we have \[\cF_{\geq y}b_1\oplus \cF_{\geq y-x}b_2=\cF_{\geq y}b_1'\oplus \cF_{\geq y-x}b_2',\] for all $y\in \R$. 
Applying this equation with $y=0$ (resp. $y=x$), we get $b_1\in \cF_{\geq 0}b_1'\oplus \cF_{\geq -x}b_2'$ (resp. $b_2\in \cF_{\geq x}b_1'\oplus \cF_{\geq 0}b_2'$) and hence $\tilde{\cF}_{\geq y}b_1\subset \tilde{\cF}_{\geq y}b_1'\oplus \tilde{\cF}_{\geq y-x}b_2'$ (resp. $\tilde{\cF}_{\geq y-x}b_2\subset \tilde{\cF}_{y}b_1'\oplus \tilde{\cF}_{\geq y-x}b_2'$) for every $y\in \R$. Therefore $\tilde{B}_y\subset\tilde{B'_y}$ for every $y\in \R$. By symmetry, we deduce $\tilde{B}_y\subset\tilde{B'_y}$ for every $y\in \R$.

 Therefore for  $x\in \tilde{\cF}$, we have $-\ln(\tilde{\eta}(x))=\max\{y\in \R\mid x\in \tilde{B}_y\}=-\ln(\tilde{\eta'}(x))$ and hence $\tilde{\eta}=\tilde{\eta}'$.
\end{proof}

This lemma implies the following proposition.

\begin{Proposition}\label{p_emd_ext_nrm}
Let $(\tilde{\cF},\tilde{\omega})$ be an extension of $(\cF,\omega)$. Then the map $\I_{\sN}^{(\cF,\omega)}\rightarrow \I_{\sN}^{(\tilde{\cF},\tilde{\omega})}$ defined by $[\eta_{b,x}^{(\cF,\omega)}]\mapsto [\eta_{b,x}^{(\tilde{\cF},\tilde{\omega})}]$ for $b$ a basis of $\cF^2$ and $x\in \R$ is well-defined and is a $\mathrm{GL}_2(\cF)$-equivariant embedding . Its restriction to $\I^{(\cF,\omega)}_{\sL}$ is the map $[\cO b_1\oplus \cO b_2]\mapsto [\tilde{\cO} b_1\oplus \tilde{\cO} b_2]$ where $\tilde{\cO}$ is the valuation ring of $(\tilde{\cF},\tilde{\omega})$. 
\end{Proposition}

\begin{Remark}\label{r_fld_ext_vrt}
\begin{enumerate}
\item Let $a\in \I_\sN$. Then there exists an extension $(\tilde{\cF},\tilde{\omega})$ of $(\cF,\omega)$ such that $a$ is in the orbit of $0$, for the action of $\tilde{G}:=\mathrm{SL}_2(\tilde{\cF})$ on $\I_{\sN}^{(\tilde{\cF},\tilde{\omega})}$. Indeed, translating by  an element of $G$, we can assume that $a\in \A$. By \cite[Theorem 2.2.1]{engler2005valued}, there exists an extension $(\tilde{\cF},\tilde{\omega})$ of $(\cF,\omega)$ such that $a\in 2\tilde{\omega}(\tilde{\cF})\subset \R$. Let $x\in \tilde{\cF}$ be such that $2\omega(x)=a$. Then by Lemma~\ref{l_act_torus},  $a=\bt_x.0\in \tilde{G}.0$. This fact generalizes to Bruhat-Tits buildings, see \cite[Proposition 4.3]{remy2015BTbuildings}. We use it to determine the fixator of a point in $G$ (see Proposition~\ref{p_fix_pt}).

\item One also uses fields extensions in the theory of masures, in order to define the fixator of a point: for that purpose, it is convenient to consider a field extension $(\tilde{\cF},\tilde{\omega})$ of $(\cF,\omega)$ such that $\tilde{\omega}(\tilde{\cF})$ is dense in $\R$, see Section~\ref{s_parahoric}.
\end{enumerate}
\end{Remark}

\section{Fixators of intervals in $\A$}

If $\Omega$ is a non-empty subset of $\I_{\sN}$, we denote by $G_\Omega$\index[notation]{g@$G_\Omega$} the (pointwise) fixator of $\Omega$ in $G$: $G_\Omega=\{g\in G\mid g.x=x,\forall x\in \Omega\}$. If $x\in \I_{\sN}$, we simply write $G_x$ instead of $G_{\{x\}}$. 

\begin{Proposition}\label{p_fix_pt}
Let $x\in \A\simeq \R$. Then $G_x=\begin{psmallmatrix} \cO & \cF_{\geq -x}\\\cF_{\geq x} & \cO\end{psmallmatrix}\cap G$. In particular, the fixator of $0_\I$ is \[G_{0_\I}=\mathrm{SL}_2(\cO)=\{\begin{psmallmatrix} a& b\\ c& d\end{psmallmatrix}\in \mathrm{SL}_2(\cF)\mid a,b,c,d\in \cO\}.\]
\end{Proposition}

\begin{proof}
Set $\eta=\eta_{e,0}$. Let $g\in \mathrm{GL}_2(\cF)$ be such that $g.\eta=\eta$. Write $g^{-1}=\begin{psmallmatrix} a & b \\ c & d\end{psmallmatrix}$. Then $g.\eta (e_1)=\eta(e_1)=1=\eta(g^{-1}.e_1)=\max(|a|,|c|)$. Similarly, $g.\eta(e_2)=1=\max(|b|,|d|)$. Therefore $a,b,c,d\in \cO$. As $g^{-1}.\eta=\eta$, we deduce $g\in \mathrm{GL}_2(\cO)$. Conversely, if $g\in \mathrm{GL}_2(\cO)$, then $g.\eta\leq \eta$. As $g^{-1}\in \mathrm{GL}_2(\cO)$, we deduce $g.\eta=\eta$. Therefore \begin{equation}\label{e_fix_eta}
\{g\in \mathrm{GL}_2(\cF)\mid g.\eta=\eta\}=\mathrm{GL}_2(\cO).
\end{equation}

For $\mu \in \Lambda$, choose $x_\mu\in \cF$ such $\omega(x_\mu)=\mu$. Let now $g\in \mathrm{GL}_2(\cF)$ be such that $g.[\eta]=[\eta]$. Write $g.\eta=k\eta$, with $k\in \R_{>0}$. We have $\eta(\cF^2)\setminus\{0\}=\exp(\Lambda)=k\exp(\Lambda)$ and hence $k\in \exp(\Lambda)$.  Write $k=e^\lambda$, with $\lambda\in \Lambda$. Let $t=\begin{psmallmatrix} x_\lambda & 0 \\ 0 & 1\end{psmallmatrix}$ and $h=t^{-1}g$. Then we have $h.\eta=\eta$, according to Lemma~\ref{l_act_torus}. Therefore $g\in E\cdot \mathrm{GL}_2(\cO)$, where $E=\{\begin{psmallmatrix} x_\lambda & 0 \\ 0 & 1\end{psmallmatrix}\mid\lambda\in \Lambda\}$. By \eqref{e_fix_eta} and Lemma~\ref{l_act_torus}, we have the converse inclusion and thus \begin{equation}
\Fix_{\mathrm{GL}_2(\cF)}(\{0_\I\})=E\cdot  \mathrm{GL}_2(\cO).
\end{equation}

We deduce that $G_{0_\I}=\mathrm{SL}_2(\cO)$. 

Let $x\in \A$. By \cite[Theorem 2.2.1]{engler2005valued}, there exists an extension $(\tilde{\cF},\tilde{\omega})$ of $(\cF,\omega)$ such that $x\in 2\tilde{\omega}(\tilde{\cF})\subset \R$. Let $b\in \tilde{\cF}$ be such that $2\omega(b)=a$. Then by Lemma~\ref{l_act_torus},  $a=\bt_b.0$.   Let $\tilde{G}=\mathrm{SL}_2(\tilde{\cF})$. We have: \[\Fix_{\tilde{G}}(x)=t\Fix_{\tilde{G}}(\{0_\I\})t^{-1}=\begin{psmallmatrix} \tilde{\cO} & a^{-2}\tilde{\cO}\\ a^2 \tilde{\cO} & \tilde{\cO}\end{psmallmatrix}\cap \tilde{G}=\begin{psmallmatrix} \tilde{\cO} & \tilde{\cF}_{\geq -x}\\ \tilde{\cF}_{\geq x} & \tilde{\cO}\end{psmallmatrix} \cap \tilde{G},\] where $\tilde{\cO}$ is the valuation ring of $(\tilde{\cF},\tilde{\omega})$. Using Proposition~\ref{p_emd_ext_nrm}, we deduce that $G_x=G\cap \Fix_{\tilde{G}}(\{x\})$. Proposition follows.
\end{proof}

\begin{Corollary}\label{c_pt_fxd_g}
Let $g=\begin{psmallmatrix} a & b \\ c & d\end{psmallmatrix} \in G$. Then $g$ fixes a point of $\A$ if and only if $a,d\in \cO$ and $-\omega(b)\leq \omega(c)$. In this case $g$ fixes the segment $[-\omega(b),\omega(c)]\subset \A\simeq \R$ (or $]-\infty,\omega(c)]$ or $[-\omega(b),+\infty[$ if $b=0$ or $c=0$ respectively). 
\end{Corollary}

\begin{proof}
Let $x\in \R$.  Then $g$ fixes $x$ if and only if $g\in \Fix_G(\{x\})=\begin{psmallmatrix} \cO & \cF_{\geq -x}\\\cF_{\geq x} & \cO\end{psmallmatrix}\cap G$, if and only if $a,d\in \cO$ and  $\omega(b)\geq -x$ and $\omega (c)\geq x$. 
\end{proof}

\begin{Corollary}\label{c_fx_apt_end}
\begin{enumerate}
\item The fixator of $\A$ is $T_\cO=T\cap \mathrm{SL}_2(\cO)$.

\item Let $x\in \R$. Then $G_{[x,+\infty[}=\begin{psmallmatrix}
\cO & \cF_{\geq -x}\\
0 & \cO
\end{psmallmatrix}\cap G$ and $G_{]-\infty,x]}=\begin{psmallmatrix}
\cO & 0\\
\cF_{\geq x} & \cO
\end{psmallmatrix}\cap G$
\end{enumerate}
\end{Corollary}

\begin{Corollary}\label{c_Iwa_sbgp}
Let $x,y\in \R$ be such that $x\leq y$. Then $G_{[x,y]}=G_x\cap G_y=G\cap \begin{psmallmatrix} \cO & \cF_{\geq-x}\\ \cF_{\geq y} & \cO\end{psmallmatrix}$,

\end{Corollary}

\begin{proof}
We have $G_{[x,y]}\subset G_x\cap G_y$. By Proposition~\ref{p_fix_pt}, $G_x\cap G_y=G\cap \begin{psmallmatrix} \cO & \cF_{\geq-x}\\ \cF_{\geq y} & \cO\end{psmallmatrix}$ and by Proposition~\ref{p_fix_pt}, this subgroup fixes $[x,y]$.
\end{proof}

\section{Intersection of two apartments and Weyl groups}\label{ss_inter_apt}

We admit the following lemma (see \cite[Corollaire 3.3]{parreau2000immeubles}). 

\begin{Lemma}\label{l_parreau_3.3}
Let $b=(b_1,b_2),b'=(b_1',b_2')$ be two bases of $\cF^2$ and   $g\in  \mathrm{GL}_2(\cF)$ be such that $b'=g.b$. Write $g=(g_{i,j})_{i,j\in \llbracket 1,2\rrbracket}$.  Let $\fS_2$ be the group of permutations of $\llbracket 1,2\rrbracket$. Then there exists $\sigma\in  \fS_2$ such that for every ultrametric  norm adapted to both $b$ and $b'$, we have: \begin{enumerate}
\item $|\det(g)|=\frac{\eta(b_1)\eta(b_2)}{\eta(b_1')\eta(b_2')},$

\item $\eta(b'_j)=|g_{\sigma(j),j}|\eta(b_j)$, for both $j\in \llbracket1,2\rrbracket$. 
\end{enumerate} 
\end{Lemma}

The following proposition implies that $\I_{\sN}$ satisfies   axiom (A2) (see Theorem~\ref{t_building_SL2}), since $G$ acts transitively on the set of apartments:

\begin{Proposition}\label{p_A2}
\begin{enumerate}
\item Let $A$ be an apartment of $\IN$. Then there exists $h\in G$ such that $h.A=\A$ and $h$ fixes $\A\cap A$ pointwise. More precisely, let $g\in G$.  Then there exists  $n$ in $N$ such that $gn$ fixes $\A\cap g.\A$ pointwise.

\item Let $A$ be an apartment of $\IN$. Then $A\cap \A$ is enclosed i.e $A\cap \A\in \{ \emptyset, \A,]-\infty,b], [a,+\infty[, [a,b]\mid a,b\in \Lambda, a\leq b\}$.  

\item The stabilizer of $\A$ in $G$ is $N$.
\end{enumerate}

\end{Proposition}

\begin{proof}
1) We follow \cite[3.4]{parreau2000immeubles}. Let $g\in G$. If $\A\cap g.\A=\emptyset$, there is nothing to prove. We now assume $A\cap g.\A\neq \emptyset$.     Set $b=g^{-1}.e$.  By Lemma~\ref{l_parreau_3.3}, there exists $\sigma\in \mathfrak{S}_2$ such that for every ultrametric norm $\eta$ adapted to both $e$ and $b$, we have $\eta(b_j)=|g_{\sigma(j),j}|\eta(e_{\sigma(j)})$. 

For $j\in \llbracket 1,2\rrbracket$, set $\tilde{b_j}=b_{\sigma(j)}/g_{j,\sigma(j)}$. Let $\tilde{b}=(\tilde{b}_1,\tilde{b}_2)$.  Let $\tilde{n}\in \mathrm{GL}_2(\cF)$ be such that $\tilde{b}=\tilde{n}^{-1}.b$. Then $\tilde{n}$ is a monomial matrix. In particular, it stabilizes $\A$. Moreover, a norm is adapted to ($b$ and $e$) if and only if it is adapted to ($\tilde{b}$ and $e$). 

Let $\eta$ be a norm adapted to both $\tilde{b}$ and $e$. Let $j\in \{1,2\}$.  Then: \begin{equation}\label{e_n_adapt}
\eta(\tilde{b}_j)=\eta(b_{\sigma(j)})/|g_{j,\sigma(j)}|=|g_{\sigma^2(j),\sigma(j)}| \eta(e_{\sigma^2(j)})/g_{j,\sigma(j)}=\eta(e_j).
\end{equation}

Let $x\in  \A\cap g.\A$. Write $x=g\tilde{n}.[\eta_{e,t}]$, with $t\in \R$ and set $\eta=\eta_{e,t}$. Then $\eta$ is adapted to $e$. Write $x=\eta_{e,t'}$, with $t'\in \R$.  Then $[\eta]=\tilde{n}^{-1}g^{-1}.x=[\eta_{\tilde{n}^{-1}g^{-1}.e,t'}]=[\eta_{\tilde{b},t'}]$, by Lemma~\ref{l_act_G_norms}. Therefore $\eta$ is also adapted to $\tilde{b}$.

Let $v\in \cF^2$. Write $v=\lambda_1 e_1+\lambda_2 e_2$, with $\lambda_1,\lambda_2\in \cF$. As $\eta$ is adapted to both $\tilde{b}$ and $e$,  we have:  \begin{align*}g\tilde{n}.\eta(v)&=\eta(\lambda_1 \tilde{n}^{-1}g^{-1}.e_1+\lambda_2 \tilde{n}^{-1}g^{-1}. e_2)\\ &=\max(|\lambda_1|\eta(\tilde{b}_1),|\lambda_2|\eta(\tilde{b}_2))= \max(|\lambda_1|\eta(e_1),|\lambda_2|\eta(e_2))=\eta(v),
\end{align*} by  \eqref{e_n_adapt} and hence $g\tilde{n}$ fixes $\eta$. Therefore $g\tilde{n}$ fixes $\A\cap g\tilde{n}.\A=\A\cap A$.

 By Lemma~\ref{l_parreau_3.3}, we have $|\det(g\tilde{n})|=1=|\det(\tilde{n})|$. Let now $t=\begin{psmallmatrix} 1/\det(\tilde{n}) & 0 \\ 0 & 1\end{psmallmatrix}$ and $h=g\tilde{n}t$. Then $h\in G$ and by Lemma~\ref{l_act_torus}, $t$ fixes $\A$ pointwise so $h$ fixes $g.\A\cap \A$ pointwise. Therefore $h$ satisfies the condition of 1.

2) By 1., there exists $h\in G$ such that $h.\A=A$ and $h$ fixes $\A\cap A$. We then conclude with Corollary~\ref{c_pt_fxd_g}. 

3) By Subsection~\ref{ss_Action_N_A}, we already know that $N$ stabilizes $\A$.  Let $g\in \Stab_{G}(\A)$. Then by 1., there exists $n\in N$ such that $gn$ fixes $\A$. By Corollary~\ref{c_fx_apt_end}, $gn\in T_\cO$ and thus $g=gn.n^{-1}\in N$, which proves (3).
\end{proof}

Note that \begin{equation}\label{e_stab_A_GL_2}
\Stab_{\mathrm{GL}_2(\cF)}(\A)=N \cdot T_{\mathrm{GL}_2(\cF)}, \text{ where }T_{\mathrm{GL_2}(\cF)}=\begin{psmallmatrix}\cF^\times & 0 \\ 0 & \cF^\times\end{psmallmatrix}.
\end{equation}

Indeed let $g\in \mathrm{GL}_2(\cF)$ be such that $g.\A=\A$. Then if $t\in T_{\mathrm{GL}_2(\cF)}$ is such that $\det(t)=\det(g)$, we have $gt^{-1}.\A=\A$ and thus we conclude with Proposition~\ref{p_A2}.

Let $W^v=N/T$\index[notation]{w@$W^v$} be the \textbf{vectorial Weyl group} of $(G,T)$. Let $\tilde{s}=\begin{psmallmatrix} 0 & 1\\ -1 & 0\end{psmallmatrix}\in G$. The group $W^v=\{T, \tilde{s}T\}$ is isomorphic to $\Z/2\Z$. Let $\vec{\nu}:N\rightarrow \{\pm \Id\}$ be defined by $\vec{\nu}^{-1}(\{\Id\})=T$. 

The fixator of $\A$ is $T_\cO=\{\begin{psmallmatrix}
    x & 0 \\ 0 & x^{-1}
\end{psmallmatrix}\mid x\in \cO^*\}$\index[notation]{t@$T_{\cO}$}. The \textbf{affine Weyl group} is $\tilde{W}=N/T_\cO$\index[notation]{w@$\tilde{W}$}. This is the set of affine automorphisms of $\A$ induced by an element of $N$. By Subsection~\ref{ss_Action_N_A}, the map $\nu:N\rightarrow W^v\ltimes 2\Lambda$ defined by $\nu(n)=(\vec{\nu}(n),n.0)$ induces an isomorphism $\tilde{W}\simeq W^v\ltimes 2\Lambda$.

\begin{Corollary}\label{c_W_gp_int_apt}
\begin{enumerate}
\item Let $g\in G$. Then   there exist $n\in N$ and  $w\in  \tilde{W}$ such that for all $x\in \A\cap g^{-1}.\A$, we have $g.x=w.x$.

\item Let $g\in G$ be such that $g$ fixes at least two points of $\A$. Then $g$ fixes $\A\cap g^{-1}.\A$. 
\end{enumerate}
\end{Corollary}

\begin{proof}
1) By Proposition~\ref{p_A2} applied with $g^{-1}$, there exists $n\in N$ such that $g^{-1}n.x=x$ for all $x\in \A\cap g^{-1}.\A$. Then $g.x=n.x$ for all $x\in \A\cap g^{-1}.\A$. Let $w\in \tilde{W}$ be the map induced by $n$ on $\A$. Then  $w.x=n.x$ for all $x\in \A$, which proves 1).

2) Let $g\in G$ fixing at least two points. Let $h\in G$ be such that $h.\A=g.\A$ and $h$ fixes $\A\cap g.\A$. Then $h^{-1}g$ stabilizes $\A$ and thus it belongs to $N$. Let $w\in \tilde{W}$ be such that $h^{-1}g.x=w.x$, for all $x\in \A$. Then $w$ is an affine automorphism of $\A\simeq \R$ fixing two points: it fixes $\A$ and thus $g.x=h.x=x$, for all $x\in \A$. 
\end{proof}

\section{Filters and faces of $\I_{\sN}$}

\subsection{Motivation}

First assume that $\Lambda=\omega(\cF^\times)=\Z$. Let $x\in \A\simeq \R$. Its fixator $G_x$ is  $\begin{psmallmatrix} \cO & \cF_{\geq -x}\\\cF_{\geq x} & \cO\end{psmallmatrix}\cap G$. If $x\in \R\setminus \Z$, we have $G_x=\begin{psmallmatrix} \cO & \cF_{\geq \lceil -x\rceil }\\\cF_{\geq \lceil x\rceil } & \cO\end{psmallmatrix}\cap G=G_y$ for all $y\in ]\lfloor x\rfloor,\lfloor x\rfloor+1[$. 

For $x,y\in \IN$, write $x\sim y$ if $G_x=G_y$ and denote by $F_x$ the equivalence class of $x$.

\begin{Proposition}
Assume that $\omega(\cF^\times)=\Z$.

\begin{enumerate}
\item Let $x\in \A\simeq \R$. Then $F_x=]\lfloor x\rfloor,\lfloor x\rfloor+1[$ if $x\in \R\setminus \Z$ and $F_x=\{x\}$ if $x\in  \Z$.

\item  Let $x\in \IN$ and $g\in G$. Then $F_{g.x}=g.F_{x}$. In particular, the  $F_{x}$ such that $x\in \IN$ are the $g.F_x$ such that $g\in G$ and $x\in \A$. 
\end{enumerate}
\end{Proposition}

\begin{proof}
Let $x\in \A$. Then by the paragraph before the proposition, we have $F_x\cap \A=\begin{cases} &=\{x\}\\
&= ]\lfloor x\rfloor,\lfloor x\rfloor+1[\end{cases}$. Let $y\in F_x$. Then by Corollary~\ref{c_frndl_pr_tr} below, there exists an apartment $A$ containing $x$ and $y$. By Proposition~\ref{p_A2}, there exists $g\in G$ such that $g.A=\A$ and $g$ fixes $A\cap \A$. Then $G_{g.x}=G_x$ and $G_{g.y}=G_{g.x}=G_x$. Therefore $g.y\in F_x\cap \A$. By Proposition~\ref{p_A2} again, $A\cap \A\supset F_x\cap \A$ and thus $y=g^{-1}.g.y=g.y$. Consequently, $y\in \A$, hence $F_x=F_x\cap \A$ which proves (1). (2) is clear. 
\end{proof}

The $F_x$, $x\in \I_{\sN}$ are usually called the faces of $\I_{\sN}$.  As we will use a slightly different notion of faces, we will call them the simplicial faces here. We have $G_{]0,1[}=G_{]0,\epsilon[}=\begin{psmallmatrix} \cO & \cO\\ \fm & \cO\end{psmallmatrix}\cap G:=K_I$, for any $\epsilon\in \R_{>0}$, where $\fm=\{x\in \cO\mid \omega(x)>0\}$ is the maximal ideal of $\cO$.  As we shall see, we have the affine Bruhat decomposition $G=K_I\cdot  N\cdot  K_I$ (Proposition~\ref{p_Br_dec}). This is actually equivalent to the fact that for every pair of simplicial faces, there exists an apartment containing them, by Proposition~\ref{p_frnd_pr_T}.

Assume now that the valuation $\omega$ is dense, which means that $\Lambda$ is dense in $\R$. Then for every $\epsilon,\epsilon'\in \R_{>0}$ such that $\epsilon\neq \epsilon'$, we have: \[G_{]0,\epsilon'[}\neq G_{]0,\epsilon[}\subsetneq G_{]0,0^+[}:=\bigcup_{\epsilon''\in \R_{>0}} G_{]0,\epsilon''[}=\begin{psmallmatrix} \cO & \cO\\ \fm & \cO\end{psmallmatrix}\cap G.\]

As we shall see (see Lemma~\ref{l_fcs_dense}), for every $\epsilon\in \R_{>0}$, we have $G_{]0,\epsilon[}\cdot N\cdot G_{]0,\epsilon[}\neq G$.  However, we have $G=G_{]0,0^+[}\cdot N\cdot  G_{]0,0^+[}$. Equivalently, if $J_1,J_2$ are two intervals of $\R\simeq \A$ not reduced to single points, there exist $g_1,g_2\in G$ such that $g.J_1$ and $g.J_2$ are not contained in a common apartment. Therefore defining the faces as translates of intervals of $\A$ by elements of $G$ is not satisfactory in the non-discrete case.  In order to define faces uniformly (in discrete and dense cases), a solution is to use the notion of filters.

\subsection{Filters and faces of $\I_{\sN}$}\label{ss_filters}

 We now give a definition of faces which works for any valuation $\omega$ with values in $\R$. For this we use the notion of filters. For example, the filter $]0,0^+[$, consisting of all the subsets of $\I_{\sN}$ containing $]0,\epsilon[\subset \A$, for some $\epsilon\in \R_{>0}$, will be a face.

\subsubsection{Filters}

\begin{Definition}
A filter  on a set $\mathcal{E}$ is a nonempty set $\cV$ of nonempty subsets of $\mathcal{E}$ such that, for all subsets $\Omega$, $\Omega'$ of $\mathcal{E}$, we have:\begin{itemize}
\item $\Omega$, $\Omega'\in \cV$ implies $\Omega\cap \Omega'\in \cV$ 

\item $\Omega'\subset \Omega$ and $\Omega'\in \cV$ implies $\Omega\in \cV$.

\end{itemize}

\end{Definition}

If $\mathcal{E}$ is a set and $\cV,\cV'$ are filters on $\mathcal{E}$, we define $\cV\cup \cV'$ to be the filter $\{\Omega\cup \Omega'\mid (\Omega,\Omega')\in \cV\times \cV'\}$.

\begin{Definition}[filter generated by a prefilter]

Let $\cE$ be a non-empty set. Let $\sB$ be a set of subsets of $\cE$. We say that $\sB$ is a prefilter if $\emptyset\notin \sB$ and for all $\Omega_1,\Omega_2\in \sB$, there exists $\Omega\in \sB$ such that $\Omega\subset \Omega_1\cap \Omega_2$.

If $\sB$ is a prefilter, the \textbf{filter generated by $\sB$}  is the filter on $\cE$ consisting of the subsets of $\cE$ containing an element of $\sB$. 

\end{Definition}

Note that if $\cV$ is a filter on a set $\cE$, then it is a prefilter on $\cE$. Therefore if $\cE'$ is a set containing $\cE$, we can consider the filter  on $\cE'$ generated by $\cV$. By abuse of language, we often identify these two filters, when there is no ambiguity. For example, we often define filters on $\A$ and then consider them as filters on $\IN$.

Let $\mathcal{E},\mathcal{E}'$ be sets, $f:\mathcal{E}\rightarrow\mathcal{E}'$ be a map and $\cV$ be a filter on $\mathcal{E}$. Then $f(\cV)$ is the filter on $\cE'$ generated by the prefilter $\{f(\Omega)\mid \Omega\in \cV\}$. 

If $\cV$ is a filter on a set $\mathcal{E}$, and $\Omega$ is a subset of $\mathcal{E}$, we say  that $\cV$ contains $\Omega$ if every element of $\cV$ contains $\Omega$. We denote it $\cV\supset \Omega$. If $\Omega$ is nonempty, the \textbf{principal filter} on $\cE$ associated with $\Omega$ is the filter  $\sF_{\Omega,\cE}$ of subsets of $\cE$ containing $\Omega$.

 A filter $\cV$ is said to be contained in another filter $\cV'$: $\cV\subset  \cV'$ (resp. in a subset $Z$ in $\mathcal{E}$: $\cV\subset Z$) if every set in $\cV'$  is in $\cV$ (resp. if $Z\in \cV$). 
 
 These definitions of containment are inspired by the following facts. Let $\cE$ be a set, $\cV$ be a filter on $\cE$ and $\Omega,\Omega'\subset \cE$. Then :\begin{itemize}
\item $\Omega\subset \Omega'$ if and only if $\sF_{\Omega,\cE}\subset \sF_{\Omega',\cE}$,

\item $\Omega\subset \cV$ if and only if $\sF_{\Omega,\cE}\subset \cV$,

\item $\Omega\supset \cV$ if and only if $\sF_{\Omega,\cE}\supset \cV$.  
\end{itemize}

If $\Omega$ is a non-empty subset of a non-empty topological space $\cE$,  and $x\in \overline{\Omega}$, then the \textbf{germ  of $\Omega$ at $x$}\index{germ} is the filter $\germ_x(\Omega)$\index[notation]{g@$\germ_x(\Omega)$} on $\cE$  composed with the subsets of $\cE$ containing a neighborhood of $x$ in $\Omega$. If $\cV$ is a filter on $\cE$, then $\overline{\cV}$ is the filter on $\cE$ generated by the prefilter $\{\overline{\Omega}\mid \Omega\in \cV\}$. If $\cV$, $\cV'$ are filters on $\cE$, we say that $\cV$ \textbf{dominates}\index{d@domination} $\cV'$ if $\overline{\cV}\supset \cV'$.

Let $\cE$ be a set and let $f:\cE\rightarrow \cE$. Let $\cV$ be a filter on $\cE$. We say that $f$ \textbf{fixes} $\cV$ if it fixes pointwise at least one element of $\cV$. Assume now that $f$ is a  bijection.   We say that $f$ \textbf{stabilizes $\cV$} if for every subset  $\Omega$ of $\cE$, we have $f(\Omega)\in \cV$ if and only if $\Omega\in \cV$. In particular, if $H$ is a group acting on $\cE$, we can talk of the fixator and the stabilizer of the filter $\cV$\index{stabilizer of a filter}. These are subgroups of $H$.

\subsubsection{Faces of $\A$}

We set $\alpha=\Id:\A\simeq \R\rightarrow \R$ and $\Phi=\{\pm \alpha\}$. 

The \textbf{fundamental vectorial chamber} is $C^v_f=\{x\in \A\mid \alpha(x)>0\}=\R_{>0}$. The subfaces of $\overline{C^v_f}$ are $\{0\}$ and $C^v_f$ itself.  The \textbf{vectorial faces} (resp. \textbf{vectorial chambers}) of $\A$ are the transforms the subfaces of $\overline{C^v_f}$ by $W^v$. Explicitly, these are $\R_{<0},\R_{>0}$  and $\{0\}$ (resp. $\R_{>0}$ and $\R_{<0}$). 

A \textbf{sector} $\fq$ of $\A$ is a set of the form $a+C^v$, for some vectorial chamber $C^v$. Its \textbf{germ at infinity}, denoted $\fQ=\germ_\infty(\fq)$ is the filter on $\A$ consisting of the subsets of $\A$ containing $b+C^v$, for some $b\in \A$. As $\A$ is one-dimensional, its sector-germs at infinity are usually called the \textbf{ends} of $\A$. For $\epsilon\in \{-,+\}$,   we denote by $\fQ_{\epsilon\infty}$ the germ at infinity of the sector $\epsilon C^v_f$. 

A \textbf{(closed) half-apartment} (resp. \textbf{open half-apartment}, \textbf{wall}) of $\A$ is a set of the form $D(\beta,k):=\{x\in \A\mid \beta(x)+k\geq 0\}$\index[notation]{d@$D(\beta,k)$} (resp. $\mathring{D}(\beta,k)=\{x\in \A\mid \beta(x)+k>0\}$, $M(\alpha,k)=\{k\}$ or $M(-\alpha,k)=-k$) where $\beta\in \{\alpha,-\alpha\}$ and $k\in \Lambda$. We set $D(\beta,+\infty)=\A=\mathring{D}(\beta,+\infty)$. 

If $\Omega$ is a filter on $\A$ (in particular it can be a non-empty set), its \textbf{enclosure} $\cl(\Omega)$\index[notation]{c@$\cl$} is the filter on $\A$ consisting of the subsets $E$ of $\A$ containing an intersection of the form $D(\alpha,k_\alpha)\cap D(-\alpha,k_{-\alpha})$ containing $\Omega$, with $k_\alpha,k_{-\alpha}\in \Lambda\cup\{+\infty\}$.  A subset of $\A$ is called \textbf{enclosed}\index{enclosed} if it is the intersection of finitely many (possibly zero) half-apartments of $\A$.

Let $x\in \A\simeq \R$. We denote by $]x,x^+[$\index[notation]{z@$]x,x^+[,]x^-,x[$}  (resp. $]x^-,x[$) the filter on $\A$ composed with the subsets of $\A$ containing a neighbourhood of $x$ in $]x,+\infty[$ (resp. $]-\infty,x[$). An \textbf{alcove} of $\A$ is a filter of $\A$ of the form $]x^-,x[$ or $]x,x^+[$, for some $x\in \R$. A \textbf{(local) face} of $\A$ is an alcove or a singleton of $\A$.

We extend to $\IN$ all the notions of $\A$ which are stable under $\tilde{W}$. For example the faces, sectors, sectors-germs of $\IN$ are the filters of the form $g.\cV$, where $\cV$ is a face, a sector, a sector-germ of $\A$ respectively.

\section{Bruhat decompositions of $G$ and consequences}

In this section, we expose the vectorial and the Bruhat decompositions of $G$ and derive building theoretic consequences. 

If $\cV$ is a filter on $\A$, we denote by $G_{\cV}$\index[notation]{g@$G_\cV$} its fixator in $G$. We have $G_{\cV}=\bigcup_{\Omega\in \cV} G_{\Omega}$.

\subsection{Decompositions of $G$ and existence of apartments}

We start this section by  relating certain decompositions of $G$ to the existence of apartments containing pairs of filters.

\begin{Proposition}\label{p_frnd_pr_T}
Let $\cV_1,\cV_2$ be two filters on $\A$. Then the following are equivalent: \begin{enumerate}
\item $G=G_{\cV_1}\cdot N \cdot  G_{\cV_2}$,

\item for all $\cV_1'\in G.\cV_1$, $\cV_2'\in G.\cV_2$, there exists an apartment containing $\cV_1'\cup \cV_2'$. 
\end{enumerate}
\end{Proposition}

\begin{proof}
(2) $\Rightarrow$ (1): assume (2). Let $g\in G$ and consider $\cV_1'=\cV_1$ and $\cV_2'=g.\cV_2$. Let $A$ be an apartment containing $\cV_1$ and $\cV_2'$. By Proposition~\ref{p_A2}, there exists $h\in G$ such that $h.A=\A$ and $h$ fixes $\A\cap A$. We have $h\in G_{\cV_1}$. By Corollary~\ref{c_W_gp_int_apt}, there exists $n\in N$ such that $n.x=hg.x$ for all $x\in \A\cap (hg)^{-1}\A$.  We have $hg.\cV_2\subset \A$ and thus $n^{-1}hg\in G_{\cV_2}$.  Consequently, $g=h^{-1} n (n^{-1}hg)\in G_{\cV_1}\cdot N\cdot G_{\cV_2}$, which proves (1). 

(1) $\Rightarrow$ (2):  we follow \cite[Théorème 4.7.18]{bruhat1972groupes}. Assume (1). Let $\cV_1'\in G.\cV_1$ and $\cV_2'\in G.\cV_2'$. As $G$ acts transitively on the set of apartments, we can assume $\cV_1'=\cV_1$. Let $g\in G$ be such that $g^{-1}.\cV_2'\subset \A$. Write $g=b_1nb_2$, with $b_1\in G_{\cV_1}$, $b_2\in  G_{g^{-1}.\cV_2'}$ and $n\in N$.  Let $A=b_1.\A$. Then $A\supset b_1.\cV_1=\cV_1$ and $\cV_2'=gg^{-1}.\cV_2 '=b_1nb_2.g^{-1}.\cV_2'=b_1n g^{-1}.\cV_2'\subset b_1N.\cV_2'\subset b_1.\A$, which proves (2).  
\end{proof}

\subsection{Vectorial Bruhat decomposition}

Let $\alpha:T\rightarrow \cF^*$ be defined by $\alpha(\begin{psmallmatrix}
    x & 0\\
    0 & x^{-1}
\end{psmallmatrix})=x$, for $x\in \cF^*$. The group of algebraic characters of $T$ is $\Z \alpha$, with the additive notation. For $a\in \cF$, we set $x_\alpha(a)=\begin{psmallmatrix}
    1 & a\\ 0 & 1
\end{psmallmatrix}$ and $x_{-\alpha}(a)=\begin{psmallmatrix}
    1 & 0 \\ -a & 1
\end{psmallmatrix}$. Let $U=U^+=x_{\alpha}(\cF)$ and $U^-=x_{-\alpha}(\cF)$\index[notation]{u@$U^{\pm}$}. These are the   \textbf{positive and negative standard unipotent subgroups} of $G$. Then $x_{\pm\alpha}$ is an isomorphism between $(\cF,+)$ and $U^{\pm}$. The root system of $(G,T)$ is $\{\alpha,-\alpha\}$, where $\alpha:\A\simeq \R\rightarrow \R$ is the identity, with the identification of \eqref{e_identification_A_R}. 

By Corollary~\ref{c_pt_fxd_g} and Corollary~\ref{c_W_gp_int_apt}, we have for all $a\in \cF$, $\beta\in \{\alpha,-\alpha\}$:  \begin{equation}\label{e_action_unipotent}
x_{\beta}(a).\A\cap \A=\{x\in \A\mid x_\beta(a).x=x\}=D_{\beta,\omega(a)}:=\{x\in \A\mid \beta(x)+\omega(a)\geq 0\}.
\end{equation}

Set \[\alpha^\vee(a)=\begin{psmallmatrix} a & 0 \\ 0 & a^{-1}\end{psmallmatrix}, (-\alpha)^\vee(a)=\begin{psmallmatrix} a^{-1} & 0 \\ 0 & a\end{psmallmatrix}\text{ and }\tilde{s}=\begin{psmallmatrix} 0 & 1\\ -1 & 0\end{psmallmatrix}.\]

By an easy calculation, we have: \begin{equation}\label{e_commutaion_relation}
x_{-\beta}(a^{-1})x_{\beta}(a)x_{-\beta}(a^{-1})=\beta^\vee(a)\tilde{s},
\end{equation} for $a\in  \cF^\times$ and $\beta\in \{-\alpha,\alpha\}$.

For $\beta\in \{-\alpha,\alpha\}$ and $\Omega\subset \A$, we set $U_{\beta,\Omega}=\{u\in U_\beta\mid u.z=z,\forall z\in \Omega\}$\index[notation]{u@$U_{\beta,\Omega}, U_{\beta,\cV}$}. For $\cV$ a filter on $\A$, we set $U_{\beta,\cV}=\bigcup_{\Omega\in \cV} U_{\beta,\Omega}$. If $x\in \R$, we write $U_{\beta,x}$ instead of $U_{\beta,\{x\}}$.

Set   \[B=B^{+}=\begin{psmallmatrix}
\cF^\times  & \cF\\
0 & \cF^\times
\end{psmallmatrix}\cap G \text{ and } B^-=\begin{psmallmatrix} \cF^\times & 0\\ \cF & \cF^\times\end{psmallmatrix}\cap G .\]\index[notation]{b@$B$, $B^{\pm}$}These are the \textbf{positive and negative standard Borel subgroups} of $G$. We have $T\cdot U^\pm=B^\pm$. 

\begin{Proposition}\label{p_fix_infty}
Let $\epsilon\in \{-,+\}$. The stabilizer of $\fQ_{\epsilon \infty}$ is $B^\epsilon$ and the fixator of $\fQ_{\epsilon \infty}$ is $T_\cO\cdot  U^\epsilon$, where $T_\cO=\begin{psmallmatrix} \cO & 0\\ 0 & \cO \end{psmallmatrix}\cap G.$

\end{Proposition}

\begin{proof}
 By Lemma~\ref{l_act_torus} and \eqref{e_action_unipotent}, we have $T_\cO \cdot U^\epsilon\subset \Fix_G(\epsilon \infty)$ and $B^\epsilon\subset \Stab_{G}(\fQ_{\epsilon \infty})$. By Corollary~\ref{c_fx_apt_end}, $\Fix_G(\fQ_{\epsilon \infty})\subset T_\cO\cdot  U^\epsilon$, which proves the second equality. 
 
 Take $g\in \Stab_{G}(\fQ_{\epsilon \infty})$. By definition, there exist $x,x'\in \R$  such that $g.[x,\epsilon\infty[\subset [x',\epsilon\infty[$. By Corollary~\ref{c_W_gp_int_apt}, there exists $w\in \tilde{W}$ such that for all $y\in [x,\epsilon\infty[$, we have $g.y=w.y$. Then $w$ is necessarily  a translation. Let $t\in T$ be such that $t.y=w.y$, for all $y\in [x,\epsilon\infty[$. Then $gt^{-1}$ fixes $[x,+\infty[$ and hence $g\in T_\cO\cdot  U^\epsilon \cdot T=B^\epsilon$. 
\end{proof}

\begin{Proposition}\label{p_Bruhat} (Bruhat decomposition)
We have $G=B\cdot N\cdot B$.
\end{Proposition}

\begin{proof}
This is an application of Gauss elimination. Let $g=\begin{psmallmatrix} a & b \\ c & d\end{psmallmatrix}$. If $c=0$, then $g\in B$ so we assume $c\neq 0$. Then $x_\alpha(-a/c)gx_{-\alpha}(-d/c)=\begin{psmallmatrix} 0 & -1/c\\ c & 0\end{psmallmatrix}\in N$ and thus $g\in B\cdot N\cdot B$. 
\end{proof}

\begin{figure}\label{f_bruhat_vec}
\centering
\includegraphics[scale=0.3]{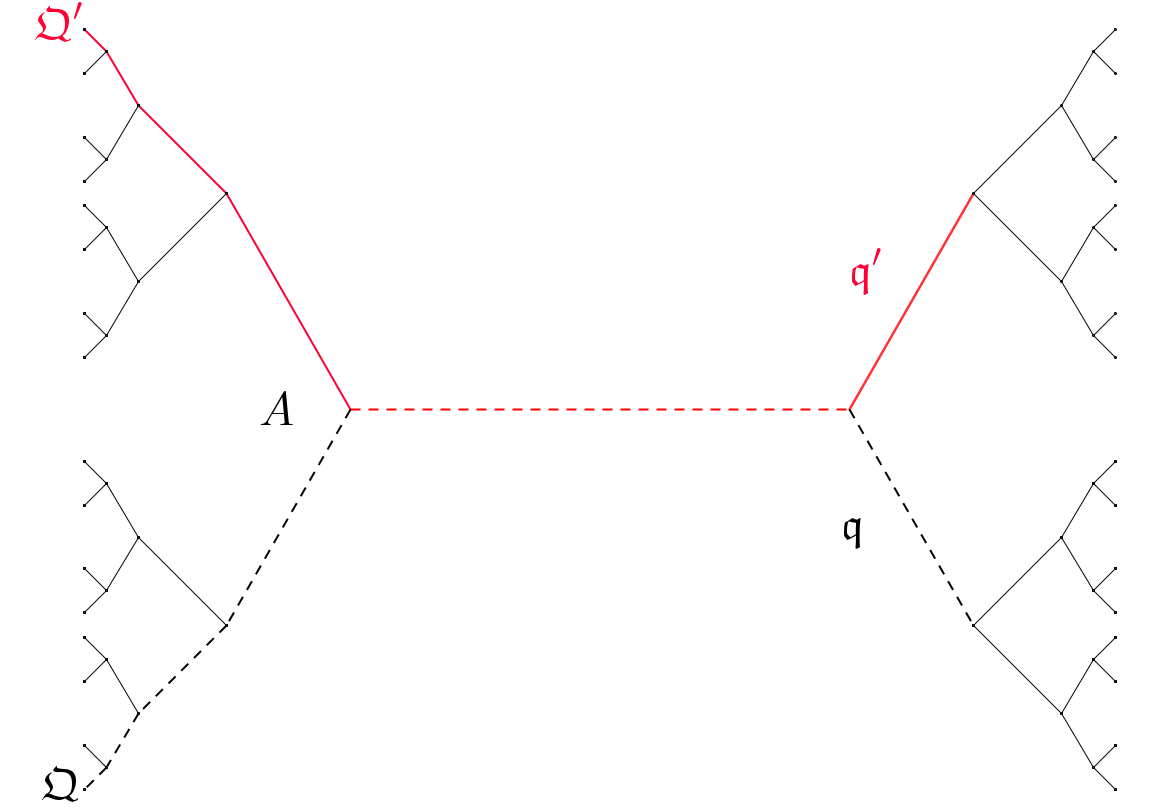}
\caption{The  vectorial Bruhat decomposition (or axiom (A4)) for the tree: $\fq$ (in dotted line) and $\fq'$ (in red) are sectors of $\IN$. No apartment contains $\fq\cup \fq'$, but the apartment $A$ (in continuous red and dotted black) contains the associated sector-germs $\fQ$ and $\fQ'$.}\label{f_A4}
\end{figure}

\begin{Corollary}\label{c_ends_frndly}(see Figure~\ref{f_A4})
Let $\fQ_1,\fQ_2$ be two sector-germs of $\I_{\sN}$. Then there exists an apartment containing $\fQ_1\cup \fQ_2$.  
\end{Corollary}

\begin{proof}
Since $G$ acts transitively on the set of sector-germs of $\I_{\sN}$, we can assume that $\fQ_1=\fQ_{+\infty}$. Let $g\in G$ be such that $g^{-1}.\fQ_2\subset \A$. Using Proposition~\ref{p_Bruhat}, we write $g=bnb'$, with $b,b'\in B$ and $n\in N$. By Proposition~\ref{p_fix_infty}, we have $b.\A\supset b.\fQ_{+\infty}=\fQ_{+\infty}=\fQ_1$ and $b.\A=bn.\A\supset bn.\fQ_{+\infty}=bnb'.\fQ_{+\infty}=\fQ_2$, which proves the corollary.
\end{proof}

\subsection{Affine Bruhat decomposition}

Recall that if $x\in \A\simeq \R$, then  $]x,x^+[$ (resp. $]x^-,x[$) denotes the filter on $\A$ composed with the subsets of $\A$ containing a neighbourhood of $x$ in $]x,+\infty[$ (resp. $]-\infty,x[$). We set $[x,x^+[=]x,x^+[\cup \{x\}$ and $]x^-,x]=\{x\}\cup ]x^-,x[$. By Corollary~\ref{c_Iwa_sbgp}, we have: \begin{equation}
\forall x\in \A, G_{]x,x^+[}=G\cap \begin{psmallmatrix}\label{e_Iwahori_sugbroup}
\cO & \cF_{\geq -x}\\
\cF_{>x} & \cO
\end{psmallmatrix}\text{ and }G_{]x^-,x[}=G\cap \begin{psmallmatrix}
\cO & \cF_{> -x}\\
\cF_{\geq x} & \cO
\end{psmallmatrix}. 
\end{equation}

In particular if $x\in \A\setminus \Lambda$, we have $G_{x}=G_{]x^-,x^+[}$, where $]x^-,x^+[=]x^-,x]\cup [x,x^+[$ and if $x\in \Lambda$, then $G_{[x,x^+[},G_{]x^-,x]}\subsetneq G_{x}$. 

The fixator of $[0,0^+[$ in $G$ is the \textbf{Iwahori subgroup} $K_I=G\cap \begin{psmallmatrix} \cO & \cO\\ \fm & \cO\end{psmallmatrix}$, where $\fm=\{x\in \cO\mid \omega(x)>0\}$\index[notation]{m@$\fm$} is the maximal ideal of $\cO$. When $\Lambda=\Z$ (which is equivalent to $\Lambda$ discrete, up to renormalizing $\omega$), then $K_I$ is the fixator of $[0,1]$, since $\fm=\cF_{\geq 1}$. However when $\Lambda$ is dense in $\R$,  $K_I$ does not fix any set of the form $[0,x[$, for any $x\in \R_{>0}$. This motivates the use of filters rather than sets in the  definition of faces in this case.

Using the affine Bruhat decomposition for $\mathrm{GL}_2(\cF)$ (\cite[Proposition 3.10]{parreau2000immeubles}), we get the affine Bruhat decomposition:

\begin{Proposition}\label{p_Br_dec}(affine Bruhat decomposition)
Let $x,y\in \A$. Then we have $G=G_{]x,x^+[}\cdot N\cdot G_{]y,y^+[}$.
\end{Proposition}

The following Corollary is axiom (A3') (see Theorem~\ref{t_building_SL2}).

\begin{figure}[h]
\centering
\includegraphics[scale=0.3]{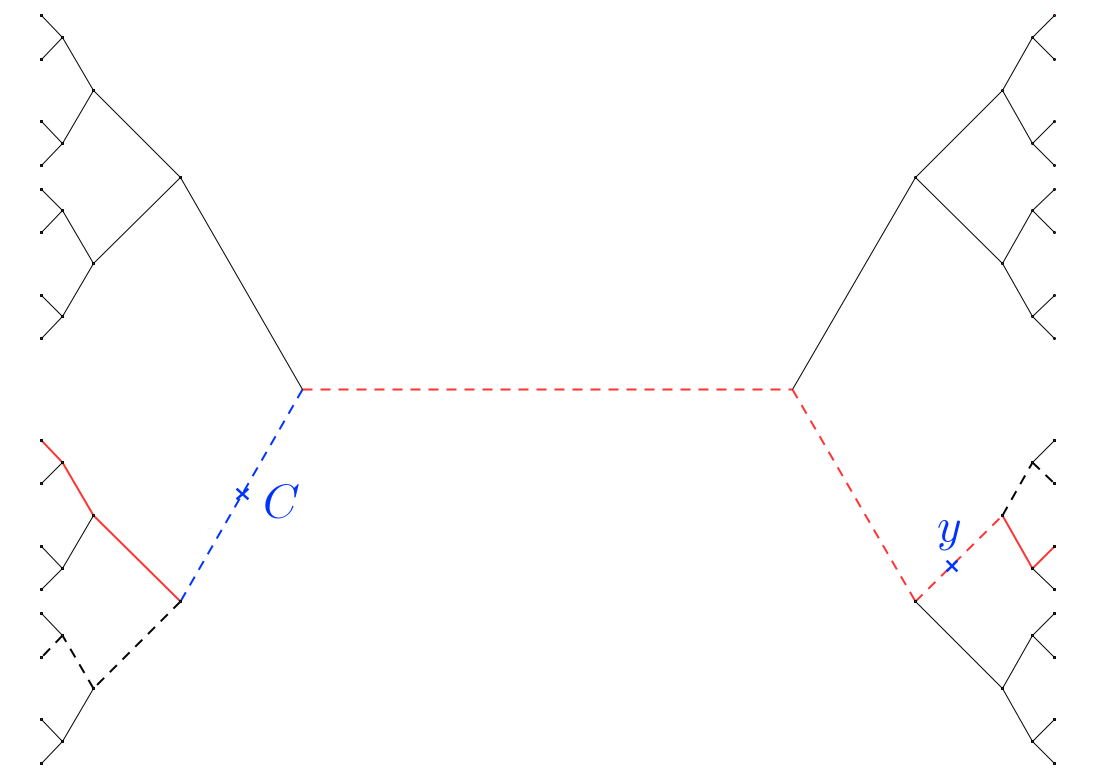}
\caption{The affine Bruhat decomposition (or axiom (A3')) for the tree: $C$ (in dotted blue) and $y$ (in blue) are an alcove and a point of $\IN$ respectively. The apartment in dotted line and the apartment consisting of ($C$ union the red part of $\IN$) are two apartments containing $x\cup C$.}\label{f_bruhat_aff}
\end{figure}

\begin{Corollary}\label{c_frndl_pr_tr}(see Figure~\ref{f_bruhat_aff})
Let $F_1,F_2$ be two faces of $\I_{\sN}$. Then there exists an apartment containing $\overline{F_1}$ and $\overline{F_2}$. 
\end{Corollary}
\begin{proof}
Let $C_1,C_2$ be two alcoves of $\I_{\sN}$ such that $F_1\subset \overline{C_1}$ and $F_2\subset \overline{C_2}$. Let $i\in \{1,2\}$. By definition, there exist $g_1,g_2\in G$ such that $g_i.C_i\subset \A$. Replacing $g_i.C_i$ by $-g_i.C_i=\tilde{s}g_i.C_i$  if necessary, we can assume that $g_i.C_i=]x_i,x_i^+[$, for some $x_i\in \R$. By Proposition~\ref{p_Br_dec} and Proposition~\ref{p_frnd_pr_T}, there exists an apartment $A$ containing $g_1^{-1}.g_1.C_1=C_1$ and $g_2^{-1}.g_2.C_2=C_2$. By Proposition~\ref{p_A2}, $A$ contains $\overline{C_1}$ and $\overline{C_2}$, which proves the lemma.
 \end{proof}

\begin{Lemma}\label{l_fcs_dense}
Assume that $\omega$ is a dense valuation. Let $\epsilon_1,\epsilon_2\in \R_{>0}$. Then $G\neq G_{]0,\epsilon_1[}\cdot N\cdot G_{]0,\epsilon_2[}$. 
\end{Lemma}

\begin{figure}[h]
\centering
\includegraphics[scale=0.3]{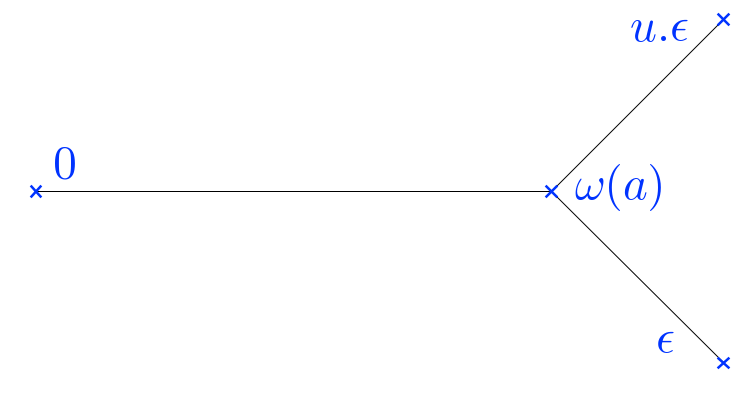}
\caption{No apartment can contain $0$, $u.\epsilon$ and $\epsilon$ (note that apartments branch almost everywhere but are not represented).}\label{f_non_dec}
\end{figure}

\begin{proof}

By Proposition~\ref{p_frnd_pr_T}, it suffices to find $\Omega_1'\in G.]0,\epsilon_1[$ and $\Omega_2'\in G.]0,\epsilon_2[$ such that $\Omega_1'\cup \Omega_2'$ is not contained in any apartment.

Let $\epsilon=\min(\epsilon_1,\epsilon_2)$.  Let $a\in \cF$ be such that $\omega(a)\in ]0,\epsilon[$. Set $\Omega_1'=]0,\epsilon_1[$ and $\Omega'_2=u.]0,\epsilon_2[$, where $u=\begin{psmallmatrix} 1 & 0\\ a & 1\end{psmallmatrix}$. By Corollary~\ref{c_fx_apt_end} and \eqref{e_action_unipotent}, we have \[u.\A\cap \A=]-\infty,\omega(a)]=\{x\in \A\mid u.x=x\}.\]

We have $\omega(a)\in ]0,\epsilon[\cap u.]0,\epsilon[$. But $\epsilon\notin u.[0,\epsilon]$ and $u.\epsilon\notin [0,\epsilon]$. Therefore, no apartment contains $0,\epsilon$ and $u.\epsilon$ (see Figure~\ref{f_non_dec}) and thus no apartment contains $\Omega_1'\cup  \Omega_2'$, which proves the lemma.
\end{proof}

\section{Building axiomatic satisfied by $\I_\sN$}

In this section, we introduce certain axioms of affine $\R$-buildings and prove that they are satisfied by $\I_{\sN}$. They are similar to the axioms of masures, but an important difference  is axiom (A3'), which implies that two points are always contained in a common apartment, and which is false in general for masures.

\subsection{Axioms (A2) to (A4)}\label{ss_axioms}

An \textbf{apartment of type $(\A,\tilde{W})$} is a set $A$ equipped with a non-empty set $\mathrm{Isom}(\A,A)$ of bijections from $\A$ to $A$   such that if $f_0\in \mathrm{Isom}(\A,A)$, then $\mathrm{Isom}(\A,A)=\{f_0\circ w\mid w\in \tilde{W}\}$. In an apartment of type $(\A,\tilde{W})$, we can define any notion defined in $\A$ and preserved by $\tilde{W}$. For example, we can define faces, sectors, walls, etc.

Recall that an apartment of $\I_{\sN}$ is a set of the form $A_{b}=\{[\eta_{b,t}]\mid t\in \R\}$. Let $A$ be an apartment of $\I_{\sN}$. We define $\mathrm{Isom}(\A,A)$ as the set of $h|_{\A}^{A}$ such that $h\in G$ is such that $h.\A=A$. As $G$ acts transitively on the set of apartments (by Lemma~\ref{l_act_G_norms}), $\mathrm{Isom}(\A,A)$ is non-empty. As $N$ is the stabilizer of $\A$ and as the set of automorphisms of $\A$ induced $N$ is $\tilde{W}$ (see Subsection~\ref{ss_inter_apt}), this equips $A$ with the structure of an apartment of type $(\A,\tilde{W})$. We have the following theorem:

\begin{Theorem}\label{t_building_SL2}
The set $\I_\sN$ is an affine building of type $(\A,\tilde{W})$ in the sense of \cite[Définition 1.1]{parreau2000immeubles}. In particular, it satisfies the following properties. \begin{enumerate}
\item[\textbf{(A2)}] For every two apartments $A$ and $B$, $A\cap B$ is enclosed and there exists $g\in G$ such that $g.A=B$ and $g$ fixes $A\cap B$.

\item[\textbf{(A3')}] For every two faces $F_1,F_2$ of $\I_{\sN}$, there exists an apartment containing $F_1\cup F_2$. This applies in particular when $F_1$ and $F_2$ are points.

\item[\textbf{(A4)}] For every two sectors $\fq,\fq'$ of $\I_{\sN}$, there exists an apartment containing their germs at infinity $\fQ$ and $\fQ'$. 

\end{enumerate}

\end{Theorem}

\begin{proof}
By \cite[Théorème 1.21]{parreau2000immeubles}, the definition given here is equivalent to \cite[Définition 1.1]{parreau2000immeubles}. Then (A2) is implied by Proposition~\ref{p_A2}, using the transitivity of the action of $G$ on the set of apartments, (A3') is Corollary~\ref{c_frndl_pr_tr} and (A4) is Corollary~\ref{c_ends_frndly}.
\end{proof}

Let $x,y\in \I_\sN$. Then by (A3'), there exists an apartment $A$ containing $x$ and $y$. We define $[x,y]$ as the segment joining $x$ to $y$ in $A$. By (A2), this segment does not depend on the choice of $A$ and is contained in any apartment containing $x$ and $y$. If $x,y\in \I_{\sN}$, we define: \begin{equation}\label{e_definition_fs}
\fs_{x,y}:[0,1]\rightarrow [x,y], t\mapsto (1-t)x+ty.
\end{equation}\index[notation]{s@$\fs_{[x,y]}$}

Let $g\in G$ be such that $g.A=\A$. Then we set $d(x,y)=|g.x-g.y|$\index[notation]{d@$d(x,y)$}, where $\A$ is identified with $\R$. This is well-defined, independently of the choices of $A$ and of $g$. Indeed, take an apartment $A'$ and $g'\in G$ such that $g'.A'=\A$. Let $h\in G$ be such that $h.A=A'$ and $h$ fixes $A\cap A'$, which exists by (A2). Then $n:=gh^{-1}g'^{-1}$ stabilizes $\A$ so it belongs to $N$. Let $w\in \tilde{W}$ be the element induced by $n$ on $\A$. For $z\in A\cap A'$, we have $g.z=ng'h.z=ng'.z=w.g'.z$ and thus $|g'.y-g'.x|=|w^{-1}(g.y-g.x)|=|g.y-g.x|$, since $\tilde{W}$ acts by isometry on $\A$. This proves that $d$ is well-defined. 
 
 By definition, $d$ is invariant under the action of $G$. Actually, since $\Stab_{\mathrm{GL}_2(\cF)}(\A)=N\cdot T_{\mathrm{GL}_2(\cF)}$ acts by isometry on $\A$, $d$ is also  $\mathrm{GL}_2(\cF)$-invariant. It is obviously symmetric and its restriction to every apartment is isometric to $(\R,|\cdot|)$. It also  satisfies the triangle inequality and thus it is a distance, but this requires a bit more work, see Corollary~\ref{c_d_distance}.

Let $A$ and $A'$ be two apartments such that $A\cap A'$ has non-empty interior. By (A2), there exists an apartment isomorphism $f:A\rightarrow A'$. If $f':A\rightarrow A'$ is an other apartment isomorphism fixing $A\cap A'$, then $f'^{-1}\circ f:A\rightarrow A$ is an apartment isomorphism fixing $A\cap A'$. In particular, it is an affine automorphism of $A$ fixing a subset with non-empty interior and thus it is the identity. Therefore $f=f'$: there exists a unique apartment isomorphism $A\rightarrow A'$ fixing $A\cap A'$.

\subsection{Exchange condition and sundial configuration}

In this subsection, we introduce the exchange condition and the sundial configuration. These two properties will be very useful in our study of masures. 

\begin{figure}[h]
\centering
\includegraphics[scale=0.3]{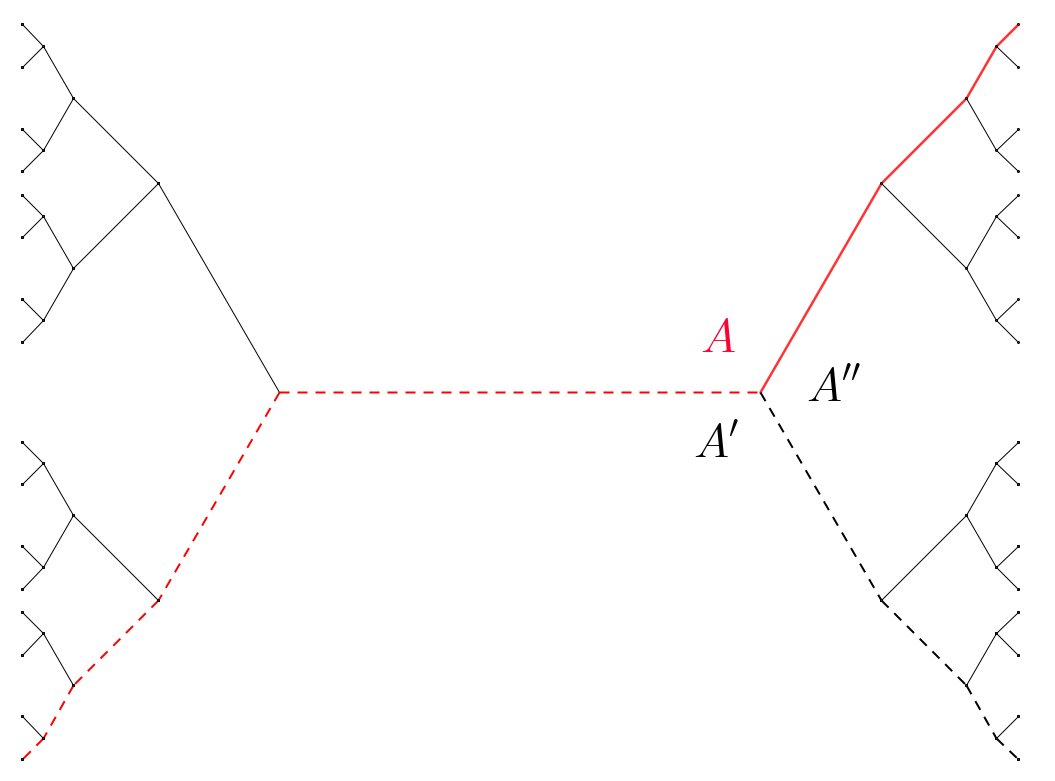}
\caption{The exchange condition for the tree. The apartment $A$ is in red, $A'$ is in dotted red and dotted black and $A''$ is in continuous red and dotted black}\label{f_EC}
\end{figure}

The following proposition is called ``\textbf{exchange condition}'' in \cite[2]{bennett2014axiomatic} (see Figure~\ref{f_EC}).

\begin{Proposition}\label{p_EC_T}
Let $A$ and $A'$ be two distinct apartments of $\I_\sN$ such that $A\cap A'$ is a half-apartment. Let  $D_1=\overline{A'\setminus A}$ and $D_2=\overline{A\setminus A'}$. Then $A'':=D_1\cup D_2$ is an apartment of $\I_{\sN}$ and $A\cap A'\cap A''=D_1\cap D_2$ is point, which is moreover a vertex of $\I_{\sN}$. 
\end{Proposition}

\begin{proof}
Using apartment isomorphisms, we can assume that $A=\A$ and by symmetry, we can assume that $A'$ contains $\fQ_{-\infty}$. Let $g\in G$ be such that $g.\A=A'$ and $g$ fixes $\A\cap A'$, which exists by (A2). By Proposition~\ref{p_fix_infty}, we can write $g=x_{-\alpha}(a)t$, with $a\in \cF$ and $t\in T_\cO$. Then we have $A'=x_{-\alpha}(a).\A$. By \eqref{e_action_unipotent}, we have $D_1=x_{-\alpha}(a).D_+$, where $D_+=\{y\in \A\mid  y\geq \omega(a)\}=D_2$. Let $D_-=\{y\in \A\mid  y\leq \omega(a)\}$. By \eqref{e_commutaion_relation} and Lemma~\ref{l_act_torus}, we have: \[x_{-\alpha}(a).D_+=x_\alpha(-a^{-1})(-\alpha)^\vee(a)\tilde{s} x_{\alpha}(-a^{-1}).D_+=x_\alpha(-a^{-1})(-\alpha)^\vee(a)\tilde{s}.D_+=x_{\alpha}(-a^{-1}).D_-.\] Therefore $D_1\cup D_2=x_{\alpha}(-a^{-1}).\A$. We set $A''=x_\alpha(-a^{-1}.\A)$. Then $A''$ is an apartment and $\A\cap A'\cap A''=\{\omega(a)\}$, which proves the result.
\end{proof}

We deduce the ``Sundial configuration property'', which is the first point of the next proposition (see \cite[2]{bennett2014axiomatic}).

\begin{Proposition}\label{p_SC_tr}
\begin{enumerate}
\item Let $C$ be an alcove or sector-germ (at infinity) of $\I_{\sN}$ and $A$ be an apartment. Then we can write $A=D_1\cup D_2$ where $D_1,D_2$ are two half-apartments of $\I_{\sN}$ such that $D_1\cap D_2$ is a point and such that for both $i\in \{1,2\}$, there exists an apartment $A_i$ containing $C$ and $D_i$. 

\item In particular, if $F$ is a face and $\fQ$ is a sector-germ of $\IN$, then there exists an apartment containing $F$ and $\fQ$.
\end{enumerate} 
\end{Proposition}

\begin{proof}
1) We first assume that $\fQ:=C$ is a sector-germ of $\IN$. Using an element of $G$, we may assume that $A=\A$. If $\fQ$ is contained in $\A$, there is nothing to prove so we assume that $\fQ$ is not contained in $\A$. Then by (A4), there exists an apartment $A_1$ containing $\fQ$ and $\fQ_{-\infty}$.  By (A2), $D_1:=A_1\cap \A$ is  a half-apartment of $\A$. Let $D'=\overline{A_1\setminus \A}$. Then $D'$ is a half-apartment of $A_1$ and thus by Proposition~\ref{p_EC_T}, if $D_2=\overline{\A\setminus D_1}$, then $A_2:=D'\cup D_2$ is an apartment of $\I_{\sN}$. By assumption, $\fQ\subset D'$ and thus $D_1$ and $D_2$ satisfy the condition of  1). This proves 1) when $C=\fQ$ is a sector-germ of $\IN$. 

2) Let $F$ be a face of $\IN$ and $\fQ$ be a sector-germ of $\IN$. Let $A$ be an apartment containing $F$. Then by 1) applied with $C=\fQ$, $F$ is contained in $A_1$ or $A_2$, with the notation of 1) This proves 2). Now to prove 1) when C is an alcove, one can replace (A4) by 2) the proof of 1).
\end{proof}

Note that we will prove that if $C$ is not contained in $A$, then $\{D_1,D_2\}$ is unique, see Remark~\ref{r_SC_tr}.

Using Proposition~\ref{p_frnd_pr_T} and Proposition~\ref{p_fix_infty}, Proposition~\ref{p_SC_tr} (2) implies the Iwasawa decomposition (see Corollary~\ref{c_Iwa_T}). Using retractions, we will describe this decomposition more precisely.

\subsection{Retractions and Hecke paths}

Retractions are an important tool in the study of buildings and masures. We introduce here the retractions centered at infinity and deduce the Iwasawa decomposition of $G$. We then study the image of a  line segment of $\I_{\sN}$ by a retraction.

\begin{Proposition}\label{p_retractions}
Let $C$ be an alcove or a sector-germ of $\I_{\sN}$ and $A$ be an apartment containing $C$. Then there exists a unique retraction $\rho=\rho_{C,A}:\I_{\sN}\rightarrow A$ such that for each apartment $A'$ containing $C$, $\rho|_{A'}^{A}$ is the apartment isomorphism from $A'$ to $A$ fixing $A\cap A'$. 
\end{Proposition}

 \begin{proof}
 Let $x\in \I_{\sN}$. Then by (A3') or Proposition~\ref{p_SC_tr}, there exists an apartment $A_x$ containing $C$ and $x$. Let $f:A_x\rightarrow A$ be the apartment isomorphism fixing $A\cap A_x$. Then if $\rho$ exists, we necessarily have $\rho(x)=f(x)$. Now to ensure the existence of $\rho$ it suffices to prove that $f(x)$ is independent of the choice of $A_x$. Let $A'$ be an apartment containing $x$ and $C$ and let $f':A'\rightarrow A$ be the apartment isomorphism fixing $A\cap A'$. Let $\phi:A'\rightarrow A_x$ be the apartment isomorphism fixing $A'\cap A_x$. Then $f',f\circ \phi:A'\rightarrow A$ are apartment isomorphism fixing $C$, and thus they fix a subset of $A'$ with non-empty interior. Therefore $f\circ \phi=f'$. In particular, $f'(x)=f\circ\phi(x)=f(x)$, which proves that $\rho$ is well-defined.
 \end{proof}

\begin{figure}[h]
\centering
\includegraphics[scale=0.3]{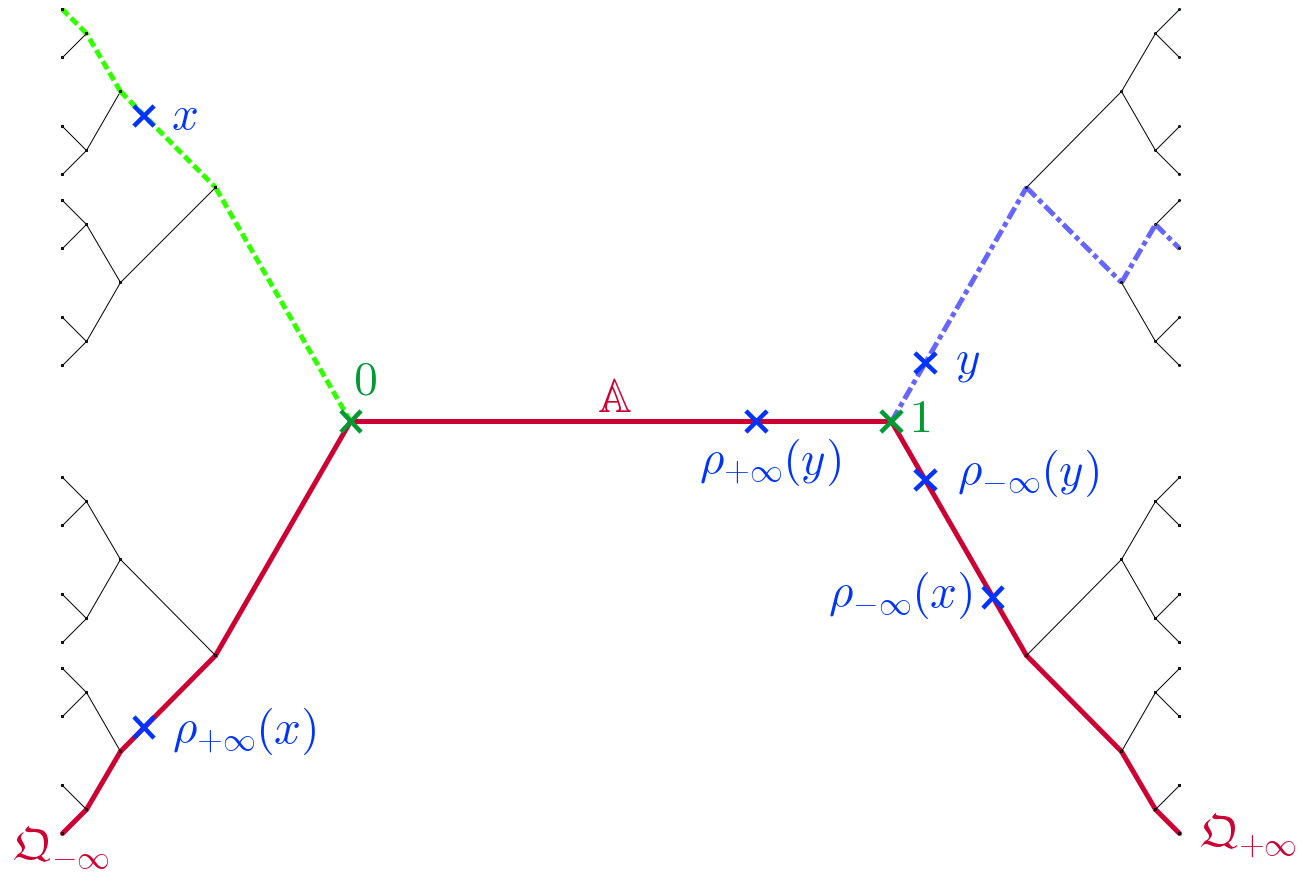}
\caption{The standard apartment is the red one.  The distance is distorted on the picture: the distance (defined in Subsection~\ref{ss_distance}) between two consecutive vertices is always $1$. The point $x$ is contained in the apartment consisting  of ($]\fQ_{-\infty},0]$ (in red) union the green doted half-apartment), and $y$ is contained in the apartment consisting of ($[1,\fQ_{+\infty}[$ union the purple doted half-apartment). The point $x$ is at distance $7/4$ from $0$ and thus $\rho_{+\infty}(x)$ and $\rho_{-\infty}(x)$ are at distance $7/4$ of $0$. The point $y$ is at distance $1/4$ from $1$ and thus $\rho_{+\infty}(y)$ and $\rho_{-\infty}(y)$ are at distance $1/4$ from $1$. Note that $d(\rho_{-\infty}(x),\rho_{-\infty}(y))=1/2<d(x,y)=3$. }\label{f_retractions}
\end{figure}

 \begin{Notation}
Recall that $\fQ_{+\infty}$  (resp. $\fQ_{-\infty}$) is the sector-germ of $\A$ consisting of the subsets of $\A$ containing a neighborhood of $+\infty$   (resp. $-\infty$). We denote by $\rho_{+\infty}$\index[notation]{r@$\rho_{+\infty}$, $\rho_{-\infty}$} and $\rho_{-\infty}$ the retractions $\rho_{\fQ_{+\infty,\A}}$ and $\rho_{\fQ_{-\infty},\A}$, respectively, with the notation of Proposition~\ref{p_retractions}. These retractions are illustrated in figure~\ref{f_retractions}.
\end{Notation} 
 
 \begin{Remark}\label{r_char_ret}
 Let $C$ be an alcove or sector-germ of $\IN$ and $x\in \IN$. Let $A$ be an apartment containing  $C$ Then $\rho_{A,C}$ is the unique $G_C$-invariant map $\I\rightarrow A$ whose restriction to $A$ is the identity. For $x\in \IN$, $\rho_{A,C}(x)$ is the unique point of $G_C.x\cap A$. 
 \end{Remark}

\begin{Corollary}\label{c_Iwa_T}(Iwasawa decomposition)
For each $\bw\in \tilde{W}$, choose $n_\bw\in N$ such that $n_\bw$ acts by $\bw$ on $\A$. Let $\epsilon\in \{-,+\}$.  Then we have $G=\bigsqcup_{\bw\in \tilde{W}} U^{\epsilon} n_{\bw} G_{[0,0^+[}$. More precisely, if $g\in G$, then  $g\in U^{\epsilon} n_{\bw} G_{]0,0^+[}$, where $\bw$ is the unique element of $\tilde{W}$  such that $\rho_{\epsilon\infty}(g.]0,0^+[)=n_{\bw}.]0,0^+[$. 
\end{Corollary}

\begin{proof}

Let $g\in G$. By Remark~\ref{r_char_ret} and Proposition~\ref{p_fix_infty}, we can write $\rho_{\epsilon \infty}(g.]0,0^+[)=tug.]0,0^+[=ug.]0,0^+[$, for some $u\in U^\epsilon$ and $t\in T_\cO$. By Corollary~\ref{c_W_gp_int_apt}, there exists $\bw\in \tilde{W}$ such that $ug.]0,0^+[=n_{\bw}.]0,0^+[$. 
Then $g\in u^{-1} n_\bw G_{]0,0^+[}$, which proves that
 $G=\bigcup_{\bw\in \tilde{W}}U^\epsilon n_{\bw} G_{]0,0^+[}$. Now if $\bv,\bw\in \tilde{W}$ are such that $n_{\bv}\in U^{\epsilon}n_\bw G_{]0,0^+[}$. Then $\rho_{\epsilon\infty}(n_\bv.]0,0^+[)=\rho_{+\infty}(n_{\bw}.]0,0^+[)=\bv.]0,0^+[=\bw.]0,0^+[$, which proves that $\bv=\bw$ and completes the proof of the corollary.
\end{proof}

If $C$ is an alcove  of $\A$, we say that $C$ is positive if it is of the form $]x,x^+[$ for some $x\in \A$ and negative otherwise. 

We now  study the image of a segment-line by a retraction. We prove that this image, which is called a Hecke path, satisfies very restrictive properties (for example they can fold only once). Hecke paths are crucial in our study of masures: we use it to prove that the intersection of two apartments is convex. These paths are also used in many applications of the theory of masures to the study of Kac--Moody groups.

\begin{Proposition}\label{p_Hecke_paths}
Let $C$ be an alcove or a sector-germ of $\A$ and $\rho=\rho_{\A,C}:\I\rightarrow \A$ be the retraction defined in Proposition~\ref{p_retractions}. Let $x,y\in \I_{\sN}$ and $\pi=\rho\circ \fs_{x,y}:[0,1]\rightarrow \A$. Then $\pi$ is a Hecke path of shape $d(x,y)$ with respect to $C$ in the sense that:\begin{enumerate}
\item $\pi$ is a piecewise affine path,

\item for all $t\in [0,1]$ for which it makes sense (i.e $t\neq 0$ for $\pi'_-(t)$ and $t\neq 1$ for $\pi'_+(t)$), we have $\pi'_{+}(t),\pi'_-(t)\in \{\pm d(x,y)\}$,

\item $\{t\in ]0,1[\mid \pi'_+(t)\neq \pi'_-(t)\}$ contains at most one point $t_0$. When $t_0$ exists, we have $\pi(t_0)\in \A_\sL=\I_\sL\cap \A$ and $\pi'_+(t_0)$ is of the same sign as $C$.
\end{enumerate}
\end{Proposition}

\begin{proof}
Let $A$ be an apartment containing $[x,y]$. If $A$ contains $C$, then if $f:A\rightarrow \A$ is the apartment isomorphism fixing $C$, we have $\pi=f\circ \fs_{x,y}$. Thus $\pi$ is a segment  and the result is clear. We now assume that $A$ does not contain $C$. Using Proposition~\ref{p_SC_tr}, we write $A=D_1\cup D_2$, where $D_1$ and $D_2$ are two half-apartments such that $D_1\cap D_2$ is a point and $D_i$ is contained in an apartment $A_i$ containing $C$ for both $i\in \{1,2\}$. If $[x,y]$ is contained in $D_i$, for $i\in \{1,2\}$, then we can replace $A$ by $A_i$ and we get the result, by the previous case. We now assume that $[x,y]$  is neither contained in $D_1$ nor in $D_2$.

 For $i\in \{1,2\}$, we denote by $f_i:A_i\rightarrow \A$ the apartment isomorphism fixing $A_i\cap \A$.  Let $\phi_i:A_i\rightarrow A$ be the apartment isomorphism fixing $A_i\cap A$ and $\psi_i=f_i\circ \phi_i^{-1}:A\rightarrow \A$. Let $a\in A$ be the point of $D_1\cap D_2$. Then $\psi_1(a)=\psi_2(a)=\rho(a)$. Moreover, $\psi_1(D_1)$ and $\psi_2(D_2)$ are two half-apartments of $\A$ delimited by $\rho(a)$ and not containing $C$ and thus $\psi_1(D_1)=\psi_2(D_2)$. Therefore $\psi_1(D_2)\neq \psi_2(D_2)$ and $\psi_1\neq \psi_2$. Therefore $r:=\psi_{1}\circ \psi_2^{-1}:\A\rightarrow \A$ is an element of $\tilde{W}$ fixing $\rho(a)$ and different from the identity: it is the reflection with respect to $\rho(a)$. Now for $z\in D_1$, we have $\rho(z)=f_1(z)=\psi_1(z)$ and for $z\in D_2$, we have $\rho(z)=f_2(z)=\psi_2(z)$. Let $t_0$ be such that $\fs_{x,y}(t_0)=a$. Then $\pi|_{[0,t_0]}$ is a segment from $\rho(x)$ to $\rho(a)$ and $\pi|_{[t_0,1]}$ is a segment from $\rho(a)$ to $\rho(y)$. As these two segments are contained in $\rho(D_1)=\rho(D_2)$ and as this half-apartment does not contain $C$, we deduce the sign condition. 
\end{proof}

\section{Study of the distance $d$ on $\I_{\sN}$}\label{ss_distance}
We now study the ``distance'' $d$ on $\IN$ defined in Subsection~\ref{ss_axioms}. We prove that it is a distance and  introduce an equivalent definition of $d$. We then use it to describe the $G$-orbits of $\I_{\sL}$. 

\subsection{Definition of the distance}
 We define a map $d:\I_{\sN}\times \I_{\sN}\rightarrow \R_{\geq 0}$\index[notation]{d@$d$} as follows. Let $x,y\in \I_{\sN}$. By Corollary~\ref{c_frndl_pr_tr}, there exists an apartment $A$ containing $x$ and $y$. Let $g\in G$ be such that $g.A=\A\simeq \R$. We set $d(x,y)=|g^{-1}.x-g^{-1}.y|$. This does not depend on the choices of $A$  and $g$.

\begin{Proposition}
If $C$ is an alcove or a sector-germ of $\A$ and $A$ is an apartment containing $C$, then $\rho=\rho_{A,C}:(\I_\sN,d)\rightarrow (A,d)$ is $1$-Lipschitz continuous. In particular, $\I_\sN$ satisfies axiom (A5')  (see \cite{parreau2000immeubles} or \cite[2.4.4.1 Definition]{rousseau2023euclidean}): for all apartment $A$ and for all $x\in \IN$, there exists a retraction $\rho:\IN\rightarrow A$ such that $\rho^{-1}(\{x\})=\{x\}$, which is $1$-Lipschitz continuous in restriction to any apartment. 
\end{Proposition}

\begin{proof}
Using an element of $G$, we may assume that $A=\A$. Then the fact that $\rho$ is $1$-Lipschitz continuous follows from Proposition~\ref{p_Hecke_paths}. Let us prove axiom (A5'). Let $x\in \I_{\sN}$ and $A$ be an apartment containing $x$. Let $C$ be an alcove of $A$ dominating $x$ and $\rho=\rho_{A,C}$. Let $y\in \rho^{-1}(\{x\})$ and $A'$ be an apartment containing $C$ and $y$. Then  the restriction $f$ of $\rho$ to $A'$ is the apartment isomorphism fixing $A\cap A'$. Therefore $f(y)=f(x)=x$ and thus $x=y$, which proves that $\rho^{-1}(\{x\})=\{x\}$. We deduce that (A5') is satisfied.

\end{proof}

 \begin{Proposition}\label{p_pr_apt}
 Let $A$ be an apartment of $\I_\sN$ and $z\in \I_\sN$. Then: \begin{enumerate}
 \item there exists a unique $\pr_A(z)\in A$ such that $d(z,A)=d(z,\pr_A(z))$.
 
 \item For all $y\in A$, we have $[y,z]=[y,\pr_A(z)]\sqcup ]\pr_A(z),z]$ and $d(y,z)=d(y,\pr_A(z))+d(\pr_A(z),z)$. 
 
 \item If $z\notin A$, then $\pr_A(z)$ is a vertex of $\I_{\sN}$.
 
 \item Assume that $z\notin A$ and that $D_+,D_-$ are two half-apartments of $\I_{\sN}$ such that $D_+\cup D_-=A$ and such that  $\{z\}\cup D_i$ is contained in an apartment of $\I_{\sN}$, for both $i\in  \{-,+\}$. Then $D_-\cap D_+=\{\pr_{A}(z)\}$. 
\end{enumerate}   
 \end{Proposition}
 
In figure~\ref{f_retractions}, we have $\pr_{\A}(z)=0$, for all $z$ in the dotted green half-apartment and $\pr_{\A}(z)=1$, for all $z$ in the dotted purple half-apartment.

 \begin{proof}
If $z\in A$, the result is clear, with $z=\pr_A(z)$. Assume now  $z\notin A$.  Using apartment isomorphisms, we may assume that $A=\A$. By Proposition~\ref{p_SC_tr}, there exist two half-apartments $D_-$ and $D_+$ of $\A$ such that $D_-\cap D_+$ is non-empty and such that for both $\epsilon\in \{-,+\}$, there exists an apartment $A_\epsilon$ containing $D_\epsilon$ and $z$. Take $a\in D_-\cap D_+$. Up to renumbering, we can assume that $D_\epsilon\supset \fQ_{\epsilon\infty}$, for $\epsilon\in \{-,+\}$. Replacing $D_+$ by $[a,+\infty[$ if necessary, we can assume that $D_-\cap D_+=\{a\}$.  Let $D=\overline{A_-\setminus \A}$.  Using the exchange property (Proposition~\ref{p_EC_T}), we can assume that $A_+=D\cup D_+$. 
  
Let $y\in \A$. Let us prove: \begin{equation}\label{e_proj_tr_apt}
a\in [y,z].
\end{equation}  Let $\epsilon\in \{-,+\}$. Consider the apartment isomorphism $f_\epsilon:A_{\epsilon}\rightarrow \A$ fixing $A_{\epsilon}\cap  \A$. Then $f_\epsilon(D)$ is the half-apartment of $\A$ opposite to $D_\epsilon$ and thus $f_\epsilon(D)=D_{-\epsilon}$. Let $y\in \A$. Set $\epsilon=-$ if  $y\in ]-\infty,a]$ and  $\epsilon=+$ if $y\in ]a,+\infty[$. Then $f_{\epsilon}([y,z])=[f_{\epsilon}(y),f_{\epsilon}(z)]=[y,f_{\epsilon}(z)]$. As $y\in D_\epsilon$ and $f_{\epsilon}(z)\in D_{-\epsilon}$, we have $a\in [y,f_{\epsilon}(z)]$, thus $f_{\epsilon}([y,z])\supset [y,a]$ and hence $[y,z]\supset [y,a]=f_{\epsilon}^{-1}([y,a])$. This proves \eqref{e_proj_tr_apt}. Therefore we have $[y,z]=[y,a]\sqcup ]a,z]$. Moreover, $f_\epsilon([y,z])=f_\epsilon([y,a])\sqcup f_{\epsilon}(]a,z])=[y,a]\sqcup f_{\epsilon}(]a,z])$ and thus $d(y,z)=d(y,a)+d(a,z)$. 

We deduce $d(z,y)\geq d(z,a)$, for all $y\in \A$, with equality only if $a=y$. This proves (1) and (2), with $\pr_A(z)=a$.

 As $a$ was any point of $D_+\cap D_-$, we deduce that we have $D_+\cap D_-=\{\pr_A(z)\}$, for  any pair $D_-,D_+$ of half-apartments satisfying the assumption of (4), which proves (4). 

We have $\pr_A(z)=a\in A_-\cap A_+\cap \A$. As $z\notin \A$, we deduce that $a$ is a vertex of $\I_{\sN}$, by Proposition~\ref{p_EC_T}.
 \end{proof}
 
 \begin{Remark}\label{r_SC_tr}
 Let $A$ be an apartment of $\IN$ and $z\in \IN$. Assume $z\notin A$ and  set $a=\pr_A(z)$. Denote by $D_1$ and $D_2$ the two half-apartments of $A$ delimited by $a$. Then Propositions~\ref{p_pr_apt} and Proposition~\ref{p_SC_tr} actually prove that $D_i\cup \{z\}$ is contained in an apartment of $\IN$ for both $i\in \{1,2\}$, and that if $D$ is a half-apartment of $A$ such that $D\cup \{z\}$ is contained in an apartment, then $D\subset D_1$ or $D\subset D_2$. 
 \end{Remark}

 \begin{Corollary}\label{c_d_distance}
 The map $d:\I_\sN\rightarrow \I_{\sN}$ is  a distance.
 \end{Corollary}
 
\begin{proof}
It remains to prove that $d$ satisfies the triangle inequality. Let $x,y,z\in \I_{\sN}$. Using apartment isomorphisms, we can assume that $x,y\in \A\simeq \R$ and that $x\leq y$, when we identify $\A$ with $\R$. Let $a=\pr_\A(z)$.  There are three cases: $a\in ]-\infty,x]$, $a\in ]x,y[$ and $a\in [y,+\infty[$. 

If $a\in ]-\infty,x]$, then $d(x,y)+d(y,z)=d(x,y)+d(x,y)+d(a,x)+d(z,a)=d(x,z)+2d(x,y)\geq d(x,z)$.

If $a\in ]x,y[$, then $d(x,z)+d(z,y)=d(x,a)+d(a,z)+d(a,y)+d(a,z)=d(x,z)+2d(a,y)\geq d(x,z)$.

If $a\in [y,+\infty[$, then $d(x,z)+d(z,y)=d(x,y)+d(y,a)+d(a,z)=d(x,z)$. In the three cases, we have the triangle inequality, which proves that $d$ is a distance.
\end{proof}

 \begin{Corollary}\label{c_geodesic}
 For every $x,y\in \I_{\sN}$, we have \[[x,y]=\{z\in \I_{\sN}\mid d(x,z)+d(z,y)=d(x,y)\}.\]
 \end{Corollary}
 
 \begin{proof}
 Let $x,y\in \I_{\sN}$ and $A$ be an apartment containing them. Then $A$ contains $[x,y]$ and  by definition of $d$, we have $d(x,z)+d(z,y)=d(x,y)$, for all $z\in [x,y]$. Let now $z\in \I_{\sN}$ and $a=\pr_A(z)$, with the notation of Proposition~\ref{p_pr_apt}. Then $d(x,z)+d(y,z)=d(x,a)+d(a,y)+2d(a,z)\geq d(x,z)+2d(a,z)$. Therefore $d(x,z)+d(y,z)=d(x,y)$ implies that $a=z$, which implies $z\in A$ and $z\in [x,y]$, which proves the other inclusion.
 \end{proof}
 
\subsection{Interpretation of $d$ in terms of norms and $G$-orbits of $\I_{\sL}$}

 \begin{Proposition}\label{p_d_infty}
 Let $x,y\in \I_{\sN}$. By (A3'), there exists a basis $b$ of $\cF^2$ and $t_x,t_y\in \R$ such that $x=[\eta_{b,t_x}]$ and $y=[\eta_{b,t_y}]$. We set $d_\infty(x,y)=|t_x-t_y|$. Then $d_\infty:\I_{\sN}\times \I_{\sN}\rightarrow \R_{\geq 0}$ is well-defined, independently of the choices we made. Moreover, we have $d_\infty=d$. 
\end{Proposition}

\begin{proof}
Let $b'$ be a basis of $\cF^2$ such that $x,y\in A_{b'}$. Write $x=[\eta_{b',t_x'}]$ and $y=[\eta_{b',t_y'}]$. Let $g\in \mathrm{GL}_2(\cF)$ be such that $g.b=b'$. Let $\eta_x=\eta_{b,t_x}$ and $\eta_y=\eta_{b,t_y}$. Then $\eta_x$ and $\eta_y$ are adapted to both $b$ and $b'$. Thus by Lemma~\ref{l_parreau_3.3}, there exists $\sigma\in \fS_2$ such that $\frac{\eta_x(b_2')}{\eta_x(b_2)}=|g_{\sigma(2),2}|=\frac{\eta_y(b_2')}{\eta_y(b_2)}$. We have $\eta_x(b_2')=e^{-t_x'}$, $\eta_x(b_2)=e^{-t_x}$, $\eta_y(b_2')=e^{-t_x'}$ and $\eta_y(b_2)=e^{-t_y}$. Therefore $e^{t_y-t_x}=e^{t_y'-t_x'}$ and thus $t_y-t_x=t_y'-t_x'$. This proves that $d_\infty$ is well-defined.

The map $d_\infty$ is $G$-invariant by Lemma~\ref{l_act_G_norms} and its restriction to $\A$ is $d$. Therefore $d_\infty=d$. 
\end{proof}

\begin{Remark}\label{r_d_serre}
\begin{enumerate}

\item  Let $x,y\in \IL$. Then by Lemma~\ref{l_description_A_L} and (A3'), we can write $x=[\eta_{b,t_x}]$ and $y=[\eta_{b,t_y}]$, for some basis $b$ and $t_x,t_y\in \Lambda$. Let $a_x,a_y\in \cF$ be such that $\omega(a_x)=t_x$ and $\omega(a_y)=t_y$. The closed balls of center $0$ and radius $1$ for $\eta_{b,t_x}$ and $\eta_{b,t_y}$ are respectively $\cO b_1\oplus \cO a_x b_2$ and $\cO b_1\oplus \cO a_y b_2$. Therefore $d_\infty$ is the distance $d$ introduced in \cite[II 1.1]{serre1977arbres} and \cite[Chapter 4; 3]{chiswell2001introduction}.

\item If we consider $\mathrm{SL}_n$, with $n\geq 3$ instead of $\mathrm{SL}_2$, then $\A$ is isomorphic to $\R^{n-1}$. Then  $d_\infty$ restricts to a uniform norm on $\A$ (i.e there exists  a basis $(e_i)_{i\in \llbracket 1,n-1\rrbracket}$ of $\A$  such that for all $(x_i)\in \R^{n-1}$, we have $d_\infty(\sum_{i=1}^{n-1} x_i e_i,0)=\max \{|x_i|\mid i\in \llbracket 1,n-1\rrbracket\}$). For the definition of $d$ on $\IN$, we usually choose a Euclidean norm on $\A$, which is invariant under $W^v$ and we apply the process of the beginning Subsection~\ref{ss_distance}.  Therefore $d_\infty$ and $d$ differ if $n\geq 3$. 
\end{enumerate}
\end{Remark}

We now  use this interpretation of $d$ to describe the $G$-orbits of $\I_{\sL}$.

For $g=\begin{psmallmatrix}
    a&b\\
    c& d
\end{psmallmatrix}\in \mathrm{GL}_2(\cF)$, we set  \[\omega(g)=\min\left(\omega(a),\omega(b),\omega(c),\omega(d)\right).\]\index[notation]{o@$\omega$}
 Recall that when we identify $\I_{\sL}$ and $\ve(\I_{\sN})$, we have  $0_\I=[\cO e_1\oplus \cO e_2]\in \I_{\sL}$.

For $g\in \mathrm{GL}_2(\cF)$, we have (\cite[4 Lemma 3.5]{chiswell2001introduction}) \begin{equation}\label{e_action_on_0}
    d(g.0_\I,0_\I)=\omega(\det(g))-2\omega(g).
\end{equation}

\begin{Proposition}\label{p_G_orbits}
    The group $\mathrm{GL}_2(\cF)$ acts transitively on $\I_{\sL}$. For  every $x,y\in \I_{\sL}$, we have $G.x=G.y$ if and only if $d(x,y)\in 2\Lambda$. 
    
      The set of $G$-orbits of $\I_{\sL}$ is in bijection with $\Lambda/2\Lambda$. More precisely, choose,  for each $a\in \Lambda/2\Lambda$, an element $\xi_{a}\in \cF$ such that $a=\omega(\xi_a)+2\Lambda$ and set $L_a=\cO e_1\oplus \cO \xi_a e_2$, where $(e_1,e_2)$ is the canonical base of $\cF^2$. Then $\I_{\sL}=\bigsqcup_{a\in \Lambda/2\Lambda} G.[L_a]$ and if $x\in \I_{\sL}$, then $x\in G.[L_a]$, where $a=d(0,x)+2\Lambda$. 
\end{Proposition}

\begin{proof}
    As the action of $\mathrm{GL}_2(\cF)$ is transitive on the set of bases of $\cF^2$, $\mathrm{GL}_2(\cF)$ acts transitively on $\I_{\sL}$. 

    Let $x,y\in \I_\sL$ and assume that $y=g.x$, for some $g\in G$. Write $x=h.0_\I$, with $h\in \mathrm{GL}_2(\cF)$. We have $d(x,y)=d(h.0_\I,gh.0_\I)=d(0_\I,h^{-1}gh.0_\I)=-2\omega(h^{-1}gh)\in 2\Lambda$, by \eqref{e_action_on_0}.
    
     Conversely, let $x,y\in \I_{\sL}$ be such that $D:=d(x,y)\in 2\Lambda$.     By Proposition~\ref{p_d_infty} and Remark~\ref{r_d_serre}, there exist a basis $(b_1,b_2)$ of $\cF^2$ and $z\in \cO$ such that $\omega(z)=D$, $x=[\cO b_1\oplus \cO b_2]$ and $y=[\cO b_1\oplus \cO z b_2]$. Let $z_1\in \cO$ be such that $\omega(z_1)=D/2$. Let $g\in G$ whose matrix in $(b_1,b_2)$ is $\begin{psmallmatrix}
        z_1^{-1} & 0 \\0 & z_1
    \end{psmallmatrix}$. 
        Then $g.x=[\cO b_1\oplus \cO (z_1)^2 b_2]=y$ since $\cO z_1^2=\{y\in \cF\mid \omega(y)\geq \omega(z_1^2)\}=\cO z$. Consequently, $x$ and $y$ belong to the same $G$-orbit.

        Let $a,a'\in \Lambda/2\Lambda$.  Then by Remark~\ref{r_d_serre},  $d([L_a],[L_{a'}])=\omega(\xi_a\xi_{a'}^{-1})=\omega(\xi_a)-\omega(\xi_{a'}\in \Lambda$. Therefore $[L_a]$ and $[L_{a'}]$ are in the same $G$-orbit if and only if $a=a'$. 

        Let now $x\in \I_{\sL}$. Set $a=d(0_\I,x)+2\Lambda\in \Lambda/2\Lambda$. Let $A$ be an apartment containing $0_\I$ and $x$, which exists by Corollary~\ref{c_frndl_pr_tr}. According to Proposition~\ref{p_A2}, there exists $g\in G$ such that $g.A=\A$ and $g$ fixes $A\cap \A$. Then $g.0_\I=0_\I$ and $d(0_\I,x)=d(0_\I,g.x)$. Up to considering $\tilde{s}g.x$ instead of $g.x$, we can assume that $g.x$ and $[L_a]$ have the same sign,  when $(\A,0_\I)$ is identified with $(\R,0)$. Then $d(g.x,[L_a])=d(0,[L_a])-d(0,g.x)=\omega(\xi_a)-d(0,x)\in 2\Lambda$. Therefore $g.x\in G.[L_a]$ and $x\in G.[L_a]$ which concludes the proof of the proposition.
        
\end{proof}

 \section{Tree structures on $\I_\sN$ and $\I_\sL$}
 
 In this section, we prove that $\I_{\sN}$ and $\I_{\sL}$, equipped with their distance, are trees. Recall that $\Lambda=\omega(\cF^\times)$.

\begin{Definition}
  A $\Lambda$-\textbf{metric space} is a metric space $(X,d)$ where $X$ is a set and $d$ is a distance on $X$ with values in $\Lambda\cap \R_{\geq 0}$.   Let $X$ be a set and $d:X\rightarrow \Lambda\cap \R_+$ be a distance. A $\Lambda$-\textbf{segment} of $X$ is a subset of $X$ of the form $\alpha([a,b]\cap \Lambda)$, where $a,b\in \Lambda$ and $\alpha:[a,b]\rightarrow X$ is an isometry. The \textbf{endpoints} of the segment are $\alpha(a)$ and $\alpha(b)$ and we say that the segment joints $\alpha(a)$ and $\alpha(b)$. 

    A \textbf{$\Lambda$-tree} is a $\Lambda$-metric space $(X,d)$ such that: \begin{enumerate}
        \item $(X,d)$ is geodesic, i.e for every $x,y\in X$, there exists a segment in $X$ with endpoints $x$ and $y$,

        \item if two segments of $X$ intersect in a single point, which is an endpoint of both, then their union is a segment of $X$. 

        \item the intersection of two segments with a common endpoint is also a segment.
    \end{enumerate}
\end{Definition}

Note that the $\Z$-trees are the graphs which are connected and without circuit, by \cite[2 Example 2]{chiswell2001introduction}.

For $x,y\in \I_{\sL}$, we denote by $[x,y]_{\sL}$ the segment of $\I_{\sL}$ between $x$ and $y$, i.e $[x,y]_{\sL}=[x,y]\cap \I_{\sL}$. If $x,y\in \Lambda\subset \R$, we set $[x,y]_{\Lambda}=[x,y]\cap \Lambda$. 
 
 \begin{Theorem}(see \cite[Chapitre 2, 1.1 Théorème 1]{serre1977arbres} and \cite[Chapter 4, Theorem 3.4]{chiswell2001introduction})
 \begin{enumerate}

 \item The metric space  $(\I_\sN,d)$ is an $\R$-tree.

\item  The metric space $(\I_\sL,d)$ is a $\Lambda$-tree, where $\Lambda=\omega(\cF^\times)$.  
 
 \end{enumerate}
 \end{Theorem}
 
 \begin{proof}
(1) By definition of $d$, if $x,y\in \I_{\sN}$ and if $g\in G$ is such that $g.x,g.y\in \A$, then $[x,y]$ is isometric to $[g.x,g.y]$, which proves that $(\I_\sN,d)$ is geodesic. Moreover if $x,y\in \I_{\sL}$, then by Lemma~\ref{l_description_A_L}, $g.[x,y]_{\sL}=[g.x,g.y]_{\sL}\simeq [g.x,g.y]_{\Lambda}$.  Thus $\I_{\sL}$ is a geodesic $\Lambda$-metric space in the sense of \cite[Chapter 1; 2]{chiswell2001introduction}. 

 Let $J_1,J_2$ be two segments of $\I_{\sN}$, which intersect in a single point, which is an end point of both. We can write $J_1=[x,y]$ and $J_2=[x,z]$, with $x,y,z\in \I_{\sN}$. Using apartment isomorphisms, we can assume that $[x,y]\subset \A$.  If $z\in \A$, then $[z,x]\cup [x,y]$ is the segment $[z,y]$. Assume now $z\notin \A$. 
 Let $a=\pr_\A(z)$, with the notation of Proposition~\ref{p_pr_apt}. Then by Proposition~\ref{p_pr_apt}, we have $[z,y]=[z,a]\sqcup ]a,y]$ and $[z,x]=[z,a]\cup ]a,x]$. Identify $\A$ and $\R$. Up to changing the identification, we can assume that $x\geq a$ and that $x\leq y$, if $a=x$. As $[z,x]\cap [x,y]=\{x\}$, we have $y\geq x$. Therefore $[a,y]=[a,x]\cup [x,y]$ and $[z,y]=[z,a]\cup [a,x]\cup [x,y]=[z,x]\cup [x,y]$, which proves that $J_1\cup J_2$ is a segment of $\IN$.

Assume moreover $x,y,z\in \I_{\sL}$. As $J_1\cup J_2$ is  a segment, it is contained in an apartment of $\IN$. Thus we can assume that $x,y,z\in \A$ and then $[z,y]_\sL=[z,x]_{\sL}\cup [x,y]_{\sL}$ is a segment of $\I_{\sL}$.

Let  now $J_1,J_2$ be two segments with a common endpoint. Write $J_1=[x,y]$ and $J_2=[x,z]$, with $x,y,z\in \I_{\sN}$. Let $A$ be an apartment containing $x$ and $y$. Set $a=\pr_A(z)$. Then $[x,y]\cap [x,z]\subset [x,z]\cap A=[x,a]$. Therefore $[x,y]\cap [x,z]=[x,y]\cap [x,a]$. As $x,y$ and $a$ belong to $A\simeq \R$, we deduce that $[x,y]\cap [x,z]$ is a segment of $\I_{\sN}$. Assume moreover $x,y,z\in \I_{\sL}$. If $x,y$ and $z$ are contained in an apartment $A_{\sL}$, then $[x,y]_\sL\cap [x,z]_{\sL}$ is a segment of this apartment and thus it is a segment of $\I_{\sL}$. If $x,y,z$ are not contained in any apartment,  then $a\in \I_{\sL}$, by Proposition~\ref{p_pr_apt} and thus $[x,y]_{\sL}\cap [x,z]_{\sL}=[x,a]_{\sL}$. This proves the theorem.
 \end{proof}
 
 Note that when $\Lambda=\Z$, then $\I_{\sL}$ has no leaf (or terminal vertex) since $\A_{\sL}$ has no leaf. In this case, $\IN$ can be obtained from $\IL$ by attaching an edge to each adjacent vertices, (where two vertices $x,y$ are  adjacent if $d(x,y)=1$). For $x\in \IN$, we denote by $\sC_x$ the set of alcoves of $\IN$ dominating $x$ (i.e the set of alcoves $C$ such that $x\in \overline{C}$). 
 
\begin{Proposition}
\begin{enumerate}
\item Let $x\in \A$. Then $U_{-\alpha,x}$ acts transitively on $\sC_x\setminus \{]x^-,x[\}$. The map $U_{-\alpha,x}\rightarrow \sC_x$ defined by $u\mapsto u.]x,x^{+}[$ induces a bijection $U_{-\alpha,x}/U_{-\alpha,]x,x^+[}\simeq \sC_x\setminus \{]x^-,x[\}$. 

\item Let $x\in \IN$. Then $ |\sC_x|=\left\{\begin{aligned} &2 &\text{ if } x \notin \IL&\\ &|\kk|+1 &\text{ if }x\in \IL \end{aligned}\right.,$ where $\kk=\cO/\fm$ is the residue field of $\cF$. 

\end{enumerate}
\end{Proposition} 

\begin{proof}
(1) Let $x\in \A$. Let $C\in \sC_x\setminus \{]x^-,x[\}$. Then by Proposition~\ref{p_SC_tr}, there exists an apartment $A$ containing $\fQ_{-\infty}$ and $C$. Let $g\in G$ be such that $g.A=\A$ and $g$ fixes $A\cap \A$. Then by Proposition~\ref{p_fix_infty}, we can write $g=tu$, with $u\in U_{-\alpha}$ and $t\in T_{\cO}$. As $t$ fixes $\A$ pointwise, we can replace $g$ by $u$. Then $A$ contains $C$ and thus $A$ contains $x$. Therefore $u\in U_{-\alpha,x}$.  Moreover $u.C$ is an alcove of $\A$ dominating $x$ and thus $u.C\in \{]x^-,x[,]x,x^+[\}$. As $A\cap \A$ contains $]-\infty,x]$, we have  $u.C\neq ]x^-,x[$, which proves that $u.C=]x,x^+[$ and thus $U_{-\alpha,x}$ acts transitively on $\sC_x\setminus \{]x^-,x[\}$. 

Let $u\in U_{-\alpha,x}$. Then $u$ fixes $]x,x^+[$ if and only if there exists $\epsilon\in \R_{>0}$ such that $u$ fixes $x+\epsilon$ if and only if $u\in \bigcup_{\epsilon\in \R_{>0}} U_{-\alpha,x+\epsilon}=U_{-\alpha,]x,x^+[}$. We deduce (1).

(2) Let $x\in \A$.  The group $U_{-\alpha}$ is isomorphic to $(\cF,+)$ via the isomorphism $x_{-\alpha}$. Moreover $x_{-\alpha}^{-1}(U_{-\alpha,x})=\cF_{\geq x}$ and $U_{-\alpha,]x,x^+[}=\cF_{>x}$. Therefore $U_{-\alpha,x}/U_{-\alpha,]x,x^+[}\simeq \cF_{\geq x}/\cF_{>x}$. If $x\notin \Lambda$, then $\cF_{\geq x}=\cF_{>x}$. If $x\in \Lambda$, write $x=\omega(a)$, with $a\in \cF$. Then $\cF_{\geq x}/\cF_{>x}=(a\cF_{\geq 0})/(a\cF_{>0})\simeq \cF_{\geq 0}/\cF_{>0}=\cO/\fm=\kk$. Using Lemma~\ref{l_description_A_L}, we deduce  (2), for $x\in \A$. As $\IL$ is stabilized by $G$, we deduce the result. 
\end{proof}
 
When $\Lambda=\Z$, we recover the fact that  $\IN$ is  a homogeneous tree with valency $|\kk|+1$.

\begin{Remark}
Note that when $\Lambda=\Z$, our definition of alcoves slightly differs from the usual one. Indeed, in this case, the alcoves of $\A$ are usually the sets of the form $]n,n+1[$, with $n\in \Z$. In our case, we have two alcoves of $\A$ associated with each $x\in \R$: $C_x^+:=]x,x^+[$ and $C_x^-:=]x^-,x[$. If $x\in \R\setminus \Z$, then $C_x:=] \lfloor x \rfloor,\lceil x\rceil[$ can be recovered from $C_x^-$ and $C_x^+$ by taking  the interior of the intersection of all the half-apartments containing $C_x^-$ or $C_x^+$. Any apartment containing $C_x^\pm$ contains $C_x$ and $\Fix_G(C_x)=\Fix_G(C_x^{\pm})$.
\end{Remark}

\section{Parohoric viewpoint on the tree}\label{s_parah_tree}
We defined the tree of $\mathrm{SL}_2(\cF)$ as a quotient of a  set of ultrametric norms on $\cF^2$. This approach has several advantages: it is explicit, relatively elementary, it enables to make computations, etc. However, we do not know if this approach enables to define the Bruhat--Tits building of any split reductive group or the masure of a Kac--Moody group. 

For arbitrary split reductive groups and Kac--Moody groups, Bruhat--Tits buildings and masures are defined using the ``parahoric subgroups''. In this section, we make this construction in the special case of $\mathrm{SL}_2$  and relate it to the construction we made.

For $x\in \A$, recall that we denoted by $G_x$ the fixator of $x$ in $G$.  The group $G_x$ is the \textbf{parahoric subgroup associated with $x$.}

We saw that: \begin{enumerate}
\item $\Stab_G(\A)=N$,

\item The restriction to $\A$ of the action of $N$ is given as follows: $\begin{psmallmatrix}0 & 1 \\ -1 & 0\end{psmallmatrix}$ acts on $\A$ by $-\Id$ and for $x\in \A$,  $\bt_x=\begin{psmallmatrix} x & 0\\ 0 & x^{-1}\end{psmallmatrix}$ acts by the translation of vector $2x$ on $\A$, for $x\in \R\simeq \A$. 

\item For every $x\in \A$, we have $G_x=\begin{psmallmatrix} \cO & \cO_{\geq -x}\\ \cO_{\geq x} & \cO\end{psmallmatrix}$ by Proposition~\ref{p_Br_dec}.
\end{enumerate}

These three properties actually enable to construct the Bruhat-Tits tree of $G$ with the following more abstract construction.

Let $\fT=\begin{psmallmatrix} * & 0\\0& *\end{psmallmatrix}\cap \mathrm{SL}_2$.
The root system of $(\mathrm{SL}_2,\fT)$ is $\{\pm \alpha\}$, where $\alpha:\fT\rightarrow \mathbb{G}_m$ is defined by $\alpha(\begin{psmallmatrix} x& 0 \\0 & x^{-1}\end{psmallmatrix})=x^2$, for $x$ a non-zero element of a field. The coroot system of  $(\mathrm{SL}_2,\fT)$ is $\{\pm \alpha^\vee\}$, where $\alpha^\vee:\mathbb{G}_m\rightarrow \fT$ is defined by $x\mapsto \begin{psmallmatrix} x& 0 \\0 & x^{-1}\end{psmallmatrix}$. 
Let $X$ and $Y$ be the character and cocharacter lattices of $(\mathrm{SL}_2,\fT)$. We have $X=1/2\Z\alpha$ and $Y=\Z\alpha^\vee$. Set $\A=Y\otimes \R$.  We have $\A\simeq \R$ and  $N$ acts on $\A$ via the  action defined above. Let $\sim$ be the equivalence relation on $G\times \A$ defined as follows. For $(g,x),(h,y)\in G\times \A$, \[(g,x)\sim (h,y) \Leftrightarrow (\exists n\in N\mid n.x=y\text{ and } g^{-1}hn\in G_x).\]

One then sets $\I_{\cP}=\I(G,\cF,\omega)=(G\times \A)/\sim$. For $g,x\in G\times \A$, we denote by $[g,x]$ the class of $(g,x)$ for $\sim$.  We embed $\A$ in $\I_{\cP}$ via $a\mapsto [a,1]$, for $a\in \A$. Then $G$ naturally acts on $\I_{\cP}$ via $g.[h,x]=[gh.x]$, for $g,h\in G$, $x\in \A$. It is easy to check that if $x\in \A$, then the fixator of $x$ in $G$ for this action is $G_x$.

\begin{Proposition}
Let $f:\I_{\cP}\rightarrow \I_{\sN}$ be defined by $f([g,x])=g.x$, for $g\in G$ and $x\in \A$. Then $f$ is well-defined and is a $G$-equivariant bijection between $\I_\cP$ and $\I_{\sN}$. 
\end{Proposition}

\begin{proof}
Let $g,h\in G$ and $x,y\in \A$ be such that $[g,x]=[h,y]$. Then there exists $n\in N$ such that $n.x=y$ and $g^{-1}hn\in G_x$.  Therefore $g^{-1}hn.x=x$ and $h.y=g.x$, which proves that $f$ is well-defined. 

As we have $\I_{\sN}=\bigcup_{g\in G}g.\A$, $f$ is surjective. Let $g,h\in G$ and $x,y\in \A$ be such that $f([g,x])=f([h,y])$. Then $h^{-1}g.x=y$. As $x,y\in \A$, there exists $n\in N$ such that $h^{-1}g.x=n.x$, by Corollary~\ref{c_W_gp_int_apt}. Then $g^{-1}hn\in G_x$ and thus $[g,x]=[h,y]$. Therefore $f$ is injective, which proves the  proposition.
\end{proof}

Note that a notion of morphism of buildings is defined in \cite{appenzeller2026morphisms}. We can then compare the two buildings for this definition, see \cite[Theorem 5.19]{appenzeller2026morphisms}.

\chapter{Kac--Moody algebras and vectorial apartment of a masure}\label{C_KM_alg}

The aim of this paper is to study Kac--Moody groups over valued fields. There are several possible definitions of Kac--Moody groups, but for all of them, a Kac--Moody group is supposed to  ``integrate'' a  Kac--Moody algebra, in some sense. In this chapter, we recall some facts on Kac--Moody algebras.

Since the works of Cartan and Killing at the end of the nineteenth
century, the classification of finite dimensional semi-simple Lie algebras over $\C$ is known. One way to describe
it  is to use Cartan matrices: to each such Lie algebra, one  associates a
Cartan matrix, which is a matrix whose coefficients are  integers and satisfy some conditions. Serre used Cartan matrices to give a presentation of such Lie algebras by generators and relations (see
\cite[Chapitre 6]{serre1966algebres}).  

In 1967, Kac and Moody independently introduced a new class of Lie algebras, the
Kac--Moody algebras (or more precisely, the Kac--Moody Lie algebras).  A Kac--Moody algebra is a Lie algebra defined by  Serre’s presentation, but for a Kac--Moody matrix (often called a generalized Cartan matrix) instead of a Cartan matrix.

Kac--Moody algebras are infinite dimensional, unless they are associated with a Cartan matrix. 
Despite this difference, Kac--Moody algebras
share many properties with finite dimensional semi-simple Lie algebras, which we now sketch. A Kac--Moody algebra $\ffg$ comes together with a Cartan Lie subalgebra $\fh$, which is a finite dimensional commutative Lie algebra. Then $\ffg$ decomposes as a direct sum of its eigenspaces with respect to $\fh$: $\ffg=\fh\oplus\bigoplus_{\alpha\in \Delta} \ffg_{\alpha}$, where $\Delta$ is the set of $\alpha\in \fh^*\setminus \{0\}$ whose associated eigenspace $\ffg_\alpha$ is non-zero. One can also define a vectorial Weyl group $W^v$, which is infinite in general, which acts on $\fh^*$ and stabilizes $\Delta$. The elements of $\Delta$ which are the transform of a simple roots by an element of $W^v$ are called real roots and the other roots are called imaginary. Then $\ffg$ is finite dimensional if and only if all the roots are real. 

In this chapter, we introduce briefly Kac--Moody algebras and their basic properties. We also introduce the vectorial apartment of a masure. It is a finite dimensional vector space - more precisely, it is an $\R$ form of the Cartan subalgebra of $\ffg$ - equipped with the hyperplane arrangement defined by the real roots of $(\ffg,\fh)$.

This chapter is organized as follows. In Section~\ref{s_example_sln}, we describe the Kac--Moody algebra $\mathfrak{sl}_n(\C)$. In Section~\ref{s_KM_algebras}, we define the Kac--Moody algebras. In Section~\ref{s_vectorial_apartment_KM_algebra}, we introduce the vectorial apartment associated with the Kac--Moody algebra. In Section~\ref{s_Vectorial apartments}, we describe the vectorial apartment in examples: in the affine case and in the indefinite case for particular choices of size 2 Kac--Moody matrices.

\section{Example of $\mathfrak{sl}_n(\C)$}\label{s_example_sln}

We start with the example of $\mathfrak{sl}_n(\C)$. We explicitly describe the roots, coroots, Weyl group, root spaces decompositions etc. 

Let $n\in \Z_{\geq 2}$. Let $\ffg=\mathfrak{sl}_n(\C)=\{a\in \cM_n(\C)\mid \mathrm{tr}(a)=0\}$, equipped with the bracket $[a,b]=ab-ba$, for $a,b\in \ffg$. The standard Cartan subalgebra $\fh$ is the set of diagonal matrices whose trace is $0$. For $\alpha\in \fh^*$, set $\ffg_{\alpha}=\{x\in \ffg\mid \forall h\in \fh, [h,x]=\alpha(h) x\}$. Let $\Delta=\{\alpha\in \fh^*\setminus \{0\}\mid \ffg_\alpha\neq 0\}$ be the \textbf{root system} of $(\ffg,\fh)$. 

 For $i,j\in \llbracket 1,n\rrbracket$, we denote by $E_{i,j}$ the matrix having $1$ in $i,j$ and $0$ elsewhere. 

For $i\in \llbracket 1,n-1\rrbracket$, we define $\alpha_i:\fh\rightarrow \C$ by  $\alpha_i(\sum_{j=1}^n h_j E_{j,j})=h_i-h_{i+1}$ for all $(h_j)\in \C^n$ such that $\sum_{j=1}^n h_j=0$. For $(i,j)\in \llbracket 1,n\rrbracket^2$ such that $i\neq j$, we define $\alpha_{i,j}:\fh\rightarrow \C$ by  $\alpha_{i,j}(\sum_{k=1}^n h_kE_{k,k})=h_i-h_j$ (we have $\alpha_{i,j}=\sum_{k=i}^{j-1} \alpha_k$ if $i<j$). We have $\alpha_{i,j}=-\alpha_{j,i}$.

 We have $\ffg_{\alpha_{i,j}}=\C E_{i,j}$ and $\ffg_0=\fh$.

If $\Delta_+=\{\alpha_{i,j}\mid (i,j)\in \llbracket 1,n\rrbracket^2, i<j\}$, we have $\Delta=\Delta_+\cup -\Delta_+$ and $\ffg=\fh\oplus \fn^+\oplus \fn^-$, where  $\fn^-=\bigoplus_{\alpha\in \Delta_-}\ffg_\alpha$ and $\fn^+=\bigoplus_{\alpha\in \Delta_+}$.   Note that $\fn^+$ (resp. $\fn^-$) is the set of strictly upper (resp. lower) triangular matrices of $\ffg$.

For all $\alpha\in \Delta$, there exists a unique $\alpha^\vee\in \fh_{\alpha}:=[\ffg_{\alpha},\ffg_{-\alpha}]$ such that $\alpha(\alpha^\vee)=2$. The set $\Delta^\vee$ is the \textbf{dual root system} of $(\ffg,\fh)$.  For $i,j\in \llbracket 1,n\rrbracket$ such that $i\neq j$, we have $\alpha_{i,j}^\vee=E_{ii}-E_{jj}$. The Cartan matrix $(\alpha_{i,i+1}(\alpha_{j,j+1}^\vee))_{i,j\in \llbracket 1,n-1\rrbracket}$ is the matrix with $2$ on the diagonal, $-1$ on the lines above and below the diagonal and $0$ elsewhere.

\begin{figure}[h]
\centering
\includegraphics[scale=0.3]{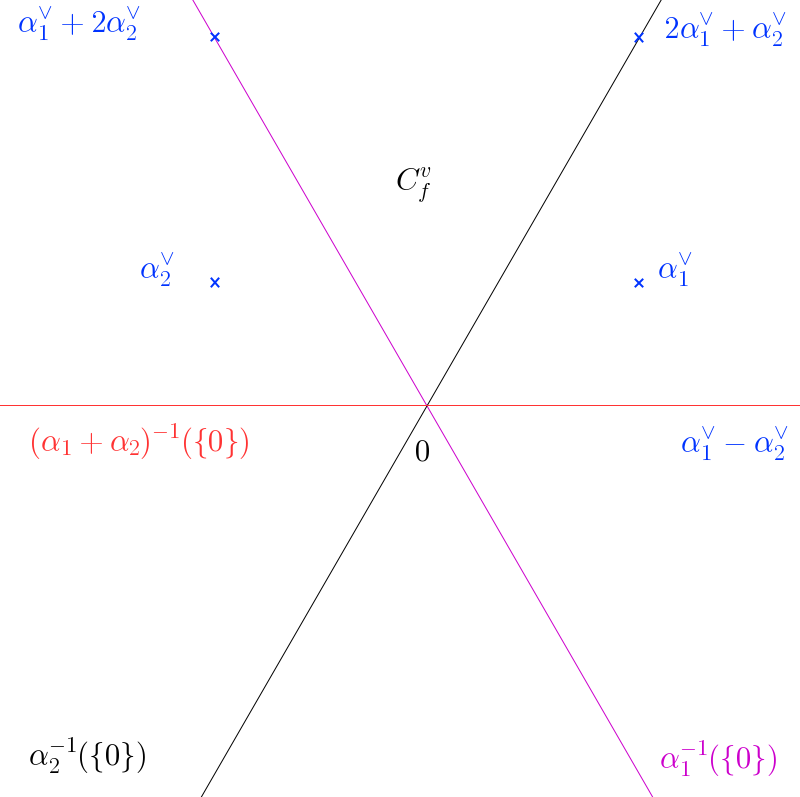}
\caption{Vectorial apartment of $\mathfrak{sl}_3$. We have $\Delta=\pm\{\alpha_1,\alpha_2,\alpha_1+\alpha_2\}$ and $\Delta^\vee=\pm \{\alpha_1^\vee,\alpha_2^\vee, \alpha_1^\vee+\alpha_2^\vee\}$. The fundamental chamber $C^v_f$ is a Weyl chamber and $W^v$ acts simply transitively on the Weyl chambers.}\label{figSysteme de racine SL3}
\end{figure}

Let $\A=\fh_\R:=\fh\cap \cM_n(\R)$ (see Figure~\ref{figSysteme de racine SL3} for $n=3$). Then by restriction we can regard $\Delta$ as a subset of $\A^*$   and we have $\Delta^\vee\subset \A$. Let $W^v=\fS_n$ be the set of permutations of $\llbracket 1,n\rrbracket$. Then $W^v$ acts on $\A$ by setting $w.\sum_{i\in \llbracket 1,n\rrbracket} x_i E_{i,i}=\sum_{i\in \llbracket 1,n\rrbracket} x_i E_{w(i),w(i)}$, for $w\in W^v$ and $\sum_{i\in \llbracket 1,n\rrbracket} x_i E_{i,i} \in \A$.

Let $\cH=\{\ker(\alpha)\mid \alpha\in \Delta\}$. The elements of $\cH$ are called vectorial walls. The hyperplane arrangement $\cH$ enables to define the Weyl chambers and a simplicial structure on $\A$. The (open) Weyl chambers are the connected components of $\A\setminus \bigcup_{\alpha\in \Delta} \ker(\alpha)$ and the faces are the relative interior of the sets of the form $\overline{C^v}\cap \bigcap_{H\in \cH'} H$, where $C^v$ is a Weyl chamber and $\cH'\subset \cH$.

\section{Kac--Moody algebras}\label{s_KM_algebras}

In this section, we define Kac--Moody algebras by generators and relations. We first define Kac--Moody data, which are the data to which Kac--Moody algebras and groups are associated. Then we recall the root decomposition of these Lie algebras. Our references are mainly~\cite{kac1994infinite},   \cite[Chapter I]{kumar2002kac} and \cite[II]{marquis2018introduction}.

\subsection{Kac--Moody data}\label{ss_KM_data}

In this subsection, we define Kac--Moody data. These are the data to which  Kac--Moody groups are associated. 

\subsubsection{Kac--Moody matrices}

\begin{Definition}
 A \textbf{Kac--Moody matrix}\index{Kac--Moody matrix} (also known as a generalized Cartan matrix) is a square matrix $A=(a_{i,j})_{i,j\in I}$ with integers coefficients, indexed by  a finite set $I$ and satisfying: 
\begin{enumerate}
\item $\forall i\in I,\ a_{i,i}=2$

\item $\forall (i,j)\in I^2\mid i \neq j,\ a_{i,j}\leq 0$

\item $\forall (i,j)\in I^2,\ a_{i,j}=0 \Leftrightarrow a_{j,i}=0$.
\end{enumerate}

The matrix $A$ is said to be \textbf{decomposable}\index{decomposable Kac--Moody matrix} if for some reordering of the indices, one can write $A$ as a non trivial block diagonal matrix $\begin{pmatrix}
A_1 & 0\\ 0& A_2
\end{pmatrix}$. We say that $A$ is \textbf{indecomposable}\index{indecomposable Kac--Moody matrix} if it is not decomposable.
\end{Definition}

\subsubsection{Kac--Moody data}

 A  \textbf{Kac--Moody datum}\index{Kac--Moody datum} (also called ``root generating system'')   of type $A$ is a $5$-tuple $\mathcal{S}=(A,X,Y,(\alpha_i)_{i\in I},(\alpha_i^\vee)_{i\in I})$ made of a Kac--Moody matrix $A$ indexed by $I$, of two dual free $\Z$-modules $X$ (of \textbf{characters}\index{character}) and $Y$ (of \textbf{cocharacters}\index{cocharacter}) of finite rank $\mathrm{rk}(X)$, a family $(\alpha_i)_{i\in I}$ (of \textbf{simple roots}\index{simple roots}) in $X$ and a family $(\alpha_i^\vee)_{i\in I}$ (of \textbf{simple coroots}\index{simple coroots}) in $Y$. They have to satisfy the following compatibility condition: $a_{i,j}=\alpha_j(\alpha_i^\vee)$ for all $i,j\in I$. We say that $\cS$ is \textbf{free} (resp. \textbf{cofree}) if the family $(\alpha_i)_{i\in I}$ (resp. $(\alpha_i^\vee)_{i\in I}$) is free.

 In this paper, we only consider free Kac--Moody data. However, we consider cofree and non-cofree Kac--Moody data.
 
\begin{Remark}
 Note that non-free and non-cofree roots system appear naturally. For example, the functor $\cK\mapsto \mathrm{SL}_2(\cK[t,t^{-1}])$ (where $\cK$ is a field) is a Kac--Moody group associated with the Kac--Moody data $(\begin{psmallmatrix} 2 & -2\\ -2& 2\end{psmallmatrix},\Z,\Z,2,2)$ (see \cite[Example 7.85]{marquis2018introduction}). The functor $\cK\mapsto \mathrm{SL}_2(\cK[t,t^{-1}])\rtimes \cK^*$ is associated with a free but non cofree Kac--Moody data.

 In this paper, we only consider free Kac--Moody data. Without this assumption, the theory could be very different, since, for example, the vectorial chambers of the standard apartment could be reduced to $\{0\}$ (this is the case of the Kac--Moody datum of affine $\mathrm{SL}_2$ above). However, if $\cS$ is a non-free Kac--Moody datum, the associated Kac--Moody group $\fG$ embeds in a Kac--Moody group $\fG'$ associated  with a free Kac--Moody datum, by \cite[Proposition 7.66]{marquis2018introduction}. If $\cF$ is a valued field, we can then study $\fG(\cF)$ via its action on the masure of $\fG'(\cF)$.  
\end{Remark}

  Let $Q$\index[notation]{q@$Q$} be the free $\Z$-module $\bigoplus_{i\in I} \Z\alpha_i$. It is called the \textbf{root lattice} of $\cS$. Let $Q^\vee=\sum_{i\in I}\Z\alpha_i^\vee$. We define $\htt:Q\rightarrow \Z$ by $\htt(\sum_{i\in I} n_i\alpha_i)=\sum_{i\in I} n_i$\index[notation]{h@$\htt$}, for $(n_i)\in \Z^I$.

\subsection{Kac--Moody algebras}\label{ss_KM_alg}

\subsubsection{Definition of the Kac--Moody algebra associated with a Kac--Moody datum} Let $A=(a_{i,j})_{i,j\in I}$ be a Kac--Moody matrix and  $\mathcal{S}=(A,X,Y,(\alpha_i)_{i\in I},(\alpha_i^\vee)_{i\in I})$ be a free Kac--Moody datum. Let $\fh=Y\otimes \C$.   The Kac--Moody algebra $\ffg=\ffg_\cS$ is the Lie algebra over $\C$ generated by $\fh$ and the symbols $e_i,f_i$ for $i\in I$ subject to the following relations for all $i,j\in I$ and $h\in \fh$: 

$(R_1)\ [\fh,\fh]=0$

$(R_2)\ [h,e_i]=\alpha_i(h)e_i$; $[h,f_i]=-\alpha_i(h)f_i$,

$(R_3)\ [e_i,f_j]=-\delta_{i,j} \alpha_i^\vee$

$(R_4)\ (\mathrm{ad}\ e_i)^{1-a_{i,j}}(e_j)=0$, $i\neq j$

$(R_5)\ (\mathrm{ad}\ f_i)^{1-a_{i,j}}(f_j)=0$, $i\neq j$.

The $e_i,f_i$, $i\in I$ are called the \textbf{Chevalley generators}\index{Chevalley generators of a Kac--Moody algebra} of $\ffg$. Note that when $\ffg=\mathfrak{sl}_n(\C)$, for $n\in \Z_{\geq 2}$, we have $e_i=E_{i,i+1}$ and $f_i=E_{i+1,i}$, for $i\in \llbracket 1,n-1\rrbracket$.

\begin{Remark}
\begin{enumerate} 
\item In general (for example in \cite{kac1994infinite}), Kac-Moody algebras are associated with free and cofree root data (and usually there is no minus in (R3)). We follow \cite[Definition 7.13]{marquis2018introduction}. Note also  that the definition we take slightly differs from the one of \cite{kac1994infinite}, where $(R_4)$ and $(R_5)$ are not used but where $\ffg$ is defined as a quotient of a Lie algebra satisfying $(R_1)$ to $(R_3)$. In the case where the Kac--Moody matrix is symmetrizable, the two definitions coincide, but in general, we do not know, see \cite[Remark 3.19]{marquis2018introduction}.

\item The Kac--Moody algebra $\ffg_{\cS}$ depends only on $\cS_\C:=(A,Y\otimes \C, (\alpha_i)_{i\in I}, (\alpha_i^\vee)_{i\in I})$. 

\end{enumerate}
\end{Remark}

\subsubsection{Root space decomposition}\label{sss_rt_space_dec}

\begin{Proposition}\label{p_rt_sp_dec}
 The natural map $\fh\rightarrow \ffg$ is injective.  Let $Q=\bigoplus_{i\in I} \Z\alpha_i\subset \fh^*$ be  the root lattice of $\ffg$. For  $\alpha\in \fh^*$,
   set $\ffg_\alpha=\{x\in \ffg\mid  [h,x]=\alpha(h)x,\forall h\in \fh\}$. Let $\Delta=\{\alpha\in X\setminus \{0\}\mid \ffg_\alpha\neq 0\}$\index[notation]{d@$\Delta$}.  
  Let $Q_\N=\bigoplus_{i\in I}\N\alpha_i\subset Q$. 
   Let $\fn= \fn^+$\index[notation]{n@$\fn=\fn^+$, $\fn^-$} (resp. $\fn^-$) be the  sub-Lie algebra of $\ffg$ generated by $\{e_i\mid i\in I\}$ (resp. $\{f_i\mid i\in I\}$). Then: \begin{enumerate}
   \item $\ffg=\bigoplus_{\alpha\in \Delta\cup \{0\}}\ffg_{\alpha}$ and $\ffg_\alpha$ is finite dimensional, for $\alpha\in \Delta$.
   
   \item Let $\Delta_+=\Delta\cap Q_\N$\index[notation]{d@$\Delta_+,\Delta_-$} and $\Delta_-=\Delta\cap (-Q_\N)$. Then $\Delta=\Delta_+\cup \Delta_-$.

   \item We have $\fn=\bigoplus_{\alpha\in \Delta_+} \ffg_\alpha$, $\fh=\ffg_0$ and $\fn^-=\bigoplus_{\alpha\in \Delta_-}\ffg_\alpha$. 
   
   \end{enumerate}
 \end{Proposition}
 
 \begin{proof}
First assume that  $\cS$ is cofree. Let $\mathfrak{i}$  be the ideal of $\tilde{\ffg}(A)$ (in the notation of \cite[Definition 3.13]{marquis2018introduction})
 generated by $\mathrm{ad}(e_i)^{1-{a_{i,j}}}e_j$, 
 $\mathrm{ad}(f_i)^{1-a_{i,j}}f_j$, $i,j\in I$. By (6)
  of the proof of  \cite[Proposition 3.14]{marquis2018introduction}, 
  $\mathfrak{i}$ intersects $\fh$ (regarded as a subspace of 
$\tilde{\ffg}(A)$) trivially and hence $\fh\subset \tilde{\ffg}(A)$ embeds  in $\tilde{\ffg}(A)/\mathfrak{i}=\ffg$. Points (1) and (2) are proved in \cite[3.5, after Remark 3.19]{marquis2018introduction} and (3) follows from \cite[Proposition 3.14 (5)]{marquis2018introduction}. This proves the proposition when $\cS$ is cofree.

We now longer assume that $\cS$ is cofree. By \cite[7.3.2]{marquis2018introduction}, there exists an ``extension'' $\tilde{\cS}=(A,\tilde{X},\tilde{Y},(\alpha_i)_{i\in I},(\tilde{\alpha_i^{\vee}})_{i\in I})$, which is both free and cofree. We have $\tilde{Y}=Y\oplus \overline{Y}$, for some lattice $\overline{Y}$. Then if $\overline{\fh}=\overline{Y}\otimes  \C$, we have $\ffg_{\tilde{\cS}}=\ffg_{\cS}\oplus_{\omega} \overline{\fh}$, for the notation of \cite[(7.12) and (7.13)]{marquis2018introduction}. In particular, we have a projection (of Lie algebras) $\pi:\ffg_{\tilde{\cS}}\rightarrow \ffg_{\cS}$ such that for all $\alpha\in \Delta$, $\pi(\ffg_{\tilde{\cS},\alpha})\subset \ffg_{\cS,\alpha}$, with obvious notations.  Therefore we get (1) and (3) for $\ffg_{\cS}$  and if $\alpha\in \Delta$, $\pi$ induces an isomorphism $\ffg_{\tilde{\cS},\alpha}\simeq \ffg_{\cS,\alpha}$ and $\pi(\tilde{\fh})=\fh$, where $\fh=Y\otimes \C$ and $\tilde{\fh}=\tilde{Y}\otimes \C$. As $\Delta$ is the same for $\ffg_{\cS}$ and $\ffg_{\tilde{\cS}}$, we already know (2), which proves the proposition. 
 \end{proof}

We write explicitly the root decomposition of an untwisted affine Kac--Moody algebra in \ref{sss_aff_la}.

\begin{Remark}\label{rqueNonnullité des ei}
 For $k\in \N$ and $x_1,\ldots,x_k\in \ffg$, we define $[x_1,\ldots,x_k]\in \ffg$ as follows. We set $[x_k]=x_k$ and if $\ell\in \llbracket 1,k-1\rrbracket$, $[x_\ell,\ldots,x_k]=[x_\ell,[x_{\ell+1},\ldots,x_k]]$.  If $k\in \N$ and  $i_1,\ldots,i_k\in I$, then  $[e_{i_1},\ldots,e_{i_k}]\in \ffg_{\alpha_{i_1}+\ldots +\alpha_{i_k}}$ and $[f_{i_1},\ldots,f_{i_k}]\in \ffg_{-\alpha_{i_1}-\ldots-\alpha_{i_k}}$.  Then by \cite[Proposition 3.14 (5)]{marquis2018introduction},  for all $\alpha\in Q_\N$ (resp. $-Q_\N$), $\ffg_\alpha$ is the vector space spanned by the $[e_{i_1},\ldots,e_{i_k}]$ (resp. $[f_{i_1},\ldots,f_{i_k}]$) such that $i_1,\ldots,i_k\in I$ and $\alpha_{i_1}+\ldots+ \alpha_{i_k}=\alpha$ (resp. $\alpha_{i_1}+\ldots+\alpha_{i_k}=-\alpha$). 

\end{Remark}

\section{Vectorial apartment of a Kac--Moody algebra}\label{s_vectorial_apartment_KM_algebra}

The masure is a union of subsets called apartments, which are all isomorphic to one of them called the standard apartment $\A$. As a set, $\A$ is $Y\otimes \R$, where $Y$ is the cocharacter lattice of $\cS$. It is then equipped with additionnal structures: the structure of a vectorial apartment, which we describe now and the structure of an affine apartment, which will be described in Subsection~\ref{ss_aff_apt}. 

In a first approximation, the vectorial apartment is $\A$, equipped with the hyperplane arrangement $\cH_0:=\{\alpha^{-1}(\{0\})\mid \alpha\in \Phi\}$, whereas the affine apartment is $(A,\cH_{\Lambda})$, where $\cH_{\Lambda}:=\{\alpha^{-1}(\{k\})\mid k\in \Lambda,\alpha\in \Phi\}$, where $\Phi\subset \Delta$ is the set of real roots defined in Subsection~\ref{ss_Weyl_gp_rts} and $\Lambda=\omega(\cF^\times)$, when we study a Kac--Moody group over a valued field $(\cF,\omega)$. Note that the vectorial apartment does not take into account the valuation on the field over which we study the Kac--Moody group.

The apartment $\A$ is naturally equipped with a Weyl group $W^v$, which is generated by reflections with respect to the elements of $\cH_0$. The set of simple roots $(\alpha_i)_{i\in I}$ enables to define a fundamental chamber $C^v_f$, which enables to define the Tits cone $\cT=W^v.\overline{C^v_f}$ and its opposite $-\cT$. The Tits cone naturally decomposes into vectorial faces, whose maximal ones are called vectorial chambers. The face decomposition of the Tits cone enables to define the ``face structure at infinity'' of the masures: it enables to define sector-germs and, more generally, sector-face-germs.

In the reductive case, $\cT=\A=-\cT$. In the (non-reductive indecomposable) Kac--Moody case however, $\cT$ is a strict cone of $\A$ and $-C^v_f\cap \cT$ is almost reduced to zero in general. It corresponds to the fact that in the reductive case, there is only one conjugacy class of Borel subgroups whereas in the Kac--Moody case, there are two.

We start by giving the definitions and then we describe the vectorial apartment when the Kac--Moody matrix has size $2$.

\subsection{Weyl group, real and imaginary roots}\label{ss_Weyl_gp_rts}

 Let $\mathcal{S}=(A,X,Y,(\alpha_i)_{i\in I},(\alpha_i^\vee))$ be a root generating system and $\ffg$ be the Kac--Moody algebra associated. Let $\A=Y\otimes \R$. If $i\in I$, we define the \textbf{simple reflection}\index{simple reflexion} $r_i$ of the space $\A$ by $r_i(v)=v-\alpha_i(v)\alpha_i^\vee$, for $v\in \A$. Then $r_i(\alpha_i^\vee)=-\alpha_i^\vee$ and $r_i$ fixes $\{x\in \A\mid\alpha_i(x)=0\}$. The \textbf{Weyl group of $\ffg$}\index{Weyl group of $\cS$} is the subgroup of $\mathrm{GL}(\A)$ spanned by the fundamental reflections. We denote it $W^v$\index[notation]{w@$W^v$}. By \cite[Proposition~1.3.21]{kumar2002kac}, $(W^v,\{r_i\mid i\in I\})$ is a Coxeter system. The action of $W^v$ on $\A$ induces an action of $W^v$ on $\A^*$ by the formula $w.f (x):=f(w^{-1}.x)$ for all $f\in \A^*$, $x\in \fh$ and  $w\in W^v$.

By \cite[Corollary~1.3.6]{kumar2002kac}, the set of roots $\Delta\subset \bigoplus_{i\in I} \Z\alpha_i\subset \A^*$ is stable under the action of $W^v$. Moreover for all $\alpha\in \Delta$ and $w\in W^v$, $\dim \ffg_{\alpha}=\dim \ffg_{w.\alpha}$. 
A root  is said to be \textbf{real}\index{real root} if it is of the form $w.\alpha_i$, for some $i\in I$ and $w\in W^v$. We denote by $\Phi$\index[notation]{p@$\Phi$} the set of \textbf{real roots}. Let $\Delta^{im}=\Delta\backslash \Phi$. An element of $\Delta^{im}$\index[notation]{d@$\Delta^{im}$} is called an \textbf{imaginary root}\index{imaginary root}. By Remark~\ref{rqueNonnullité des ei}, $ \ffg_{\alpha_i}=\C e_i\neq \{0\}$ and thus $\dim \ffg_\alpha=1$ for all $\alpha\in \Phi$. By Remark~\ref{rqueNonnullité des ei}, if $i\in I$ and $n\in \C$, $n\alpha_i\in \Delta$ if and only if $n\in \{-1,1\}$. Thus we get that if $\alpha\in \Phi$ and $n\in \C$ such that $n\alpha\in \Delta$, then $n\in \{-1,1\}$. On the contrary, by \cite[Proposition~5.4]{kac1994infinite}, if $\alpha\in \Delta^{im}$ and $r\in \Q^*$ are  such that $r\alpha\in Q$, then $r\alpha\in \Delta^{im}$.

Let $\alpha\in \Phi$. By definition, we can write $\alpha=w.\alpha_i$, with $w\in W^v$ and $i\in I$. We then set $\alpha^\vee=w.\alpha_i^\vee\in Y$\index[notation]{$\alpha$} and $r_\alpha=wr_iw^{-1}\in W^v$\index[notation]{r@$r_\alpha$}. By \cite[Lemma 4.19]{marquis2018introduction},  these elements are well-defined, independently of the choices of $w$ and $i$. We call $\alpha^\vee$ the \textbf{coroot associated with $\alpha$} and $r_\alpha$ the \textbf{reflection} associated with $\alpha$. We then have $r_\alpha(x)=x-\alpha(x)\alpha^\vee$, for $x\in \A$.

For $w\in W^v$, we set \[\Inv(w)=\{\alpha\in \Delta_+\mid w.\alpha\in \Delta_-\}.\] By \cite[1.3.14 Lemma]{kumar2002kac}, if $r_{i_k}\ldots r_{i_1}$ is a reduced expression of $w$, with $i_1,\ldots,i_k\in I$, then we have \begin{equation}\label{e_Kmr_1.3.14}
\Inv(w)=\{\alpha_{i_1},r_{i_1}.\alpha_{i_2},\ldots,r_{i_1}\ldots r_{i_{k-1}}.\alpha_{i_k}\}.
\end{equation}\index[notation]{i@$\Inv$}

In particular, $\Inv(w)\subset \Phi$.

\subsection{Vectorial apartment}\label{subsubVectorial apartment}

 Let $C^v_f:=\{x\in \A\mid \alpha_i(x)>0, \forall i\in I\}$\index[notation]{c@$C^v_f$} be the \textbf{vectorial fundamental chamber}\index{fundamental vectorial chamber}. For $J\subset I$, we set  $F^v(J)=\{v\in \A\mid\alpha_i(v)=0, \forall i\in J,\alpha_i(v)>0, \forall i\in J\backslash I\}$\index[notation]{f@$F^v(J)$}. Then the closure $\overline{C_f^v}$ of $C_f^v$ is the union of the $F^v(J)$ for $J\subset I$. The \textbf{positive} (resp. \textbf{negative}) \textbf{vectorial faces}\index{face (vectorial)} are the sets $w.F^v(J)$ (resp. $-w.F^v(J)$) for $w\in W^v$  and $J\subset I$. A \textbf{vectorial face}\index{vectorial face} is either a positive vectorial face or a negative vectorial face. Note that in certain cases (for example if $A$ is a Cartan matrix), non-zero vectorial faces can be positive and negative at the same times. We call \textbf{positive chamber} (resp. \textbf{negative}) every cone  of the form $w.C_f^v$ for some $w\in W^v$ (resp. $-w.C_f^v$). A vectorial face of codimension $1$ is called a \textbf{vectorial panel}.  
 
A vectorial face $F^v$ is said to be \textbf{spherical}\index{spherical vectorial face} if its fixator in $W^v$ is finite. Vectorial chambers and panels are examples of spherical vectorial faces.

 The \textbf{Tits cone}\index{Tits cone} $\mathcal T$ is defined by $\mathcal{T}=\bigcup_{w\in W^v} w.\overline{C^v_f}$. The \textbf{integral Tits cone} is $Y^+:=Y\cap \cT$\index[notation]{y@$Y^+$}
 
 If $J\subset I$, the fixator $W^v_{F^v(J)}$\index[notation]{w@$W^v_{F^v},W^v_{F^v(J)}$} of $F^v(J)$ in $W^v$ is \begin{equation}\label{e_Wv_FJ}
W^v_{F^v(J)}= \langle r_j\mid j\in J\rangle,
 \end{equation} by  \cite[Proposition 3.12 a]{kac1994infinite}.

\begin{Proposition}\label{p_T_cone} (\cite[Proposition 3.12 b,c, h]{kac1994infinite})
\begin{enumerate}
\item The fundamental chamber $\overline{C^v_f}$ is a fundamental domain for the action of $W^v$ on $\cT$, i.e., any orbit $W^v.h$ of $h\in \cT$ intersects $\overline{C^v_f}$ in exactly one point. In particular, the positive faces form a partition of $\cT$.

\item We have $\cT=\{h\in \A\mid\{\alpha\in \Delta_+ \mid  \alpha(h)<0\}\mathrm{\ is\ finite }\}$. In particular, $\cT$ is a convex cone of $\A$.

\item Let $\mathring{\cT}$ be the interior of $\cT$. Then $\mathring{\cT}$ is the union of the positive spherical vectorial faces of $\A$. In other words, $\mathrm{\cT}$ is the set of elements of $\cT$ which have a finite fixator in $W^v$.
\end{enumerate}

\end{Proposition}

For $x,y\in \A$, we write $x\leq_{Q^\vee}y$\index[notation]{z@$\leq_{Q^\vee}$} (resp. $x\leq_{Q^\vee_\R} y$\index[notation]{z@$\leq_{Q^\vee_\R}$})  if $y-x\in \bigoplus_{i\in I} \N\alpha_i^\vee$ (resp. $y-x\in \bigoplus_{i\in I} \R_{\geq 0} \alpha_i^\vee$).

By Proposition~\ref{p_T_cone}  and \cite[Lemma 2.4]{gaussent2014spherical}, we have: \begin{equation}\label{e_GR2.4}
    \lambda \leq_{Q^\vee} \lambda^{++}, x\leq_{Q^\vee_\R}x^{++}
\end{equation} for every $\lambda\in Y^{+}$ and $x\in \cT$.

\paragraph{Tits preorder on $\A$} We define the \textbf{Tits preorder}\index{Tits preorder} on $\A$ as follows. If $x,y$ in $\cT$, then one writes $x\leq y$ if $y-x\in \cT$. This is indeed a preorder by Proposition~\ref{p_T_cone}. By definition, it is $W^v$-invariant.

A hyperplane of the form $\alpha^{-1}(\{0\})$ for some $\alpha\in \Phi$ is called a \textbf{vectorial wall}\index{vectorial wall} of $\A$.

\subsection{Essential apartment}~\label{ss_ess_apt}

Let $\A_{\ines}=\bigcap_{i\in I} \ker(\alpha_i)$\index[notation]{a@$\A_{\ines}$}. Let $\A_{\mathrm{es}}=\A/\A_{\ines}$ be the \textbf{essentialization of $\A$}\index{essentialization of $\A$}. 

If $X$ is an affine space, we denote by $\mathrm{Aut}(X)$ its group of affine automorphisms. 
 
\begin{Lemma}\label{l_res_ess} (see \cite[Lemma 2.2.10]{hebert2018study})
Let $w\in W^v$. Then the map $\overline{w}:\A_{\mathrm{es}}\rightarrow \A_{\mathrm{es}}$ defined by $\overline{w}(x+\A_{\ines})=w(x)+ \A_{\ines}$ is well-defined. Moreover, the morphism $\Gamma: W^v\rightarrow \mathrm{Aut}(\A_{\mathrm{es}})$ sending each $v\in W^v$ on $\overline{v}$ is injective.
\end{Lemma}

\subsection{Classification of indecomposable Kac--Moody matrices}
We now recall the classification of indecomposable Kac--Moody matrices in three types: Cartan matrices (also called Kac--Moody matrix of finite type), affine Kac--Moody matrices and indefinite matrices.  If $(x_i)\in \R^I$, we write $x\geq 0$ if $x_i\geq 0$ for all $i\in I$ and $x>0$ if $x_i>0$, for all $i\in I$.

 \begin{Theorem}\label{t_clas_KM_mat}(\cite[Theorem~4.3]{kac1994infinite})
Let  $A$ be an indecomposable Kac--Moody matrix. Then, one of the following three mutually exclusive possibilities holds for both $A$ and ${^tA}$:\begin{itemize}
\item[\textbf{(Fin)}] $\det(A)\neq 0$, there exists $u\in \R^I$ such that $u>0$ and $A u >0$. 

\item[\textbf{(Aff)}] $\mathrm{corank}(A)=1$; there exists $u\in \R^I$ such that $u>0$ and $Au=0$. For all $v\in \R^I$, $Av\geq 0$ implies $Av=0$.

\item[\textbf{(Ind)}] there exists $u\in \R^I$ such that $u>0$, $Au<0$ and for all $v\in (\R_{\geq 0})^I$, $Av\geq 0$ implies $v=0$. 
\end{itemize}
\end{Theorem}

\begin{Remark}\label{r_clas_KM_mat}
 
 Using Theorem~\ref{t_clas_KM_mat}, we deduce that if $A$ is an indecomposable Kac--Moody matrix, then we have the following properties (and the same hold for $(\alpha_i)_{i\in I}$ with the correct adjustments) \begin{itemize}
\item If $A$ is of finite type, then    $ \sum_{i\in I}\R_{>0} \alpha_i^\vee\cap C^v_f\neq \emptyset$ and $Q^\vee\cap \overline{C^v_f}\subset (\bigoplus_{i\in I}\R_{>0} \alpha_i^\vee)\cup \{0\}$.

\item If $A$ is of affine type, $\mathrm{corank}(A)=1$; $ \A_{\ines}\cap \sum_{i\in I}\R_{>0}\alpha_i^\vee \neq \emptyset$  and $Q^\vee\cap \overline{C^v_f}\subset \A_{\ines}$.

\item If $A$ is of indefinite type,  $(\sum_{i\in I}\R_{>0}\alpha_i^\vee)\cap -C^v_f\neq \emptyset$ and $(\sum_{i\in I}\R_{\geq 0} \alpha_i^\vee)\cap \overline{C^v_f}=\{0\}$.
\end{itemize}

\end{Remark}

We say that $A$ if of \textbf{finite type}\index{finite type KM matrix} or is a \textbf{Cartan matrix}\index{Cartan matrix} if $A$ satisfies \textbf{(Fin)}. We say that $A$ is of \textbf{affine type}\index{affine type KM matrix}  (resp. \textbf{indefinite type}\index{indefinite KM matrix}) if $A$ satisfies \textbf{(Aff)} (resp. \textbf{(Ind)}).

\begin{Proposition}\label{p_char_fd} (see \cite[Proposition 5.9 and Proposition 4.36]{marquis2018introduction})
Suppose that $A$ is indecomposable. The following conditions are equivalent: 

\begin{enumerate}

\item The Kac--Moody matrix $A$ is a Cartan matrix.

\item The vectorial Weyl group $W^v$ is finite.

\item  $\A=\cT$.

\item The set of roots $\Delta$ is finite.

\item The set of real roots $\Phi$ is finite.

\item The set $\Delta^{im}$ is empty.

\item $\ffg(A)$ is a finite dimensional Lie algebra.

\end{enumerate} 
\end{Proposition}

 \begin{Proposition}\label{propClassification des matrices de Kac--Moody de taille 2} (see \cite[Proposition 2.2.12]{hebert2018study})
Let $A=\begin{pmatrix}
2 & -a\\ -b & 2
\end{pmatrix}$, with $a,b\in \Ne$ be an indecomposable Kac--Moody matrix of size $2$. Then: \begin{itemize}

\item if $ab\leq 3$, then  $A$ is of finite type,

\item if $ab=4$, then $A$ is of affine type,

\item otherwise $A$ is of indefinite type.
\end{itemize}
\end{Proposition}

\begin{Definition}\label{d_pos_free}
Let $(x_j)_{j\in J}$ be a family in $\A$. We say that $(x_j)$ is \textbf{positively free}\index{positively free} if we have: \[ \forall (t_j)\in (\R_{\geq 0})^J, \sum_{j\in J} t_j x_j=0 \Rightarrow (t_j)=0.\] A free family of $\A$ is positively free.

We say that $\cS$ is \textbf{positively cofree}\index{positively cofree}  if $(\alpha_i^\vee)_{i\in I}$ is positively cofree.
\end{Definition}

\begin{Lemma}\label{l_char_pf}
Let $V$ be a finite dimensional vector space and $(b_j)_{j\in J}\in V^J$, where $J$ is a finite set. Then $(b_j)_{j\in J}$ is positively-free if and only if there exists $f\in V^*$ such that $f(b_j)>0$, for all $j\in J$, where $V^*$ denotes the dual of $V$. 
\end{Lemma}

\begin{proof}
The sense $\Leftarrow$ is clear. Assume that $(b_j)_{j\in J}$ is positively free. Let $\cC=\sum_{j\in J}\R_{\geq 0} b_j$. By Farkas's Lemma (see for example \cite[Lemma 4.3.3]{hiriart2012fundamentals}), $\cC$ is closed in $V$. Let $\cC'=\conv(\{-b_j\mid j\in J\})$. By Carathéodory's Theorem (see for example \cite[Theorem 1.43]{hiriart2012fundamentals}), $\cC'$ is compact. As $(b_j)$ is positively-free, we have $\cC\cap \cC'=\emptyset$.  By Hahn-Banach's geometric theorem (see for example \cite[Corollary 4.1.3]{hiriart2012fundamentals}), there exists $g\in V^*$ such that $m:=\sup\{g(x)\mid x\in \cC\}<\inf\{g(x)\mid x\in \cC'\}=:M$. As $\cC$ is stable under multiplication by $\R_{\geq 0}$,  we necessarily have $m=0$. We then set $f=-g$ and we have $f(b_j)\geq M>0$ for all $j\in J$, which proves the lemma. 

\end{proof}

\begin{Lemma}\label{l_pos_free}
 \begin{enumerate}
\item If the Kac--Moody matrix $A$ is invertible, then $(\alpha_i^\vee)_{i\in I}$ is free.

\item If all the indecomposable components of  $A$ are of finite type or of indefinite type, then $(\alpha_i^\vee)_{i\in I}$ is positively free.
\end{enumerate}
\end{Lemma}

\begin{proof}
Let $(t_i)\in \R^I$ and  $x=\sum_{i\in I} t_i\alpha_i^\vee$. Then \begin{equation}\label{e_KM_mat}
{^t\!}A.(t_i)_{i\in I}=(\alpha_i(x))_{i\in I}. 
\end{equation} 

(1) Assume that $A$ is invertible. Then if $x=0$, we have $\alpha_i(x)=0$ for all $i\in I$ and thus  ${^t\!}A.(t_i)=0$. Therefore  $(t_i)=0$, which proves that $(\alpha_i^\vee)$ is free.

(2) Let $A_1,\ldots,A_k$ be the indecomposable components of $A$ and denote by $I_j$ the index set of $A_j$, for $j\in \llbracket 1,k\rrbracket$. Let $(t_i)_{i\in I}\in (\R_{\geq 0})^I$ be such that $\sum_{i\in I} t_i\alpha_i^\vee=0$. Let $j\in \llbracket 1,k\rrbracket$. Let $\underline{t}_j=(t_i)_{i\in I_j}$.  Then $A_j \underline{t_j}=(\alpha_i(\sum_{k\in I_j} t_k \alpha_k^\vee))_{i\in I_j}=0$. By Theorem~\ref{t_clas_KM_mat}, we deduce that $(t_k)_{k\in I_j}=0$. Lemma follows.

\end{proof}

\section{Vectorial apartments in the non-reductive case}\label{s_Vectorial apartments}
In this section, we describe the vectorial apartment in the affine case and the indefinite case. 

In~\ref{ss_cT_aff}, we describe the link between Kac--Moody matrices associated to untwisted affine Kac--Moody algebras and loop Lie algebras. We then describe the vectorial apartment of such algebras by using the affine apartment of the underlying semi-simple finite dimensional matrix.

In~\ref{ss_vecApt_indf}, we describe the vectorial apartment in the indefinite case. We give a more precise description in the two dimensional case.

\subsection{The untwisted affine case}\label{ss_cT_aff}

The affine Kac--Moody matrices fall into two types: the untwisted ones, which appear in \cite[5.1 Table Aff 1]{marquis2018introduction} and the twisted ones which appear in \cite[5.1 Table Aff 2 and 3]{marquis2018introduction}. We recall some facts on untwisted ones in this section, and refer  to \cite[5.3]{marquis2018introduction} for the untwisted ones.

Let $\mathring{A}$ be a Cartan matrix. We can associate to it an untwisted affine Kac--Moody matrix $A$ by adding a line and a column. If $\mathring{\ffg}$ and $\ffg$ are the Kac--Moody matrices of $\mathring{A}$ and $A$ respectively (and if the Kac--Moody datum are well-chosen) then $\ffg$ is isomorphic to a double extension of  $\mathring{\ffg}\otimes \C[t,t^{-1}]$. We detail it here and describe the vectorial apartment in this case. We mainly use \cite[XIII]{kumar2002kac}. 

\subsubsection{Underlying finite dimensional simple Lie algebra}\label{sss_undr_fd_la}  Let $\mathring{A}=(a_{i,j})_{1\leq i,j\leq \ell}$ be an indecomposable Cartan matrix. By   \cite[Proposition~8 of V.11]{serre2001complex}, this determines a reduced root system uniquely up to isomorphism. By Theorem~8 and Theorem of the appendix of VI of~\cite{serre2001complex}, there exists a unique semi-simple Lie algebra $\ffgo(\mathring{A})$ having $\mathring{A}$ as a Cartan matrix and this algebra is the Kac--Moody algebra $\ffgo(\mathring{A})$.

Let $\mathring{\mathcal{S}}_\C=(\fho, (\alpha_i)_{i\in \llbracket 1,\ell \rrbracket}, (\alpha_i^\vee)_{i\in \llbracket 1,\ell \rrbracket})$ be the complex minimal free realization of $\mathring{A}$, i.e, $\fho$ is a complex space of dimension $\ell$, $(\alpha_i^\vee)_{i\in \llbracket 1,\ell\rrbracket}$ is a basis of $\fho$, $(\alpha_i)_{i\in \llbracket 1,\ell\rrbracket}$ is a basis of $\fho^*$ and $\alpha_i(\alpha_j^\vee)=a_{j,i}$, for all $i,j\in I$. Let $\ffgo=\ffgo_{\mathring{\mathcal{S}}_\C}$: $\ffgo$ is the Kac--Moody algebra associated with $\mathring{\cS}$, for any Kac--Moody datum $\mathring{\cS}$ associated with $\mathring{A}$ whose character lattice has rank the size of $\mathring{A}$. We denote by $\mathring{e_i},\mathring{f_i}$, $i\in \llbracket 1,\ell \rrbracket$ the Chevalley generators of $\ffgo$. Let $\mathring{\Phi}\subset \fho^*$ be the root system of $(\ffgo,\mathring{\mathcal{S}}_\C)$.

  Let $\mathring{X}=\bigoplus_{i=1}^\ell \Z\alpha_i$ and $\mathring{Y}=\Hom_\Z(\mathring{X},\Z)$ (we regard the $\alpha_i^\vee$, $i\in \llbracket 1,\ell\rrbracket$ as elements of $\mathring{Y}$). Let $c,d$ be symbols and $Y=\mathring{Y}\oplus \Z c \oplus \Z d$.  Let $X=\Hom_\Z(Y,\Z)$.  Let $\fh=Y\otimes_\Z \C$ and $\fh^*=X\otimes_\Z \C$.
We embed $\mathring{X}$ in $X$  by requiring $\lambda(\Z c+\Z d)=\{0\}$ for all $\lambda\in \mathring{Y}$. Let $\delta\in X$ defined by $\delta(d)=1$ and $\delta(Y\oplus \Z c)=0$.

 By  \cite[VI, Proposition 25 (i)]{bourbaki1981elements}, $\mathring{\Phi}^+$ admits a unique highest root: this is a root $\theta=\sum_{i=1}^\ell a_i \alpha_i\in \mathring{\Phi}^+$ such that for all $\alpha=\sum_{i=1}^\ell n_i\alpha_i\in \mathring{\Phi}^+$, we have $n_i\leq a_i$. In particular $a_i\in \N^*$ for all $i\in \llbracket 1,\ell \rrbracket$. 

\subsubsection{Realization of $\ffg$ as a loop Lie algebra}\label{sss_aff_la}

The system $\mathring{\mathcal{S}^\vee}:=({^t\mathring{A}},\mathring{Y},\mathring{X},(\alpha_i^\vee)_{i\in \llbracket 1,\ell \rrbracket},(\alpha_i)_{i\in \llbracket 1,\ell \rrbracket})$ is a Kac--Moody datum.
 Let $\theta^\vee\in \fho_\Z^*$ be the highest root of $\mathring{\Phi}^{\vee+}$ (where $\mathring{\Phi^\vee}$ is the root system of $\mathring{\mathcal{S}^\vee}$). We set $\alpha_0=\delta-\theta$ and $\alpha_0^\vee=c-\theta^\vee$. 

We set $a_{0,0}=2=\theta(\theta^\vee)$ and for $j\in \llbracket 1,\ell \rrbracket$, $a_{0,j}=-\alpha_j(\theta^\vee)$ and $a_{j,0}=-\theta(\alpha_j^\vee)$. Let $A=(a_{i,j})_{0\leq i,j\leq l}$.

 The matrix $A$ is then  a Kac--Moody matrix of untwisted affine type.

Let $I=\llbracket 0,\ell \rrbracket$.  The system $\mathcal{S}=(A,X,Y,(\alpha_i)_{i\in I},(\alpha_i^\vee)_{i\in I})$ is a Kac--Moody datum.

Let $\ffg=(\C[t,t^{-1}]\otimes_\C \ffgo)\oplus \C c\oplus \C d$, where $t$ is an indeterminate. Let $\langle,\rangle$ be the invariant (symmetric nondegenerate) bilinear form on $\ffgo$ defined in  \cite[Theorem~1.5.4]{kumar2002kac} (the properties of $\langle,\rangle$ are important to prove that $\ffg$ is isomorphic to the Kac--Moody algebra of $\mathcal{S}$ but as we admit this result, we will not explicitly use them). We equip $\ffg$ with the bracket: \[[t^m\otimes x+\lambda c+ \mu d,t^{m'}\otimes x'+\lambda' c+\mu' d]=t^{m+m'}\otimes [x,x']+\mu m't^{m'}\otimes x'-\mu'mt^m \otimes x +m\delta_{m,-m'}\langle x,x'\rangle c,\] 
for $\lambda,\mu,\lambda',\mu'\in \C$, $m,m'\in \Z$, $x,x'\in \ffgo$.

Let $\Delta=\{j\delta\mid j\in \Z\backslash \{0\}\}\cup \{j\delta+\beta\mid j\in \Z,\ \beta\in \mathring{\Phi}\}$. If $\alpha\in \fh^*$, we set $\ffg_\alpha=\{x\in \ffg\mid [h,x]=\alpha(h)x\ \forall h\in \fh\}$. By \cite[(5.9)]{marquis2018introduction}, we have:

\begin{Lemma}\label{l_rt_dec_aff}
We have $\ffg=\fh\oplus \bigoplus_{\alpha\in \Delta}\ffg_\alpha$. Moreover if $j\in \Z\backslash \{0\}$, $\ffg_{j\delta}=t^j\otimes \fho$ and if $j\in \Z$ and $\beta\in \mathring{\Phi}$, $\ffg_{\beta+j\delta}=t^j\otimes \ffgo_\beta$, where $\ffgo_\beta=\{x\in \ffgo\mid[h,x]=\beta(h)x\ \forall h\in \fho\}$.
\end{Lemma}

Let $\mathring{\omega}$ be the Chevalley involution of $\ffgo$ (see for example \cite[after Lemma 3.22]{marquis2018introduction}). Let $x_0\in \ffgo_\Theta$ such that $\langle x_0,\mathring{\omega}(x_0)\rangle=-1$. Let $E_0=-t\otimes \mathring{\omega}(x_0)$ and $F_0=-t\otimes x_0$.

\begin{Theorem}\label{t_af_la_loop} (\cite[Theorem~1.3.1.3]{kumar2002kac})
There exists a unique Lie algebra isomorphism $\Upsilon:\ffg_{\mathcal{S}} \overset{\sim}{\rightarrow} \ffg $ fixing $\fh$ and mapping $e_0$ to $E_0$, $f_0$ to $F_0$, $e_i$ to $\mathring{e_i}$ and $f_i$ to $\mathring{f_i}$ for all $i\in \llbracket 1,\ell \rrbracket$.
\end{Theorem}

\begin{Corollary}\label{c_rt_aff}
The set $\Delta$ is the set of  roots of $\ffg$. We have $\Delta_+=\{j\delta+\beta\mid j\in \Ne,\beta\in \mathring{\Phi}\cup \{0\}\}\cup \mathring{\Phi}^+$ and $\Delta^{im}=(\Z\backslash\{0\}) \delta$. Moreover, if $j\in \Z\backslash\{0\}$, $\dim \ffg_{j\delta}=\ell$.
\end{Corollary}

\subsubsection{Apartment of $\mathring{A}$ viewed as a subspace of $\A_{\mathrm{es}}$}
In this subsection, we show that we can see the affine apartment of $\mathring{A}$ as an affine hyperplane of the essentialization $\A_{\mathrm{es}}$ of $\A$. This is inspired by  \cite[Section 6]{kac1994infinite}.

Let $\A=Y\otimes \R$. Let $\A^1=\{x\in \A\mid\delta(x)=1\}$.  Then $\A^1$ is an affine subspace of $\A$ whose direction contains $\R c=\A_{\ines}$.
 
We set $\widetilde{\A}=\A^1/ \R c\subset \A_{\mathrm{es}}$.  We  can consider each element of $\Phi$ as a linear form on $\A_{\mathrm{es}}$ and we consider $W^v$ as a subgroup of $\mathrm{Aut}(\A_{\mathrm{es}})$, which is possible by Lemma~\ref{l_res_ess}. Let $\mathring{\A}=\bigoplus_{i=1}^\ell \R\alpha_i^\vee\subset \A$. The map $\iota:\mathring{\A}\rightarrow \A_{\mathrm{es}}$ defined by $\iota(\sum_{i=1}^\ell  x_i \alpha_i^\vee)=\sum_{i=1}^\ell  x_i\alpha_i^\vee +\R c$ is injective and  we consider $\mathring{\A}$ as a subspace of $\A_{\mathrm{es}}$ through this embedding. 

Let $\psi:\mathring{\A}\rightarrow \widetilde{\A}$ be defined by $\psi(x)=x+d$, for $x\in \mathring{\A}$. Then $\psi$ is an isomorphism of affine spaces. 

We equip $\mathring{\A}$ with its structure of affine apartment. The walls of $\mathring{\A}$ are the $\mathring{\alpha}^{-1}(\{k\})$ such that $\mathring{\alpha}\in \mathring{\Phi}$ and $k\in \Z$. Let $\mathring{Q}^\vee_\Z=\bigoplus_{i=1}^\ell  \Z\alpha_i^\vee$ and $\mathring{W}=\mathring{W^v}\ltimes \mathring{Q}^\vee_\Z$  be the affine Weyl group of $\mathring{\A}$. Let $\mathcal{M}_0=\{\alpha^{-1}(\{0\})\mid\alpha\in \Phi\}$ be the set of vectorial walls of $\A_{\mathrm{es}}$ and $\widetilde{\cM}=\{M\cap \widetilde{\A}\mid M\in \mathcal{M}\}$. Then by Corollary~\ref{c_rt_aff}, if  $\mathring{\cM}:=\{\mathring{\alpha}^{-1}(\{k\})\mid (\mathring{\alpha},k)\in \mathring{\Phi}\times \Z\}$ is the set of walls of $\mathring{\A}$, then $\mathring{\cM}=\psi^{-1}(\widetilde{\cM})$. 

 Let $\mathring{C^v_f}=\{x\in \mathring{\A}\mid  \alpha_i(x)>0,\ \forall i\in \llbracket 1,\ell \rrbracket\}$ and $\mathring{C}_0^+=\mathring{C}^v_f\cap \theta^{-1}(]0,1[)$. By  \cite[VI Proposition~5]{bourbaki1981elements},  $\overline{\mathring{C}_0^+}$ is a fundamental domain for the action of $\mathring{W}$ on $\mathring{\A}$ and is a connected component of $\mathring{\A}\backslash \bigcup_{M\in \mathcal{M}_0}M$. 

Let $C^v_{f,\mathrm{es}}\subset \A_{\mathrm{es}}$ be the image of $C^v_f$ in $\A_{\mathrm{es}}$. Let $\mathcal{F}(C^v_{f,\mathrm{es}})$   (resp. $\mathcal{F}(\mathring{C}_0^+)$) be  the poset consisting of the  faces of $\overline{C^v_{f,\mathrm{es}}}$ (resp. $\overline{\mathring{C}_0^+}$).

\begin{Lemma}\label{l_f_Ch_af_red} (see \cite[Lemma 2.3.5]{hebert2018study})
The map $\Lambda:\sF(C^v_{f,\mathrm{es}})\rightarrow \sF(\mathring{C}_0^+)$ defined by $\Lambda(F^v)=\psi^{-1}\big(F^v\cap \widetilde{\A}\big)$ is well-defined and is an isomorphism of posets.
\end{Lemma}

As $\delta(Q^\vee_\R)=\{0\}$, we have $\delta\circ w=\delta$ for all $w\in W^v$. Therefore, $W^v$ acts on $\widetilde{\A}$. We deduce an action of $W^v$ on $\mathring{\A}$  through $\psi$. If $w\in W^v$, we denote by $\mathring{w}=\psi^{-1}\circ w\circ \psi:\mathring{\A}\rightarrow \mathring{\A}$ the corresponding affine automorphism.

\begin{Lemma}\label{l_ref_af_red} (see \cite[Lemma 2.3.6]{hebert2018study})\begin{enumerate}
\item Let $i\in \llbracket 1,\ell \rrbracket$. Then for all $x\in \mathring{\A}$, $\mathring{r}_i(x)=x-\alpha_i(x)\alpha_i^\vee$.
\item For all $x\in \mathring{\A}$, $\mathring{r}_0(x)=x-(\theta(x)-1)\theta^\vee$.
\end{enumerate}
\end{Lemma}

Let $\mathcal{F}^+(\A_{\mathrm{es}})$ (resp. $\sF(\mathring{\A})$) be the simplicial complex consisting of the positive vectorial faces of $\A_{\mathrm{es}}$ (resp. of the faces of $\mathring{\A}$).
 
\begin{Proposition}\label{p_red_apt_af}(see for example \cite[Proposition 2.3.7]{hebert2018study})
\begin{enumerate}
\item  The morphism $\Xi:W^v\rightarrow \mathrm{Aut}(\mathring{\A})$ mapping each $w\in W^v$ to $\mathring{w}$ is injective and its image is $\mathring{W}=\mathring{W^v}\ltimes \mathring{Q^\vee_\Z}$.

\item  The map $\Lambda:\mathcal{F}^+(\A_{\mathrm{es}})\rightarrow \sF(\mathring{\A})$ sending each $F^v$ on $\psi^{-1}(F^v\cap \widetilde{\A})$ is well-defined and is an isomorphism compatible with the actions of $W^v$ and $\mathring{W}$: for all $w\in W^v$ and $F^v\in \mathcal{F}^+(\A)$ and $w\in W^v$, $\Lambda(w.F^v)=\mathring{w}. \Lambda (F^v)$.
\end{enumerate}
\end{Proposition}

\begin{Corollary}(see for example \cite[Corollary 2.3.8]{hebert2018study})\label{c_TC_aff}
The Tits cone $\cT$ of $\A$ is $\delta^{-1}(\R^*_+)\cup \A_{\ines}$. In particular, $\overline{\cT}=\delta^{-1}(\R_+)$.
\end{Corollary}

\begin{figure}[h]
\centering
\includegraphics[scale=0.4]{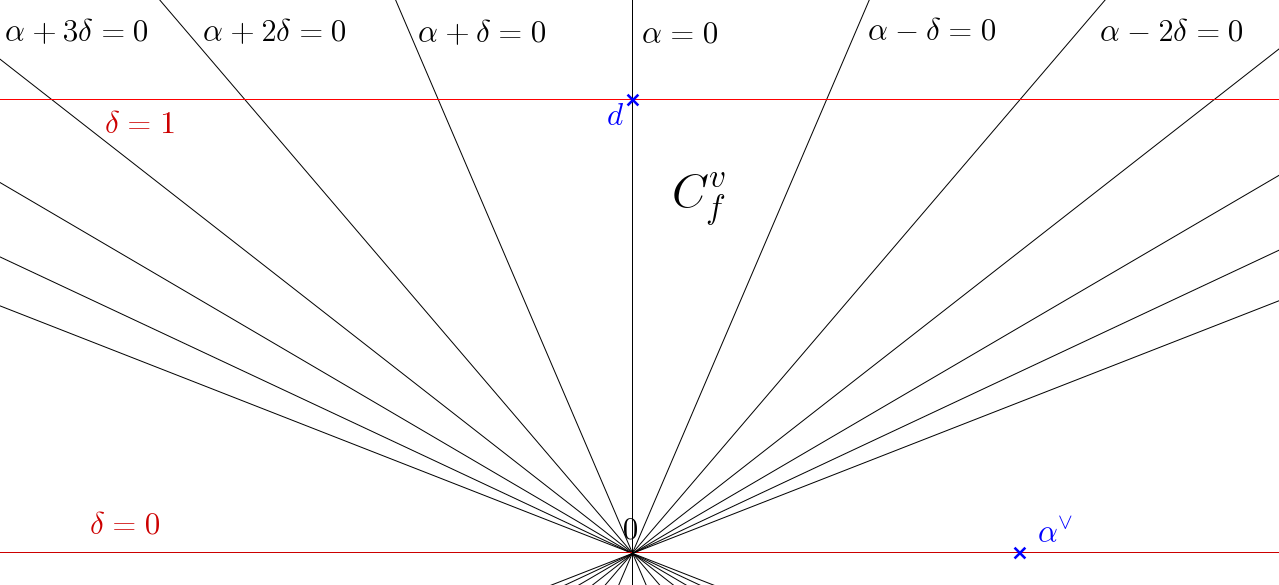}
\caption{Vectorial apartment of $\widetilde{\mathrm{SL}_2}$ with thirteen walls}
\end{figure}

\begin{figure}[h]
\centering
\includegraphics[scale=0.4]{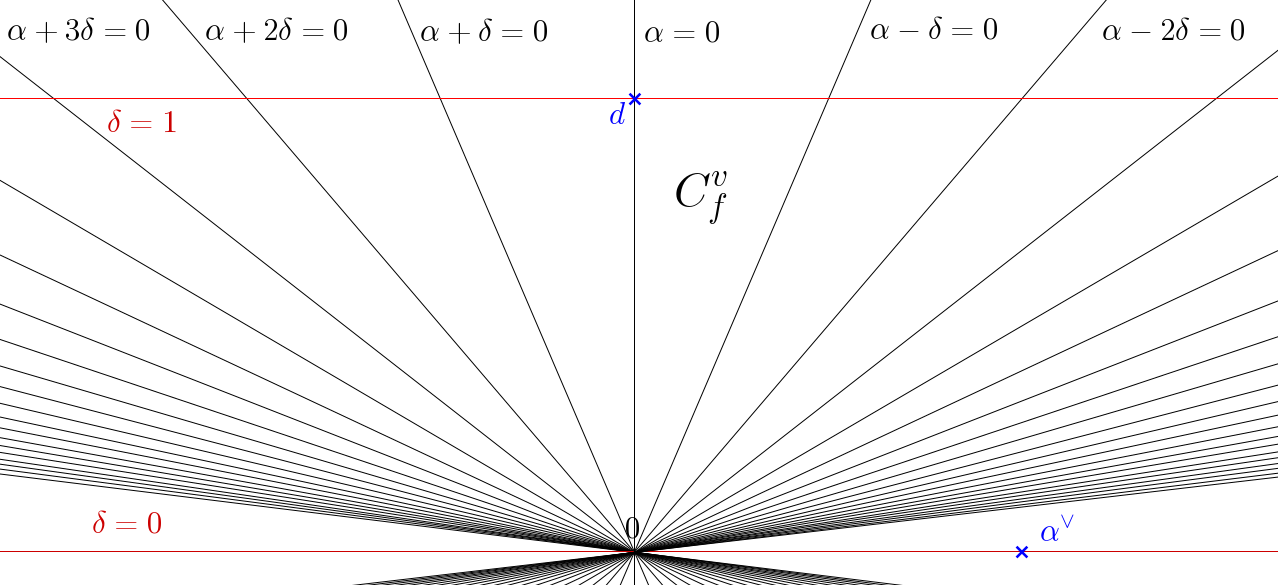}
\caption{Vectorial apartment of $\widetilde{\mathrm{SL}_2}$ with forty-one walls}
\end{figure}

\subsection{Vectorial apartment in the indefinite case}\label{ss_vecApt_indf}
We now describe the vectorial apartment in the indefinite case. Let $\cS$ be a Kac--Moody datum associated to an indefinite Kac--Moody matrix. The following proposition gathers well-known facts on imaginary roots.

\begin{Proposition}(\cite[Proposition 2.3.9]{hebert2018study})\label{p_T_cne_indef}
Suppose that $A$ is an indefinite Kac--Moody matrix. Then: \begin{enumerate}
\item $\overline{\cT}=\{x\in \A\mid\alpha(x)\geq 0,\ \forall \alpha\in \Delta^{im}_+\}$, and in particular, $\mathring{\cT}\subset \{x\in \A\mid\alpha(x)> 0,\ \forall \alpha\in \Delta^{im}_+\}$,

\item for all $i\in I$, $\alpha_i^\vee\in \A\backslash (\overline{\cT}\cup \overline{-\cT})$,

\item $\overline{\cT}$ is not a half-space of $\A$,

\item for all $\alpha\in \Delta^{im}_+$, $\Z_{>0} \alpha \subsetneq \Delta^{im}_+$. 
\end{enumerate}
\end{Proposition}

We now study the case where $A$ is of size $2$. We set $I=\{1,2\}$. We write  $A=A(a,b)=\begin{pmatrix}
2 & -a\\ -b &2
\end{pmatrix}$ with $a,b\in \Ne$. We have $ab\geq 5$ by Proposition~\ref{propClassification des matrices de Kac--Moody de taille 2}. In particular, $A$ is invertible and the minimal free realization $\A$ of $A$ over $\R$ has dimension $2$. The vectorial faces of $\A$ are $\{0\}$, the vectorial panels and the vectorial chambers. Except $\{0\}$, these faces are spherical and by Proposition~\ref{p_T_cone}, $\mathring{\cT}=\cT\backslash \{0\}$.

\begin{Proposition} (see Figure~\ref{figVectorial apartment of A(1,5)} and Figure~\ref{figVectorial apartment of A(3,3)})\label{propDescription de l'appart standard en indéf de dim 2}
The set  $\A\backslash (\cT\cup-\cT)$ has two connected components. We denote them by $\Gamma_1$ and $\Gamma_2$. We have $\Gamma_1=-\Gamma_2$. If $i\in \{1,2\}$, then $\Gamma_i$ is a convex cone, has nonempty interior,   does not meet any vectorial wall of $\A$ and $\Gamma_i\cup\{0\}$ is closed. Moreover, $\cT\cup(-\cT)$ is delimited by the eigenspaces of $r_1r_2$. 
\end{Proposition}

\begin{proof}
The first part is  \cite[Proposition 2.3.10]{hebert2018study}. Let $L_1$ and $L_2$ be the two lines delimiting $\cT\cup (-\cT)$.  Let $i\in \{1,2\}$. Then  $r_i$ fixes pointwise a line $L_i'$ such that $L_i'\setminus \{0\}$ is contained in $\mathring{\cT}$. As $r_i$ stabilizes exactly two lines, $r_i$ stabilizes neither $L_1$ nor $L_2$. As $r_i$ stabilizes $\cT$, we deduce that $(r_i(L_1),r_i(L_2))=(L_{2},L_1)$. Consequently $r_1r_2$ stabilizes $L_1$ and $L_2$, which are thus exactly the eigenspaces of $r_1r_2$.
\end{proof}

\begin{figure}
\centering
\includegraphics[scale=0.4]{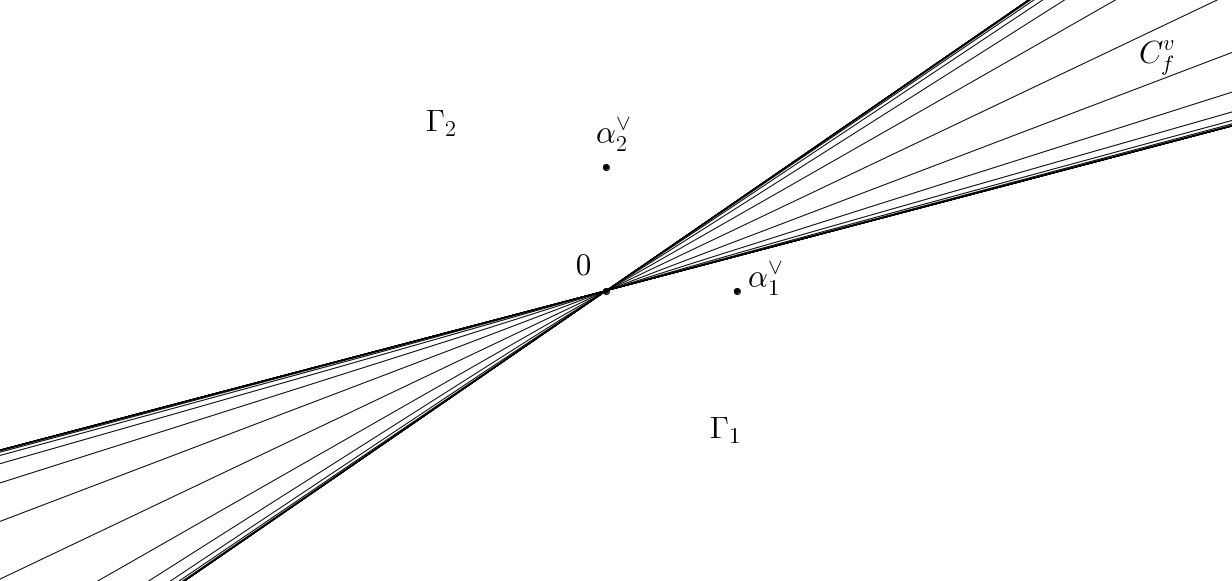}
\caption{Vectorial apartment of A(1,5) }\label{figVectorial apartment of A(1,5)}
\end{figure}

\begin{figure}
\centering
\includegraphics[scale=0.4]{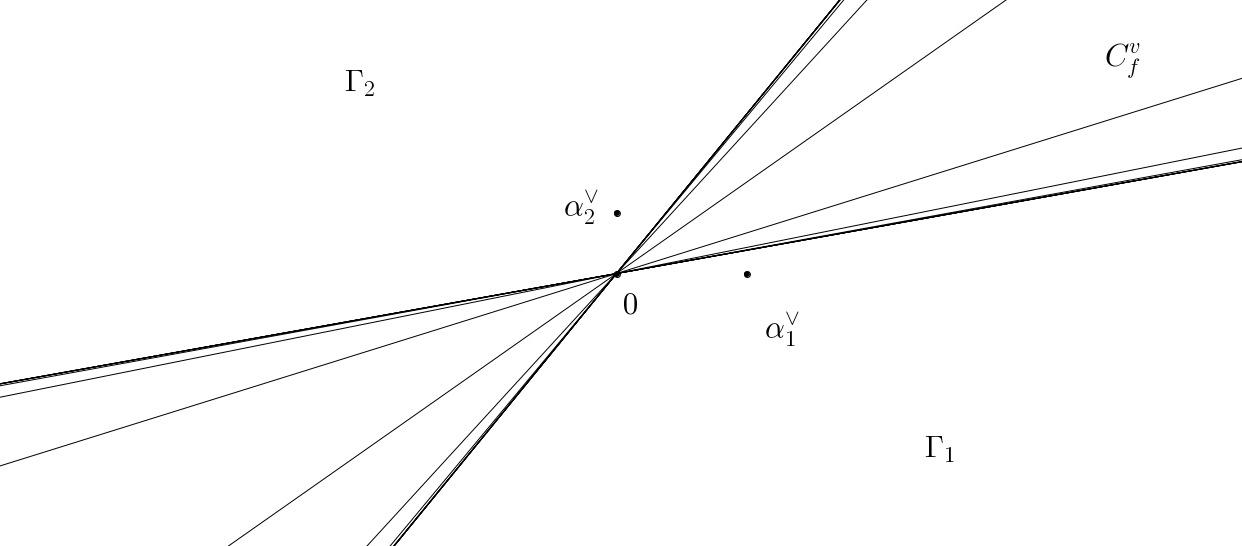}
\caption{Vectorial apartment of A(3,3)}\label{figVectorial apartment of A(3,3)}
\end{figure}

\chapter{Split Kac--Moody groups, à la Tits and à la Mathieu}\label{C_splt_KM_grps}

Let $\cS$ be a free  Kac--Moody datum. In Chapter~\ref{C_KM_alg}, we introduced a Kac--Moody algebra to $\cS$.  In \cite{tits1987uniqueness}, Tits defined the constructive Tits functor associated   to this datum. It is a group functor  $\fG_\cS$ from the category of rings to the category of groups, which integrates $\ffg$ in some sense. Moreover, Tits also defined a series of axioms that a group functor needs to satisfy to be called a ``minimal Kac--Moody group functor'' - minimal meaning that it integrates only the real roots of $(\ffg,\fh)$. He proved that $\fG_\cS$ satisfies its axioms over the fields and that if  $\cK$, the value of a functor satisfying its axioms is unique and isomorphic to $\fG_\cS(\cK)$.  If $\cR$ is a ring, then  $\fG_{\cS}(\cR)$ is generated by a torus $\fT_\cS(\cR)$, isomorphic to $(\cR^\times)^{\dim \fh}$, and by the root groups $x_\alpha(\cR)$, $\alpha\in \Phi$, which are all isomorphic to $(\cR,+)$.   

On the other hand, Mathieu defined in \cite{mathieu1989construction} a completed Kac--Moody group functor $\fG_\cS^{ma+}=\fG^{ma+}$ from the category of rings to the category of groups. It integrates the real roots and the positive imaginary roots of $(\ffg,\fh)$ and there is a group functor morphism from $\fG_\cS$ to $\fG^{ma+}$, which is an injection on every field. It is used by Marquis to define a ``minimal group functor'' which satisfies Tits's axioms (on all rings and not only on fields).

Our objective is to study minimal Kac--Moody groups. However, certain commutation relations are missing in them, compared to the reductive case. More precisely, assume that $\ffg$ is finite dimensional and that $\cK$ is a field.  Then by \cite[Lemma 15]{steinberg2016lectures}, for all $\alpha,\beta\in \Delta=\Phi$ such that $\alpha\neq -\beta$, the commutator $[x_\alpha(\cK),x_{\beta}(\cK)]$ is contained in $\prod_{\gamma\in (\Z_{\geq 1}\alpha+\Z_{\geq 1} \beta)\cap \Phi} x_{\gamma}(\cK)$. Consider now the Kac--Moody group $G=\mathrm{SL}_2(\cK[t,t^{-1}])$, where $\cK$ is a field. Its root system is $\{\pm \mathring{\alpha}+k\delta\mid k\in \Z\}$, for the notation of Subsection~\ref{ss_cT_aff}. Then by \cite[Example 7.85 \& (7.47)]{marquis2018introduction}, we can take $x_{\mathring{\alpha}+k\delta}(r)=\begin{psmallmatrix} 1 & rt^k\\ 0 & 1\end{psmallmatrix}$ and $x_{-\mathring{\alpha}+k\delta}(r)=\begin{psmallmatrix} 1 & 0\\ r t^k & 1\end{psmallmatrix}$ for $r\in \cK$. Then for $r,s\in \cK$ and $k,\ell\in \Z$, we have:  \[\begin{psmallmatrix} 1 & 0 \\ rt^k & 0\end{psmallmatrix}\begin{psmallmatrix} 1 & s t^\ell \\ 0 & 1\end{psmallmatrix}=\begin{psmallmatrix} 1  &st^\ell/(1+rst^{k+\ell})\\ 0 & 1\end{psmallmatrix}\begin{psmallmatrix} 1/(1+rs t^{k+\ell})  & 0 \\ 0 & 1+rs t^{k+\ell}\end{psmallmatrix}\begin{psmallmatrix} 1 & 0 \\ rt^k/(1+rst^{k+\ell}) & 1\end{psmallmatrix}.\] If $(k+\ell)rs\neq 0$, then the three matrices on the right hand side do not belong to $G$. 

In order to handle this kind of elements, it is convenient to embed $G$ in $\mathrm{SL}_2(\cK(\!(t)\!))$ or $\mathrm{SL}_2(\cK(\!(t^{-1})\!))$. More generally, if $\fG$ is a general Kac--Moody group, we use  its Mathieu's completions $\fG^{ma+}$ and $\fG^{ma-}$ (where $p$ and $n$ stand for positive and negative), which are $\mathrm{SL}_2(\cK(\!(t)\!))$ and $\mathrm{SL}_2(\cK(\!(t^{-1})\!))$ in the case of $\mathrm{SL}_2(\cK[t,t^{-1}])$. 

The definition of $\fG^{ma+}$ is complicated over arbitrary rings. However, over fields, a simpler description, by generators and relations was provided by Marquis in \cite[8.5]{marquis2018introduction}. Let $\fU^+=\langle x_\alpha \mid \alpha\in \Phi_+\rangle$ and $\fU^-=\langle x_\alpha\mid \alpha\in \Phi_-\rangle$. The group functor $\fU^+$ is completed in a group functor $\fU^{ma+}$ and then if $\cK$ is a field, $\fG^{ma+}(\cK)$ is defined   as a quotient of the free product $\fG(\cK)*\fU^{ma+}(\cK)$ by certain relations (see Definition~\ref{d_sch_th_KM}). 

An advantage of $\fU^{ma+}$ compared to $\fU^+$ is that we can use the ``coordinates''  $X_\alpha, \alpha\in \Delta_+$ introduced by Rousseau in \cite{rousseau2016groupes} to describe its elements (Rousseau does not define the $X_\alpha$ there, but he defines the ``twisted exponentials'', from which the  $X_\alpha$ are defined). The $X_\alpha$  generalize the $x_\alpha$ used in the reductive case.   Let $\cR$ be a ring. For $\alpha\in \Delta_+$, we have a map $X_\alpha:\ffg_{\alpha,\Z}\otimes \cR\rightarrow \fU^{ma+}(\cR)$ (depending on many choices). Then if we fix an order on $\Delta_+$ compatible with the height, we can express every element of $\fU^{ma+}(\cR)$  uniquely as a product $\prod_{\alpha\in \Delta_+} X_\alpha(\underline{u_\alpha})$, where $(\underline{u_\alpha})\in \prod_{\alpha\in \Delta_+} \ffg_{\alpha,\Z}\otimes \cR$. These coordinates are probably useful in the study of Kac--Moody groups over arbitrary rings or fields in general. For our purpose,  where $\cR=\cF$ is a valued field and $G=\fG_\cS(\cF)$, they enable to define the fixator $U^+$ and $U^-$-components of the fixator of a subset of $\A$, for the action of $G$ on the masure.

In Section~\ref{s_min_KM_gp}, we define the constructive Tits functor associated to a Kac--Moody datum. Its values on any field $\cK$ is the ``minimal Kac--Moody group over $\cK$``.

In Section~\ref{s_gp_Upma}, we define and study $\fU^{ma+}$, which is a completion of $\fU^+$, and which is a standard unipotent subgroup of Mathieu's completed group.

In Section~\ref{s_cplt_KM_gp}, we introduce Marquis's scheme theoretic and Mathieu's completion of the minimal Kac--Moody group.

\section{Minimal Kac--Moody group}\label{s_min_KM_gp}

In this section, we introduce and study the constructive Tits functor $\fG_\cS$ associated by Tits with each Kac--Moody datum $\cS$. It is defined by generators and relations.  For a field $\cK$, $\fG_\cS(\cK)$ is what we call the minimal split Kac--Mody group over $\cK$ associated with $\cS$. For an arbitrary ring $\cR$, the ``right'' definition of the minimal Kac--Moody group over $\cR$ associated with $\cS$ is not $\fG_{\cS}(\cR)$. We give Marquis's definition in Subsection~\ref{ss_min_KM_gps}, using Mathieu's completion.

\subsection{Nilpotent sets and prenilpotent pairs of roots}\label{subsubNilpotent sets}
In this subsection, we define the notion of nilpotency for sets of roots. This notion is important to define the constructive Tits functor: it is used to define certain commutation relations satisfied by this group.

\begin{Definition}\label{d_closed}
Let $\Psi\subset \Delta_+$. We say that $\Psi$ is \textbf{closed} (resp. \textbf{strongly closed}) if for every $\alpha,\beta\in \Psi$ such that $\alpha+\beta\in \Delta_+$ (resp. for all $p,q\in \Z_{\geq 1}$ such that $p\alpha+q\beta\in \Delta$), we have $\alpha+\beta\in \Psi$ (resp. $p\alpha+q\beta\in \Psi$).

 Let $\Psi'\subset \Psi\subset \Delta_+$ be closed subset of roots. We say that $\Psi'$ is \textbf{coclosed} in $\Psi$ (resp. \textbf{strongly coclosed}) if $\Psi\setminus \Psi'$ is closed (resp. if $\Psi$ is strongly closed and $\Psi\setminus \Psi'$ is strongly closed). 
 
 We also  call $\Psi'$ an \textbf{ideal in $\Psi$} (resp. a \textbf{strong ideal in $\Psi$}) if $\alpha+\beta\in \Psi'$ (resp. $p\alpha+q\beta\in \Psi'$) whenever $\alpha\in \Psi$, $\beta\in \Psi'$ and $\alpha+\beta\in \Psi$ (resp. $p,q\in \Z_{\geq 1}$, $\alpha\in \Psi$, $\beta\in \Psi'$ and $p\alpha+q\beta\in \Delta_+$).
\end{Definition}

Rousseau and Marquis use  different definitions of closeness and ideals in \cite{rousseau2016groupes} and \cite{marquis2018introduction}. We add the adjective ``strong'' to refer to the notion used by Rousseau. We do not know whether it is really different from the one of Marquis.

\begin{Definition}
Let  $\Psi$  be a subset of $\Phi$ (or $\Delta$).  The set $\Psi$ is said to be \textbf{prenilpotent}\index{prenilpotent} if there exists $w,w'\in W^v$ such that $w.\Psi\subset \Delta_+$ and $w'.\Psi\subset \Delta_-$. Then $\Psi$ is finite since $w.\Psi\subset \Inv(w'w^{-1})$. We say that $\Psi$ is \textbf{nilpotent} if it is prenilpotent and closed. A set $\Psi$ is prenilpotent if and only if $-\Psi$ is prenilpotent.
\end{Definition}

If $\{\alpha,\beta\}\subset \Phi$ is a prenilpotent pair, we denote by $]\alpha,\beta[_\N$\index[notation]{z@$]\alpha,\beta[_{\N}$} the finite set $(\Ne \alpha+\Ne\beta) \cap \Phi$ and we set  $[\alpha,\beta]_\N=]\alpha,\beta[_\N\cup \{\alpha,\beta\}$\index[notation]{z@$[\alpha,\beta]_{\N}$}. For all $\alpha\in \Phi$, the pair $\{\alpha,-\alpha\}$ is not prenilpotent.

In the remaining part of this subsection, we give criteria for pairs of roots to be prenilpotent in order to illustrate this notion. 

\begin{Lemma}\label{lem5.4.2 de remy} (see  \cite[Lemme 5.4.2]{remy2002groupes})
\begin{enumerate}
\item Let $\Psi\subset \Phi$. Then $\Psi$ is prenilpotent if and only if $\cT\cap\bigcap_{\alpha\in \Psi}\alpha^{-1}(\R_+)$ and $-\cT\cap\bigcap_{\alpha\in \Psi}\alpha^{-1}(\R_+)$ has nonempty interior.
\item Let $\alpha,\beta\in \Phi$ be such that $\alpha\neq -\beta$. Then at least one pair of $\{\{\alpha,\beta\},\{\alpha,-\beta\}\}$ is prenilpotent.
\end{enumerate}

\end{Lemma}

\begin{Proposition}\label{propPaires prénilpotentes en réductif} (see \cite[Proposition 3.3.2]{hebert2018study})
Suppose that $A$ is of finite type. Then a pair $\{\alpha,\beta\}\subset \Phi$ is prenilpotent if and only if $\alpha\neq -\beta$.
\end{Proposition}

\begin{Proposition}\label{propPaires prénilpotentes en affine} (see \cite[Proposition 3.3.3]{hebert2018study})
Suppose that $A$ is of untwisted affine type. Let $\alpha,\beta\in \Phi$ be  such that $\alpha\neq \beta$. We write $\alpha=\mathring{\alpha}+k\delta$, $\beta=\mathring{\beta}+ \ell \delta$, with $\mathring{\alpha},\mathring{\beta}\in \mathring{\Phi}$, $k,\ell \in \Z$, which is possible by Corollary~\ref{c_rt_aff}. Then $\{\alpha,\beta\}$ is prenilpotent if and only if $\mathring{\alpha}\neq \mathring{-\beta}$.
\end{Proposition}

For the next proposition, we assume that the matrix $A$ is indefinite of size $2$. Let $\Gamma_1$, $\Gamma_2$ be the connected components of $\A\backslash (\cT\cup -\cT)$ (see Proposition~\ref{propDescription de l'appart standard en indéf de dim 2}).

\begin{Proposition}\label{propPaires prénilpotentes en dim 2} (see \cite[Proposition 3.3.5]{hebert2018study})
Suppose that $A$ is indefinite of size $2$. Let $\alpha,\beta\in \Phi$. Let $\Gamma_1$, $\Gamma_2$ be as in Proposition~\ref{propDescription de l'appart standard en indéf de dim 2}. Then: \begin{enumerate}
\item~\label{itCaractérisation des paires prénilp, indéfini} The pair $\{\alpha,\beta\}$ is prenilpotent if and only if there exists $i\in \{1,2\}$ such that $\alpha^{-1}(\R_+)\cap \beta^{-1}(\R_+)\supset \Gamma_i$. 

\item Exactly one pair of $\{\{\alpha,\beta\},\{\alpha,-\beta\}\}$ is prenilpotent.
 \end{enumerate}
\end{Proposition}

\begin{figure}
\centering
\includegraphics[scale=0.4]{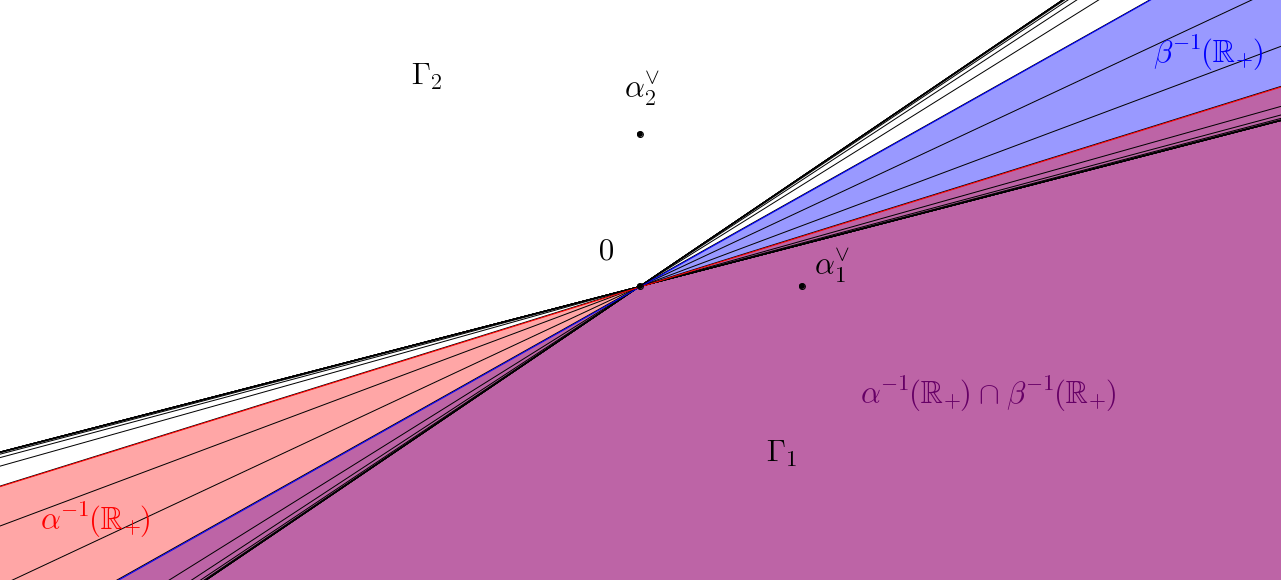}
\caption{Exemple of a prenilpotent pair $\{\alpha,\beta\}$}
\end{figure}

\begin{figure}
\centering
\includegraphics[scale=0.4]{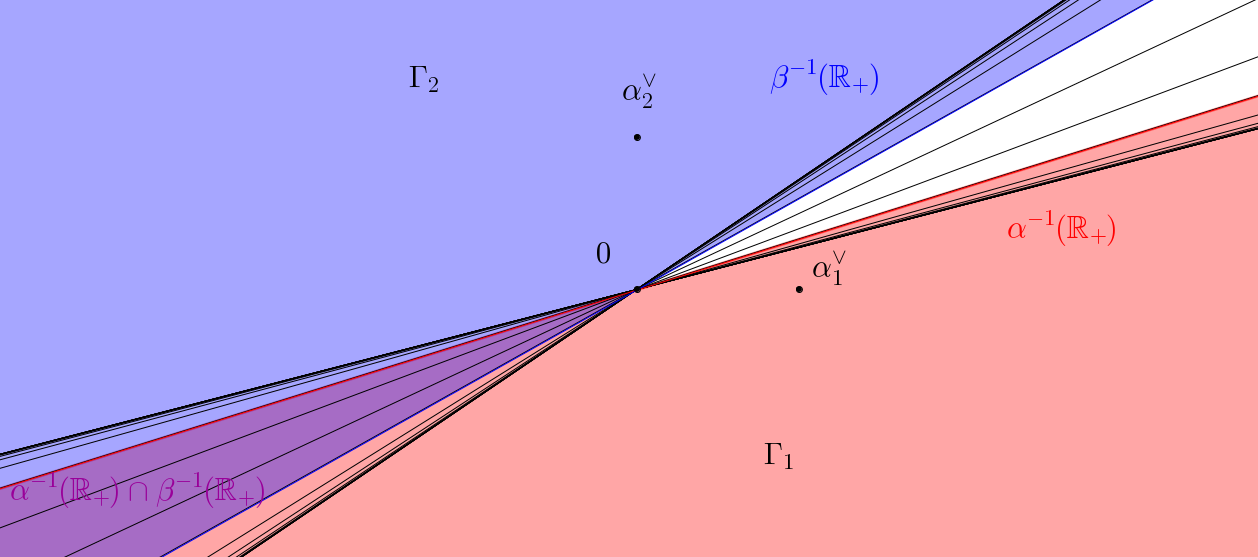}
\caption{Exemple of a non prenilpotent pair $\{\alpha,\beta\}$}
\end{figure}

\subsection{Enveloping algebras of $\ffg$}\label{ss_env_alg}
Let $\ffg=\ffg_\cS$\index[notation]{g@$\ffg$} be the Kac--Moody Lie algebra over $\C$ associated with $\cS$ (see  Subsection~\ref{ss_KM_alg}) and $\cU_\C(\ffg)$\index[notation]{U@$\cU_\C(\ffg)$} be its enveloping algebra.
In this subsection, we study   $\cU_\C(\ffg)$.  The group $\fU^{ma+}(\C)$,  will be defined as a subgroup of a completion of $\cU_\C(\ffg)$. As we want to define $\fU^{ma+}(\cR)$, for any ring $\cR$,  we will also consider $\Z$-forms of $\ffg$ and $\cU_\C(\ffg)$. 

 If $\ffg'$ is a Lie-subalgebra of $\ffg$, we denote by $\cU_\C(\ffg')$ the envelopping algebra of $\ffg'$. We regard it as a subalgebra of $\cU_\C(\ffg)$.

  Following \cite[4]{tits1987uniqueness} (see \cite[7.2]{marquis2018introduction}), we introduce a $\Z$-form of the enveloping algebra $\cU_\C(\ffg)$. Set $\fh'_\Z=\bigoplus_{i\in I} \Z\alpha_i^\vee\subset \fh$ and $\fh'=\bigoplus_{i\in I} \C \alpha_i^\vee$.

  Let $x\in \cU_\C(\ffg)$. We set $x^{(s)}=x^s/s!$ and $\tbinom{x}{s}=x(x-1)\ldots (x-s+1)/s!$. 
 
 \begin{Definition}
We define $\cU^+$, $\cU^-$ and $\cU^0$ as the $\Z$-subalgebras of $\cU_\C(\ffg)$ generated respectively by $\{e_i^{(s)}\mid i\in I,s\in \N\}$, $\{f_i^{(s)}\mid i\in I, s\in \N\}$ and $\{\tbinom{h}{s}\mid s\in \N, h\in \fh_\Z'\}$. We define $\cU$\index[notation]{U@$\cU$} as the $\Z$-subalgebra of $\cU_\C(\ffg)$ generated $\cU^+$, $\cU^-$ and $\cU^0$.  
 \end{Definition}
 
By \cite[Proposition 7.4]{marquis2018introduction}, $\cU^+, \cU^-$, $\cU^0$ and $\cU$ are $\Z$-forms of $\cU_\C(\fn^+)$, $\cU_\C(\fn^-)$, $\cU_\C(\fh')$ and $\cU_\C(\ffg)$ respectively and the product map $\cU^-\otimes \cU^0\otimes \cU^+\rightarrow \cU$ is bijective.

We graduate $\cU_\C(\ffg)$ by $Q$ as follows (see \cite[3.4]{marquis2018introduction} for the definition of a gradation). For $i\in I$,  $\deg(e_i)=\alpha_i$, $\deg(f_i)=-\alpha_i$, $\deg(h)=0$, for $h\in Y$, $\deg(xy)=\deg(x)+\deg(y)$ for all $x,y\in \cU_\C(\ffg)$ which can be written as a product of non-zero elements of $\ffg$. The algebra $\cU_\C(\ffg)$  decomposes as $\cU_\C(\ffg)=\bigoplus_{\alpha\in Q}\cU_\C(\ffg)_\alpha$ where $\cU_\C(\ffg)_\alpha$ is the set of elements of $\cU_\C(\ffg)$ whose degree is well-defined and is $\alpha$, for $\alpha\in Q$.  For $\alpha\in Q$, we set $\ffg_\alpha=\ffg_\alpha\cap \cU$\index[notation]{g@$\ffg_{\alpha,\Z}$}, $\cU_\alpha=\cU_\C(\ffg)_\alpha\cap \cU$\index[notation]{U@$\cU_{\alpha}$} and $\cU_{\alpha,\cR}=\cU_\alpha\otimes \cR$\index[notation]{U@$\cU_{\alpha,\cR}$}.

  For a ring $\cR$, we set $\cU_\cR=\cU\otimes_\Z \cR$\index[notation]{U@$\cU_{\cR}$}, $\cU^+_{\cR}=\cU^+\otimes_\Z \cR$ and $\cU^-=\cU^-\otimes_\Z \cR$.

Let $h\in \fh$ and $\alpha\in \Delta$. Let $e\in g_{\alpha,\Z}$.  We have: \begin{equation}\label{e_cmt}
he=he-eh+eh=[h,e]+eh=\alpha(h)e+eh=e(\alpha(h)+h).
\end{equation}

\begin{Lemma}\label{l_wghts_cU}
Let $h\in \fh$, $\lambda\in Q_+$ and $u\in \cU_\lambda$. Then $hu=uh+\lambda(h)u$.
\end{Lemma}

\begin{proof}
 Write $\lambda=\sum_{i\in I} n_i\alpha_i$, with $(n_i)\in \Z^I$. An element of $\cU_\lambda$ can be written as a sum of elements of the form \[u=h_1\ldots h_k x_{i_1}^{|s_1|}\ldots x_{i_\ell}^{|s_\ell|},\] where $k,\ell\in \N$, $h_1,\ldots,h_k\in \fh_\Z$, $s_1,\ldots,s_\ell\in \Z$ are such that $\sum_{i=1}^\ell  s_i\alpha_i=\lambda$ and $x_i=e_i$ if $s_i>0$ and $x_i=f_i$ is $s_i<0$. As $\fh$ is commutative, we have $hh_1\ldots h_k=h_1\ldots h_k h$ and by an induction on $|s_{i_j}|$ using \eqref{e_cmt}, we have $hx_{i_j}^{|s_{i_j}|}=x_{i_j}^{|s_{i_j}|} (s_{i_j}\alpha_{i_j}(h)+h)$, for $j\in \llbracket 1,k\rrbracket$. Lemma follows.
 \end{proof}

The algebra $\cU_\C(\ffg)$ and $\cU$ are naturally equipped with the structure of cocommutative Hopf algebras (see \cite[Definition 8.42]{marquis2018introduction} or \cite{bourbaki2006groupes2}).  By  \cite[§ 1.4 Proposition 7]{bourbaki2006groupes2}, there exists a unique comultiplication $\nabla:\cU_\cR\rightarrow \cU_\cR\otimes \cU_\cR$\index[notation]{n@$\nabla$} such that: \begin{equation}\label{e_nabla}
\nabla(x)=x\otimes 1+1\otimes x, \forall x\in \ffg_\cR.
\end{equation}
 The counit $\epsilon:\cU_\cR\rightarrow \cR$\index[notation]{e@$\epsilon$} associates to each $x\in \cU_\cR$ its \textbf{constant term} (see \cite[§ 2.1]{bourbaki2007groupes1}).  Explicitly, denote by  $\tau:\cU_\cR\rightarrow \cU_\cR$\index[notation]{t@$\tau$} the   coinverse. Then $\nabla$, $\epsilon$ and $\tau$ are the unique algebra morphisms such that for $h\in \fh$, $n\in \N$ and $i\in I$, we have: \begin{align}
\nabla (\tbinom{h}{n})&=\sum_{k,\ell\in \N\mid k+\ell=n}\tbinom{h}{k}\otimes \tbinom{h}{\ell}, &\epsilon (\tbinom{h}{n})=0,\text{for }n>0\text{ and }&\tau(\tbinom{h}{n})=\tbinom{-h}{n},\label{e_hopf_h}\\
\nabla (e_i^{(n)})&=\sum_{k,\ell\in \N\mid k+\ell=n} e_i^{(k)}\otimes e_i^{(\ell)}, &\epsilon(e_i^{(n)})=0,\text{ for }n>0\text{ and }&\tau(e_i^{(n)})=(-1)^n e_i^{(n)},\label{e_hopf_e}\\
\nabla(f_i^{(n)})&=\sum_{k,\ell\in \N\mid k+\ell=n} f_i^{(k)}\otimes f_i^{(\ell)}, &\epsilon(f_i^{(n)})=0, \text{ for }n>0\text{ and }&\tau(f_i^{(n)})=(-1)^n f_i^{(n)}.\label{e_hopf_f}
\end{align}

\subsection{The split torus scheme}
We now introduce the split torus scheme (see \cite[7.3.3]{marquis2018introduction}). We denote by $\Z-\mathrm{alg}$ the category of $\Z$-algebras and by $\mathrm{Grp}$ the category of groups.

\begin{Definition}
The \textbf{split torus scheme} associated with $\cS$ is the group functor $\fT=\fT_\cS/\Z-\alg\rightarrow \mathrm{Grp}$ defined by \[\fT(\cR)=Y\otimes \cR^\times\text{ for any ring }\cR.\]
\end{Definition}

The functor $\fT$ is isomorphic to $(\GG_m)^r$, where $r$ is the rank of $X$ over $\Z$. Alternatively, $\fT=\mathrm{Spec}(\Z[X])$. For a given ring $\cR$,  $\fT(\cR)=\Hom_{\Z-\alg}(\Z[X],\cR)\simeq \Hom_{\mathrm{Grp}}(X,\cR^\times)$, the isomorphism is given by: \[Y\otimes_\Z \cR^\times\rightarrow \Hom_{\mathrm{Grp}}(X,\cR^\times): h\otimes r\mapsto  [r^h:X\rightarrow \cR^\times: \gamma\mapsto r^h(\gamma):=r^{\gamma(h)}].\]

For $\alpha\in \Delta$ and $t\in \fT(\cR)$, we often write $\alpha(t)$  instead of $t(\alpha)$.

\subsection{The Steinberg functor}\label{ss_Tits_fnctr}

We now define the Steinberg functor $\mathfrak{St}_A$ associated with the Kac--Moody matrix $A$. The constructive Tits functor will be defined as a quotient of the amalgamated product $\mathfrak{St}_A*\fT_{\cS}$. Let $\ffg$ be the Kac--Moody algebra of $\cS$. We use the same notation as in Section~\ref{s_KM_algebras}. Our reference is \cite[7.4]{marquis2018introduction}.

If $\alpha\in \Phi$. Let $E_\alpha=\{e_\alpha,e_{-\alpha}\}$\index[notation]{e@$e_\alpha$} be the double basis of  \cite[Remark 7.6]{marquis2018introduction} ($e_\alpha$ is a particular choice of an element of $\ffg_{\alpha,\Z}$). Let $\fU_\alpha$ be the group scheme over $\Z$ isomorphic to $\GG_a$ (where $\GG_a(\cR)=(\cR,+)$ for all ring $\cR$) whose Lie algebra is the $\Z$-Lie subalgebra of $\ffg$ generated by $E_\alpha$. Let $x_\alpha:\GG_a\rightarrow \fU_\alpha$ denote the corresponding isomorphism, which we write for any ring $\cR$ as $x_\alpha:\GG_a(\cR)\overset{\sim}{\rightarrow}\fU_{\alpha}(\cR):r\mapsto \exp(r e_\alpha)$. Note that this is just a suggestive notation.

Let $\Psi\subset \Phi$ be a nilpotent set of roots.  The sum $\ffg_{\Psi}=\bigoplus_{\alpha\in \Psi}\ffg_\alpha$ is a nilpotent Lie algebra. Let $\cU_{\Psi}=\cU\cap \cU_{\C}(\ffg_{\Psi})$ (where $\cU_\C(\ffg_{\Psi}\subset \cU_\C(\ffg)$ is the envelopping algebra of $\ffg_\Psi$).
 This is a subalgebra of $\cU$. 
 
 \begin{Proposition}\label{p_Mar7.43}(see \cite[Proposition 7.43]{marquis2018introduction} or \cite[§ 9.2.2]{remy2002groupes}) Let $\{\alpha,\beta\}\in \Delta$ be a prenilpotent pair of distinct roots. Fix an arbitrary order on $]\alpha,\beta[_{\N}$. Then there exist integers $C_{i,j}^{\alpha,\beta}$ depending only on $\alpha,\beta$ and the chosen order, such that in the ring $\cU_{[\alpha,\beta]_{\N}}[\![t_1,t_2]\!]$ 
 (where $t_1,t_2$ are indeterminates), we have: \[[\exp(t_1 e_\alpha),\exp(t_2 e_\beta)]=\prod_{\gamma=p\alpha+q\beta\in ]\alpha,\beta[_{\N}} \exp(t_1^p t_2^q C_{p,q}^{\alpha,\beta} e_\gamma). \]
 
 \end{Proposition}

\begin{Definition}\label{d_St_fnctr}
The \textbf{Steinberg functor} associated to the Kac--Moody matrix $A$ is the group functor $\fS\ft_A$\index[notation]{s@$\mathfrak{St}_A$} from the category of rings to the category of groups defined as follows. If $\cR$ is a ring, then $\fS\ft_A(\cR)$ is the quotient of the free product of all the $\fU_{\gamma}(\cR)$, for $\gamma\in \Phi$, by the relations \begin{equation}
[x_{\alpha}(t_1),x_\beta(t_2)]=\prod_{\gamma=p\alpha+q\beta\in ]\alpha,\beta[_{\N}} x_\gamma(C_{p,q}^{\alpha,\beta} t_1^p t_2^q),
\end{equation}

for all prenilpotent pairs $\{\alpha,\beta\}$ and $t_1,t_2\in \cR$, where the $C_{p,q}^{\alpha,\beta}$ are as in Proposition~\ref{p_Mar7.43}. 
\end{Definition}

Let $\cR$ be a ring. For $i\in I$ and $r\in \cR^\times $, we set: \begin{equation}\label{e_tilde_si}
\tilde{s_i}(r)=x_{\alpha_i}(r)x_{-\alpha_i}(r^{-1})x_{\alpha_i}(r), \tilde{s_i}=\tilde{s_i}(1).
\end{equation}\index[notation]{s@$\tilde{s_i}$}

Let $i\in I$. In order to be coherent with the notation $\tilde{s_i}$ and $s_i^*$ below, we sometimes write $s_i$ instead of $r_i$. The group $W^v$ acts on $\fT(\cR)$ by $w.t(\lambda)=t(w^{-1}.\lambda)$, for $w\in W^v$ and $\lambda\in X$.

Denote by $\mathrm{Aut}(\ffg)$ the group of automorphisms of the Kac--Moody algebra $\ffg$.  By   \cite[Lemma 4.1]{marquis2018introduction}, if $x\in \{e_i\mid i\in I\}\cup \{f_i\mid i\in I\}$, $\mathrm{ad}(x):\ffg\rightarrow \ffg$ is \textbf{locally nilpotent}\index{locally nilpotent}, which means that for all $v\in \ffg$, there exists $n(v)\in \N$ such that $(\mathrm{ad}(x))^{n(v)}(v)=0$. This enables to define $\exp \mathrm{ad} (x):\ffg\rightarrow \ffg$ by the usual formula. Let  $i\in I$. As in  \cite{marquis2018introduction}, we set: \begin{equation}\label{e_si*}
 s_i^*=r_i^{\mathrm{ad}}=\exp(\mathrm{ad}(f_i))\exp(\ad(e_i))\exp(\ad(f_i))\in \mathrm{Aut}(\ffg),
 \end{equation}\index[notation]{s@$s_i^*$} which is well-defined  by \cite[Lemma 4.1]{marquis2018introduction}. By \cite[(7.33)]{marquis2018introduction}, we have \[s_i^*(x_\beta(r))=x_{s_i(\beta)}(\pm r),\] for $i\in I$, $\beta\in \Phi$ and $r\in \cR$, where the sign is given by $s_i^*(e_\beta)=\pm e_{r_i(\gamma)}$.

Let $w\in W^v$, $i\in I$ and $\alpha=w.\alpha_i$. Let $w=r_{i_1}\ldots r_{i_k}$ be  a reduced decomposition of $w$ (i.e $k=\ell(w)$ and $i_1,\ldots,i_k\in I$).  We set: \begin{equation}\label{e_star_w}
w^*=s_{i_1}^*\ldots s_{i_k}^*\in \mathrm{Aut}(\ffg).
\end{equation}\index[notation]{w@$w^*$}

This is well-defined, independently of the choice of a reduced decomposition, by \cite[Proposition 7.57]{marquis2018introduction}.

 \subsection{The constructive Tits functor}\label{ss_Tits_minKM}
We now introduce the constructive Tits functor, which was introduced in \cite{tits1987uniqueness}.

\begin{Definition}
The \textbf{constructive Tits functor of type $\cS$} is the group functor $\fG_{\cS}$ from the category of rings to the category of groups such that, for each ring $\cR$, the group $\fG_{\cS}(\cR)$ is the quotient of the free product $\fS \ft_A(\cR)*\fT(\cR)$ by the following relations, where $i\in I$, $r\in \cR$, $\beta\in \Phi$ and $t\in \fT(\cR)$:

\begin{enumerate}[label=\blue{(KMT4)}]
\item\label{a_KMT4} $t. x_{\alpha_i}(r). t^{-1} =x_{\alpha_i}(\alpha_i(t)r)$\axiom{KMT4@\ref{a_KMT4}}
\end{enumerate}
\vspace{-6mm}
\begin{enumerate}[label=\blue{(KMT5)}]
\item\label{a_KMT5} $\tilde{s}_i .t .\tilde{s_i}^{-1} = s_i(t)$\axiom{KMT5@\ref{a_KMT5}}
\end{enumerate}
\vspace{-6mm}
\begin{enumerate}[label=\blue{(KMT6)}]
\item\label{a_KMT6} $\tilde{s}_i(r^{-1}) = \tilde{s}_i \alpha_i^{\vee}(r)$\axiom{KMT6@\ref{a_KMT6}}
\end{enumerate}
\vspace{-6mm}
\begin{enumerate}[label=\blue{(KMT7)}]
\item\label{a_KMT7} $\tilde{s}_{i}.x_{\beta}(r).\tilde{s}_{i}^{-1} = x_{r_i.\beta}(\epsilon r)\text{, where }\epsilon\in \{-1,1\}\text{ is such that }s_i^*(e_\beta)=\epsilon e_{r_i.\beta}.$\axiom{KMT7@\ref{a_KMT7}}
\end{enumerate}

\end{Definition}

By~\ref{a_KMT4}, \ref{a_KMT5} and \ref{a_KMT7}, we have for all $t\in \fT(\cR)$ and $\alpha\in \Phi$: 
\begin{enumerate}[label=\blue{(KMT4')}]
\item\label{a_KMT4'}  $tx_\alpha(r)t^{-1}=x_\alpha(\alpha(t)r).$\axiom{KMT4'@\ref{a_KMT4'}}
\end{enumerate}

 Let $\cR$ be a ring. By \cite[Remark 7.52]{marquis2018introduction}, if $\alpha\in \Phi$, then the natural maps $\fU_\alpha(\cR)\rightarrow \fG_{\cS}(\cR)$ and $\fT(\cR)\rightarrow \fG_{\cS}(\cR)$ are injective. One then identifies these groups and their images in $\fG_\cS(\cR)$.

For $i\in I$, and $r\in \cR$, we have $(\tilde{s_i}(r))^{-1}=\tilde{s_i}(-r)$, by defintion. By \ref{a_KMT6} applied with $r=1$, we deduce: \begin{equation}\label{e_s_i_sq}
\tilde{s_i}^2=\alpha_i^\vee(-1)\in \fT(\cR).
\end{equation}

Let $w\in W^v$, $i\in I$ and $\alpha=w.\alpha_i$. Let $w=r_{i_1}\ldots r_{i_k}$ be  a reduced decomposition of $w$.  We set: \begin{equation}\label{e_tilde_w}
\tilde{w}=\tilde{s}_{i_1}\ldots \tilde{s}_{i_k}.
\end{equation}\index[notation]{w@$\tilde{w}$}

This is well-defined, independently of the choice of a reduced decomposition, by \cite[Proposition 7.57]{marquis2018introduction}.

If $\Psi\subset \Phi$, we set $\fU_{\Psi}(\cR)=\langle \fU_{\alpha}(\cR)\mid \alpha\in \Psi\rangle\subset \fG_{\cS}(\cR)$\index[notation]{u@$\fU_{\Psi}$}. This is normalized by $\fT(\cR)$ by ~\ref{a_KMT4}. We set $\fU^{\pm}=\fU_{\Phi_{\pm}}(\cR)$ and we write $\fB^\pm(\cR):=\fT(\cR)\cdot \fU^{\pm}(\cR)=\langle \fT(\cR),\fU^{\pm}(\cR)\rangle$. These are respectively the positive and negative standard \textbf{unipotent subgroups}\index{standard unipotent subgroup} and \textbf{Borel subgroups}\index{standard Borel subgroup} of $\fG_{\cS}$.  We also let $\fN(\cR)$ be the subgroup of $\fG_\cS(\cR)$ generated by $\fT(\cR)$ and the elements $\tilde{s_i}$, for $i\in I$. By \cite[8.4.1]{remy2002groupes}, if $\cR$ is a field with at least four elements, then $\fN(\cR)$ is the normalizer of $\fT(\cR)$ in $\fG(\cR)$.

 By \cite[9.5 \& 9.6]{remy2002groupes}, the following axiom is satisfied:

\begin{enumerate}
\item[(KMG5)] There exists a morphism $\mathrm{Ad}:\fG_\cS(\C)\rightarrow \mathrm{Aut}(\ffg)$ whose kernel is contained in $\fT(\C)$, such that for all $c\in \C$, $t\in \fT(\C)$ and $i\in I$, \[\mathrm{Ad}[\varphi_i(\begin{psmallmatrix}
1 & c \\ 0 & 1 \end{psmallmatrix})]=\exp \mathrm{ad} (ce_i),\quad  \mathrm{Ad}[\varphi_i(\begin{psmallmatrix}
1 & 0 \\ c & 1 \end{psmallmatrix})]= \exp \mathrm{ad} (-cf_i),\]  \[\mathrm{Ad}(\eta(t))(e_i)=t(\alpha_i)(e_i),\quad \mathrm{Ad}(\eta(t))(f_i)=-t(\alpha_i)f_i.\]
\end{enumerate}

Note that this axiom (KMG5) corresponds to the fact that $\fG$ ``integrates'' $\ffg$.

We set $\fN(\cR)=\langle \tilde{s_i}\mid i\in I,\fT(\cR)\rangle\subset \fG(\cR)$\index[notation]{n@$\fN$}. If $\cR$ is a field with at least four elements, then $\fN(\cR)$ is the normalizer of $\fT(\cR)$ in $\fG(\cR)$.

By \cite[Corollary 7.54]{marquis2018introduction}, there exists a natural transformation $\Ad:\fG_\cS\rightarrow \mathfrak{Aut}(\ffg_\cS)$ \index[notation]{a@$\Ad$} (resp. $\Ad:\fG_{\cS}\rightarrow \mathfrak{Aut}(\cU_{\cS})$)  satisfying, for any ring $\cR$: \begin{align}
\Ad_{\cR}(x_\alpha(r))&=\exp(\mathrm{ad}\  e_\alpha\otimes r)=\sum_{n\geq 0} \frac{(\mathrm{ad}\ e_\alpha)^n}{n!}\otimes r^n,\forall\alpha\in \Phi,r\in \cR\label{e_Ad_U}\\
\Ad_{\cR}(t)(u_\beta)&=t(\beta)\cdot u_\beta,\forall t\in \fT(\cR),\beta\in Q, u_\beta\in \cU_{\beta,\cR}\label{e_Ad_T}.
\end{align}

\subsection{Root datum for a split minimal Kac--Moody group}\label{ss_rt_data}
Let $\cK$ be a field and $G=\fG_\cS(\cK)$. We use gothic  letters to denote the functors $\fT$, $\fU_\alpha$, ...  and roman letters to denote their evaluation in $\cK$: $T=\fT(\cK)$, $U_\alpha=\fU_\alpha(\cK)$, ... 

The following proposition asserts in particular that $(G,(U_\alpha)_{\alpha\in \Phi},T)$ is a ``donnée radicielle jumelée entière'' for \cite[Définition 6.2.5]{remy2002groupes}.

\begin{Proposition}\label{p_DRJ}

\begin{enumerate}
\item We have $T\simeq (\cK^\times)^{\mathrm{rk}_\Z(X)}$ and  for all $\alpha\in \Phi$, $U_\alpha$ is isomorphic to $(\cK,+)$.

\item The group $G$ is generated by $T$ and the $U_\alpha$, $\alpha\in \Phi$.

\item For all $\alpha\in\Phi$, $U_\alpha$ is normalized by $T$. More precisely, $tx_\alpha(u)t^{-1}=x_\alpha(\alpha(t)u)$ for all $\alpha\in \Phi$, $t\in T$ and $u\in \cK$.

\item\label{i_Stnbrg_pres} Let $\{\alpha,\beta\}$, be a prenilpotent  pair of roots.  We fix an order on $[\alpha,\beta]$. Then the product map $\prod_{\gamma\in [\alpha,\beta] } U_\gamma \rightarrow \langle U_\gamma\mid\gamma\in [\alpha,\beta]\rangle$ is a bijection. Moreover, there exist integers $C_{p,q}$ such that for all $u,v\in \cK$, $[x_\alpha(u),x_\beta(v)]=\prod_{\gamma\in ]\alpha,\beta[, \gamma=p\alpha+q\beta} x_\gamma (C_{p,q} u^p v^q)$. These integers a priori depend on the order choosen on $[\alpha,\beta]$ but not on the field $\cK$.

\item\label{i_prenilp_pairs} Let $\{\alpha,\beta\}$ be a non-prenilpotent pair, with $\alpha+\beta\neq 0$. Then the canonical morphism $U_\alpha * U_\beta\rightarrow \langle U_\alpha, U_\beta\rangle$ (where $U_\alpha * U_\beta$ is the free product of $U_\alpha$ and $U_\beta$) is an isomorphism.

\item For all $\alpha\in \Phi$ and all $u\in U_\alpha\backslash\{1\}$, there exist $u',u''\in U_{-\alpha}$ such that $m(u):=u'uu''$ conjugates $U_\beta$ in $U_{r_\alpha(\beta)}$ for all $\beta\in \Phi$. Moreover, one can require that for all $u,v\in U_\alpha\backslash\{1\}$, $m(u)T=m(v)T$.

\item If $U^+$ (resp. $U^-$) is the subgroup of $G$ spanned by the $U_\alpha$ for $\alpha\in \Phi_+$ (resp. $\alpha\in \Phi_-$), we have $T\cdot U^+\cap U^-=\{1\}$.

\end{enumerate}
\end{Proposition}

\begin{proof}
(1) follows from the definition of $\fT$ and of  the $\fU_\alpha$. (5) is Proposition~5 of  Cours 1989-1990 of \cite{tits2013resume}. The other assertions are a combination of \cite[8.4.1]{remy2002groupes} or \cite[Proposition 1.5]{rousseau2006groupes} and of  \cite[1.6]{rousseau2006groupes}.
\end{proof}

Let $N=\fN(\cK)$.

Let $\alpha\in \Phi$ and $a\in \cK^*$. By \cite[1.5]{rousseau2006groupes} one can take: \begin{equation}\label{e_m}
m(x_\alpha(a))=x_{-\alpha}(a^{-1})x_{\alpha}(a)x_{-\alpha}(a^{-1})
\end{equation} in Proposition~\ref{p_DRJ} (6). By \eqref{e_tilde_si} and \ref{a_KMT7} and \ref{a_KMT6}, if we write $\alpha=w.\alpha$, for $w\in W^v$ and $i\in I$, we have  \begin{equation}\label{e_m_s}
m(x_{\alpha}(a))=\tilde{w}\tilde{s_i} \alpha_i^\vee(\epsilon a^{-1})\tilde{w}^{-1}\in N, 
\end{equation}

for some $\epsilon\in \{-,+\}$.

By \cite[Corollary 7.70]{marquis2018introduction} and \cite[Proposition 3.2]{kac1985defining}, we have the Bruhat decomposition of $G$ for both $\epsilon\in \{-,+\}$: \begin{equation}\label{e_Bruhat}
G=\bigsqcup_{n\in N} U^\epsilon n U^\epsilon
\end{equation}\index{d@Bruhat decomposition}

and the Birkhoff decomposition of $G$ (\cite[Corollary 7.70 and Proposition B.33]{marquis2018introduction}): \begin{equation}\label{e_Birkhoff}
G=\bigsqcup_{n\in N} U^+ nU^-.
\end{equation}\index{d@Birkhoff decomposition}

\begin{Lemma}\label{l_utu}
Let $(u_+,u_-,t),(u_+',u_-',t')\in U^+\times U^-\times T$ be such that $u_+ u_- t=u'_+u'_- t'$ or $u_+t u_-=u'_+ t' u'_-$. Then $u_+=u'_+$, $u_-=u'_-$ and $t=t'$. 
\end{Lemma}

\begin{proof}
By uniqueness in the  Birkhoff decomposition \eqref{e_Birkhoff}, we get the result in the second case. In the first case ($u_+ u_- t=u'_+u'_- t'$), \ref{a_KMT4'} provides the existence of $v_-$, $v'_-\in U^-$ such that $u_+t v_-=u'_+ t' v'_-$.  Using uniqueness in the  Birkhoff decomposition again, we deduce $t=t'$ and $u_+=u'_+$. Therefore $u_-=u'_-$. 
\end{proof}

\subsection{Untwisted affine Kac--Moody groups and loop groups}\label{subsubAffine KM groups and loop groups}

Let $\mathring{A}$ be an indecomposable Cartan matrix and $A$ be the untwisted associated affine Kac--Moody matrix (see Subsection~\ref{ss_cT_aff}). Let $\ffg$ be the Kac--Moody algebra of $A$ and $\mathring{\ffg}$ be the (finite dimensional) Kac--Moody algebra of $\mathring{\ffg}$. We saw in Subsubsection~\ref{sss_aff_la} an explicit description of $\ffg$ using the loop algebra $\mathring{\ffg} \otimes \C[t,t^{-1}]$ (with $t$ an indeterminate). This description can be ``integrated''. More precisely, if $\cS$ is a Kac--Moody datum associated with $A$, then  we can describe explicitly $\fG_\cS(\cK)$, for  $\cK$ a field, up to a central extension (if we assume  $(\alpha_i)_{i\in I}$ and $(\alpha_i^\vee)_{i\in I}$ to be free).

Let $\cS=(A,X,Y,(\alpha_i^\vee)_{i\in I},(\alpha_i^\vee)_{i\in I})$ be a Kac--Moody datum such that $(\alpha_i)$ and $(\alpha_i^\vee)$ are free and $Y/\sum_{i\in I} \Z\alpha_i^\vee$ is torsion-free (this is condition of 6.1.16 of \cite{kumar2002kac}). Let $\mathring{\fG}$ be the connected, simply-connected algebraic group such that $\mathring{G}:=\mathring{\fG}(\C)$ admits $\mathring{\ffg}$ as a Lie algebra. One defines an action of $\C^*$ on $\C[t,t^{-1}]$ by setting $a.P(t)=P(at)$ for all $a\in \C^*$ and $P(t)\in\C[t,t^{-1}]$. This action defines an action of $\C^*$ on $\mathring{\fG}(\C[t,t^{-1}])$. Let $C\subset \mathring{G}\subset \C^*\ltimes \mathring{\fG}(\C[t,t^{-1}])$ be the center of $\mathring{G}$.

Let $G=\fG_\cS(\C)$. By \cite[3.20]{rousseau2016groupes}, $G$ is the minimal Kac--Moody group of  \cite[7.4]{kumar2002kac} (the Kac-Peterson group). By  \cite[Corollary 13.2.9]{kumar2002kac}, there exists an isomorphism $\psi:G/\mathcal{Z}(G)\overset{\sim}{\rightarrow} (\C^*\ltimes \mathring{G}(\C[t,t^{-1}])/C$, where $\mathcal{Z}(G)$ is the center of $G$. Moreover, by  \cite[Theorem~13.2.8 c)]{kumar2002kac},  if the $\fU_\alpha$, $\alpha\in \Phi$ and $\mathring{\fU}_{\mathring{\alpha}}$, $\mathring{\alpha}\in \mathring{\Phi}$ (with the same notation as in Subsection~\ref{ss_cT_aff}) are the  root subgroups of $\widetilde{\fG}_\cS$ and $\mathring{\fG}$, we have  $\fU_{\mathring{\alpha}+k\delta}(\C)=\{1\}\ltimes \mathring{\fU}_{\mathring{\alpha}}(t^k\C)\subset \C^* \ltimes  \mathring{\fG}(\C[t,t^{-1}])$. 
Let $\mathring{\fT}:=\mathrm{Hom}_{\Z-\mathrm{alg}}(\Z[\mathring{X}],\GG_m)$ be the standard maximal torus of $\mathring{\fG}$ and $T=\C^*\ltimes \mathring{T}\subset  \C^* \ltimes  \mathring{\fG}(\C[t,t^{-1}])$ (here the semi-direct product is the direct product because $\C^*$ acts trivially on $\mathring{G}$). Then by  \cite[Theorem 13.2.8 b)]{kumar2002kac}, $T/C$ is the image by $\psi$ of the standard maximal torus of $G$. The adjoint action of $\C^* \ltimes  \mathring{\fG}(\C[t,t^{-1}])$ on $\mathring{\ffg}\otimes \C[t,t^{-1}]\oplus \C c\oplus \C d$ can also be described explicitly, see \cite[page 493]{kumar2002kac}. We also refer to \cite[7.6]{marquis2018introduction} for the description of $\fG_{\cS}(\cK)$, for an arbitrary field $\cK$.

\section{The affine group scheme $\fU^{ma+}$}\label{s_gp_Upma}

We now introduce the ``affine group scheme`` $\fU^{ma+}$. It is a ``completion'' of the functor $\fU^+$. If $\cR$ is a ring, then $\fU^{ma+}(\cR)$ is  defined as a subgroup of a completion of the Hopf enveloping  algebra $\cU_\cR$. 

We start by defining the completion of a graded vector space (Subsection~\ref{ss_gded_v_s}) and apply this to define the completed envelopping algebra $\widehat{\cU}^+_\cR$. We  define the 
twisted exponentials (Subsection~\ref{ss_twstd_exp}), which enable to give a normal form to the elements of $\fU^{ma+}(\cR)$. We study decompositions of $\fU^{ma+}$ (see Subsections~\ref{ss_dec_fU_ma+} and \ref{ss_Lv_dec}). We then describe explicitly $\fU^{ma+}$ in the case of affine $\mathrm{SL}_2$, using Garland-Mitzman polynomials (see Subsection~\ref{ss_Gr_Mitz_pl}).

\subsection{Graded vector spaces}\label{ss_gded_v_s}

Let $V=\bigoplus_{n\in \Z} V_n$ be a $\Z$-graded free module over $\cR$. Let $\widehat{V}=\widehat{V}^p:=\bigoplus_{n\in \Z_{<0}} V_n \oplus \prod_{n\in \N} V_n$ be its \textbf{positive completion with respect to the graduation}.  For $n\in \Z$, set $V_{\geq n}=\bigoplus_{m\in \Z_{\geq n}} V_m$ and $\widehat{V}_{\geq n}=\prod_{m\in \Z_{\geq n}} V_m$. We equip $\widehat{V}$ with the metric $d$ defined by \[d(x,y)=\exp(-\max\{n\in \Z\mid x-y\in \widehat{V}_{\geq n}\}).\] Then by \cite[Exercise 8.7]{marquis2018introduction}: \begin{enumerate}
\item $d$ is a translation-invariant ultrametric on $\widehat{V}$,

\item  $\widehat{V}$ is the completion of $V$ with respect to $d$,

\item  a basis of neighborhood of $0$ is given by $\{\widehat{V}_{\geq n}\mid n\in \N\}$

\item every bounded subset of $\widehat{V}$ is contained in some $\widehat{V}_{\geq n}$, for some $n\in \Z$. 
\end{enumerate} 

Let $\cA=\bigoplus_{n\in \Z} \cA_n$ be a $\Z$-graded free module over $\cR$. We say that $\cA$ is a \textbf{graded (associative) $\cR$ algebra}, if $\cA$ is equipped with the structure of an $\cR$-algebra, compatible with the grading, i.e  $\cA_k\cdot \cA_n\subset \cA_{k+n}$, for $k,n\in \Z$. Let $\widehat{\cA}^p:=\bigoplus_{n\in -\N^*}\cA_n\oplus \prod_{n\in \N} \cA_n$.  For $(a_k), (b_k)\in \prod_{k\in \Z}\cA_k$ such that $a_k=b_k=0$ for $k\ll 0$, , we define \[(\sum_{k\in \N} a_k)( \sum_{\ell\in \N} b_\ell)=\sum_{k,\ell\in \N}a_k b_\ell\in \widehat{\cA}^p.\] This is well-defined and this equips $\widehat{\cA}^p$ with the structure of an associative $\cR$-algebra, whose multiplication is continuous. 

Let $M=\bigoplus_{n\in \Z}M_n$ be a $\Z$-graded free module over $\cR$. We assume that $M$ is equipped with an action $\pi$ of $\cA$ on $M$, compatible with the grading, i.e $\cA_k. M_n\subset M_{k+n}$, for $k,n\in \Z$. Let $\widehat{M}^p$ be its positive completition (as a vector space) with respect to the graduation. Then $\pi$ naturally extends to an action $\hat{\pi}$ of $\widehat{\cA}^p$ on $\widehat{M}^p$ defined by:  \begin{equation}\label{e_cplted_act}
\hat{\pi}(\sum_{k\in \Z} a_k )(\sum_{n\in \Z} m_n)=\sum_{k,n\in \Z} \pi(a_k)(m_n),
\end{equation} for $(a_k)\in \prod_{k\in \Z} \cA_k$ and $(m_k)\in \bigoplus_{k\in \Z} M_k$ such that $a_k=0$ and $m_k=0$ for $k\ll 0$.

\subsection{The completed enveloping algebra $\widehat{\cU}_\cR$}

Recall that if $\alpha=\sum_{i\in I} n_i\alpha_i\in Q$, we set $\htt(\alpha)=\sum_{i\in I} n_i$. 

Let $\cR$ be a ring. We set  $\widehat{\cU}^+_\cR=\prod_{\alpha\in Q_+} \cU_{\alpha,\cR}$\index[notation]{U@$\widehat{\cU}_\cR$}. This is the completion of $\cU^+_\cR$ with respect to the graduation $Q_+$. We write   $\widehat{\cU}^+=\widehat{\cU}^+_\Z$\index[notation]{U@$\widehat{\cU}$}.  We apply the results of Subsection~\ref{ss_gded_v_s}, by setting $\cU_n=\bigoplus_{\alpha\in Q_+\mid \htt(\alpha)=n} \cU_\alpha$ (and similarly for $\cU_\cR$). 

If $(u_\alpha)\in \prod_{\alpha\in Q_+} \cU_{\alpha,\cR}$, we write $\sum_{\alpha\in Q_+} u_\alpha$ the corresponding element of $\widehat{\cU}^+_{\cR}$. A sequence $(\sum_{\alpha\in Q_+} u_\alpha^{(n)})_{n\in \N}$, with $\deg(u_\alpha^{(n)})=\alpha$ for all $\alpha\in Q_+$ and $n\in \N$,  converges in $\widehat{\cU}_\cR$ if and only if for every $\alpha\in \Delta_+$, the sequence $(u_\alpha^{(n)})_{n\in \N}$ is stationary.

Let $(E,\leq)$ be a totally ordered set. Let $(u^{(e)})\in  (\widehat{\cU}_{\cR}^+)^E$. For $e\in E$, write  $u=\sum_{\alpha\in Q_+} u_\alpha^{(e)}$, with $u_\alpha^{(e)}\in \cU_{\alpha,\cR}$, for $\alpha\in Q_+$. We assume that for every $\alpha\in Q_+$, $\{e\in E\mid u_\alpha^{(e)}\neq 0\}$ is finite. Then we set $\prod_{e\in E} u^{(e)}=\sum_{\alpha\in Q_+} u_\alpha$,  where \[u_\alpha=\sum_{\stackrel{(\beta_1,\ldots,\beta_k)\in Q_+^{(\N)},}{\beta_1+\ldots+\beta_k=\alpha}}\sum_{\stackrel{(e_1,\ldots,e_k)\in E,}{e_1<\ldots <e_k}} u_{\beta_1}^{(e_1)}\ldots u_{\beta_k}^{(e_k)}\in  \cU_{\alpha,\cR},\]  for $\alpha\in \Delta_+.$ This is well-defined since in the sum defining $u_\alpha$, only finitely many non-zero terms appear.

 The counit $\epsilon$ and comultiplication $\nabla$ both extend by continuity to $\widehat{\cU}^+_\cR$. Note however that the extension of the comultiplication $\nabla$ is a map $\widehat{\cU}^+_\cR\rightarrow \widehat{\cU}^+_\cR\hat{\otimes} \widehat{\cU}^+_\cR$. Where $\widehat{\cU}^+_\cR\hat{\otimes} \widehat{\cU}^+_\cR=\prod_{m\in \N}\sum_{\alpha,\beta\in Q_+\mid \htt(\alpha+\beta)=m} \cU_{\alpha,\cR}\otimes \cU_{\beta,\cR}$ expresses the fact that the $\cR$-bilinear property of the tensor product $\otimes$ is replaced by the stronger property \[(\sum_{\alpha\in Q_+}a_\alpha u_\alpha)\hat{\otimes}( \sum_{\beta\in Q_+} b_{\beta} v_\beta)=\sum_{\alpha,\beta\in Q_+} a_\alpha b_\beta u_\alpha\hat{\otimes}v_\beta,\] for $(u_\alpha),(v_\beta)\in \cR^{Q_+}$ and $(u_\alpha),(v_\alpha)\in \prod_{\alpha\in Q_+} \cU_{\alpha,\cR}$ (see \cite[Remark A.18]{marquis2018introduction} for more details).  For simplicity, we write $\otimes$ instead of $\hat{\otimes}$. 

We extend the adjoint action of $\cU^+_\cR$ on itself via \eqref{e_cplted_act}. Explicitly, we have $\widehat{\Ad}:\widehat{\cU}_\cR^{+}\rightarrow \End (\widehat{\cU}_\cR^+)$ defined by: \begin{equation}\label{e_ad_hat}
\widehat{\Ad}(\sum_{\lambda\in Q_+} u_\lambda)(\sum_{\mu\in Q_+} v_\mu)=\sum_{\lambda,\mu\in Q_+} \Ad(u_\lambda)(v_\mu),
\end{equation}
for $(u_\lambda),(v_\lambda)\in \prod_{\lambda\in Q_+} \cU_{\lambda,\cR}.$ This representation restricts to \eqref{e_Ad_U} on $\fU^+(\cR)$.

\subsection{The group scheme  $\fU^{ma+}$ and twisted exponentials}\label{ss_twstd_exp}  We now introduce the group-scheme $\fU^{ma+}$. We then define the ``twisted exponentials``, which enable to describe the elements of $\fU^{ma+}$.

\begin{Definition/Proposition}\label{dp_grp_lk}
We call an element $x$ of $\widehat{\cU^+_\cR}\setminus\{0\}$ \textbf{group-like}\index{group-like} if $\nabla(x)=x\otimes x$. In this case, $x$ is invertible, $x^{-1}=\tau(x)$ and $\epsilon(x)=1$. Moreover, the set of group-like elements of $\widehat{\cU^+_\cR}$ forms a subgroup of $(\widehat{\cU^+_\cR})^{\times}$.  
\end{Definition/Proposition}

\begin{proof}
Let $x\in \widehat{\cU^+_\cR}$ be such that $\nabla(x)=x\otimes x$. Then by definition of a Hopf algebra (see the two diagrams preceding \cite[Definition A.12]{marquis2018introduction}), we have $(\epsilon\otimes \Id)\circ \nabla(x)=1\otimes x=(\epsilon\otimes \Id) (x\otimes x)$ and hence $\epsilon(x)=1$. We also have $\tau(x)x=\epsilon(x)=1$ and $x\tau(x)=\epsilon(x)=1$, which proves that $x$ is invertible. Moreover, as $\nabla$ is an algebra morphism, we have $\nabla(x^{-1})=\nabla(x)^{-1}=x^{-1}\otimes x^{-1}$. 
\end{proof}

\begin{Definition}
For $\cR$ a ring and $A$ a Kac--Moody matrix, we define $\fU^{ma+}(\cR)=\fU^{ma+}_A(\cR)$ as the set of group-like elements of $\widehat{\cU^+}$ having $1$ as a constant term (i.e if $u=\sum_{\alpha\in Q_+} u_\alpha$, with $(u_\alpha)\in \cU_{\alpha,\cR}$ for all $\alpha\in Q_+$, then $u_0=1$).
\end{Definition}

In \cite[8.5.1]{marquis2018introduction}, $\fU^{ma+}$ is defined as the set of algebra morphisms of some Hopf algebra. By \cite[Theorem 8.51 (3)]{marquis2018introduction}, it is equivalent to the definition we gave.

The \textbf{standard filtration on $\cU_{\C}(\ffg)$} is defined as follows. For $m\in \N$, $\cU_\C(\ffg)_m$ is the image of $\bigoplus_{s=0}^m \ffg^{\otimes s}$ in $\cU_\C(\ffg)$. The filtration  of an element $x\in \cU_\C(\ffg)$  is the minimum $m$ such that $x\in \cU_{\C}(\ffg)_m$.

\begin{Definition}
Let $\alpha\in \Delta\cup \{0\}$ and $x\in \ffg_{\alpha,\Z}$. An \textbf{exponential sequence} for $x$ in $\cU$ or simply an exponential sequence for $x$ (resp. an exponential sequence for $x$  in $\cU_{\Q}$) is a sequence $(x^{[n]})_{n\in \N}$ of elements of $\cU$ (resp. $\cU_{\Q}$) satisfying the following three conditions:

\begin{enumerate}[label=\blue{(ES1)}]
\item\label{a_es1} $x^{[0]}=1$, $x^{[1]}=x$ and $x^{[n]}\in \cU_{n\alpha}$, for $n\in \Z_{\geq 1}$.\axiom{es1@\ref{a_es1}}
\end{enumerate}

\begin{enumerate}[label=\blue{(ES2)}]
\item\label{a_es2} $x^{[n]}-x^{(n)}$ has filtration  strictly less that $n$ in $\cU_{\C}(\ffg)$, for all $n\in \Z_{\geq 1}$.\axiom{es2@\ref{a_es2}}
\end{enumerate}

\begin{enumerate}[label=\blue{(ES3)}]
\item\label{a_es3} $\nabla(x^{[n]})=\sum_{p,q\in \N\mid p+q=n} x^{[p]}\otimes x^{[q]}$ and $\epsilon(x^{[n]})=0$, for all $n\in \Z_{\geq 1}$. \axiom{es3@\ref{a_es3}}
\end{enumerate}

 For $r\in \cR$ and $(x^{[n]})$ an exponential sequence in $\cU$ for $x$, we then set: \[[\exp](rx):=\sum_{n\in \N} x^{[n]}\otimes r^n\in \widehat{\cU}_\cR^+.\]\index[notation]{e@$[\exp]$} This is the \textbf{twisted exponential} of $rx$ associated with the sequence $(x^{[n]})_{n\in \N}$.

\end{Definition}

By definition of exponential sequences, we have the following lemma.

\begin{Lemma}\label{l_twst_exp}
Let $\alpha\in \Delta\cup \{0\}$ and  $x\in \ffg_{\alpha,\Z}$. Let $r\in \cR$. Write $[\exp](rx)=\sum_{\beta\in Q_+} u_\beta\in \widehat{\cU}_{\cR}^+$, with $u_\beta\in \cU_{\beta,\cR}$, for $\beta\in Q_+$. Then we have $u_0=1$ and $u_\alpha=x$.
\end{Lemma}

Let $\alpha\in \Delta$.   Let $x\in \ffg_{\alpha,\Z}$.  By \cite[Proposition 2.7]{rousseau2016groupes} or \cite[Proposition 8.50]{marquis2018introduction}, there always exists an exponential sequence for $x$. Note that it is not unique in general. However, if $\alpha\in \Phi_+$, the unique exponential sequence for $x$ is $(x^{(n)})_{n\in \N}=(\frac{1}{n!}x^n)$ by \cite[2.9 2)]{rousseau2016groupes}.

For $n\in \N$, set $\fU_n^{ma+}(\cR)=\fU^{ma+}(\cR)\cap (1+\widehat{\cU}^+_{\geq n})$, where $\widehat{\cU}^+_{\geq n}=\prod_{\beta\in \Delta_+\mid \htt(\beta)\geq n}\cU_{\beta}$. Then by \cite[Lemma 8.67 and Definition 8.119]{marquis2018introduction}  the sequence $(\fU_{n}^{ma+}(\cR))$ is a separated conjugation-invariant filtration of $\fU^{ma+}(\cR)$ in the sense of \cite[Exercise 8.5]{marquis2018introduction}. It equips $\fU^{ma+}(\cR)$ with the structure of a topological group, where a basis of neighborhood of the identity is given by the $\fU_n(\cR)$, $n\in \N$.

We fix for every $\alpha\in \Delta_+$ a $\Z$-basis $\cB_\alpha$ of $\ffg_{\alpha,\Z}$\index[notation]{b@$\cB_{\alpha}$}.  Set $\cB=\bigcup_{\alpha\in \Delta_+}\cB_\alpha$\index[notation]{B@$\cB$}. We fix an order $\leq$ on each  $\cB_\alpha$ and on $\Delta_+$.  
Let $\alpha\in \Delta_+$. We define: \[X_\alpha=X_\alpha^{(\cB_\alpha,\leq)}=X_\alpha^{\cB_\alpha}:\ffg_{\alpha,\Z}\otimes \cR\rightarrow \fU^{ma+}(\cR)\]\index[notation]{X@$X_\alpha$} by \[X_\alpha(\sum_{x\in \cB_\alpha} \lambda_x.x)=\prod_{x\in \cB_\alpha} [\exp]\lambda_x.x, \forall (\lambda_x)\in \cR^{\cB_\alpha}.\]   When $\alpha\in \Phi_+$, we have $\ffg_{\alpha,\Z}=\Z e_\alpha$, where $e_\alpha$ is defined in \cite[Remark 7.6]{marquis2018introduction}. We set $x_\alpha(r)=[\exp](r e_\alpha)$, for $r\in \cR$. We have: \begin{equation}\label{e_rt_sbgrp}
X_\alpha(\ffg_{\alpha,\Z}\otimes_\Z \cR)=x_\alpha(\cR):=\fU_{\alpha}(\cR).
\end{equation}

Note that by \ref{a_es3}, if $\alpha\in \Delta_+$, $x\in \ffg_{\alpha,\Z}$ and $r\in \cR$, then  $[\exp](rx)$ is group-like and by \ref{a_es1}, its constant term is $1$. Thus $[\exp](rx)\in \fU^{ma+}(\cR)$. Conversely, by \cite[Theorem 8.51]{marquis2018introduction}, every $g\in \fU^{ma+}(\cR)$ can be written in a unique way as a   product \begin{equation}\label{e_normal_form_Upma}
g=\prod_{\alpha\in \Delta_+}X_\alpha(c_\alpha),
\end{equation}  where $c_\alpha\in \ffg_{\alpha,\Z}\otimes \cR$, for $\alpha\in \Delta_+$, where the product is taken in the given order  on $\Delta_+$.

Let  $\Psi\subset \Delta_+$ be a closed subset. We set: \[\fU_{\Psi}^{ma}(\cR)=\prod_{\alpha\in \Psi} X_\alpha(\ffg_{\alpha,\Z}\otimes \cR)\subset \fU^{ma+}(\cR).\]\index[notation]{u@$\fU_{\Psi}^{ma}(\cR)$}

This is a subgroup of $\fU^{ma+}$, which does not depend on the chosen order on $\Delta_+$ (for the product). This is not the definition given in \cite{rousseau2016groupes} or \cite[page 210]{marquis2018introduction}, but it is equivalent by \cite[Theorem 8.51]{marquis2018introduction}.

\subsection{Decompositions of $\fU^{ma+}$}\label{ss_dec_fU_ma+}

We equip $\fU^{ma+}(\cR)$ with the topology induced by the topology of $\widehat{\cU}^+_\cR$. This is the same as the topology used by Marquis in \cite{marquis2018introduction}, by \cite[page 234]{marquis2018introduction} and it equips $\fU^{ma+}(\cR)$ with the structure of a topological group.

\begin{Proposition}\label{p_Mar8.58}(\cite[8.58]{marquis2018introduction}) Let $\cR$ be a ring and let $\Psi''\subset \Psi'\subset \Psi\subset \Delta_+$ be closed sets of roots. Then we have: \begin{enumerate}
\item $\fU_{\Psi}^{ma}(\cR)$ is a closed subgroup  of $\fU^{ma+}_A(\cR)
$ coinciding with $\fU^{ma+}_A(\cR)\cap \widehat{\cU}_{\cR}(\Psi)\subset \widehat{\cU}^+_\cR$.  

\item If $\Psi'$ is coclosed in $\Psi$, then we have a unique decomposition  $\fU_{\Psi}^{ma}(\cR)=\fU_{\Psi'}^{ma}(\cR)\cdot \fU_{\Psi\setminus \Psi'}^{ma}(\cR)$.

\item If $\Psi'$ is an ideal in $\Psi$, then $\fU^{ma}_{\Psi'}(\cR)$ is a normal subgroup of $\fU^{ma}_\Psi(\cR)$.

\item If $\Psi'$ is a coclosed ideal in $\Psi$, then $\fU^{ma}_{\Psi}(\cR)=\fU^{ma}_{\Psi'}(\cR)\rtimes \fU_{\Psi\setminus \Psi'}^{ma}(\cR)$.

\item If $\alpha+\beta\in \Psi''$ whenever $\alpha\in \Psi,\beta\in \Psi'$ and $\alpha+\beta\in \Delta$, then $[\fU^{ma}_{\Psi}(\cR),\fU^{ma}_{\Psi'}(\cR)]\subset \fU^{ma}_{\Psi''}(\cR)$. 

\item If $\alpha+\beta\in\Psi'$ whenever $\alpha,\beta,\alpha+\beta\in \Psi$, then $\fU^{ma}_{\Psi}(\cR)/\fU^{ma}_{\Psi'}(\cR)$ is isomorphic to the additive group of the $\cR$-module $\prod_{\alpha\in \Psi\setminus \Psi'} \ffg_{\alpha,\cR}$. 
\end{enumerate}

\end{Proposition}

Let $i\in I$. Consider $\Psi=\Delta_+$ and $\Psi'=\Delta_+\setminus\{\alpha_i\}$. Then by Proposition~\ref{p_Mar8.58} (4) and \eqref{e_rt_sbgrp}, we have \begin{equation}\label{e_smd_pdct_U}
\fU^{ma+}(\cR)=U^{ma}_{\Delta_+\setminus\{\alpha_i\}}(\cR)\rtimes U_{\alpha_i}(\cR). 
\end{equation}

\begin{Lemma}\label{l_inv_cclsd}
Let $w\in W^v$. Then $\Delta_+\setminus \Inv(w)$ is strongly closed and coclosed in $\Delta_+$. 
\end{Lemma}

\begin{proof}
Let $\alpha,\beta\in \Delta_+\setminus \Inv(w)$ (resp. $\alpha,\beta\in \Inv(w)$) and $p,q\in \Z_{\geq 1}$ be such that $p\alpha+q\beta\in \Delta_+$. Then $w.(p\alpha+q\beta)=pw.\alpha+qw.\beta\in \Delta_+$ (resp. $w.(p\alpha+q\beta)\in \Delta_-$) and thus $p\alpha+q\beta\in  \Delta_+\setminus \{\Inv(w)\}$ (resp. $p\alpha+q\beta\in \Inv(w)$). 
\end{proof}

\begin{Remark}
Let $w\in W^v$. In general, $\Delta_+\setminus \Inv(w)$ is not a strong ideal nor an ideal of $\Delta_+$. Indeed, choose $i,j\in I$ such that $i\neq j$ and set $s=r_i$, $t=r_j$. Assume that $\langle s,t\rangle$ is infinite. Set $w=st\in W^v$. Then $\Delta_+\setminus \Inv(w^{-1})=\Delta_+\cap w.\Delta_+=\Delta_+\setminus\{\alpha_s,t.\alpha_s\}$, with $\alpha_s=\alpha_i$ and $\alpha_t=\alpha_j$. Let $\alpha=\alpha_s\in \Inv(w^{-1})$ and $\beta=\alpha_t\in \Delta_+\cap w.\Delta_+$. Let $p=1$ and $q=-\alpha_t(\alpha_s^\vee)$. Then $p\alpha+q\beta=t.\alpha_s\in \Inv(w^{-1})$.  If moreover  we assume $\alpha_t(\alpha_s^\vee)=-1$, (which implies $\alpha_s(\alpha_t^\vee)\leq -4$, in order to have $\langle s,t\rangle$ infinite), then $\Delta_+\setminus \Inv(w^{-1})$ is not an ideal of $\Delta_+$. 

\end{Remark}

\begin{Corollary}\label{c_dec_inv}
Let $\cR$ be a ring. Let $w\in W^v$. Then we have a unique decomposition $\fU^{ma+}(\cR)=\fU^{ma}_{\Delta_+\setminus \Inv(w)}(\cR)\cdot \fU_{\Inv(w)}^{ma}(\cR)$. Moreover, fix an order on $\Inv(w)$. Then: \begin{equation}\label{e_U_Inv_w}
\fU_{\Inv(w)}(\cR):=\fU^{ma}_{\Inv(w)}(\cR)=\prod_{\beta\in \Inv(w)} \fU_{\beta}(\cR)\subset \fG(\cR).
\end{equation} More precisely, write $\Inv(w)=\{\beta_1,\ldots,\beta_{\ell(w)}\}$ where $(\beta_j)$ increasing for the chosen order on $\Inv(w)$.  Then every element $u$ of $\fU_{\Inv(w)}^{ma}(\cR)$ has a unique expression $u=u_1\ldots u_{\ell(w)}$, 
with $u_j\in \fU_{\beta_j}(\cR)$, for all $j\in  \llbracket 1,\ell(w)\rrbracket$.
\end{Corollary}

\begin{proof}
The first statement follows from Proposition~\ref{p_Mar8.58} and Lemma~\ref{l_inv_cclsd}. By definition, $\Inv(w)\subset \Delta_+$ and $w.\Inv(w)\subset \Delta_-$ and thus $\Inv(w)$ is prenilpotent. By Lemma~\ref{l_inv_cclsd}, $\Inv(w)$ is nilpotent. Thus the statement after ``moreover'' is a particular case of \cite[Exercise 7.62]{marquis2018introduction}. 
\end{proof}

\subsection{On the  Levi decomposition of $U^+$}\label{ss_Lv_dec}

In this subsection, we obtain a kind of Levi decomposition of $\fU^+(\cR)$. This will be useful in Section~\ref{s_chimney}, to describe the fixator of a chimney in the masure.

Let $F^v$ be a vectorial subface of $\overline{C^v_f}$ and let $\cR$ be a ring. We consider  $U^+=\fU^+(\cR)$. 
 
We set $\Delta^u_+(F^v)=\{\alpha\in \Delta_+\mid \alpha(F^v)>0\}$\index[notation]{d@$\Delta^u_+(F^v), \Delta^m_+(F^v)$} and $\Delta^m_+(F^v)=\{\alpha\in \Delta_+\mid \alpha(F^v)=\{0\}\}=\Delta_+\setminus \Delta_+^u(F^v)$.

For $\Psi$ a closed subset of $\Delta_+$, we set $U^+_{\infty}(\Psi)=U^+\cap\fU^{ma}_{\Delta^m_+(F^v)}(\cR)$\index[notation]{u@$U^+_\infty(\Psi)$}. Note that $\Delta^u_+(F^v)$ is strongly closed, strongly  co-closed and it is a strong ideal in $\Delta_+$. Therefore by Proposition~\ref{p_Mar8.58}, we have: \begin{equation}\label{e_smd_pdct_Del_u}
\fU^{ma+}(\cR)=\fU^{ma}_{\Delta^u_+(F^v)}(\cR)\rtimes \fU^{ma}_{\Delta^m_+(F^v)}(\cR).
\end{equation}
Set $U^+(\Phi^m_+(F^v))=\langle U_\alpha\mid \alpha\in \Phi^m_+(F^v)\rangle\subset U^+$. Denote by $U^+(F^v)$  the smallest normal subgroup of $U^+$ containing $\langle U_\alpha\mid \alpha\in \Phi^u_+(F^v)\rangle$, where $\Phi^m_+(F^v)=\{\alpha\in \Phi_+\mid \alpha(F^v)=\{0\}\}$ and  $\Phi^u_+(F^v)=\{\alpha\in \Phi_+\mid \alpha(F^v)=\R_{>0}\}$\index[notation]{p@$\Phi^m_+(F^v)$, $ \Phi^u_+(F^v)$}. 
 
\begin{Proposition}\label{p_Levi_dec_U}
\begin{enumerate}
\item We have $U^+=U^+(\Phi^m_+(F^v))\ltimes U^+(F^v)$. 
\item We have $U^+(\Phi^m_+(F^v))=U^+_{\infty}(\Delta^m_+(F^v))$ and  $U^+(F^v)=U^+_{\infty}(\Delta^u_+(F^v))$. 
\item Let $\pr_u:\fU^{ma+}(\cR)\twoheadrightarrow \fU^{ma}_{\Delta^u(F^v)}(\cR)$ and $\pr_m:\fU^{ma+}(\cR)\twoheadrightarrow \fU^{ma}_{\Delta^m_+(F^v)}(\cR)$ be the canonical projections from the decomposition~\eqref{e_smd_pdct_Del_u}. Then $\pr_u$ and $\pr_m$ stabilize $U^+$. 
\end{enumerate}
\end{Proposition}
\begin{proof}
(1) For $u\in U^+$, we denote by $\ell(u)$ the minimal $k\in \N$ such that there exist $\beta_1,\ldots,\beta_k\in \Phi_+$ such that $u\in U_{\beta_1}\cdot\ldots \cdot U_{\beta_k}$. We prove the result by induction on $\ell(u)$. Let $k\in \N$ and assume that for every element of $U^+$ such that $\ell(u)\leq k$, we can write $u=v_1v_2$, with $v_1\in U^+(\Phi^m_+(F^v))$ and $v_2\in U^+(F^v)$. Let $u\in U^+$ be  such that $\ell(u)\leq k+1$. Write $u=x_{\beta_1}(a_1)\ldots x_{\beta_{k+1}}(a_{k+1})$, with $\beta_1,\ldots,\beta_{k+1}\in \Phi_+$ and $a_1,\ldots,a_{k+1}\in \cR$. By assumption, we can write $x_{\beta_1}(a_1)\ldots x_{\beta_k}(a_k)=v_1v_2$, with  $v_1\in U^+(\Phi^m_+(F^v))$ and $v_2\in U^+(F^v)$. If $\beta_{k+1}\in \Phi^u(F^v)$, then $v_2x_{\beta_{k+1}}(a_{k+1})\in U^+(F^v)$. Otherwise, $\beta_{k+1}\in \Phi^m_+(F^v)$ and $u=\underbrace{v_1 x_{\beta_{k+1}}(a_{k+1})}_{\in U^+(\Phi^m_+(F^v))}\underbrace{x_{\beta_{k+1}}(-a_{k+1})v_2 x_{\beta_{k+1}}(a_{k+1})}_{\in U^+(F^v)}$. By induction we deduce that $U^+=U^+(\Phi^m_+(F^v))\cdot  U^+(F^v)$. We have  $U^+(\Phi^m_+(F^v))\cap  U^+(F^v)\subset \fU^{ma}_{\Delta^m_+(F^v)}(\cR)\cap  \fU^{ma}_{\Delta_+^u(F^v)}(\cR)=\{1\}$, by \eqref{e_smd_pdct_Del_u}, which proves (1).

(2) For $\alpha\in \Phi^m_+(F^v)$, we have $U_\alpha\subset U^+_{\infty}(\Delta^m_+(F^v))$ and thus:  \[U^+(\Phi^m_+(F^v))\subset U^+_{\infty}(\Delta^m_+(F^v)).\]

By \eqref{e_smd_pdct_Del_u}, $U^+_{\infty}(\Delta^u_+(F^v))$ is normal in $U^+$. Moreover, $U_\alpha\subset U^+_{\infty}(\Delta^u_+(F^v))$, for all $\alpha\in \Phi^u_+(F^v)$. As $U^+(F^v)$ is the smallest normal subgroup of $U^+$ containing the $\langle U_\alpha,\alpha\in \Phi^u_+(F^v)\rangle$, we deduce that: \[U^+(F^v)\subset U^+_{\infty}(\Delta^u_+(F^v)).\]

 Let $u\in U^+_{\infty}(\Delta^u_+(F^v))$. By (1), we can write $u=u_1 u_2$, with $u_1\in U^+(\Phi^m_+(F^v))$ and $u_2\in U^+(F^v)$. Then $u_2=u_1^{-1}u\in  U^+_{\infty}(\Delta^u_+(F^v))$. Therefore $u_2\in \fU^{ma}_{\Delta^u_+(F^v)}(\cR)\cap \fU^{ma}_{\Delta^m_+(F^v)}(\cR)$ which implies that $u_2=1$, by \eqref{e_smd_pdct_Del_u}. Therefore $u\in U^+(\Phi^m_+(F^v))$ and $U^+(\Phi^m_+(F^v))=U^+_{\infty}(\Delta^m_+(F^v))$. With the same reasoning,  we have $U^+(F^v)=U^+_{\infty}(\Delta^u_+(F^v))$. Using (1), we deduce (2). Then (3) follows from (1) and (2). 
 
\end{proof}

\subsection{Description of $\fU^{ma+}(\cR)$ and of the twisted exponentials in the untwisted  affine case}\label{ss_Gr_Mitz_pl}

\subsubsection{Computation of twisted exponentials in the case of affine $\mathrm{SL}_2$}\label{sss_twisted}

We now  describe a choice of  $\bx=(x_s)$ in the $\tilde{A}$ case, following \cite{marquis2018introduction}. For simplicity, we only work in the case of $A^{(1)}_1$ (i.e $A=\begin{psmallmatrix}
2 & -2 \\ -2 & 2
\end{psmallmatrix}$), although the case of $A^{(1)}_{\ell}$ is similar, see \cite[Example 8.60]{marquis2018introduction}. As $\fU^{ma+}$ depends only on $A$ (and not on $\cS$), we can choose $\cS$ to be $\mathrm{coad}(\cD^A_{\min})$ (see \cite[Example 7.12 and 7.14]{marquis2018introduction}) for definitions. In this case, $\ffg_{\cS}=\mathfrak{sl}_2(\C[t,t^{-1}])$. We have $\fh_\Z\simeq \Z$. Denote by $h$ the element of $\fh_\Z$ corresponding to $1$. We have $\Delta^{im}_+=\Z_{\geq 1} \delta$. For $s\in \Z_{\geq 1}$, the space $\cB_{s\delta}$ admits $h_s=h\otimes t^{s}$ as a basis. We can then choose, for $n\in \N$ and $s\in \Z_{\geq 1}$, the exponential sequence  \[(h_s)^{[n]}:= \Lambda_n(h_s,h_{2s},h_{3s},\ldots).\]

By \cite[Example 8.60]{marquis2018introduction}, we have a natural representation $\tilde{\pi}_\cR:\widehat{\cU}^+_\cR\rightarrow \End(\cR(\!(t)\!)^2)$ and this representation restricts and corestricts to a faithful representation $\hat{\pi}_\cR:\fU^{ma+}(\cR)\rightarrow \mathrm{GL}_2(\cR(\!(t)\!))$. Via this identification, we have: \begin{equation}
[\exp](r h_{s})=\begin{pmatrix} 1/(1-rt^s) & 0 \\ 0 & 1-rt^s\end{pmatrix},
\end{equation} for $s\in \Z_{\geq 1}$. 

\begin{Lemma}\label{l_gd_pdcts}
The map $\psi:\cR^{\Z_{\geq 1}}\rightarrow 1+t\cR[\![t]\!]$ defined by $\psi((r_n))=\prod_{n\in \Z_{\geq 1}} (1-r_nt^n)$, for $(r_n)\in \cR^{\Z_{\geq 1}}$ is well-defined and is a bijection.
\end{Lemma}

\begin{proof}
The map $\psi$ is clearly well-defined. For $n\in \Z_{\geq 1}$, denote by $\pi_n:\cR[\![t]\!]\rightarrow  \cR[\![t]\!]/t^n \cR[\![t]\!]\simeq \cR\oplus  t\cR\oplus \ldots \oplus t^{n-1}\cR$ the natural projection.  Let $n\in \Z_{\geq 1}$.

Let $x\in 1+t\cR[\![t]\!]$.  Assume that there exists a unique $(r_1,\ldots,r_{n-1})\in \cR^{n-1}$ such that $\pi_n(x)=\pi_n(\prod_{j=1}^{n-1} (1-r_j t^j))$. By assumption, we can write $\prod_{j=1}^{n-1} (1-r_j t^j)=x+t^n f_n$, for some $f_n\in \cR[\![t]\!]$. Write $x=\sum_{j=0}^\infty x_n t^n$, with 
$x_n\in \cR$ for all $n\in \N$. Let $r_n\in \cR$.  Then: \[\pi_n(\prod_{j=1}^{n} (1-r_j t^j))=\pi_n((1-r_nt^n) (x+t^n f_n))=\pi_n(x_0+x_1 t+\ldots + x_nt^n-r_n t^n+t^n f_n^{(0)}),\] where we write $f_n=\sum_{m=0}^\infty f_n^{(m)} t^m$, with $f_n^{(m)}\in \cR$ for all $m\in \N$.  Therefore $\pi_n(\prod_{j=1}^{n} (1-r_j t^j))=\pi_n(x)$ if and only if $r_n=f_n^{(0)}$, and the lemma follows.
\end{proof}

Using Lemma~\ref{l_gd_pdcts} and \eqref{e_normal_form_Upma}, we deduce that: \begin{equation}
\label{e_U_Zdelta}\fU^{ma+}_{\Z_{\geq 1}\delta}(\cR)=\fU_{\Delta^{im}_+}(\cR)=\prod_{n\in \Z_{\geq 1}} X_{n\delta}(\ffg_{n\delta,\Z}\otimes \cR)=\begin{pmatrix}
1+t\cR[\![t]\!] & 0\\
0 	 & 1+t\cR[\![t]\!]
\end{pmatrix}\cap \mathrm{SL}_2(\cR[\![t]\!]). \end{equation}

\subsubsection{Description of $\fU^{ma+}(\cR)$}

By \cite[1.3.1.4 Corollary]{kumar2002kac}, we have: \[\Delta_+=\Z_{\geq 1}\delta\sqcup\{j\delta+\mathring{\beta}\mid j\in \Z_{\geq 1},\mathring{\beta}\in \mathring{\Phi}\}\sqcup \mathring{\Phi}_+.\] Set $\Delta_+(\mathring{\Phi}_+)=\Delta_+\cap (\mathring{\Phi}_+ +\Z \delta)$, $\Delta_+(\mathring{\Phi}_-)=\Delta_+\cap (\mathring{\Phi}_-+\Z \delta)$ and $\Delta^{im}_+=\Z_{\geq 1}\delta$. 

\begin{Proposition}\label{p_dec_U_aff}
Let $\cR$ be a ring. Let $A$ be a Kac--Moody matrix of untwisted affine type. Then $\fU^{ma+}(\cR)$ decomposes uniquely as:  \[\fU^{ma+}(\cR)=\fU^{ma}_{\Delta_+(\mathring{\Phi}_-)}(\cR)\cdot \fU^{ma}_{\Delta^{im}_+}(\cR) \cdot \fU^{ma}_{\Delta_+(\mathring{\Phi}_+)}(\cR).\]  
\end{Proposition}

\begin{proof}
The set $\Delta_+(\mathring{\Phi}_-)$ is strongly closed and strongly coclosed in $\Delta_+$. Therefore by Proposition~\ref{p_Mar8.58},  we have $\fU^{ma+}(\cR)=\fU^{ma}_{\Delta_+(\mathring{\Phi}_-)}(\cR)\cdot \fU^{ma}_{\Delta_+\setminus \Delta_+(\mathring{\Phi}_-)}(\cR)$. 

We have $\Delta_+\setminus  \Delta_+(\mathring{\Phi}_-)=\Delta^{im}_+(\cR)\sqcup \Delta_+(\mathring{\Phi}_+)$. Moreover,   $\Delta_+^{im}$ is closed and coclosed in $\Delta_+(\mathring{\Phi}_+)$, which proves that $\fU^{ma}_{\Delta_+\setminus \Delta_+}(\cR)=\fU^{ma}_{\Delta^{im}_+}(\cR)\cdot \fU^{ma}_{\Delta_+(\mathring{\Phi}_+)}(\cR)$, by Proposition~\ref{p_Mar8.58}. Proposition follows.
\end{proof}

We now assume that $A$ has type $A_1^{(1)}$: $A=\begin{psmallmatrix} 2 & -2\\ -2 & 2\end{psmallmatrix}$. Then by \cite[Example 7.7]{marquis2018introduction}, we can take $\mathring{\ffg}=\mathfrak{sl}_2(\C)$. We have $\mathring{\Phi}=\{\pm \mathring{\alpha}\}$ for a root $\mathring{\alpha}$ such that $\ffg_{\mathring{\alpha},\Z}=\begin{psmallmatrix} 0 & \Z \\ 0 	& 0\end{psmallmatrix}$ and $\ffg_{-\mathring{\alpha},\Z}=\begin{psmallmatrix} 0 & 0 \\ \Z 	& 0\end{psmallmatrix}$.  Then by Lemma~\ref{l_rt_dec_aff}, we have: \begin{equation}\label{e_root_SL2}
U_{\mathring{\alpha}+k\delta}(\cR)=\begin{psmallmatrix} 1 & t^k \cR\\ 0 & 1\end{psmallmatrix},
\end{equation} if $k\in \Z$. Similarly, if $k\in \Z$, we have: \begin{equation}\label{e_rt_-_SL2}
U_{-\mathring{\alpha}+k\delta}(\cR)=\begin{psmallmatrix} 1 & 0\\ t^k \cR & 0\end{psmallmatrix}. 
\end{equation}

\begin{Proposition}\label{p_desc_U_aff}
Assume that the Kac--Moody matrix $A$ is $\begin{psmallmatrix} 2 & -2 \\ -2 & 2\end{psmallmatrix}$. Let $\cR$ be a ring. Then under the identification $\fU^{ma+}(\cR)\subset \mathrm{GL}_2(\cR(\!(t)\!))$ of \cite[Example 8.59]{marquis2018introduction}, we have $\fU^{ma}_{\Delta_+(\mathring{\Phi}_-)}(\cR)=\begin{psmallmatrix} 1 & 0 \\ t\cR[\![t]\!] & 1\end{psmallmatrix}$, $\fU^{ma}_{\Delta_+(\mathring{\Phi}_+)}(\cR)=\begin{psmallmatrix} 1 & \cR[\![t]\!] \\ 0& 1\end{psmallmatrix}$ and $\fU^{ma}_{\Delta^{im}_+}(\cR)=\begin{psmallmatrix} 1+t\cR[\![t]\!] & 0 \\ 0 & 1+t\cR[\![t]\!] \end{psmallmatrix}\cap \mathrm{SL}_2(\cR[\![t]\!])$. 

We have \begin{equation}\label{e_des_U_aff}
\fU^{ma+}(\cR)=\begin{psmallmatrix} 1+t\cR[\![t]\!] & \cR[\![t]\!] \\
 t\cR[\![t]\!] & 1+t\cR[\![t]\!] \end{psmallmatrix}\cap \mathrm{SL}_2(\cR[\![t]\!]).
 \end{equation}
\end{Proposition}

\begin{proof}
By definition,  $\fU^{ma}_{\Delta_+(\mathring{\Phi}_-)}(\cR)$ consists of the elements that  can be written as an infinite product $\prod_{n\in \Z_{\geq 1}} u_n$, where $u_n\in U_{-\mathring{\alpha}+n\delta}(\cR)$, for $n\in \Z_{\geq 1}$. This proves that $\fU^{ma}_{\Delta_+(\mathring{\Phi}_-)}(\cR)=\begin{psmallmatrix} 1 & 0 \\ t\cR[\![t]\!] & 1\end{psmallmatrix}$. Similarly, we have $\fU^{ma+}_{\Delta_+(\mathring{\Phi}_+)}(\cR)=\begin{psmallmatrix} 1 & \cR[\![t]\!] \\ 0& 1\end{psmallmatrix}$ and the description of $\fU^{ma+}_{\Delta^{im}_+}(\cR)$ is given in \eqref{e_U_Zdelta}.

The right hand side of \eqref{e_des_U_aff} is a group and thus by Proposition~\ref{p_dec_U_aff}, we have the inclusion ``$\subset$''. Let $M=\begin{psmallmatrix} a & b\\ c & d \end{psmallmatrix}\in \begin{psmallmatrix} 1+t\cR[\![t]\!] & \cR[\![t]\!] \\
 t\cR[\![t]\!] & 1+t\cR[\![t]\!] \end{psmallmatrix}\cap \mathrm{SL}_2(\cR[\![t]\!])$.  Then $a\in (\cR[\![t]\!])^\times$ and $\begin{psmallmatrix} a' & b'\\ c' & d'\end{psmallmatrix}:=\begin{psmallmatrix} 1 & 0 \\ -c/a & 1\end{psmallmatrix} M\in \begin{psmallmatrix} 1+t\cR[\![t]\!] & \cR[\![t]\!] \\
 0 & 1+t\cR[\![t]\!] \end{psmallmatrix}\cap \mathrm{SL}_2(\cR[\![t]\!])$. Then  $\begin{psmallmatrix} a' & b'\\ c' & d'\end{psmallmatrix}\begin{psmallmatrix} 1 & -b'/a' \\ 0 & 1\end{psmallmatrix}=\begin{psmallmatrix} a' & 0 \\ 0 & d'\end{psmallmatrix}$.  Therefore: \[\begin{psmallmatrix} 1+t\cR[\![t]\!] & \cR[\![t]\!] \\
 t\cR[\![t]\!] & 1+t\cR[\![t]\!] \end{psmallmatrix}\cap \mathrm{SL}_2(\cR[\![t]\!])=\begin{psmallmatrix} 1 & 0 \\ t\cR[\![t]\!] & 1\end{psmallmatrix} \cdot \left(\begin{psmallmatrix} 1+t\cR[\![t]\!] & 0\\ 0 & 1+t\cR[\![t]\!]\end{psmallmatrix}\cap  \mathrm{SL}_2(\cR[\![t]\!])\right)\cdot \begin{psmallmatrix} 1 & \cR[\![t]\!] \\ 0 & 1\end{psmallmatrix}. \] Using the first part of the proposition, we deduce the sense $\supset$ of \eqref{e_des_U_aff}, which proves the proposition.
\end{proof}

\section{Completed Kac--Moody groups}\label{s_cplt_KM_gp}

Now that we defined the constructive Tits functor $\fG$ and the affine group scheme $\fU^{ma+}$, we can define completions of $\fG$: the scheme theoretic completion $\fG^{sch+}$ and Mathieu's completion $\fG^{ma+}$. The values of these functors are isomorphic over fields. We then use Mathieu's completed group to define  a minimal Tits functor.

\subsection{Marquis's scheme theoretic completion}

In \cite[8.5.4]{marquis2018introduction}, Marquis introduce the (positive) scheme theoretic completion $\fG^{\sch+}$ of $\fG$. Its is defined by generators and relations. For every field $\cK$, $\fG^{\sch+}(\cK)$ is isomorphic, as a topological group, to $\fG^{ma+}(\cK)$ (\cite[Proposition 8.120]{marquis2018introduction}).   The ``advantage'' of $\fG^{\sch+}$ is that its definition is simpler than the one of $\fG^{ma+}$.

Let $\cR$ be a ring. Let $i\in I$. Then by \cite[(8.55)]{marquis2018introduction}, we have two representations by continuous automorphisms $\Ad_\cR:\fU_{-\alpha_i}(\cR)\rightarrow \mathrm{Aut}(\fU_{\Delta_+\setminus\{\alpha_i\}}(\cR))$ and $\Ad_\cR:T\rightarrow \Aut(\fU^{ma+}(\cR))$.

\begin{Definition}\label{d_sch_th_KM}(see \cite[Definition 8.65]{marquis2018introduction}) Let $\cR$ be a ring. The \textbf{scheme-theoretic completion} of $\fG_{\cS}(\cR)$, denoted $\fG^{\sch+}_{\cS}(\cR)=\fG^{\sch+}(\cR)$ or simply $\fG^{\sch}_{\cS}(\cR)=\fG^{\sch}(\cR)$ is the quotient of the free product of $\fG_\cS(\cR)*\fU^{ma+}(\cR)$ by the following relations, where $i\in I$, $r\in \cR$ $t\in \fT(\cR)$ and $u\in \fU^{ma+}(\cR)\subset \widehat{\cU}_{\cR}^+$: \begin{align}
&\textbf{(R5)} &x_{\alpha_i}(r)=\exp(r e_i)\label{i_R5}\\
&\textbf{(R6)} &tut^{-1}=\mathrm{Ad}_{\cR}(t)(u)\label{i_R6}\\
&\textbf{(R7)} &x_{-\alpha_i}(r)u x_{-\alpha_i}(-r)=\mathrm{Ad}_{\cR}(x_{-\alpha_i}(r))(u).\label{i_R7}
\end{align} 
\end{Definition}

By \cite[Corollary 8.121]{marquis2018introduction}, if $\cR$ is a ring, then the canonical map $\fT\ltimes \fU^{ma+}(\cR)\rightarrow \fG^{\sch}$ is injective.

\subsection{Mathieu's group $\fG^{ma+}$}

\subsubsection{Sketch of the construction}

We now briefly recall the definition of Mathieu's Kac--Moody group. We refer to \cite{rousseau2016groupes} and \cite{mathieu1989construction} for more details.

The \textbf{Borel subgroup} (it will be a subgroup of $\fG^{ma+}$) is $\fB=\fB_{\cS}=\fT_{\cS}\ltimes \fU^{ma+}$\index[notation]{B@$\fB$}, where $\fT= \fT_{\cS}$ acts on $\fU^{ma+}$ as follows.  Let  $\cR$ be a ring, $\alpha\in \Delta_+$, $t\in \fT(\cR)$, $r\in \cR$ and $x\in \ffg_{\alpha,\cR}$. Then by \cite[3.5]{rousseau2016groupes}, we have:  \begin{equation}\label{e_KMT4}
t[\exp](rx)t^{-1}=[\exp](\alpha(t)rx).
\end{equation}
In particular, if $\alpha\in \Phi_+$, we have:   \begin{equation}\label{e_KMT4'}
\  tx_\alpha(r)t^{-1}=x_\alpha(\alpha(t)r).
\end{equation} More generally, if $\alpha\in \Delta_+$, $\underline{x}\in \ffg_{\alpha,\cR}$  and $t\in \cR$, we have \begin{equation}\label{e_com_T_X}
t X_\alpha(\underline{x})t^{-1}=X_{\alpha}(\alpha(t) \underline{x}).
\end{equation} Indeed, if $\underline{x}=\sum_{b\in \cB_\alpha}\lambda_b b$, with $(\lambda_b)\in \cR^{\cB_\alpha}$, then $X_\alpha(\underline{x})=\prod_{b\in\cB_\alpha}[\exp] \lambda_b b$. Then $tX_\alpha(\underline{x})t^{-1}=\prod_{b\in \cB_\alpha} t[\exp] (\lambda_b b]) t^{-1}=\prod_{b\in \cB_\alpha} [\exp](\alpha(t)\lambda_b b),$ by \eqref{e_KMT4}.

Let $i\in I$. We denote by $\fA^X_i$ the affine group scheme over $\Z$ associated to  the Kac--Moody datum $((2),X,Y,(\alpha_i),(\alpha_i^\vee))$. In other words, if $X_i$ is the $\Z$ dual of a complement of $\Z\alpha_i^\vee$ in $\Z$, then: \begin{equation}\label{e_fUi_X}
\fA_i^{X}=\fT_{X_i}\ltimes \mathrm{SL}_2
\end{equation} (see \cite[Exercise 7.33 (5)]{marquis2018introduction}). For each $i\in I$, Mathieu defines an (infinite dimensional) affine group scheme $\fP_{i}=\fU_{\alpha_i}^Y \ltimes \fU^{ma}_{\Delta_+\setminus \{\alpha_i\}}$ (see \cite[Definition 8.65]{marquis2018introduction} for the definition of the action of $\fU_{-\alpha_i}$ on $\fU^{ma}_{\Delta_+\setminus \{\alpha_i\}}$).

 We do not detail the definition of $\fG^{ma+}$\index[notation]{G@$\fG^{ma+}$} and we refer to \cite{mathieu1989construction}, \cite[8.7]{marquis2018introduction} or \cite[3.6]{rousseau2016groupes}. This is an ind-group scheme containing the $\fP_{i}$ for every $i\in I$. Let $w\in W^v$ and write $w=r_{i_1}\ldots r_{i_k}$, with $k=\ell(w)$ and $i_1,\ldots,i_k\in I$. Then the multiplication map $\fP_{i_1}\times \ldots  \times \fP_{i_k}\rightarrow \fG^{ma+}$ is  a  scheme morphism, and we have $\fG^{ma+}(\cR)=\bigcup_{(i_1,\ldots,i_n)\in \mathrm{Red}(W^v)} \fP_{i_1}(\cR)\times \ldots \times \fP_{i_k}(\cR)$, where $\mathrm{Red}(W^v)$ is the set of reduced words of $W^v$ (i.e $\mathrm{Red}(W^v)=\{(i_1,\ldots,i_k)\in I^{(\N)}\mid \ell(r_{i_1}\ldots r_{i_k})=k\}$). 

Let $w\in W^v$, $i\in I$ and $\alpha=w.\alpha_i$. We set $\fU_{\alpha}=\tilde{w}.\fU_{\alpha_i} .\tilde{w}^{-1}$, where $\tilde{w}$ is defined in \eqref{e_tilde_w}. There is an isomorphism of group schemes $x_\alpha: \mathbb{G}_a\rightarrow \fU_{\alpha}$ (see \cite[page 262]{marquis2018introduction}). The group $\fG^{ma+}$ is generated by the $\fP_i$, $i\in I$. Moreover, if $i\in I$, then $\fP_i$ is generated by $\fT$, $\fU_{\pm \alpha_i}$ and $\tilde{r}_{i}=x_{\alpha_i}(1)x_{-\alpha_i}(1)x_{\alpha_i}(1)$. Thus $\fG^{ma+}$ is generated by $\fU^{ma+}$, $\fT$, $\fU_{-\alpha_i}$ and the $\tilde{r}_{i}$, for $i\in I$  and thus we have: \begin{equation}\label{e_generators_Gpma}
\fG^{ma+}=\langle \fU^{ma+},\fT,\fU_\alpha,\alpha\in \Phi_-\rangle.
\end{equation}

There is a group functor morphism $\iota:\fG\rightarrow \fG^{ma+}$ such that for any ring $\cR$, $\iota_\cR$ maps $x_\alpha(r)$ to $x_\alpha(r)$ and $t$ to $t$, for each $\alpha\in \Phi$, $r\in \cR$, $t\in \fT(\cR)$. When $\cR$ is a field, this morphism is injective (see \cite[3.12]{rousseau2016groupes} or \cite[Proposition 8.117]{marquis2018introduction}).

Note that the relations (R5) to (R7) (\eqref{i_R5} to \eqref{i_R7}) are satisfied in $\fG^{ma+}(\cR)$, by \cite[Proposition 8.120]{marquis2018introduction}.

\subsubsection{Conjugation by elements of the Weyl group}

Recall that for $n\in \N$, $\fU_{n}^{ma+}(\cR)=\fU_{\Psi(n)}^{ma+}(\cR)$ where $\Psi(n)=\{\alpha\in \Delta_+\mid  \htt(\alpha)\geq n\}$.  Then by \cite[Lemma 8.67 and Definition 8.119]{marquis2018introduction}  the sequence $(\fU_{n}^{ma+}(\cR))$ is a separated conjugation-invariant filtration of $\fG^{\sch+}(\cR)$ and $\fG^{ma+}(\cR)$, in the sense of \cite[Exercise 8.5]{marquis2018introduction}. It equips $\fG^{\sch+}(\cR)$ and $\fG^{ma+}(\cR)$ with the structure of a topological group, where a basis of neighborhood of the identity is given by the $\fU_n(\cR)$, $n\in \N$. 

Recall the definitions of $w^*$ and $\tilde{w}$ in \eqref{e_star_w} and \eqref{e_tilde_w}.

\begin{Lemma}\label{l_cnj_X_w}
Let $\alpha \in \Delta_+$. Let $\underline{u_\alpha}\in \ffg_{\alpha}\otimes \cF$ and $w\in W^v$. Then $\tilde{w} X_{\alpha}^{\cB_{\alpha}}(\underline{u_\alpha}) \tilde{w}^{-1}=X_{w.\alpha}^{w^*.\cB_\alpha}(w^*.\underline{u_\alpha})$.
\end{Lemma}

\begin{proof}
First assume $\alpha\in \Inv(w)$. Then $\ffg_{\alpha,\Z}=\Z e_\alpha$ has rank $1$, and by \ref{a_KMT7}, $\tilde{w}X_{\alpha}(\underline{u_\alpha})\tilde{w}^{-1} =\tilde{w} x_\alpha( \underline{u_\alpha})\tilde{w}^{-1}=x_{w.\alpha}(w^*. \underline{u_\alpha})=X_{w.\alpha}(w^*. \underline{u_\alpha})$.

Assume now $\alpha\in \Delta_+\setminus \Inv(w)$.  Write $\underline{u_\alpha}=\sum_{b\in \cB_\alpha}r_b b$, with $r_b\in \cR$, for $b\in \cB_\alpha$. Then $X_\alpha^{\cB_{\alpha}}(\underline{u_\alpha})=\prod_{b\in \cB_\alpha}[\exp](r_b b)$. 

By \eqref{e_Kmr_1.3.14}, $\alpha\neq \alpha_{i_1}$, $r_{i_1}.\alpha\neq \alpha_{i_2},\ldots, r_{i_1}\ldots r_{i_{k-1}}.\alpha\neq \alpha_{i_k}$. Thus we can apply \cite[Lemma 8.77]{marquis2018introduction} inductively and we have $\tilde{w}[\exp](r_b b)\tilde{w}^{-1}=[\exp ] (r_b w^*.b)$, for $b\in \cB_\alpha$. Therefore $\tilde{w} X_{\alpha}^{\cB_\alpha}(\underline{u_\alpha})\tilde{w}^{-1}=X_{w.\alpha}^{w^*.\cB_\alpha}(w^*.\underline{u_\alpha})$.
\end{proof}

Let $(\underline{u_\alpha})\in \prod_{\alpha\in \Delta_+}(\ffg_{\alpha,\Z}\otimes \cR)$. Then by Lemma~\ref{l_cnj_X_w}, \begin{equation}\label{e_cnj_X_tilde_w}
\tilde{w}\prod_{\alpha\in \Delta_+} X_\alpha^{\cB_\alpha}(\underline{u_\alpha})\tilde{w}^{-1}=\prod_{\alpha\in \Delta_+}X_{w.\alpha}^{w^*.\cB_\alpha}(w^*.\underline{u_\alpha}),
\end{equation} and this product is well-defined in $\fG^{ma+}(\cR)$ or $\fG^{\sch+}(\cR)$ since  for all  $n\in \N$, $\{\alpha\in \Delta_+\mid \htt(w.\alpha)\leq n\}$ is finite.

The following proposition holds in $\fG^{\sch+}$ and $\fG^{ma+}$.

\begin{Proposition}\label{p_dec_inv}
Let $\cR$ be a ring and  $w\in W^v$. Then: \[\fU^{ma}_{\Inv(w)}(\cR)\subset \fU^{ma+}(\cR)\cap \tilde{w}^{-1} \fU^-(\cR)\tilde{w}\text{ and }\fU^{ma}_{\Delta_+\setminus \Inv(w)}(\cR)\subset \fU^{ma+}(\cR)\cap \tilde{w}^{-1} \fU^{ma+}(\cR)\tilde{w}.\] If moreover $\cR=\cK$ is a field, then these inclusions are actually equalities.
\end{Proposition}

\begin{proof}
Let $u\in \fU^{ma}_{\Inv(w)}(\cR)$. Then by Corollary~\ref{c_dec_inv}, we can write $u=\prod_{\alpha\in \Inv(w)}x_\alpha(a_\alpha)$, for some $(a_\alpha)\in \cR^{\Inv(w)}$. By \eqref{e_cnj_X_tilde_w} or (KM7), we have $\tilde{w} u\tilde{w}^{-1}=\prod_{\alpha\in \Inv(w)}x_{w.\alpha}(w^*.a_\alpha)\in \fU^-(\cR)$, which proves that $\fU^{ma}_{\Inv(w)}(\cR)\subset \fU^{ma+}(\cR)\cap \tilde{w}^{-1}\fU^{-}(\cR)\tilde{w}$. 

Let $u\in \fU^{ma}_{\Delta_+\setminus \Inv(w)}(\cR)$. We can write $u=\prod_{\alpha\in \Delta_+\setminus \Inv(w)} X_\alpha^{\cB_\alpha}(\underline{u_\alpha})$, where $\underline{u_\alpha}\in \ffg_{\alpha,\Z}\otimes \cR$, for $\alpha\in \Delta_+\setminus \Inv(w)$. Let $\alpha\in \Delta_+\setminus \Inv(w)$.  By \eqref{e_cnj_X_tilde_w}, we have $\tilde{w} u\tilde{w}^{-1}=\prod_{\alpha\in \Delta_+}X_{w.\alpha}^{w^*.\cB_\alpha}(w^*.\underline{u_\alpha})\in \fU^{ma+}(\cR)$,  which proves that $\fU^{ma+}_{\Delta_+\setminus \Inv(w)}(\cR)\subset \fU^{ma+}(\cR)\cap  \tilde{w}^{-1}\fU^{ma+}(\cR)\tilde{w}$. 

Assume now that $\cR=\cK$ is a field. Let $u\in \fU^{ma+}(\cK)\cap \tilde{w}^{-1}\fU^-(\cK)\tilde{w}$ (resp. $u\in \fU^{ma+}(\cK)\cap \tilde{w}^{-1}\fU^{ma+}(\cK)\tilde{w}$). 
Using Corollary~\ref{c_dec_inv}, we write $u=u_1 u_2$, with $u_1\in \fU^{ma+}_{\Inv(w)}$ and $u_2\in \fU^{ma+}_{\Delta_+\setminus \Inv(w)}(\cK)$. Then by what we just proved $\tilde{w}u_2\tilde{w}^{-1}\in \fU^{ma+}(\cK)$ (resp. $\tilde{w}u_1\tilde{w}^{-1}\in \fU^{-}(\cK)$). Moreover $\tilde{w} u_2 \tilde{w}^{-1}=(\tilde{w}u \tilde{w}^{-1})^{-1} \tilde{w}u_1\tilde{w}^{-1}\in \fU^-(\cK)$ (resp. $\tilde{w} u_1\tilde{w}^{-1}\in \fU^{ma+}(\cK)$). Therefore $\tilde{w} u_2\tilde{w}^{-1}\in \fU^-(\cK)\cap \fU^{ma+}(\cK)$ (resp. $\tilde{w} u_1\tilde{w}^{-1}\in \fU^-(\cK)\cap \fU^{ma+}(\cK)$). Using \cite[3.16 Remarque]{rousseau2016groupes}, we deduce that $\tilde{w} u_2 \tilde{w}^{-1}=1$ (resp. $\tilde{w} u_1\tilde{w}^{-1}=1$) and thus $u=u_1\in \fU^{ma+}_{\Inv(w)}(\cK)$ (resp. $u=u_2\in \fU^{ma+}_{\Delta_+\setminus\Inv(w)}(\cK)$), which prove the equalities. 
\end{proof}

Note that by (R7) (\eqref{i_R7}), if $i\in I$, then $\fU_{-\alpha_i}(\cR)$ normalizes $\fU^{ma+}_{\Delta_+\setminus\{\alpha_i\}}(\cR)$.

For $w\in W^v$, we set: \begin{equation}\label{e_fU_ma_w}
\fU^{ma}_{w.\Delta_+}(\cR)=\tilde{w}\fU^{ma+}(\cR)\tilde{w}^{-1}\subset \fG^{ma+}(\cR).
\end{equation} By Proposition~\ref{p_dec_inv}, if $\cR=\cK$ is a field, then  we have  $\fU^{ma+}(\cK)\cap \fU_{w.\Delta_+}^{ma}(\cK)=\fU_{\Delta_+\cap w.\Delta_+}(\cK)$.

Let $\Psi\subset \Delta$. We assume that for some $w\in W^v$, we have $w.\Psi\subset \Delta_+$ and $w.\Psi$ is closed. Then we set:  \begin{equation}
\fU^{ma+}_{\Psi}(\cR)=\{\prod_{\alpha\in \Psi} X_\alpha(\underline{u_\alpha})\mid (\underline{u_\alpha})\in \prod_{\alpha\in \Psi}(\ffg_{\alpha,\Z}\otimes \cR)\}\subset \fU^{ma+}(\cR).
\end{equation} 

Then by \eqref{e_cnj_X_tilde_w}, we have \begin{equation}\label{e_cnj_U_Wv}
\tilde{w} \fU_{\Psi}^{ma+}(\cR)\tilde{w}^{-1}=\fU^{ma}_{w.\Psi}(\cR), \forall w\in W^v.
\end{equation}

\subsubsection{Mathieu's group and ring morphisms}

\begin{Proposition}\label{p_mrphsm}
Let $\cR$ and $\cR'$ be two rings and $\varphi:\cR\rightarrow \cR'$ be a ring morphism. 

 Let $f_\varphi^{\widehat{\cU}^+}:\widehat{\cU}_\cR^+\rightarrow \widehat{\cU}_{\cR'}^+$ and $f_\varphi:\fG^{pma}(\cR)\rightarrow \fG^{pma}(\cR')$ be the induced morphisms.  Then $f_\varphi^{\widehat{\cU}^+}(\fU^{pma}(\cR))\subset \fU^{pma}(\cR')$ and we have: \begin{enumerate}
\item Let $\cB=\bigcup_{\alpha\in \Delta}\cB_{\alpha}$. Let $\cB'\subset \cB$. Then  for every $(r_x)\in \cR^{\cB'}$, we have $f_\varphi^{\widehat{\cU}^+}\left(\prod_{x\in \cB}[\exp](r_x x)\right)=\prod_{x\in \cB}[\exp](\varphi(r_x) x)$.

\item For $\alpha\in \Delta_+$ and $(\lambda_x)\in \cR^{\cB_\alpha}$, we have $f_\varphi\left(X_\alpha(\sum_{x\in \cB_\alpha} \lambda_x x)\right)=X_\alpha\left(\sum_{x\in \cB_\alpha}\varphi(\lambda_x) x\right)$. 

\item We have $f_{\varphi}(u)=f_\varphi^{\widehat{\cU}^+}(u)$ for $u\in \fU^{pma}(\cR)$, $f_{\varphi}(x_{\alpha}(r))=x_\alpha(\varphi(r))$, for $\alpha\in \Phi$ and $r\in \cR$, and $f_\varphi(\chi(r))=\chi(\varphi(r))$, for $\chi\in Y$ and $r\in \cR^\times$. 

\item If $\varphi$ is surjective, then $f_\varphi$ is surjective.
\end{enumerate}
\end{Proposition}

\begin{proof}
(1), (2) By definition, we have: \[f_\varphi^{\widehat{\cU}^+}(\sum_{\alpha\in Q^+} \sum_{j\in J_\alpha} u_{\alpha,j}\otimes r_j)=\sum_{\alpha\in Q^+} \sum_{j\in J_\alpha} u_{\alpha,j}\otimes \varphi(r_j)\] if $J_\alpha$ is a finite set and $(r_j)\in \cR^{J_\alpha}$ and $u_{\alpha,j}\in \cU_{\alpha,\cR}$, for every $\alpha\in Q_+$. Thus $f_{\varphi}^{\widehat{\cU}^+}$ commutes with infinite sums and product, which proves (1) and (2). 

(3) Let $i\in I$. Then the morphism $\fP_i(\cR)\rightarrow \fP_i(\cR')$ induced by $\varphi$ satisfies the formula above. Using the fact that $x_\alpha=\tilde{w} x_{-\alpha_i} \tilde{w}^{-1}$, for $\alpha=-w.\alpha_i$, with $w\in W^v$, $i\in I$ and $\tilde{w}$ defined as in \eqref{e_tilde_w}, we have (3).

(4)  Assume $\varphi$ is surjective. By \eqref{e_normal_form_Upma} and (1), the restriction of $f_\varphi$ to $\fU^{pma}(\cR)$ is surjective. By (3), $f_{\varphi}(\fU^-(\cR))=\fU^-(\cR')$ and $f_\varphi(\fT(\cR))=\fT(\cR')$. We conclude by using the fact that $\fG^{pma}$ is generated by $\fU^{pma}$, $\fU^-$ and $\fT$ (see \eqref{e_generators_Gpma}).
\end{proof}

\subsection{Minimal Kac--Moody groups over rings}\label{ss_min_KM_gps}

\subsubsection{Axiomatic properties of Tits functor}

Let $\cS$ be  a Kac--Moody root datum. In Subsection~\ref{ss_Tits_fnctr}, we introduced the constructive Tits functor $\fG_\cS$. In some sense, it integrates the Kac--Moody algebra $\ffg_{\cS}$. In order to give a precise meaning to ``integrating $\ffg_\cS$'', Tits defines a series of fives axioms ((KM1) to (KMG5)) that a group functor has to satisfy to be a ``minimal Kac--Moody group functor of type $\cS$'' (see \cite{tits1987uniqueness} and also \cite[8]{remy2002groupes}).  Tits claims that if $\fG$ is a group functor satisfying its axioms, then there exists a morphism of group functors $\pi:\fG_{\cS}\rightarrow \fG$ satisfying certain conditions and that  under some mild conditions, $\pi$ is an isomorphism for each field (see \cite[Theorem 1]{tits1987uniqueness}). However he only sketches the proofs and we do not know a reference where it is proved in details.   Marquis proposes in \cite[7.5]{marquis2018introduction} a slight modification of Tits's axioms: he defines axioms (KMG1'), (KMG2'), (KMG4') and (KMG5'), which are close to Tits's axioms. We call such a functor a \textbf{Tits functor}. He gives a detailed proof of \cite[Theorem 1]{tits1987uniqueness} in certain cases (when $\cS$ is adjoint). He also proves  the existence of a Tits functor (\cite[8.8]{marquis2018introduction}) $\fG^{\min}_\cS=\fG^{\min}$. Over arbitrary rings, its definition use Mathieu's group and we recall it here. This functor has the  properties (Gmin1) to (Gmin4) and (KMG5) below  (see \cite[2.6]{marquis2025presentation}).

Let $\rE_2$ be the \textbf{elementary subgroup functor} of $\mathrm{SL}_2$. It is defined by $\rE_2(\cR)=\langle \begin{psmallmatrix} 1 & \cR\\ 0 & 1\end{psmallmatrix} , \begin{psmallmatrix} 1 & 0\\ \cR & 1\end{psmallmatrix}\rangle$, if $\cR$ is  a ring. We have $\rE_2(\cR)=\mathrm{SL}_2(\cR)$, if $\cR$ is Euclidean, by \cite[Exercise 7.2 (3)]{marquis2018introduction}.
\begin{enumerate}
\item[(Gmin1)]
$\fG_{\cS}^{\min}$ comes equipped with group functor morphisms $\varphi_i: \mathrm{SL}_2\to\fG_{\cS}^{\min}$ ($i\in I$) and $\eta:\fT_{\Lambda}\to \fG_{\cS}^{\min}$ such that for each ring $\cR$, the group morphism $\eta_\cR:\fT(\cR)\to \fG_{\cS}^{\min}(\cR)$ is injective. We then identify $\fT(\cR)$ with a subgroup of $G^{\min}_\cR:=\fG_{\cS}^{\min}(\cR)$. 
\item[(Gmin2)]
There is a (unique) group functor morphism $\varphi:\fG_{\cS}\to \fG_{\cS}^{\min}$ such that for each ring $\cR$, the restriction of $\varphi_\cR: \fG_{\cS}(\cR)\to \fG_{\cS}^{\min}(\cR)$ to $\fT(\cR)$ is the identity, and for each $i\in I$,
$$\varphi_{i\cR}\begin{psmallmatrix}1&r\\ 0&1\end{psmallmatrix}=\varphi_\cR(x_i(r))\quad\textrm{and}\quad \varphi_{i\cR}\begin{psmallmatrix}1&0\\ -r&1\end{psmallmatrix}=\varphi_\cR(x_{-i}(r)) \quad\textrm{for all $r\in \cR$.}$$ 
The morphism $\varphi_\cR$ is injective on each $\fU_{\gamma}( \cR)$ ($\gamma\in\Phi$), and we keep the notations $\fU_{\gamma}(\cR)$, $x_{\gamma}(r)$, $x_{\pm i}(r)$ and $\widetilde{s}_i$ ($i\in I$) for the corresponding objects in $G^{\min}_\cR$. We also set $\overline{G}_{i\cR}:=\varphi_\cR(G_{i\cR})$. Note that $$\varphi_{i\cR}\begin{psmallmatrix}0&1\\ -1&0\end{psmallmatrix}=\widetilde{s}_i\quad\textrm{and}\quad \overline{G}_{i\cR}=\varphi_{i\cR}(\mathrm{SL}_2(\cR))\quad\textrm{for all $i\in I$.}$$
\item[(Gmin3)]
If $\cR_1\to \cR_2$ is an injective ring morphism, then the group morphism $\fG_{\cS}^{\min}(\cR_1)\to\fG_{\cS}^{\min}(\cR_2)$ is injective.
\item[(Gmin4)]
If $\cK$ is a field, then $\varphi_{\cK}: \fG_{\cS}(\cK)\to \fG_{\cS}^{\min}(\cK)$ is an isomorphism. We then identify $\fG_\cS(\cK)$ and $\fG^{\min}_\cS(\cK)$ and we simply write $\fG_\cS(\cK)$.

\end{enumerate}

\subsubsection{Definition of a Tits functor}

Let $i\in I$. We denote by $\fA^X_i$ the affine group scheme over $\Z$ associated to  the Kac--Moody datum $((2),X,Y,(\alpha_i),(\alpha_i^\vee))$ (see \eqref{e_fUi_X}). Let $\varphi_i:\mathrm{SL}_2\rightarrow \fA_{i}^X$\index[notation]{P@$\varphi_i$} be the natural inclusion morphism. 

\begin{Definition}\label{d_G_min}
For a ring $\cR$, we set \[\fG^{\min}_{\cS}(\cR)=\fG^{\min}(\cR)=\langle \varphi_i\left(\mathrm{SL}_2(\cR)\right),\fT(\cR)\rangle\subset \fG^{ma+}(\cR).\]\index[notation]{G@$\fG^{\min}$} This is the \textbf{minimal group functor of type 
$\cS$}.
\end{Definition}

This group functor is introduced by Marquis in  \cite[Definition 8.126]{marquis2018introduction}. By \cite[Proposition 8.129]{marquis2018introduction}, it is a nondegenerate Tits functor in the sense of \cite[Definition 7.83]{marquis2018introduction} and we have $\fG^{\min}(\cR)\simeq \fG(\cR)$, for any field $\cR$. 

Note that if $\varphi$ is a ring morphism between two rings $\cR$ and $\cR'$, the induced morphism $\fG^{ma+}(\cR)\rightarrow \fG^{ma+}(\cR')$ restricts to a morphism $\fG^{\min}(\cR)\rightarrow \fG^{\min}(\cR')$. 

Let $\cR$ be a semilocal ring, i.e a ring with finitely many maximal ideals, then by \cite[4.3.9 Theorem]{hahn1989classical}, $\mathrm{SL}_2(\cR)$ is generated by $\begin{psmallmatrix} 1 & \cR\\ 0 & 1\end{psmallmatrix}$ and  $\begin{psmallmatrix} 1 & 0\\ \cR & 1\end{psmallmatrix}$. Therefore, \begin{equation}\label{e_minimal_group_semilocal_ring}
\fG^{\min}(\cR)=\langle \fU_{\pm\alpha_i}(\cR),\fT(\cR)\mid i\in I\rangle\subset \fG^{ma+}(\cR).
\end{equation}

\begin{Example}
If $n\in\N_{\geq 2}$, then the functor $k\mapsto \mathrm{SL}_n(k)$ is an example of Tits functor, where $\fT(k)$ is the group of diagonal matrices, $\fU(k)$ (resp. $\fU^-(k)$) is the group of upper (resp. lower) triangular matrices having $1$ on the diagonal. It is associated to the Lie algebra $\mathfrak{sl}_n(\C)$. We keep the same notation as in Section~\ref{s_example_sln}. If $i,j\in \llbracket 1,n-1\rrbracket$ is such that $i<j$ and $\alpha=\alpha_{i,j}$, one can choose $\fU_\alpha(k)=I+k E_{i,j}$ and $\fU_{-\alpha}(k)=I+kE_{j,i}$.
\end{Example}

\chapter{Abstract masure: a convexity result}\label{C_w_mas}

In Chapter~\ref{C_prop}, we will define the masure associated with a Kac--Moody group $G$, as the quotient $\I=(G\times \A)/\sim$, where $\A$ is the apartment associated with $G$ and $\sim$ is an equivalence relation. This definition is complicated, but we can recover many properties of $\I$ by considering it as an abstract masure. Roughly speaking (see Definition~\ref{d_w_mas} for a precise definition), this means  that $\I$  is a union of subsets of $\I$ called apartments  satisfying the following properties: \begin{enumerate}
                                                                                                                                                                                                                                                                                                                                                                                          \item[(MA1)] each apartment is ``isomorphic'' to the standard apartment,

                                                                                                                                                                                                                                                                                                                                                                                          \item[(MA2)] for every two apartments $A,B$ of $\I$, $A\cap B$ is a finite intersection of half-spaces of $A$, delimited by real roots of $G$, and there exists an apartment isomorphism from $A$ to $B$ fixing $A\cap B$,

                                                                                                                                                                                                                                                                                                                                                                                          \item[(MA3)] some pairs of filters $\cV_1,\cV_2$ are contained in  a common apartment.                                                                                                                                                                                                                                                                                                                                                                                       \end{enumerate}

Note that when $\I$ is associated with a split Kac--Moody group $G$ over a valued field with split maximal torus $T$, the standard apartment is $Y\otimes \R$ equipped with some set of walls, where $Y$ is the cocharacter lattice of $(G,T)$. The apartments are then the $g.\A$ such that $g\in G$. 

                                                                                                                                                                                                                                                                                                                                                                                                                                                                                                                                                                                                                                                                                                                                                                                                                                                                                                                                                                                                                                                                                                                    While in the reductive case, axiom~\ref{a_ma2} is rather easy to obtain, new difficulties appear in the Kac--Moody case. The main ones are the fact that the Tits cone is a proper subcone of $\A$ and  the fact that the construction of $\I$ involves completions of $G$, which have imaginary roots. Actually in  \cite{rousseau2011masures} and \cite{rousseau2016groupes}, Rousseau proved that $\I$ satisfies weak versions of~\ref{a_ma2}. We then obtained that these weak versions actually imply~\ref{a_ma2} in  \cite{hebert2020new} and \cite{hebert2022new}.

In this chapter, we define abstract masures. We then prove that a certain series of axioms imply~\ref{a_ma2}  (see Theorem~\ref{t_MA2}). We work abstractly, which means that we consider a set covered by apartments and satisfying certain axioms, but we do not assume that this set is attached to a Kac--Moody group. We assume that  the Kac--Moody data has a cofree family of simple roots to obtain Theorem~\ref{t_MA2}. However, in the case of a masure associated with a Kac-Moody group, we will drop this assumption.

While a part of this chapter is purely building theoretic (from Subsection~\ref{s_chrctz_apt}), we introduce objects which are particularly important in the theory of Kac--Moody groups over valued fields, namely retractions and Hecke paths.

\section{Axioms}\label{s_axioms}

\subsection{Apartment of type $\underline{\A}$}\label{ss_aff_apt}

In this subsection, we define the appartment of a masure. Let $A=(a_{i,j})_{i,j\in I}$ be a Kac--Moody matrix and $\cS=(A,X,Y,(\alpha_i)_{i\in I},(\alpha_i^\vee))$ be a free Kac--Moody datum, as defined in Subsection~\ref{ss_KM_data}.

\subsubsection{Affine Weyl group and real walls}\label{sss_af_Wl_gp}   Let $\A=Y\otimes \R$.  As the results of this chapter also apply in the almost split case (see \cite{rousseau2017almost}), 
we use a slightly more general framework here. This introduces a bit of technicality. The reader willing to avoid it will not lose much in making the following assumptions, which are true in the split case: \begin{enumerate}
\item There exists a subgroup $\Lambda$ of $\R$  such that $\Lambda_\alpha=\Lambda$, for all $\alpha\in \Phi$ (and $\Lambda_\alpha=\R$ for all $\alpha\in \Delta\setminus \Phi$). The subgroup $\Lambda$ is then $\omega(\cF^\times)$ if the masure is associated with a split group over the valued field $(\cF,\omega)$.

\item The set of real walls is $\cM=\{\alpha^{-1}(\{k\})\mid \alpha\in \Phi, k\in \Lambda\}$,

\item The affine Weyl group is $W^a:=W^v\ltimes Q^\vee$, where $Q^\vee=\bigoplus_{i\in I} \Lambda\alpha_i^\vee$, and the group $\widehat{W^a}$ is either $W^a$ itself or the extended affine Weyl group $\tilde{W^a}=W^v\ltimes (\Lambda\otimes Y)$.
\end{enumerate}

For $\alpha\in \A^*$ and $k\in \R$, set $M(\alpha,k)=\{v\in \A\mid  \alpha(v)+k=0\}$\index[notation]{m@$M(\alpha,k)$}, $D(\alpha,k)=\{v\in \A\mid \alpha(v)+k\geq 0\}$\index[notation]{d@$D(\alpha,k)$} and $D(\alpha,+\infty)=\A$. 

Recall that if $\alpha\in \Phi$, then $r_\alpha$ is the reflection of $W^v$ fixing $\alpha^{-1}(\{0\})$. We set $W^a_\R=W^v \ltimes \A$\index[notation]{w@$W^a$}.  If $a\in \A$, we denote by $\bt_a$ the translation on $\A$ by the vector $a$: $\bt_a(x)=x+a$, for $x\in \A$. We have $\bt_a. M(\alpha,k)=M(\alpha,k-\alpha(a))$, for $\alpha\in \Phi$  and $a\in \A$. We set  $\cM_\R=\{\alpha^{-1}(\{\lambda\})\mid \lambda\in \R\}$. If $M=M(\alpha,k)\in \cM_\R$, the associated reflection $r_M$\index[notation]{r@$r_M$} is the reflection of $W^a_\R$ fixing $M$. Explicitly, $r_M=\bt_{k\alpha^\vee} r_\alpha$.

 Let $\sL$ be the set of families $(\Lambda_\alpha)_{\alpha\in \Phi}$  consisting of the infinite subsets of $\R$ containing $0$ and  satisfying: \begin{equation}\label{e_cond_walls}
 \Lambda_\beta=\Lambda_\beta-\beta(\alpha^\vee)\Lambda_\alpha,\forall \alpha,\beta\in \Phi.
 \end{equation}

\begin{Definition}
Let $(\Lambda_\alpha)\in \sL$. The set  $\cM=\cM((\Lambda_\alpha))=\{\alpha^{-1}(\lambda)\mid \alpha\in \Phi,\lambda\in \Lambda_\alpha\}$\index[notation]{m@$\cM$} is  the set of \textbf{affine walls} of $\A$ associated with $(\Lambda_\alpha)$. The group $W^a=W^a((\Lambda_\alpha))=\{r_M\mid M\in \cM\}$\index[notation]{w@$W^a$} is the \textbf{affine Weyl group} associated with $(\Lambda_\alpha)$. 
\end{Definition}

Let $(\Lambda_\alpha)\in \sL$. Condition~\ref{e_cond_walls} is equivalent to the fact that $r_M.\cM=\cM$ for all $M\in \cM$. Therefore $W^a$ stabilizes $M$.

If $\alpha\in \Phi$, we denote by  $\tilde{\Lambda}_\alpha$\index[notation]{l@$\tilde{\Lambda_\alpha}$} the subgroup of $\R$ generated by $\Lambda_\alpha$.  By~\eqref{e_cond_walls}, $\Lambda_\alpha=\Lambda_{\alpha}+2\tilde{\Lambda}_\alpha$ for all $\alpha\in \Phi$. In particular, $\Lambda_\alpha=-\Lambda_{\alpha}$ and when $\Lambda_\alpha$ is discrete, $\tilde{\Lambda}_\alpha=\Lambda_\alpha$ is isomorphic to $\Z$.

We set  $Q^\vee= \sum_{\alpha\in \Phi} \tilde{\Lambda}_\alpha \alpha^\vee$. This is a subgroup of $\A$ stable under the action of $W^v$. Then we have: \[W^a=W^v\ltimes Q^\vee.\] An element of $W^a$ is called a \textbf{Weyl automorphism}\index{Weyl automorphism} of $(\A,(\Lambda_\alpha))$.

We set: \begin{equation}\label{e_Aut_vw}
\Aut^{\vw}(\A,(\Lambda_\alpha))=\{\bw \in W^a_\R\mid \bw.\cM((\Lambda_\alpha))=\cM((\Lambda_\alpha))\}.
\end{equation}\index[notation]{a@$\Aut^{\vw}(\A,(\Lambda_\alpha))$}
 We call it the \textbf{group of vectorially Weyl automorphisms} of $(\A,(\Lambda_\alpha))$ and an element of this group is called a \textbf{vectorially Weyl automorphism}\index{vectorially Weyl automorphism}.

\begin{Definition}\label{d_apt}
An \textbf{affine apartment} $\underline{\A}$ (or simply an \textbf{apartment}) is a quadruple $\underline{\A}=(\A,\cS,(\Lambda_\alpha),\widehat{W^a})$ such that $\cS$ is a  Kac--Moody datum, $(\Lambda_\alpha)_{\alpha\in \Phi}\in \sL$ and  $\widehat{W^a}$ is a group of vectorially Weyl automorphisms  such that $W^a((\Lambda_\alpha))\subset \widehat{W^a}\subset  \Aut^{\vw}(\A,(\Lambda_\alpha))$.  

 We usually denote it $\underline{\A}=(\A,\cS,(\Lambda_\alpha),\widehat{W^a})$, where $\A=Y\otimes \R$ with $Y$ the cocharacter lattice of $\cS$. 

In the split case, we simply take $(\Lambda_\alpha)=(\Lambda)_{\alpha\in \Delta}$ and $\widehat{W^a}\in \{W^a,W^a\ltimes (Y\otimes \Lambda)\}$.

\end{Definition}

 A set of the form $M(\alpha,k)$, with $\alpha\in \Phi$ and $k\in \Lambda_\alpha$  is called a \textbf{wall}\index{wall} of $\A$ and a set of the form $D(\alpha,k)$, with $\alpha\in \Phi$ and $k\in \Lambda_\alpha$ is called a \textbf{half-apartment}\index{half-apartment} of $\A$.

\begin{Lemma}\label{l_reflction}
Let $\bw\in W^a_\R$ be such that $\bw$ fixes a wall $M$ of $\A$. Then $\bw\in \{\Id,r_M\}$. 
\end{Lemma}

\begin{proof}
We write $\bw=\tau\circ w$, with $\tau$ a translation of $\A$ and $w\in W^v$. Then $w.M$ is a wall parallel to $M$. Let $M_0$ be the wall parallel to $M$ containing $0$. Then $w.M_0$ is a wall parallel to $M_0$ and containing $0$. Therefore $w.M_0=M_0$. Let $C^v$ be a vectorial chamber dominating a vector panel of $M_0$. Then $w.C^v$ is a chamber adjacent to $C^v$. Consequently, $w.C^v\in \{C^v,r_{M_0}.C^v\}$, where $r_{M_0}$ is the reflection of $W^v$ with respect to $M_0$. As the action of $W^v$ on the set of vectorial chambers is simply transitive, we deduce that $w\in \{\Id,r_{M_0}\}$, which proves the lemma.
\end{proof}

 \subsubsection{Enclosure maps}\label{sss_enclosure}
Let $\underline{\A}=(\A,\mathcal{S},(\Lambda_\alpha),\widehat{W^a})$ be an affine apartment and $\A$ be the underlying affine space.  A subset $X$ of $\A$ is said to be \textbf{enclosed} if there exist $k\in \N$, $\beta_1,\ldots,\beta_k\in \Phi$ and $(\lambda_1,\ldots,\lambda_k)\in\prod_{i=1}^k \Lambda_{\beta_i}$ such that $X=\bigcap_{i=1}^k D(\beta_i,\lambda_i)$.

    Let  $\cV$ be a filter on  $\A$.  Its \textbf{enclosure}\index{enclosure}  is the filter $\cl(\cV)$\index[notation]{c@$\cl$}  on $\A$ generated by the enclosed subsets of $\A$ containing $\cV$.

 For $\alpha \in X\setminus \Phi$, we set $\Lambda_\alpha=\R$.  Let $\cP\subset \Delta$. The \textbf{$\cP$-enclosure} of $\cV$  is the filter $\cl^{\cP}(\cV)$\index[notation]{c@$\cl^{\cP},\cl^{\Delta}$} on $\A$ generated by the sets containing $\cV$  of the form $\bigcap_{\alpha\in \cP} D(\alpha,k_\alpha)$, with $(k_\alpha)\in \prod_{\alpha\in \cP} \Lambda_\alpha$. We often use the $\Delta$-enclosure $\cl^{\Delta}$. Note that when $\Phi$ is infinite,  the enclosure $\cl$ differs from $\cl^{\Phi}$, since $\cl$ only allows finite intersections whereas $\cl^{\Phi}$ allows infinite ones.

   Note that in \cite{gaussent2008kac} and \cite{rousseau2016groupes}, the authors use the $\Delta$-enclosure. This is natural since this is the enclosure that naturally appears when we consider the completed Kac--Moody group. We actually proved in \cite{hebert2022new} that we can replace the $\Delta$-enclosure by the  enclosure in the definition of masures.

\subsubsection{Definitions of  faces, chimneys and related notions}\label{ss_d_faces}
Recall the definition of filters given in Subsection~\ref{ss_filters}.

Let $\underline{\A}=(\A,\mathcal{S},(\Lambda_\alpha),\widehat{W^a})$ be an affine apartment and $\A$ be the underlying affine space.

Let $x\in \A$ and $F^v$ be a vectorial face.  The \textbf{local-face}\index{face (local)} or simply the \textbf{face} $F^\ell(x,F^v)=germ_x(x+F^v)$  is the filter on $\A$ defined as the  intersection of $x+F^v$ with the filter of neighborhoods of $x$ in $\A$. We say that $F$ is \textbf{positive} (or \textbf{negative}) if $F^v$ is. 

Recall that  $[0,0^+[$ (resp. $]0^-,0^+[$, $]0,0^+[$) is the filter on $\R$ generated by the sets of the form $[0,\epsilon]$ (resp. $[-\epsilon,\epsilon]$, $]0,\epsilon[$), with $\epsilon\in \R_{>0}$. A \textbf{segment-germ}\index{sector-germ} $\fs:[0,0^+[\rightarrow \A$ is a map $\fs$ defined on an element $\Omega\in [0,0^+[$ such that for some $a,b\in \A$, we have $\fs(t)=(1-t)a+tb$, for all $t\in \Omega$.

Let $F,F'$ be two faces. We say that $F'$ \textbf{dominates}\index{dominate} $F$ if $F\subset \overline{F'}$. The dimension of a face $F$ is the smallest dimension of an affine space generated by some $S\in F$. Such an affine space is unique and is called its \textbf{support}. A face is said to be \textbf{spherical}\index{spherical} if the direction of its support meets the open Tits cone $\mathring \cT$; then its pointwise stabilizer $W_F$ in $W^a$ is finite.

An \textbf{alcove}\index{alcove} is a face of the form $F^\ell(x,C^v)$ where $x\in \A$ and $C^v$ is a vectorial chamber of $\A$. We set $C_0^+=\germ_0(C^v_f)$\index[notation]{c@$C_0^+$} and $C_0^-=-C_0^+=\germ_0(-C^v_f)$. These are the \textbf{positive} and the \textbf{negative fundamental alcoves}\index{fundamental alcoves} respectively.

A \textbf{panel} is a face whose associated vectorial face is maximal among the faces contained in at least one wall.

\begin{Remark}
Let $x\in \A$ and $F^v$ be a vectorial face. In \cite{rousseau2011masures}, Rousseau defines the face $F(x,F^v)$ as the filter consisting of the subsets containing an intersection  of half-apartments $D(\alpha,\lambda_\alpha)$ or open half-apartments $\mathring{D}(\alpha,\lambda_\alpha)$, with $\lambda
_\alpha\in \Lambda_\alpha\cup\{+\infty\}$ for all $\alpha\in \Phi$  (at most one $\lambda_\alpha\in \Lambda_\alpha$ for each $\alpha\in \Phi$).  In this paper, we will only use the local-faces and not these faces. We make this choice in order to ease the use of galleries of alcoves. Note that by axiom~\ref{a_ma2} (see \ref{sss_axioms}) an apartment contains a local face if and only if it contains the associated face. 
\end{Remark}

Let $F=F^\ell(x,F_1^v)$ be a face ($x\in \A$ and $F_1^v$ is a vectorial face) and $F^v$ be a vectorial face. The \textbf{chimney}\index{chimney}  $\mathfrak{r}(F,F^v)$ is the filter consisting of the sets containing an enclosed set containing $F+F^v$. The vectorial face $F^v$ is the \textbf{direction} of $\fr$. A face $F'$ such that $\fr=\fr(F',F^v)$ is a  \textbf{basis} of the chimney. A chimney is \textbf{splayed}\index{splayed (chimney)} if $F^v$ is spherical. A \textbf{shortening} of a chimney $\mathfrak{r}(F,F^v)$, with $F=F^\ell(x,F_0^v)$ is a chimney of the form $\mathfrak{r} (F^\ell(x+\xi,F_0^v),F^v)$ for some $\xi\in \overline{F^v}$. The \textbf{germ}\index{germ of a chimney} (at infinity) of a chimney $\mathfrak{r}$, denoted $germ_\infty(\fr)$ is the filter of subsets of $\A$ containing a shortening of $\mathfrak{r}$ (this definition of shortening is slightly different from the one of \cite{rousseau2011masures} 1.12 but follows \cite{rousseau2017almost} 3.6) and we obtain the same germs with these two definitions). We often use a capital letter (like $\fR$) to denote the germ of a chimney. 

A \textbf{sector}\index{sector} $\fq$ is a set of the form $x+C^v$, for some $x\in \A$ and $C^v$ a vectorial chamber. The point $x$ is the \textbf{basis} of $\fq$ and $C^v$ is its \textbf{direction}. Its \textbf{germ}\index{germ of a sector} (at infinity) $\germ_\infty(\fq)$, often denoted with a capital letter (like $\fQ$) is the filter on $\A$ generated by the $y+C^v$ such that $y\in \A$. If $\fQ$ is a sector-germ, we write $y+\fQ=y+C^v$\index[notation]{q@$y+\fQ$}, if $C^v$ is the direction of $\fq$. We denote by $\fQ_{\infty}$\index[notation]{q@$\fQ_{+\infty}$, $\fQ_{-\infty}$} and $\fQ_{-\infty}$ the germs at infinity of $C^v_f$ and $-C^v_f$ respectively.

 A \textbf{sector-face}\index{sector-face} $\ff$ is a set of the form $x+F^v$, where $x\in \A$ and $F^v$ is a vectorial face. The point $x$ is the \textbf{basis} of $\ff$ and $F^v$ is its \textbf{direction}. Note that a sector-germ is a particular case of a chimney germ: when the direction of a chimney is a vectorial chamber $C^v$, its germ at infinity coincides with the germ at infinity of any sector of the direction $C^v$.
 
 A face $F=F^\ell(x,F^v)$ or a sector-face $\ff=x+F^v$ is called \textbf{spherical}\index{spherical face of sector-face} if the fixator of $F^v$ in $W^v$ is finite.

\begin{figure}[h]
\centering
\includegraphics[scale=0.25]{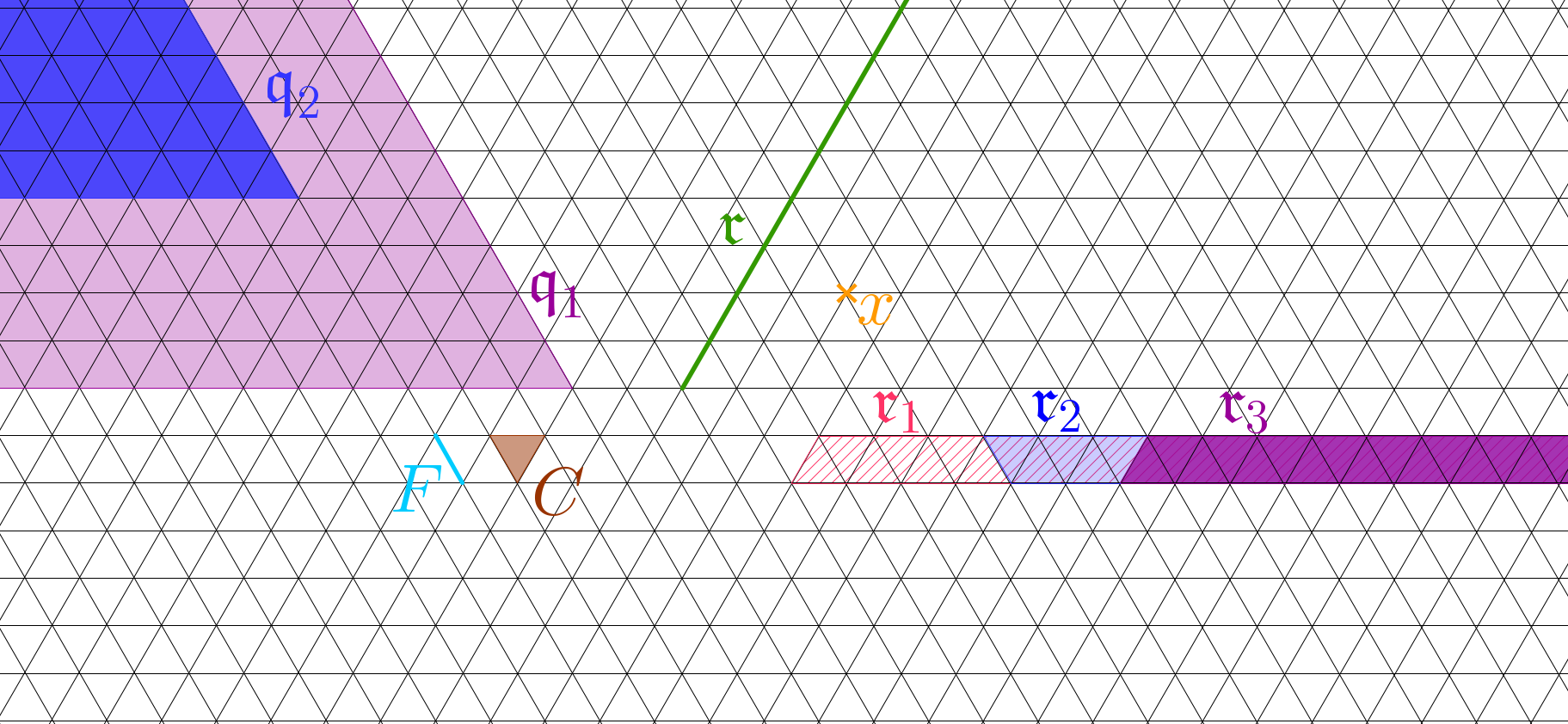}
\caption{Affine apartment of $\mathrm{SL}_3$ when $\Lambda=\Z$. All the sets represented are chimneys: $\fq_1$ and $\fq_2$ have the same direction, which is a vectorial chamber, $\fr_1$, $\fr_2$ and $\fr_3$ have the same direction, which is a ray, the direction of $\fr$ is a ray and $x$, $F$ and $C$ have direction $\{0\}$.  The sectors $\fq_1$ and $\fq_2$ define the same germ at infinity, which is up and left of the picture, $\fr_1,\fr_2$ and $\fr_3$ define the same germ, which is a ``strip at infinity'' at the right of the picture, and the germ of $\fr$ at infinity is a ``ray at infinity'' contained in $\fr$.}\label{f_apt_aff_A2}
\end{figure}

\begin{figure}[h]
\centering
\includegraphics[scale=0.32]{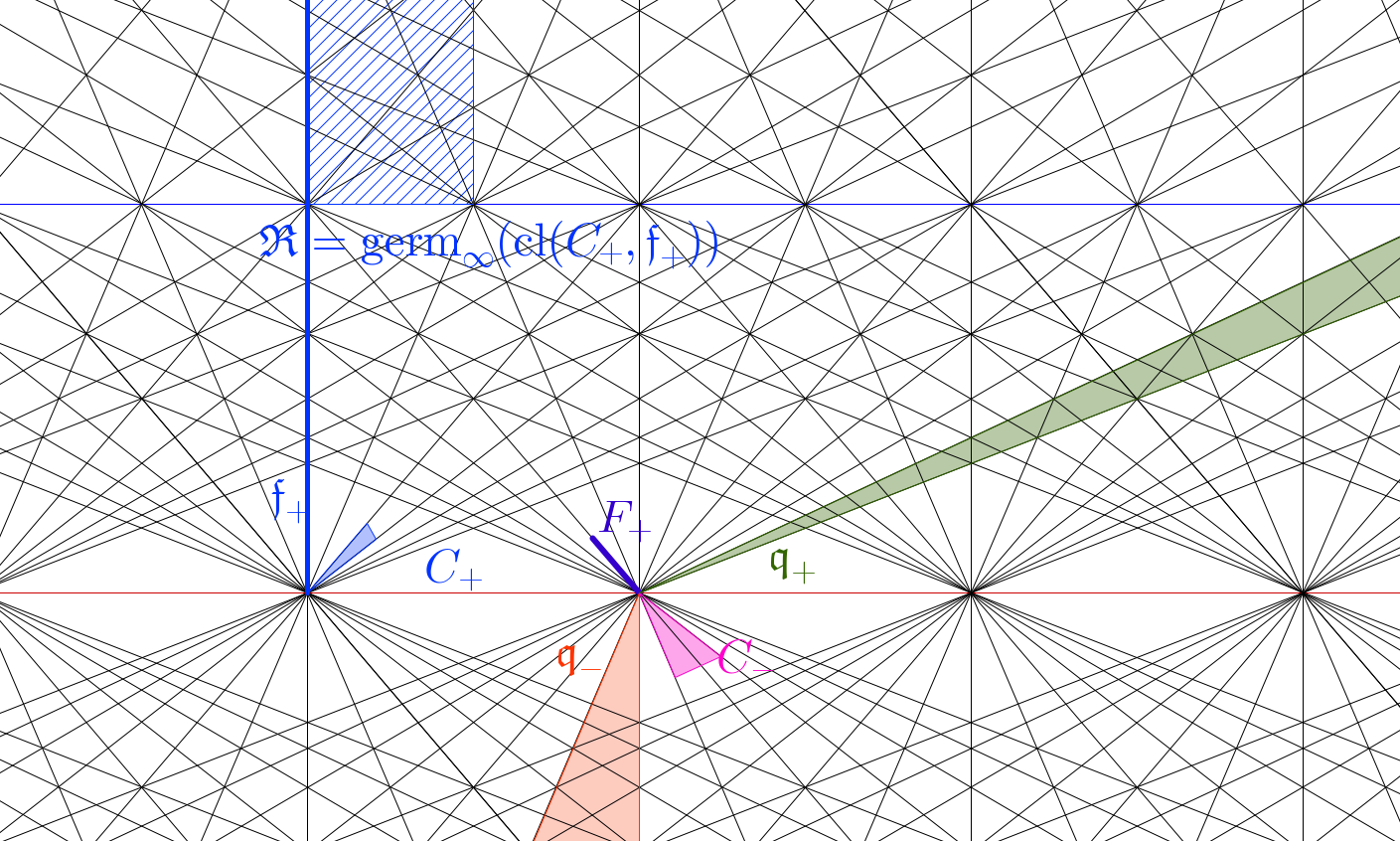}
\caption{Affine apartment of affine  $\mathrm{SL}_2$ when $\Lambda=\Z$.  All the filters represented are chimneys. The set $F^v_+$ is a vectorial panel, $C_+$ is an alcove and the hatched strip represents the germ at infinity of the chimney $\fr=\cl(C_+,\ff_+)$ (actually it is the ``strip at infinity'' in the direction of the hatched strip). The filter $F_+$ is a local panel and $C_-$ is an alcove. The sets  $\fq_+$ and $\fq_-$ are sectors and $\ff_+$ is a sector-panel. For simplicity, all the local faces and sector-faces are based at vertices, but this is arbitrary.}\label{f_apt_aff_SL2}
\end{figure}

\subsubsection{Galleries}\label{sss_gal}

Let $(\fQ,\fQ')$ be a pair of either:\begin{enumerate}
\item alcoves of $\A$ of the same sign based at the same point $x\in \A$,

\item sector-germs of $\A$ (at infinity) of the same sign,

\item vectorial chambers of $\A$ of the same sign.
\end{enumerate}

By symmetry, we assume that $\fQ$ and $\fQ'$ are positive. In the first case, we write $\fQ=germ_x(x+w.C^v_f)$, $\fQ'=\germ_x(x+w'.C^v_f)$, where $w,w'\in W^v$.  We then set $d(\fQ,\fQ')=w^{-1}w'$. In the second case, we write $\fQ=w.\fQ_{+\infty}$ and $\fQ'=w'.\fQ_{+\infty}$, with $w,w'\in W^v$. We set $d(\fQ,\fQ')=w^{-1}w'$. In the third case, we write $\fQ=w.C^v_f$ and $\fQ'=w'.C^v_f$ and we set $d(\fQ,\fQ')=w^{-1}w'$. The ``distance'' $d$ is $W^v$ invariant: we have $d(v.\fQ,v.\fQ')=d(\fQ,\fQ')$ for all $v\in W^v$, in the three situations.

In the three situations, we say that $\fQ$ and $\fQ'$ are \textbf{adjacent} if $d(\fQ,\fQ')\in \{1\}\cup \{r_i\mid i\in I\}$.  If $\fQ$ and $\fQ'$ are alcoves, this means that they dominate a common panel. If they are vectorial chambers, it means that they dominate a common vectorial panel. If they are sector-germs, it means that for all $x\in \A$, $x+\fQ$ and $x+\fQ'$ dominate a common sector-face panel, i.e a sector face of codimension $1$. 

A \textbf{gallery} $\Gamma=(\fQ_1,\ldots,\fQ_k)$ is a finite sequence of either alcoves, sector-germs or vectorial chambers (all the $\fQ_i$ have to be of the same type) such that for each $i\in \llbracket 1,k-1\rrbracket$, $\fQ_i$ and $\fQ_{i+1}$ are adjacent or equal. The \textbf{length} of $\Gamma$ is $k-1$. 

\begin{Lemma}\label{l_gal_apt}
Let $(\fQ, \fQ')$ be as in the beginning of this subsubsection. Then there exists a gallery between $\fQ$ and $\fQ'$. Moreover, the minimal possible length of such a gallery is $\ell(d(\fQ,\fQ'))$.
\end{Lemma}

\begin{proof}
As the three situations are similar, we treat the case of a pair of  vectorial chambers. Write $\fQ=w.C^v_f$ and $\fQ'=w'.C^v_f$, with $w,w'\in W^v$. Let $v=w^{-1}w'$. Let $k=\ell(v)$ and $r_{i_1}\ldots r_{i_k}$ be a reduced expression of $v$ (i.e $i_1,\ldots,i_k\in I$ and $v=r_{i_1}\ldots r_{i_k}=v$). Set $\fQ_0=\fQ=w.C^v_f$, $\fQ_1=wr_{i_1} C^v_f$, $\ldots$, $\fQ_k=wr_{i_1}\ldots r_{i_k}.C^v_f=w'.C^v_f$. Then $\Gamma$ is a gallery from $\fQ$ to $\fQ'$. 
\end{proof}

\subsubsection{Apartment of type $\underline{\A}$}

Let $\underline{\A}=(\A,\cS,(\Lambda_\alpha),\widehat{W^a})$ be an affine apartment as in Definition~\ref{d_apt}.  An \textbf{apartment of type }$\underline{\A}$ is a set $A$ equipped with a nonempty set $\mathrm{Isom}(\A,A)$ of bijections (called \textbf{isomorphisms}) such that if $f_0\in \mathrm{Isom}(\A,A)$ then $f\in \mathrm{Isom}(\A,A)$ if and only if there exists $w\in \widehat{W^a}$ satisfying $f=f_0\circ w$.  Let $A$ and $A'$ be two apartments of type $\underline{\A}$. An (apartment) isomorphism between $A$ and $A'$ is an element of the form $g\circ f^{-1}$, where $f\in \mathrm{Isom}(\A,A)$ and $g\in \mathrm{Isom}(\A,A')$. We extend all the notions that are preserved by $\widehat{W^a}$ to each apartment. Thus sectors, enclosures, faces and chimneys are well-defined in any apartment of type $\underline{\A}$ (they actually only depend on the set of walls of $\underline{\A}$, which is stabilized by $\widehat{W^a}$).

\subsection{Weak masures and masures}\label{ss_wma}

In this subsection, we introduce several axioms of masures. We then give the definition of masures.

\subsubsection{Axioms}\label{sss_axioms}

 We consider a space $\cI$ endowed with a covering $\mathcal{A}$ of subsets called \textbf{apartments} such that:  \begin{enumerate}[label=\blue{(MA1)}]
                                                                                                                    \item\label{a_ma1} Any $A\in \mathcal{A}$ is equipped with the  structure of an apartment of type $\underline{\A}$.\axiom{ma1@\ref{a_ma1}}

\end{enumerate}

    The elements of $\cA$ are called \textbf{apartments}. We can define the notion of walls, faces, enclosure $\cl_A$, $\Delta$-enclosure $\cl^{\Delta}_A$, segments $[x,y]_A$ between two points $x,y$ of $A$, ... in any apartment $A$ of type $A$. We also consider the following axioms:

\begin{enumerate}[label=\blue{(MA2)}]
\item\label{a_ma2} Let $A$ and $B$ be two apartments. Then $A\cap B$ is enclosed, i.e there exist $k\in \N$ and  half-apartments $D_1,\ldots,D_k$ of $A$ such that $A\cap B=\bigcap_{i=1}^k D_i$ and there exists an apartment isomorphism $\phi:A\rightarrow B$ fixing $A\cap B$.\axiom{ma2@\ref{a_ma2}}

\end{enumerate}

\begin{enumerate}[label=\blue{(MA3)}]
\item\label{a_ma3} Let $\fQ$ be a splayed chimney germ and $\fR$ be a local face or a chimney germ. Then there exists an apartment containing $\fQ$ and $\fR$.\axiom{MA3@\ref{a_ma3}}
\end{enumerate}

\begin{enumerate}[label=\blue{(wMA3)}]
\item\label{a_wma3} (``weak MA3'') : Let $\fQ, \fQ'$  be two sector-germs at infinity. Then  there exists an apartment containing them.\axiom{MA3@\ref{a_wma3}}
\end{enumerate}

\begin{enumerate}[label=\blue{(OCO)}]
 \item\label{a_oco} (``ordered convexity''): Let $\cV$ be a filter on an apartment $A$ such that there exists a  sector $\fq$ based at some $x\in A$ such that $\{x\}\subset \overline{\cV}\subset \overline{\fq}$ (i.e $\overline{\fq}\in \cV$ and $x\in \overline{\Omega}$,  for all $\Omega\in \cV$). Then for every apartment $A'$ containing $\cV$, we have $\cl^{\Delta}_A(\cV)\subset A'$ and there exists an apartment isomorphism from $A$ to $A'$ fixing $\cl^{\Delta}_A(\cV)$.\axiom{OCO@\ref{a_oco}}
\end{enumerate}

                                                                                                                    In an space satisfying \ref{a_ma1} and \ref{a_wma3}, we say that two sector-germs at infinity are \textbf{adjacent} if there exists an apartment containing them in which they are adjacent sector-germs. We now introduce two last axioms.

                                                                                                                    \begin{enumerate}[label=\blue{(SC)}]
\item\label{a_sc} (``sundial configuration'', see Figure~\ref{f_SC}) Let $\fQ$ and $\fQ'$ be two adjacent sector-germs at infinity and $A$ be an apartment containing $\fQ$. Then there exist two opposite half-apartments $D_1$ and $D_2$ of $A$ (i.e $D_1\cap D_2$ is the wall of both $D_1$ and $D_2$) such that for both $i\in \{1,2\}$, there exists an apartment containing $D_i\cup \fQ'$.\axiom{SC@\ref{a_sc}}
\end{enumerate}

\begin{enumerate}[label=\blue{(EC)}]
\item\label{a_ec} (``exchange condition'', see Figure~\ref{f_EC}) Let $A$ and $B$ be two apartments such that $A\cap B$ is a half-apartment. Let $M$ be the wall of $A\cap B$. Then $(A\setminus B)\cup (B\setminus A)\cup M$ is an apartment.\axiom{EC@\ref{a_ec}}
\end{enumerate}

Note that axiom~\ref{a_oco} is a strong version of axiom (MAO) of \cite{rousseau2011masures}. It is clearly implied by~\ref{a_ma2}. 

Axiom \ref{a_wma3} is a weak version of axiom~\ref{a_ma3}, since sector-germs at infinity are chimney germs. These two axioms are  building theoretic translations of decompositions of $G$, when $\I$ is associated to a Kac--Moody group $G$. Axiom \ref{a_wma3} is the translation of  the Bruhat and the Birkhoff decompositions in $G$ (see \eqref{e_Bruhat} and \eqref{e_Birkhoff}) and~\ref{a_ma3} is a consequence of a decomposition of $G$ proved in Subsection~\ref{ss_MA3}.

\begin{Definition}\label{d_w_mas}
Let $\cS$ be a free Kac--Moody datum. Let $\underline{\A}=(\cS,(\Lambda_\alpha)_{\alpha\in \Phi},\widehat{W^a})$ be an apartment in the sense of Definition~\ref{d_apt}. \begin{enumerate}

\item A \textbf{weak masure} of type $\underline{\A}$ is a set $\mathcal{I}$ endowed with a covering $\mathcal{A}$ of subsets called  apartments  satisfying~\ref{a_ma1}, \ref{a_wma3},~\ref{a_oco} and \ref{a_sc}.  

\item  A \textbf{masure} of type $\underline{\A}$ is a set $\mathcal{I}$ endowed with a covering $\mathcal{A}$ of subsets called  apartments satisfying \ref{a_ma1},~\ref{a_ma2} and~\ref{a_ma3}. 

\end{enumerate}
\end{Definition}

It is clear that~\ref{a_ma3} implies \ref{a_wma3} and that~\ref{a_ma2} implies~\ref{a_oco}. By  \cite[Lemma 8]{hebert2021distances}, a masure (associated with a free and cofree Kac--Moody datum $\cS$) satisfies \ref{a_sc}. Therefore (at least when $\cS$ is cofree), a masure is a weak masure.

The main theorem of this chapter is Theoreme~\ref{t_MA2}, which asserts that a weak masure (with an additional assumption on the Kac--Moody datum) satisfies~\ref{a_ma2}. Along the proof of Theorem~\ref{t_MA2}, we obtain properties of masures which are useful for the applications of the theory of masures. From Section~\ref{s_rtrct_infty} to Section~\ref{s_spltng_apt}, we study basic properties of (weak) masures. In Section~\ref{s_H_paths}, we introduce and study Hecke paths, which are a fundamental tool in the study of Kac--Moody groups, for example  they enable  to define Hecke algebras of Kac--Moody groups (\cite{gaussent2014spherical}  and \cite{bardy2016iwahori}). Sections~\ref{s_cvx_intrsctn} and \ref{s_encls_int} might be less useful for the reader ready to admit that $\I$ satisfies~\ref{a_ma2}.

\begin{Remark}\label{r_MA1}
Let $(\Lambda_\alpha)\in \sL$ and $\widehat{W^a_1}, \widehat{W^a_2}$ be two groups such that $W^a((\Lambda_\alpha))\subset \widehat{W^a_1}\subset \widehat{W^a_2}\subset \Aut^{\vw}(\A,(\Lambda_\alpha))$. Set $\underline{\A_i}=(\A,\cS,(\Lambda_\alpha),\widehat{W^a_i})$, for $i\in \{1,2\}$. Then if $\I$ is a masure of type $\underline{\A_1}$, it can naturally be regarded as a masure of type $\underline{\A_2}$. Indeed, it suffices to set $\mathrm{Isom}_{\underline{\A_2}}(\A,A)=\{f\circ \bw\mid f\in \mathrm{Isom}_{\underline{\A_1}}(\A,A), \bw\in \widehat{W^a_2}\}$, for any apartment $A$ of $\I$, where $\mathrm{Isom}_{\underline{\A_i}}(\A,A)$ is the set of isomorphisms corresponding to the structure of type $\underline{\A_i}$, for $i\in \{1,2\}$. 
\end{Remark}

\subsection{Apartment isomorphisms fixing filters}

We equip $\A$  with its topology of finite dimensional affine space. As the elements of $W^a$ are homeomorphisms for this topology, this equips any apartment $A$ of type $\underline{\A}$ with a topology. 

\begin{Remark}\label{r_iso_fix_open}
Let $A,B$ be two apartments and $f,f':A\rightarrow B$ be apartment isomorphisms. Assume that $f$ and $f'$ coincide on a subset of $A\cap B$ with non-empty interior. Then $f=f'$. Indeed, identify $A$ and $\A$. Then  $f'^{-1}\circ f:\A\rightarrow \A$  is an apartment isomorphism and thus it belongs to $W^a$. Therefore $f'^{-1}\circ f$ is an affine automorphism fixing an open subset of $\A$ and thus $f'^{-1}\circ f$ is the identity.
\end{Remark}

\begin{Lemma}\label{l_iso_sect}
Let $(\I,\cA)$ be a space satisfying~\ref{a_ma1} and ~\ref{a_oco}. Let $A,B$ be two apartments containing a sector $\fq$. Then there exists a unique apartment isomorphism $f:A\rightarrow B$ fixing $A\cap B$. Moreover, if $y\in A\cap B$, then the translate of $\fq$ at $y$ (in $A$) is contained in $A\cap B$. 
\end{Lemma}

\begin{proof}
We identify $A$ and $\A$. Let $\fq$ be a sector contained in $\A\cap B$.  Write $\fq=x+C^v$, with $C^v$ a vectorial chamber of $\A$. By~\ref{a_oco}, there exists an apartment isomorphism $f:\A\rightarrow B$ fixing $\fq$.  Let $y\in \A\cap B$. Then $(y+C^v)\cap (x+C^v)$ is a sector of $\A$ and we can write it $z+C^v$, for some $z\in \A$. Let $\Omega=y\cup (z+C^v)\subset y+C^v$. Then by~\ref{a_oco}, $\overline{y+C^v}\subset \cl^{\Delta}(\Omega)\subset \A\cap B$ and there exists an apartment isomorphism $f_y:\A\rightarrow B$ fixing $y+C^v$. But $f_y$ and $f$ coincide on $z+C^v$ and thus by Remark~\ref{r_iso_fix_open}, $f=f_y$ and $f(y)=f_y(y)=y$.  Therefore $f$ fixes $\A\cap B$, which proves the lemma. 
\end{proof}

\section{Retractions centered at infinity and inseparable filters}\label{s_rtrct_infty}

From now on in this chapter, $\I$ is a weak masure: it satisfies~\ref{a_ma1}, \ref{a_wma3},~\ref{a_oco} and \ref{a_sc}.

In this section, we prove that an inseparable filter (see Definition~\ref{d_insep_filter} below) and a sector-germ at infinity are always contained in a common apartment. This enables to define retractions centered at infinity, since points are particular cases of inseparable filters.

Recall that if $\alpha\in \Phi$ and $x\in \R$, $D(\alpha,x)=\{a\in \A\mid \alpha(a)+x\geq 0\}$. 
\begin{Definition}\label{d_insep_filter}
Let $\cV$ be a filter on $\A$. We say that $\cV$ is \textbf{inseparable}\index{inseparable} if for every $\alpha\in \Phi$ and $x\in \Lambda_\alpha$, we have $\cV\subset D(\alpha,x)$ or $\cV\subset D(-\alpha,-x)$. 
\end{Definition}

\begin{Example}
\begin{enumerate}
\item A face or a sector-face-germ at infinity is inseparable.

\item The germ at infinity of a ray is inseparable.

\item A chimney germ is inseparable.

\item Let $\cV$ be an inseparable filter. Then $\cl(\cV)$ is inseparable. 
\end{enumerate}

\end{Example}

Let $\fQ$, $\fQ'$ be two sector-germs of $\A$ of the same sign $\epsilon$. We say that they are \textbf{adjacent} if for some apartment $A$ containing $\fQ$ and $\fQ'$ (which exists by \ref{a_wma3}) and some apartment isomorphism $f:A\rightarrow \A$, $f(\fQ)$ and $f(\fQ')$ are adjacent. This does not depend on the choice of $f$ and by Lemma~\ref{l_iso_sect}, this does not depend on the choice of $A$. This enables us to define the notion of galleries between two sector-germs of the same sign, similarly as in \eqref{sss_gal}.

\begin{Proposition}\label{p_mas_th_iwa}
Let $\fQ$ be a sector-germ at infinity and $\cV$ be an inseparable filter. Then there exists an apartment containing $\fQ$ and $\cV$.  
\end{Proposition}

\begin{proof}
Let $A$ be an apartment containing $\cV$. Let $\fR$ be a sector-germ of $A$ of the same sign as $\fQ$. Let $B$ be an apartment containing $\fQ$ and $\fR$.  Let $\Gamma=(\fQ_1,\ldots,\fQ_n)$ be a gallery  of sector-germs from $\fR=\fQ_1$ to $\fQ=\fQ_n$ (in $B$).  By the sundial configuration \ref{a_sc}, we can write $A=D_1\cup D_2$, with $D_i\cup \fQ_2$ contained in some apartment $A_i$, for both $i\in \{1,2\}$. As $\cV$ is inseparable, either $D_1$ or $D_2$ contains $\cV$. Up to exchanging $D_1$ and $D_2$, we may assume $D_1\supset \cV$. Then we can replace $\Gamma$ by $(\fQ_2,\ldots,\fQ_n)$ and by induction, we deduce the result. 
\end{proof}

\begin{Corollary}
Let $F$  be a face (resp. $\fs:[0,0^+[\rightarrow \I$ be a segment-germ) and $\fQ$ be a sector-germ. Then there exists an apartment containing $F$ and $\fQ$ (resp. $\fs$ and $\fQ$). 
\end{Corollary}

\begin{Corollary}\label{c_fcs_pnt_frnd}
Let $F_1,F_2$ be two faces of $\I$ based at the same point. Then there exists an apartment containing $F_1$ and $F_2$. 
\end{Corollary}

\begin{proof}
Let $C_1$ be an alcove dominating $F_1$. By~\ref{a_oco}, any apartment containing $C_1$ contains $F_1$. Therefore we can assume that $F_1=C_1$. Let $A$ be an apartment containing $C_1$ and $\fq$ be the sector based at the origin $x$ of $C_1$ and containing $C_1$. Let $\fQ$ be the sector-germ at infinity of $\fq$. By Proposition~\ref{p_mas_th_iwa}, there exists an apartment $A'$ containing $C_2$ and $\fQ$. Then $A$ contains $x$ and $\fQ$ and thus by~\ref{a_oco}, $A\cap A'$ contains $\overline{\fq}$ and hence  $C_1$, which proves the corollary.
\end{proof}

\begin{Definition}
Let $x\in \I$ and $\fQ$ be a sector-germ at infinity. Then by Proposition~\ref{p_mas_th_iwa} there exists an apartment $A$ containing $x$ and $\fQ$. We denote by $x+\fQ$ the sector of $A$ based at $x$  and whose germ at infinity is $\fQ$. This sector is well-defined, independently of the choice of $A$, by~\ref{a_oco}. 

\end{Definition}

Let $\fQ$ be a sector-germ at infinity and $A$,  $B$ be two apartments containing $\fQ$, then we have \begin{equation}\label{e_apat_sctr}
y+\fQ\subset A\cap B,  \forall y\in A\cap B,
\end{equation} by Lemma~\ref{l_iso_sect}.

\begin{Definition}\label{d_rtrct_infty}
Let $A$ be an apartment and $\fQ$ be a sector-germ. We define the \textbf{retraction $\rho_{A,\fQ}:\I\rightarrow A$ onto $A$ and centered at $\fQ$} as follows. Let $x\in \I$. By Proposition~\ref{p_mas_th_iwa}, there exists an apartment $A_x$ containing $x$ and $\fQ$. By Lemma~\ref{l_iso_sect}, there exists an apartment isomorphism $f:A_x\rightarrow A$ fixing $A\cap A_x$. We set $\rho_{A,\fQ }(x)=f(x)$\index[notation]{r@$\rho_{A,\fQ}$}. This is well-defined, independently on the choice of $A_x$. Indeed, if $A'$ is an other apartment containing $x$ and $\fQ$, then there exists an apartment isomorphism $f':A'\rightarrow A$ fixing $A'\cap A$, by Lemma~\ref{l_iso_sect}. Let $h:A'\rightarrow A_x$ be the apartment isomorphism fixing $A\cap A'$. Then $f'^{-1}\circ f\circ h:A'\rightarrow A'$ is an apartment isomorphism fixing $\fQ$ pointwise. As it is an affine automorphism of $A'$, it is the identity, which proves that $f'(x)=f(x)$. 

We denote by $\rho_{+\infty}:\I\rightarrow \A$\index[notation]{r@$\rho_{+\infty}, \rho_{-\infty}$} and $\rho_{-\infty}:\I\rightarrow \A$ the retractions centered at $\fQ_{+\infty}$ and $\fQ_{-\infty}$ respectively. 
\end{Definition}

\begin{Remark}
When $\I$ is associated with a Kac--Moody group $G$, Proposition~\ref{p_mas_th_iwa} corresponds to an Iwasawa decomposition of $G$ (using Lemma~\ref{l_frndly_prs}). However, in order to prove that $\I$ is a weak masure, we will use the  Iwasawa decomposition of $G$ (Proposition~\ref{p_Iwa}): our proof that ~\ref{a_oco} is satisfied uses Theorem~\ref{t_G_x}, which uses the Iwasawa decomposition. Therefore we do not get a fully building theoretic proof of the Iwasawa decomposition of $G$.
\end{Remark}

\section{Exchange condition}

In this section, we prove that a weak masure satisfies the exchange condition.

Recall that $\I$ is a weak masure. Recall that a half-apartment is the image by an apartment isomorphism of  a set of the form $D(\alpha,k)\subset \A$, with $\alpha\in \Phi$ and $k\in \Lambda_\alpha$. If $D_1$ and  $D_2$ are two half-apartments of an apartment $A$, we say that they are \textbf{opposite}\index{opposite half-apartments} if $D_1\cap D_2$ is a wall (and thus $D_1\cup D_2=A$). 

\begin{Lemma}\label{l_half_apt}
Let $A$ and $B$ be two distinct apartments such that $A\cap B$ contains a half-apartment. Then $A\cap B$ is a  half-apartment.
\end{Lemma}

\begin{proof}
We can assume that $A=\A$. Let $D$ be a half-apartment contained in $\A\cap B$ and $M$ be its boundary. Let $\alpha\in \Phi$ be such that $D=\alpha^{-1}(\R_{\geq \lambda})$, for some $\lambda\in \Lambda_\alpha$. Let $M=\alpha^{-1}(\{\lambda\})$. Let $\fq$ be a sector contained in $D$ and whose boundary contains a sector-panel $P$ of $M$. Let $\tilde{P}$ be the sector-panel of $M$ opposite to $P$. Let $\tilde{\fq}$ be the sector of $D$ dominating $\tilde{P}$. Let $\fQ$ and $\tilde{\fQ}$ be the germs at infinity of $\fq$ and $\tilde{\fq}$ respectively. Then if $x\in \A\cap B$, we have $\overline{(x+\fQ)}\cup \overline{(x+\tilde{\fQ})}\subset \A\cap B$. Therefore $\bigcup_{y\in \overline{x+\tilde{\fQ}}}\overline{y+\fQ}=\alpha^{-1}(\R_{\geq \alpha(x)})\subset \A\cap B$, by \eqref{e_apat_sctr}. Consequently, as $B\neq \A$ by assumption, we deduce the existence of $a\in \R$ such that $\A\cap B\in \{\alpha^{-1}(\R_{>a}),\alpha^{-1}(\R_{\geq a})\}$. Choose $x\in \A$ such that $\alpha(x)=a$. Then $x+\fQ\subset \A\cap B$ and hence by~\ref{a_oco}, $x\in \overline{x+\fQ}\subset \A\cap B$ and $\cl^{\Delta}(x+\fQ)\subset \A\cap B$. Therefore $a\in \Lambda_\alpha$ and $\A\cap B=\alpha^{-1}(\R_{\geq a})$, which proves the lemma.  
\end{proof}

\begin{figure}[h]
\centering
\includegraphics[scale=0.3]{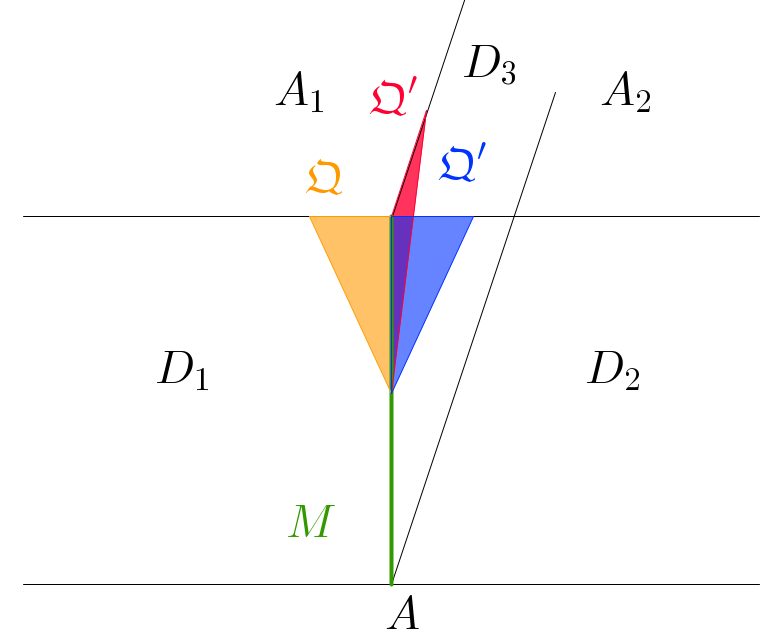}
\caption{Sundial configuration. We use the notation of Lemma~\ref{l_Y_prop}}\label{f_SC}
\end{figure}

\begin{Lemma}\label{l_Y_prop}(see Figure~\ref{f_SC})
Let $A$ be an apartment. Let $\fQ$, $\fQ'$ be two adjacent sector-germs such that $\fQ\subset A$ and $\fQ'\not\subset A$. Let $D_1,D_2$ be two opposite half-apartments of $A$ such that $D_i\cup \fQ'$ is contained in an apartment $A_i$ for both $i\in \{1,2\}$. Let $M=D_1\cap D_2$. Up to renumbering, we may assume that $\fQ\subset D_1$.   Then:\begin{enumerate}

\item There exists a unique sector-germ $\fQ''$ of $A$, not contained in $D_1$ and adjacent to $\fQ$. This sector-germ is adjacent to $\fQ'$. 

\item  $A_1\cap A_2$ is a half-apartment of $\I$,

\item $A_i\cap A=D_i$ for both $i\in \{1,2\}$,

\item  $A\cap A_1\cap A_2=M$,

\item Let $D_3=A_1\cap A_2$. Then $A_i=D_i\cup D_3$, for both $i\in \{1,2\}$. In particular, if $i,j\in \{1,2\}$ is such that $\{i,j\}=\{1,2\}$, we have  $A_j=(A_i\setminus A)\cup (A\setminus A_i)\cup M$. 
\end{enumerate} 
\end{Lemma}

\begin{proof}
The apartment $A_1$ contain $\fQ$ and $\fQ'$. Let $x\in \In{(D_1)}$. Set  $\fq=x+\fQ$ and $\fq'=x+\fQ'$. Then $\fp=\overline{\fq}\cap \overline{\fq'}$ is a sector-panel of $A_1$, contained in $A\cap A_1$. Let $\fQ''$ be the sector-germ of $A$ such that $x+\fQ''$ dominates $P$ and $\fQ''\neq \fQ$. If $\fQ''$ was contained in $D_1$, then $\fQ,\fQ'$ and $\fQ''$ would be contained in $A_1$ and would dominate $P$: this is impossible. Therefore $\fQ''$ is not contained in $D_1$. Moreover, $\fQ''$ is adjacent to $\fQ$. Up to changing the choice of $\A$ and of $C^v_f$, we can assume that $A=\A$ and that $\fQ=\fQ_{\pm \infty}$. By symmetry, we can assume that $\fQ=\fQ_{+\infty}$. Then there exists $i\in I$ such that $\fQ''=r_i.\fQ_{+\infty}$. As $\fQ''$ is not contained in $D_1$, we have $D_1=D(\alpha_i,\lambda)$, for some $\lambda\in \R$.  Then the other sector-germs adjacent to $\fQ$ are of the form $r_j.\fQ_{+\infty}$, for some $j\in I\setminus \{i\}$ and thus they are all contained in $D_1$. This proves the uniqueness of $\fQ''$.

The apartments $A_1$ and $A_2$ contain $\fQ'$ and $M$. Hence by Lemma~\ref{l_iso_sect}, $A_1\cap A_2$ contains $\bigcup_{x\in M} x+\fQ'$, which is  a half-apartment of both $A_1$ and $A_2$. Therefore $A_1\cap A_2$ is a half-apartment, by Lemma~\ref{l_half_apt}.

Let $i\in  \{1,2\}$. By assumption $D_i\subset A_i\cap A$. Let $j\in \{1,2\}\setminus \{i\}$. Let $x\in \In(D_j)$. Assume by contradiction that $x\in A\cap A_i$. Then $A\cap A_i$ is a half-apartment containing $x$ and $D_i$. Let $y\in M$. Then $A_i$ contains $\germ_y(y+\fQ), \germ_y(y+\fQ')$ and $\germ_y(y+\fQ'')$. Moreover, these three alcoves are distinct: indeed, $y\in A_1\cap A_2\cap A$. In $A$, $\germ_y(y+\fQ)$ and $\germ_y(y+\fQ'')$ are distinct, since $A$ contains $\fQ$ and $\fQ''$. Similarly, as $y\in A_1$,  $\germ_y(y+\fQ)\neq \germ_y(y+\fQ')$ and as $y\in A_2$, $\germ_y(y+\fQ')\neq \germ_y(\fQ'')$.  Therefore the  three alcoves  $\germ_y(y+\fQ)$,  $\germ_y(y+\fQ')$ and $\germ_y(\fQ'')$ dominate a single panel of $M$: a contradiction. Therefore $x\notin A\cap A_i$ and $A\cap A_i=D_i$. This proves that $A\cap A_1\cap A_2=M$. 

Let $i\in \{1,2\}$. Then $D_i$ and $D_3$ are two half-apartments of $A_i$. Moreover, $D_i\cap D_3\subset A\cap A_1\cap A_2=M$ and hence $D_i\cap D_3=M$. Let $\psi:A_i\rightarrow \A$ be an apartment isomorphism. Then $\psi(D_i)$ and $\psi(D_3)$ are two half-apartments whose intersection is a wall. Therefore $\psi(D_i)\cup \psi(D_3)=\A$ and hence $D_i\cup D_3=A_i$, which completes the proof of the lemma.
\end{proof}

The weak masure $\I$ satisfies the \textbf{exchange condition}\index{e@(EC), exchange condition} \ref{a_ec} of  \cite{bennett2014axiomatic}:
\begin{Proposition}\label{p_EC_MT} (see Figure~\ref{f_SC})
Let $A$ and $B$ be two apartments such that $A\cap B$ is a half-apartment. Let $M$ be the wall of $A\cap B$. Then $(A\setminus B)\cup (B\setminus A)\cup M$ is an apartment.
\end{Proposition}

\begin{proof}
Let $\fR$, $\fR'$ be two  sector-germs of $B$ of the same sign such  that $\fR\subset A\cap B$ and $\fR'\subset B\setminus A$. Let $\Gamma=(\fQ_0,\ldots,\fQ_k)$ be a gallery of sector-germs of $B$ from $\fR$ to $\fR'$. Let $i\in \llbracket 0,k-1\rrbracket$ be minimal such that $\fQ_i\subset A\cap B$ and $\fQ_{i+1}\not \subset A\cap B$. We set $\fQ=\fQ_i$ and $\fQ'=\fQ_{i+1}$. Then by the sundial configuration \ref{a_sc}, there exist two half-apartments $D_1$ and $D_2$ of $A$ such that for both $j\in \{1,2\}$, there exists an apartment $A_j$ containing $D_j\cup \fQ'$. Let $j\in \{1,2\}$ be such that $D_j\supset \fQ$. Let $\fQ''$ be the unique sector-germ of $A$, not contained in $D_j$ and adjacent to $\fQ$.  Then $A\cap A_j$ is a half-apartment delimited containing $\fQ$ but not $\fQ''$ and hence $A\cap A_j$ and $B\cap A$ are two parallel half-apartments of $A$. Therefore $B\cap A_j$ contains a half-apartment of $A$. Moreover $B\cap A_j$ contains $\fQ'$ and hence $B=A_j$, by Lemma~\ref{l_half_apt}. Therefore the proposition follows from Lemma~\ref{l_Y_prop}. 
\end{proof}

\section{Local sundial configuration}

In this section, we prove that the weak masure $\I$ satisfies the local sundial configuration. We also study the existence of apartments containing pairs of local faces.

\begin{figure}[h]
\centering
\includegraphics[scale=0.35]{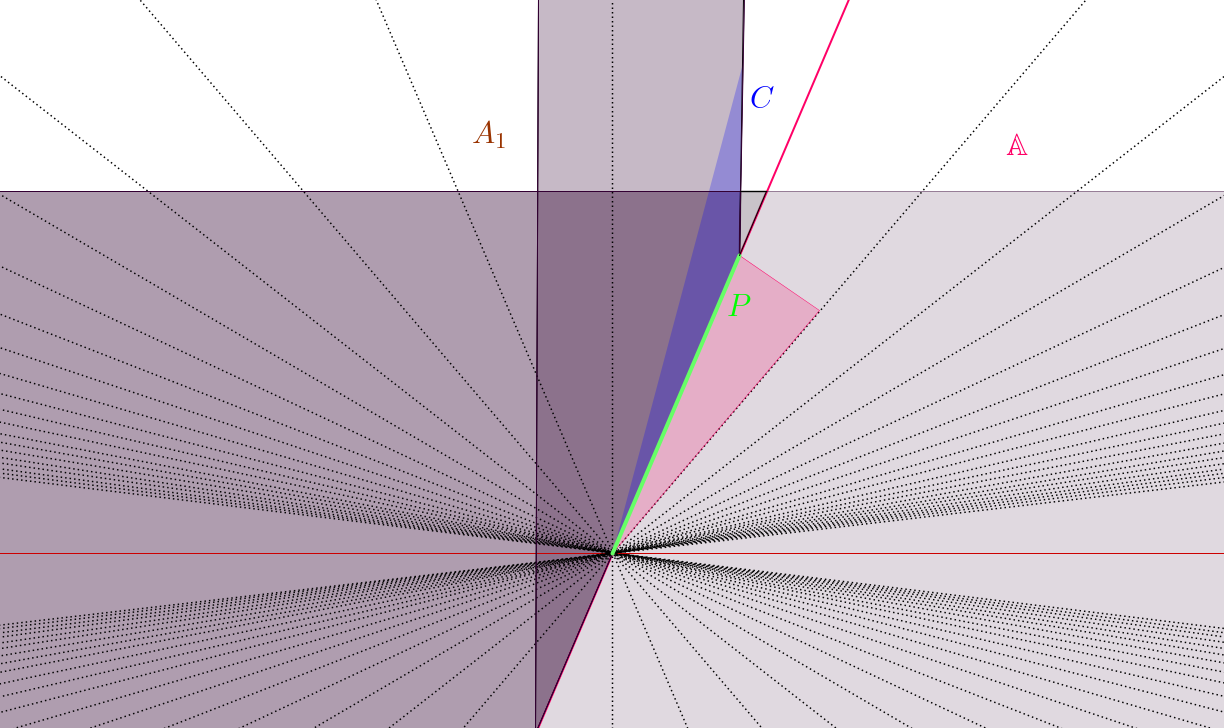}
\caption{Local sundial configuration when the masure is associated with affine $\mathrm{SL}_2$. The panel $P$ is dominated by $C$. The red alcove of $\A$ and $C$ dominate the panel $P$. The apartment $A_1$ contains $C$ and is such that $\A\cap A_1$ is a half-apartment.}\label{f_l_SC}
\end{figure}

\begin{Proposition}\label{p_loc_SC}(\textbf{local sundial configuration}, see Figure~\ref{f_l_SC}) Let $D$ be a  half-apartment of $\I$, $P$ be a (local) panel contained in the wall of $D$ and $C$ be an alcove of $\I$ dominating $P$. Then there exists an apartment $A_1$ containing $D\cup C$. 
\end{Proposition}

\begin{proof}
Let $A$ be an apartment containing $D$. If $C\subset A$, it suffices to take $A_1=A$. We now assume $C\nsubseteq A$.  Let $M$ be the wall of $D$. Let $x$ be the base-point of $P$. Let $\fq$ be the sector of $D$ based at $x$ and  such that $\germ_x(\fq)$ is the alcove of $D$ dominating $P$. Let $\fQ=\germ_\infty(\fq)$.  By Proposition~\ref{p_mas_th_iwa}, there exists an apartment $B$ containing $\fQ$ and $C$. Let $\fq'$ be the sector of $B$ such that $\germ_x(\fq')=C$. The apartment $B$ contains $x$ and $\fQ$ and thus it contains $\fq$, by~\ref{a_oco}. Moreover $\fq'$ and $\fq$ dominate $P$ and thus $\fQ':=\germ_{\infty}(\fq')$ and $\fQ$ are adjacent. 

By \ref{a_sc}, we can write $A=D_1\cup D_2$, where $D_1,D_2$ are opposite half-apartments such that $D_i\cup \fQ'$ is contained in an apartment $A_i$ for both $i\in \{1,2\}$. Up to renumbering, we may assume that $P\subset A_1$. Then $P\subset A\cap A_1$. Moreover $x\in \overline{P}\subset A_1$ and $A_1$ contains $\fQ'$. Therefore $A_1$ contains $\fq'$ and in particular, $A_1$ contains $C$. Using Lemma~\ref{l_half_apt}, we deduce that  $A\cap A_1$ is a half-apartment of $A$ delimited by $P$. As $A\cap A_1$ contains $\fQ$, it contains $D$ and thus $A\cap A_1=D$. Therefore $A_1$ satisfies the condition of the proposition. 
\end{proof}

\begin{Remark}
 Proposition~\ref{p_loc_SC} is \cite[Proposition 2.9 1)]{rousseau2011masures}. However our proofs differ since it is based on the existence of retractions and on the sundial configuration. Rousseau uses the existence of apartments containing chimneys germs. We do not know whether our axioms imply the existence of such apartments. 
\end{Remark}

Let $(\cV_j)_{j\in J}$ be a family of filters. We say that $(\cV_j)$ is \textbf{friendly} if there exists an apartment containing $\bigcup_{j\in J}\cV_j$.

\begin{Proposition}(see \cite[Proposition 5.17]{hebert2020new})\label{p_frndly_fac}
Let $a,b\in \I$ and $F_a$, $F_b$ be two faces of $\I$ based at $a$ and $b$ respectively. Then $F_a\cup F_b$ is friendly  if and only if $\{a,b\}$ is friendly.
\end{Proposition}

\begin{proof}

By~\ref{a_oco}, any apartment containing $F_a\cup F_b$ contains $\cl^\Delta(F_a)\cup \cl^\Delta(F_b)\supset \overline{F_a}\cup \overline{F_b}\supset \{a,b\}$. This proves the sense $\Rightarrow$. 

Let $a,b\in \I$ be friendly. By~\ref{a_oco}, we can assume that $C_a:=F_a$ and $C_b:=F_b$ are two alcoves. Let $A$ be an apartment containing $\{a,b\}$. Let $C_{a,A}$ be an alcove of $A$ based  at $a$ and of the same sign as $C_a$. Let $\Gamma_a=(C_1,\ldots, C_n)$ be a gallery of alcoves (in the sense of \ref{sss_gal}) from $C_{a,A}=C_1$ to $C_a=C_n$. By Corollary~\ref{c_fcs_pnt_frnd}, there exists an apartment $B_1$ containing $C_1$ and $C_2$. Let $P$ be the panel of $B_1$ dominated by both $C_1$ and $C_2$. Then by~\ref{a_oco}, $P$ is contained in $A$. Let $D_1,D_2$ be the two half-apartments of $A$ delimited by $P$. Then either $D_1$ or $D_2$ contains $b$. Up to renumbering, we can assume $D_1$ contains $b$. By Proposition~\ref{p_loc_SC}, there exists an apartment $A_1$ containing $C_2$ and $D_1$ and then $A_1$ contains $C_2$ and $b$. By induction, we deduce the existence of an apartment $B$ containing $C_n=C_a$ and $b$. 

Let $C_{b,B}$ be an alcove based at $b$ and  of the same sign as $C_b$. Let $\Gamma_b=(C_1',\ldots,C_{n'}')$ be a gallery of alcoves from $C_{b,B}$ to $C_b$. Let $P'$ be the panel between $C_1'$ and $C_2'$. Let $D_1',D_2'$ be the two half-apartments delimited by $P'$, with $D_1'$ containing $C_a$. Then by Proposition~\ref{p_loc_SC}, there exists an apartment containing $D_1'\cup C_2'$. By induction, we obtain the existence of an apartment containing $C_a$ and $C_{n'}'=C_b$, which proves the proposition. 
\end{proof}

\section{Splitting of apartments}\label{s_spltng_apt}
Recall that $\I$ is a weak masure and  that a subset of $\A$ is said to be enclosed if it can be written as a finite intersection of half-apartments of $\A$.
In this section, we use the sundial configuration to prove that if $A$ is an apartment and $\fQ$ is a sector-germ, then we can decompose $A$ as a union $A=\bigcup_{i=1}^k P_i$, where $k\in \N$  and the $P_i$ are enclosed subsets such that $P_i\cup \fQ$ is contained in some apartment for all $i$. This result will be important in our proof of~\ref{a_ma2}. First, it explains why the half-spaces considered in~\ref{a_ma2} are only related to real roots and not to the imaginary roots, and second, it will enable us to prove that the image of a segment by a retraction is a Hecke path, which is one of the ingredients of this proof. These decompositions first appeared in \cite{hebert2021distances}.

Let $A$ be an apartment and $\fQ$ be a sector-germ of $\I$. Let $\epsilon$ be the sign of $\fQ$. We define $d(\fQ,A)\in \N$ as the minimal possible length of a gallery of sector-germs of $\I$ of sign $\epsilon$ from a sector-germ of $A$ to $\fQ$.

\subsection{Splitting of apartments}

\begin{Lemma}\label{l_nnpty_int}

Let $A$ be a finite dimensional affine space and  $P_1,\ldots,P_k$ be polyhedra  of $A$ such that $A=\bigcup_{i=1}^k P_i$. Let $\Omega\subset A$ be closed, non-empty and such that $\Omega=\overline{\mathring{\Omega}}$. Let $J=\{i\in \llbracket 1,k\rrbracket\mid \mathring{P_i}\cap \mathring{\Omega}\neq \emptyset\}$. Then $\Omega=\bigcup_{j\in J}P_j\cap \Omega$. 
\end{Lemma}
\begin{proof}
 Let $\cM_i$ be the union of the supports of the facets of $P_i$. Let $\cM=\bigcup_{i=1}^k \cM_i$.  As $\cM$ is finite, the set $\cN=\bigcup_{M\in \cM} M$ has measure $0$, for any Lebesgue measure on $A$. Therefore $\mathring{\Omega}\setminus \cN$ has full measure in $\mathring{\Omega}$ and in particular it is dense in $\mathring{\Omega}$.

Let $|\cdot |$ be  a norm on $A$. Let $x\in  \Omega$. By assumption, there exists $(x_n)\in \mathring{\Omega}^\N$ such that $|x-x_n|\leq 1/(n+1)$, for all $n\in \N$. Then there exists $y_n\in \mathring{\Omega}\setminus \cN$ such that $|y_n-x_n|\leq \frac{1}{n+1}$. Then $(y_n)$ converges to $x$. Let $i\in \llbracket 1,k\rrbracket$ be such that $E:=\{n\in \N\mid y_n\in P_i\}$ is infinite. Then $x\in P_i$. Moreover, for $n\in E$, $y_n\in \mathring{\Omega}\cap \mathring{P_i}$, which proves that this set is non-empty. Therefore  $\Omega\subset \bigcup_{j\in J}P_j\cap \Omega$, and the lemma follows.
\end{proof}

Recall that a subset of an apartment  $A$ is called enclosed if it is a finite intersection of half-apartments of $A$.

\begin{Definition}\label{d_dec}
Let $A$ be an apartment of $\I$ and $\fQ$ be a sector-germ of $\I$. A \textbf{decomposition of $A$ with respect to $\fQ$}\index{decomposition of an apartment} is a triple $(P_i,A_i,f_i)_{i\in \llbracket 1,k\rrbracket}$ such that $k\in \N$ and for $i\in \llbracket 1,k\rrbracket$:\begin{enumerate}
\item $P_i$ is contained in $A$,  enclosed in $A$ and has non-empty interior (in $A$),

\item $A_i$ is an apartment containing $P_i\cup \fQ$,

\item $f_i:A\rightarrow A_i$ is an apartment isomorphism fixing $P_i$,

\item $\bigcup_{j\in \llbracket 1,k\rrbracket} P_j=A$.  
\end{enumerate}
\end{Definition}

\begin{figure}[h]
\centering
\includegraphics[scale=0.5]{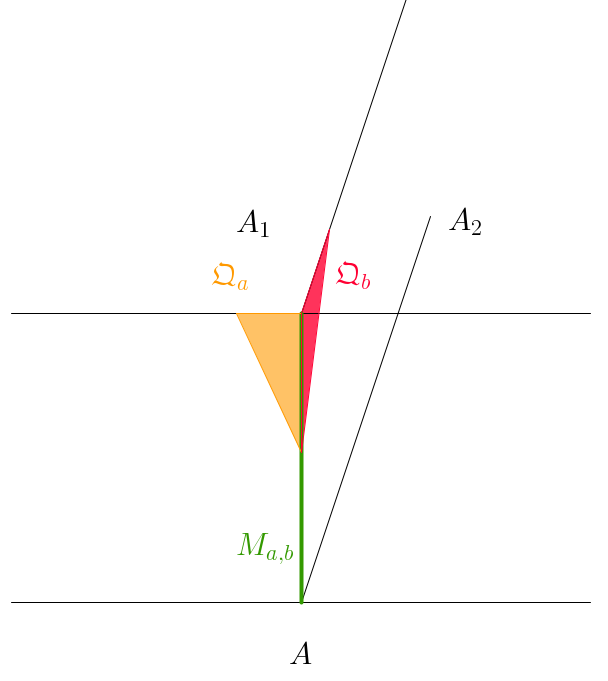}
\caption{In  the figure, $\fQ_a$ is a sector-germ of an apartment $A$. The sector-germ $\fQ_b$ is adjacent to $\fQ_a$ and not contained in $A$.  Let $D_1$ (resp. $D_2$) be the half-apartments delimited by $M_{a,b}$ and on the left (resp. right) of $M_{a,b}$. The apartments $A_i$  contains $\fQ_b$ and $D_i$, for $i\in \{1,2\}$.  The wall $M_{a,b}$ separates $\fQ_a$ and $\fQ_b$ in $A_1$. }\label{f_splt_2}
\end{figure}

\begin{figure}[h]
\centering
\includegraphics[scale=0.5]{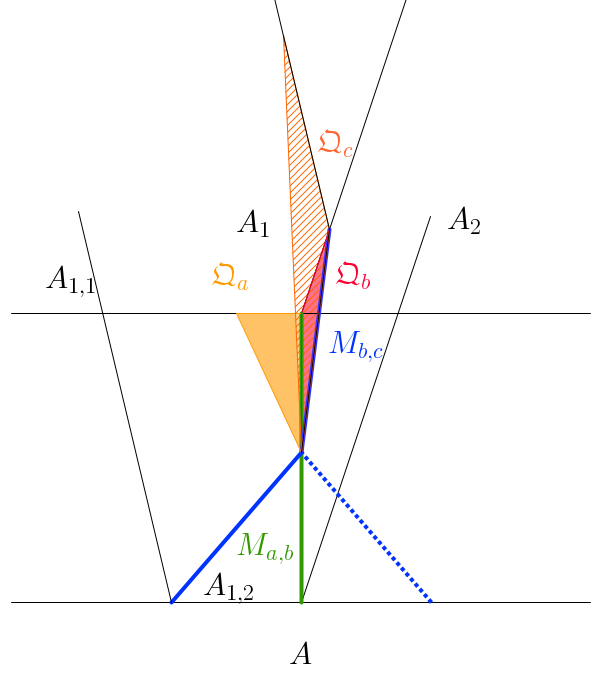}
\caption{We keep the same notation as in Figure~\ref{f_splt_2}. The gallery $(\fQ_a,\fQ_b,\fQ_c)$ is minimal. The apartment $A_1$ is divided in two half-apartments: $D_{1,1}$, on the left of $M_{b,c}$ and $D_{1,2}$ on the right of $M_{b,c}$. Then there exists an apartment $A_{1,i}$ containing $D_{1,i}$ and $\fQ_c$, for both $i\in \{1,2\}$. }\label{f_splt_4}
\end{figure}

\begin{figure}[hh]
\centering
\includegraphics[scale=0.3]{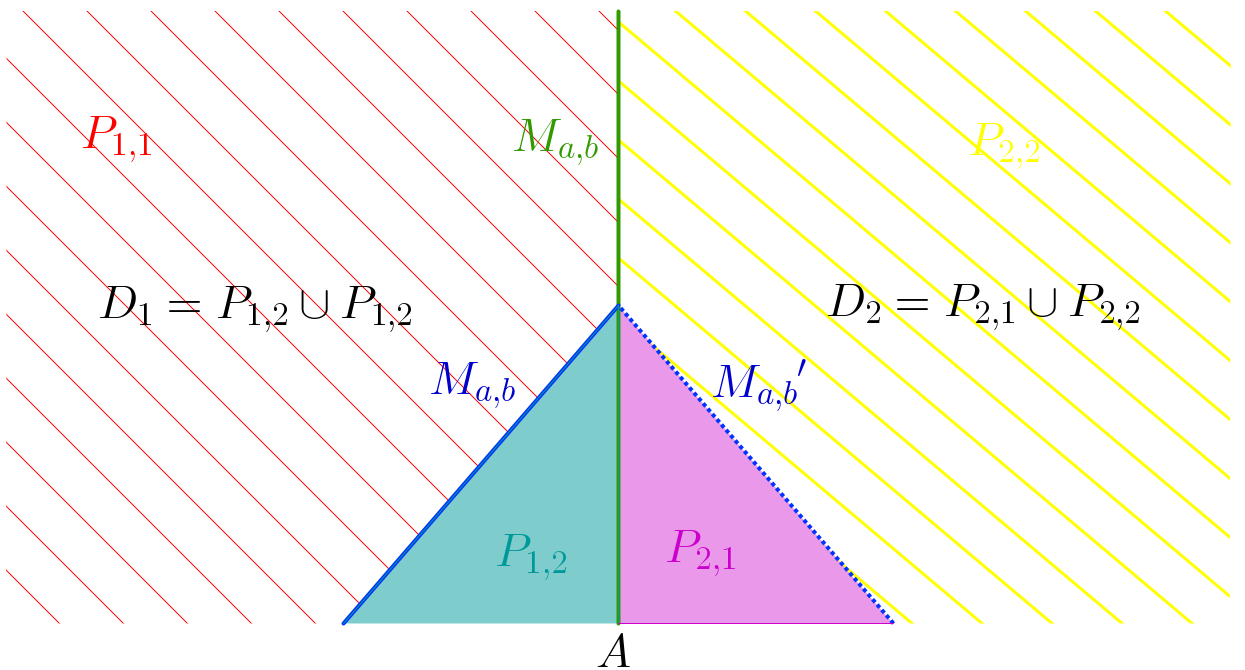}
\caption{We keep the same notation as in Figure~\ref{f_splt_2} and Figure~\ref{f_splt_4}. Similarly as we divided $A_1$ in two half-apartments, we can divide $A_2$ in two half-apartments $D_{2,1}$ and $D_{2,2}$, in such a way that $D_{2,i}\cup \fQ_c$ is contained in an apartment $A_{2,i}$, for $i\in \{1,2\}$. The wall $M'_{a,b}$ is $D_{2,1}\cap D_{2,2}$. Let $i,j\in \{1,2\}$. We set $P_{i,j}=D_i\cap D_{i,j}$. Then $D_{i,j}\subset A_{i,j}$. The apartment $A_{i,j}$ contains $D_{i,j}\cup \fQ_c$. Therefore $A$ can be divided in four polyhedra, each of them being contained in an apartment containing $\fQ_c$, which is at distance $2$ from $\fQ_a\subset A$.}\label{f_splt_4_apt}
\end{figure}

\begin{Proposition}\label{p_splt_apt}(see Figures~\ref{f_splt_2}, \ref{f_splt_4} and \ref{f_splt_4_apt})
Let $A$ be an apartment of $\I$, $\fQ$ be a sector-germ of $\I$ and $n=d(\fQ,A)$. Then: \begin{enumerate}
\item We can write $A=D_1\cup D_2$, where $D_1$ and $D_2$ are opposite half-apartments of $A$ such that for both $i\in \{1,2\}$, there exists an apartment $A_i$ containing $D_i$ and such that $d(A_i,\fQ)=n-1$. 

\item There exists $k\in \llbracket 1, 2^n\rrbracket$ such that $A$ admits a decomposition $(P_i,A_i,f_i)_{i\in \llbracket 1,k\rrbracket}$ with respect to $\fQ$, in the sense of Definition~\ref{d_dec}.
\end{enumerate}
\end{Proposition}

\begin{proof}
(1) Let $\epsilon$ be the sign of $\fQ$.  Let $\fQ'$ be a sector-germ of $A$ of sign $\epsilon$ and such that $d(\fQ,\fQ')=n$. Let $\Gamma=(\fQ_1,\ldots,\fQ_n)$ be a gallery of sector-germs from $\fQ$ to $\fQ'$.  Let $\fQ''=\fQ_{2}$. We have  $d(\fQ,\fQ'')=1$, $d(\fQ'',\fQ')=n-1$. Then by \ref{a_sc}, we can write $A=D_1\cup D_2$, where for both $i\in \{1,2\}$, there exists an apartment $A_i$ containing $D_i\cup \fQ''$, which proves (1).

(2) We prove it by induction on $n=d(\fQ,A)$. We first  prove the existence of the $P_i$, possibly with empty interior.  If $n=0$, we can take $P_1=A$ and there is nothing to prove. We now assume that $n\geq 1$ and that the result is true for all apartment $A'$ such that $d(\fQ,A')\leq n-1$. Using (1), we write $A=D_1\cup D_2$, where $D_1$ and $D_2$ are opposite half-apartments of $A$ such that for both $i\in \{1,2\}$, there exists an apartment $B^{(i)}$ containing $D_i$ and such that $d(B^{(i)},\fQ)=n-1$.

Let $j\in \{1,2\}$. We apply the induction assumption to $B^{(j)}$. Let $(Q_i,A_i^{(j)},f_i^{(j)})_{i\in \llbracket 1,k_j\rrbracket}$  be a decomposition of $B^{(j)}$ with respect to $\fQ$, where $k_j\in \llbracket 1,2^{n-1}\rrbracket$. By Lemma~\ref{l_iso_sect}, there exists an apartment isomorphism $f_j:A\rightarrow B^{(j)}$ fixing $A\cap B^{(j)}$. For $i\in \llbracket 1,k_j\rrbracket$, we set $P_i^{(j)}=D_j\cap Q_i^{(j)}$. Then $A_j^{(i)}$ contains $P_i^{(j)}\cup \fQ$ and $f_i^{(j)}\circ f_j:A\rightarrow A_i^{(j)}$ fixes $P_i^{(j)}$. Moreover, $P_i^{(j)}=Q_i^{(j)}\cap D_j=f_j(Q_i^{(j)})\cap D_j$ is enclosed. We have $A=D_1\cup D_2=(\bigcup_{i\in \llbracket 1,k_1\rrbracket} P_i^{(1)})\cup (\bigcup_{i\in \llbracket 1,k_1\rrbracket} P_i^{(2)})$. Using Lemma~\ref{l_nnpty_int} (with $\Omega=A=\A$), we can keep only the sets having non-empty interior among the elements of $\{P_i^{(j)}, j\in \{1,2\},i\in \llbracket 1,k_j\rrbracket\}$ and thus we obtain (2) by induction.
\end{proof}

A \textbf{parametrized line} of $\I$ is a  map $\fs:\R\rightarrow \I$ such that there exist $a,b\in \A$, with $a\neq b$ and $\psi:\A\rightarrow A$ an apartment isomorphism, with $A$ an apartment, such that $\fs(t)=\psi(ta +(1-t)b)$, for all $t\in \R$.

We deduce the following corollary, which is often proved using a compactness argument:

\begin{Corollary}\label{c_splt_seg}
Let $A$ be an apartment of $\I$ and $\fs:\R\rightarrow A$ be a parametrized line. Then there exist $k\in \N$ and $t_0,\ldots,t_k\in \R$ such that  $t_0<t_1<\ldots < t_k$  and  such that:\begin{enumerate}
\item for all $i\in \llbracket 0,k-1\rrbracket$, $\fs([t_i,t_{i+1}])\cup \fQ$ is contained in an apartment,

\item  $\fs(]-\infty,t_0])\cup \fQ$ (resp.  $\fs([t_k,+\infty[)\cup \fQ$) is contained in an apartment.
\end{enumerate}  
\end{Corollary}

\subsection{Topological consequences}

We equip $\A$ with its topology of a finite dimensional affine space. Using \ref{a_ma1}, this equips every apartment with a topology.

\begin{Lemma}\label{l_clsd_inter}
Let $A$ and $B$ be two apartments of $\I$. Then $A\cap B$ is closed (in $A$ or $B$). More precisely, if $(x_n)\in (A\cap B)^{\N}$, then $(x_n)$ converges in $A$ if and only if $(x_n)$ converges in $B$, and when $(x_n)$ converges, the limit in $A$ equals the limit in $B$.
\end{Lemma}

\begin{proof}
We begin by the case where $A\cap B$ contains a sector. By symmetry, we can assume that $A\cap B$ contains a positive sector. Up to changing the choice of $(\A,C^v_f)$, we can assume that $A=\A$ and $\fQ=\fQ_{+\infty}$. 

Let $(x_n)\in (\A\cap B)^\N$ be a sequence converging in $\A$. Let $x$ be the limit of $(x_n)$, regarded as a sequence of $\A$. Let us prove the existence of a sequence $(y_n)\in (\A\cap B)^\N$ such that $y_n\to x$ (in $\A$) and $y_n\in x+\overline{C^v_f}$ for all $n\in \N$. 

Let $(f_i)_{i\in J}$ be a basis of $\A^*$ such that $J\supset I$ and $f_i=\alpha_i$, for all $i\in I$. Let $(b_j)_{j\in J}$ be the dual basis of $(f_i)$, i.e $f_i(b_j)=\delta_{i,j}$ for $i,j\in J$. Then $\overline{C^v_f}=\sum_{i\in I}\R_{\geq 0} b_i+\sum_{i\in J\setminus I} \R b_i$. Let $|\cdot|_\infty$ be the norm on $\A$ defined by $|\sum_{j\in J} t_j b_j|_\infty=\max_{j\in J} |t_j|$, for $(t_j)\in \R^J$. 

For $n\in \N$ and $j\in J$, set $t_j(n)=f_j(x_n)$ and $t_j=f_j(x)$.  Let  $n\in \N$. Set \[y_n=\sum_{i\in I} \max(t_i,t_i(n))b_i+\sum_{i\in J\setminus I} t_i b_i.\]
 Then $y_n\in (x+\overline{C^v_f})\cap (x_n+\overline{C^v_f})\subset \A\cap B$ and $y_n\to x$. Let $\Omega=\{y_n\mid n\in \N\}$. By \ref{a_oco}, $\cl^\Delta(\Omega)\subset \A\cap B$ and there exists an apartment isomorphism $\psi:A\rightarrow B$ fixing $\cl^{\Delta}(\Omega)$. Therefore $x\in \A\cap B$ and $(x_n)=(\psi(x_n))$ converges to $x$ in $B$ and thus $x\in \A\cap B$. By symmetry, we deduce the result in the case where $\A\cap B$ contains a sector.

We no longer assume that $A\cap B$ contains a sector. Let $\fQ_B$ be a sector germ of $B$. Using Proposition~\ref{p_splt_apt}, we find a decomposition $(P_i,A_i,\psi_i)_{i\in \llbracket 1,k\rrbracket}$ of $A$ with respect to $\fQ_B$. Assume that $A\cap B$ is non-empty. Let $(x_n)\in (A\cap B)^{\N}$ be a sequence converging in $A$. Let $x$ be the limit of $(x_n)$ in $A$. 

For $i\in \llbracket 1,k\rrbracket$, set $E_i=\{n\in \N\mid x_n\in P_i\}$.   Let $i\in \llbracket 1,k\rrbracket$ be such that $E_i$ is infinite. Let $\phi_i:A_i\rightarrow B$ be the apartment isomorphism fixing $A_i\cap B$, which exists by Lemma~\ref{l_iso_sect}. Then $\phi_i\circ \psi_i$ fixes $\{x_\ell\mid \ell\in \N\}\cap P_i$. As $P_i$ is enclosed, $x\in P_i$. Therefore $(x_\ell)_{\ell\in E_i}=(\psi_i(x_\ell))$ converges to $x$ in $A_i$. We have $\phi_i(x_\ell)=x_\ell\in A_i\cap B$, for all $\ell\in E_i$. By the first case we treated, $x\in A_i\cap B$ and $(x_\ell)_{\ell\in E_i}$ converges to $x$ in $B$. As this is true for every $i\in \llbracket 1,k\rrbracket$ such that $E_i$ is infinite, we deduce that $(x_\ell)$ converges to $x$ in $B$. By symmetry, we obtain the lemma. 
\end{proof}

\begin{Proposition}\label{p_cont_retrct}
Let $\fQ$ be a sector-germ of $\I$ and $A$ be an apartment of $\I$ containing $\fQ$. Let $\rho=\rho_{A,\fQ}:\I\rightarrow A$ be the retraction onto $A$ centered at $\fQ$. Let $B$ be an apartment of $\I$. Then $\rho_{|B}:B\rightarrow \A$ is continuous (for the affine topologies on $B$ and $A$).
\end{Proposition}

\begin{proof}
By Proposition~\ref{p_splt_apt}, we have a decomposition $(P_i,B_i,\psi_i)_{i\in \llbracket 1,n\rrbracket}$ of $B$ with respect to $\fQ$.
 For all $i\in \llbracket 1,n\rrbracket$, we denote by $\phi_i$ the apartment isomorphism $B_i\rightarrow A$ fixing $A\cap B_i$, which exists by Lemma~\ref{l_iso_sect}. Then  $\rho|_{P_i}=\phi_i\circ \psi_i|_{P_i}$ for all $i\in \llbracket 1,n\rrbracket$. 

Let $(x_k)\in B^\N$ be a converging sequence and $x=\lim x_k$. For $i\in \llbracket 1,n\rrbracket$, set $E_i=\{k\in \N\mid x_k\in P_i\}$. Let $i\in \llbracket 1,n\rrbracket$ be such that $E_i$ is infinite. Then $x\in P_i$  and for all $k\in E_i$, we have  $\rho(x_k)=\phi_i\circ \psi_i(x_k)\underset{k\to +\infty}{\rightarrow} \phi_i\circ \psi_i(x)=\rho(x)$. Therefore $(\rho(x_k))_{k\in\N}$ converges to $\rho(x)$, which proves the proposition since $A$ and $B$ are metric spaces.
\end{proof}

\section{Axiom \ref{a_co}}

In this section, we assume that $\I$ is a weak masure: it satisfies \ref{a_ma1}, \ref{a_wma3}, \ref{a_oco} and \ref{a_sc}.

Let $\fQ,\fQ'$ be two sector-germs (at infinity) of $\I$. By \ref{a_wma3}, there exists an apartment containing $\fQ$ and $\fQ'$. However,  in general, there exists $x\in \I$ such that $x+\fQ$ and $x+\fQ'$ is not contained in any apartment. For example if $\I$ is  a tree (not isometric to a line), then it is easy to find $x\in \I$ and $\fQ,\fQ'$ such that $(x+\fQ)\cup (x+\fQ')$ is not contained in any apartment.

In this section, we prove axiom \ref{a_co} (see below), which yields a sufficient condition for the existence of an apartment containing $x+\fQ$ and $x+\fQ'$. This criterion is based on the relative position of the (local) germs of $x+\fQ$ and $x+\fQ'$ at $x$. 

We then deduce a characterization of $\A$ (see Corollary~\ref{c_char_apt_alc}). We will use this characterization later in order to prove that the intersection of two apartments is enclosed (see Section~\ref{s_encls_int}).

\subsection{Definitions of opposition and of axiom \ref{a_co}}\label{ss_co}

Let $\fQ$ be a sector-germ (at infinity) of $\I$ and $x\in \I$. We set $\fQ_x=\germ_x(x+\fQ)$\index[notation]{q@$\fQ_x$}. This is an alcove of $\I$. Note that by Proposition~\ref{p_frndly_fac}, if $\fQ,\fQ'$ are two sector-germs and $x\in \I$, then there exists an apartment containing $\fQ_x$ and $\fQ'_x$.

\begin{Definition}
\begin{enumerate}
\item Let $\fQ,\fQ'$ be two sector-germs of $\A$ and $x\in \A$. We say that $\fQ$ and $\fQ'$ are opposite at $x$ if there exists $\epsilon\in \{-,+\}$  and $w\in W^v$ such that $\fQ_x=\germ_x(x+\epsilon w.C^v_f)$ and $\fQ'_x=\germ_x(x-\epsilon w.C^v_f)$. 

\item Let $\fQ,\fQ'$ be two sector-germs of $\I$ and $x\in \I$. We say that $\fQ_x$ and $\fQ_x'$ are \textbf{opposite}\index{opposite, opposite at $x$} or that $\fQ$ and $\fQ'$ are \textbf{opposite at $x$} if there exists an apartment $A$ containing $\fQ_x,\fQ'_x$ and such that $\fQ_x$ and $\fQ'_x$ are opposite in $A$, i.e for any apartment isomorphism $f: A\rightarrow \A$, $f(\fQ)$ and $f(\fQ')$ are opposite at $x$. We denote it $\fQ_x\opp \fQ'_x$\index[notation]{o@$\opp$}.

\item Let $\fQ,\fQ'$ be two sector-germs of $\I$. We say that $\fQ$ and $\fQ'$ are \textbf{opposite} if there exists an apartment $A$ containing $\fQ\cup \fQ'$ in which $\fQ$ and $\fQ'$ are opposite, i.e for some apartment isomorphism $f:A\rightarrow \A$, there exist $w\in W^v$ and $\epsilon\in \{-,+\}$ such that $f(\fQ)=\germ_\infty(\epsilon w.C^v_f)$ and $f(\fQ')=\germ_\infty(-\epsilon w.C^v_f)$. 
\end{enumerate}
\end{Definition}

Note that if $\fQ$ and $\fQ'$ are two sector-germs such that $\fQ_x \opp \fQ'_x$, for some $x\in \I$, then $\fQ$ and $\fQ'$ are opposite. 

We now introduce axiom \ref{a_co} (see \cite[1.5.3]{parreau2000immeubles} or \cite[3.1.1 Proposition]{rousseau2023euclidean} for the case of  affine buildings). 

\begin{enumerate}[label=\blue{(CO)}]
 \item\label{a_co} (``chambres opposées''?) Let $\fQ$ and $\fQ'$ be two sector-germs at infinity of $\I$ and let $x\in \I$. If $\fQ$ and $\fQ'$ are opposite at $x$, then there exists an apartment containing $(x+\fQ)\cup (x+\fQ')$.   \axiom{co@\ref{a_co}}
\end{enumerate}

\begin{Lemma}\label{l_cl_seg}
Let $x,x',x''\in \A$. We assume that $x\in ]x',x''[$ and that there exists a vectorial chamber $C^v$ of $\A$ such that $x'\in x+C^v$.  Then $\cl^{\Delta}([x',x''])$ contains $(x'-C^v)\cap (x''+C^v)$. 
\end{Lemma}

\begin{proof}
Let $\Omega\in \cl^{\Delta}([x',x''])$. By definition, there exists $(k_\alpha)\in (\R\cup \{+\infty\})^{\Delta}$ such that $\Omega\supset \bigcap_{\alpha\in \Delta} D(\alpha,k_\alpha)\supset [x',x'']$. Up to exchanging $x'$ and $x''$, we can assume that $C^v$ is positive. Let $w\in W^v$ be such that $C^v_f=w.C^v$. Let $\alpha\in \Delta$. First assume $w^{-1}.\alpha\in \Delta_+$. We have $\alpha(x'')+k_\alpha\geq 0$ and thus $\alpha(x''+w.C^v_f)+k_\alpha=\alpha(x'')+\R_{>0}+k_\alpha> 0$. Consequently, $x''+C^v\subset D(\alpha,k_\alpha)$. Similarly, if $w^{-1}.\alpha\in \Delta_-$, then we have $x'-C^v\subset D(\alpha,k_\alpha)$. Therefore $(x'-C^v)\cap (x''+C^v)\subset \bigcap_{\alpha\in w.\Delta_+}D(\alpha,k_\alpha)\cap \bigcap_{\alpha\in w.\Delta_-}D(\alpha,k_\alpha)=\bigcap_{\alpha\in \Delta} D(\alpha,k_\alpha)\subset \Omega$, which proves the lemma.
\end{proof}

\begin{Lemma}\label{l_opp}
Let $\fQ,\fQ'$ be two sector-germs of $\I$. Then:
\begin{enumerate}
\item Let $x\in \I$. We assume that $\fQ$ and $\fQ'$ are opposite at $x$. Then for any apartment $B$ containing $\fQ_x\cup \fQ'_x$, the alcoves $\fQ_x$ and $\fQ'_x$ are opposite at $x$. Moreover, if $B$ and $B'$ are two apartments containing $\fQ_x\cup \fQ'_x$, then $B\cap B'$ contains $x$ in its interior and there exists an apartment isomorphism from $B$ to $B'$ fixing a neighborhood of $x$ in $B$.

\item We assume that $\fQ$ and $\fQ'$ are opposite. Then there exists a unique apartment containing $\fQ$ and $\fQ'$. 
\end{enumerate}
 
\end{Lemma}

\begin{proof}
(1) Let $A$ be an apartment containing $\fQ_x$ and $\fQ'_x$ be such that $\fQ_x$ and $\fQ'_x$ are opposite in $A$. Let $\Omega\in \fQ_x$ and  $\Omega'\in \fQ_x'$ be such that $\Omega\cup \Omega'\subset A\cap B$. We assume that $\Omega$ and $\Omega'$ are convex, which is possible, up to reducing them. We identify $A$ and $\A$. Let $C^v$ be the vectorial chamber of $\A$ such that $\fQ_x=\germ_x(x+C^v)$. Then $\fQ'=\germ_x(x-C^v)$. 

Let $\lambda\in C^v$. Let $t\in \R_{>0}$ be such that $x+ t\lambda\in \Omega$ and  $x-t\lambda\in \Omega'$. Then by \ref{a_oco}, there exists an apartment isomorphism $f:\A\rightarrow B$ fixing $\cV:=\cl^\Delta([x-t\lambda,x+t\lambda])$. By Lemma~\ref{l_cl_seg}, $\cV$ contains $x$ in its interior (which means that for all $\Omega''\in \cV$, we have $x\in \mathring{\Omega}''$). Therefore $f$ fixes an open neighborhood $\Omega''$ of $x$. Then $\fQ_x$ and $\fQ_x'$ are contained in $\Omega''$ and they are opposite  in $B$.

(2) By assumption, there exists an apartment $A$ containing $\fQ\cup \fQ'$ and such that $\fQ$ and $\fQ'$ are opposite in $A$.  Let $B$ be an apartment containing $\fQ\cup \fQ'$.  Up to changing the choice of $(A,C^v_f)$ and to exchanging the roles of $\fQ$ and $\fQ'$, we can assume that $A=\A$ and $C^v=C^v_f$. Let $x_+\in C^v_f$ and $x_-\in (-C^v_f)$ be such that $\A\cap B\supset (x_++C^v_f)\cup (x_--C^v_f)$. Then by \ref{a_oco} and Lemmas~\ref{l_clsd_inter} and \ref{l_cl_seg}, $\A\cap B\supset \bigcup_{x'\in x_++C^v_f,x''\in x_--C^v_f}(x'-C^v_f)\cap (x''+C^v_f)=\A$. Consequently, $\A\subset \A\cap B$. By Lemma~\ref{l_iso_sect}, there exists an apartment isomorphism $f:\A\rightarrow B$ fixing $\A\cap B=\A$. Then $f(\A)=B=\A$, which proves (2). 
\end{proof}

\subsection{Axiom~\ref{a_co} and characterization of $\A$}

We now introduce the Tits preorder $\leq$ on $\I$.
\begin{Definition}
Let $x,y\in \I$. We write $x\leq y$ if there exists an apartment $A$ containing $x$ and $y$ and an apartment isomorphism $f:A\rightarrow \A$ such that $f(y)-f(x)\in \cT$ (where $\cT$ is the Tits cone of $\A$).

If $x,y\in \I$ we say that $\{x,y\}$ is preordered if $x\leq y$ or $y\leq x$. 
\end{Definition}

\begin{Remark}
\begin{enumerate}
\item The Tits preorder is a crucial tool in the theory of masures. It is  transitive and thus it is indeed a preorder (see \cite[Théorème 5.9]{rousseau2011masures}). Its transitivity enables to obtain the existence of apartments containing certain pairs of points. It is also very used in the theory of Kac--Moody groups over valued fields: if $\I$ is associated with a Kac--Moody group $G$, it enables to define the semi-group $G^+$, which is used in the representation theory of $G$ (for example to define Hecke algebras, see \cite{gaussent2014spherical} and \cite{bardy2016iwahori}).

\item If $\I$ is an affine building (i.e if $W^v$ is finite), then $\leq$ is trivial: for all $x,y\in \I$, we have $x\leq y$. 
\end{enumerate}

\end{Remark}

The following proposition is proved (but not stated) in the proof of \cite[Proposition 5.4]{rousseau2011masures}. We reproduce the proof here.

\begin{Proposition}\label{p_line}
Let $x,y\in \I$ be preordered and distinct. Let $\delta_x$ (resp. $\delta_y$) be a ray of $\I$ based at $x$ (resp. $y$) and containing $y$ (resp. $x$). Let $L=\delta_x\cup \delta_y$. Then there exists an apartment containing $L$ and $L$ is a line of any apartment containing $L$.
\end{Proposition}

\begin{proof}
Let $A_x$ (resp. $A_y$) be an apartment containing $\delta_x$ (resp. $\delta_y$). Let $\fQ_y$ be a sector-germ of $A_y$ such that $y+\overline{\fQ_y}$ contains $\delta_y$. By Proposition~\ref{p_splt_apt}, we have a decomposition $(P_i,B_i,f_i)_{i\in \llbracket 1,k\rrbracket}$ of $A_x$ with respect to $\fQ_y$.
 Let $\fs:\R_{\geq 0}\rightarrow A_x$ be an affine parametrization of $\delta_x$, i.e $\fs(0)=x$, $\fs$ is an affine map and $\fs(\R_{\geq 0})=\delta_x$. Set $t_0=0$. There exist $t_1,\ldots,t_n\in  \R_{\geq 0}$ such that $t_0<t_1<\ldots t_n$ and such that  for all $i\in \llbracket 0, n-1\rrbracket$ (resp. $i=n$), there exists $j(i)\in \llbracket 1,k\rrbracket$ such that $P_{j(i)}$ contains $\fs([t_i,t_{i+1}])$ (resp. $\fs([t_n,+\infty[)$). For $i\in \llbracket 0,n-1\rrbracket$, we denote by  $\cH_i$ the assumption \[\fs(t_i)+\overline{\fQ_y}\supset \fs([t_0,t_i]).\] Then $\cH_0$ is clear. Let $i\in \llbracket 1,n\rrbracket$ be  such that $\cH_{i-1}$ is true. The apartment $A_{j(i)}$ contains $\fs([t_{i-1},t_i])\cup \fQ_y$. Thus it contains $\fs(t_{i-1})+\overline{\fQ}_y$ and hence it contains $\fs([t_0,t_{i-1}])$, by $\cH_{i-1}$.  Consequently, $A_{j(i)}$ contains $\fs([0,t_i])$. By applying \ref{a_oco} with $\cV=\fs([t_0,t_i])$, $A=A_x$ and $A'=A_{j(i)}$, we obtain that $\fs|_{[t_0,t_i]}$ is a segment of $A_{j(i)}$. In $A_{j(i)}$, $\fs(t_i)+\overline{\fQ_{y}}$ is a translate of $\fs(t_{i-1})+\overline{\fQ_{y}}$. Therefore $\fs(t_i)+\overline{\fQ_{y}}$ contains $\fs([t_0,t_i])$, which proves that $\cH_i$ is true. In particular, $\cH_{n-1}$ is true.

Now $A_{j(n)}$ contains $\fs([t_n,+\infty[)\cup \fQ_y$. Therefore it contains $\fs(t_n)+\overline{\fQ_y}$ and by $\cH_{n-1}$, it contains $\fs([t_0,t_n])$. We deduce that $A_{j(n)}$ contains $\fs([t_0,+\infty[)=\delta_x$. In particular, $A_{j(n)}$ contains $\fs(t_0)+\overline{\fQ_y}$ and hence it contains $\delta_y$. Therefore $A_{j(n)}$ contains $L=\delta_x\cup \delta_y$.

Let now $B$ be an apartment containing $L$. By \ref{a_oco}, $\delta_x$ and $\delta_y$ are rays of $B$. As $\delta_x\cap \delta_y\supset [x,y]$, we deduce that $L$ is the unique line of $B$ containing $[x,y]$. 
\end{proof}

\begin{Corollary}\label{c_opp}
Let $\I$ be a weak masure. Then $\I$ satisfies \ref{a_co}. More precisely,  Let $\fQ,\fQ'$ be two sector-germs at infinity of $\I$ and $x\in \I$. We assume that $\fQ_x$ and $\fQ_x'$ are opposite. Then there exists an apartment $A$ containing $(x+\fQ)\cup (x+\fQ')$. Moreover, $\fQ$ and $\fQ'$ are opposite and $A$ is the unique apartment containing them.
\end{Corollary}

\begin{proof}
Let $B$ be an apartment containing $x+\fQ$ and $B'$ be an apartment containing $x+\fQ'$. Let $B''$ be an apartment containing $\fQ_x$ and $\fQ'$, which exists by Proposition~\ref{p_mas_th_iwa}. Then $B''$ contains $x$ and $\fQ'$, thus it contains $x+\fQ'$ and hence  $\fQ'_x$. By Lemma~\ref{l_opp}, $\fQ_x$ and $\fQ'_x$ are opposite at $x$. We identify $B''$ and $\A$. Write $\fQ_x=\germ_x(x+C^v)$, where $C^v$ is a vectorial chamber of $\A$. Then $x+\fQ'=x-C^v$.  Let $\lambda\in C^v$. Let $t\in \R_{>0}$ be  such that  $x+t\lambda\in \A\cap B$. Let $\delta_1$ be the ray of $B$ based at $x+t\lambda$ and containing $x$. Let $\delta_2$ be the ray of $A$ based at $x$ and containing $x+t\lambda$. Then by Proposition~\ref{p_line}, $\delta_1\cup \delta_2$ is contained in an apartment $A$ of $\I$ and it is a line of $A$. 

We have $\delta_2\subset x+\fQ$. By \ref{a_oco}, $A\cap B\supset \cl^{\Delta}(\delta_1)$. Let $\fQ_B'$ be the sector-germ of $B$ opposite to $\fQ$. By Lemma~\ref{l_cl_seg}, $\cl^\Delta(\delta_1)\supset \bigcup_{x'\in \delta_1} (x+\fQ)\cap (x+\fQ_B')=x+\fQ$. Therefore $A\supset x+\fQ'$. Similarly, $A\supset x+\fQ$. Then $A$ contains $\fQ_x$ and $\fQ'_x$ and by Lemma~\ref{l_opp}, $\fQ_x$ and $\fQ'_x$ are opposite in $A$. Therefore $\fQ$ and $\fQ'$ are opposite and we conclude with Lemma~\ref{l_opp}. 
\end{proof}

\begin{Corollary}\label{c_char_apt_alc}
Let $B$ be an apartment and $\fQ,\fQ'$ be two opposite sector-germs of $B$. Then $B=\{x\in \I\mid \fQ_x\opp \fQ'_x\}$.  In particular, $\A=\{x\in \I\mid \fQ_{+\infty,x}\opp \fQ_{-\infty,x}\}$. 
\end{Corollary}

\begin{proof}
We identify $B$ and $\A$. Write $\fQ=\germ_\infty(\epsilon w.C^v_f)$, with $w\in W^v$ and $\epsilon\in \{-,+\}$. Then $\fQ'=\germ_\infty(-\epsilon w.C^v_f)$ and for all $x\in \A$, we have $\fQ_x=\germ_x(x+\epsilon w.C^v_f)$ and $\fQ'_x=\germ_x(x-\epsilon w.C^v_f)$ and thus $\fQ_x$ and $\fQ'_x$ are opposite at $x$. Conversely, let $x\in \I$ be such that $\fQ_x$ and $\fQ'_x$ are opposite at $x$. Then by \ref{a_co}, there exists an apartment $A$ containing $(x+\fQ)\cup (x+\fQ')$. Then $A$ contains $\fQ\cup \fQ'$.  By Lemma~\ref{l_opp}, $A=\A$, which proves that $x\in \A$ and hence the lemma.
\end{proof}

\section{Hecke paths}\label{s_H_paths}

Hecke paths were introduced by Kapovich and Milson in \cite{kapovich2008path} in the context of Euclidean buildings,  in order to study the representation theory of a reductve group over $\C$. There were then generalized in the frameworks of masures by Gaussent and Rousseau in \cite{gaussent2008kac}. They are a crucial tool in the study of the representation theory of Kac--Moody groups over local fields.

Recall that $\I$ is a weak masure: it satisfies~\ref{a_ma1}, \ref{a_wma3}, \ref{a_sc} and~\ref{a_oco}. 

\subsection{Definition}

We now give the definition of Hecke paths. Under some additional assumption involving the Tits preorder, Hecke paths are exactly the images of line segments by retractions (see \cite[Theorem 6.2 and 6.3]{gaussent2008kac}.

\begin{Definition}{($(W^v_a,\fQ)$-chain)}\label{d_chain}
Let $a\in \A$ and $\xi,\xi'\in \A$. Let $C^v$ be a vectorial chamber of $\A$ and $\fQ$ be its germ at infinity. A $(W^v_a,\fQ)$-\textbf{chain} from $\xi$ to $\xi'$ with respect to $C^v$ is a pair $((\xi_0,\ldots,\xi_k),(\beta_1,\ldots,\beta_k))\in \A^{k+1}\times \Phi^{k}$ such that    $\xi_0=\xi$, $\xi_k=\xi'$ and for all $i\in \llbracket 1,k\rrbracket$, we have: \begin{enumerate}
\item $r_{\beta_i}(\xi_{i-1})=\xi_i$,

\item $\beta_i(\xi_{i-1})<0$,

\item $\beta_i(a)\in \Lambda_{\beta_i}$: $a$ belongs to  a wall of direction $\ker(\beta_i)$,

\item $\beta_i(C^v)<0$. 
\end{enumerate}
\end{Definition}

\begin{Definition}{(Hecke path)}\label{d_H_path}
 A  \textbf{piecewise affine path} (resp. piecewise affine ray, resp. piecewise affine line) $\gamma$ (resp. $\gamma:\R_{\geq 0}\rightarrow \A$ , resp. $\gamma:\R\rightarrow \A$) is a map $\gamma:[0,1]\rightarrow \A$ (resp. $\gamma:\R_{\geq 0}\rightarrow \A$, resp. $\gamma:\R\rightarrow \A$) such that $\gamma$ is continuous and there exist $n\in \N$  and  $0=t_0<t_1\ldots <t_n=1$ (resp. $0=t_0<\ldots <t_n\in \R_{\geq 0}$, resp. $t_0<\ldots <t_n$) such that $\gamma$ is affine on $[t_i,t_{i+1}]$ for all $i\in \llbracket 0,n-1\rrbracket$ (resp. and on $[t_n,+\infty[$, resp. and on $]-\infty,t_0]$ and $[t_n,+\infty[$).
 
 A \textbf{billiard path} (resp. a billiard ray, resp.  a billiard line) $\gamma:[0,1]\rightarrow \A$ (resp. $\gamma:\R_{\geq 0}\rightarrow \A$, resp. $\gamma:\R\rightarrow \R$) is a piecewise affine path such that there exists $\lambda\in \A$ such that we have $\gamma'_+(t_+),\gamma'_-(t_-)\in W^v\cdot \lambda$, for all $t_+\in (0,1],t_-\in [0,1)$ (resp. $t_+\in \R_{\geq 0}$ and $t_-\in \R_{>0}$, resp. $t_+,t_-\in \R$), where $\gamma'_+(t_+),\gamma'_-(t_-)$ are respectively the right and left derivative of $\gamma$ at $t_+$ and $t_-$. We then say that $\gamma$ is of shape $\lambda$. The \textbf{folding times} of $\gamma$ are the $t\in ]0,1[$ such that $\gamma'_{-}(t)\neq \gamma'_+(t)$. 
 
Let $C^v$ be a vectorial chamber of $\A$. A \textbf{Hecke path} (resp. Hecke ray, Hecke line) with respect to $C^v$ is a billiard path (resp. billiard ray, resp. billiard line) such that for all $t\in ]0,1[$ (resp. $t\in \R_>0$, resp. $t\in \R$), there exists a $(W^v_{\gamma(t)},\germ_\infty(C^v))$-chain from $\gamma'_-(t)$ to $\gamma'_+(t)$. 
\end{Definition}

Note that Hecke paths can also be defined as the images of line segments after successive foldings, see \cite[3]{hebert2024affine} for this approach, in a slightly different context. In \cite[1.8]{gaussent2014spherical}, the authors define the paths with respect to $-C^v_f$ and thus condition 4. is replaced by $\beta_i(C^v_f)>0$.

\subsection{Hecke paths and dominance order}

Recall that if $x,y\in \A$, we write $x\leq_{Q^\vee_\R} y$ to mean that $y-x\in \bigoplus_{i\in I}\R_{\geq 0} \alpha_i^\vee$.

\begin{Proposition}\label{p_H_paths}
Let $\gamma:[0,1]\rightarrow \A$ be a Hecke path with respect to $C^v_f$. Then: \begin{enumerate}
\item $\gamma$ is a billiard path of shape $\gamma'_+(0)$.

\item for all $t,t'\in ]0,1[$ such that $t<t'$, we have: \begin{equation}\label{e_ineq_H_pth}\gamma'_+(0)\geq_{Q^\vee_\R} \gamma'_-(t)\geq_{Q^\vee_\R} \gamma'_+(t)\geq_{Q^\vee_\R} \gamma'_-(t')\geq_{Q^\vee_\R} \gamma'_+(t')\geq_{Q^\vee_\R} \gamma'_-(1).\end{equation}

\item  We have $\gamma'_+(0)\geq_{Q^\vee_\R} \gamma(1)-\gamma(0)$. If moreover $(\alpha_i^\vee)_{i\in I}$ is positively free, then  this inequality is strict, unless $\gamma$ is a line segment.
\end{enumerate}
\end{Proposition}

\begin{proof}
(1) is clear. 

(2) By definition, if $\xi,\xi',a\in \A$ are such that there exists a $(W^v_a,\fQ_{+\infty})$-chain from $\xi$ to $\xi'$, then $\xi \geq_{Q^\vee_\R} \xi'$ (since all the $\beta_i$ appearing in Definition~\ref{d_chain} are negative roots). Therefore for all $t\in (0,1)$, we have $\gamma'_+(t)\leq_{Q^\vee_\R} \gamma'_-(t)$. Now if $t,t'\in (0,1)$ are such that $\gamma$ is affine on $[t,t']$ and $t<t'$, we have $\gamma'_+(t)=\gamma'_-(t')$ and hence we deduce \eqref{e_ineq_H_pth}.

(3) By integrating \eqref{e_ineq_H_pth} between $0$ and $1$, we obtain that $\gamma'_+(0)\geq_{Q^\vee_\R} \gamma(1)-\gamma(0)$. Moreover, if this inequality is an equality and $(\alpha_i^\vee)$ is positively free, then we necessarily have $\gamma'_+(0)= \gamma'_-(t)= \gamma'_+(t)= \gamma'_-(t')=\gamma'_+(t')= \gamma'_-(1)$ for all $t,t'\in (0,1)$ such that $t<t'$ and hence $\gamma$ is a line segment. 
\end{proof}

\begin{Remark}\label{r_H_pth}

\begin{enumerate}
\item If we take a Hecke path with respect to $-C^v_f$ in Proposition~\ref{p_H_paths}, then we get the same inequalities, but with $\leq_{Q^\vee_\R}$ replaced by $\geq_{Q^\vee_\R}$. 

\item Let $\lambda\in\overline{C^v_f}$ and $\gamma:[0,1]\rightarrow \A$ be a Hecke path of shape $\lambda$ with respect to $C^v_f$. Then by \cite[Lemma 1.3.13]{kumar2002kac}, we have the following property. For $t_-\in [0,1]$ where it makes sense, we write $\gamma'_+(t)=w'_+(t).\lambda$ and $\gamma'_-(t)=w'_-(t).\lambda$, where $w'_-(t), w'_+(t)\in W^v$ have minimal lengths for these properties. Then for all $t,t'\in [0,1]$ such that $0\leq t< t'\leq 1$, we have $w'_-(t)\leq w'_+(t) \leq w'_-(t')\leq w'_+(t')$, where we delete the derivatives that do not make sense (for $t=0$ or $t'=1$). Roughly speaking it means that each times that the path folds, its derivative moves away from $C^v_f$.
\end{enumerate}

\end{Remark}

\subsection{Simple foldings}\label{ss_smp_fld}

\begin{figure}[h]
\centering
\includegraphics[scale=0.4]{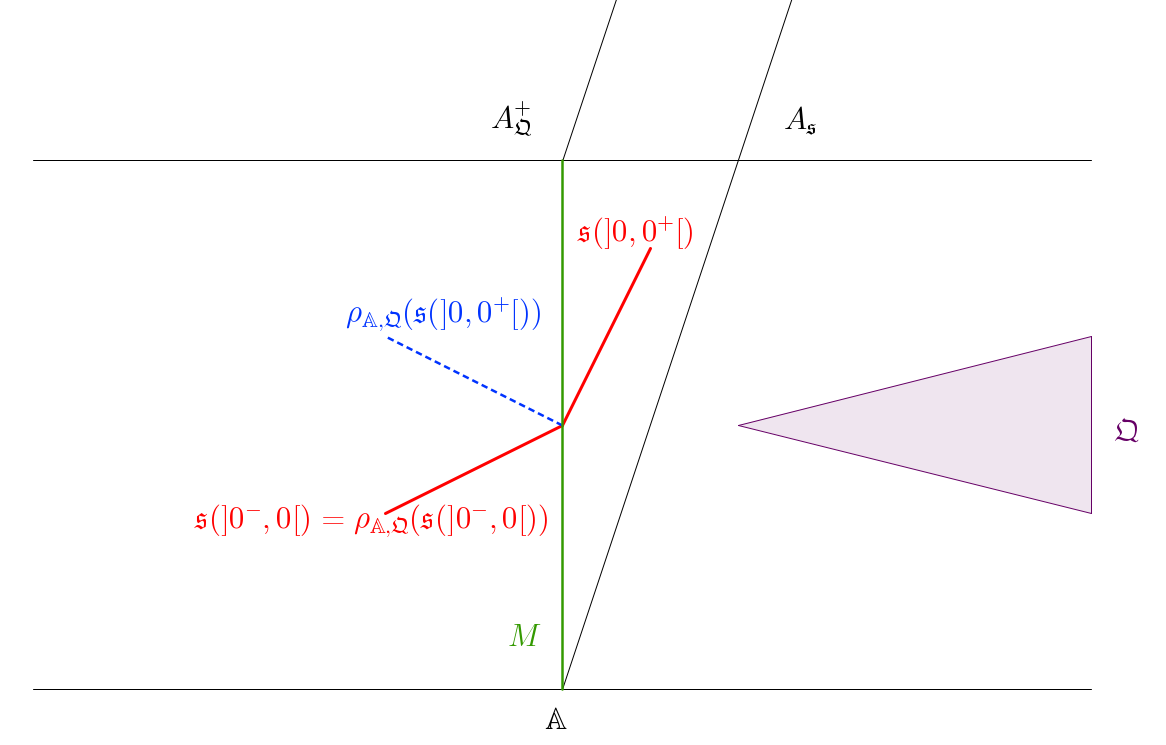}
\caption{Folding of a the retraction of a segment-germ $\fs:]0^-,0^+[\rightarrow \I$. The apartment $A_{\fs}$ contains $\fs$ but not $\fQ$. The apartment $\A$ (resp. $A_{\fQ}^+$) contain $\fQ$ and $\fs(]0^-,0[)$ (resp. $\fs(]0,0^+[)$). We have $\rho_{\A,\fQ}\circ \fs'_+(0)=r_M \circ\fs'_-(0)$, where $r_M$ is the reflection of $\A$ with respect to $M$. }\label{f_folding}
\end{figure}

In this subsection, we study the image of a segment-germ via a retraction, when the folding is ``simple'', i.e, when the length (the number $k$) of the chain is $1$ in Definition~\ref{d_chain} (see Figure~\ref{f_folding}).

\begin{Lemma}\label{l_loc_iso}
Let $x\in \I$ and $A$, $B$ be two apartments containing $x$. We assume that $A\cap B$ contains an alcove $C$ based at $x$. Let $f:A\rightarrow B$ be the apartment isomorphism fixing $C$. Then there exists a neighborhood $\Omega$ of $x$ in $A$ such that $f$ fixes $A\cap B\cap \Omega$. 
\end{Lemma}

\begin{proof}
Using apartment isomorphisms, we can assume that $A=\A$.  Let $C^v$ be the vectorial chamber of $\A$ such that $C$ is contained in $x+C^v$.  Let   $\Omega'\in C$ be such that $\Omega'\subset (A\cap B)\cap (x+C^v)$. Let $y\in \Omega'$.  Then for $z\in \Omega'$ close enough to $x$, we have $z\in y-C^v$  and hence up to reducing $\Omega$,~\ref{a_oco} implies that $\Omega:=\conv_\A(y,\Omega\cap \A\cap B)\subset \A\cap B$ and the existence of $g:\A\rightarrow B$ fixing $\Omega$, which proves the lemma.
\end{proof}

\begin{Lemma}\label{l_lc_hf_apt_alc}
Let $A$ and $B$ be two apartments such that $A\cap B\neq \emptyset$.  Let $x\in A\cap B$. We assume that there exists a half-apartment $D$ of $A$ such that  for some neighborhood $\Omega$ of $x$, we have $D\cap \Omega\subset A\cap B$. Let $P$ be a panel  based at $x$ and contained in the wall  of $D$. Let $C$ be the alcove of $A$ dominating $P$ and not contained in $D$. Then if $C\subset B$, we have  $x\in \In(A\cap B)$. 
\end{Lemma}

\begin{proof}
Using apartment isomorphisms, we may assume that $A=\A$.  Let $C^v$ be the vectorial chamber of $\A$ such that $C$ is contained in $x+C^v$.  Let   $\Omega\in C$ be such that $\Omega\subset (A\cap B)\cap (x+C^v)$. Let $y\in \Omega$.  Then for $z\in \Omega$ close enough to $x$, we have $z\in y-C^v$  and hence up to reducing $\Omega$, we have $\conv_\A(\Omega\cap D,z)\subset \A\cap B$, by~\ref{a_oco}. As $\conv_\A(\Omega\cap D,z)$ contains $x$ in its interior, we deduce the lemma.
\end{proof}

\begin{Lemma}\label{l_locInt_hlfpt}
Let $x\in \I$ and $A,B$ be two apartments containing $x$. We assume that there exists a   neighborhood $\Omega$ of $x$ in $A$ and a half-apartment $D$ of $A$ such that $A\cap B\supset D\cap \Omega$. Then either $x\in \In(A\cap B)$, or  there exists an apartment $A'$ such that $A\cap A'$ is a half-apartment and  $x\in \In (B\cap A')$. 
\end{Lemma}

\begin{proof}
Reducing $\Omega$ if necessary, there exists an apartment isomorphism $f:A\rightarrow B$ fixing $\Omega\cap A\cap B$, by Lemma~\ref{l_loc_iso}. Let $D_B=f(D)$. Let $P$ be a panel based at $x$ and contained in the wall $M$ of $D$. Then $P\subset A\cap B$. Let $C$ be the alcove of $B$ dominating $P$ and not contained in $D_B$. By the local sundial configuration (Proposition~\ref{p_loc_SC}), there exists an apartment $A'$ containing $D$ and $C$. Then $A'\cap B$ contains $(\Omega\cap D_B)\cup C$. Therefore by Lemma~\ref{l_lc_hf_apt_alc}, $x\in \In(B\cap A')$. If $C\subset A$, then we can take $A=A'$. Otherwise, we have $A\cap A'=D$,  by Lemma~\ref{l_half_apt}.
\end{proof}

A \textbf{polyhedron}\index{polyhedron} $\cP$ of $\A$ is a non-empty finite intersection of closed half-spaces of $\A$. We refer to  \cite[Definition 1.3]{bruns2009polytopes} for   a definition of the \textbf{faces}\index{face of a polyhedron} of $\cP$. A \textbf{facet}\index{facet of a polyhedron} of $\cP$ is  a  proper face of $\cP$ of maximal dimension. By \cite[Theorem 1.10]{bruns2009polytopes}, $\cP$ has finitely many faces.

\begin{Lemma}\label{l_U_Pi_cnt_pt}
Let $k\in \N$ and $P_1,\ldots,P_k$ be closed subsets of a non-empty metric space $(E,d)$ such that $\bigcup_{i=1}^k P_i=E$. Let $y\in E$ and $J_y=\{i\in \llbracket 1,k\rrbracket\mid y\in P_i\}$. Then $\bigcup_{i\in J_y} P_i$ contains $y$ in its interior.
\end{Lemma}

\begin{proof}
Let $J\subset \llbracket 1,k\rrbracket$ be such that $ \bigcup_{j\in J} P_j$ does not contain $y$ in its interior. Then there exists a sequence $(y_m)\in (E\setminus \bigcup_{i\in J} P_i)^\N$ converging to $y$. Let $i\in \llbracket 1,k\rrbracket$ be such that $\{n\in \N\mid y_n\in P_i\}$ is infinite. Then $i\notin J$.  Up to extracting a subsequence of $(y_n)$, we can assume that $y_n\in P_i$, for all $n\in \N$. But then $y\in P_i$, which proves that $i\in J_y$. As $i\notin J$, we obtain that $J\neq J_y$, which proves the result.
\end{proof}

Let $B$ be an apartment,  $\fQ$ be a sector-germ of $\A$ and $\rho=\rho_{\A,\fQ}:\I\rightarrow \A$ be the retraction onto $\A$ centered at $\fQ$. Using Proposition~\ref{p_splt_apt}, we have a decomposition $(P_i,B_i,f_i)_{i\in \llbracket 1,k\rrbracket}$ of $B$ with respect to $\fQ$¨. Let $\cE=\{P_1,\ldots,P_k\}$. If $P\in \cE$, then  $P$ is a finite intersection of half-spaces and thus it is a polyhedron. We define a wall of $P$ to be the support of one facet of $P$. Let $\cM$ be the set of all the walls of the elements of $\cE$. We set $\sS=\bigcup_{M\in \cM}M$ and $\sS^2=\bigcup_{M,M'\in \cM\mid M\neq M'} M\cap M'$. 

\begin{Lemma}\label{l_ptp_loc_hsp}
Let $x\in \sS\setminus \sS^2$. We assume that $x\notin \bigcup_{i=1}^k \In(P_i)$.  Then there exists a neighborhood $\Omega$ of $x$, a wall $M$ of $\A$ and $i_1,i_2\in \llbracket 1,k\rrbracket$ such that if $D_1,D_2$ are the two half-apartments delimited by $M$, we have $\Omega\cap P_{i_j}=\Omega\cap D_j$ for $j\in \{1,2\}$. 
\end{Lemma}

\begin{proof}
 As $x\in \sS\setminus \sS^2$, there exists a unique $M\in \cM$ such that $x\in M$. Let $D_1,D_2$ be the two half-apartments of $B$ delimited by $M$. Let $J=\{i\in \llbracket 1,k\rrbracket \mid x\in P_i\}$. Then by Lemma~\ref{l_U_Pi_cnt_pt}, $\bigcup_{i\in J} P_i$ contains $x$ in its interior. Choose a convex open neighborhood  $\Omega$ of $x$ contained in $\bigcup_{i\in J} P_i$. Reducing $\Omega$ if necessary, we may assume that $\Omega$ does not meet any wall of $\cM\setminus \{M\}$. For $j\in \{1,2\}$, choose $y_j\in \Omega\cap \In(D_j)$. Choose $i_j\in \llbracket 1,k\rrbracket$ such that $y_j\in P_{i_j}$. Then as one of the walls of $P_{i_j}$ is $M$, we have $P_{i_j}\subset D_{j}$ and hence $\Omega\cap P_{i_j}\subset \Omega\cap D_j$. Let us prove that these sets are equal. Let $z\in \Omega\cap D_j$. Assume by contradiction that $z\notin  P_{i_j}$. Then $]z,y_j[\cap \partial P_{i_j}=\{a\}$, for some $a\in \Omega$. By choice of $\Omega$, we have $a\in M$ and hence $[z,y_{j}]\subset M$: a contradiction. Therefore $\Omega\cap D_j=\Omega\cap P_{i_j}$.
\end{proof}

We say that ``$\fs:]0^-,0^+[\rightarrow B$ is a local line-germ'' to  mean that there exists $\epsilon>0$ such that $\fs:[-\epsilon,\epsilon]\rightarrow B$ is an affine map. We denote by   $\Phi_+(\fQ)$  the set of $\alpha\in \Phi$ such that $\alpha(x+\fQ)$ contains a neighborhood of $+\infty$, for any $x\in \A$. 

\begin{Lemma}\label{l_H_pth_1_stp}
Let $x\in \sS\setminus \sS^2$. We assume that $x\notin \bigcup_{i=1}^k \In(P_i)$.  Let $\fs:]0^-,0^+[\rightarrow B$ be a local line-germ such that $\fs(0)=x$. Let $\gamma=\rho\circ \fs$. Then  either $\gamma$ does not fold at $0$ or there exists $\alpha\in \Phi_{+}(\fQ)$ such that $\alpha(\gamma_-'(0))>0$ and $\gamma'_+(0)=r_\alpha.\gamma'_-(0)$. In particular, if $\fQ=\fQ_{+\infty}$ (resp. $\fQ=\fQ_{-\infty}$), we have $\gamma'_+(0)\leq_{Q^\vee_\R} \gamma'_-(0)$ (resp. $\gamma'_+(0)\geq_{Q^\vee_\R} \gamma'_-(0)$). 
\end{Lemma}

\begin{proof}
If $\gamma$ does not fold, there is nothing to prove so we assume that $\gamma$ folds at $0$.  Let $H$ be the wall as in Lemma~\ref{l_ptp_loc_hsp}. Let $D_1$ and $D_2$ be the two half-apartments delimited by $H$. Let $\Omega$ be a neighborhood of $x$ such that $P_1\cap \Omega=D_1\cap \Omega$ and $P_2\cap \Omega=D_2\cap \Omega$. Let $B_1$ be an apartment containing $P_1\cup \fQ$.

 If $B_1\cap B$ contained $x$ in its interior, then $\gamma$ would not fold so $x\notin \In (B_1\cap B)$. By Lemma~\ref{l_locInt_hlfpt}, there exists an apartment  $B_3$ such that $B_1\cap B_3$ is a half-apartment and $x\in \In(B_1\cap B_3)$.

As $\gamma$  folds, $B_3$ does not contain $\fQ$.  Let $B_2=(B_3\setminus B_1)\cup (B_1\setminus B_3)\cup M$, where $M$ is the wall of $B_1\cap B_3$. Then $B_2$ is an apartment by the exchange condition (Proposition~\ref{p_EC_MT}). Moreover $B_2$ contains $\fQ$ since $B_3$ does not. We have the following (non-commutative) diagram:

\[\xymatrix{ B_3\ar[d]^{f_3}\ar[r]^{f_2} & B_2\ar[d]^{g_2} \ar[ld]^{h} \\ B_1\ar[r]^{g_1}&\A}, \] where  $f_3$ and $f_2$ fix $B_3\cap B_1$ and $B_3\cap B_2$ respectively, and the other maps fix $\fQ$. We have $\rho\circ \fs|_{[0,0^+[}=g_1\circ f_3\circ \fs|_{[0,0^+[}$ and $\rho\circ \fs|_{]0^-,0]}=g_2\circ f_2\circ \fs|_{]0^-,0]}$. We have $g_2=g_1\circ h$ (since all these isomorphisms fix $\fQ$) and $f_3=r'\circ h\circ f_2$, by Lemma~\ref{l_reflction}, where $r'$ is the reflection on $B_1$ fixing $B_1\cap B_2\cap B_3$ and contained in  the affine Weyl group of $B$. We thus have the following commutative diagram: \[\xymatrix{ B_3\ar[d]^{f_3}\ar[r]^{f_2} & B_2\ar[d]^{h}  \\ B_1\ar[r]^{r'}\ar[d]^{g_1} &B_1 \ar[d]^{g_1}\\
\A \ar[r]^{r} & \A }, \] where $r:\A\rightarrow \A$ is the reflection with respect to $g_1(B_1\cap B_2 \cap B_3)$ contained the affine Weyl group. 

Let $\alpha':B_1\rightarrow \R$ be the affine root such that $\alpha'(B_1\cap B_2\cap B_3)=\{0\}$ and $\alpha'(\fQ)$ is the filter on $\R$ consisting of the sets containing a neighborhood of $+\infty$. We have $\alpha'(\fs(t))=\alpha'(x)+t\vec{\alpha'}(fs'_-(0))$, for $t\in \R_{<0}$ close enough to $0$, where $\vec{\alpha'}$ is the linear part of $\alpha'$. Moreover $\alpha'(\fs(t))<0$, for $t\in \R_{<0}$ small enough. Lemma follows.
\end{proof}

\subsection{Image of a line segment by a retraction}

In this subsection, we conclude the proof of Theorem~\ref{t_H_pth}, which asserts that the image of a segment by a retraction is a Hecke path.

\begin{Lemma}\label{l_lcl_H_pth}
Let $\fQ$ be a sector-germ of $\A$ and $\rho=\rho_{\A,\fQ}:\I\rightarrow \A$. Let $\fs:]0^-,0^+[\rightarrow B$ be a local line-germ and $a=\rho(\fs(0))$. Then $\pi:=\rho\circ \fs$ is a local Hecke path, i.e there exists a $(W^v_a,\fQ)$-chain from $\pi'_-(0)$ to $\pi'_+(0)$. 
\end{Lemma}

\begin{proof}

\begin{figure}[h]
\centering
\includegraphics[scale=0.3]{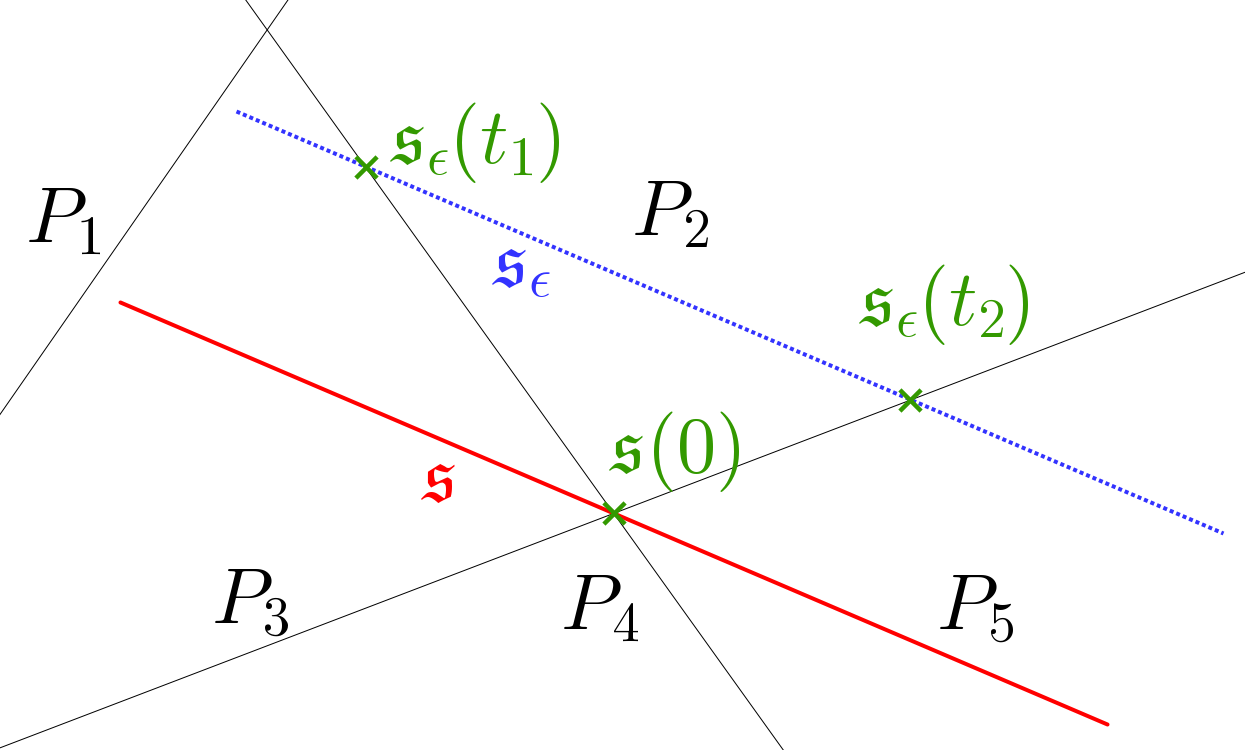}
\caption{Proof of Lemma~\ref{l_lcl_H_pth}. The point $x=\fs(0)$ belongs to $P_2\cap P_3\cap P_4\cap P_5$. We consider the translate $\fs_{\epsilon}$ of $\fs$ by a vector $\epsilon \vec{u}$, for $\epsilon\in \R_{>0}$ and some well-chosen vector $\vec{u}\in \A$.  We choose $\epsilon$ and $\vec{u}$ so that for all $t\in \R$,  $\fs_{\epsilon}(t)$ is well-defined and  $\fs_\epsilon(t)$ belongs to at most two $P_i$. We can then apply Lemma~\ref{l_H_pth_1_stp}.}\label{f_H_path}
\end{figure}

The proof is illustrated by Figure~\ref{f_H_path}. Let $x=\fs(0)$. By Proposition~\ref{p_splt_apt}, there exists a decomposition $(P_i,A_i,f_i)_{i\in \llbracket 1,n\rrbracket}$ of $B$ with respect to $\fQ$. For all $i\in \llbracket 1,n\rrbracket$, $P_i$ is polyhedron delimited by walls. If $x$ belongs to $\In(P_i)$, for some $i\in \llbracket 1,n\rrbracket$, then $\pi$ does not fold, and there is nothing to prove. Thus we assume that $x\notin \bigcup_{i=1}^n \In(P_i)$.

If $i\in \llbracket 1,n\rrbracket$, we define the  walls of $P_i$ to be the supports of the facets of the $P_i$ (the facets of $P_i$ are its faces of codimension $1$ in $B$). Let $\cM$ be the set of all the  walls of the $P_i$, for $i\in \llbracket 1,k\rrbracket$. Let $\cM_x=\{M\in \cM\mid x\in M\}$.   Let $\Omega$ be a convex open neighborhood of $x$ in $B$ such that for all $M\in \cM\setminus \cM_x$, we have $\Omega\cap M=\emptyset$. Let $m=|\cM_x|$.  For each $M\in \cM_x$, choose an affine root $\beta_M$ on $B$ such that $M=\beta_{M}^{-1}(\{0\})$. Let $V$ be a segment of $\R$ containing $0$ in its interior on which $\fs$ is defined. We assume that $\fs(V)\subset \Omega$, which is possible, up to reducing $V$.

Let $\vec{B}$ be the direction of $B$ (for example one can take $\vec{B}=\A$). If $\beta:B\rightarrow \R$ is an affine map, we denote by $\vec{\beta}$ its linear part. Let $\vec{\lambda}=\fs'(0)\in \vec{B}$. Let \[\vec{u}\in \vec{B}\setminus \left(\bigcup_{M\in \cM_x}\ker \vec{\beta}_M\cup  \bigcup_{M,M'\in \cM_x\mid M\neq M'} \ker(\vec{\beta}_M-\frac{\vec{\beta}_{M'}(\vec{\lambda})}{\vec{\beta}_{M}(\vec{\lambda})}\vec{\beta}_{M'})\right).\]

Let $t_-=\min(V)$ and $t_+=\max(V)$. Let $x_-=\fs(t_-)$ and $x_+=\fs(t_+)$. For $\epsilon\in \R_{>0}$, we set $\fs_{\epsilon}=\fs+\epsilon\vec{u}$.

Let $t\in V$ and $M\in \cM_x$.   We have:
\begin{equation}\label{e_beta_M}
 \beta_M(\fs_\epsilon(t))=\beta_M(x+t\epsilon \vec{u})=\beta_M(x)+t\epsilon \vec{\beta}_M(\vec{u}).
\end{equation}  Let $\eta\in \{-,+\}$.  Then for $\epsilon$ small enough, we have: \begin{equation}\label{e_non_zero}
\beta_{M}(\fs_\epsilon(t_\eta))\neq 0, \forall \eta\in \{-,+\}, \forall M\in \cM_x.  
 \end{equation}
 
 Let $\cM'_x=\{M\in \cM_x\mid \exists t\in [t_-,t_+], \beta_M(\fs_\epsilon(t))=0\}$. By \eqref{e_beta_M}, for every $M\in \cM'_x$, there exists a unique $t_M\in ]t_-,t_+[$ such that $\beta_M(\fs_{\epsilon}(t_M))=0$. More precisely, we have 
 $t_M=t_M(\epsilon)=-\epsilon \frac{\vec{\beta}_M(\vec{u})}{\vec{\beta}_M(\vec{\lambda})}$. By our assumption on  $\vec{u}$, the family $(t_M)_{M\in \cM'_x}$ is injective. Let $k=|\cM'_x|$. We write $\cM'_x=\{M_1,\ldots,M_k\}$ in such a way that if we set $t_i=t_{M_i}$ for $i\in \llbracket 1,k\rrbracket$, we have  $t_{1}<t_{2}<\ldots < t_{k}$.

Up to renumbering, we can assume that for some $n_x\in \llbracket 1,n\rrbracket$, we have  $\{i\in \llbracket 1,n\rrbracket\mid x\in P_i\}=\llbracket 1,n_x\rrbracket$. Set $t_0=t_-$ and $t_{k+1}=t_+$. Then for all $i\in \llbracket 1,k\rrbracket$, there exists $j_i(\epsilon)\in \llbracket 1,n_x\rrbracket$ such that: 

 \begin{equation}\label{e_fs_cnted_P_i}
\fs_\epsilon([t_i,t_{i+1}])\subset P_{j_i(\epsilon)}.
\end{equation} 

  As $\llbracket 1,n_x\rrbracket$ is finite, we can find a sequence $(\epsilon_m)\in (\R_{>0})^\N$ converging to $0$ and such that $(j_i(\epsilon_m))_{m\in \N}$ is constant, for all $i\in \llbracket 1,k\rrbracket$. We write $j_i$ the value $j_i(\epsilon_m)$, for $i\in \llbracket 1,k\rrbracket$ and $m\in \N$.

 For $i\in \llbracket 1,n\rrbracket$, we set $g_i=h_i\circ f_i:B\rightarrow \A$, where $h_i:A_i\rightarrow \A$ is the apartment isomorphism fixing $A_i\cap \A$. For $i\in \llbracket 1,n\rrbracket$ and $y\in P_i$, we have $\rho(y)=g_i(y)$. 

 We set $g=g_{j_1}$.  Then for all $i\in \llbracket 1,k \rrbracket$, there exists $w_i\in W^a$ such that $g_{j_i(\epsilon)}=w_i\circ g$.  Let $\pi_\epsilon=\rho\circ \fs_\epsilon$.  Let $i\in \llbracket 1,k\rrbracket$. Then by Lemma~\ref{l_H_pth_1_stp}, we have the following two possibilities.  Either $\pi'_{\epsilon,-}(t_i)=\pi'_{\epsilon,+}(t_i)$ or  there exists $\gamma_i\in \Phi_+(\fQ)$ such that $\gamma_i(\pi'_{\epsilon,-}(t_i))>0$ and $\pi'_{\epsilon,+}(t_i)=r_{\gamma_i}(\pi'_{\epsilon,+}(t_i))$. Let $i\in \llbracket 0,k\rrbracket$. By \eqref{e_fs_cnted_P_i}, $\pi_\epsilon$ is a segment on $[t_i,t_{i+1}]$. In particular, $\pi_{\epsilon,-}'(t_{i+1})=\pi_{\epsilon,+}'(t_i)$.  Thus we proved the existence of a $(W^v_a,\fQ)$-chain from $\pi'_{\epsilon,+}(t_-)$ and $\pi'_{\epsilon,-}(t_+)$. 
 
 By \eqref{e_fs_cnted_P_i}, $\fs_{\epsilon}(t_+)\in P_{j_k}$. Therefore \eqref{e_non_zero} implies that $\fs_{\epsilon}(t_+)\in \In(P_{j_k(\epsilon)})$. Similarly, $\fs_{\epsilon}(t_-)\in \In( P_{j_1(\epsilon)})$.   Therefore $\pi(t)=g(\fs(t))$, for $t\in [t_-,0]$ and $\pi(t)=g_{j_+}(\fs(t))$ for $t\in [0,t_+]$. We deduce that $\pi$ and $\pi_\epsilon$ are translates of each other on $[t_-,t_-^+]$ and on $[t_+^-,t_+]$. In particular $\pi'_-(0)=\pi'_+(t_-)=\pi'_{\epsilon,+}(t_-)$ and $\pi'_+(0)=\pi'_-(t_+)=\pi'_{\epsilon,-}(t_+)$, which proves the existence of a $(W^v_a,\fQ)$-chain from $\pi'_-(0)$ to $\pi'_+(0)$.

\end{proof}

The following theorem generalizes \cite[Theorem 6.2]{gaussent2008kac}. In comparison, we drop the preorder assumption on $\fs$. Our proof is also different and relies on the sundial configuration.

\begin{Theorem}\label{t_H_pth}
Let $\fs:[0,1]\rightarrow \I$ be a segment (resp. $\fs:\R_{\geq 0}\rightarrow \I$ be a ray, resp. $\fs:\R\rightarrow \I$  be a  line). Let $\fQ$ be a sector-germ at infinity of $\A$. Let $\rho=\rho_{\A,\fQ}$ be the retraction onto $\A$ centered at $\fQ$. Then $\gamma:=\rho\circ \fs$   is a Hecke path (resp. a Hecke ray, resp. a Hecke line) of shape $\gamma'_+(0)$ with respect to $\fQ$. 
\end{Theorem}

\begin{proof}
By Corollary~\ref{c_splt_seg}, $\gamma$  are  piecewise affine paths, with finitely many folding times. By Lemma~\ref{l_lcl_H_pth}, we deduce the result. 
\end{proof}

\section{A characterization of the points of $\A$ via retractions}\label{s_chrctz_apt}

Recall that $\I$ is a weak masure. In this section,  we make the additional assumption that $(\alpha_i^\vee)_{i\in I}$ is positively free. Under this assumption, we give a characterization of the points of $\A$: we prove that an element $x\in \I$ belongs to $\A$ if and only if $\rho_{+\infty}(x)=\rho_{-\infty}(x)=x$ (see Corollary~\ref{c_char_aprtmt}). 

We follow \cite{hebert2017gindikin}. We use Hecke paths to majorize $\htt(\rho_{+\infty}(x)-\rho_{-\infty}(x))$ by the ``distance'' between $x$ and $\A$, for $x\in \I$ (see Proposition~\ref{p_bnd_tnu}). Note that these methods also enable to prove the Gindikin-Karpelevich finiteness (i.e if $\I$ is semi-discrete with finite thickness, then $\rho_{+\infty}^{-1}(\{\lambda\})\cap \rho_{-\infty}^{-1}(\{\mu\})$ is finite, for every $\lambda,\mu\in Y$), which is the main goal of \cite{hebert2017gindikin}.

\begin{Definition}\label{d_pseudo_hgt}
Let $\htt':Q^\vee\otimes \R \rightarrow \R$ be a linear map. We say that $\htt'$ is a \textbf{pseudo-height}\index{pseudo-height} if $\htt'(\alpha_i^\vee)>0$, for all $i\in I$. 
\end{Definition}

By Lemma~\ref{l_char_pf},  $(\alpha_i^\vee)_{i\in I}$ is positively free if and only if there exits a pseudo-height. Note that by \cite[Theorem 5.6 c)]{kac1994infinite}, if the Kac--Moody matrix $A$ is indecomposable of indefinite type, then we can choose a positive imaginary root as a pseudo-height. 

\begin{Notation}
We fix a pseudo height $\htt'$\index[notation]{h@$\htt'$} on $Q^\vee\otimes \R$ such that $\min_{i\in I} \htt'(\alpha_i^\vee)=1$. If $(\alpha_i^\vee)_{i\in I}$ is free, we choose $\htt'=\htt$. 
\end{Notation}

\begin{Lemma}\label{l_H_pth_htt}
Let $\dd\in \R_{\geq 0}$, $\mu\in \A$,  $a\in \mathbb{A}$ and  $\nu\in Y^{++}:=Y\cap \overline{C^v_f}$. We assume the existence of a Hecke path $\gamma$ of shape $\dd\nu$ with respect to $-C^v_f$ from $a$ to $a+\dd\nu-\mu$. Then: \begin{enumerate}
\item we have $\mu\in Q^\vee_{\R_+}$. Consequently $\htt'(\mu)$ is well-defined.

\item  If $\dd>\htt'(\mu)$, there exists $t'\in [0,1[$ such that $\gamma$ is differentiable on $]t',1]$ and $\gamma'_{|]t',1]}=d\nu$. Furthermore, let $t^*$ be the smallest $t' \in [0,1]$ having this property, then $t^*\leq \frac{\htt'(\mu)}{d}$.

\end{enumerate} 

\end{Lemma}

\begin{proof} The main idea of 2) is to use the fact that during the time when $\gamma'(t')\neq \dd\nu$, we have $\gamma'(t')=\dd\nu-\dd\lambda(t')$ with $\lambda(t')\in Q^\vee_+\backslash \{0\}$. Hence for $\dd$ large, $\gamma$ decreases quickly for the $Q^\vee$ order, but it cannot decrease too much because $\mu$ is fixed.

Let $n\in \N$ and $0=t_0<\ldots<t_n=1$ be a subdivision of $[0,1]$ such that for all $i\in \llbracket 0,n-1\rrbracket $, $\gamma_{|[t_i,t_{i+1}]}$ is a line segment. For $i\in \llbracket 0,n-1 \rrbracket$, let $w_i\in W^v$ be such that $\gamma'_{|]t_i,t_{i+1}[}=w_i.\dd\nu$ (which exists by Proposition~\ref{p_H_paths}). If $w_i.\nu=\nu$, we choose $w_i=1$.

For $i\in \llbracket 0,n-1\rrbracket$, we have $w_i.\nu =\nu -\lambda_i$, for some  $\lambda_i \in Q^\vee_+$, according to \eqref{e_GR2.4}. Moreover if $w_i\neq 1$, then  $\lambda_i\neq 0$. We have: \begin{equation}\label{e_diff_gamma}
\gamma(1)-\gamma(0)=\dd\nu-\sum_{i=0, w_i\neq 1}^{n-1}(t_{i+1}-t_i)\dd\lambda_i=\dd\nu -\mu
\end{equation}

and we deduce 1).

Assume  now $\dd>\htt'(\mu)$. Let us prove the existence of  $i\in \llbracket 0,n-1\rrbracket$ such that $w_i=1$. Let $i\in \llbracket 0,n-1\rrbracket$. For all $i$ such that $w_i\neq 1$, we have $\htt'(\lambda_i)\geq 1$. Hence: \[\sum_{i=0, w_i\neq 1}^{n-1}(t_{i+1}-t_i)\leq \frac{\htt'(\mu)}{\dd}<1=\sum_{i=0}^{n-1}(t_{i+1}-t_i).\] Thus there exists $i\in \llbracket 0,n-1\rrbracket$ such that $w_i=1$. But then by Remark~\ref{r_H_pth}, we obtain that $w_j=1$ for all $j\in \llbracket i,n-1\rrbracket$. 

 This shows the existence of $t'\in [0,1[$ such that $\gamma$ is a line-segment on $[t',1]$. As $\gamma$ has finitely many folding times,   $t^*$ is well-defined and there exists $i^*\in \llbracket 0,n-1\rrbracket$ such that $t^*=t_{i^*}$. Then \eqref{e_diff_gamma} becomes $\sum_{i=0}^{i^*-1} (t_{i+1}-t_i)\dd\lambda_i=\mu$. As $\htt'(\lambda_i)\geq 1$ for all $i\in \llbracket 0,i^*-1\rrbracket$, we deduce $\htt'(\mu)\leq \sum_{i=0}^{i^*-1}t_{i+1}-t_i=t^*$, which proves 2). 
\end{proof}

Let $\fQ$ be a sector-germ of $\A$ and  $\fq=0+\fQ$. Write $\fq=\epsilon w.C^v_f$, with $\epsilon\in \{-,+\}$ and $w\in W^v$. The expression ``for $\nu\in \fq$ \textbf{sufficiently dominant}''\index{sufficiently dominant} means that $\epsilon w.\alpha_i(\nu)$ can be as large as we want, for any $i\in I$. 

\begin{figure}[h]
\centering
\includegraphics[scale=0.3]{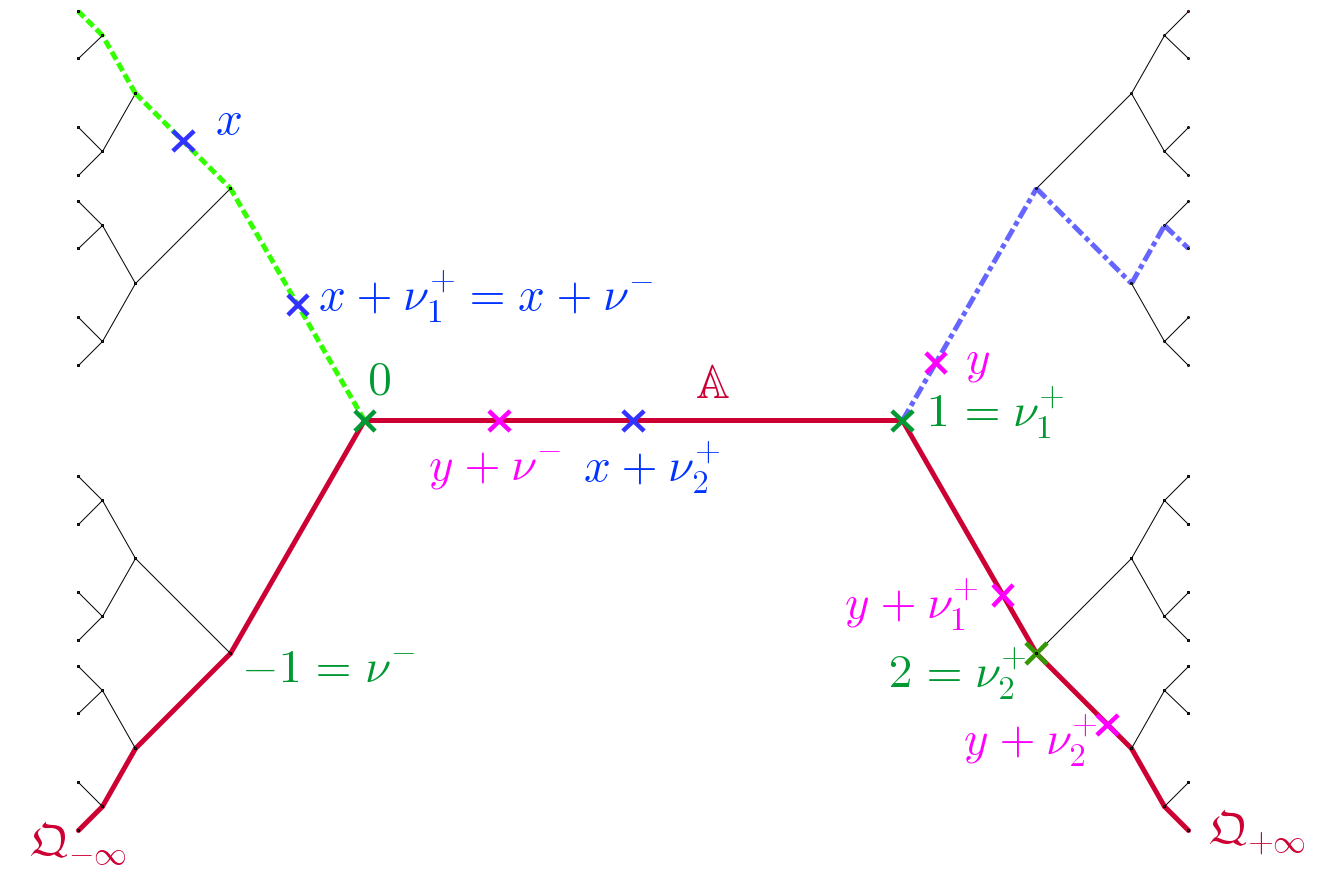}
\caption{We drop the indices $\fQ_{+\infty}$ and $\fQ_{-\infty}$ in the figure. Let  $\delta_x$ be the green dotted ray and $\delta_y$ be the purple dotted ray. The point $x$ is contained in the apartment $A_x^-=]\fQ_{-\infty},0]\cup \delta_x$ (resp. $A_x^+=[0,\fQ_{+\infty}[\cup \delta_x$), which contains $\fQ_{-\infty}$ (resp. $\fQ_{+\infty}$). Then $x+_{\fQ_{-\infty}}\nu^-$ (resp. $x+_{\fQ_{+\infty}}\nu_1^+$) is the translate of $x$ by  $1$ in the direction $\fQ_{-\infty}$ (resp. $\fQ_{+\infty}$). The apartment $A_y^+=\delta_y\cup [1,\fQ_{+\infty}[$ (resp. $A_y^-=\delta_y\cup ]\fQ_{-\infty},1]$) contains $y$ and $\fQ_{+\infty}$  and thus we obtain $y+_{\fQ_{+\infty}}\nu_1^+$,  $y+_{\fQ_{+\infty}}\nu_2^+$ and $y+_{\fQ_{-\infty}}\nu^-$ by translating by $1$, $2$ and $1$ in the direction $\fQ_{+\infty}$, $\fQ_{+\infty}$ and $\fQ_{-\infty}$ respectively.}\label{f_translations}
\end{figure}

If $x\in \I$, we define $x+\overline{\fQ}$ as the closure of $x+\fQ$ in any apartment containing $x$ and $\fQ$. This is well-defined by Lemma~\ref{l_iso_sect}.  Let $\rho=\rho_{\fQ,\A}$.  For $x\in \I$ and $\nu\in \overline{\fq}$, we now define $x+_{\fQ}\nu\in x+\overline{\fQ}$ as the unique element $y$ of $x+\overline{\fQ}$ such that $\rho(y)=x+_{\fQ}\nu$. This is well-defined by Lemma~\ref{l_def_x+tnu} 1). Heuristically, $x+_{\fQ}\nu$ is the translate of $x$ by the vector $\nu$, see Figure~\ref{f_translations}.

\begin{Lemma}\label{l_def_x+tnu}
\begin{enumerate}
\item Let $x\in \I$ and $\nu\in \overline{\fq}$. Then there exists a unique $y\in x+\overline{\fQ}$ such that $\rho(y)=\rho(x)+\nu$. More precisely, if $A$ is an apartment containing $x$ and $\fQ$ and $f:A\rightarrow \A$ is the apartment isomorphism fixing $\fQ$, then $x+_{\fQ}\nu=f^{-1}(\rho(x)+\nu)=f^{-1}(f(x)+\nu)$.  

\item Let $x\in \I$ and $\nu,\nu'\in \overline{\fq}$. Then $x+_{\fQ}(\nu+\nu')=(x+_{\fQ}\nu)+_{\fQ}\nu'$.

\item Let $x\in \A$ and $\nu \in \overline{\fq}$. Then $x+_{\fQ}\nu$ coincides with $x+\nu$, when we regard $\A$ as a vector space. 

\item Let $x\in \I$. Then for $\nu\in \fq$ sufficiently dominant, $x+_{\fQ}\nu=\rho(x)+\nu\in \A$. 
\end{enumerate}
\end{Lemma}

\begin{proof}
1) Let $A$ be an apartment containing $x$ and $\fQ$.   Let $f:A\rightarrow \A$ be the apartment isomorphism fixing $A\cap \A$, which exists by Lemma~\ref{l_iso_sect}.  Let $z\in \overline{x+\fQ}$. Assume that $\rho(z)=\rho(x)+\nu$. Then as $z\in A$, we have $f(z)=\rho(x)+\nu$ and hence $z=f^{-1}(\rho(x)+\nu)$. 

Conversely, $\rho(f^{-1}(\rho(x)+\nu))=f(f^{-1}(\rho(x)+\nu))=\rho(x)+\nu$. Moreover $f(x)+\nu\in \rho(x)+\overline{\fQ}$ and thus $f^{-1}(f(x)+\nu)\in f^{-1}(\rho(x)+\overline{\fQ})=f^{-1}(f(x))+\overline{\fQ}=x+\overline{\fQ}$. 

2) Let $A$ be an apartment containing $x$ and $\fQ$. Let $f:A\rightarrow  \A$ be the apartment isomorphism fixing $\fQ$. Then $f(x+_{\fQ}(\nu+\nu'))=x+\nu+\nu'$ and $f((x+_{\fQ}\nu)+_{\fQ}\nu')=f(x+_{\fQ}\nu)+\nu'=x+\nu+\nu'$.

3) If $x\in \A$, then we can take $A=\A$ and then $f$ is the identity.

4) We use  the same notation as in 1). Then  $A$ contains a subsector $\fq'$ of $\fq$. For $\nu\in \fq$ sufficiently dominant, $f(x)+\nu\in \fq'$ and as $f$ fixes $A\cap \A$, we have  $x+_{\fQ}\nu=f^{-1}(f(x)+\nu))=f(x)+\nu=\rho(x)+\nu\in \A$. 
\end{proof}

\begin{Lemma}\label{l_df_tnu}
Let $x\in \I$ and $\nu\in C^v_f$. Then there exists $\dd_{\nu}=\dd_{\nu}(x)\in \R_{\geq 0}$\index[notation]{d@$d_\nu$} such that $\{t\in \R_{\geq 0}\mid x+_{\fQ_{+\infty}}t\nu\in \A\}=[\dd_{\nu},+\infty[$. 
\end{Lemma}

\begin{proof}
By Lemma~\ref{l_def_x+tnu} 4), there exists $t_0\in \R_{\geq 0}$ such that $x+_{\fQ_{+\infty}} t_0\nu\in \A$. By Lemma~\ref{l_def_x+tnu} 2), we have $x+_{\fQ_{+\infty}} t\nu \in \A$, for all $t\in [t_0,+\infty[$. Therefore if $E=\{t\in \R\mid x+_{\fQ_{+\infty}}\nu\in \A\}$ and $\dd_{\nu}=\inf E$, we have either $E=]\dd_{\nu},+\infty[$ or $E=[\dd_{\nu},+\infty[$.  Let $A$ be an apartment containing $x$ and $\fQ_{+\infty}$. By Lemma~\ref{l_iso_sect}, \[A\cap \A\supset\{x+_{\fQ_{\infty}}t\nu \mid t\in ]\dd_{\nu},+\infty[\}=\{\rho_{+\infty}(x)+t\nu\mid t\in ]\dd_{\nu},+\infty[\}:=\delta.\]

Let $f:A\rightarrow \A$ be the apartment isomorphism fixing $A\cap \A$, which exists by Lemma~\ref{l_iso_sect}. As $\delta$ is a preordered ray, the closure  $\overline{\delta}$ of $\delta$ in $A$ is contained in $\cl^\Delta(\delta)\subset A\cap \A$ (by~\ref{a_oco}), and hence it is fixed by $f$. Therefore $x+_{\fQ_{+\infty}}\dd_{\nu} \nu\in A$, which proves the lemma. 
\end{proof}

With the notation of Lemma~\ref{l_df_tnu}, $x+_{\fQ_{+\infty}} \dd_\nu(x) \nu=x+_{\fQ_{+\infty}} \dd_{\nu'}(x) \nu'$, for all $\nu'\in \R^*_+.\nu$. This is the unique element $z$ of $\A$ such that $(x+_{\fQ_{+\infty}}\R_{\geq 0}\nu)\cap \A=z+_{\fQ_{+\infty}}\R_{\geq 0} \nu=z+\R_{\geq 0} \nu$: this is a kind of projection of $x$ onto $\A$ along $\R_{>0}\nu$. We set: \begin{equation}\label{e_pr_nu}
\pr_\nu(x)=x+_{\fQ_{+\infty}}\dd_\nu(x) \nu\in \A.
\end{equation}\index[notation]{p@$\pr_{\nu}$} This is the projection of $x$ onto $\A$ along $\R_{>0}\nu$.

\begin{Remark}

Assume  that $\I$ is a tree. In figure~\ref{f_retractions}, we have $x+_{\fQ_{+\infty}} \dd_\nu(x) \nu=0$ and $y+_{\fQ_{+\infty}} \dd_\nu(y) \nu=1$,  for any $\nu\in [0,\fQ_{+\infty}[\subset \A$. 

 For all $x'\in \I$, we have: \[\rho_{+\infty}(x')=\underbrace{x'+_{\fQ_{+\infty}}\dd_\nu(x')\nu}_{\in \A\simeq \R} -\dd_\nu(x') \nu\text{ and }\rho_{-\infty}(x')=\underbrace{x'+_{\fQ_{-\infty}}\dd_\nu(x')\nu}_{\in \A\simeq \R} +\dd_\nu(x') \nu.\] Therefore: \[\rho_{+\infty}(x')-\rho_{-\infty}(x')=-2\dd_\nu(x') \nu:\] we obtain   $\rho_{+\infty}(x')-\rho_{-\infty}(x')$ as a function of the distance $\dd_\nu(x')\nu$ from  $x'$ to $\A$. In general (when $\I$ is no longer assumed to be a tree), we obtain that $\htt'(\rho_{+\infty}(x)-\rho_{-\infty}(x))$ is majorized by a function of $\dd_\nu(x)\nu$ (see Proposition~\ref{p_bnd_tnu}). 
\end{Remark}

\begin{Lemma}\label{l_diff_H_pth}
Let $x\in \I\setminus \A$, $\nu\in C^v_f$ and $\dd_{\nu}=\dd_{\nu}(x)$, with the notation of Lemma~\ref{l_df_tnu}. Let $\fs:[0,1]\rightarrow \I$ be defined by $\fs(t)=x+_{\fQ_{+\infty}} t\dd_{\nu}\nu$, for $t\in [0,1]$. Let $\gamma=\rho_{-\infty}\circ \fs$. Then $\gamma$ is a Hecke path of shape $\dd_\nu \nu$ with respect to $C^v_f$ such that  for all $t\in [0,1]$ such that $\gamma$ is differentiable at $t$, we have $\gamma'(t)\neq \dd_{\nu} \nu$.

\end{Lemma}

\begin{proof}
By Proposition~\ref{p_H_paths} and \eqref{e_GR2.4}, it suffices to prove that $\gamma'_-(1)\neq \dd_{\nu} \nu$. Suppose by contradiction that it is the case. By Proposition~\ref{p_mas_th_iwa}, there exists an apartment $A$ containing $\fs(]1^-,1])\cup \fQ_{-\infty}$. By Lemma~\ref{l_iso_sect}, there exists an apartment isomorphism  $f:A\rightarrow \A$ fixing $A\cap \A$.  Let $y=\pr_\nu(x)=x+_{\fQ_{+\infty}}\dd_\nu \nu\in \A$.  Then by assumption, \[\gamma(]1^-,1])=\gamma(1)+]0^-,0] \nu=y+]0^-,0]\nu\subset y-C^v_f.\] Moreover $A\cap \A$ contains $y$ and $\fQ_{-\infty}$ and thus it contains $y+\fQ_{-\infty}=y-C^v_f$. Therefore $f$ fixes $y-C^v_f$ and $\fs([1^-,1])=f^{-1}(\gamma([1^-,1]))=\gamma([1^-,1])\subset \A$. We reach a contradiction with the definition of $\dd_{\nu}$ and thus $\gamma'_-(1)\neq \dd_{\nu} \nu$. 
\end{proof}

\begin{Proposition}\label{p_bnd_tnu}
Assume that $\cS$ is positively cofree. Let $x\in \I$ and $\nu\in Y\cap C^v_f$. Let $\htt'$ be a pseudo-height such that $\min_{i\in I}\htt'(\alpha_i^\vee)=1$.  Then:

\begin{enumerate}
\item We have $\rho_{-\infty}(x)-\rho_{+\infty}(x)\in Q^\vee_{\R_{\geq 0}}=\sum_{i\in I}\R_{\geq 0} \alpha_i^\vee$. 

\item Let  $\nu\in Y\cap C^v_f$ and $d_\nu=\dd_{\nu}(x)$ be as in Lemma~\ref{l_df_tnu}. Then $\dd_{\nu}(x)\leq \htt'(\rho_{-\infty}(x)-\rho_{+\infty}(x))$. 
\end{enumerate}
\end{Proposition}

\begin{proof}
Let $y=\pr_\nu(x)=x+_{\fQ_{+\infty}} \dd_{\nu} \nu\in \A$. Let $\fs:[0,1]\rightarrow \I$ be defined by $\fs(t)=x+_{\fQ_{+\infty}} td_\nu \nu$, for $t\in [0,1]$. Then $\fs$ is a segment and by Theorem~\ref{t_H_pth}, $\gamma:=\rho_{-\infty}\circ \fs$ is a Hecke path of shape $\dd_{\nu} \nu$ with respect to $\fQ_{-\infty}$. We have $\gamma(0)=\rho_{-\infty}(x)$ and $\gamma(1)=\rho_{-\infty}(y)=y$ since $y\in \A$. By Lemma~\ref{l_def_x+tnu}, $y=\rho_{+\infty}(x)+\dd_{\nu} \nu=\rho_{-\infty}(x)+\dd_{\nu} \nu +(\rho_{+\infty}(x)-\rho_{-\infty}(x))$. By Lemma~\ref{l_diff_H_pth}, $\gamma'_-(1)\neq \dd_{\nu} \nu$ and  thus by Lemma~\ref{l_H_pth_htt} applied with $\mu=\rho_{-\infty}(x)-\rho_{+\infty}(x)$, we get the result. 
\end{proof}

\begin{Corollary}\label{c_char_aprtmt}
Assume that $\cS$ is positively cofree. We have $\A=\{x\in \I\mid \rho_{+\infty}(x)=\rho_{-\infty}(x)\}$. 
\end{Corollary}

\begin{proof}
If $x\in\A$, then $\rho_{-\infty}(x)=\rho_{+\infty}(x)=x$. Let $x\in \I$ be such that $\rho_{-\infty}(x)=\rho_{-\infty}(x)$. Then by Proposition~\ref{p_bnd_tnu}, we have $\dd_{\nu}(x)=0$, which proves that $x\in \A$. 
\end{proof}

\section{Convexity of the intersection of two apartments}\label{s_cvx_intrsctn}

Let $\I$ be a space satisfying \ref{a_ma1}. We define the following two  axioms, which are weak versions of \ref{a_ma2}.

\begin{enumerate}[label=\blue{(MA2')}]
\item\label{a_ma2'} Let $A$ and $B$ be two apartments. Then $A\cap B$ is convex  and there exists an apartment isomorphism $\phi:A\rightarrow B$ fixing $A\cap B$.\axiom{ma2'@\ref{a_ma2'}}
\end{enumerate}

\begin{enumerate}[label=\blue{(MA2'')}]
\item\label{a_ma2''} Let $A$ and $B$ be two apartments whose intersection has non-empty interior in $A$. Then $A\cap B$ is convex.\axiom{ma2''@\ref{a_ma2''}}
\end{enumerate}

In this section, we prove that if  $\I$ is a weak masure associated with a positively cofree Kac--Moody datum or an indecomposable Kac--Moody matrix of affine type, then $\I$ satisfies \ref{a_ma2'}.
 We proceed as follows:\begin{enumerate}
\item We prove that \ref{a_ma2''} implies \ref{a_ma2'}.

\item We prove \ref{a_ma2''}. We separate the positively free case and the indecomposable affine case. In the positively free case, we consider two apartments $A$ and $B$ such that $A\cap B$ has non-empty interior and  $x,y\in \In(\A\cap B)$. We want to prove that $[x,y]_{\A}=[x,y]_B$. We consider an affine parametrization $\fs:[0,1]\rightarrow [x,y]_B$ of $[x,y]_B$ and we study $\gamma:=\rho_{-\infty}\circ \fs$ and $\tilde{\gamma}:=\rho_{+\infty}\circ \fs$. Using the fact that these are Hecke paths, we obtain inequalities (with respect to $\leq_{Q^\vee_\R}$) satisfied by these paths. We deduce that these two paths are necessarily the line segment of $\A$ between $x$ and $y$. Using the characterization of $\A$ given by Corollary~\ref{c_char_aprtmt}, we deduce that $[x,y]_B\subset \A$. 
\end{enumerate}

\subsection{\ref{a_ma2''} implies \ref{a_ma2'}}

\begin{Lemma}\label{l_fix_nonem}
Let $B$ be an apartment and $\Omega\subset \A\cap B$. We assume that $\Omega$ is convex and has non-empty interior. Then there exists an apartment isomorphism $\psi:\A\rightarrow B$ fixing $\overline{\Omega}$. In particular, the interior of $\A\cap B$ is the same if we regard $\A\cap B$ as a subset of $\A$ or as a subset of $B$. 
\end{Lemma}

\begin{proof}
Let $x\in \mathring{\Omega}$. By Lemma~\ref{l_loc_iso},  for every $x\in \mathring{\Omega}$, there exists a unique apartment isomorphism $\psi_x:\A\rightarrow B$ fixing an open neighborhood $V_x$ of $x$ in $\mathring{\Omega}$.

Fix $x\in \mathring{\Omega}$ and set $\psi=\psi_x$. Let $\Omega_1=\{y\in \mathring{\Omega} \mid \psi_y=\psi\}$. Let $y\in \Omega_1$.  Then for all $z\in V_y$, we have $\psi_z=\psi_y=\psi$ and thus $V_y\subset \Omega_1$. Therefore $\Omega_1$ is open in $\mathring{\Omega}$.

 Let $(y_k)\in (\Omega_1)^\N$ be a sequence converging in $\mathring{\Omega}$ and $y=\lim y_k$. Then for $n$ large enough, $y_n\in V_y$, thus  $\psi_y$ fixes a neighborhood of $y_n$ and hence $\psi_y=\psi_{y_k}=\psi$. This proves that $\Omega_1$ is open and closed in $\mathring{\Omega}$ and thus $\Omega_1=\mathring{\Omega}$. For $y\in \mathring{\Omega}$, we have $\psi(y)=\psi_y(y)=y$, which proves that $\psi$ fixes $\mathring{\Omega}$.

Let $x\in \overline{\Omega}$. As $\Omega$ is convex, there exits $(x_n)\in \In(\Omega)^\N$ such that $(x_n)$ converges to $x$ in $\A$. Then $(\psi(x_n))$ converges to $\psi(x)$ in $B$, by definition of the topology on $B$. But by Lemma~\ref{l_clsd_inter}, $(x_n)=(\psi(x_n))$ converges to $x$ in $B$, so we have $\psi(x)=x$. Therefore $\psi$ fixes $\overline{\Omega}$, which proves the lemma.
\end{proof}

\begin{Lemma}\label{l_cvxty_inter_imp}
Let $\I$ be a weak masure in the sense of Definition~\ref{d_w_mas}. Assume that $\I$ satisfies \ref{a_ma2''}. Then  $\I$ satisfies \ref{a_ma2'}.
\end{Lemma}

\begin{proof}
Let $A$ and $B$ be two apartments. We assume that $\Omega:=A\cap B$ is non-empty. Let $a,b\in \Omega$. Let $C_a$ (resp. $C_b)$ be an alcove of $A$ (resp. of $B$) based at $a$ (resp. $b$). By Proposition~\ref{p_frndly_fac}, there exists an apartment $\tilde{A}$ containing $C_a$ and $C_b$. By assumption, $A\cap \tilde{A}$ and $\tilde{A}\cap B$ are convex.  By Lemma~\ref{l_fix_nonem}, there exist $\phi_A:A\rightarrow \tilde{A}$ and $\phi_B:\tilde{A}\rightarrow B$ fixing pointwise $A\cap \tilde{A}$ and $B\cap \tilde{A}$ respectively. Set $\phi_{a,b}=\phi_B\circ\phi_A$. Then $[a,b]_A=[a,b]_{\tilde{A}}=[a,b]_B$ and  $\phi_{a,b}$ fixes $[a,b]_A$.  Therefore $\Omega$ is convex.

 Let $H$ be the support of $\Omega$. If $E\subset H$, we denote by $\In_H(E)$ its interior in $H$. Let $\fQ_A$ be a sector-germ at infinity of $A$.  By Proposition~\ref{p_splt_apt}, there exists a decomposition $(P_i, B_i, f_i)_{i\in \llbracket 1,n\rrbracket }$ of $B$ with respect to $\fQ_A$.
 
 We have $\Omega=\bigcup_{i=1}^n \Omega\cap P_i$. As the $\Omega\cap P_i$ are convex and $\In_H(\Omega)\neq \emptyset$, there exists $i\in \llbracket 1,n\rrbracket$ such that $\In_H(\Omega\cap P_i)\neq \emptyset$ (otherwise the measure of $\Omega$ would be $0$ for any Lebesgue measure on $H$). Let $\phi_i=g_i\circ f_i$, where $g_i:B_i\rightarrow A$ is the apartment isomorphism fixing $A\cap B_i$. 

Let us prove that $\phi_i$ fixes $A\cap B$. Let $a\in \In_H(\Omega\cap P_i)$. Let $b\in \Omega\setminus \{a\}$. Let $f=\phi_i^{-1}\circ \phi_{a,b}:A\rightarrow A$. Then $\phi$ fixes a neighborhood of $a$ in $[a,b]$ and as it is an affine isomorphism, it fixes the line containing $\{a,b\}$. In particular $f(b)=b$ and hence $\phi_i(b)=\phi_{a,b}(b)=b$, which proves that $\phi_i$ fixes $\Omega$.
\end{proof}

\subsection{The ``positively free'' case}\label{ss_pos_free}

In this subsection, we prove that a weak masure satisfies \ref{a_ma2''} when  $(\alpha_i^\vee)$ is positively free (see Definition~\ref{d_pos_free}), following \cite{hebert2022new}.

\begin{Lemma}\label{l_ma2'_pf}
Assume that $(\alpha_i^\vee)_{i\in I}$ is positively free in $\A$. Then $\I$ satisfies \ref{a_ma2''}: if $A$ and $B$ are two apartments such that $A\cap B$ has non-empty interior, then $A\cap B$ is convex. 
\end{Lemma}

\begin{proof}
Using apartment isomorphisms, we can assume that $A=\A$.  Let $a\in \In(B\cap \A)$ and $b\in B\cap \A$. We begin by proving that $[a,b]_B\subset \A$.   Let $\fs:[0,1]\mapsto [a,b]_B$ be the affine parametrization of $[a,b]$ such that $\fs(0)=a$ and $\fs(1)=b$. Let $\gamma=\rho_{-\infty}\circ \fs$ and $\tilde{\gamma}=\rho_{+\infty}\circ \fs$. Then by Proposition~\ref{p_H_paths} (using Remark~\ref{r_H_pth}), $\gamma'_{+}(0)\leq_{Q^\vee_\R} b-a$ and $\tilde{\gamma}_{+}(0)\geq_{Q^\vee_\R} b-a$. Moreover $\fs([0,0^+[)\subset \A$. Therefore $\gamma|_{[0,0^+[}=\tilde{\gamma}|_{[0,0^+[}=\fs|_{[0,0^+[}$. Consequently: \[\gamma'_+(0)=\tilde{\gamma}'_+(0)\leq_{Q^\vee_\R} b-a\leq_{Q^\vee_\R} \tilde{\gamma}'_+(0)=\gamma'_+(0).\]
As $(\alpha_i^\vee)_{i\in I}$ is positively free, we deduce that $\gamma'_+(0)=\tilde{\gamma}_+(0)=b-a$. The path $\gamma$   is thus a line segment, since otherwise  we would have $b-a=\gamma(1)-\gamma(0)>_{Q^\vee_\R} b-a$ by Proposition~\ref{p_H_paths}. By symmetry, $\tilde{\gamma}$ is a line segment too and thus $\gamma=\tilde{\gamma}$. Therefore   $\rho_{-\infty}\circ\fs(t)=\rho_{+\infty}\circ \fs(t)$ for $t\in [0,1]$. Using Corollary~\ref{c_char_aprtmt}, we deduce that $\fs([0,1])\subset \A$ and that  $[a,b]_B\subset \A$. By symmetry, $[a,b]_{\A}\subset \A\cap B$. Therefore: \[\forall a\in \In(\A\cap B), \forall b\in \A\cap B, [a,b]_\A\subset \A\cap B.\]

Let us prove that $\In(\A\cap B)$ is dense in $\A\cap B$. Let $b\in \A\cap B$. Let $a\in \In(\A\cap B)$.  Let $\Omega'$ be an open  neighborhood of $a$ in $\A$ contained in $\A\cap B$.  If $n\in \N^*$, then $(1-\frac{1}{n})b+\frac{1}{n}\Omega'$ is a neighborhood of $(1-\frac{1}{n})b+\frac{1}{n}a$ contained in $\bigcup_{a'\in\Omega'} [a',b]_{\A}\subset  \A\cap B$. Therefore $b\in \overline{\In(B\cap \A)}$ and thus we proved: \[\A\cap B\subset \overline{\In(\A\cap B)}.\]

Let now $a,b\in \A\cap B$ and $c\in [a,b]_\A$. Write $c=(1-t)a+tb$, with $t\in [0,1]$. Let $(a_n),(b_n)\in (\In(\A\cap B))^{\N}$ such that $a_n\to a$ and $b_n\to b$. Then $c_n:=(1-t)a_n+tb_n\in \A\cap B$ for all $n\in \N$ and $c_n\to c$. As $\A\cap B$ is closed (by Lemma~\ref{l_clsd_inter}), we deduce that $c\in \A\cap B$, which proves the lemma. 
\end{proof}

\subsection{The indecomposable affine case}

We saw in Lemma~\ref{l_pos_free} that if $\cS$ is not positively free, then the underlying Kac--Moody matrix has an indecomposable affine  component.

Let $\I$ be a weak masure of type $\cS$. We assume that the Kac--Moody matrix $A$ to which $\cS$ is associated is indecomposable and of affine type. We give a  proof of \ref{a_ma2}   without assuming that the family $(\alpha_i^\vee)_{i\in I}$ is free, following  \cite[5.3]{hebert2020new}. Note that this case will not be used elsewhere since we will prove that the masure of a Kac--Moody group   satisfies \ref{a_ma2} by deducing it  from the fact that a weak masure whose Kac--Moody datum is positively free satisfies \ref{a_ma2}.

Recall that $\A_{\ines}=\bigcap_{i\in I}\ker(\alpha_i)$. By \cite[Proposition 5.8]{kac1994infinite}, if $\delta\in \Delta^{im}_+$ denotes the smallest positive imaginary root, then we have: \begin{equation}\label{e_T_cn_aff}
\cT=\delta^{-1}(\R_{>0})\sqcup \A_{\ines}.
\end{equation}

\begin{Lemma}\label{l_ma2'_aff}
Let $\I$ be a weak masure associated with an indecomposable Kac-Moody matrix of affine type. Then $\I$ satisfies \ref{a_ma2''}: for every two apartments $A$ and $B$ whose intersection has non-empty interior, the set $A\cap B$ is convex.
\end{Lemma}

\begin{proof}
Let $x\in \In(\A\cap B)$ and let $C$ be an alcove of $\A$ based at $x$. By~\ref{a_oco}, there exists $f:\A\rightarrow B$ fixing $C$. Let $\Omega\in C$ be a convex open subset of $\A\cap B$ fixed by $f$. Let $y\in \A\cap B$. By \eqref{e_T_cn_aff}, there exists $z\in \Omega$ such that $z\leq y$ or $z\geq y$. Then by~\ref{a_oco}, there exists $f':\A\rightarrow B$ fixing $[z,y]$. Then $(f')^{-1}\circ f:\A\rightarrow \A$  is an affine automorphism fixing a neighborhood of $z$ in $[z,y]$ and thus it fixes the line $L$ containing $x$ and $y$. Therefore  $(f')^{-1}(f(y))=y=f(y)$. Therefore $f$ fixes $\A\cap B$.

Take $a\in [x,y]$. Write $a=tx+(1-t)y$, with $t\in [0,1]$. Let $(z_n)\in \Omega^\N$ be such that $z_n\to x$ and $z_n\leq y$ or $z_n\geq y$ for all $n\in \N$. Then $t z_n+(1-t)y\in \A\cap B$ for all $n\in \N$ and thus $a\in \A\cap B$, by Lemma~\ref{l_clsd_inter}.  Therefore $\A\cap B$ is convex. 
\end{proof}

\begin{Remark}\label{r_decompo_I}
 In a previous version of this paper, we claimed that when $A$ is a decomposable Kac--Moody matrix, then a weak masure of type $A$ decomposes as a product of weak masures whose type is an indecomposable Kac--Moody matrix. We then deduced that any weak masure satisfies \ref{a_ma2}, combining the indecomposable  affine case and the positively cofree case. Actually, we do not know if our claim was correct. In the reductive case, the fact that an affine building associated with a decomposable Cartan matrix decomposes as a product of affine buildings is true (see \cite[Proposition 2.1]{parreau2000immeubles} or \cite[2.4.8.1 Theorem]{rousseau2023euclidean}) but requires some work, and the proofs given in these references do not obviously generalize to the case of weak masures. However  we  will prove that the masure of a Kac--Moody group satisfies~\ref{a_ma2} (without assuming  $(\alpha_i^\vee)$ to be positively free), using the fact that it can be obtained as a quotient of a masure whose family of simple coroots is free. 
\end{Remark}

\section{Encloseness properties of the intersection of two apartments}\label{s_encls_int}

We proved in Section~\ref{s_cvx_intrsctn} that a weak masure satisfies \ref{a_ma2'} in many cases (see Lemmas~\ref{l_ma2'_pf} and \ref{l_ma2'_aff}). 

In this section, we prove that a weak masure satisfies  \ref{a_ma2'} if and only if it satisfies \ref{a_ma2}. We then gather the results of this chapter to prove that if $\cS$ is positively cofree or is  associated with an indecomposable affine Kac--Moody matrix, then a weak masure of type $\cS$ satisfies \ref{a_ma2}.

In order to prove that \ref{a_ma2'} implies \ref{a_ma2}, we proceed as follows: we begin by proving that the intersection of two apartments can be written as a finite union of enclosed subsets, using the characterizations of $\A$ obtained in Corollary~\ref{c_char_apt_alc} or in  Corollary~\ref{c_char_aprtmt} (see Lemma~\ref{l_un_encl_gen}). We then prove that a convex set that can be written as a finite union of polyhedra is actually a polyhedra (see Lemma~\ref{l_un_enclsd_cvx1}), which enables to conclude.

\subsection{The intersection of two apartments is a finite union of enclosed subsets}

In this subsection, we prove that the intersection of two apartments $A$ and $B$ can be written as a finite union of enclosed subsets. We start with the case where $A$ and $B$ share a common sector and then use \ref{a_sc} to deduce the general case. For the first case, we provide two independent proofs based on the two characterizations of $\A$ that we obtained (Corollaries~\ref{c_char_apt_alc} and  \ref{c_char_aprtmt}). Note that the second proof is less general since it requires $\cS$ to be positively cofree.

\begin{Lemma}\label{l_U_encl_sct1}
Let $B$ be an apartment containing $\fQ_{+\infty}$. Let $x\in \Int(\A\cap B)$. Let $\fQ^-$ be the sector-germ of $B$ such that $\fQ^-_x$ is opposite $\fQ_{+\infty,x}$. Using Proposition~\ref{p_splt_apt}, we get a decomposition $(P_i,B_i,f_i)_{i\in \llbracket 1,n\rrbracket}$ with respect to $\fQ^-$. Let $J=\{j\in \llbracket 1,n\rrbracket\mid \Int(P_j)\cap \Int(\A\cap B)\neq \emptyset\}$. Then $\A\cap B=\bigcup_{j\in J} P_j$. 
\end{Lemma}

\begin{proof}
As $\fQ_x^-$ and $\fQ_{+\infty,x}$ are opposite, the sector-germs $\fQ^-$ and $\fQ_{+\infty}$ are opposite. Therefore if $j\in J$ and $y\in \In(P_j)\cap \Int(\A\cap B)$, then  $\fQ_y^-$ and $\fQ_{+\infty,y}$ are opposite. Let $z\in \Int(P_j)$. Then  $y+\fQ^-$  and $z+\fQ^-$  are parallel in $B_j$, therefore $\fQ^-_y$ and $\fQ_{z}^-$ are parallel in $B_j$. Using $f_j^{-1}$, we deduce that $\fQ^-_y$ and $\fQ_{z}^-$ are parallel in $\A$. As $\A$ contains $y,z$ and $\fQ_{+\infty}$, the alcoves $\fQ_{+\infty,y}$ and $\fQ_{+\infty,z}$ are parallel. Therefore $\fQ_{+\infty,z}\opp \fQ^{-}_z$ and thus $z\in \A\cap B$, by Corollary~\ref{c_char_apt_alc}. Therefore $\Int(P_j)\subset \A\cap B$. Using Lemma~\ref{l_clsd_inter}, we deduce that $P_j=\overline{\Int(P_j)}\subset \A\cap B$, which proves that $\bigcup_{j\in J} P_j\subset \A\cap B$. 

Conversely, let $z\in \A\cap B$.  Then $z+\fQ_{+\infty}\subset \A\cap B$. As  the $P_i$ are polyhedra, the union $E$ of the boundaries of the $P_i$ is contained in a finite union of affine hyperplanes of $\A$.  Therefore there exists $(z_m)\in (\A\setminus E)\cap (z+\fQ_{+\infty})$ such that $z_m\to z$. Let $j\in \llbracket 1,n\rrbracket$ be such that $\{m\in \N\mid z_m\in P_j\}$ is infinite. Then $z=\lim z_n\in P_j$. Moreover, $z_m\in \Int(\A\cap B)\cap \Int(P_j)$ for $m\in \N$. In particular, $j\in J$, which proves that $\A\cap B\subset \bigcup_{j\in J} P_j$, which proves the lemma.
\end{proof}

\begin{Lemma}\label{l_U_encl_pf}
Assume that $\cS$ is positively cofree. Let $B$ be an apartment containing $\fQ_{+\infty}$.  Using Proposition~\ref{p_splt_apt}, we find a decomposition $(P_i,B_i,f_i)_{i\in \llbracket 1,n\rrbracket}$ with respect to $\fQ_{-\infty}$. Then: \[\A\cap B=\bigcup_{i\in \llbracket 1,n\rrbracket\mid (\A\cap B)\cap \mathring{P_i})\neq \emptyset}P_i.\]
\end{Lemma}

\begin{proof}
Let $\Omega=\A\cap B$.  Let $\phi:B\rightarrow \A$ be the isomorphism fixing $\Omega$, which exists by Lemma~\ref{l_iso_sect}. For $i\in \llbracket 1,n\rrbracket$, let $\phi_i=g_i\circ f_i$, where $g_i:B_i\rightarrow \A$ is the apartment isomorphism fixing $\A\cap B_i$. For $i\in \llbracket 1,n\rrbracket$, $f_i$ fixes $\Omega\cap P_i$. Let $i\in \llbracket 1,n\rrbracket$ be such that $\mathring{P_i}\cap \Omega\neq \emptyset$ and let $x\in \Omega\cap \mathring{P_i}$. Then $x+\fQ_{+\infty}=x+C^v_f\subset \Omega$. Therefore $P_i\cap \Omega$ has non-empty interior. For all $x\in P_i\cap \Omega$, $\phi_i(x)=\rho_{-\infty}(x)=x=\rho_{+\infty}(x)=\phi(x)$. Therefore $\phi_i=\phi$, by Remark~\ref{r_iso_fix_open}.

 For all $x\in P_i$, we have $\rho_{-\infty}(x)=\phi_i(x)=\phi(x)=\rho_{+\infty}(x)$ and thus \begin{equation}\label{e_P_i_sub_Om} P_i\subset \Omega,
\end{equation} by Corollary~\ref{c_char_aprtmt}.

 As $\Omega$ is convex, closed and with non-empty interior, we have $\Omega=\overline{\mathring{\Omega}}$. Therefore by Lemma~\ref{l_nnpty_int}, we have $\Omega=\bigcup_{i\in \llbracket 1,n\rrbracket\mid \Omega\cap \mathring{P_i}\neq \emptyset}\Omega\cap P_i=\bigcup_{i\in \llbracket 1,n\rrbracket\mid \Omega\cap \mathring{P_i}\neq \emptyset} P_i,$ which proves the lemma.
\end{proof}

\begin{Lemma}\label{l_un_encl_gen}
Let $B$ be an apartment. Then there exist $k\in \N$ and enclosed subsets $P_1,\ldots,P_k$ of $\A$ such that  $\A\cap B=\bigcup_{i=1}^k P_i$. 
\end{Lemma}

\begin{proof}
Let $d=\min \{d(\fQ,\fQ_{+\infty})\mid \fQ\text{ is a positive sector-germ of }B\}$. We prove the result by induction on $d$. If $d=0$, then $B$ contains $\fQ_{+\infty}$ and thus this follows from Lemma~\ref{l_U_encl_sct1} (or Lemma~\ref{l_U_encl_pf} if $\cS$ is positively cofree). Assume now that $d>0$. Let $\fQ$ be a sector-germ at infinity of $B$ such that $d(\fQ,\fQ_{+\infty})=d$. Let $\fQ'$ be a sector-germ at infinity adjacent to $\fQ$ and such that $d(\fQ_{+\infty},\fQ')=d-1$. Using \ref{a_sc}, we write $B=D_1\cup D_2$, where $D_1$ and  $D_2$ are opposite half-apartments such that for all $i\in \{1,2\}$, $D_i$ is contained in a apartment $B_i$ containing $\fQ'$.  We assume that for $i\in \{1,2\}$,  we can write $B_i\cap \A=\bigcup_{j=1}^{k_i} P_j^i$, where $k_i\in \N$ and  the $P_j^i$ are enclosed. Then $B\cap \A=\bigcup_{j=1}^{k_1} (D_1\cap P_j^1) \cup \bigcup_{j=1}^{k_2} (D_2\cap P_j^2)$, which proves the lemma. 
\end{proof}

\subsection{Convex union of polyhedra}

In this subsection, we prove that a convex set that can be written as a finite union of polyhedra is a polyhedron. 

\begin{Lemma}\label{l_un_enclsd_cvx1}
Let $A$ be a finite dimensional affine space. Let $n\in \N$ and $P_1,\ldots,P_n$ be polyhedra of $A$ having non-empty interior. We assume that $P:=\bigcup_{i=1}^n P_i$ is convex. Then $P$ is a polyhedron. More precisely, for $i\in \llbracket 1,n\rrbracket$, write $P_i=\bigcap_{D\in E_i} D$, where $E_i$ is a finite set of half-spaces of $A$. Let $E=\{D\in \bigcup_{i=1}^n E_i\mid P\subset D\}$. Then $P=\bigcap_{D\in E} D$. 
\end{Lemma}

\begin{proof}
By definition, we have $P\subset \bigcap_{D\in E} D$. If $P=A$ there is nothing to prove so we assume that there exists $x\in A\setminus P$.  Let $H_1,\ldots, H_k$ be the hyperplanes delimiting the elements of $\bigcup_{i=1}^n E_i$. We assume that the $H_i$ are all distinct. Then $\bigcup_{i,j\in \llbracket 1,k\rrbracket\mid i\neq j}H_i\cap H_j$ is a finite union of codimensional $2$ subspaces and thus \[\cN=\bigcup_{i=1}^k H_i\cup \bigcup_{i,j\in \llbracket 1,k\rrbracket\mid i\neq j} \supp(H_i\cap H_j\cup \{x\})\] is  a finite union of hyperplanes. Therefore $A\setminus \cN$ is dense in $A$. We fix $y\in \cN\cap \mathring{P}$. Then by choice of $\cN$,we have: \begin{equation}\label{e_cond_cN}
\forall z\in [x,y], |\{i\in \llbracket 1,k\rrbracket\mid z\in H_i\}|\leq 1.
\end{equation} As $P\cap [x,y]$ is a closed convex subset of $[x,y]$, there exists $z\in P$ such that $[x,y]\cap P=[y,z]$. We have $z\in \partial P$ and thus $z\in \bigcup_{i\in \llbracket 1,k\rrbracket} H_i$. By \eqref{e_cond_cN}, there exists a unique $i\in \llbracket 1,k\rrbracket$ such that $z\in H_i$. 

Let $V$ be a convex open neighborhood of $z$ such that $V\cap H_j=\emptyset$, for all $j\in \llbracket 1,k\rrbracket\setminus\{i\}$. Let $D_i$ be the half-space of $E$ delimited by $H_i$ and containing $P_i$. Let us prove: \begin{equation}\label{e_V_cap_Di}
V\cap P=V\cap D_i.
\end{equation}

We have $]z,y[\subset \mathring{D_i}\cap \mathring{P}$. Fix $x_0\in V\cap \mathring{P}\cap \mathring{D_i}$. Let $x'\in V\cap D_i$. If $x'\notin P$, then $]x_0,x'[$ meets $\partial P$ and thus $[x_0,x']\subset \partial D_i$: a contradiction. Therefore $V\cap D_i\subset V\cap P$. Assume for contradiction that there exists $a\in V\cap P \setminus D_i$. Then $\conv(a,D_i\cap V)\subset P$ contains $z$ in its interior: a contradiction. Therefore we have \eqref{e_V_cap_Di}.

 Let $a\in P$. Then $]z,a[\cap P$ contains an element $a'$ of $V\cap P$. For $a'\in A\setminus \{z\}$, denote by $\delta_{[z,a')}$ the closed ray based at $z$ and containing $a'$. Then $P\subset \bigcup_{a'\in V\cap P\setminus \{z\}} \delta_{[z,a')}=\bigcup_{a'\in V\cap D_i\setminus \{z\}}\delta_{[z,a')}\subset D_i$. Therefore $D_i\in E$ and $x\notin D_i$. Therefore $x\notin \bigcap_{D\in E} D$, which proves that $\bigcap_{D\in E}\subset P$. Therefore $P=\bigcap_{D\in E} D$. 
\end{proof}

\begin{Lemma}\label{l_supp_polyhd}
Let $A$ be a finite dimensional affine space, $k\in \Ne$, $D_1,\ldots,D_k$ be half-spaces of $A$ and $M_1,\ldots, M_k$ be their hyperplanes. Then their exists $J\subset \llbracket 1,k\rrbracket$ (maybe empty) such that $\mathrm{supp}(\bigcap_{i=1}^k D_i) =\bigcap_{j\in J} M_j$
\end{Lemma}

\begin{proof}
Let $d\in \Ne$ and $\ell\in \N$. Let $\mathcal{P}_{d,\ell}$:``for all affine space $A'$ such that $\dim A'\leq d$ and for all half-spaces $E_1,\ldots, E_\ell$ of $A'$, there exists $J\subset \llbracket 1,\ell \rrbracket $ such that $\mathrm{supp}( \bigcap_{i=1}^\ell E_i)=\bigcap_{j\in J} H_j$, where for all $j\in J$, $H_j$ is the hyperplane of $E_j$''.

It is clear that for all $\ell\in \N$, $\mathcal{P}_{1,\ell}$ is true and that for all $d\in \N$, $\mathcal{P}_{d,0}$ and $\mathcal{P}_{d,1}$ are true. Let $d\in \Z_{\geq 2}$ and $\ell\in \N$ and assume that (for all $d'\leq d-1$ and $\ell'\in \N$, $\mathcal{P}_{d',\ell'}$ is true) and that (for all $\ell'\in \llbracket 0,\ell\rrbracket$, $\mathcal{P}_{d,\ell'}$ is true). 

Let $A'$ be a $d$ dimensional affine space, $E_1,\ldots,E_{\ell+1}$ be half-spaces of $A'$ and $H_1,\ldots, H_{\ell+1}$ be their hyperplanes. Let $L=\bigcap_{j=1}^{\ell} E_j$ and $S=\mathrm{supp}(L)$. Then $E_{\ell+1}\cap S$ is either $S$ or a half-space of $S$. In the first case, $E_{\ell+1}\supset S\supset L$, thus $\bigcap_{i=1}^{\ell+1} E_i=L$ and thus by $\mathcal{P}_{d,\ell}$, $\mathrm{supp}(\bigcap_{i=1}^{\ell+1} E_i)=\bigcap_{j\in J} H_j$  for some $J\subset \llbracket 1,\ell\rrbracket$.

Assume that $E_{\ell+1}\cap S$ is a half-space of $S$. Then either $\mathring{E}_{\ell+1}\cap L\neq \emptyset$ or $\mathring{E}_{\ell+1}\cap L=\emptyset$. In the first case, we choose $x\in \mathring{E}_{\ell+1}\cap L$ and  a sequence $(x_n)\in (\In_r (L))^\N$  converging to $x$ (where $\In_r$ denotes the relative interior). Then for $n\gg 0$, $x_n\in  \mathring{E}_{\ell+1}\cap \In_r(L)$. Consequently,  $L\cap E_{\ell+1}$ has nonempty interior in $S$. Thus $\mathrm{supp}(\bigcap_{i=1}^{\ell+1} E_i)=S$ and  by $\mathcal{P}_{d,\ell}$, $\mathrm{supp}(\bigcap_{i=1}^{\ell+1} E_i)=\bigcap_{j\in J} H_j$  for some $J\subset \llbracket 1,\ell\rrbracket$.

Assume now that $\mathring{E}_{\ell+1} \cap L$ is empty. Then $L\cap E_{\ell+1}\subset H_{\ell+1}$, where $H_{\ell+1}$ is the hyperplane of $E_{\ell+1}$. Therefore $\bigcap_{i=1}^{\ell+1} E_i= \bigcap_{i=1}^{\ell+1} (E_i\cap H_{\ell+1})$  and thus by $\mathcal{P}_{d-1,\ell+1}$, $\mathrm{supp}(\bigcap_{i=1}^{\ell+1} E_i)=\bigcap_{j\in J} H_j$  for some $J\subset \llbracket 1,\ell+1\rrbracket$. 

Therefore $\cP_{d,\ell+1}$ is true. We deduce that $\cP_{d',\ell'}$ is true for every $d',\ell'\in \N$ and thus we get the lemma.

\end{proof}

\begin{Lemma}\label{l_un_enclsd_cvx2}
Let $A$ be a finite dimensional affine space.  Let $n\in \N$ and $P_1,\ldots, P_n$ be polyhedra such that $P:=\bigcup_{i=1}^n P_i$ is convex. Then $P$ is a polyhedron. More precisely, if $P_i=\bigcap_{D\in E_i} D$, for each $i\in \llbracket 1,n\rrbracket$, where $E_i$ is a set of half-spaces of $A$, then there exists $E\subset \bigcup_{i\in \llbracket 1,n\rrbracket}E_i$ such that  $P=\bigcap_{D\in E} D$. 
\end{Lemma}

\begin{proof}
Let $A'=\supp(P)$. As $P$ is convex and has non-empty interior in $A'$, we have $P=\overline{\In_{A'}(P)}$, where $\In_{A'}(P)$ is the interior of $P$ in $A'$. By Lemma~\ref{l_nnpty_int}, we can assume that all the $P_i$ have non-empty interior in $A'$. By Lemma~\ref{l_supp_polyhd}, there exists $F'\subset \bigcup_{i=1}^n E_i$ such that $A'=\bigcap_{D\in F'} D$. For $i\in \llbracket 1,n\rrbracket$, set $E_i'=\{D\cap A'\mid D\in E_i, D\nsupseteq A'\}$. Then $P_i=\bigcap_{D'\in E_i'}D'$, for $i\in \llbracket 1,n\rrbracket$. By Lemma~\ref{l_un_enclsd_cvx1}, there exists $E'\subset \bigcup_{i=1}^n E_i'$ such that $\bigcup_{i=1}^n P_i=\bigcap_{D'\in E'} D'$. For $D'\in E'$, choose $D_{D'}\in \bigcup_{i=1}^n E_i$ such that $D_{D'}\cap A'=D'$. Let $E=\{ D_{D'}\mid D'\in E'\}$. Then $\bigcap_{D\in E\cup F'}D=\bigcap_{D\in E}(D\cap A')=\bigcap_{D'\in E'} (D_{D'}\cap A')=A'\cap \bigcap_{D'\in E'} D'=P$. 
\end{proof}

\subsection{Conclusion}\label{ss_cnclsin_cl}

Combining Lemmas~\ref{l_un_encl_gen} and Lemma~\ref{l_un_enclsd_cvx2} (and using the fact that an enclosed subset of $\A$ is a polyhedron), we obtain:

\begin{Lemma}\label{l_ma2'_imp_ma2}
Let $\I$ be a weak masure satisfying \ref{a_ma2'}. Then $\I$ satisfies \ref{a_ma2}.
\end{Lemma}

\begin{Theorem}\label{t_MA2}
Let $\cS=(A,X,Y,(\alpha_i)_{i\in I}, (\alpha_i^\vee)_{i\in I})$ be a free Kac--Moody datum, which is either positively cofree or  associated with an indecomposable affine Kac--Moody matrix. Let $\I$ be  a space satisfying~\ref{a_ma1} \ref{a_sc},~\ref{a_oco} and \ref{a_wma3}. Then $\I$ satisfies~\ref{a_ma2}: for all apartments $A,B$, $A\cap B$ is enclosed and there exists an apartment isomorphism from $A$ to $B$  fixing $A\cap B$. 
\end{Theorem}

\begin{proof}
This is a combination of Lemmas~\ref{l_ma2'_aff}, \ref{l_ma2'_pf} and Lemma~\ref{l_ma2'_imp_ma2}.
\end{proof}

\begin{Remark}
\begin{enumerate}
\item Note that in \cite[Corollary 3.7]{hebert2022new}, we proved that definition of an abstract masure of \cite[Définition 2.1]{rousseau2011masures} is equivalent to Definition~\ref{d_w_mas}, when $(\alpha_i^\vee)$ is free. Actually, by adapting  the proof we gave here, we obtain that the two definitions are equivalent, when $(\alpha_i^\vee)$ is positively cofree. Indeed, the arguments of \cite{hebert2022new} only use the positively freeness.

\item In \cite[Proposition 4.25 and Remark 4.26]{hebert2020new}, it is proved that if $\I$ is associated with a split Kac--Moody group, then every  enclosed subset of $\A$ with non-empty interior is the intersection of two apartments of $\I$. This seems to indicated that~\ref{a_ma2} is optimal. Note that this question of which subsets of $\A$ are the intersections of two apartments is studied in \cite{abramenko2010intersection} in the case where $\I$ is  a building (spherical or affine).

\end{enumerate}
\end{Remark}

\chapter{Parahoric subgroup associated with a point}\label{C_d_mas}

Let $\cS=(A,X,Y,(\alpha_i)_{i\in I},(\alpha_i^\vee))$ be a free Kac--Moody datum and $\fG=\fG_\cS$ be the split minimal Kac--Moody group associated with $\cS$ in Chapter~\ref{C_splt_KM_grps}. Let $\cF$ be a field equipped with a non-trivial valuation $\omega:\cF\rightarrow \R\cup \{+\infty\}$.    We assume that $\omega(\cF)\supset \Z$, which is possible up to renormalization and if $\omega(\cF)$ is discrete, we assume that $\omega(\cF)=\Z$. We set $G=\fG(\cF)$.      We mainly use gothic letters to denote the subfunctors of  $\fG$: $\fB^{\pm}$, $\fU^{\pm}$, $\fT$, $\fN$, ... and we use roman letters to denote their evaluations at $\cF$: $B^{\pm}=\fB^{\pm}(\cF)$, $U^{\pm}=\fU^{\pm}(\cF)$, $T=\fT(\cF)$, $N=\fN(\cF)$, ...

We want to define the masure of $(G,\omega)$.  We  use the same methods as in \cite{bruhat1972groupes} (in the reductive case), \cite{gaussent2008kac} and  \cite{rousseau2016groupes}: \begin{enumerate}
                                                                                                                                                                                          \item Let $\fT$ be a maximal split torus of $\fG$ and $T=\fT(\cF)$.   Let $Y$ be the cocharacter lattice of $(\fG,\fT)$ and let $N$ be the normalizer of $T$ in $G$. We define the standard apartment of the masure  as $\A=Y\otimes \R$.  There is a natural action $\nu^v$ of $N$ on $\A$ by linear maps (and independant of the valuation on $\cF$). We define an action $\nu$ of $N$ on $\A$ by affine automorphisms such that for each $n\in N$, $\nu^v(n)$ is the linear part of $\nu(n)$.

                                                                                                                                                                                          \item We define for each $x\in \A$ its fixator $G_x$ in the (not yet defined) masure. If $x\in \A$, $G_x$ is defined as the subgroup of $G$ generated by $U_{x,\infty}^+$, $U_{x,\infty}^-$ and $N_x$, where $N_x$ is the fixator of $x$ in $N$ and $U_{x,\infty}^+$ is defined as the ``fixator'' of $x$ in $U^+$ (and similarly for $U_{x,\infty}^-$). The group $U_{x,\infty}^+$ is defined via the coordinates $X_\alpha$, $\alpha\in \Delta_+$ defined in Subsection~\ref{ss_twstd_exp} (see Definition~\ref{d_U_Vma} or Proposition~\ref{p_Rou4.5}).

                                                                                                                                                                                          \item The masure $\I$ is then defined as $(G\times \A)/\sim$, where $\sim$ is an equivalence relation defined in such a way that $G_x$ is the fixator of $x$ in $G$, for $x\in \A$.

Note that eventually, we obtain that  $G_x$ can be defined similarly as in the reductive case:   $G_x=\langle U_{\alpha,x}\mid \alpha\in \Phi, N_x\rangle$ (see Proposition~\ref{p_sphrcl_para}). This definition involves the minimal group only and is simpler. We could in principle take it as a definition. However, we use the masure to prove this result and we do not know how to prove that $G_x$ and the masure have ``good properties'' without using the maximal groups.

                                                                                                                                                                                         \end{enumerate}

In this chapter, we define the parahoric subgroups $G_x$, for $x\in \A$.

Note that we change the notation  compared to \cite{gaussent2008kac} and \cite{rousseau2016groupes}: if $\cV$ is a filter on $\A$,  we write $U_{\cV,\fin}$ instead of $U_{\cV}$, $U_{\cV,\fin}^{\pm}$ instead of $U_{\cV}^{\pm}$, $U_{\cV,\infty}^+$ instead of $U_{\cV}^{pm+}$, $U_{\cV,\infty}^{-}$ instead of $U_{\cV}^{nm-}$, $U_{\cV,\infty+}$ instead of $U_{\cV}^{pm}$ and $U_{\cV,\infty-}$ instead of $U_{\cV}^{nm}$. Intuitively, the index ``$\infty$'' means that infinite products (of $X_\alpha$, $\alpha\in \Delta$) are involved whereas ``$\fin$'' means that only finite products (of $x_\alpha$, with $\alpha\in \Phi$) are involved.

 In Section~\ref{s_act_N_A}, we define the vectorial and affine actions on $N$ on $\A$.

 In Section~\ref{s_unip_sbgrps}, we associate unipotent subgroups to filters on $\A$. As a particular case, we introduce the groups $U_{x,\infty}^+$ and $U_{x,\infty}^-$, for $x\in \A$. We also prove the Iwasawa decomposition.

 In Section~\ref{s_parahoric}, we introduce and  study the subgroup $G_x$, for $x\in \A$. We obtain a decomposition of $G_x$ which is crucial in our study of masures.

\section{Action of $N$ on $\A$}\label{s_act_N_A}

Recall that as a set, the standard apartment is  $\A=Y\otimes \R$, where $Y$ is the cocharacter lattice of $\cS$. We equip it with the set of walls $\cM=\{\alpha^{-1}(\{\lambda\})\mid \alpha\in \Phi, \lambda\in \Lambda\}$.

 In this section, we define an action of $N=\fN(\cF)$ on $\A$ by affine automorphisms. We start by defining an action $\nu^v$ of $\fN(\cR)$ on $\A$ by linear maps, for any ring $\cR$. We then define an action of $N$ on $\A$ by affine automorphisms, which takes the valuation $\omega$ into account  and such that the linear part of $\nu(n)$ is $\nu^v(n)$, for all $n\in N$.

\subsection{Action of $\fN(\cR)$ on $\A$ by linear automorphisms}

Let $\A=Y\otimes \R$. Let $\cR$ be a ring. Recall from Subsection~\ref{ss_Tits_minKM} that we set $\fN(\cR)=\langle \tilde{s_i}\mid i\in I,\fT(\cR)\rangle\subset \fG(\cR)$ and that if $\cR$ is a field with at least four elements, then $\fN(\cR)$ is the normalizer of $\fT(\cR)$ in $\fG(\cR)$.

We set: \begin{equation}\label{e_N_0_T_0'}
\fN_0(\cR)=\langle \tilde{s}_i\mid i\in I\rangle\subset \fN( \cR)\text{ and }\fT_0'(\cR)=\langle \alpha_i^\vee(-1)\mid i\in I\rangle\subset \fN(\cR)\cap \fT(\cR).
\end{equation}\index[notation]{n@$\fN_0$}

By \ref{a_KMT5}, we have: \begin{equation}\label{e_dec_N}
\fN(\cR)=\fT(\cR)\cdot \fN_0(\cR).
\end{equation}

\begin{Proposition}\label{p_d_nu_v}
\begin{enumerate}
\item Let $\cR$ be a ring.  There exists a unique action $\nu^v$ of $\fN(\cR)$ on $\A$ by linear automorphisms such that $\nu^v(T)=\{\Id\}$ and $\nu^v(\tilde{s_i})=r_i$, for all $i\in I$. It is induced by the restriction of $\Ad_\cR$ to $\fN(\cR)$.

\item For all  $n_0\in \fN_0(\cR)$, there exists a unique $t'\in \fT_0'(\cR)$ and $n_0'\in \fN_0(\cR)$ such that: \begin{enumerate}
\item $n_0=t'n_0'$

\item $n_0'$ can be written $\tilde{s}_{i_1}\ldots \tilde{s}_{i_k}$, with $k=\ell(\nu^v(n_0))$ and $i_1,\ldots,i_k\in I$.
\end{enumerate}  

\item We have $\ker(\nu^v)=\fT(\cR)$ and $\fT_0'(\cR)=\fN_0(\cR)\cap \ker(\nu^v)=\fN_0(\cR)\cap\fT(\cR)$.

\end{enumerate}
\end{Proposition}

\begin{proof}
(1) Let $i\in I$. By \eqref{e_Ad_U}, $\Ad_{\cR}(\tilde{s_i})=s^*_i=r_i^{\mathrm{ad}}$, for the notation of \cite[Proposition 4.18 and Definition 7.46]{marquis2018introduction}. By \cite[Proposition 4.18 (2) and (7.31)]{marquis2018introduction}, $s_i^*$ stabilizes $Y$ and acts by $r_i$ on $Y$. By \eqref{e_Ad_T}, since $Y\subset \cU_0$, we have $\Ad_\cR(t)(x)=x$, for all $x\in Y$. Therefore if $\nu^v:\fN(\cR)\rightarrow \mathrm{GL}(Y)$ is defined by $\nu^v(n)=\Ad_\cR(n)|_{Y}^{Y}$, for $n\in \fN(\cR)$, then we have $\nu^v(\tilde{s_i})=r_i$, for all $i\in I$ and $\nu^v(t)=\Id$, for all $t\in T$. This proves the existence of $\nu^v$, by considering $\mathrm{GL}(Y)$ as a subgroup of $\mathrm{GL}(\A)$.

By \eqref{e_dec_N}, we have  $\fN(\cR)=\fT(\cR)\cdot \fN_0(\cR)$, which proves the uniqueness of $\nu^v$.

(2) By definition, we have $\fT(\cR)\subset \ker(\nu^v)$. Let $n\in \fN_0(\cR)$.   Write $n=\tilde{s}_{i_1}\ldots \tilde{s}_{i_k}$, with $k\in \N$ and $i_1,\ldots,i_k\in I$. Let $w=\nu^v(n)$.  Assume that $\ell(w)<k$. By the word property in the Coxeter group $W^v$ (\cite[Theorem 3.3.1]{bjorner2005combinatorics}), up to  changing the expression $r_{i_1}\ldots r_{i_k}$ by using braid-moves, we can assume that there exists $m\in \llbracket 1,k-1\rrbracket$ such that $i_m=i_{m+1}$. By \cite[Proposition 7.57]{marquis2018introduction}, the $\tilde{s}_i$ satisfy the braid relations and thus  $\tilde{s}_{i_1}\ldots \tilde{s}_{i_k}$ is unchanged by the change of expression. Then \begin{align*}
n&=\tilde{s}_{i_1}\ldots \tilde{s}_{i_{m-1}} (\tilde{s}_{i_m})^2 \tilde{s}_{i_{m+2}}\ldots \tilde{s}_{i_k}=
\tilde{s}_{i_1}\tilde{s}_{i_2}\ldots \tilde{s}_{i_{m-1}}\alpha_{i_m}^\vee(-1) \tilde{s}_{i_{m+2}}\ldots \tilde{s}_{i_k}\text{ by \eqref{e_s_i_sq}}.
\end{align*}

Now by \ref{a_KMT5}, we have \[\tilde{s}_{i_1}\ldots \tilde{s}_{i_{m-1}}\alpha_{i_m}^\vee(-1)=(\tilde{s}_{i_1}\ldots \tilde{s}_{i_{m-1}}.\alpha_{i_m}^\vee)(-1)\in \fT_0'(\cR).\] Therefore $n\in \fT_0'(\cR)\tilde{s}_{i_1}\tilde{s}_{i_2}\ldots \tilde{s}_{i_{m-1}} \tilde{s}_{i_{m+2}}\ldots \tilde{s}_{i_k}$.  By repeating this reasoning $(k-\ell(w))/2-1$ times (if $k\geq \ell(w)+2$), we deduce that $n_0$ can be written $n_0=t' \tilde{s}_{j_1}\ldots \tilde{s}_{j_\ell(w)}$, where $t'\in \fT_0'(\cR)$ and $j_1,\ldots,j_{\ell(w)}\in I$.

Let $t''\in \fT'_0(\cR)$ and $j_1'',\ldots,j_{\ell(w)}''\in I$ be such  that $t''\tilde{s}_{j_1''}\ldots \tilde{s}_{j_{\ell(w)}''}=n_0=t'\tilde{s}_{j_1}\ldots \tilde{s}_{j_\ell(w)}$. We have $\nu^v(n_0)=r_{j_1''}\ldots r_{j_{\ell(w)}''}=r_{j_1}\ldots r_{j_{\ell(w)}}$. By the word property in the Coxeter group $W^v$ (\cite[Theorem 3.3.1]{bjorner2005combinatorics}), we can pass from the expression $r_{j_1}\ldots r_{j_{\ell(w)}}$ of $w$ to the $r_{j_1''}\ldots r_{j_{\ell(w)}''}$ by using only braid-moves.  By \cite[Proposition 7.57]{marquis2018introduction}, we then have $\tilde{s}_{j_1''}\ldots \tilde{s}_{j_{\ell(w)}''}=\tilde{s}_{j_1}\ldots \tilde{s}_{j_{\ell(w)}}$. Therefore $t''=t'$, which proves the uniqueness in (2).

(3) By definition, $\fT(\cR)\subset \ker(\nu^v)$.  Let $n\in \ker(\nu^v)$. By \eqref{e_dec_N}, we can write $n=tn_0$, where $t\in \fT(\cR)$ and $n_0\in \fN_0(\cR)$. We have $\nu^v(n)=1=\nu^v(tn_0)=\nu^v(n_0)$.  By (2), we have $n_0\in \fT_0'(\cR)\subset \fT(\cR)$ (since $\ell(\nu^v(n_0))=0$), which proves that $n_0\in \fT(\cR)$. Therefore $n\in \fT(\cR)$, which proves that $\ker(\nu^v)=\fT(\cR)$.

We have $\fT_0'(\cR)\subset \fN_0(\cR)\cap \ker(\nu^v)$. Let $n\in  \fN_0(\cR)\cap \ker(\nu^v)$. Then $\nu^v(n)=\Id$ and thus by (2), $n\in \fN_0(\cR)\cap \ker(\nu^v)$. This proves (3).
\end{proof}

\subsection{Action of $N$ on $\A$}\label{ss_act_N_A}
In this subsection, we define an action $\nu$ of $N$ on $\A$ by affine automorphisms, following \cite[2.9]{rousseau2006groupes}.  We set $\Lambda=\omega(\cF^\times)$ and for all $\alpha\in \Phi$, we set $\Lambda_\alpha=\Lambda$. We denote by $\cl$ the associated enclosure, see \ref{sss_enclosure}. 

We define an action $\nu$ of $T=\fT(\cF)$ on $\A$ by translation as follows: if $t\in T$, then  $\nu(t)$ translates $\A$ by the unique vector (still denoted) $\nu(t)$ satisfying: \begin{equation}\label{e_act_T}
\chi(\nu(t))=-\omega(\chi(t)), \forall \chi\in X.
\end{equation}

We also set \[T_0=\fT(\cO).\]\index[notation]{t@$T_0$} Then $T_0$ is the kernel of $\nu:T\rightarrow \mathrm{\Aut}(\A)$.

\begin{Lemma}\label{l_act_cort}
Let $\alpha^\vee\in \Phi^\vee$ and $r\in \cF^\times$. Then $\nu(\alpha^\vee(r))=-\omega(r)\alpha^\vee$.
\end{Lemma}

\begin{proof}
Let $\chi\in X$, $t=\alpha^\vee(r)$ and $\bt=\nu(t)$. We have $\chi(\bt)=-\omega(\chi(t))=-\omega(r^{\chi(\alpha^\vee)})=-\chi(\alpha^\vee)\omega(r)=\chi(-\omega(r)\alpha^\vee)$ and hence $\bt=-\omega(r)\alpha^\vee$.
\end{proof}

\begin{Lemma}\label{l_nu_equiv}
The action of $\nu$ of $T$ on $\A$ is $W^v$-equivariant.
\end{Lemma}

\begin{proof}
Let $i\in I$ and $t\in T$. Let $\chi\in X$. Then $\chi(\nu(r_i(t)))=-\omega(\chi(r_it))=-\omega(r_i.\chi(t))$ and $\chi(r_i.\nu(t))=(r_i.\chi)(\nu(t))=-\omega((r_i.\chi)(t))$, which proves that $r_i.\nu(t)=\nu(r_i.t)$. As $W^v$ is generated by the $r_i$, $i\in I$, we deduce the result.
\end{proof}

\begin{Proposition}\label{p_d_nu}(see \cite[Lemme 2.9.2]{rousseau2006groupes})
 There exists a unique   action $\nu$ of $N$ on $\A$ by affine automorphisms such that $\nu$ restricts to \eqref{e_act_T} on $T$ and such that for all $i\in I$, $\nu(\tilde{s_i})=\nu^v(\tilde{s_i})=r_i$. For all $n\in N$, the linear part of $\nu(n)$ is $\nu^v(n)$. 

\end{Proposition}

\begin{proof}
Recall the definitions of $\fN_0(\cF)$ and $\fT_0'(\cF)$ in \eqref{e_N_0_T_0'}. Set $N_0'=\fN_0(\cF)$ (we add a prime to avoid a confusion with $N_0$, which will be defined as the fixator of $0$ in $N$ for the action of $\nu$ on $\A$). The uniqueness of $\nu$ follows from \eqref{e_dec_N}, since $\nu^v(\tilde{s_i})=r_i$, for all $i\in I$, and the $\tilde{s_i}$ generate $N_0'$. It remains to  prove the existence of $\nu$.

 The uniqueness of $\nu$ follows from \eqref{e_dec_N}, since $\nu^v(\tilde{s_i})=r_i$, for all $i\in I$, and the $\tilde{s_i}$ generate $N_0'$. It remains to  prove the existence of $\nu$.

 Let $T_0=\fT(\cO)\subset \fG(\cF)$. By \eqref{e_act_T}, we have $\nu(T_0)=\{0\}$: $T_0$ acts trivially on $\A$. Thus it suffices to define an action of $N/T_0$ on $\A$. Let us prove: \begin{equation}\label{e_N_T_0}
 N/T_0=T/T_0\rtimes (N_0'\cdot T_0) /T_0.
 \end{equation} By \ref{a_KMT5}, $N_0'\cdot T_0$  is a subgroup of $N$.  By \eqref{e_dec_N}, we have $N=T\cdot (N_0'\cdot T_0)$ and hence $N/T_0=T/T_0\cdot (N_0'\cdot T_0)/T_0$.   As $T$ is normal in $N$, $T/T_0$ is normal in $N/T_0$. Let $x\in (N_0'\cdot T_0)/T_0\cap T/T_0$. Write $x=nT_0=tT_0$, with $n\in N_0'$,  and $t\in T_0$. Then $\nu^v(n)=1$ and hence $n\in T\cap N_0'$. Then by Proposition~\ref{p_d_nu_v}, $n\in \fT_0'(\cF)$ and by  \eqref{e_N_0_T_0'}, $n\in T_0$.  This proves \eqref{e_N_T_0}.

  Define $\nu':T/T_0\rtimes (N_0'\cdot T_0)/T_0\rightarrow \mathrm{Aut}(\A)$ by $\nu'(tT_0,n_0 T_0)=\nu^v(n_0)+\nu(t)$, for $n_0\in N_0'\cdot T_0$ and $t\in T$ (where $\nu^v(n_0)+\nu(t)$ is the map $\A\rightarrow \A$ sending each $a\in \A$ on $\nu^v(n_0)+a+\nu(t).a$). Let us prove that $\nu'$ is a group morphism. Let $i\in I$ and $t\in T$. Then by \ref{a_KMT5}, we have $\tilde{s}_i t \tilde{s}_i^{-1}=r_i(t)$. Therefore $\nu(\tilde{s}_i t\tilde{s}_i^{-1})=r_i.\nu(t)=\nu^v(\tilde{s}_i).\nu(t)$, by Lemma~\ref{l_nu_equiv}. By induction, we deduce that if $n\in N_0'$, then: \begin{equation}\label{e_cnj_N_nu}
\nu(ntn^{-1})=\nu^v(n).\nu(t).
\end{equation}

   Let $(\overline{t},\overline{n}),(\overline{t'},\overline{n'})\in T/T_0\rtimes (N_0'\cdot T_0)/T_0$. Then $\overline{t}\overline{n}\cdot \overline{t'}\overline{n'}=\overline{t}
  \overline{nt'^{-1}n^{-1}}\cdot \overline{nn'}.$ Thus \begin{align*}
  \nu'(t nt'n')&=\nu(t)+\nu(nt'n^{-1})+\nu^v(nn')\\
  &= \nu(t)+\nu^v(n).\nu(t')+\nu^v(nn') &\text{by \eqref{e_cnj_N_nu}}\\
  &= \nu'(tn)\circ \nu'(t'n').
  \end{align*}

   Then we define $\nu:N\rightarrow \mathrm{Aut}(\A)$ by $\nu(n)=\nu'(nT_0)$, for $n\in N$ and $\nu$ satisfies the conditions of the proposition.
\end{proof}

\begin{Definition}\label{d_ext_aff_W}
We define the \textbf{extended affine Weyl group} as $\widetilde{W^a}=\nu(N)$.\index[notation]{w@$\widetilde{W^a}$}
\end{Definition}

\begin{Proposition}\label{p_des_ext_aff}
We have $\widetilde{W^a}=W^v\ltimes (Y\otimes \Lambda)$ and $\widetilde{W^a}\subset \Aut^{\vw}(\A,(\Lambda)_{\alpha\in \Phi})$, where $\Aut^{\vw}(\A,(\Lambda))$ is defined in \eqref{e_Aut_vw}. 
\end{Proposition}

\begin{proof}
By \eqref{e_act_T}, we have $\nu(T)=Y\otimes \Lambda$. By Proposition~\ref{p_d_nu}, we have $\nu(\fN_0(\cF))=W^v$. Therefore by \eqref{e_dec_N}, we have $\nu(N)=\langle Y\otimes \Lambda, W^v\rangle=W^v\ltimes (Y\otimes \Lambda)$. As $W^v$ and $Y\otimes \Lambda$ stabilize $\{\alpha^{-1}(\{k\})\mid \alpha\in \Phi, k\in \Lambda\}$, we deduce the result.
\end{proof}

\section{Unipotent and pro-unipotent subgroups  associated with filters on $\A$}\label{s_unip_sbgrps}

In this section, we introduce and study unipotent and pro-unipotent subgroups associated with filters on $\A$. We introduce many definitions, but the more important definitions are the one of $U_{\cV,\infty}^+$ and $U_{\cV,\infty}^-$, for $\cV$ a filter on $\A$ (see Definition~\ref{d_U_Vma}). They will be the fixators of $\cV$ in $U^+$ and $U^-$ respectively (see Theorem~\ref{t_fltrs_gd_fix}), for the action on the masure.

\subsection{Fonction $f_\cV$ associated with a filter $\cV$}\label{ss_f_cV}

\subsubsection{Extension of the valuation on the root spaces}\label{sss_ext_val}

Let $\alpha\in \Delta$. Let $\cB_\alpha$ be a $\Z$-basis of $\ffg_{\alpha,\Z}$. Let $\underline{u_\alpha}\in \ffg_{\alpha,\Z}\otimes \cF$. Write $\underline{u_\alpha}=\sum_{b\in \cB_\alpha} x_b b$, where $(x_b)\in \cF^{\cB_\alpha}$. We set $\omega(\underline{u_\alpha})=\min_{b\in \cB_\alpha} \omega(x_b)\in \Lambda\cup \{+\infty\}$. Then $\omega:\ffg_{\alpha,\Z}\otimes \cF\rightarrow \Lambda\cup \{+\infty\}$\index[notation]{o@$\omega(\underline{u_\alpha})$} is well-defined, independently of the choice of $\cB_\alpha$.

\subsubsection{Fonction $f_\cV$}

Let $u\in U^+$. Consider $u$ as an element of $U^{ma+}$ and write $u=\prod_{\alpha\in \Delta_+} X_{\alpha}(\underline{u_\alpha})$, where $\underline{u_\alpha}\in \prod_{\alpha\in \Delta_+} \ffg_{\alpha,\Z}\otimes \cF$. We want to define the masure in such a way that $u$ fixes $\bigcap_{\alpha\in \Delta_+}\{a\in \A\mid \omega(\underline{u_\alpha})+\alpha(a)\geq 0\}$. To that end, it is convenient to associate  to each filter a function $f_\cV$ which enables to  determine when $\alpha(\cV)+\lambda\geq 0$,  for  $\alpha\in X$ and $\lambda\in \R$. We embed $\R$ in an ordered monoid $\tilde{\R}$, following \cite[6.4.1]{bruhat1972groupes} and \cite[4.2]{rousseau2016groupes}.

Let $\tilde\R$\index[notation]{r@$\tilde{\R}$} be the monoid  consisting of the elements $r,r^+$ for $r\in \R$ and of $+\infty$, satisfying, for all $r,s\in \R$: \begin{enumerate}
\item if $r<s$, then $r<r^+<s<s^+<+\infty$,

\item  $r^++s^+=r^++s=r+s^+=(r+s)^+$. 
\end{enumerate}  We denote by $\tilde{\Lambda}$ the set of the elements of $\tilde{\R}$ which are a lower bound of a bounded below subset of $\Lambda$.  We have $\tilde{\Lambda}=\Lambda\cup \{+\infty\}$ if $\Lambda$ is discrete and $\tilde{\Lambda}=\Lambda\cup \{r^+\mid r\in \R\}\cup \{+\infty\}$ otherwise. If $E$ is a bounded below subset of $\Lambda$, we denote by $\inf_{\tilde{\Lambda}}(E)$ its lower bound in $\tilde{\Lambda}$, i.e the lower bound of $E$ regarded as a subset of $\tilde{\Lambda}$.  Note that if $E$ is a subset of 
$\Lambda$ such that $\inf_{\tilde{\Lambda}}(E)\in \R$, then $E$ admits a minimum.

For $\lambda\in \tilde{\Lambda}$, we set $\cF_{\omega\geq \lambda}=\cF_{\geq \lambda}=\{x\in \cF\mid \omega(x)\geq \lambda\}$\index[notation]{f@$\cF_{\omega\geq \lambda}=\cF_{\geq \lambda}$}.

 For $\lambda\in \tilde{\Lambda}$ and $\alpha\in \Phi$,  \[U_{\alpha,\lambda}:=\{x_\alpha(a)\mid a\in \cF_{\omega\geq \lambda}\}\] is a subgroup of $U_\alpha$. For $\alpha\in X$ and $\lambda\in \R$, we set $D(\alpha,\lambda)=\{x\in \A\mid \alpha(x)+\lambda\geq 0\}$.

Let $\cV$ be a filter on $\A$. For $\alpha\in X$, we set: \[F_\cV(\alpha)= \{\lambda\in \Lambda \mid  D(\alpha,\lambda)\in \cV\}.\]\index[notation]{f@$F_{\cV}(\alpha)$} Then we have: \begin{equation}\label{e_terminal}
\forall \lambda \in \Lambda, (\lambda\in F_\cV(\alpha))\Rightarrow [\lambda,+\infty[\cap \Lambda \subset F_\cV(\alpha).
\end{equation} 

For $\alpha\in X$ and $\cV$ a filter on $\A$, we set: 
\[f_\cV^\Lambda(\alpha)=f_\cV(\alpha)=\mathrm{inf}_{\tilde{\Lambda}} \{\lambda\in \Lambda\mid \cV\subset D(\alpha,\lambda)\}=\mathrm{inf}_{\tilde{\Lambda}} \left(F_\cV(\alpha)\right)\in \tilde{\Lambda}\cup \{-\infty,+\infty\}.\]\index[notation]{f@$f_\cV(\alpha)$}

\begin{Lemma}\label{l_char_fV}
Let $\lambda\in \R$, $\alpha\in X$ and $\cV$ be a filter on $\A$. Then $f_{\cV}(\alpha)\leq \lambda$ if and only if $\cV\subset D(\alpha,\lambda)$. 
\end{Lemma}

\begin{proof}
If $\cV\subset D(\alpha,\lambda)$, then by definition of the infimum, we have $f_{\cV}(\alpha)\leq \lambda$. Conversely, let $\lambda\in \R$ be such that $f_{\cV}(\alpha)\leq \lambda$. First assume  $f_{\cV}(\alpha)<\lambda$. Then there exists $\mu\in \Lambda$ such that $\cV\subset D(\alpha,\mu)$ and $\mu<\lambda$. Then $D(\alpha,\mu)\subset D(\alpha,\lambda)$ and hence $\cV\subset D(\alpha,\lambda)$. Assume now $f_{\cV}(\alpha)=\lambda$. Then as $\lambda\in \R=\tilde{\R}\setminus \{r^+\mid r\in \R\}$, $\lambda\in F_{\cV}(\alpha)$ and hence $\cV\subset D(\alpha,\lambda)$. 
\end{proof}

\begin{Example}\label{ex_f_cV}
\begin{enumerate}
\item If $\cV=\A$ and $\alpha\in \Delta$, we have $F_{\cV}(\alpha)=\emptyset$ and $f_\A(\alpha)=+\infty$.

\item If $\cV=\fQ_{+\infty}$ and $\alpha\in \Phi_+$, then we have $F_{\fQ_{+\infty}}(\alpha)=\Lambda$ and $f_{\fQ_{+\infty}}(\alpha)=-\infty$. If $\alpha\in \Phi_-$, we have $F_{\fQ_{+\infty}}(\alpha)=\emptyset$ and thus $f_{\fQ_{+\infty}}(\alpha)=+\infty$. 

\item If $\cV=C_0^-:=\germ_0(-C^v_f)$ and $\alpha\in \Delta_+$, we have $F_{C_0^-}(\alpha)=\R_{>0}\cap \Lambda$ and $f_{C_0^-}(\alpha)=\begin{cases} &=0^+\text{ if }\overline{\Lambda}=\R\\ 
&= 1\text{ if }\overline{\Lambda}=\Z\end{cases}$. Indeed, let $\lambda\in \R$. Then $C_0^-\subset D(\alpha,\lambda)$ if and only if there exists a neighborhood $\Omega$ of $0$ in $\A$ such that  $D(\alpha,\lambda)\supset \Omega\cap -C^v_f$. Therefore if $\lambda\in \R_{>0}$, $0\in \mathring{D}(\alpha,\lambda)$ and hence $C_0^-\subset D(\alpha,\lambda)$. Conversely, if $\lambda\in \R$ is such that $C_0^-\subset D(\alpha,\lambda)$, then there exists $x\in -C^v_f$ such that $\alpha(x)+\lambda\geq 0$ and thus $\lambda>0$. 

\item If $\cV=C_0^+:=\germ_0(C^v_f)$ and $\alpha\in \Delta_+$, we have $F_{C_0^+}(\alpha)=\R_{\geq 0}$ and $f_{C_0^+}(\alpha)=0$. Indeed, if $\lambda\in \R$ is such that $D(\alpha,\lambda)\supset C_0^+$, then there exists $\Omega\in C_0^+$ such that  $\Omega\subset D(\alpha,\lambda)$. Then $0\in \overline{\Omega}\subset D(\alpha,\lambda)$ and thus $\lambda\in \R_{\geq 0}$. Conversely, $D(\alpha,0)\supset C^v_f$ and $C^v_f\in C_0^+$ and thus $F_{C_0}^+(\alpha)=\R_{\geq 0}$.

\end{enumerate}

\end{Example}

\begin{Proposition}\label{p_f_cV_alpha}
Let $\alpha\in \Delta$. Then we have: \begin{enumerate}
\item Let $\Omega$ be a non-empty subset of $\A$. Then $F_\cV(\alpha)=[\sup_\R(-\alpha(\Omega)),+\infty[\cap \Lambda$ (this set can be empty if $ \sup_\R(-\alpha(\Omega))=+\infty$).

\item Let $\cV$, $\cV'$ be two filters on $\A$ such that $\cV\subset \cV'$. Then $F_\cV(\alpha)\supset F_{\cV'}(\alpha)$ and $f_{\cV}(\alpha)\leq f_{\cV'}(\alpha)$. In particular, for all $\Omega\in \cV$, we have $f_{\cV}(\alpha)\leq f_{\Omega}(\alpha)$. 

\item Let $\cV$ be a filter on $\A$ such that $F_\cV(\alpha)\notin \{\emptyset,\Lambda\}$. Then one of the two mutually exclusive  possibilities holds: \begin{enumerate}
\item There exists $\Omega\in \cV$ such that $f_\cV(\alpha)=f_\Omega(\alpha)$.

\item There exists $(\Omega_n)\in \cV^{\N}$ such that $(\Omega_n)$ is decreasing (for $\subset$) and  $f_{\cV}(\alpha)=\inf_{\tilde{\Lambda}} \{f_{\Omega_n}(\alpha)\mid n\in \N\}\in \{r^+\mid r\in \R\}$. 
\end{enumerate}
\item Let $\cV$ be a filter on $\A$. If $\alpha\in \Delta$ (resp. $\alpha\in \Phi$), then $f_{\cV}(\alpha)=f_{\cl^\Delta(\cV)}(\alpha)$ (resp. $f_{\cV}(\alpha)=f_{\cl(\cV)}(\alpha)$). 

\item Let $\Omega$ be a subset of $\A$.  The map $f_{\Omega}$ is \textbf{concave} in the sense of \cite{bruhat1972groupes}
\begin{equation}\label{e_concavity}
\forall \alpha,\beta\in X,f_{\Omega}(\alpha+\beta)\leq f_{\Omega}(\alpha)+f_{\Omega}(\beta)\text{ and }f_{\Omega}(0)=0.
\end{equation}

\item   Let $\cV$ be a filter on $\A$ and $\alpha\in X$. Then if $f_{\cV}(\pm\alpha)=-\infty$, we have $f_{\cV}(\mp\alpha)=+\infty$. If $-\infty\notin \{f_{\cV}(\alpha),f_{\cV}(-\alpha)\}$, then $f_{\cV}(\alpha)+f_{\cV}(-\alpha)\geq 0$. 

\end{enumerate}
\end{Proposition}

\begin{proof}
(1) Let $\lambda\in \Lambda$. Then $\Omega \subset D(\alpha,\lambda)$ if and  only if $\alpha(\Omega)+\lambda\geq 0$ if and only if $\lambda\geq \sup_\R(-\alpha(\Omega))$. Therefore $F_\cV(\alpha)=\Lambda\cap [\sup_\R(-\alpha(\Omega)),+\infty[$. 

(2) is clear.

 (3) As $F_\cV(\alpha)\neq \Lambda$, the set $F_\cV(\alpha)$ is bounded below, by \eqref{e_terminal}. If $F_\cV(\alpha)$ admits a minimum, then we have (a). Otherwise, consider a sequence $(\lambda_n)\in \Lambda^\N$ such that $(\lambda_n)$ decreases to $\inf_\R(F_\cV(\alpha))$. For $n\in \N$, set $\Omega_n=D(\alpha,\lambda_n)$. Then $f_{\Omega_n}(\alpha)=\lambda_n$ and hence $\inf_{\tilde{\Lambda}}\{f_{\Omega_n}(\alpha)\mid n\in \N\} = f_{\cV}(\alpha)$.

(4) If $\alpha\in \Delta$ (resp. $\alpha\in \Phi$), for all $\lambda\in \Lambda$, we have $D(\alpha,\lambda)\in \cV$ if and only if $D(\alpha,\lambda)\in \cl^{\Delta}(\cV)$ (resp. $D(\alpha,\lambda)\in \cl(\cV)$). Therefore $F_{\cV}(\alpha)=F_{\cl^{\Delta}(\cV)}(\alpha)$ (resp. $F_{\cl(\cV)}(\alpha)=F_{\cV}(\alpha)$, which implies (4). 

(5) This follows from (1). 

(6) Assume that $f_{\cV}(\alpha)=-\infty$. Then $F_{\cV}(\alpha)=\Lambda$. Let $\lambda\in \Lambda$. If $D(-\alpha,-\lambda)\in \cV$, then  as $D(\alpha,\lambda-1)\in \cV$, we have $D(\alpha,\lambda-1)\cap D(-\alpha,-\lambda)=\emptyset\in \cV$: a contradiction. Therefore $D(\alpha,\lambda)\notin \cV$ and thus $F_{\cV}(\alpha)=\emptyset$. This proves that $f_{\cV}(-\alpha)=+\infty$. If $f_{\cV}(\alpha)= +\infty$ or $f_{\cV}(-\alpha)=+\infty$ and $f_{\cV}(\alpha),f_{\cV}(-\alpha)\neq -\infty$, the result is clear. Assume now $f_{\cV}(\alpha),f_{\cV}(-\alpha)\neq \pm \infty$. Then by (3), there exists $(\Omega_n)\in \cV^\N$ such that $(f_{\Omega_n}(\alpha))$ and $(f_{\Omega_n}(-\alpha))$ decrease and  $\inf_{\tilde{\Lambda}}f_{\Omega_n}(\alpha)=f_{\cV}(\alpha)$ and $\inf_{\tilde{\Lambda}}(f_{\Omega_n}(-\alpha))=f_{\cV}(-\alpha)$. By (5), we have $f_{\Omega_n}(\alpha)+f_{\Omega_n}(-\alpha)\geq 0$ for all $n\in \N$ and thus $f_{\cV}(\alpha)+f_{\cV}(-\alpha)\geq 0$. 
 \end{proof}

\subsection{Pro-unipotent subgroups of $U^{ma}$ associated with filters on $\A$}\label{ss_unip_sbgrps}

In this subsection, we introduce the subgroup $U_\cV^{ma\epsilon}\subset U^{ma\epsilon}$, for $\epsilon\in \{-,+\}$ and $\cV$ a filter on  $\A$. This subgroup is important because its intersections with $U^\epsilon$ will be the fixator of $\cV$ in $U^\epsilon$. 

We give a definition of $U_{\cV}^{ma\epsilon}$ using the completed enveloping algebra $\widehat{\cU^\epsilon}$. We then give a description of $U_{\cV}^{ma\epsilon}$ involving the $X_\alpha$ and the valuation $\omega$ (see Proposition~\ref{p_Rou4.5}), which is crucial for our purpose. We then study the action of $N$ on the $U_\cV^{ma}$ by conjugation.

\subsubsection{Submodules of the enveloping algebras associated with filters on $\A$}

Recall that we have $\cU=\bigoplus_{\alpha\in Q} \cU_\alpha$.

\begin{Definition}
For $\alpha\in Q$ and $\Omega$  a non-empty subset of $\A$, we set $\cU_{\alpha,\Omega}=\cU_\alpha\otimes \cF_{\geq f_\Omega(\alpha)}$\index[notation]{u@$\cU_{\alpha,\Omega}$}.  We define
$\cU_{\Omega}=\bigoplus_{\alpha\in Q} \cU_{\alpha,\Omega}$. We also set $\widehat{\cU}^+_\Omega =\prod_{\alpha\in Q_+} \cU_{\alpha,\Omega}\subset \widehat{\cU}^+_\cF$ and  $\widehat{\cU}^-_\Omega=\bigoplus_{\alpha\in Q_-} \cU_{\alpha,\Omega}\subset \widehat{\cU}^-_\cF$.  

Let now  $\cV$ be  a filter on $\A$. We set $\cU_\cV=\bigcup_{\Omega\in \cV} \cU_\Omega$\index[notation]{u@$\cU_{\cV}$}, $\widehat{\cU}^+_\cV=\bigcup_{\Omega\in \cV}\widehat{\cU}^+_\Omega$\index[notation]{u@$\widehat{\cU}^+_{\cV}$} and  $\widehat{\cU}^-_\cV=\bigcup_{\Omega\in \cV}\widehat{\cU}^-_\Omega$\index[notation]{u@$\widehat{\cU}^-_{\cV}$}. By concavity of $f_\cV$ (\eqref{e_concavity}), these are subalgebras of $\cU_\cF$, $\widehat{\cU}^+_\cF$ and $\widehat{\cU}^-_\cF$.
\end{Definition}

Note that if $\alpha\in Q$ and $\cV$ is a filter on $\A$, we have $\cU_{\alpha,\cV}=\cU_\alpha \otimes \cF_{\geq f_{\cV}(\alpha)}$. However in general, $\widehat{\cU}^+_\cV$ is different from $\prod_{\alpha\in Q_+}\cU_{\alpha,\cV}$, by Lemma~\ref{l_dif_descri} below.

\subsubsection{Definition of pro-unipotent subgroups}

\begin{Definition}\label{d_U_Vma}

Let $\Psi\subset \Delta_+$ be a closed subset. Let $\Omega$ be a non-empty subset of $\A$.  We set \[U^{ma}_{\cV}(\Psi)=\fU^{ma}_{\Psi}(\cF)\cap \widehat{\cU}^+_{\cV}.\]\index[notation]{u@$U^{ma}_{\cV}(\Psi), U^{ma+}_{\Omega}$} We write $U^{ma+}_\cV$ instead of $U^{ma}_{\cV}(\Delta_+)$.

If $\Psi\subset \Delta_-$ and $\cV$ is a filter on $\A$, then we define $U^{ma}_\cV(\Psi)\subset U^{ma-}$\index[notation]{u@$U^{ma}_{\cV}(\Psi), U^{ma-}_{\cV}$} and $U^{ma-}_{\cV}\subset U^{ma-}$ similarly.
\end{Definition}

The definition of $U^{ma}_\cV(\Psi)$ we gave has the advantage of being intrinsic: it does not depend on the choice of an order on $\Delta_+$ nor on the choice of bases of $\ffg_{\alpha,\Z}$, for $\alpha\in \Delta_+$. However, for our purpose, it will be more convenient to use the description of proposition below.

\begin{Proposition}(see \cite[4.5]{rousseau2016groupes})\label{p_Rou4.5}
\begin{enumerate}
\item Let $\cV$ be a filter on $\A$ and $\Psi$ be a closed subset of $\Delta_+$. Then $U^{ma}_{\cV}(\Psi)$ is a subgroup of $U^{ma}(\Psi)$.

\item Let $\Omega$ be a non-empty subset of $\A$. Choose any order on $\Delta_+$ compatible with the height. Then: \[U^{ma}_\Omega(\Psi)=\{\prod_{\alpha\in \Psi} X_\alpha(\underline{u_\alpha})\mid (\underline{u_\alpha})\in \prod_{\alpha\in \Psi}\ffg_{\alpha,\Z}\otimes \cF_{\geq f_{\Omega}(\alpha)}\}.\]
\end{enumerate}

\end{Proposition}

\begin{proof}
(1) Let $\Omega\in \cV$. By concavity of $f_{\Omega}$ (\eqref{e_concavity}), $U^{ma}_{\Omega}(\Psi)$ is stable by product. Let $x=\sum_{\alpha\in Q_+} u_\alpha\in U^{ma}_{\Omega}(\Psi)$. Then by Definition/Proposition~\ref{dp_grp_lk}, $x^{-1}=\tau(x)=\sum_{\alpha\in Q_+}\tau(u_\alpha)$. By \eqref{e_hopf_h} to \eqref{e_hopf_f},  $\tau(\cU_{\alpha,\Omega})=\cU_{\alpha,\Omega}$ for all $\alpha\in Q_+$. Therefore $\tau(x)=x^{-1}\in \widehat{\cU}_{\Omega}$ and hence $x^{-1}\in U_{\Omega}^{ma}(\Psi)$, which proves that $U_{\Omega}^{ma}(\Psi)$ is  a subgroup of $U^{ma}(\Psi)$.

  Let now $u_1,u_2\in U_{\cV}^{ma}(\Psi)$. Then there exist $\Omega_1,\Omega_2\in \cV$ such that $u_i\in U_{\Omega_i}^{ma}(\Psi)$, for both $i\in \{1,2\}$. Then $u_1u_2^{-1}\in U_{\Omega_1\cap \Omega_2}^{ma}(\Psi)\subset U_{\cV}^{ma}(\Psi)$, which proves that $U_{\cV}^{ma}(\Psi)$ is a subgroup of $U^{ma}(\Psi)$.

(2) Let $\Omega$ be a non-empty subset of $\A$. Let $(\underline{u_\alpha})_{\alpha\in \Psi}\in \prod_{\alpha\in \Psi} \ffg_{\alpha,\Z}\otimes \cF_{\geq f_{\Omega}(\alpha)}$. Let $\alpha\in \Psi$. Write $\underline{u_\alpha}=\sum_{b\in \cB_\alpha} \lambda_b b$, where $(\lambda_b)\in (\cF_{\geq f_{\Omega}(\alpha)})^{\cB_{\alpha}}$. Then $[\exp](\lambda_b b)=\sum_{n\in \N} \lambda_b^n b^{[n]}$. For $n\in \N$, $\omega(\lambda_b)+\alpha(\Omega)\geq 0$ and hence $\omega(\lambda_b^n)+n\alpha(\Omega)\geq 0$, which proves that $\lambda_b^n b^{[n]}\in \cU_{n\alpha,\Omega}$. Therefore $[\exp](\lambda_b b)\in U_{\Omega}^{ma}(\Psi)$, for $b\in \cB_{\alpha}$ and hence $X_{\alpha}(\underline{u_\alpha})\in U_{\Omega}^{ma}(\Psi)$, by (1). Let $\alpha\in Q_+$. The $\cU_{\alpha}$-component of $\prod_{\beta\in \Psi}X_\beta(\underline{u_\beta})$ is the $\cU_{\alpha}$-component of $\prod_{\beta\in F} X_\beta(\underline{u_\beta})$, where $F=\{\beta\in \Psi\mid \beta\leq_{Q^\vee} \alpha\}$. By (1), this proves that it belongs to $\cU_{\alpha,\Omega}$ and thus $\prod_{\beta\in \Psi}X_\beta(\underline{u_\beta})\in U_{\Omega}^{ma}(\Psi)$. Consequently, $\{\prod_{\alpha\in \Psi} X_\alpha(\underline{u_\alpha})\mid (\underline{u_\alpha})\in \prod_{\alpha\in \Psi}\ffg_{\alpha,\Z}\otimes \cF_{\geq f_{\Omega}(\alpha)}\}\subset U^{ma}_\Omega(\Psi)$.

Let now $(\underline{u_\alpha})\in (\prod_{\alpha\in\Psi} \ffg_{\alpha,\Z}\otimes \cF)\setminus \prod_{\alpha\in\Psi} \ffg_{\alpha,\Z}\otimes \cF_{\geq f_{\Omega}(\alpha)}$. Let us prove that $\prod_{\alpha\in \Psi} X_{\alpha}(\underline{u_\alpha})\in U^{ma}(\Psi)\setminus U_{\Omega}^{ma}(\Psi)$. Let $\alpha\in \Psi$ be of minimal height such that $\underline{u_\alpha}\notin \ffg_{\alpha,\Z}\otimes \cF_{\geq f_{\Omega}(\alpha)}$. By what we proved above, it suffices to prove that $u:=\prod_{\beta\in \Psi, \htt(\beta)\geq \htt(\alpha)} X_{\beta}(\underline{u_\beta})\notin U_{\Omega}^{ma}(\Psi)$. As the order on $\Psi$ is compatible with the height, the $\cU_{\alpha}$-component of $u$ is the $\cU_{\alpha}$-component of $X_\alpha(\underline{u_\alpha})$. Write $\underline{u_\alpha}=\sum_{b\in \cB_{\alpha}} \lambda_b b$, with $(\lambda_b)\in \cF^{\cB_{\alpha}}$. By \ref{a_es1}, the $\cU_{\alpha}$-component of $X_{\alpha}(\underline{u_\alpha})$ is $\sum_{b\in \cB_{\alpha}}\lambda_b b\in \cU_{\alpha}\setminus \cU_{\alpha,\Omega}$. In particular $u\notin U^{ma}_\Omega(\Psi)$. Thus we proved that $U^{ma}_{\Omega}(\Psi)\subset \{\prod_{\alpha\in \Psi} X_\alpha(\underline{u_\alpha})\mid (\underline{u_\alpha})\in \prod_{\alpha\in \Psi}\ffg_{\alpha,\Z}\otimes \cF_{\geq f_{\Omega}(\alpha)}\}$. (2) follows.
\end{proof}

If $\cV$ is a filter on $\A$ and $C$ is  a non-empty subset of $\A$, we set $\cV+C=\{\Omega+C\mid \Omega\in \Omega\in \cV\}$.

\begin{Corollary}\label{c_Uma_cl}
Let $\cV$ be a filter on $\A$ and $\Psi$ be a closed subset of $\Delta_+$. Then $U^{ma}_{\cV}(\Psi)=U^{ma}_{\cl^\Delta(\cV)}(\Psi)=U^{ma}_{\cl^\Delta(\cV)+\overline{C^v_f}}(\Psi)$. 
\end{Corollary}

\begin{proof}
We have $\cV\subset \cl^{\Delta}(\cV)$ and thus $U_{\cl^{\Delta}(\cV)}^{ma}(\Psi)\subset U_{\cV}^{ma}(\Psi)$. Let $u\in U_{\cV}^{ma}(\Psi)$. Let $\Omega\in \cV$ be such that $u\in U_{\Omega}^{ma}(\Psi)$. Write $u=\prod_{\alpha\in \Delta_+}X_\alpha(\underline{u_\alpha})$, with $(\underline{u_\alpha})\in \prod_{\alpha\in \Delta_+}\ffg_{\alpha,\Z}\otimes \cF_{\geq f_{\Omega}(\alpha)}$, which is possible by Proposition~\ref{p_Rou4.5}. Set $\Omega'=\bigcap_{\alpha\in \Delta_+}D(\alpha,\omega(\underline{u_\alpha}))$. Then $\Omega'\in \cl^{\Delta}(\cV)$ and $u\in U_{\Omega'}^{ma}(\Psi)$. Therefore $u\in U_{\cl^{\Delta}(\cV)}$ and hence $U_{\cV}^{ma}(\Psi)\subset U_{\cl^{\Delta}(\cV)}^{ma}(\Psi)$. Consequently, $U_{\cV}^{ma}(\Psi)= U_{\cl^{\Delta}(\cV)}^{ma}(\Psi)$.

Up to replacing $\cV$ by $\cl^{\Delta}(\cV)$, it suffices to prove that $U_{\cV}^{ma}(\Psi)=U_{\cV+\overline{C^v_f}}^{ma}(\Psi)$. As $\cV\subset \cV+\overline{C^v_f}$, we have $U_{\cV+\overline{C^v_f}}^{ma}(\Psi)\subset U_{\cV}^{ma}(\Psi)$. Let $u\in U_{\cV}^{ma}(\Psi)$. Let $\Omega\in \cV$ be such that $u\in U_{\Omega}^{ma}(\Psi)$. By definition, we have $f_{\Omega}(\alpha)=f_{\Omega+\overline{C^v_f}}(\alpha)$, for all $\alpha\in \Delta_+$. Therefore we have $U_{\Omega+\overline{C^v_f}}^{ma}(\Psi)=U_{\Omega}^{ma}(\Psi)$, by Proposition~\ref{p_Rou4.5}. As $\Omega+\overline{C^v_f}\in \cV+\overline{C^v_f}$, we deduce $u \in U_{\cV+\overline{C^v_f}}^{ma}(\Psi)$, and the corollary follows.
\end{proof}

\begin{Corollary}\label{c_U_ma_inter}
Let $\Omega$ be a non-empty subset of $\A$. Let $\Psi$ be a closed subset of $\Delta_+$.  Then $U_{\Omega}^{ma}(\Psi)=\bigcap_{x\in \Omega} U_x^{ma}(\Psi)$. 
\end{Corollary}

\begin{proof}
This follows from Proposition~\ref{p_Rou4.5} (2) and \eqref{e_normal_form_Upma}.
\end{proof}

\begin{Lemma}\label{l_dif_descri}
Let $\cV$ be a filter on $\A$ and $\Psi$ be a closed subset of $\Delta_+$. Then: \[U^{ma}_\cV(\Psi)\subset\{\prod_{\alpha\in \Psi} X_\alpha(\underline{u_\alpha})\mid (\underline{u_\alpha})\in \prod_{\alpha\in \Psi}\ffg_{\alpha,\Z}\otimes \cF_{\geq f_{\cV}(\alpha)}\}.\] However, if $A$ is not a Cartan matrix, then there exists a filter $\cV$ on $\A$ and $\Psi$ a closed subset of $\Delta_+$ such that the inclusion is strict.
\end{Lemma}

\begin{proof}

Assume that $\Delta_{im}$ is non-empty and choose $\alpha\in \Delta_{im}^+$. Let $C_0^-=\germ_0(-C^v_f)$.  
For $\beta\in Q_+$, define $\pr_{\beta}:\widehat{\cU}^+_{\cF}\rightarrow \cU_{\beta,\cF}$ by $\pr_{\beta}(v)=v_\beta$, for $v=\sum_{\gamma\in Q_+}v_{\gamma}\in \widehat{\cU}^+_{\cF}$. Choose $\qp\in \cF$ such that $\omega(\qp)=1$ and choose $u_{k\alpha}\in \ffg_{k\alpha,\Z}\setminus \{0\}$, for all $k\in \Z_{\geq 1}$. Set $u=\prod_{k\in \Z_{\geq 1}}X_{k\alpha}(\underline{u_{k\alpha}})$. Then by definition of $[\exp](\qp u_{k\alpha})$ and \ref{a_es1}, we have: \[\pr_{nk\alpha}(X_{k\alpha}(\qp u_{k\alpha}))=\pr_{nk\alpha}([\exp](\qp u_{k\alpha}))=\qp^{n}u_{k\alpha}^{[n]},\] for $k,n\in \Z_{\geq 1}$, where $(u_{k\alpha}^{[n]})_{n\in \N}$ is an exponential sequence for $u_{k\alpha}$. Therefore if $k\in \Z_{\geq 1}$, there exists $u_{k\alpha}'\in \cU_{k\alpha}$ satisfying  $\omega(u_{k\alpha}')>1$ and  such that:  \begin{equation}\label{e_in_val}
\pr_{k\alpha}(u)=\qp u_{k\alpha}^{[1]}+u_{k\alpha}'=\qp u_{k\alpha}+u_{k\alpha}'.
\end{equation} Therefore by Example~\ref{ex_f_cV}, $u\in  \{\prod_{\beta\in \Delta_+} X_\beta(\underline{u_\beta})\mid (\underline{u_\beta})\in \prod_{\beta\in \Delta_+}\ffg_{\beta,\Z}\otimes \cF_{\geq f_{C_0^-}(\beta)}\}$. 

However, $u\notin U_{C_0^-}^{ma+}$. Indeed, let $v\in U_{C_0^-}^{ma+}$. Consider $\Omega\in C_0^-$ such that $v\in U_{\Omega}^{ma+}$.  Write $v=\prod_{\beta\in \Delta_+} X_{\beta}(\underline{v_\beta})$, with $\underline{v_\beta}\in \prod_{\beta\in \Delta_+} \ffg_{\beta,\Z}\otimes \cF_{\geq \Omega}$.  Let $a\in -C^v_f\cap \Omega$ be such that $(\alpha_i(a))_{i\in I}$ has constant value, say $-\epsilon$.  Then for all $\beta\in \Delta_+$, we have $\omega(\underline{v_{\beta}})+\beta(a)\geq 0$ and thus $\omega(\underline{v_\beta})\geq \htt(\beta) \epsilon$, for $\beta\in \Delta_+$. Using \eqref{e_in_val}, we deduce that $u\notin U_{C_0^-}^{ma+}$, which proves the lemma.
\end{proof}

\subsubsection{Conjugation of the unipotent subgroups by $N$}

Let $\Psi\subset \Delta$. We assume that there exists $w\in W^v$ such that $w.\Psi\subset \Delta_+$ and $w.\Psi$ is closed. Let $\Omega$ be a non-empty subset of $\A$. 
We set: \begin{equation}\label{e_U_cV_wPsi} U_{\Omega}^{ma}(\Psi)=\{\prod_{\alpha\in \Psi} X_{\alpha}(\underline{u_\alpha})\mid (\underline{u_\alpha})\in \prod_{\alpha\in \Psi}\ffg_{\alpha,\Z}\otimes \cF_{\geq f_{\Omega}(\alpha)}\}\subset \fU^{ma}_{\Psi}(\cF).
\end{equation} If $\cV$ is a filter on $\A$, we set: \[U_\cV^{ma}(\Psi)=\bigcup_{\Omega\in \cV}U_\Omega^{ma}(\Psi).\]  This is a subgroup of $U_{w^{-1}.\Delta_+}^{ma}(\cF)$ and it is independent of the choice of the order (compatible with the height) on $\Delta_+$ by Proposition~\ref{p_cnj_N} below.

\begin{Definition}\label{d_U_alpha_cV}
Let $\alpha\in \Phi$ and $\cV$ be a filter on $\A$. Then $\{\alpha\}$ is a closed subset of $\Delta_+$.  We set $U_{\alpha,\cV}=U^{ma}_{\cV}(\{\alpha\})$\index[notation]{u@$U_{\alpha,\cV}$}.

For $k\in \tilde{\R}$, we set $U_{\alpha,k}=x_{\alpha}(\ffg_{\alpha,\Z}\otimes \cF_{\geq k})$.  By Proposition~\ref{p_Rou4.5}, if $\Omega\in \cV$ we have $U_{\alpha,\Omega}=U_{\alpha,f_{\Omega}(\alpha)}$. We have $U_{\alpha,\cV}=\bigcup_{\Omega\in \cV}U_{\alpha,f_{\Omega}(\alpha)}=U_{\alpha,f_{\cV}(\alpha)}$. 
\end{Definition}

\begin{Lemma}\label{l_cnj_N_x}
Let $\cV$ be a filter on $\A$, $\alpha\in\Phi$ and  $n\in N$. Set $w=\nu^v(n)$ and $\bw=\nu(n)$.  Then $n U_{\alpha,\cV} n^{-1}=U_{w.\alpha,\bw.\cV}$. 
\end{Lemma}

\begin{proof}
Let $\tilde{n}\in \fN_0(\cF)$. Then by \ref{a_KMT7}, we have $\tilde{n} U_{\alpha,\cV} \tilde{n}^{-1}=U_{\nu^v(\tilde{n}),\nu(\tilde{n})(\cV)}$.

By \ref{a_KMT4'}, if $t\in T$, then $tU_{\alpha,\cV} t^{-1}=U_{\alpha,\nu(t).\cV}$. Using \eqref{e_dec_N}, we deduce the result.
\end{proof}

\begin{Proposition}\label{p_cnj_N}
Let $n\in N$, $\cV$ be a  filter on $\A$ and $\Psi\subset \Delta_+$ be a closed set. Let $w=\nu^v(n)\in W^v$ and $\bw=\nu(n)\in W^a$. Then $nU_{\cV}^{ma}(\Psi)n^{-1}=U_{\bw.\cV}^{ma}(w.\Psi)$.
\end{Proposition}

\begin{proof}
Let $\Omega\in \cV$.  We first consider the case where $n=t\in T$. Let $\alpha\in \Delta_+$. Let $\underline{u_\alpha}\in \ffg_{\alpha,\Z}\otimes \cF_{\geq f_{\Omega}(\alpha)}$. Then by \eqref{e_com_T_X}, $tX_\alpha(\underline{u_\alpha})t^{-1}=X_{\alpha}(\alpha(t)\underline{u_\alpha})$ and \[\alpha(t)\underline{u_\alpha}\in \ffg_{\alpha,\Z}\otimes \cF_{\geq f_{\Omega}(\alpha)+\omega(\alpha(t))}=\ffg_{\alpha,\Z}\otimes \cF_{\geq f_{\nu(t).\Omega}(\alpha)}.\] By Proposition~\ref{p_Rou4.5}, we deduce that $tU_{\Omega}^{ma}(\Psi)t^{-1}=U_{\nu(t).\Omega}^{ma}(\Psi)$.

Let $w\in W^v$ and $k=\ell(w)$. Write $w=r_{i_1}\ldots r_{i_k}$, with $i_1,\ldots,i_k\in I$. Let $n_0=\tilde{s_{i_1}}\ldots \tilde{s_{i_k}}\in \fN_0(\cF)$. Let $\alpha\in \Delta_+$. 

If $\alpha\in \Phi$, we have $n_0 U_{\alpha,\Omega} n_0^{-1}=U_{w.\alpha, \bw.\Omega}$, by Lemma~\ref{l_cnj_N_x}.

Assume now $\alpha\in \Delta^{im}_+$.  By \eqref{e_Kmr_1.3.14}, $v.\alpha\in \Delta_+$ for all $\alpha\in \Delta_+$ and $v\in W^v$. Let $w^*=s_{i_1}^*\ldots s_{i_k}^*\in \mathrm{GL}(\ffg_\Z)$, with the notation of \eqref{e_si*}. Let $\cB_{w.\alpha}'=w^*.\cB_{\alpha}$. Then $\cB_{w.\alpha}'$ is a $\Z$-basis of $\ffg_{w.\alpha,\Z}$. Denote by $X_{w.\alpha}'$ the analog of $X_{w.\alpha}$ associated with $\cB'_\alpha$ instead of $\cB_{\alpha}$.   By \cite[Lemma 8.77]{marquis2018introduction} (either applied in $\fG^{\sch+}(\cF)\simeq \fG^{ma+}(\cF)$ or applied directly in $\fG^{ma+}(\cF)$, with the same proof), we have:  \[n_0.X_\alpha(\underline{u_\alpha})n_0^{-1}=X_{w.\alpha}'(w^* \underline{u_\alpha}),\] for $\underline{u_{\alpha}}\in \ffg_{\alpha,\Z}\otimes \cF$. If $\underline{u_\alpha}\in \ffg_{\alpha,\Z}\otimes \cF_{\geq f_{\Omega}(\alpha)}$, then $w^*.\underline{u_\alpha}\in \ffg_{w.\alpha,\Z}\otimes \cF_{\geq f_{\Omega}(\alpha)}=\ffg_{w.\alpha,\Z}\otimes \cF_{\geq f_{w.\Omega}(w.\alpha)}$. Therefore: \[n_0 X_{\alpha}(\ffg_{\alpha,\Z}\otimes \cF_{\geq f_{\Omega}(\alpha)})n_0^{-1}=X_{w.\alpha}'(\ffg_{w.\alpha,\Z}\otimes \cF_{\geq f_{w.\Omega}(w.\alpha)}).\]

Using Proposition~\ref{p_Rou4.5}, we deduce that $n_0 U_{\Omega}^{ma}(\Psi)n_0^{-1}=U_{w.\Omega}^{ma}(w.\Psi)$.
 
  By \eqref{e_dec_N}, we can write every element of $n$ as $n=n_0t$, with $t\in T$ and $n_0\in \fN_0(\cF)$ and thus $nU_{\Omega}^{ma}(\Psi)=U_{\bw.\cV}^{ma}(w.\Psi)$. As $U_{\cV}^{ma}(\Psi)=\bigcup_{\Omega\in \cV} U_{\Omega}^{ma}(\Psi)$, we deduce the result.
\end{proof}

\subsection{Unipotent and pro-unipotent subgroups of $U^+$ and $U^-$ associated with filters on $\A$}

In this subsection, we introduce the subgroups $U_{\cV,\fin}^\epsilon$ and $U_{\cV,\infty}^\epsilon$, for $\cV$ a filter on $\A$ and $\epsilon\in \{-,+\}$. This last subgroup will be the $U^\epsilon$ component of the fixator of $\cV$ for the action of $G$ on the masure.

\subsubsection{Unipotent subgroups $U_{\cV,\fin}$}

\begin{Definition}\label{d_U_N_cV} Let $\cV$ be a filter on $\A$.  For $\alpha\in \Phi$, we set $U_{\cV}^{(\alpha)}=\langle U_{\alpha,\cV},U_{-\alpha,\cV}\rangle\subset G$\index[notation]{u@$U_{\cV}^{(\alpha)}$}. We define $U_{\cV,\fin}=\langle U_{\alpha,\cV}\mid \alpha\in\Phi\rangle$\index[notation]{u@$U_{\cV,\fin}$}  and $U_{\cV,\fin}^{\epsilon}=U_{\cV,\fin}\cap U^{\epsilon}$\index[notation]{u@$U_{\cV,\fin}^+,U_{\cV,\fin}^-$}, for $\epsilon\in \{-,+\}$.

 We set $N_{\cV}^{(\alpha)}=N\cap U_{\cV}^{(\alpha)}$\index[notation]{n@$N_{\cV}^{(\alpha)}$} and   $N_\cV^u=\langle N_{\cV}^{(\alpha)}\mid \alpha\in \Phi\rangle$\index[notation]{n@$N_\cV^u$}.

\end{Definition}

Let $\cV,\cV'$ be two filters on $\A$.  The groups $U_{\cV}$, $U_{\cV}^{(\alpha)}$, $N_\cV^{(\alpha)}$ and  $N_{\cV}^u$ are normalized by $T_0$, by \eqref{e_KMT4'}.   If $\cV'\subset \cV$, we have $U_{\cV',\fin}\supset U_{\cV,\fin}$ etc.

By Lemma~\ref{l_cnj_N_x}, if $\cV$ is a filter on $\A$, $\alpha\in\Phi$ and  $n\in N$, we have: \begin{equation}\label{e_cnj_U_fin}
n U_{\cV,\fin} n^{-1}=U_{\bw.\cV,\fin},
\end{equation}
where $w=\nu^v(n)$ and $\bw=\nu(n)$.

If $\cV$ is a filter on $\A$ and $\alpha\in \Phi$, we define $\cl^{(\alpha)}(\cV)$ as the filter on $\A$ generated by all the sets containing $\cV$ of the form $D(\alpha,k_\alpha)\cap D(-\alpha,k_\alpha')$, where $k_\alpha,k_\alpha'\in \Lambda\cup \{+\infty\}$.

\begin{Lemma}\label{l_cnj_N_U}
Let $\cV$ be a filter on $\A$ and $\alpha\in \Phi$. Then $U_{\cV,\fin}=U_{\cl(\cV),\fin}, N_\cV^u=N_{\cl(\cV)}^u$, $U_{\cV}^{(\alpha)}=U_{\cl^{(\alpha)}(\cV)}^{(\alpha)}$ and $N_{\cV}^{(\alpha)}=N_{\cl^{(\alpha)}}^{(\alpha)}(\cV)$. 
\end{Lemma}

\begin{proof}
As $\cl(\cV)\supset \cV$, we have $U_{\cl(\cV),\fin}\subset U_{\cV,\fin}$. Now if $u\in U_{\cV,\fin}$, we can write $u=\prod_{i=1}^n x_{\beta_i}(a_i)$, where $n\in \N$, $\beta_1,\ldots,\beta_n\in \Phi$ and $a_1,\ldots,a_n\in \cF$ are such that $x_{\beta_i}(a_i)\in U_{\beta_i,\cV}$, for $i\in \llbracket 1,n\rrbracket$. Set $\Omega=\bigcap_{i=1}^n D(\beta_i,\omega(a_i))$. Then $u\in U_{\Omega,\fin}$ and $\Omega\in \cl(\cV)$. Therefore $u\in U_{\cl(\cV),\fin}$. The other equalities are obtained similarly.
\end{proof}

\subsubsection{Study of  $U_{\cV}^{(\alpha)}$}

 In this subsubsection we use computations in $\mathrm{SL}_2$ to study the subgroup $U_{\cV}^{(\alpha)}$, for $\alpha\in \Phi$ and $\cV$ a filter on $\A$.  This study will be useful to obtain decompositions of the different subgroups $U$ associated with $\cV$, see Proposition~\ref{p_GR08_3.4}.

\begin{Lemma}\label{l_f_cv}
Let $\alpha\in \Phi$ and $\cV$ be a filter on $\A$. Then ($-\infty \notin \{f_{\cV}(\alpha),f_{\cV}(-\alpha)\}$ and $f_{\cV}(\alpha)+f_{\cV}(-\alpha)=0$) if and only if there exist $\mu\in \Lambda$ and $\Omega\in \cV$ such that $\alpha(\Omega)=\{\mu\}$.
\end{Lemma}

\begin{proof}
Assume that $-\infty\notin \{f_{\cV}(\alpha),f_{\cV}(-\alpha)\}$ and that  $f_{\cV}(\alpha)+f_{\cV}(-\alpha)=0$. Then $f_{\cV}(\alpha), f_{\cV}(-\alpha)\in \Lambda$ (since otherwise $f_{\cV}(\alpha)$ or $f_{\cV}(-\alpha)$ would belong to $\{r^+\mid r\in \R\}$). By Proposition~\ref{p_f_cV_alpha} (3), this implies the existence of $\Omega_\epsilon\in \cV$ such that $f_{\cV}(\epsilon\alpha)=f_{\Omega_\epsilon}(\alpha)$, for both $\epsilon\in \{-1,1\}$. Set $\Omega=\Omega_+\cap \Omega_-$. Then by Proposition~\ref{p_f_cV_alpha}, we have $f_{\cV}(\epsilon \alpha)=f_{\Omega}(\epsilon \alpha)$, for both $\epsilon\in \{-1,1\}$. 

By Proposition~\ref{p_f_cV_alpha}, $\sup_\R(-\alpha(\Omega)),\sup_\R(\alpha(\Omega))\in \Lambda$. If $\alpha(\Omega)$ was not reduced to a singleton, it would contain at least two elements $a,b$, with $a<b$ and we would have $\sup_\R(-\alpha(\Omega))-\sup_\R(\alpha(\Omega))\geq b-a>0$. Therefore $\alpha(\Omega)$ is a singleton of $\Lambda$. 
\end{proof}

Let $x_-,x_+:\cF\rightarrow \mathrm{SL}_2(\cF)$ be defined by $x_-(r)=\begin{psmallmatrix} 1 & 0\\ r & 1\end{psmallmatrix}$ and $x_+(r)=\begin{psmallmatrix} 1 & r\\ 0 & 1\end{psmallmatrix}$, for $r\in \cF$. Let $\overline{T_0}=\begin{psmallmatrix} \cO^\times & 0 \\ 0 & \cO^\times \end{psmallmatrix}\cap \mathrm{SL}_2(\cF)$. If $a\in \tilde{\Lambda}$ and $\epsilon\in \{-,+\}$, we set $\overline{U_{\epsilon,a}}=x_\epsilon(\cF_{\geq a})$.

For $k\in \Lambda$, we set $\overline{N_k}=\overline{T_0}\cup \overline{T_0}\cdot \begin{psmallmatrix} 0 & a\\ -a^{-1} & 0\end{psmallmatrix}\subset \mathrm{SL}_2(\cF)$, for any $a\in \cF$ such that $\omega(a)= k$.  This is a well-defined subgroup of $\mathrm{SL}_2(\cF)$.

\begin{Lemma}\label{l_cmp_SL2}
Let $k_+,k_-\in \tilde{\Lambda}$ be such that $k_++k_-\geq 0$. Let $\overline{G_{k_+,k_-}}=\langle \overline{U_{+,k_+}},\overline{U_{-,k_-}}\rangle\subset \mathrm{SL}_2(\cF)$.   Then we have: \begin{align*}
\overline{G_{k_+,k_-}}&=\overline{U_{+,k_+}}\cdot \overline{U_{-,k_-}}\cdot \overline{T_0}=\overline{U_{-,k_-}}\cdot\overline{U_{+,k_+}}\cdot  \overline{T_0} &\text{if }k_++k_->0\\
 \overline{G_{k_+,k_-}}&=\overline{U_{+,k_+}}\cdot \overline{U_{-,k_-}}\cdot \overline{N_{k_-}}=\overline{U_{-,k_-}}\cdot\overline{U_{+,k_+}}\cdot  \overline{N_{k_+}} &\text{if }k_++k_-=0.
\end{align*}
\end{Lemma}

\begin{proof}
Let $a,b\in \cF$ be such that $1+ab\neq 0$. By an easy computation, we have: \begin{equation}\label{e_com_x_+_x_-}
 x_-(b)x_+(a)=x_+(a/(1+ab))x_-(b(1+ab))\begin{psmallmatrix} 1/(1+ab)& 0\\ 0 & 1+ab\end{psmallmatrix}.
 \end{equation}

Let $a,b\in \cF$ be such that $x_+(a)\in \overline{U_{+,k_+}}$ and $x_-(b)\in \overline{U_{-,k_-}}$. Then $\omega(ab)\geq k_++k_-\geq 0$. First assume $k_++k_-> 0$. Then by \eqref{e_com_x_+_x_-}, we have: \begin{equation}\label{e_U_+_U_-}
\overline{U_{-,k_-}}\cdot \overline{U_{+,k_+}}\subset \overline{U_{+,k_+}}\cdot \overline{U_{-,k_-}}\cdot \overline{T_0}.
\end{equation}

Therefore: \begin{align*} \overline{U_{-,k_-}}\cdot \overline{T_0}\cdot \overline{U_{+,k_+}}\cdot \overline{U_{-,k_-}}\cdot \overline{T_0}&\subset  \overline{U_{-,k_-}}\cdot \overline{U_{+,k_+}}\cdot \overline{U_{-,k_-}}\cdot \overline{T_0}\cdot \overline{T_0}   &\text{by \eqref{e_KMT4'}}\\
&\subset \overline{U_{+,k_+}}\cdot \overline{U_{-,k_-}}\cdot \overline{T_0} &\text{ by \eqref{e_U_+_U_-}}.
\end{align*}

Therefore $\overline{U_{+,k_+}}\cdot \overline{U_{-,k_-}}\cdot \overline{T_0}$ is stable under left multiplication by $\overline{U_{-,k_-}}$ and $\overline{T_0}$. As it is obviously stable under left multiplication by $\overline{U_{+,k_+}}$, we deduce that it is a subgroup of $\mathrm{SL}_2(\cF)$ and thus it is equal to $\overline{G_{k_+,k_-}}$.

We now assume that $k_++k_-=0$. If $\omega(1+ab)=0$, then $x_-(b)x_+(a)\in \overline{U_{+,k_+}}\cdot \overline{U_{-,k_-}}\cdot \overline{T_0}$, by \eqref{e_com_x_+_x_-}. We now assume $\omega(1+ab)>0$. Write $ab=-1+c$, with $c\in \cF_{>0}$. Then we have $b=-a^{-1}+ca^{-1}$. We have:
\begin{equation}\label{e_df_r(a)}
x_+(a)x_-(-a^{-1})x_+(a)=\begin{psmallmatrix} 0 & a \\-a^{-1} & 0\end{psmallmatrix}=:r(a).
\end{equation}

Therefore $x_-(-a^{-1})x_+(a)=x_+(-a)r(a)$ and \[x_-(b)x_+(a)=x_-(ca^{-1})x_-(a^{-1})x_+(a)=x_-(ca^{-1})x_+(a)r(a)\in \overline{U_{+,k_+}}\cdot \overline{U_{-,k_-}}\cdot \overline{N_{k_+}}.\]

We deduce: \begin{equation}\label{e_U_-_U_+_0}
\overline{U_{-,k_-}}\cdot \overline{U_{+,k_+}}\subset \overline{U_{+,k_+}}\cdot \overline{U_{-,k_-}} \cdot \overline{N_{k_+}}.
\end{equation}

We have: \begin{equation}\label{e_cnj_r}
r(a) \overline{U_{-,k_-}} r(a)^{-1}=\begin{psmallmatrix} 1 & -a^2 \cF_{\geq k_-} \\  0 & 1\end{psmallmatrix}= \begin{psmallmatrix} 1 & \cF_{\geq k_-+2k_+} \\  0 & 1\end{psmallmatrix}=\overline{U_{+,k_+}}\end{equation} and similarly: \begin{equation}\label{e_cnj_r_+}
r(a) \overline{U_{+,k_+}}r(a)^{-1}=\overline{U_{-,k_-}}
\end{equation}

Using \eqref{e_U_-_U_+_0}, \eqref{e_cnj_r} and \eqref{e_cnj_r_+}, we get that $\overline{U_{+,k_+}}\cdot \overline{U_{-,k_-}}\cdot \overline{N_{k_+}}$ is stable by left multiplication by $\overline{U_{+,k_+}}$, $\overline{U_{-,k_-}}$ and $ \overline{N_{k_+}}$, which proves that it is a group. By symmetry of the roles of $+$ and $-$, we get the other equalities.
\end{proof}

Let $\alpha\in \Phi$. We define $r_\alpha:\A\rightarrow \A$ by $r_\alpha(x)=x-\alpha(x)\alpha^\vee$, for $x\in \A$. For $k\in \Lambda$, set $r_{\alpha,k}=r_\alpha \bt_{k\alpha^\vee}$, where $\bt_{k\alpha^\vee}$ is the translation on $\A$ by the vector $k\alpha^\vee$. Then $r_{\alpha,k}$ fixes $M(\alpha,k):=\{x\in \A\mid \alpha(x)+k=0\}$.

\begin{Lemma}\label{l_nu_m}
Let $a\in \cF$ and $\alpha\in \Phi$. Then with the notation of \eqref{e_m},  \[\nu\left(m(x_\alpha(a))\right)=\nu\left(x_{-\alpha}(a^{-1})x_{\alpha}(a)x_{-\alpha}(a)\right)=r_{\alpha,\omega(a)}.\]
\end{Lemma}

\begin{proof}
Write $\alpha=w.\alpha_i$, with $i\in I$ and $w\in W^v$. Then by \eqref{e_m_s}, there exists $\epsilon\in \{-,+\}$ such that $m(x_\alpha(a))=\tilde{w}\tilde{s_i}(\epsilon a)\tilde{w}^{-1}$. Then: \[\nu(m(x_\alpha(a)))=w \nu(\tilde{s_i}(\epsilon a))w^{-1}.\] 

By \ref{a_KMT6}, \[\nu(\tilde{s_i}(\epsilon a))=\nu(\tilde{s_i} \alpha_i^\vee(\epsilon a^{-1}))=r_i\nu(\alpha_i^\vee(\epsilon a^{-1}))=r_i \bt_{\omega(a)\alpha_i^\vee},\] by Lemma~\ref{l_act_cort}. Therefore $\nu(m(x_\alpha(a)))=r_{\alpha}\bt_{\omega(a)\alpha^\vee}=r_{\alpha,\omega(a)}$. 
\end{proof}

\begin{Lemma}\label{l_GR08_3.3}(\cite[Lemma 3.3]{gaussent2008kac})
Let $\cV$ be a filter on $\A$ and $\alpha\in \Phi$. Then: \begin{enumerate}
\item $U_{\cV}^{(\alpha)}=U_{\alpha,\cV}\cdot  U_{-\alpha,\cV} \cdot N_{\cV}^{(\alpha)}=U_{-\alpha,\cV}\cdot  U_{\alpha,\cV}\cdot  N_{\cV}^{(\alpha)}$.

\item If $f_{\cV}(\alpha)+f_{\cV}(-\alpha)>0$, then $N_{\cV}^{(\alpha)}\subset T_0$. If $f_{\cV}(\alpha)=-f_{\cV}(-\alpha)=:k\in \Lambda$, then $\nu(N_{\cV}^{(\alpha)})=\{r_{\alpha,k},1\}$.

\item $N_{\cV}^{(\alpha)}$ fixes $\cV$, i.e for all $n\in N_{\cV}^{(\alpha)}$, there exists $\Omega\in \cV$ fixed pointwise by $\nu(n)$.
\end{enumerate}
\end{Lemma}

\begin{proof}
We reduce the proof to computations in $\mathrm{SL}_2(\cF)$.  Write $\alpha=w.\alpha_i$, with $i\in I$ and $w\in W^v$. Let $k=\ell(w)$ and $i_1,\ldots,i_k\in I$ be such that $w=r_{i_1}\ldots r_{i_k}$. Let $\tilde{w}=\tilde{s}_{i_1}\ldots \tilde{s}_{i_k}\in N$.  Then by \ref{a_KMT7}, we have $U_{\cV}^{(\alpha)}=\tilde{w} \cdot U_{w^{-1}.\cV}^{(\alpha_i)}\cdot \tilde{w}^{-1}$. Thus, replacing $\cV$ by $w^{-1}.\cV$ if necessary, we can assume that $\alpha=\alpha_i$.

 By \cite[Exercise 7.51]{marquis2018introduction}, the assignment $\pi:\mathrm{SL}_2(\cF)\rightarrow \fG(\cF)$ defined by $\overline{x_{\pm }}(r)\mapsto x_{\pm \alpha}(r)$, for $r\in \cF$, defines a group morphism. We have $\langle U_{\alpha,\cV},U_{-\alpha,\cV}\rangle=\langle \pi(\overline{U_{+,f_\cV}(\alpha)}),\pi(\overline{U_{-,f_\cV(-\alpha)}})\rangle=\pi(\langle \overline{U_{+,f_\cV(\alpha)}},\overline{U_{-,f_\cV(-\alpha)}}\rangle)$. By Lemma~\ref{l_cmp_SL2}, we have: \[U_{\cV}^{(\alpha)}=\pi(\overline{U_{+,f_\cV(\alpha)}})\cdot \pi(\overline{U_{-,f_\cV(-\alpha)}})\cdot (N'_{\cV})^{(\alpha)}=U_{\alpha,\cV}\cdot U_{-\alpha,\cV}\cdot (N'_{\cV})^{(\alpha)},\] where $(N'_{\cV})^{(\alpha)}=\begin{cases}&\pi(\overline{T_0})\text{ if }f_{\cV}(\alpha)+f_{\cV}(-\alpha)>0\\ 
 &\pi(\overline{N_{f_{\cV}(\alpha)}})\text{ if }f_{\cV}(\alpha)+f_{\cV}(-\alpha)=0\end{cases}$.

If $f_\cV(\alpha)+f_{\cV}(-\alpha)=0$, choose $a\in \cF^\times$ such that $\omega(a)=f_\cV(\alpha)=:k$, which is possible by Lemma~\ref{l_f_cv}.  Then $\pi(\overline{T_0})\subset T_0$ and $\pi(\begin{psmallmatrix} 0 & a\\ -a^{-1} & 0\end{psmallmatrix})=\tilde{s}_i(a)\in N$, by \eqref{e_df_r(a)} and \eqref{e_tilde_si}. Therefore $(N'_\cV)^{(\alpha)}\subset N$. 

Let $g\in N\cap U_{\cV}^{(\alpha)}$. Write $g=u_+u_-n$, with $u_+\in U_{\alpha,\cV}$, $u_-\in U_{-\alpha,\cV}$ and $n\in (N'_\cV)^{(\alpha)}$. Then $u_+^{-1} gn^{-1} u_-^{-1}=1$ and thus by the Birkhoff decomposition, we have $u_+=u_-=1$ and $g=n$. Therefore \[N\cap U_{\cV}^{(\alpha)}=(N'_\cV)^{(\alpha)}=N_{\cV}^{(\alpha)}.\] We deduce the left hand side of (1) and by symmetry, we get (1).

As $\pi(\overline{T_0})\subset T_0$, we have  $N_\cV^{(\alpha)}\subset T_0$ if $f_{\cV}(\alpha)+f_{\cV}(-\alpha)>0$ and $\nu(N_\cV^{(\alpha)})=\{1, \nu(\tilde{s_i}(a))\}=\{1,r_{\alpha,k}\}$, if $f_{\cV}( \alpha)+f_{\cV}(-\alpha)=0$, by Lemma~\ref{l_nu_m}. This proves (2).

(3) If $f_{\cV}(\alpha)+f_{\cV}(-\alpha)>0$, then $N_{\cV}^{(\alpha)}\subset  T_0$ fixes $\A$ and thus it fixes $\cV$. If $f_{\cV}(\alpha)+f_{\cV}(-\alpha)=0$, then by Lemma~\ref{l_f_cv}, there exists $\Omega\in \cV$ such that $\alpha(\Omega)=\{k\}$. Then $r_{\alpha,k}$ fixes $\Omega$ and  $N^{(\alpha)}_{\cV}$ fixes  $\Omega$, which proves (3).
\end{proof}

\subsubsection{Pro-unipotent subgroups $U_{\cV,\infty}^+$ and $U_{\cV,\infty}^-$}

\begin{Definition}\label{d_U_infty_pm}
Let $\cV$ be a filter on $\A$ and $\epsilon\in \{-,+\}$.  Let $\Psi\subset \Delta_\epsilon$ be a closed subset. We set $U_{\cV,\infty}(\Psi)=U^{ma\epsilon}_{\cV}(\Psi)\cap U^\epsilon$\index[notation]{u@$U_{\cV,\infty}(\Psi)$}, (see Definition~\ref{d_U_Vma} for the definition of  $U^{ma}_{\cV}(\Psi)$). We write $U_{\cV,\infty}^+$\index[notation]{u@$U_{\cV,\infty}^{\pm}$} instead of $U_{\cV,\infty}(\Delta_+)$ and $U_{\cV,\infty}^-$ instead of $U_{\cV,\infty}(\Delta_-)$.

We also set $U_{\cV,\infty+}=\langle  U_{\cV,\fin},U_{\cV,\infty}^+\rangle\subset G$\index[notation]{u@$U_{\cV,\infty+}, U_{\cV,\infty-}$} and $U_{\cV,\infty-}=\langle U_{\cV,\fin},U_{\cV,\infty}^-\rangle\subset G$.

For $\cV$ a filter on $\A$ and $w\in W^v$, we set $U_{\cV,\infty}(w.\Delta_-)=U_{\cV}^{ma-}(w.\Delta_-)\cap G$\index[notation]{u@$U_{\cV,\infty}(w.\Delta_-),U_{\cV,\infty}(w.\Delta_+)$} and
 $U_{\cV,\infty}(w.\Delta_+)=U_{\cV}^{ma+}(w.\Delta_+)\cap G$.
\end{Definition}

When $G$ is reductive, then $U^{ma+}=U^+$ and $U_{\cV,\fin}^+=U_{\cV,\infty}^+=\langle U_{\alpha,\cV}\mid \alpha\in \Phi_+\rangle$, by \cite[7.1.8, 7.1.11 and 7.4.8]{bruhat1972groupes}. In particular, Proposition~\ref{p_Rou4.5} provides a description of $U_\cV^{ma}$ by the ``coordinates'' $X_\alpha$, which are simply the $x_\alpha$ in this case. In the Kac--Moody case, Proposition~\ref{p_Rou4.5} provides a description of $U_{\cV}^{ma+}$ by coordinates. However, we do not know ``coordinates'' which are adapted to describe $U_{\cV,\fin}^{+}$. In general, $U_{\cV,\fin}^+\neq \langle U_{\alpha,\cV}\mid \alpha\in \Phi_+\rangle$ (see \ref{sss_cmp_Ufnty}) and if $(\alpha,\beta)\in (\Phi_+)^2$ is not a prenilpotent pair, then $\langle U_\alpha,U_\beta\rangle$ is the free product $U_\alpha*U_\beta$, by Proposition~\ref{p_DRJ} and we do not know whether it is possible to have a normal form for the elements of $U_{\cV,\fin}^+$.

\begin{Lemma}\label{l_pre_GR08}
Let $\cV$ be a filter on $\A$. Then: \[\sU_{\cV,\infty}:=U_{\cV,\infty}^+\cdot U_{\cV,\infty}^-\cdot N^u_{\cV},\ \sU_{\cV,\fin}:=U_{\cV,\fin}^+\cdot U_{\cV,\fin}^-\cdot N^u_{\cV},\]\index[notation]{u@$\sU_{\cV,\infty}, \sU_{\cV,\fin}$} \[\sU_{\cV,\infty+}:=U_{\cV,\infty}^+\cdot U_{\cV,\fin}^-\cdot N^u_\cV\text{ and }\sU_{\cV,\infty-}:=U_{\cV,\infty}^-\cdot U_{\cV,\fin}^+\cdot N^u_\cV\]\index[notation]{u@$\sU_{\infty\pm}$} are independent of the choice of a set of positive roots in it $W^v$-class.
\end{Lemma}

\begin{proof}
Let $i\in I$. We have $r_i(\Delta_-)=(\Delta_-\setminus\{-\alpha_i\})\cup\{\alpha_i\}$ (by \eqref{e_Kmr_1.3.14}). Set $U_{\cV,\infty}(\Delta_-\setminus \{-\alpha_i\})=U^{ma-}_\cV(\Delta_-\setminus \{-\alpha_i\})\cap G$.  Then by  \eqref{e_smd_pdct_U} intersected with $G$ (and by symmetry), we have: \begin{equation}\label{e_dec_Unm}
U_{\cV,\infty}^-=U_{-\alpha_i,\cV}\ltimes U_{\cV,\infty}(\Delta_-\setminus \{-\alpha_i\})\text{ and }U_{\cV,\infty}^+=U_{\alpha_i,\cV}\ltimes U_{\cV,\infty}(\Delta_+\setminus \{\alpha_i\}).
\end{equation}

Let $U_{\cV,\infty}(r_i.\Delta_-)=U_{\cV}^{ma-}(r_i.\Delta_-)\cap G$. By Proposition~\ref{p_cnj_N}, we have $U_{\cV,\infty}(r_i.\Delta_-)=\tilde{r}_i U_{r_i.\cV,\infty}^- \tilde{r}_i$. Therefore by \eqref{e_dec_Unm} and Proposition~\ref{p_cnj_N}, we have $U_{\cV,\infty}(r_i.\Delta_-)=U_{\alpha_i,\cV}\ltimes U_{\cV,\infty}(\Delta_-\setminus \{\alpha_i\})$.

We have: \begin{align*}
\sU_{\cV,\infty}&=U_{\cV,\infty}(\Delta_+\setminus\{\alpha_i\})\cdot U_{\alpha_i,\cV}\cdot U_{\cV,\infty}(\Delta_-\setminus \{-\alpha_i\})\cdot U_{-\alpha_i,\cV} \cdot N_\cV^u\\
&= U_{\cV,\infty}(\Delta_+\setminus\{\alpha_i\})\cdot U_{\cV,\infty}(\Delta_-\setminus \{-\alpha_i\})\cdot U_{\alpha_i,\cV}\cdot U_{-\alpha_i,\cV} \cdot N_\cV^u \\
&=U_{\cV,\infty}(\Delta_+\setminus\{\alpha_i\})\cdot U_{\cV,\infty}(\Delta_-\setminus \{-\alpha_i\})\cdot U_{-\alpha_i,\cV}\cdot U_{\alpha_i,\cV} \cdot N_\cV^u  &\text{by Lemma~\ref{l_GR08_3.3}}\\
&=U_{\cV,\infty}(\Delta_+\setminus\{\alpha_i\})\cdot U_{-\alpha_i,\cV}\cdot U_{\cV,\infty}(\Delta_-\setminus \{-\alpha_i\})\cdot U_{\alpha_i,\cV} \cdot N_\cV^u \\
&=U_{\cV,\infty}(r_i.\Delta_+)\cdot U_{\cV,\infty}(r_i.\Delta_-)\cdot N_\cV^u &\text{by \eqref{e_dec_Unm}.}
\end{align*}

By induction on the length in $W^v$, using \eqref{e_cnj_U_Wv} and \eqref{e_smd_pdct_U},  we  deduce that $\sU_{\cV,\infty}$ is independent of the choice of $\Delta_+$. 

Let $\Psi\subset \Delta$ be such that $w.\Psi$ is a closed subset of $\Delta_+$, for some $w\in W^v$. We set $U_{\cV,\fin}(\Psi)=U_{\cV,\fin}\cap U^{ma+}_\Psi$.
 Let $u\in U_{\cV,\fin}^+$. Then  by \eqref{e_smd_pdct_U}, we can write $u=u_1u_2$, with $u_1\in U_{\alpha_i,\cV}$ and $u_2\in U_{\cV}^{ma+}(\Delta_+\setminus\{\alpha_i\})$. Then   $u_2=u_1^{-1}u\in U_{\cV}^{ma+}(\Delta_+\setminus \{\alpha_i\})$. Therefore $u_2\in U_{\cV,\fin}(\Delta_+\setminus \{\alpha_i\})$. Consequently, \begin{equation}
U_{\cV,\fin}^+=U_{\alpha_i,\cV}\ltimes U_{\cV,\fin}(\Delta_+\setminus \{\alpha_i\}).
\end{equation}

We also have $\tilde{r}_i U_{\cV,\fin}(r_i.\Delta_+)\tilde{r}_i^{-1}=U_{r_i.\cV,\fin}^+=U_{\alpha_i,\cV}\ltimes U_{\cV}(r_i.\Delta_+\setminus \{\alpha_i\})$.
Therefore $U_{\alpha_i,\cV}$ and $U_{-\alpha_i,\cV}$ normalize $U_{\cV,\fin}(\Delta_+\setminus \{\alpha_i\})$. Considering $\sU_{\cV,\fin}:=U_{\cV,\fin}^+\cdot U_{\cV,\fin}^-\cdot N^u_\cV$ instead of $\sU_{\cV,\infty}$, and with the same reasoning as above, we have  $\sU_{\cV,\fin}=U_{\cV,\fin}(w.\Delta_+)\cdot U_{\cV,\fin}(w.\Delta_-) \cdot N^u_\cV$, for every $w\in W^v$.

Similarly, we obtain that  $\sU_{\cV,\infty+}=U_{\cV,\infty}^+\cdot U_{\cV,\fin}^-\cdot N^u_\cV$ and $\sU_{\cV,\infty-}^-=U_{\cV,\infty}^-\cdot U_{\cV,\fin}^+\cdot N^u_\cV$ are independent of the choice of a fundamental chamber.
\end{proof}

\begin{Proposition}\label{p_GR08_3.4}(see \cite[Proposition 3.4]{gaussent2008kac}) Let $\cV$ be a filter on $\A$.

Let $\epsilon\in \{-,+\}$. We have $U_{\cV,\fin}^\epsilon\subset U_{\cV,\infty}^\epsilon$ and the following decompositions hold, independently of the choice of a set of positive roots in its $W^v$-conjugacy class : \begin{align*}
U_{\cV,\fin}&=U_{\cV,\fin}^- \cdot U_{\cV,\fin}^+\cdot N_\cV^u=U_{\cV,\fin}^+\cdot U_{\cV,\fin}^-\cdot N_\cV^u\\
U_{\cV,\infty\epsilon}&=U_{\cV,\infty}^\epsilon\cdot  U_{\cV,\fin}^{-\epsilon}\cdot N_{\cV}^u.
\end{align*}

Moreover, we have:
\begin{enumerate}
\item $U_{\cV,\fin}\cap N=U_{\cV,\infty\epsilon}\cap N=N_\cV^u$.

\item $U_{\cV,\fin}\cap N\cdot U^{\epsilon}=N_\cV^u\cdot U_{\cV,\fin}^{\epsilon}$.

\item $U_{\cV,\infty\epsilon}\cap (N\cdot U^\epsilon)=N_\cV^u\cdot U_{\cV,\infty}^\epsilon$,

\item $U_{\cV,\infty\epsilon}\cap U^{\epsilon}=U_{\cV,\infty}^{\epsilon}$.
\end{enumerate}

\end{Proposition}

\begin{proof}
By symmetry, we can assume $\epsilon=+$. 
Let $\alpha\in \Phi$. Let $w\in W^v$ be such that $\alpha\in w.\Delta_+$ (it suffices to take $w\in W^v$ such that $\alpha=w.\alpha_i$, for some $i\in I$). Then $U_{\alpha,\cV}\subset U_{\cV,\infty}(w.\Delta_+)$ and by Lemma~\ref{l_pre_GR08}, we have $\sU_{\cV,\infty}=U_{\cV,\infty}(w.\Delta_+)\cdot U_{\cV,\infty}(w.\Delta_-)\cdot N^u_{\cV}$. Therefore $\sU_{\cV,\infty}$ is stable by left  multiplication by $U_{\alpha,\cV}$ and hence $U_{\cV,\fin}\subset \sU_{\cV,\infty}$. Consequently, \begin{equation}\label{e_sU_cap_U+}
U_{\cV,\fin}^+\subset \sU_{\cV,\infty}\cap U^+.
\end{equation} Let $u\in \sU_{\cV,\infty}\cap U^+$. Write $u=u_+ u_-n$, with $u_+\in U_{\cV,\infty}^+$, $u_-\in U_{\cV,\infty}^-$ and $n\in N$.  Then $u_- n (u_+^{-1}u)^{-1}=1$, which proves that $u=u_+$, by the Birkhoff decomposition (\eqref{e_Birkhoff}). Therefore $\sU_{\cV,\infty}\cap U^+=U_{\cV,\infty}^+$ and similarly $\sU_{\cV,\infty}\cap U^-=U_{\cV,\infty}^-$. Using \eqref{e_sU_cap_U+}, we deduce $U_{\cV,\fin}^+\subset U_{\cV,\infty}^+$. Symmetrically, we get $U_{\cV,\fin}^-\subset U_{\cV,\infty}^-$.

 Therefore $U_{\cV,\fin}\subset \sU_{\cV,\fin}\subset U_{\cV,\fin}$, which proves that $U_{\cV,\fin}=\sU_{\cV,\fin}=U_{\cV,\fin}^+\cdot U_{\cV,\fin}^-\cdot N^u_\cV$ and by symmetry, we have $U_{\cV,\fin}=U_{\cV,\fin}^-\cdot U_{\cV,\fin}^+\cdot N^u_\cV$. 

Similarly, $\sU_{\infty+}$ is stable under left multiplication by $U_{\cV,\fin}$. It is obviously stable under left multiplication by $U_{\cV,\infty}^+$. As $N_{\cV}^u\subset U_{\cV,\fin}$, we deduce that $U_{\cV,\infty+}$ is stable under left multiplication by $N_{\cV}^u$. Therefore $\sU_{\cV,\infty+}$ is stable by multiplication. Let $g\in \sU_{\infty+}$. Write 
$g=u_+ u_- n$, with $u_+\in U_{\cV,\infty}^+$, $u_-\in U_{\cV,\fin}^-$ and $n\in N_{\cV}^u$. Then $g^{-1}=n^{-1} u_-^{-1} u_+^{-1}\in \sU_{\cV,\infty+}$, since $n^{-1},u_-^{-1},u_+^{-1}\in \sU_{\infty+,\cV}$. Therefore $\sU_{\cV,\infty+}$ is a group and hence $\sU_{\cV,\infty+}=U_{\cV,\infty^+}$.

Let $u\in U_{\cV,\infty+}\cap N\cdot U^{+}$. As $U_{\cV,\infty+}=(U_{\cV,\infty+})^{-1}$, we can write $u=n u_-u_+=n_1 v_+$, with $n\in N_{\cV}^u$, $u_-\in U_{\cV,\fin}^-$, $u_+\in U_{\cV,\infty}^+$, $n_1\in N$ and $v_+\in U^+$. Then $u_-^{-1} n^{-1} n_1v_+ (u_+)^{-1}=1$ and hence by  uniqueness in Birkhoff decomposition (\eqref{e_Birkhoff}), we have $v_+=u_+$ and $n_1=n$, which proves  that $u\in N_\cV^u\cdot U_{\cV,\infty}^{+}$. As the other inclusion is clear, we have $U_{\cV,\infty+}\cap N\cdot U^{+}=N^u_\cV\cdot U_{\cV,\infty}^{+}$. Therefore $U_{\cV,\infty+}\cap U^+=U_{\cV,\infty+}\cap N\cdot U^{+}\cap U^{+}=N_\cV^u\cdot U_{\cV,\infty}^{+}\cap U^{+}=U_{\cV,\infty+}^{+}$ by   uniqueness in the Birkhoff decomposition. Similarly, $U_{\cV,\infty+}\cap N=N_{\cV}^u$ and $U_{\cV,\fin}\cap N\cdot U^+=N_{\cV}^u\cdot U_{\cV,\fin}^+$. This completes the proof of the proposition.
\end{proof}

\subsubsection{Comparison of $U_{\cV,\fin}^+$ and  $U_{\cV,\infty}^+$}\label{sss_cmp_Ufnty}

If $\cV$ is  a filter on $\A$ and $\epsilon\in \{-,+\}$, we introduced two subgroups $U_{\cV,\fin}^\epsilon$ and $U_{\cV,\infty}^\epsilon$. We compare them in this subsubsection.

Let $\epsilon\in \{-,+\}$. When $G$ is reductive, $\Phi=\Delta$, $X_\alpha=x_\alpha$ for all $\alpha\in \Phi$ and thus $U_{\cV,\fin}^+=U_{\cV,\infty}^+$, for all filter $\cV$ on $\A$. Moreover, by \cite[6.4.2 and 6.4.9]{bruhat1972groupes}, we have $U_{\cV,\fin}^\epsilon=\prod_{\alpha\in \Phi_\epsilon} U_{\alpha,\cV}$, for any order on $\Phi_\epsilon$. In particular, $U_{\cV,\fin}^\epsilon=U_{\cV}^{\epsilon\epsilon}:=\langle U_{\alpha,\cV}\mid \alpha\in \Phi_\epsilon\rangle$, in the notation of \cite[4.3]{rousseau2016groupes}.

If $G$ is not reductive however, we can have $U_{\cV,\fin}^+\neq U_{\cV}^{++}$. We have $U_\cV^{++}\subset U_{\cV,\fin}^+$ by definition, but Rousseau proves in \cite[4.12.3 a)]{rousseau2016groupes} that when $G$ is affine $\mathrm{SL}_2$, we have $U_0^{++}\subsetneq U_{0,\infty}^+$.

We now prove (following  \cite[5.7.3]{rousseau2016groupes}) that there can exist a filter $\cV$ such that $U_{\cV,\infty}^+\neq U_{\cV,\fin}^+$. Let $x\in \A$ and set $\fq=x+C^v_f$. Then \[U_{\fq,\fin}^+=\langle U_{\alpha,\fq}\mid \alpha\in \Phi\rangle\cap U^+=\langle U_{\alpha,\fq}\mid \alpha\in \Phi_+\rangle\cap U^+=\langle U_{\alpha,x}\mid \alpha\in \Phi_+\rangle=U_{x}^{++}.\] By Corollary~\ref{c_Uma_cl}, we have $U_{\fq,\infty}^+=U_{x,\infty}^+$. Therefore $U_{\fq,\infty}^+\neq U_{\fq,\fin}^+$, if $x$ is such that $U_x^{++}\neq U_{x,\fin}^+$. We also have $U_{\fq,\infty+}=\langle U_{\fq,\infty}^+, U_{\alpha,\fq}\mid \alpha\in \Phi\rangle=U_{x,\infty}^+$ and $U_{\fq,\infty-}=\langle U_{\fq,\infty}^-,U_{\alpha,\fq}\mid \alpha\in \Phi\rangle=\langle 1,U_{\alpha,x}\mid \alpha\in\Phi_+\rangle=U_x^{++}$. Therefore we can have $U_{\cV,\infty+}\neq U_{\cV,\infty-}$.

However, we will see in Subsection~\ref{ss_para_minGp}, that when $\cV$ is  a spherical local face or a point, then $U_{\cV,\fin}^+=U_{\cV,\infty}^+$.

\subsection{The Iwasawa decomposition}\label{ss_Iwa_dec}

In this subsection, we prove the Iwasawa decomposition. It enables to define retractions with respect to sector-germs in the masure, which are an important tool in the study of Kac--Moody groups. It also enables to prove that the $G_x$, that we define in Section~\ref{s_parahoric}, decompose as the product of their $U^+$, $U^-$ and $N$ components (see Theorem~\ref{t_G_x}), which is a crucial property in the study of masures. We then compare the notion of inseparability, that we introduced here and the notion of narrowness, introduced in \cite{gaussent2008kac}

\subsubsection{The Iwasawa decomposition}

For $\alpha\in \Phi$, we set: \[G^{(\alpha)}=\langle T, U_\alpha,U_{-\alpha}\rangle \subset G.\] By \ref{a_KMT7}, if $\alpha=w.\alpha_i$, with $w\in W^v$ and $i\in I$, we have:
\[G^{(\alpha)}=\tilde{w} G^{(\alpha_i)} \tilde{w}^{-1}.\]
Let $\cS_i=((2),X,Y,(\alpha_i^\vee)_{\{i\}},(\alpha_i)_{\{i\}})$. Then $\cS_i$ is a Kac--Moody datum. The inclusion $\cS_i\rightarrow \cS$ is a Kac--Moody datum morphism in the sense of \cite[1.1 6) \& 7)]{rousseau2016groupes} and $\cS$ is an extension of $\cS_i$. By \cite[1.10]{rousseau2016groupes}, this induces an injective functor morphism $\fG_{\cS_{i}}\rightarrow \fG$. Then $\fG_{\cS_i}(\cF)$ is canonically isomorphic to $\langle T, U_{\alpha_i}, U_{-\alpha_i}\rangle$, by \cite[Before Proposition 7.75]{marquis2018introduction} and hence we can regard $\fG^{(\alpha_i)}$ as the Kac--Moody group associated with $\cS_i$.

We say that $\cV$ is \textbf{almost-open}\index{almost-open} if: \[\forall \alpha\in \Phi, f_{\cV}(\alpha)+f_{\cV}(-\alpha)>0.\]

Let $\alpha\in \Phi$. Write $\alpha=w.\alpha_i$, with $i\in I$ and $w\in W^v$. We set $\tilde{r_\alpha}=\tilde{w} \tilde{r_i} \tilde{w}^{-1}\in N$\index[notation]{r@$\tilde{r}_\alpha$}. This could depend on the choice of $w$ and $i$. However, $\tilde{r}_\alpha T=\nu^{-1}(\{1,r_\alpha\})$ does not.

\begin{Lemma}\label{l_GR08_3.8}(\cite[Lemma 3.8]{gaussent2008kac}) Let $\cV$ be an inseparable filter on $\A$.  Let $\alpha\in \Phi$. Then $Z_\alpha:=U_\alpha\cdot \{1,\tilde{r}_\alpha\} \cdot T\cdot U_\cV^{(\alpha)}$ contains $G^{(\alpha)}$. 
\end{Lemma}

\begin{proof}
Let us prove that if $n\in N$ is such that $\nu^v(n)=r_\alpha$ and $u\in U_\alpha$, then $nu\in Z_\alpha$. Write $u=x_\alpha(a)$, with $a\in \cF$.  If $f_\cV(\alpha)\leq \omega(a)$, then  $u\in U_{\alpha,\cV}$, and $nu\in Z_\alpha$. Otherwise, $f_{\cV}(\alpha)>\omega(a)$, thus $\cV\nsubseteq D(\alpha,\omega(a))$, by Lemma~\ref{l_char_fV}. As $\cV$ is inseparable, we deduce  $\cV\subset D(-\alpha,\omega)$ and hence  $\omega(a)<f_{\cV}(\alpha)\leq -f_{\cV}(-\alpha)$, by Lemma~\ref{l_char_fV}. Set $u'=x_{-\alpha}(a^{-1})$. Then $u'\in U_{-\alpha,\cV}$ and by \eqref{e_m}, $m:=u'uu'$ satisfies $\nu^v(m)=r_\alpha$. Hence $\nu^v(m^{-1})=r_\alpha=\nu^v(n)$ and  we can write $m^{-1}=tn $, with $t\in T$.  Using Proposition~\ref{p_cnj_N}, we get: \[nu=nu'^{-1}m u'^{-1}=\underbrace{nu'^{-1}n}_{\in U_{r_{\alpha}(-\alpha)}=U_{\alpha}}tu'^{-1}\in U_\alpha\cdot T\cdot U_{-\alpha,\cV}\in Z_\alpha.\] In other words, we proved that $r_\alpha T \cdot U_\alpha\subset Z_\alpha$ and hence: \[U_\alpha\cdot r_\alpha T \cdot U_\alpha\subset Z_\alpha.\]

By the Bruhat decomposition of $G^{(\alpha)}$ (\cite[Proposition B.25]{marquis2018introduction}), we have: \[G^{(\alpha)}\subset U_\alpha\cdot \{1,r_\alpha\}\cdot T\cdot U_\alpha=(U_\alpha\cdot  T\cdot U_\alpha)\cup (U_\alpha \cdot r_\alpha T \cdot U_\alpha ).\] By ~\ref{a_KMT4}, we have $U_{\alpha}\cdot T \cdot U_\alpha=U_{\alpha}\cdot T\subset Z_\alpha$. Lemma follows.
\end{proof}

The following proposition is the Iwasawa decomposition. It slightly generalizes \cite[Proposition 3.6]{gaussent2008kac}, by Proposition~\ref{p_nar_imp_ins} below.

\begin{Proposition}\label{p_Iwa}(Iwasawa decomposition, ) Let $\cV$ be a inseparable filter on $\A$. Then: \[G=U^+\cdot  N \cdot U_{\cV,\fin}.\] If moreover, $\cV$ is almost-open, then the natural map from $W^a=N/T_0$ onto $U^+\backslash N/U_{\cV,\fin}$ is one to one. 
\end{Proposition}

\begin{proof}
Consider $G'=U^+\cdot  N\cdot U_{\cV,\fin}$. We want to prove that $G'=G$. By Proposition~\ref{p_DRJ}, we have $G=\langle U_\alpha\mid \alpha\in \Phi,T\rangle$. Moreover, $\tilde{s_i}\in \langle U_{\alpha_i},U_{-\alpha_i}\rangle$, for $i\in I$. By \ref{a_KMT7}, we deduce that $G=\langle U_{\alpha_i},U_{-\alpha_i}\mid i\in I,T\rangle$.  Therefore it suffices to prove that $G'$ is stable under left multiplication by $U_{\pm\alpha_i}$, for all $i\in I$ and by $T$.  It is obviously stable under left multiplication by $U^+$ and it is stable by left multiplication by $T$ and thus it remains to prove that if $i\in I$, then $U_{-\alpha_i}\cdot  G'\subset G'$. 

Set $U(\Phi_+\setminus \{\alpha_i\})=U^+\cap U^{ma+}(\Delta_+\setminus\{\alpha_i\})$. By Proposition~\ref{p_Mar8.58}, we have $U^+=U_{\alpha_i}\cdot U(\Phi_+\setminus \{\alpha_i\})$ and $U(r_i.\Delta_+\setminus \{-\alpha_i\})=U_{-\alpha_i}\cdot U(\Phi_+\setminus \{\alpha_i\})$ (and we can permute the order in these products by using  $u\mapsto u^{-1}$).   We have: \begin{align*}
U_{-\alpha_i}\cdot  G'&=U_{-\alpha_i}\cdot U(\Phi_+\setminus\{\alpha_i\})\cdot  U_{\alpha_i}\cdot N\cdot U_{\cV,\fin}\\
&=  U(\Phi_+\setminus \{\alpha_i\})\cdot U_{-\alpha_i}\cdot U_{\alpha_i}\cdot N\cdot U_{\cV,\fin}\\
&\subset U(\Phi_+\setminus\{\alpha_i\})\cdot G^{(\alpha_i)}\cdot N\cdot U_{\cV,\fin}\\
&\subset U(\Phi_+\setminus\{\alpha_i\})\cdot U_{\alpha_i}\cdot \{1,\tilde{r}_{i}\}\cdot T\cdot U_{\cV}^{(\alpha_i)}\cdot N\cdot U_{\cV,\fin} &\text{by Lemma~\ref{l_GR08_3.8}}\\
&= U^+ \cdot \{1,\tilde{r}_i\} \cdot T\cdot U_{\cV}^{(\alpha_i)}\cdot N\cdot U_{\cV,\fin}
\end{align*}

By Lemma~\ref{l_GR08_3.3}, we have $U_{\cV}^{(\alpha_i)}\subset U_{\alpha_i,\cV}\cdot U_{-\alpha_i,\cV}\cdot N_{\cV}^u$ and  thus: \[
U^+\cdot T\cdot U_{\cV}^{(\alpha_i)}\cdot N\cdot U_{\cV,\fin}=U^+\cdot U_{-\alpha_i,\cV}\cdot N\cdot U_{\cV,\fin}.\]

We also have: \begin{align*}
U^+\cdot r_i T \cdot U_{\cV}^{(\alpha_i)}\cdot N\cdot U_{\cV,\fin}&=U^+\cdot r_i T\cdot U_{-\alpha_i,\cV}\cdot U_{\alpha_i,\cV}\cdot N\cdot U_{\cV,\fin}\\
&\subset U^+ \cdot U_{\alpha_i}\cdot   U_{-\alpha_i}\cdot  N\cdot U_{\cV,\fin} &\text{by \ref{a_KMT4} and Lemma~\ref{l_cnj_N_U}}\\
&= U^+ \cdot U_{-\alpha_i}\cdot N\cdot U_{\cV,\fin}.
\end{align*}

Consequently, \[U_{-\alpha_i}\cdot G'\subset U^+\cdot U_{-\alpha_i}\cdot N\cdot U_{\cV,\fin}.\] Therefore it suffices to prove that $U_{-\alpha_i}\cdot N\subset G'$. Let $u\in U_{-\alpha_i}$ and $n\in N$. Let $\beta=-\nu^{v}(n^{-1})$. Then by Lemma~\ref{l_cnj_N_U} and Lemma~\ref{l_GR08_3.8}, we have: \[un=n\cdot n^{-1}un\in nU_{\beta}\subset n U_{-\beta}\cdot \{1,r_{\beta}\}\cdot T \cdot U_{\cV}^{(\beta)}.\] Therefore: \[un\in U_{\alpha_i}\cdot n\cdot \{1,r_{\beta}\}\cdot T\cdot U_{\cV}^{(\beta)}\subset U_{\alpha}\cdot N\cdot U_{\cV,\fin}\subset G'.\] Thus $G'=G$, which proves the Iwasawa decomposition.

Assume moreover that  $\cV$ is almost-open. Then $N_\cV^u=T_0$. Let $n,n'\in N$ be such that $n'\in U^+ \cdot n\cdot U_{\cV,\fin}$. Let $w=\nu^v(n^{-1})\in W^v$. Then $n^{-1}n'\in n^{-1}\cdot U^+ \cdot n \cdot U_{\cV,\fin}=U(w.\Phi_+) \cdot U_{\cV,\fin}$. Using Proposition~\ref{p_GR08_3.4},  we deduce: \[n^{-1}n'\in U(w.\Phi_+)\cdot  U_{\cV,\fin}(w.\Phi_+)\cdot  U_{\cV,\fin}(w.\Phi_-)\cdot  T_0\subset U(w.\Phi_+)\cdot U(w.\Phi_-) \cdot T_0=U(w.\Phi_+)\cdot T_0\cdot U(w.\Phi_-).\] By uniqueness  in the Birkhoff decomposition (\eqref{e_Birkhoff}), we deduce $n^{-1}n'\in T_0$, which completes the proof of the proposition.
\end{proof}

Note that the Iwasawa decomposition has a masure theoretic interpretation : it corresponds to the existence of apartments containing certain pairs of filters, see Remark~\ref{r_iwa}. This proposition is actually a consequence of Proposition~\ref{p_mas_th_iwa}, once we know that the masure of a Kac--Moody group is a masure in the sense of Definition~\ref{d_w_mas}. However, we will use the Iwasawa decomposition in order to prove that $\I$ is an abstract masure: we use it in the proof of Theorem~\ref{t_G_x}, which gives a decomposition of the fixator $G_x$ of an element $x$  of $\A$.

\subsubsection{Comparison of the notions of narrowness and of inseparability}

We introduced the notion of inseparability in Definition~\ref{d_insep_filter}. We prove that it is slightly more general than the notion of narrowness introduced in \cite{gaussent2008kac}. Note that the notion of narrowness will not be used elsewhere in this paper.

\begin{Definition}\label{d_ao_nar}
We say that $\cV$ is \textbf{narrow}\index{narrow} if  for all $\alpha\in \Phi$, we have: \begin{enumerate}
\item $f_{\cV}(\alpha)\notin \{-\infty,+\infty\}$,

\item $f_{\cV}(\alpha)+f_{\cV}(-\alpha)\in \{0,0^+\}$ or in the discrete case, if $f_{\cV}(\alpha)+f_{\cV}(-\alpha)\in \{0,1\}$,

\item $(f_{\cV}(\alpha),f_{\cV}(-\alpha))\notin  \{(\lambda^+,(-\lambda)^+)\mid \lambda\in \Lambda\}$.
\end{enumerate} 
\end{Definition}

Note that condition (1) was implicit in \cite{rousseau2016groupes}.

\begin{Proposition}\label{p_nar_imp_ins}
Let $\cV$ be a filter on $\A$. Then the following are equivalent:\begin{enumerate}
\item $\cV$ is narrow,

\item $\cV$ is bounded modulo $\A_{\ines}$ and $\cV$ is inseparable.
\end{enumerate} 
\end{Proposition}

\begin{proof}
(1) $\Rightarrow$ (2). Assume that $\cV$ is narrow. First assume that we are in the non-discrete case. Let $\alpha\in \Phi$. By definition of narrowness, $f_\cV(\alpha)\in \R$ and $f_\cV(-\alpha)\in \R$. Then we have $\cV\subset D(\alpha,a)$ and $\cV\subset D(-\alpha,-b)$, for all $b\in \R$ such that $b<f_\cV(-\alpha)$. Therefore $\cV\subset D(-\alpha,-b)$ for all $b\in ]-\infty,a]$ and $\cV\subset D(\alpha,b)$ for all $b\in [a,+\infty[$. Therefore $\cV$ is inseparable. 

Assume now that $\Lambda=\Z$. If $f_\cV(\alpha)=-f_{\cV}(-\alpha)$, we have $\cV\subset D(\alpha,b)$ for all $b\in [f_\cV(\alpha),+\infty[$ and $\cV\subset D(-\alpha,-b)$ for all $b\in ]-\infty,f_\cV(\alpha)]$. Assume $f_\cV(\alpha)=-f_{\cV}(-\alpha)+1$. Then $\cV\subset D(\alpha,b)$, for all $b\in [f_\cV(\alpha),+\infty[$ and for all $b\in \Z$ such that $b<f_{\cV}(\alpha)$, we have $b\leq f_{\cV}(\alpha)-1$ and hence $\cV\subset D(-\alpha,1-f_{\cV}(\alpha))=D(-\alpha,f_{\cV}(-\alpha))\subset D(-\alpha,-b)$. Therefore $\cV$ is inseparable.

Let us prove that the image of $\cV$ in $\A/\A_{\ines}$ is bounded. Let $i\in I$. By definition, $f_{\cV}(\alpha)\neq +\infty$ and thus there exists $\lambda_i\in \Lambda$ such that $D(\alpha_i,\lambda_i)\in \cV$. Let $x\in \A$ be such that $\alpha_i(x)=-\lambda_i$ for all $i\in I$. Then $\bigcap_{i\in I}D(\alpha_i,\lambda_i)=x+\overline{C^v_f}\in \cV$. Considering the $-\alpha_i$ instead of the $\alpha_i$, we get the existence of $y\in \A$ such that $y-\overline{C^v_f}\in \cV$. Then $(x+\overline{C^v_f})\cap (y-\overline{C^v_f})\in \cV$. As $(\alpha_i)_{i\in I}$ induces a basis of $(\A/\A_{\ines})^*$, we deduce that the image of $(x+\overline{C^v_f})\cap (y-\overline{C^v_f})$ in $\A/\A_{\ines}$ is bounded and thus $\cV$ is bounded modulo $\A_{\ines}$. Therefore we have (2).

(2) $\Rightarrow$ (1). Assume that $\cV$ is inseparable and that $\cV$ is bounded modulo $\A_{\ines}$.  Let $\Omega\in \cV$ be bounded modulo $\A_{\ines}$. Then for all $\alpha\in \Lambda$,  there exist $\lambda,\mu\in \Lambda$ such that $\Omega\subset D(\alpha,\lambda)$ and $\Omega\cap  D(\alpha,\mu)=\emptyset$. Then $\cV\subset D(\alpha,\lambda)$ and  $D(\alpha,\mu)\notin \cV$ (otherwise, $\emptyset$ would belong to $\cV$) and thus $\cV\nsubseteq D(\alpha,\mu)$. Therefore $f_{\cV}(\alpha)\notin \{\pm\infty\}$, for all $\alpha\in \Phi$. 

By Proposition~\ref{p_f_cV_alpha}, we have $f_{\cV}(\alpha)+f_{\cV}(-\alpha)\geq 0$, for all $\alpha\in \Phi$.  First assume that $\Lambda$ is discrete. Assume by contradiction that $f_{\cV}(\alpha)+f_{\cV}(-\alpha)\geq 2$, for some $\alpha\in \Phi$. Then there exists $\lambda\in \Lambda$ such that $f_{\cV}(\alpha)=\lambda$ and $f_{\cV}(-\alpha)\geq -\lambda+2$. Then $\cV\nsubseteq D(\alpha,\lambda-1)$ and $\cV\nsubseteq D(-\alpha,-(\lambda-1))$: this contradicts the inseparability of $\cV$. Thus $f_{\cV}(\alpha)+f_{\cV}(-\alpha)\leq 1$ and hence $\cV$ is narrow.

Assume now that $\Lambda$ is dense in $\R$. Let $\alpha\in \Phi$. If $(f_{\cV}(\alpha),f_{\cV}(-\alpha))=(\lambda^+,(-\lambda)^+)$, for some $\lambda\in \Lambda$, then $\cV\nsubseteq D(\alpha,\lambda)$ and 
$\cV\nsubseteq D(-\alpha,-\lambda)$: a contradiction. Assume  $f_{\cV}(\alpha)+f_{\cV}(-\alpha)\notin\{0,0^+\}$. There exist $\lambda,\mu\in \Lambda$ such that  $f_{\cV}(\alpha)>\lambda$, $f_{\cV}(-\alpha)>\mu$ and $\lambda+\mu>0$.
 Then by Lemma~\ref{l_char_fV}, we have  $\cV\nsubseteq D(\alpha,\lambda)$ and $\cV\nsubseteq D(-\alpha,\mu)$. As $-\lambda<\mu$, we have $D(-\alpha,-\lambda)\subset D(-\alpha,\mu)$ and thus $\cV\nsubseteq D(-\alpha,-\lambda)$: a contradiction. Thus $\cV$ is narrow.
 
\end{proof}

\section{Parahoric subgroup associated with a point}\label{s_parahoric}

For $x\in \A$, we set $G_x=\langle U_{x,\infty}^+, U_{x,\infty}^-,N_x\rangle$\index[notation]{g@$G_x$}. A  crucial property of $G_x$ is that\begin{equation}\label{e_G_x}
G_x=U_{x,\infty}^+\cdot U_{x,\infty}^-\cdot N_x=U_{x,\infty}^-\cdot U_{x,\infty}^+\cdot N_x\end{equation}(see Theorem~\ref{t_G_x}). The aim of this subsection is to prove this property.

Compared to the reductive case, a difficulty is that the definition of $G_x$ involves $G^{ma+}$ and $G^{ma-}$. Moreover, in the reductive case, $\Delta$ is finite and thus we can make inductions, using commutation relations between $U_{\alpha,x}$ and $U_{\beta,x}$, for $\alpha,\beta\in \Phi$. In the Kac--Moody case, $\Delta$ is infinite and thus it is difficult to apply these methods. In order to  study $G_x$, we interpret it via the action of $G$ on modules. This interpretation will not be used elsewhere in the paper and the reader ready to admit \eqref{e_G_x}   can skip this section.

\subsection{Study of the adjoint action of $U^+$ on $\cU$ and $\widehat{\cU}$}\label{ss_ad_U_cU}

For $x\in \A$, the fixator $G_x$ wil be defined as the fixator of spaces $M_x$, for $M$ a highest-weight module or $M=\cU$ (see Definition~\ref{d_parah}). In this subsection, we study the action of $G$ on $(\cU,\Ad)$. We obtain an interpretation of $U_{\cV,\infty}^+$ as the fixator of some space (see Lemma~\ref{l_Rou4.11_U}).

Recall the definition of  $\widehat{\Ad}:\widehat{\cU}^+_\cF\rightarrow \End(\widehat{\cU^+})$ in \eqref{e_ad_hat}.

\begin{Lemma}\label{l_UVp_stable}

Let $\cV$ be a filter on $\A$.
\begin{enumerate}
\item Let $\lambda,\mu\in Q$. Let $u_\lambda\in \cU_{\lambda,\cV}$ and $v_\mu\in \cU_{\mu,\cV}$ Then $\Ad(u_\lambda)(v_\mu)\in \cU_{\lambda+\mu,\cV}$.

\item We have $\Ad(\cU_\cV)(\cU_\cV)\subset \cU_\cV$ and $\widehat{\Ad}(\widehat{\cU}_\cV^+)(\widehat{\cU}_\cV^+)\subset \widehat{\cU}_\cV^+$.
\end{enumerate}
\end{Lemma}

\begin{proof}
(1) Let $\Omega\in \cV$ be such that $u_\lambda\in \cU_{\lambda,\Omega}$ and $v_\mu\in \cU_{\mu,\Omega}$.  By additivity of $\Ad$, we can assume that $u_\lambda=t x_1\ldots x_k$, for some $t\in \cF_{\geq f_\Omega(\lambda)}$, $k\in \N$, $\beta_1,\ldots,\beta_k\in \Delta$ such that $\sum_{i=1}^k \beta_i=\lambda$ and $x_i\in \ffg_{\beta_i,\Z}$ for all $i\in \llbracket 1,k\rrbracket$.

If $\nu\in Q$,  $y\in \cU_{\nu}\otimes \cF_{\geq f_\Omega(\mu)}$ and $i\in \llbracket 1,k\rrbracket$, we have $\Ad(x_i)(y)\in \cU_{\nu+\beta_i}\otimes \cF_{\geq f_{\Omega}(\mu)}$, which implies that $t\Ad(x_1\ldots x_k)(v_\mu)\in \cU_{\mu+\lambda}\otimes t\cF_{\geq f_\Omega(\mu)}\subset \cU_{\mu+\lambda}\otimes \cF_{\geq f_\Omega(\mu)+f_\Omega(\lambda)}\subset \cU_{\mu+\lambda}\otimes \cF_{f_{\Omega}(\lambda+\mu)}$, by concavity of $f_\Omega$ (see \eqref{e_concavity}). Therefore $\Ad(u_\lambda)(v_\lambda)\in \cU_{\lambda+\mu,\Omega}\subset \cU_{\lambda+\mu,\cV}$, which proves (1). 

(2) Let $\Omega\in \cV$. Then by (1), we have $\Ad(\cU_\Omega)(\cU_\Omega)\subset \cU_{\Omega}$ and $\widehat{\Ad}(\widehat{\cU}^+_\Omega)(\widehat{\cU}^+_\Omega)\subset \widehat{\cU}^+_\Omega$, which implies (2).

\end{proof}

The following lemma is a part of \cite[Lemme 4.11]{rousseau2016groupes}.

\begin{Lemma}\label{l_Rou4.11_U}
Let $\cV$ be a filter on $\A$. Then $\{u\in U^+\mid \Ad(u)(\cU_{\cV})=\cU_\cV\}=U_{\cV,\infty}^+$.
\end{Lemma}

\begin{proof}
Let $H=\{u\in U^+\mid \Ad(u)(\cU_{\cV})=\cU_\cV\}$.

By Lemma~\ref{l_UVp_stable}, $\widehat{\Ad}(U_{\cV,\infty}^+)(\widehat{\cU}^+_\cV)\subset \widehat{\Ad}(\widehat{\cU}^{+}_\cV)(\widehat{\cU}^+_\cV)\subset \widehat{\cU}^+_\cV$. Moreover, $U_{\cV,\infty}^+\subset U^+$ and  thus $\widehat{\Ad}(U_{\cV,\infty}^+)(\cU_\cV) \subset \widehat{\cU}^+_\cV\cap \cU=\cU_\cV$.  Therefore $U_{\cV,\infty}^+\subset H$. We need to prove the converse inclusion.

We fix a total order on $\Delta_+$ compatible with the height. Let $u\in H$. Then we have $\Ad(u)|_{\cU_\cV}=\widehat{\Ad}(u)|_{\cU_\cV}$  and thus $\widehat{\Ad}(u)$ stabilizes the closure of $\cU_\cV$ in $\widehat{\cU}^+$. 
Therefore, $\widehat{\Ad}(u)$ stabilizes $\widehat{\cU}^+_\cV$. Write $u=\prod_{\alpha\in \Delta_+}X_\alpha(\underline{u_\alpha})$, with $\underline{u_\alpha}\in \ffg_{\alpha,\Z}\otimes \cF$, for $\alpha\in \Delta_+$.

Let $\alpha\in \Delta_+$. By induction, we assume that for all $\beta\in \Delta_+$ such that $\beta<\alpha$, we have $\underline{u_\beta}\in \ffg_{\beta,\Z}\otimes \cF_{\geq f_{\cV}(\beta)}$. In particular, we assume that $\prod_{\beta\in \Delta_+\mid \beta<\alpha} X_\beta(\underline{u_\beta})\in U_{\cV}^{ma+}$. Thus we can replace $u$ by $\prod_{\beta\in \Delta_+\mid \beta\geq \alpha} X_\beta(\underline{u_{\beta}})$, which also stabilizes $\widehat{\cU}^{+}_\cV$, since $U_\cV^{ma+}$ stabilizes $\widehat{\cU}^+_{\cV}$, by Lemma~\ref{l_UVp_stable}.

 Let $k=\mathrm{rk}(\ffg_{\alpha,\Z})$ and  $b_1,\ldots,b_k$ be a $\Z$-basis of $\ffg_{\alpha,\Z}$. Write $\underline{u_\alpha}=\sum_{i=1}^k t_i b_i$, with $t_i\in \cF$, for $i\in \llbracket 1,k\rrbracket$. We have $X_\alpha(\underline{u_\alpha})=\prod_{i=1}^k [\exp](t_i b_i)$. For $i\in \llbracket 1,k\rrbracket$, we have  $[\exp](t_i b_i)=\sum_{p=0}^{+\infty} t_i^p b_i^{[p]}$ and $b_i^{[p]}\in \cU_{p\alpha}$, by definition of exponential sequences.

 Let $h\in \cB_0$ (a $\Z$-basis of $\ffg_{0,\Z}=\fh\cap \ffg_\Z$) be such that $m:=\alpha(h)\neq 0$. Up to changing $h$ in $-h$ (and thus to changing the $\Z$-basis  of $\ffg_{0,\Z}$), we can assume that $m>0$.  We have $\widehat{\Ad}([\exp] (t_i b_i))=\sum_{p=0}^{\infty} t_i^p \Ad(b_i^{[p]})$. We have $X_\alpha(\underline{u_\alpha})=\prod_{i=1}^k [\exp](t_i b_i)$ and hence: \begin{equation}\label{e_X_alp}
 X_\alpha(\underline{u_\alpha})=1+\sum_{i=1}^k t_i b_i+v',
\end{equation}
 where $v'\in \prod_{\nu\in 2\alpha+Q^+} \cU_\nu$.

  Let $n\in \N$. We have $\tbinom{h}{n}\in \cU_0$ and by \eqref{e_X_alp}:  \[\widehat{\Ad}\left(X_\alpha(\underline{u_\alpha})\right)(\tbinom{h}{n})=\tbinom{h}{n}+\sum_{i=1}^k t_i \Ad(b_i)(\tbinom{h}{n})+u',\] for some $u'\in \widehat{\cU}^+$ such that $\pr_\nu(u')\neq 0\Rightarrow \nu \in \Z_{\geq 2} \alpha$, for $\nu\in Q_+$, where $\pr_\nu:\widehat{\cU}^+_\cF\rightarrow \cU_{\nu}\otimes \cF$ denotes the projection on the $\nu$-th component. Therefore: \[\pr_\alpha\left((u)(\tbinom{h}{n})\right)=\sum_{i=1}^k t_i\Ad(b_i)(\tbinom{h}{n}).\] Let $i\in \llbracket 1,k\rrbracket$. We have: \begin{align*}
 t_i \Ad(b_i)(\tbinom{h}{n})=t_i (b_i\tbinom{h}{n}-\tbinom{h}{n} b_i)=-t_i (\tbinom{h}{n} b_i-b_i\tbinom{h}{n}).
\end{align*}

By \eqref{e_cmt}, we have:  \begin{align*}n! \tbinom{h}{n} b_i&=h(h-1)\ldots (h-n+1) b_i\\ &= h(h-1)\ldots (h-n+2)b_i(h-n+1+\alpha(h))\\ &=\ldots =n! b_i \tbinom{h+\alpha(h)}{n}.\end{align*}

Therefore:  \[t_i \Ad(b_i)(\tbinom{h}{n})=-t_ib_i(\tbinom{h+\alpha(h)}{n}-\tbinom{h}{n})=-t_ib_i\sum_{q=0}^{n-1}\tbinom{h}{q}\tbinom{\alpha(h)}{n-q}\] and \[\sum_{i=1}^k  -t_ib_i\sum_{q=0}^{n-1}\tbinom{h}{q}\tbinom{\alpha(h)}{n-q}\in \cU_{\alpha,\cV}.\] As the $b_i\tbinom{h}{q}$, $i\in \llbracket 1,k\rrbracket$, $q\in \N$ belong to a $\Z$-basis of $\cU$, we deduce that \[\omega(t_i \tbinom{\alpha(h)}{n-q})\geq f_{\cV}(\alpha), \forall i\in \llbracket 1,k\rrbracket, q\in \llbracket 1,n-1\rrbracket.\] Applying it with $q=0$ and $n=m$, we obtain $\omega(t_i)\geq f_{\cV}(\alpha)$, for all $i\in \llbracket 1,k\rrbracket$, which proves that $X_\alpha(\underline{u_\alpha})\in U_{\cV}^{ma+}$. By induction, we deduce that $u\in U_\cV^{ma+}\cap U=U_{\cV,\infty}^+$. Therefore $H\subset U_{\cV,\infty}^+$, which proves that $H=U_{\cV,\infty}^+$.
\end{proof}

\subsection{Highest weight modules}

In this subsection, we recollect some facts on the highest weight modules that we use in the definition of the $G_x$. 

Recall the definition of $\fn^+$ in \ref{sss_rt_space_dec}.

\begin{Definition}
A $\ffg$-module $V$ is called a \textbf{highest-weight module} with \textbf{highest weight} $\lambda\in \fh^*$ if there exists a non-zero vector $v_\lambda\in V$ (called a \textbf{highest-weight vector}) such that \begin{enumerate}
\item $\fn^+.v_\lambda=0$ and $h.v_\lambda=\lambda(h) v_\lambda$, for all $h\in \fh$,

\item $\cU_\C(\ffg).v_\lambda=V$.
\end{enumerate}
\end{Definition}

Let $V$ be a highest-weight module with highest-weight $\lambda$ and highest weight vector $v_\lambda$.  Note that by \eqref{e_cmt}, if  $\alpha\in Q_+$ and $u\in \cU_{-\alpha,\C}$, we have $u.v_\lambda\in V_{\lambda- \alpha}$. Therefore   a highest-weight module $V$ is $\fh$-diagonalizable: $V=\bigoplus_{\alpha\in Q_+} V_{\lambda-\alpha}$ and $V_\lambda=\cU_{0,\C}.v_\lambda=\fh.v_\lambda=\C v_\lambda$ is one-dimensional.  More generally, if $\nu\in Q$ and $\mu\in \fh^*$, we have $\cU_{\nu,\C}.V_\mu\subset V_{\mu+\nu}$.

A highest-weight $\ffg$-module $M(\lambda)$ with highest-weight $\lambda$ is called a \textbf{Verma module} if every highest-weight $\ffg$-module with highest weight $\lambda$ is a quotient of $M(\lambda)$. By \cite[Lemma 4.9]{marquis2018introduction}, for every $\lambda\in \fh^*$, there exists a unique (up to isomorphism) Verma module $M(\lambda)$. Moreover, $M(\lambda)$ contains a unique proper maximal submodule $M'(\lambda)$. One then sets $L(\lambda)=M(\lambda)/M'(\lambda)$. By \cite[Proposition 4.11]{marquis2018introduction}, $L(\lambda)$ is an integrable $\ffg$-module (in the sense of  \cite[4.1]{marquis2018introduction}) if and only if $\lambda$ is \textbf{dominant integral}, which means that $\lambda(\alpha_i^\vee)\in \N$ for all $i\in I$.

Let $X^+=\{\lambda\in X\mid \lambda(\alpha_i^\vee)\geq 0, \forall i\in I\}$\index[notation]{x@$X^+$} be the set of dominant integral characters.  Let $\lambda\in X^+$.  We choose $v_\lambda\in L(\lambda)_\lambda\setminus\{0\}$ and we set $L_\Z(\lambda)=\cU.v_\lambda$. This is a $\cU$-module, which does not depend on the choice of $v_\lambda$ (up to isomorphism) since $L(\lambda)_\lambda$ is one-dimensional.

Let $\alpha\in Q$ and $u_\alpha \in \cU_{\alpha}$. Then by Lemma~\ref{l_wghts_cU}, we have:  \begin{equation}\label{e_wghts_L}
u_\alpha.v_\lambda\in L_\Z(\lambda)_{\lambda+\alpha}
\end{equation}  for all $\alpha\in Q$ and thus: \begin{equation}\label{e_dec_L_lmbd}
 L_\Z(\lambda)=\bigoplus_{\alpha\in Q_-} L_\Z(\lambda)_{\lambda+\alpha}.
 \end{equation}
 
 Note that in \cite[2.14]{rousseau2016groupes}, Rousseau takes an other definition of $L_\Z(\lambda)$. We get however the same ``fixator'' $G_x$, for $x\in \A$, by Theorem~\ref{t_G_x}.

\subsection{Highest weight modules associated with a filter}
Let $\lambda\in X^+$ and $M=L_\Z(\lambda)$. Then $M$ is naturally a $G$-module, by \cite[Theorem 7.48]{marquis2018introduction}.

In this subsection, we associate, to   each filter $\cV$ on $\A$ a  subspace $M_\cV$ of $M$ and we study the action of $U^+$ on $M$. 

 We want to define the fixator $G_x$, $x\in \A$ via the use of the $M_x$.  When $\omega(\cF^\times)=\Z$, these spaces do not characterize $x+\A_{\ines}$ in general (see Example~\ref{Ex_SL2_F_F'}). Therefore it is convenient to embed $(\cF,\omega)$ in an extension $(\cF',\omega')$ such that $\omega'(\cF')$ is dense in $\R$.

\begin{Notation}\label{n_modules}
 If $\omega(\cF^\times)$ is dense in $\R$, we set $\cF'=\cF$ and $\omega'=\omega$. Otherwise, we can assume that $\omega(\cF^\times)=\Z$, up to renormalization.  We set $\cF'=\cF(\bx)$, where $\bx$ is an indeterminate.  By \cite[Theorem 2.2.1]{engler2005valued}, for every $x\in \R$, there exists a valuation $\omega_x:\cF'\rightarrow \R\cup \{\infty\}$ such that $\omega_x(\bx)=x$ and $\omega_x|_\cF=\omega$. We choose $x_0\in \R$ such that $\Z+\Z x_0$ is dense in $\R$ and we set $\omega'=\omega_{x_0}$ and $\Lambda'=\omega'((\cF')^\times )$.
\end{Notation}

If $r\in \R$, we set $\cF'_{\geq r}=\{x\in \cF'\mid \omega'(x)\geq r\}$.

If $\cV$ is a filter on $\A$, we denote by $f_\cV':\Delta\rightarrow \tilde{\R}\cup \{+\infty\}$ the analogue of $f_\cV$ (defined in Subsection~\ref{ss_f_cV}) but with $(\cF,\omega)$ replaced by $(\cF',\omega')$.

Let $\lambda\in X^+$ and $M=L_\Z(\lambda)\otimes \cF'$.  For $\mu\in \lambda-Q_+$ and $r\in \R$, we set  $M_{\mu,r}=M_{\mu,r}(\lambda)=L_\Z(\lambda)_\mu \otimes_\Z \cF'_{\geq r}$.  If $\Omega$ is a non-empty subset of $\A$, we set $M_{\mu,\Omega}=M_{\mu,\Omega}(\lambda)=M_{\mu,f'_\Omega(\mu)}$. We set $M_{\Omega}=\bigoplus_{\alpha\in Q_-}M_{\lambda-\alpha,\Omega}$. Now if $\cV$ is a filter on $\A$, we set $M_\cV=\bigcup_{\Omega\in \cV} M_\Omega$\index[notation]{m@$M_{\cV}, M_{\Omega}$}.

\begin{Example}\label{Ex_SL2_F_F'}
Consider $G=\mathrm{SL}_2(\cF)$. Let $\qp:T\rightarrow \cF^\times$ be defined by $\qp(\begin{psmallmatrix} h & 0\\ 0 & h^{-1}\end{psmallmatrix})=h$, for $h\in \cF^\times$. Let $\alpha^\vee:\cF^\times \rightarrow T$ be defined by $\alpha^\vee(h)=\begin{psmallmatrix} h & 0\\ 0 & h^{-1}\end{psmallmatrix}$, for $h\in \cF^\times$. Let $\tilde{M}=\cF^2$ be the $G$-module defined by the canonical action of $\mathrm{SL}_2(\cF)$ on $\cF^2$. Set $v_{\qp}=(1,0)$ and $v_{-\qp}=(0,1)$. We have $\tilde{M}=\cF v_\qp \oplus \cF v_{-\qp}$. Let $x\in \A=\R \alpha^\vee$. We have: \[\tilde{M}_x=\tilde{M}_{\{x\}}= \cF_{\geq x} v_{\qp}\oplus \cF_{\geq -x} v_{-\qp}=\cF_{\geq \lceil x\rceil}v_{\qp} \oplus \cF_{\geq \lceil -x\rceil} v_{-\qp}\] (where $\tilde{M}$ is ``equipped'' with $\omega$ instead of $\omega'$).

Let now $M=\cF' v_\qp \oplus \cF'v_{-\qp}$. Then  $M_{x}=\cF'_{\geq x} v_\qp \oplus   \cF'_{\geq -x} v_{-\qp}$. As $\omega'(\cF')$ is dense in $\R$, the $M_{\{x\}}$ are all distinct, for $x\in \A$, whereas if $\omega(\cF^\times )=\Z$, then $\tilde{M}_x=\tilde{M}_y$ for all $x,y\in \R$ such that $x,y\in ]n,n+1[$ for some $n\in \Z$.
\end{Example}

\begin{Lemma}\label{l_cpltd_actn}
Let $\lambda\in X^+$ and $M=L_\Z(\lambda)\otimes \cF'$. Then the action of $\cU_\cF^+$ on $M$ naturally extends to an action of $\widehat{\cU}^+_\cF$ on $M$, defined by \begin{equation}\label{e_actn_Uhat_M}
\sum_{\nu\in Q_+} u_\nu. \sum_{\mu\in Q_+}m_{\lambda-\mu}=\sum_{\nu\in Q_+, \mu\in Q_+} u_\nu.m_{\lambda-\mu},
\end{equation}  for $(m_{\lambda-\mu})\in \prod_{\mu\in Q_+}M_{\lambda-\mu}$ with at most finitely many non-zero terms and  $(u_\nu)\in \prod_{\nu\in Q_+}\cU_{\nu,\cF}$.
\end{Lemma}

\begin{proof}
For $\mu,\nu\in Q_+$, $u_\nu.m_{\lambda-\mu}\in M_{\lambda-\mu+\nu}$, by \eqref{e_wghts_L}. Therefore $u_\nu.m_{\lambda-\mu}\neq 0$ implies $\nu\leq_{Q^\vee} \mu$, which proves that \eqref{e_actn_Uhat_M}  is well-defined.
\end{proof}

\begin{Lemma}\label{l_U_v_stb_M_v}
Let $\lambda\in X^+$,  $\cV$ be a filter on $\A$ and $M=L_\Z(\lambda)\otimes \cF$. Then $\widehat{\cU}^+_\cV$ stabilizes $M_\cV$. In particular,  $U_{\cV,\infty}^+$ stablizes $M_\cV$.
\end{Lemma}

\begin{proof}
Let $\Omega\in \cV$.   Let $(m_{\lambda-\mu})\in \prod_{\mu\in Q_+} M_{\lambda-\mu,\Omega}$ with at most finitely many non-zero terms  and $(u_\nu)\in \prod_{\nu\in Q_+} \cU_{\nu,\Omega}$. Let $m=\sum_{\mu\in Q_+}m_{\lambda-\mu}\in M_\Omega$ and $u=\sum_{\nu\in Q_+} u_\nu\in \widehat{\cU}^+_\Omega$. Let $\mu,\nu\in Q_+$. Write $m_{\lambda-\mu}=\sum_{j\in J} a_j m_j$, with $J$ a finite set, $(a_j)\in (\cF'_{\geq f_{\Omega}(\lambda-\mu)})^J$  and $(m_j)\in (L_\Z(\lambda)_{\lambda-\mu})^J$. Write $u_\nu=\sum_{k\in K} b_k u_k$, with $K$ a finite set,  $(u_k)\in (\cU_{\nu})^K$ and $(b_k)\in (\cF'_{\geq f_{\Omega}(\nu)})^K$. Let $(j,k)\in J\times K$ and $x\in \Omega$. Then  we have \[\omega'(a_jb_k)+(\lambda-\mu+\nu)(x)=\omega'(a_j)+(\lambda-\mu)(x)+\omega'(b_k)+\nu(x)\geq 0,\] since $\omega'(a_j)+(\lambda-\mu)(\Omega),\omega'(b_k)+\nu(\Omega)\geq 0$, by assumption. Therefore $a_j b_k m_j.u_k\in M_{\lambda-\mu+\nu,\Omega}$, which proves that $m_{\lambda-\mu}.u_\nu\in M_{\lambda-\mu+\nu,\Omega}$. Lemma follows, since $U_{\cV,\infty}^+$ is contained in $\widehat{\cU}_{\cV}^+$ (see Definition~\ref{d_U_Vma}).
\end{proof}

\begin{Remark}\label{r_U_v_stb_M_v}
By symmetry, we obtain that if $\cV$ is a filter on $\A$, then $U_{\cV,\infty}^-$ stabilizes $M_\cV$. 
\end{Remark}

\subsection{Parahoric subgroup associated with a point}\label{ss_para_sbgp}

Let $\cM=\{L_\Z(\lambda)\otimes \cF'\mid \lambda\in X^+\}\cup \{(\cU_\cF,\Ad)\}$.  If $M\in \cM$, then $M$ is naturally a $G$-module, by \cite[Theorem 7.48 or Corollary 7.54]{marquis2018introduction}.

\begin{Definition}\label{d_parah}
Let $x\in \A$. The \textbf{parahoric subgroup associated with $x$} is \[G_x=P_x:=\{g\in G\mid g.M_x=M_x, \forall M\in \cM\}.\]
\end{Definition}

As we shall see (Theorem~\ref{t_G_x}), this definition is independent of the choice of $\omega'$.

 For $M\in \cM$, we denote by $\sP(M)$ the set of weights of $M$, i.e $\sP(M)$ is the set of $\lambda\in X$ such that there exists $m\in M$ such that $h.m=\lambda(h) m$, for all $h\in \fh$. Equivalently, by \cite[(7.23) or Corollary 7.54]{marquis2018introduction}, this is the set of $\lambda\in X$ such that there exists $m\in M$ such that $t.m=\lambda(t)$, for $t\in T$.
  Let $\mu\in \sP(M)$ and   $r,r'\in \R$. We have $M_{\mu,r}=M_{\mu,r'}$ if and only $\cF'_{\omega'\geq r}=\cF'_{\omega'\geq r'}$, if and only if $r=r'$.

Let $\cV$ be a filter on $\A$. Let $M\in \cM$ and $\mu\in \sP(M)$.  We define $\cD(\mu,\cV,\omega
')$ as the set of half-spaces delimited by $\mu^{-1}(\{a\})$, for $a\in \omega'(\cF')$, and containing $\cV$.

\begin{Lemma}\label{l_char_mod}
Let $M\in \cM$, $\mu\in \sP(M)$ and $\Omega,\Omega'$ be two non-empty subsets of $\A$. Then $M_{\mu,\Omega}=M_{\mu,\Omega'}$ if and only if $\cD(\mu,\Omega,\omega')=\cD(\mu,\Omega',\omega')$.
\end{Lemma}

\begin{proof}
By definition, if $\cD(\mu,\Omega,\omega')=\cD(\mu,\Omega',\omega')$, then $M_{\mu,\Omega}=M_{\mu,\Omega'}$. Conversely, assume that $M_{\mu,\Omega}=M_{\mu,\Omega'}$. Then $M_{\mu,\Z}\otimes \cF'_{\geq f'_{\Omega}(\mu)}=M_{\mu,\Z}\otimes \cF'_{\geq f'_{\Omega'}(\mu)}$ and thus $\cF'_{\geq f_\Omega'(\mu)}=\cF'_{\geq f'_{\Omega'}(\mu)}$.  If $f'_{\Omega}(\mu)=+\infty$, then $f'_{\Omega'}(\mu)=+\infty$ and hence $\cD(\mu,\Omega,\omega')=\cD(\mu,\Omega',\omega')$. Assume now that $f'_{\Omega}(\mu)<+\infty$. Then $\cF'_{\geq f'_{\Omega}(\mu)}=\cF'_{\geq f'_{\Omega'}(\mu)}$ and hence as $\omega'(\cF')$ is dense in $\R$, we deduce that  $f'_\Omega(\mu)=f'_{\Omega'}(\mu)$. Hence by definition of $f'_\Omega(\mu)$ and $f'_{\Omega'}(\mu)$, we deduce that $\cD(\mu,\Omega,\omega')=\cD(\mu,\Omega',\omega')$, which proves the lemma.
\end{proof}

\begin{Lemma}\label{l_act_N_mdles}
Let $n\in N$ and $\cV$ be a filter on $\A$. Let $M\in \cM$ and $\mu\in X$. Let $v=\nu^v(n)\in W^v$ and $\bv=\nu(n)\in W^a$.  Then $n.M_{\mu,\cV}=M_{v(\mu),\bv.\cV}$.
\end{Lemma}

\begin{proof}
Let $t\in T$. Let $\Omega\in \cV$.  Then $M_{\mu,\Omega}=M_{\mu,\Z}\otimes \cF'_{\omega'\geq f'_\Omega(\mu)}$. By \cite[(7.23)]{marquis2018introduction} or \eqref{e_Ad_T}, we have: \begin{align*}t.M_{\mu,\Omega}=\mu(t)(M_{\mu,\Z}\otimes \cF'_{\geq f'_\Omega(\mu)})
&=M_{\mu,\Z}\otimes \cF'_{\geq f'_\Omega(\mu)+\omega(\mu(t))}
\\
&=M_{\mu,\Z}\otimes \cF'_{\geq f'_{\nu(t).\Omega}(\mu)}=M_{\mu,\nu(t).\Omega}.\end{align*}

 Let $n_0\in N_0'=\fN_0(\cF)$. By Proposition~\ref{p_d_nu}, we have $\nu(n_0)=\nu^v(n_0)$.  Then by \cite[Proposition 4.18 (1)]{marquis2018introduction}, we have: \[n_0.M_{\mu,\Omega}=n_0.M_{\mu,\Z}\otimes \cF'_{\geq f'_\Omega(\mu)}=M_{\nu^v(n_0)(\mu),\Z}\otimes \cF'_{\geq f'_{\Omega}(\mu)}=M_{\nu^v(n_0)(\mu),\nu(n_0).\Omega}.\] As $N=N_0'\cdot T$ (by \eqref{e_dec_N}), we deduce the result when $\cV=\Omega$  is a non-empty set. If $\cV$ is an arbitrary filter on $\A$, then $n.M_{\mu,\cV}=n.\bigcup_{\Omega\in \cV} M_{\mu,\Omega}=\bigcup_{\Omega\in \cV} M_{v(\mu),\bv.\Omega}=M_{v.\mu,\bv.\cV}$.
\end{proof}

Let $\Lambda'=\omega'(\cF')$. For $\sP\subset X$, we define $\cl^{\sP}_{\Lambda'}$ as follows. Let $\cV$ be a filter on $\A$. Then $\cl^\sP_{\Lambda'}(\cV)$ is the filter on $\A$ generated by the set of sets of the form $\bigcap_{\alpha\in \sP} D(\alpha,\lambda_\alpha)$, where $(\lambda_\alpha)\in (\Lambda')^{\sP}$ is such that $D(\alpha,\lambda_\alpha)\supset \cV$, for all $\alpha\in \sP$.

\begin{Lemma}\label{l_Rou4.11_N}
\begin{enumerate}
\item Let $M\in \cM$ and $\Omega$ be a non-empty subset of $\A$.  Let $N_\Omega(M)=\{n\in N\mid n.M_{\Omega}=M_\Omega\}$. Then $N_\Omega(M)$ is the stabilizer of  $\cl^{\sP(M)}_{\Lambda'}(\Omega)$ in $N$.

\item Let $x\in \A$. Then $\bigcap_{M\in \cM} N_x(M)$ is the fixator of $x$ in $N$.

\end{enumerate}
\end{Lemma}

\begin{proof}
(1) Let $\sP=\sP(M)$.  We have $M_\Omega=\bigoplus_{\mu\in \sP}M_{\mu,\Omega}$. Let $n\in N$. Let $v=\nu^v(n)$ and $\bv=\nu(n)$. By Lemma~\ref{l_act_N_mdles}, we have $n.M_{\Omega}=\bigoplus_{\mu\in \sP}M_{v.\mu,\bv.\Omega}$. Then: \begin{align}n.M_{\Omega}=M_{\Omega}&\Leftrightarrow M_{v.\mu,\Omega}=M_{v.\mu,\bv.\Omega}, \forall \mu\in \sP\\
&\Leftrightarrow M_{\mu,\Omega}=M_{\mu,\bv.\Omega},\forall \mu\in \sP &\text{ by \cite[Lemma 4.16]{marquis2018introduction}}\\
&\Leftrightarrow \cD(\mu,\Omega,\omega')=\cD(\mu,\bv.\Omega,\omega'),\forall \mu\in \sP\label{e_st_D} &\text{by Lemma~\ref{l_char_mod}}\\
& \Leftrightarrow \cl^{\sP}_{\Lambda'}(\Omega)=\cl^{\sP}_{\Lambda'}(\bv.\Omega),
\end{align} which proves (1).

(2) Let $n\in \bigcap_{M\in \cM} N_x(M)$.  By  \eqref{e_st_D} applied with $\Omega=\{x\}$, if $n\in \bigcap_{M\in \cM}N_x(M)$, then \begin{equation}\label{e_cond_D}
\cD(\mu,\{x\},\omega')=\cD(\mu,\{n.x\},\omega'), \forall \mu\in \bigcup_{M\in \cM} \sP(M).
\end{equation}

Let $\mu\in X$ and $y\in \A$. Then $\cD(\mu,\{y\},\omega')=\{D(\mu,a)\mid a\in \Lambda', \mu(y)+a\geq 0\}$. Therefore $\cD(\mu,\{x\},\omega')=\cD(\mu,\{n.x\},\omega')$ if and only if $\{a\in \Lambda'\mid \mu(x)+a\geq 0\}=\{a\in \Lambda'\mid \mu(n.x)+a\geq 0\}$ if and only if $\mu(x)=\mu(n.x)$, since $\Lambda'$ is dense in $\R$.

As $\mu\in \sP(L_\Z(\mu)\otimes \cF')$, for every $\mu\in X^+$, \eqref{e_cond_D} implies that $\mu(x)=\mu(n.x)$, for every $\mu\in X^+$. As $X^+$ generates $X$ as a $\Z$-module, we deduce that $\mu(x)=\mu(n.x)$, for every $\mu\in X$ and thus that $x=n.x$.

Conversely, let  $n\in N_x$. Let $M\in \cM$. Then by Lemma~\ref{l_act_N_mdles}, we have $n.M_{\mu,x}=M_{\nu(n).\mu,x}$, for every $\mu\in X$. Therefore $n.M_x=n.\bigoplus_{\mu\in X} M_{\mu,x}=\bigoplus_{\mu\in \nu^v(n).X} M_{\mu,x}=M_x$ and hence $n\in \bigcap_{M\in \cM} N_x(M)$, which proves (2). \end{proof}

Let $\cM_1\subset \cM$ be non-empty and $\cV$ be a filter on $\A$. We set: \begin{equation}
P_\cV(\cM_1)=\{g\in G\mid g.M_\cV=M_\cV,\forall M\in \cM_1\}.
\end{equation}

\begin{Lemma}\label{l_Rou4.11_2}
 Let $\cV$ be a filter on $\A$ and $\cM_1\subset \cM$ be non-empty. Then $P_\cV(\cM_1)\cap (U^+\cdot  N)=(P_\cV(\cM_1)\cap U^+)\cdot (P_\cV(\cM_1)\cap N)$.

\end{Lemma}

\begin{proof}
We begin by proving the case where $\cM_1=\{M\}$, for some $M\in \cM$.   Let $g\in P_\cV\cap U^+\cdot N$. Write $g=un$, with $u\in U^+$ and $n\in N$. Let us prove that $n\in P_\cV$. Let $\Omega\in \cV$ be such that $g\in P_\Omega$. Let $\lambda\in X$ and  $m\in M_{\lambda,\Omega}$. Set $m'=n.m$ and $w=\nu^v(n)\in W^v$.  By Lemma~\ref{l_act_N_mdles}, $m'\in M_{w.\lambda,\nu(n).x}\subset M_{w.\lambda}$.

We have $U^+\subset \widehat{\cU}^+_\cF$. Write $u=\sum_{\mu\in Q_+} u_\nu$, with $u_\mu\in \cU_{\mu,\cF}$, for $\mu\in Q_+$.

First assume that $M$ is a highest weight module with highest weight $\lambda$. Then by Lemma~\ref{l_cpltd_actn}, we have: \begin{equation}\label{e_dec_unm'}
un.m=u.m'=\sum_{\mu\in Q_+} u_\mu.m'\in M_\Omega
\end{equation} and by \eqref{e_wghts_L}, we have $u_\mu.m'\in M_{\mu+w.\lambda}$, for $\mu \in Q_+$. Therefore in  \eqref{e_dec_unm'}, $u_\mu.m'$ is the $\mu+w.\lambda$ component of $un.m$ in the decomposition $M_\Omega=\bigoplus_{\nu\in Q_+}M_{\nu,\Omega}$, for $\mu\in X$. Consequently, $u_\nu.m'\in M_{\nu+w.\lambda,\Omega}$, for $\nu\in Q_+$. In particular, $u_0.m'\in M_{w.\lambda,\Omega}$. But $u_0=1$ (by Lemma~\ref{l_twst_exp}) and thus $m'\in M_{w.\lambda,\Omega}$. Therefore $n.M_{\lambda,\Omega}\subset M_\Omega$ for $\lambda\in X$ and hence $n.M_\Omega\subset M_\Omega$, which proves that $n\in N\cap P_\Omega\subset N\cap P_\cV$.

Assume now $M=\cU_\cF$. Then as $u_0=1$, we have $u.m'=m'+m''$, where $m''\in \bigoplus_{\nu\in Q_+\setminus \{0\}} \cU_{w.\lambda+\nu}$. Therefore as $u.m'=un.m\in \cU_\Omega$, we have $m'\in \cU_{w.\lambda,\Omega}$. Therefore $n.\cU_{\lambda,\Omega}\subset \cU_\Omega$, which prove  that $n\in N\cap P_\Omega\subset N\cap P_\cV$.

In both cases, $u=un.n^{-1}\in U^+\cap P_\cV$ so  we proved that $P_\cV\cap U^+\cdot N\subset (P_\cV\cap U^+)\cdot (P_\cV\cap N)$. The converse inclusion is obvious and thus we have the result when $\cM_1=\{M\}$.

Let now $\cM_1\subset \cM$ be non-empty.
 We have $P_\cV(\cM_1)\cap (U^+\cdot  N)=\bigcap_{M\in \cM_1} (P_\cV(M)\cap U^+\cdot N)=\bigcap_{M\in \cM_1} (P_\cV(M)\cap U^+)\cdot (P_\cV(M)\cap N)\subset (P_\cV(\cM_1)\cap U^+)\cdot (P_\cV(\cM_1)\cap N)$, by uniqueness in the Birkhoff decomposition (\eqref{e_Birkhoff})). As the converse inclusion is clear, we deduce the result.
\end{proof}

\begin{Theorem}\label{t_G_x}(see \cite[Proposition 4.14 (P3)]{rousseau2016groupes})
Let $x\in \A$. Let $G_x=\langle U_{x,\infty}^+,U_{x,\infty}^-, N_x\rangle$. Then $G_x=P_x(\cM):=\{g\in G\mid g.M_x=M_x,\forall M\in \cM\}$, (for the notation of Notation~\ref{n_modules}) and we have:  \[G_x=U_{x,\infty}^+\cdot U_{x,\infty}^-\cdot N_x=U_{x,\infty}^-\cdot U_{x,\infty}^+\cdot N_x.\]
\end{Theorem}

\begin{proof}
Let $P=P_x(\cM)$. Then $U_{x,\infty}^+, U_{x,\infty}^-, N_x\subset P$, by Lemmas~\ref{l_UVp_stable} and \ref{l_U_v_stb_M_v} and Remark~\ref{r_U_v_stb_M_v}. Thus $U_{x,\infty}^+\cdot U_{x,\infty}^-\cdot N_x\subset P$.

Let $g\in P$. By Proposition~\ref{p_Iwa}, we can write $g=u_+n v$, with $u_+\in U^+$, $n\in N$ and $v\in U_{x,\fin}$. By Proposition~\ref{p_GR08_3.4},  $U_{x,\fin}\subset \langle U_{x,\infty}^+,U_{x,\infty}^-\rangle \subset P$, and thus $u_+n\in P$. By Lemma~\ref{l_Rou4.11_2} (2), we have $u_+\in U^+\cap P$ and $n\in N\cap P$. By Lemma~\ref{l_Rou4.11_U}, $u_+\in U_{x,\infty}^+$  and by  Lemma~\ref{l_Rou4.11_N} (2),  $n\in N_x$.  As $N_x$ normalizes $U_{x,\fin}$, we can write $nv=v'n$, with $v'\in U_{x,\fin}$.

By Proposition~\ref{p_GR08_3.4}, we can write $v'=v_+v_-n'$, with $v_+\in U_{x,\infty}^+$, $v_-\in U_{x,\infty}^-$ and $n'\in N_x$. Then $g=u_+v_+v_- n'n\in \cdot U_{x,\infty}^+\cdot U_{x,\infty}^-\cdot N_x$.  Therefore $P\subset U_{x,\infty}^+\cdot U_{x,\infty}^-\cdot N_x\subset P$, which proves that $P=U_{x,\infty}^+\cdot U_{x,\infty}^-\cdot N_x$. By symmetry of the roles of $+$ and $-$, we deduce the theorem.
\end{proof}

\chapter{Study of the masure associated with a split Kac--Moody group over a valued field}\label{C_prop}

Recall that $G=\fG(\cF)$, where $\fG$ is  a split minimal Kac--Moody group and $\cF$ is a field equipped with a non-trivial valuation $\omega:\cF\rightarrow \R\cup \{+\infty\}$.  In this chapter, we define the masure of $G$ and we prove that it is an abstract masure in the sense of Definition~\ref{d_w_mas}. To that end, we proceed as follows:\begin{enumerate}

\item We define $\I$ as $G\times \A/\sim$, for some equivalence relation $\sim$ on $G\times \A$ defined in such a way that the fixator of every element $x$ of $\A$ is $G_x$, as defined in Section~\ref{s_parahoric}.

\item We prove that $\I$ satisfies~\ref{a_oco}. To that end, we use the notion of ``good fixators'' introduced in \cite{gaussent2008kac}:  a filter is said to have a good fixator if its fixator satisfies certain decompositions (see Definition~\ref{d_good_fixators}). Actually,    we will see that every filter has a good fixator (see Theorem~\ref{t_fltrs_gd_fix}), but this notion will be useful to prove~\ref{a_oco}.

\item We prove that $\I$ satisfies \ref{a_sc}. 

	\item We prove that $\I$ satisfies \ref{a_wma3} and~\ref{a_ma3}: \ref{a_wma3} is obtained by translating the (vectorial) Bruhat and the Birkhoff decompositions of $G$ in terms of masures and~\ref{a_ma3} is obtained by using Levi decompositions and the affine Bruhat decomposition in certain reductive groups. 

\item We obtain that $\I$ satisfies~\ref{a_ma2}, when $\cS$ is positively cofree by using Theorem~\ref{t_MA2}.

\item When $\cS$ is not positively cofree, we consider an extension $\tilde{\cS}$ of $\cS$ which is cofree and then prove that the central extension $\fG_{\tilde{\cS}}\rightarrow  \fG_{\cS}$ induces an extension of masures $\tilde{\I}\twoheadrightarrow \I$, where $\tilde{\I}$ is the masure of $\fG_{\tilde{\cS}}(\cF)$ (see Section~\ref{s_ext_mas}). We then deduce \ref{a_ma2} for $\I$ from \ref{a_ma2} for $\tilde{\I}$.

\end{enumerate} 

\section{Definition of the masure associated with $G$}\label{ss_d_mas}

We define an equivalence relation $\sim$ on $G\times \A$ as follows: $(g,x)\sim (h,y)$ if and only if there exists $n\in N$ such that $y=\nu(n).x$ and $g^{-1}.h.n\in G_x$. The \textbf{masure} $\I=\I(\fG,\fT,\cF,\omega)$  is the set $G\times \A /\sim$. For $(g,x)\in G\times \A$, we denote by $[g,x]$ or by $g.x$ its image in $\I$, depending on the context. We define an action of $G$ on $\I$ by $g.[g',x]=[gg',x]$ if $g,g'\in G$ and $x\in \A$.

\begin{Lemma}\label{l_emb_A_I}
\begin{enumerate}
\item Let $x\in \A$. Then $G_x\cap N=N_x$. 

\item Let $\iota:\A\rightarrow \I$ be defined by $\iota(x)=[1,x]$, for $x\in \A$. Then $\iota$ is injective.
\end{enumerate}

\end{Lemma}

\begin{proof}
(1) By definition of $G_x$, we have $N_x\subset G_x\cap N$. Let $n\in G_x\cap N$. Using Theorem~\ref{t_G_x}, we write $n=u_+u_- n'$ with $u_+\in U^+$, $u_-\in U^-$ and $n'\in N_x$. Then $u_-n'n^{-1}u_+^{-1}=1$ and by uniqueness in the Birkhoff decomposition (\eqref{e_Birkhoff}), we have $n=n'$. Therefore $n\in N_x$ and we have (1).

(2) Let $x,y\in \A$ be such that $\iota(x)=\iota(y)$. Then $[1,x]=[1,y]$ and by definition of $\sim$, there exists $n\in N$ such that $\nu(n).x=y$ and $n\in G_x$. By (1), we have $n\in N_x$ and thus $x=y$, which proves (2).
\end{proof}

 Using Lemma~\ref{l_emb_A_I}, we regard $\A$ as a subset of $\I$, via the embedding $\iota$. By definition of $\sim$, we have $[n,x]=[1,\nu(n).x]$, for all $x\in \A$ and $n\in N$. Therefore $N$ stabilizes $\A$ and the action on $\A$ induced by the action on $\I$ coincides with the one of $\nu$. By definition of $\sim$, we have: \begin{equation}\label{e_def_sim}
\forall x,y\in \A,\forall g\in G, (g.x=y)\Rightarrow (\exists n\in N\mid g.x=n.x=y).
\end{equation}

 An \textbf{apartment}\index{apartment} of $\I$ is then  a set of the form $g.\A$, for some $g\in G$.

  If $\cV$ is a filter on $\A$, we define: \[G(\cV\curvearrowright \A)=\{g\in G\mid g.\cV\subset \A\}.\]\index[notation]{g@$G(\cV\curvearrowright \A$} This subgroup is denoted $G(\cV\subset  \A)$ in \cite{rousseau2016groupes}.

 \begin{Lemma}\label{l_TF_pt}
 \begin{enumerate}
 \item Let $x\in \A$. Then $G(\{x\} \curvearrowright\A)=N \cdot G_x$.

 \item Let $\cV$ be a filter on $\A$. Then  we have: $G(\cV\curvearrowright \A)=\bigcup_{\Omega\in \cV}(\bigcap_{x\in \Omega} N \cdot G_x)$.
 \end{enumerate}
 \end{Lemma}

 \begin{proof}
(1)  As $N$ stabilizes $\A$, we have  $N\cdot G_x\subset G(\{x\}\cur \A)$. Let $g\in G(\{x\}\cur \A)$. Let $y=g.x$. By definition, there exists $n\in N$ such that $y=n.x$ and $n^{-1}g\in G_x$ and hence $g\in N\cdot  G_x$, which proves (1).

 (2) If $\Omega$ is  a non-empty subset of $\A$, we have $G(\Omega\curvearrowright \A)=\bigcap_{x\in \Omega} G(\{x\}\cur \A)$. Let $g\in G$. We have:
\[ g\in G(\cV\curvearrowright \A)\Leftrightarrow \exists \Omega\in \cV\mid g.\Omega\subset \A \Leftrightarrow
 \exists \Omega\in \cV\mid g\in G(\Omega\curvearrowright \A)\] and hence $G(\cV\curvearrowright \A)=\bigcup_{\Omega\in \cV} G(\Omega\curvearrowright \A)$

 \end{proof}

\begin{Lemma}\label{l_Fx_Ain}
 Let $\A_{\ines}=\bigcap_{\alpha\in \Phi} \ker(\alpha)= \bigcap_{i\in I} \ker(\alpha_i)$ and $\cV$ be a filter on $\A$. Then we have $G_\cV=G_{\cV+\A_{\ines}}$.

 \end{Lemma}

\begin{proof}
We begin by the case where $\cV=\{x\}$, for some $x\in \A$. Let $y\in x+\A_{\ines}$.   Then we have $U_{\alpha,y}=U_{\alpha,x}$ for all $\alpha\in \Phi$, since $\alpha(x)=\alpha(y)$.  Let $n\in N_x$ and $\bw\in W$ be the automorphism of $\A$ induced by $n$. Write $\bw=a+w$, where $a\in Y\otimes \Lambda$ and $w\in W^v$. Then we have $a=x-w.x$. As $W^v$ fixes $\A_{\ines}$, we deduce $y-w.y=a$ and hence $\bw$ fixes $y$. In other words, $n$ fixes $y$ and we have  $N_x\subset N_y$. By symmetry, $N_x=N_y$. Moreover, by definition, $U_{x,\infty}^+=U_{x+\A_{\ines},\infty}^+$ and $U_{x,\infty}^-=U_{x+\A_{\ines},\infty}^-$ and thus $G_x=G_y$, by Theorem~\ref{t_G_x}. Let now $\Omega$ be a nonempty set. Then $G_\Omega=\bigcap_{x\in \Omega}G_x=\bigcap_{x\in\Omega}\bigcap_{y\in x+\A_{\ines}} G_y=G_{\Omega+\A_{\ines}}$.

Let now $\cV$ be  a filter on $\A$.  Let $\Omega$ be a subset of $\A$. Then $\Omega\in \cV+\A_{\ines}$ if and only if there exists $\Omega'\in \cV$ such that  $\Omega=\Omega'+\A_{\ines}$. Therefore: \[G_{\cV+\A_{\ines}}=\bigcup_{\Omega\in \cV+\A_{\ines}}G_\Omega=\bigcup_{\Omega'\in \cV} G_{\Omega'+\A_{\ines}}=G_\cV.\]
\end{proof}

\section{Filters with good fixators}

In this section, we introduce the notion of filters having a good fixator. This will be useful to prove axiom~\ref{a_oco}. We will eventually prove that every filter on $\A$ has a good fixator (see Theorem~\ref{t_fltrs_gd_fix}).

\subsection{Definition and first consequences}

Assume that $G$ is reductive. Let $\cV$ be a filter on $\A$. Then the fixator  of $\cV$ in $G$ is \[G_{\cV}=N_{\cV}\cdot U_{\cV,\fin}^+\cdot U_{\cV,\fin}^-=N_{\cV}\cdot U_{\cV,\fin}^-\cdot U_{\cV,\fin}^+, \] where $N_\cV$ (resp. $U_{\cV,\fin}^+$, $U_{\cV,\fin}^-$) is the fixator of $\cV$ in $N$ (resp. in $U^+$, in $U^-$), by \cite[7.1.8, 7.1.11 and 7.4.8]{bruhat1972groupes} and by symmetry of the roles of $+$ and $-$. Moreover by \cite[6.4.2 and 6.4.9]{bruhat1972groupes}, we have $U_{\cV,\fin}^\epsilon=\prod_{\alpha\in \Phi_\epsilon} U_{\alpha,\cV}$, for any  $\epsilon\in \{-,+\}$ and any order on $\Phi_\epsilon$, and  by \cite[Corollaire 7.4.9]{bruhat1972groupes}, we have $G(\cV\cur \A)=N\cdot G_{\cV}$.

Inspired by these properties, Gaussent and Rousseau introduced in \cite[4.1]{gaussent2008kac} the notion of ``good fixator''.

\begin{Definition}\label{d_good_fixators}
Let $\cV$ be a filter on $\A$. We consider the following conditions:

\begin{enumerate}[label=\blue{(GF$+$)}]
\item\label{a_GF+} $G_{\cV}=U_{\cV,\infty}^+\cdot U_{\cV,\infty}^-\cdot N_\cV$.\axiom{GF+@\ref{a_GF+}}
\end{enumerate}
\vspace{-5mm}
\begin{enumerate}[label=\blue{(GF$-$)}]
\item\label{a_GF-} $G_{\cV}=U_{\cV,\infty}^-\cdot U_{\cV,\infty}^+\cdot N_\cV$.\axiom{GF-@\ref{a_GF-}}
\end{enumerate}
\vspace{-5mm}
\begin{enumerate}[label=\blue{(TF)}]
\item\label{a_TF} $G(\cV\cur \A)=N\cdot  G_{\cV}$.\axiom{TF@\ref{a_TF}}
\end{enumerate}

We say that  $\cV$ has a \textbf{good fixator}\index{good fixator} (resp. \textbf{half-good fixator}\index{half-good fixator}, resp. \textbf{transitive fixator}\index{transitive fixator}) if it satisfies \ref{a_GF+}, \ref{a_GF-} and \ref{a_TF} (resp \ref{a_TF} and \ref{a_GF+} or \ref{a_GF-}, resp. \ref{a_TF}). 
\end{Definition}

As we shall see (see Lemma~\ref{l_TF}), a filter has a transitive fixator if and only if $G_{\cV}$ acts transitively on the set of apartments containing $\cV$. 

Let $x\in \A$. Then by Theorem~\ref{t_G_x}, $\{x\}$ satisfies \ref{a_GF+} and \ref{a_GF-}, and by Lemma~\ref{l_TF_pt}, $\{x\}$ satisfies \ref{a_TF}. Therefore every singleton of $\A$ has a good fixator.

\begin{Lemma}\label{l_u_fix}
\begin{enumerate}
\item Let $x\in \A$ and $\epsilon\in \{-,+\}$.  Then $G_x\cap U^\epsilon=U_{x,\infty}^{\epsilon}$.

\item Let $\cV$ be a filter on $\A$ and $\epsilon\in \{-,+\}$. Then $U_{\cV,\infty}^{\epsilon}=U_{\cl^\Delta(\cV),\infty}^{\epsilon}=U_{\cl^\Delta(\cV)+\epsilon\overline{C^v_f},\infty}^{\epsilon}$ and this group fixes $\cl^\Delta(\cV)+\epsilon\overline{C^v_f}$. 
\end{enumerate}

\end{Lemma}

\begin{proof}
1) By symmetry, we can assume $\epsilon=+$. By definition of $G_x$, we have $U_{x,\infty}^+\subset G_x$. Conversely, take $u\in U^{+}\cap G_x$. By Theorem~\ref{t_G_x}, we can write $u=u_+u_-n$, with $u_+\in U_{x,\infty}^+$, $u_-\in U_{x,\infty}^-$ and $n\in N_x$. Then we have $u_-n u^{-1}u_+=1$ and hence by uniqueness in the Birkhoff decomposition (\eqref{e_Birkhoff}), we deduce $u=u_+$, which proves that $u\in U_{x,\infty}^{+}$. Therefore $G_x\cap U^+$.

2) The equality of the groups follows from Corollary~\ref{c_Uma_cl} applied with $\Psi=\Delta_+$ (and its analogue for $\Delta_-$) since $U_{\cV',\infty}^{\epsilon}=U^{ma \epsilon}_{\cV'}\cap U^{\epsilon}$, for any filter $\cV'$ on $\A$. 

Let $u\in U_{\cl^\Delta(\cV)+\overline{C^v_f},\infty}^{\epsilon}$. Let $\Omega\in \cl^\Delta(\cV)+\overline{C^v_f}$ be such that $u \in U_ {\Omega,\infty}^{\epsilon}$. Let $x\in \Omega$. Then by definition of $G_x$, $u\in G_x$ and thus $u$ fixes $x$. This proves that $u$ fixes $\Omega$ and proves the lemma.
\end{proof}

\begin{Remark}\label{r_chng_ord_fix}
 Let $\cV$ be a filter on $\A$.

\begin{enumerate}
\item  As $G_\cV$ is  a subgroup of $G$, $G_{\cV}=G_{\cV}^{-1}$ and  \ref{a_GF+} can be rewritten $G_{\cV}=N_\cV \cdot U_{\cV,\infty}^-\cdot  U_{\cV,\infty}^+$ and similarly for \ref{a_GF-}. However, $G(\cV\curvearrowright \A)$ is not a subgroup of $G$ and it can differ from $G_\cV\cdot  N$  in general.

\item The filter $\cV$ satisfies \ref{a_GF+} (resp. \ref{a_GF-}, resp. \ref{a_TF}) if and only if $G_{\cV}\subset U_{\cV,\infty}^+\cdot  U_{\cV,\infty}^-\cdot N_{\cV}$ (resp. $G_{\cV}\subset U_{\cV,\infty}^-\cdot U_{\cV,\infty}^+ \cdot N_{\cV}$, resp. $G(\cV\curvearrowright \A)\subset N\cdot G_\cV$). 

\end{enumerate}
\end{Remark}

\section{Axioms \ref{a_oco} and \ref{a_ma1}}

In this section, we prove that $\I$ satisfies axioms \ref{a_oco} and \ref{a_ma1}. These two axioms are intricate since on one hand, axiom~\ref{a_oco} requires to define enclosures, which uses \ref{a_ma1}. On the other hand, we prove \ref{a_ma1} by transporting the structure of $\A$ to $g.\A$ via $g$, for $g\in G$. However to prove that it is well-defined, we use the fact that the stabilizer of $\A$ in $G$ is $N$ (see Proposition~\ref{p_N_stab_A}), which we show by using a weak form of \ref{a_oco}.

\subsection{Weak version of \ref{a_oco} and \ref{a_ma1}}

\begin{Lemma}\label{l_wd_retrac_1}
Let $\epsilon\in \{+,-\}$, $u\in U^\epsilon$ and $x\in \A$. Assume that $u.x\in \A$. Then $u.x=x$ and $u\in U^{\epsilon}_{x,\infty}$. 
\end{Lemma}

\begin{proof}
We assume that $\epsilon =+$, by symmetry. By \eqref{e_def_sim}, we can write $u.x=n.x$, for some $n\in N$. Then $u^{-1}n\in G_x$. By Theorem~\ref{t_G_x}, we can write $u^{-1}n=v_{+} v_{-}n'$, with $v_+\in U_{x,\infty}^+$, $v_{-}\in U_{x,\infty}^-$ and $n'\in N_x$. Then $uv_+ v_- n'n^{-1}=1$,  hence $v_-n'n^{-1}=(uv_+)^{-1}$. Consequently $v_-n'n^{-1} uv_+=1$. Using the Birkhoff decomposition (\eqref{e_Birkhoff}), we deduce $uv_+=1$ and thus $u\in U_{x,\infty}^+$, which proves the lemma. 
\end{proof}

\begin{Lemma}\label{l_sfct_cnd_GF}
Let $\cV$ be a filter on $\A$. We assume that $G(\cV\cur\A)\subset N \cdot  U_{\cV,\infty}^-\cdot  U_{\cV,\infty}^+$. Then $\cV$ satisfies \ref{a_GF+} and \ref{a_TF}.
\end{Lemma}

\begin{proof}
Let $g\in G_\cV$. Then $g\in G(\cV\curvearrowright \A)$ and thus we can write $g=nu_- u_+$, with $n\in N$, $u_+\in U_{\cV,\infty}^+$ and $u_-\in U_{\cV,\infty}^-$. Let $\Omega_1\in \cV$ be such that $g\in G_{\Omega_1}$. Let $\Omega_+,\Omega_-\in \cV$ be such that $u_+\in U_{\Omega_+,\infty}^+$ and $u_-\in U_{\Omega_-,\infty}^-$. Let $\Omega=\Omega_-\cap \Omega_+\cap \Omega_1\in \cV$. Let $x\in \Omega$. Then $g.x=x=nu_-u_+.x=n.x$. Therefore $n\in N_\Omega$ and thus $g\in N_{\Omega}\cdot U_{\Omega,\infty}^- \cdot U_{\Omega,\infty}^+\subset N_{\cV}\cdot U_{\cV,\infty}^-\cdot U_{\cV,\infty}^+$. Using Remark~\ref{r_chng_ord_fix} we deduce that $\cV$ satisfies \ref{a_GF+}. Therefore $G(\cV\cur \A)\subset N\cdot  U_{\cV,\infty}^-\cdot U_{\cV,\infty}^+\subset N\cdot G_\cV$, which proves that $\cV$ satisfies \ref{a_TF} by Remark~\ref{r_chng_ord_fix}. 
\end{proof}

The following lemma is a weak version of \cite[Proposition 4.3 4)]{gaussent2008kac}. Actually, we kept only the point that is useful to prove that $\I$ satisfies~\ref{a_oco}. As we will prove that every filter on $\A$ has a good fixator, we do not need the level of details of this proposition.

\begin{Lemma}\label{l_sctr_gd_fx}
Let $x\in \A$, $\epsilon\in \{-,+\}$ and $\cV$ be a filter on $\A$ such that $\{x\}\subset \cV\subset x+\epsilon \overline{C^v_f}$ (i.e $x+\epsilon \overline{C^v_f}\in \cV$ and for all $\Omega\in \cV$, $\Omega$ contains $x$). Then $\cV$ satisfies (GF$\epsilon$) and \ref{a_TF}. 
\end{Lemma}

\begin{proof}
By symmetry, we assume $\epsilon=+$. Let $g\in G(\cV\curvearrowright \A)$.  Let $\Omega\in \cV$ be such that $g\in G(\Omega\curvearrowright \A)$. Up to replacing $\Omega$ by $\Omega\cap (x+\epsilon \overline{C^v_f})$ (which belongs to $\cV$, since $\cV\subset x+\overline{C^v_f}$), we can assume that $\Omega$ is contained in $x+\overline{C^v_f}$.   

As $x\in \Omega$, we have $g\in G(\{x\}\cur \A)=N \cdot G_x=N \cdot U_{x,\infty}^- \cdot U_{x,\infty}^+$ (by Theorem~\ref{t_G_x}). Write $g=nu_-u_+$, with $n\in N$, $u_+\in U_{x,\infty}^+$ and $u_-\in U_{x,\infty}^-$. By Corollary~\ref{c_Uma_cl}, we have $U_{x,\infty}^+=U_{\Omega,\infty}^+$.  For $y\in \Omega$, we have $g.y=n u_-u_+.y=nu_-.y\in \A$.
 Hence $u_-.y\in n^{-1}.\A=\A$. Using Lemma~\ref{l_wd_retrac_1} we deduce $u_-.y=y$ and hence $u_-\in U_{y,\infty}^-$, by Lemma~\ref{l_wd_retrac_1}. Consequently $u_-\in \bigcap_{y\in \Omega} U_{y,\infty}^-$. By Corollary~\ref{c_U_ma_inter}, we have  $\bigcap_{y\in \Omega} U_{y,\infty}^-=U_{\Omega,\infty}^-$ and hence $u_-\in U_{\Omega,\infty}^-$.

 We deduce $g=nu_- u_+\in N\cdot U_{\Omega,\infty}^- \cdot  U_{\Omega,\infty}^+\subset N\cdot U_{\cV,\infty}^-\cdot U_{\cV,\infty}^-$. Using Lemma~\ref{l_sfct_cnd_GF}, we deduce that $\cV$ satisfies \ref{a_GF+} and \ref{a_TF}.  
\end{proof}

\begin{Proposition}\label{p_N_stab_A}
The stabilizer $G(\A\curvearrowright \A)$ of $\A$ in $G$ is $N$. 
\end{Proposition}

\begin{proof}
We already know that $N\subset G(\A\curvearrowright\A)$. Let $g\in G(\A\cur \A)$. Then $g\in G(C^v_f\cur \A)$.   By Lemma~\ref{l_sctr_gd_fx}, we deduce $g\in N \cdot U_{C^v_f,\infty}^- \cdot U_{C^v_f,\infty}^+=N \cdot U_{C^v_f,\infty}^+\subset N\cdot U^+$. We also have $g\in G(-C^v_f\cur \A)$ and thus $g\in N\cdot U^-$. Using the Birkhoff decomposition (\eqref{e_Birkhoff}), we deduce $g\in N$.
\end{proof}

\begin{Corollary}\label{c_Fix_A}
The fixator $G_\A$ of $\A$ in $G$ is $T_0=\fT(\cO)$. 
\end{Corollary}

\begin{proof}
By definition of the affine action of $T$ on $\A$ (see \eqref{e_act_T}), we have $T_0\subset G_\A$. Conversely, take $g\in G_\A$. By Proposition~\ref{p_N_stab_A}, $g\in N$. As $g$ stabilizes $C^v_f$ and as $W^v$ acts simply transitively on the set of positive vectorial chambers (by Proposition~\ref{p_T_cone}), we have $\nu^v(g)=\Id$ and hence $g\in T$, by Proposition~\ref{p_d_nu_v}.  But then by \eqref{e_act_T}, we have $\omega(\chi(g))=0$, for all $\chi\in X$ and hence $g\in T_0$, which proves the corollary.
\end{proof}

\begin{Definition}\label{d_aff_apt_G}
The \textbf{affine apartment} (resp. \textbf{extended affine apartment}) associated with $(\fG,\cF,\omega)$ is $\underline{\A}=(\A,\cS,(\Lambda)_{\alpha\in \Phi}, W^a)$ \index[notation]{a@$\underline{\A}$} (resp. $\underline{\widetilde{\A}}=(\A,\cS,(\Lambda)_{\alpha\in \Phi}, \widetilde{W^a})$\index[notation]{a@$\underline{\widetilde{\A}}$}), where $\cS$ is the Kac--Moody datum defining $\fG$, $\Lambda=\omega(\cF^\times)$, $W^a=W^a((\Lambda)_{\alpha\in \Phi})$ and  $\widetilde{W^a}$ is the extended affine Weyl group defined in Definition~\ref{d_ext_aff_W}. 
\end{Definition}

\begin{Definition/Proposition}\label{dp_MA1_tildeA}
\begin{enumerate}

\item Let $A$ be an apartment. We define $\Isom_{\underline{\widetilde{\A}}}(\A,A)$ as the set of $h|_\A^{A}$ such that $h\in G^u:=\langle U^+,U^-\rangle$ and $h.\A=A$. This equips $A$ with the structure of an apartment of type $\underline{\A}$ and thus $\I$ satisfies \ref{a_ma1} for $\underline{\widetilde{\A}}$. 

\item Let $A$ and $B$ be two apartments of $\I$.  An \textbf{apartment isomorphism} $f$ from $A$ to $B$ is a map of the form $f=f_2\circ f_1^{-1}$, where $f_1\in \Isom_{\underline{\widetilde{\A}}}(\A,A)$ and $f_2\in \Isom_{\underline{\widetilde{\A}}}(\A,B)$. Then for every $g\in G$ such that $g.A=B$, $g|_{A}^{B}$ is an apartment isomorphism and every apartment isomorphism from $A$ to $B$ is of this form.
\end{enumerate}
\end{Definition/Proposition}

\begin{proof}
Let $f_0\in\mathrm{Isom}_{\underline{\widetilde{\A}}}(\A,A)$. Write $f_0=g|_{\A}^{A}$, where $g\in G$. Let $\bw\in \widetilde{W^a}$. Let $n\in N$ be such that $\bw=\nu(n)$, which exists by definition of $\widetilde{W^a}$. Then $f_0\circ \bw=(gn)|_{\A}^{A}\in \mathrm{Isom}_{\underline{\widetilde{\A}}}(\A,A)$. Conversely, let $f\in \mathrm{Isom}_{\underline{\widetilde{\A}}}(\A,A)$. Write $f=h|_{\A}^A$. Then by Proposition~\ref{p_N_stab_A}, there exists $n\in N$ such that $g=hn$ and then $f_0=f\circ \nu(n)$. This proves the proposition.

(2) It suffices to prove it when $A=\A$. This is then a consequence of (1).
\end{proof}

We can equip each apartment with its structure of an apartment of type $\underline{\widetilde{\A}}$ provided by  Definition/Proposition~\ref{dp_MA1_tildeA}. Then $\I$ satisfies \ref{a_ma1}, for $\underline{\widetilde{\A}}$. 

We prove in the next  subsection that $\I$ actually satisfies \ref{a_ma1} for $\underline{\A}$, where $\underline{\A}=(\A,\cS,(\Lambda)_{\alpha\in \Phi},W^a)$, where $W^a=W^a((\Lambda)_{\alpha\in \Phi})$. This is slightly more precise than Definition/Proposition~\ref{dp_MA1_tildeA} by Remark~\ref{r_MA1}. However, this refinement will not be used elsewhere and the reader not interested by this subtlety can skip the next subsection and replace $\underline{\A}$ by $\underline{\widetilde{\A}}$, $G^u$ by $G$ and ``Weyl-isomorphism'' by ``isomorphism'' in the sequel.

\subsection{\ref{a_ma1} for $\underline{\A}$}\label{ss_MA1_A}

In this subsection, we prove that $\I$ satisfies \ref{a_ma1} for $\underline{\A}$, following \cite[5.11]{rousseau2016groupes}. For this, we use the simply connected Kac--Moody group $\fG^A$ associated with $A$.  Let $\cS^{A}$ be the Kac--Moody datum $(A,X^A,Y^A,(\alpha_i^A)_{i\in I},(\alpha_i^{\vee,A})_{i\in I})$, where the $\alpha_i^{\vee,A}$, $i\in I$ are symbols, $Y^A=\bigoplus_{i\in I}\Z \alpha_i^{A,\vee}$, $X^A$ is the $\Z$-dual of $Y^A$ and for all $i\in I$, $\alpha_i^A$ is the unique element of $X^A$ such that $(\alpha_i(\alpha_j^\vee))_{j\in I}=(a_{j,i})_{j\in I}$. Note that contrary to the rest of the paper,  $\cS^A$ need not be free: $\cS^A$ is free if and only $A$ is invertible. Let $\fG^A$ be the Kac--Moody group functor associated to $\cS^A$ in \cite[7.4]{marquis2018introduction} (or in Chapter~\ref{C_splt_KM_grps} if $A$ is invertible). We denote by $\fU^{A,+}$, $G^A$, etc. the standard positive unipotent subgroup of $\fG^A$, $\fG^A(\cF)$, etc.

Let  \[G^u=\langle U^+,U^-\rangle\subset G.\]\index[notation]{g@$G^u$}

\begin{Lemma}\label{l_Gu}

\begin{enumerate}
\item We have $G^u\cap T=\langle \alpha^\vee_i(\cF^\times)\mid i\in I\rangle\subset G$ and $G^u\cap N=\fN_0(\cF)\cdot \langle \alpha^\vee_i(\cF^\times)\mid i\in I\rangle$.

\item We have $\nu(G^u\cap T)=Q^\vee:=\bigoplus_{i\in I} \Lambda \alpha_i^\vee$ and $\nu(G^u \cap N)=W^a$. 
\end{enumerate}

\end{Lemma}

\begin{proof}
(1) Let $i\in I$. By \eqref{e_tilde_si}, \ref{a_KMT6} and \eqref{e_s_i_sq}, $\tilde{s_i}\tilde{s_i}(-r^{-1})=\alpha_i^\vee(r)\in G^u\cap T$, for all $r\in \cF^\times$. Therefore  \begin{equation}\label{e_simp_con}
\langle \alpha^\vee_i(\cF^\times)\mid i\in I\rangle\subset G^u\cap T.
\end{equation}

The map $\psi:(A,Y^A,(\alpha_i^A)_{i\in I},(\alpha_i^{\vee,A})_{i\in I})\rightarrow (A,Y,(\alpha_i)_{i\in I},(\alpha_i^\vee)_{i\in I})$ defined by $\psi(\alpha_i^{\vee,A})=\alpha_i^\vee$ and $\psi(i)=i$ for $i\in I$  is a morphism of Kac--Moody data in the sense of \cite[1.1 6)]{rousseau2016groupes}. By \cite[1.10]{rousseau2016groupes},  $\psi$ induces a group functor morphism $f_\psi:\fG^A\rightarrow \fG$. We denote by $f$ the group morphism $G^A\rightarrow G$ induced by $f_\psi$. We have \[f(\alpha_i^{\vee,A}(r))=\alpha_i^\vee(r), \forall i\in I, \forall r\in  \cF^\times\] and $f$ induces an isomorphism $U^{\epsilon,A}\rightarrow U^{\epsilon}$, for $\epsilon\in \{-,+\}$. Moreover, $f^{-1}(T)=T^A:=\fT^A(\cF)$. We deduce $G^u \cap T=f(\langle U^{+,A},U^{-,A}\rangle \cap T^A)$. By \eqref{e_simp_con} applied with $\fG^A$ instead of $\fG$, we have $\langle U^{+,A},U^{-,A}\rangle \cap T^A=T^A$. Therefore  $f(T^A)=\langle f(\alpha_i^{\vee,A}(\cF^\times))\mid i\in I\rangle)=\langle \alpha_i^\vee(\cF^\times)\mid i\in I\rangle$, which proves the first equality of  (1). 

By \eqref{e_tilde_si}, $\fN_0(\cF)\subset G^u$. Let $n\in N\cap G^u$. By \eqref{e_dec_N}, we can write $n=tn_0$, with $t\in T$ and $n_0\in \fN_0(\cF)$. Then $t=nn_0^{-1}\in G^u\cap T$, which proves (1).

(2) is a consequence of (1) and of Lemma~\ref{l_act_cort}.
\end{proof}

\begin{Definition/Proposition}\label{dp_MA1_A}
\begin{enumerate}

\item Let $A$ be an apartment. We define $\Isom_{\underline{\A}}(\A,A)$ as the set of $h|_\A^{A}$ such that $h\in G^u$ satisfies $h.\A=A$ (where $G^u=\langle U^+,U^-\rangle$). This equips $A$ with the structure of an apartment of type $\underline{\A}$ and thus $\I$ satisfies \ref{a_ma1} for $\underline{\A}$. 

\item Let $A$ and $B$ be two apartments of $\I$. An apartment isomorphism $f$ from $A$ to $B$ is called a \textbf{Weyl-isomorphism} if it can be written $f=f_2\circ f_1^{-1}$, where $f_1\in \Isom_{\underline{\A}}(\A,A)$ and $f_2\in \Isom_{\underline{\A}}(\A,B)$. Then for every $g\in G^u$ such that $g.A=B$, $g|_{A}^{B}$ is a Weyl-isomorphism and every Weyl-isomorphism from $A$ to $B$ is of this form.

\end{enumerate}

\end{Definition/Proposition}

\begin{proof}
(1) By definition, there exists $g\in G$ such that $A=g.\A$. By Proposition~\ref{p_DRJ} and \ref{a_KMT4'}, we can write $g=ht$, with $h \in G^u$ and $t\in T$. Then $A=h.\A$. Set $f_0=h|_{\A}^A$.  Let $\bw\in W^a$. By Lemma~\ref{l_Gu}, there exists $n\in G^u\cap N$ such that $\nu(n)=\bw$. Then $(hn)|_{\A}^A=f_0\circ \bw$. Conversely, let $f\in \Isom_{\underline{\A}}(\A,A)$. Then $f=h'|_{\A}^A$, for some $h'\in G^u$. Then by Proposition~\ref{p_N_stab_A}, $n:=h'^{-1}\circ h\in N\cap G^u$, hence $h'=h\circ \nu(n)\in f_0\circ W^a$, which proves (1).

(2) It suffices to prove it when $A=\A$. This is then a consequence of (1).
\end{proof}

\subsection{Axiom~\ref{a_oco}}

We proved that $\I$ satisfies \ref{a_ma1} for $\underline{\widetilde{\A}}$ and for $\underline{\A}$. The only difference between these two structures lies in the maps that are considered as ``apartment isomorphisms''. In order to distinguish these two structures, we call ``Weyl-isomorphism\index{Weyl-isomorphism}'' an apartment isomorphism for $\underline{\A}$ and simply ``apartment isomorphism'' an apartment isomorphism for $\underline{\widetilde{\A}}$. 

Using \ref{a_ma1}, we can define half-apartments in any apartment of $\I$ (this does not change if we consider $\underline{\A}$ or $\underline{\widetilde{\A}}$ since $W^a$ and $\widetilde{W^a}$ preserve the set of real walls).

\begin{Lemma}\label{l_csq_gd_fix}
Let $\cV$ be a filter on an apartment of $\I$. We assume that for some $g\in G$, $g.\cV$ is contained in  $\A$ and  has a half-good fixator. Let $A, B$ be two apartments containing $\cV$. Then $\cl^{\Delta}_A(\cV)\subset A\cap B$ and  there exists $h\in G^u$ such that $h.A=B$ and $h$ fixes $\cl^{\Delta}_A(\cV)$ pointwise. In particular, $\cl^{\Delta}_A(\cV)$ is independent of the choice of an apartment containing $\cV$ and we simply denote it $\cl^{\Delta}(\cV)$. 
\end{Lemma}

\begin{proof}
Let $\cV'=g.\cV$, $A'=g.A$ and $B'=g.B$. If we find $h'\in G$ such that $h'.A'=B'$ and $h'$ fixes $\cl^{\Delta}_{A'}(\cV')$, then $g^{-1}h'g.A=B$ and $g^{-1}h'g$ fixes $\cl^\Delta_{A}(\cV)$ pointwise. Therefore we can assume $\cV\subset \A$. 

It suffices to take $A=\A$. By symmetry, we can assume that $\cV$ satisfies \ref{a_GF+}. Let $j\in G$ be such that $j.B=\A$. Then $j.\cV\subset \A$ and thus $j\in G(\cV\cur \A)$. By \ref{a_TF}, we can write $j=nj'$, with $j'\in G_{\cV}$.  Using \ref{a_GF+} and the fact that $G_{\cV}=(G_{\cV})^{-1}$, we can write $j'=n' u_-u_+$, with $n'\in N_{\cV}$, $u_-\in U_{\cV,\infty}^-$ and $u_+\in U_{\cV,\infty}^+$. We have $B=j^{-1}.\A=u_+^{-1}u_{-}^{-1}(n')^{-1} n^{-1}.\A=u_+^{-1} u_-^{-1}.\A$, since $N$ stabilizes $\A$. Let $h=u_{+}^{-1}u_-^{-1}$. Then by Lemma~\ref{l_u_fix}, $h$ fixes $\cl^{\Delta}_{\A}(\cV)$ pointwise and thus we get the lemma. 
\end{proof}

\begin{Lemma}\label{l_exstnce_sqnc}
Let $F^v$ be a  vectorial face of $\overline{C^v_f}$. Let $x\in \A$. Let $\Omega$ be a subset of $\A$ such that $x\in \overline{\Omega}$ and $\Omega\subset x+F^v$. Then there exists $(x_k)\in \Omega^\N$ such that $x_k\to x$ and for all $k\in \N$, $\{x_i\mid i\in \llbracket 0,k\rrbracket\}\subset x_{k+1}+F^v$. 
\end{Lemma}

\begin{proof}
Let $S=\supp(x+F^v)$. Then $\Omega\subset S$. Let $(y_k)\in \Omega^\N$ be such that $y_k\to x$.   Set $x_0=y_0$. Let $k\in \N$ and assume that we have constructed $(x_i)_{i\in \llbracket 0,k\rrbracket}$. We have  $x_k\in x+F^v$, or equivalently, $x\in x_k-F^v$. As $x_k-F^v$ is open in $S$ and contains $x$, we have $y_m\in x_k-F^v$, for $m\in \N$ large enough. We choose $m\in \llbracket k+1,+\infty\llbracket$ such that $y_m\in x_k-F^v$ and we set $x_{k+1}=y_m$. Then $(x_i)_{i\in \N}$ satisfies the condition of the lemma.
\end{proof}

\begin{Lemma}\label{l_OCO}
Let $x\in \A$, $\epsilon\in \{-,+\}$ and $\cV$ be a filter on $\A$ such that $\{x\}\subset \overline{\cV}\subset x+\epsilon\overline{C^v_f}$. Then $\cV$ satisfies (GF$\epsilon$) and \ref{a_TF}. 
\end{Lemma}

\begin{proof}
By symmetry, we can assume that $\epsilon=+$. Let $g\in G(\cV\curvearrowright \A)$. Let $\Omega\in \cV$ be such that $g\in G(\Omega\curvearrowright \A)$. Up to replacing $\Omega$ by  $\Omega\cap (x+\overline{C^v_f})$ (which belongs to $\cV$ since $(x+\overline{C^v_f})\in \cV$ by assumption), we can assume that $\Omega$ is contained in $x+\overline{C^v_f}$. 

Let us prove that $g.x\in \A$. Let $(y_k)\in \Omega^{\N}$ be such that $y_k\to x$. Let $F^v$ be a vectorial face of $\overline{C^v_f}$ such that $\{k\in \N\mid y_k\in x+F^v\}$ is infinite. Up to extracting a subsequence of $(y_k)$, we can assume that $y_k\in x+F^v$, for all $k\in \N$. Let $\Omega'=\{y_k\mid k\in \N\}$.  Using Lemma~\ref{l_exstnce_sqnc}, we get a sequence $(x_k)\in (\Omega')^{\N}$ such that for all $k\in \N$, $\{x_i\mid i\in \llbracket 0,k\rrbracket\}\subset x_{k+1}+F^v$.

 For $k\in \N$, we set $\Omega_k=\{x_i\mid i\in \llbracket 1,k\rrbracket\}$. The sequence $(\supp(\Omega_k))$ is increasing and thus it is stationary. Let $k_0\in \N$ be such that $\supp(\Omega_\ell)=\supp(\Omega_{k_0})$, for all $\ell\in \Z_{\geq k_0}$.

  Let  $\ell\in \Z_{\geq k_0}$. Then $g\in G(\Omega_\ell\curvearrowright \A)$ and thus by Lemma~\ref{l_sctr_gd_fx},  we can write $g=n^{(\ell)} u_-^{(\ell)} u_+^{(\ell)}$, 
  where $n^{(\ell)}\in N$, $u_-^{(\ell)}\in U_{\Omega_\ell,\infty}^{-}$ 
  and $u_{+}^{(\ell)}\in U_{\Omega_\ell,\infty}^+$. Then: \[\tilde{n}^{(\ell)}:=(n^{(k_0)})^{-1}n^{(\ell)}=u_-^{(k_0)}u_+^{(k_0)} (u_+^{(\ell)})^{-1} (u_-^{(\ell)})^{-1}\] fixes $\Omega_{k_0}$ pointwise. As it induces an affine automorphism of $\A$, it fixes $\supp(\Omega_{k_0})=\supp(\Omega_\ell)$ pointwise. 
  
  Set $\tilde{g}=(n^{(k_0)})^{-1} g$. Then $\tilde{g}=\tilde{n}^{(\ell)} u_{-}^{(\ell)} u_{+}^{(\ell)}$ and thus $\tilde{g}\in G_{\Omega_\ell}$, for all $\ell\in \Z_{\geq k_0}$. Therefore $\tilde{g}\in G_{\tilde{\Omega}}$, where $\tilde{\Omega}=\{x_{\ell}\mid \ell\in \Z_{\geq k_0}\}$.  We have \[\tilde{g}^{-1}=(u_{+}^{(k_0)})^{-1}(u_-^{(k_0)})^{-1}.\] By Lemma~\ref{l_u_fix}, $u_-^{(k_0)}$ fixes $x_{k_0}-\overline{C^v_f} $ and thus $u_-^{(k_0)}$ fixes $\{x_\ell\mid \ell \in \Z_{\geq k_0}\}\cup \{x\}$. Therefore $u_+^{(k_0)}=(u_-^{(k_0)})^{-1} \tilde{g}$ fixes $\{x_\ell\mid \ell\in \Z_{\geq k_0}\}$ and hence $u_{+}^{(k_0)}$ fixes $x$  by Lemma~\ref{l_u_fix}. Therefore $\tilde{g}.x=u_-^{(k_0)} u_+^{(k_0)}.x=x$ and hence $g.x=n^{(k_0)}\tilde{g}.x=n^{(k_0)}.x\in \A$. This proves that $g\in G(\Omega\cup\{x\} \curvearrowright \A)$.

  By Lemma~\ref{l_sctr_gd_fx}, $\Omega\cup \{x\}$ satisfies (GF$+$) and \ref{a_TF}. Therefore $g\in 
N\cdot U_{\Omega\cup \{x\},\infty}^-\cdot U_{\Omega\cup \{x\},\infty}^+ \subset    N\cdot U_{\Omega,\infty}^{-}\cdot U_{\Omega,\infty}^+\subset N\cdot U_{\cV,\infty}^-\cdot U_{\cV,\infty}^+$. Using Lemma~\ref{l_sfct_cnd_GF}, we deduce the lemma.
\end{proof}

\begin{Proposition}\label{p_SOMA}
The masure $\I$ satisfies~\ref{a_oco}: let $\cV$ be a filter on an apartment $A$ such that there exists a sector $\fq$ based at some $x\in A$ satisfying $\{x\}\subset \overline{\cV}\subset \overline{\fq}$ (i.e $\overline{\fq}\in \cV$ and $x\in \overline{\Omega}$,  for all $\Omega\in \cV$). Then for every apartment $A'$ containing $\cV$, $\cl^{\Delta}_A(\cV)\subset A'$ and there exists a Weyl-  isomorphism from $A$ to $A'$ fixing $\cl^{\Delta}_A(\cV)$. Moreover, this isomorphism is induced by an element of $G^u$. 
\end{Proposition}

\begin{proof}
By definition, there exist $g\in G$ and $\epsilon\in \{-,+\}$ such that $g.\fq=g.x+\epsilon \overline{C^v_f}$. Using Lemma~\ref{l_OCO}, we deduce that $g.\cV$ has a half-good fixator and we conclude with Lemma~\ref{l_csq_gd_fix}. 
\end{proof}

\begin{Lemma}\label{l_TF}
Let $\cV$ be a filter on $\A$. Then $\cV$ has a transitive fixator (i.e satisfies \ref{a_TF}) if and only if $G_\cV$ acts transitively on the set of apartments containing $\cV$. 
\end{Lemma}

\begin{proof}
($\Rightarrow$) Assume that $\cV$ has a transitive fixator. Let $B$ be an apartment containing $\cV$ and let $h\in G$ be such that $h.B=\A$. Then $h\in G(\cV\curvearrowright \A)$ and hence we can write $h=ng$, with $n\in N$ and $g\in G_{\cV}$. Then $g.B=n^{-1}h.B=n^{-1}.\A=\A$ and $g$ fixes $\cV$, which proves that $G_\cV$ acts transitively on the apartments containing $\cV$. 

($\Leftarrow$) Assume that $G_\cV$ acts transitively on the set of apartments containing  $\cV$ and take an apartment $B$ containing $\cV$.  Let $g\in G(\cV\curvearrowright \A)$. Let $B=g.\A$. Then $B$ is an apartment containing $\cV$ and thus there exists $h\in G$ such that $h.\A=B$ and $h\in G_{\cV}$. Then $h^{-1}g.\A=\A$ and hence $h^{-1}g\in N$, by Proposition~\ref{p_N_stab_A}. Therefore $g\in N \cdot  G_{\cV}$, which proves that $\cV$ has a transitive fixator. 
\end{proof}

\section{Fixator of a sector-germ and axiom \ref{a_wma3}}

In this section, we prove that $\I$ satisfies axiom \ref{a_wma3}, by translating the Birkhoff and the Bruhat decompositions in terms of masures.

\begin{Lemma}\label{l_frndly_prs}
Let $\cV_1,\cV_2$ be two filters on $\A$. We assume that $G=G_{\cV_1}\cdot N\cdot G_{\cV_2}$. Then for all $\cV_1'\in G.\cV_1$, $\cV_2'\in G.\cV_2$, there exists an apartment containing $\cV_1'\cup \cV_2'$.
\end{Lemma}

\begin{proof}
We follow \cite[Théorème 4.7.18]{bruhat1972groupes}.  Let $\cV_1'\in G.\cV_1$ and $\cV_2'\in G.\cV_2$. As $G$ acts transitively on the set of apartments, we can assume $\cV_1'=\cV_1$. Let $g\in G$ be such that $\cV_2:=g^{-1}.\cV_2'\subset \A$. Write $g=b_1nb_2$, with $b_1\in G_{\cV_1}$, $b_2\in  G_{\cV_2}$ and $n\in N$.  Let $A=b_1.\A$. Then $A\supset b_1.\cV_1=\cV_1$ and $\cV_2'=gg^{-1}.\cV_2 '=b_1nb_2.\cV_2=b_1n.\cV_2\subset b_1N.\cV_2\subset b_1.\A$, which proves the result.
\end{proof}

\begin{Remark}
\begin{enumerate}
\item Note that the converse of this lemma is also true (see Proposition~\ref{p_frnd_pr_T} for a precise statement). However, to prove it, we need to know that $\I$ satisfies~\ref{a_ma2}, which is not proven yet. Once we know it, we can proceed exactly as in Proposition~\ref{p_frnd_pr_T}, using~\ref{a_ma2}.  

\item Note that two points are not always contained in a common apartment, by \cite[Remark 6.10]{gaussent2008kac}. This corresponds to the fact that the Cartan decomposition does not hold for $G$ (we have $G_0\cdot N\cdot G_0\neq G$ in general).
\end{enumerate}

\end{Remark}

For $\epsilon\in \{-,+\}$, we set  $B^\epsilon=\fB^{\epsilon}(\cF)$\index[notation]{b@$B^\epsilon$}, where $\fB^{\epsilon}=\fT\cdot \fU^{\epsilon}$\index[notation]{b@$\fB^{\epsilon}$}. 

\begin{Proposition}\label{p_wMA3}
 Let $\epsilon\in \{-,+\}$. Then the fixator of $\fQ_{\epsilon\infty}$ is $T_0\cdot U^\epsilon$ and the stabilizer of $\fQ_{\epsilon\infty}$ is $B^\epsilon$.
\end{Proposition}

\begin{proof}
 We assume $\epsilon=+$ by symmetry. Let $g\in G_{\fQ_{+\infty}}$. Let $\Omega\in \fQ_{+\infty}$ be such that $g\in G_{\Omega}$. Up to reducing $\Omega$, we can assume that $\Omega=x+\overline{C^v_f}$, for some $x\in \A$. Then by Lemma~\ref{l_sctr_gd_fx}, we have $g\in U_{\Omega,\infty}^+\cdot U_{\Omega,\infty}^-\cdot N_{\Omega}=U_{\Omega,\infty}^+\cdot T_0$. Therefore $G_{\fQ_{+\infty}}\subset U^+\cdot T_0$. As the converse inclusion is clear, we have $G_{\fQ_{+\infty}}=T_0\cdot U^+$.

 Let $g\in G(\fQ_{+\infty}\curvearrowright \A)$. Let $\Omega\in \fQ_{+\infty}$ be such that $g\in G(\Omega\curvearrowright \A)$. Reducing $\Omega$ if necessary, we can assume that $\Omega=x+\overline{C^v_f}$, for some $x\in \A$. By Lemma~\ref{l_sctr_gd_fx}, we have $g\in N \cdot G_{\Omega}\subset N\cdot G_{\fQ_{+\infty}}$. Therefore $G(\fQ_{+\infty}\curvearrowright \A)=N\cdot G_{\fQ_{+\infty}}.$

  Let $g\in \mathrm{Stab}_G(\fQ_{+\infty})\subset G(\fQ_{+\infty}\cur \A)$. Write $g=nu$, with $n\in N$ and $u\in U^+$. Then $n.\fQ_{+\infty}=\fQ_{+\infty}$ and thus the image $\nu^v(n)$ of $n$ in $W^v$ is $1$. Therefore $n\in T$ and $g\in B^+$. Conversely, $B^+\subset \mathrm{Stab}_{G}(\fQ_{+\infty})$, which proves the proposition.
\end{proof}

\begin{Corollary}\label{c_wMAIII}
The masure $\I$ satisfies \ref{a_wma3}:  let $\fQ,\fQ'$ be two sector-germs at infinity. Then there exists an apartment containing $\fQ\cup \fQ'$.
\end{Corollary}

\begin{proof}
If $\fQ$, $\fQ'$ have the same sign $\epsilon$, we have $\fQ,\fQ'\in G.\fQ_{\epsilon\infty}$. By the Bruhat decomposition (\eqref{e_Bruhat}) and by  Proposition~\ref{p_wMA3}, $G=B^\epsilon\cdot  N\cdot  B^{\epsilon}=U^{\epsilon}\cdot  T\cdot N\cdot T\cdot  U^{\epsilon}= G_{\fQ_{\epsilon\infty}}\cdot  N \cdot  G_{\fQ_{\epsilon\infty}}$ and by  Lemma~\ref{l_frndly_prs}, we deduce the existence of an apartment containing $\fQ\cup \fQ'$.

 Now if $\fQ$ and $\fQ'$ have opposite signs, we can assume that $\fQ$ is positive and $\fQ'$ is negative. Then by the Birkhoff decomposition (see \eqref{e_Birkhoff}) and by Proposition~\ref{p_wMA3}, we have $G=B^+ \cdot  N \cdot  B^-=G_{\fQ_{+\infty}}\cdot  N\cdot  G_{\fQ_{-\infty}}$  and by Lemma~\ref{l_frndly_prs}, there exists an apartment containing $\fQ\cup \fQ'$.
\end{proof}

\begin{Remark}\label{r_iwa}
Using the Iwasawa decomposition (Proposition~\ref{p_Iwa}) and Lemma~\ref{l_frndly_prs}, we can prove that if $\fQ$ is a sector-germ at infinity and $\cV=g.\cV'$, where $g\in G$ and $\cV'$ is an inseparable filter on $\A$, then there exists an apartment containing $\cV$ and $\fQ$. This is actually a particular case of Proposition~\ref{p_mas_th_iwa}, once we know that $\I$ is a weak masure.
\end{Remark}

\section{Exchange condition and sundial configuration}

In this section, we prove that $\I$ satisfies the exchange condition (axiom \ref{a_ec}) and the sundial configuration (axiom \ref{a_sc}). 

For $\alpha\in \Phi$ and $k\in \R$, we set $D(\alpha,k)=\{x\in \A\mid \alpha(x)+k\geq 0\}$. Recall that a half-apartment of $\A$ is a set of the form $D(\alpha,k)$, where $\alpha\in \Phi$ and $k\in \Lambda$.

\begin{Lemma}\label{l_int_hlf_apt_gp}
Let $D$ be a half-apartment of $\A$. Let $B$ be an apartment containing $D$. Then either $B=\A$ or there exists $a\in \cF$ such that $B=x_{\alpha}(a).\A$. We then have $B\cap \A=D(\alpha,\omega(a))$.
\end{Lemma}

\begin{proof}
By symmetry, we can assume that $D$ contains $\fQ_{+\infty}$. By Lemma~\ref{l_iso_sect}, there exists $g\in G$ such that $g.\A=B$ and $g$ fixes $\A\cap B$. As $g$ fixes $\fQ_{+\infty}$, we have $g\in U^+\cdot T_0$, by Proposition~\ref{p_wMA3} and hence by Corollary~\ref{c_Fix_A}, we can assume $g\in U^+$. 

Let $x\in D$. We have $g\in G_x\cap U^+$ and hence by Lemma~\ref{l_u_fix}, we have $g\in U_{x,\infty}^+$. By Corollary~\ref{c_U_ma_inter}, we deduce $g\in U_{D,\infty}^+$. Write $D=D(\alpha,k)$, with $\alpha\in \Phi$ and $k\in \Lambda$. For $\beta\in \Delta_+\setminus \{\alpha\}$, we have $\beta(D)=\R$. Thus by Proposition~\ref{p_Rou4.5}, we have $g\in U_\alpha$. Write $g=x_\alpha(a)$, with $a\in \cF$. If $a=0$, then $B=\A$. If $a\neq 0$, then \[\A\cap B=\{x\in \A\mid g.x= x\}=\{x\in \A\mid g\in U_{x,\infty}^+\}=\{x\in \A\mid \alpha(x)+\omega(a)\geq 0\}=D(\alpha,\omega(a)),\] by Lemma~\ref{l_u_fix} and Proposition~\ref{p_Rou4.5}. \end{proof}

\begin{Proposition}\label{p_EC}(Exchange condition)
The set $\I$ satisfies the exchange condition, i.e:  for all appartments $A,B$ of $\I$ such that $A\cap B$ is  a half-apartment, $((A\cup B)\setminus (A\cap B))\cup M$  is  an apartment, where $M$ is the wall of $A\cap B$. 

\end{Proposition}

\begin{proof}
Translating by an element of $G$ if necessary, we may assume that $A=\A$. By Lemma~\ref{l_int_hlf_apt_gp}, we have $B=x_\alpha(a).\A$, for some $a\in \cF$ and $\alpha\in \Phi$. Let $\epsilon\in \{-,+\}$ be the sign of $\alpha$. Let $u_\epsilon=x_\alpha(a)$ and $u_{-\epsilon}=x_{-\alpha}(a^{-1})$. Let $r\in W^a$ be the reflection with respect to $\alpha^{-1}(\{-\omega(a)\})$. Let $n=m(x_\alpha(a)):=u_{-\epsilon} u_{\epsilon} u_{-\epsilon}$. Then by \eqref{e_m}, \eqref{e_m_s} and Lemma~\ref{l_nu_m}, $n\in N$  and $\nu(n)=r$. Let $D=\A\cap B$ and $D'=r.D$. Then we have $u_{-\epsilon}.D'=D'$ and thus $u_\epsilon u_{-\epsilon}.D'=u_{\epsilon}.D'=u_{-\epsilon}^{-1}r.D'=u_{-\epsilon}^{-1}.D$.  Therefore: \[((\A\cup B)\setminus (\A\cap B))\cup M=u_\epsilon.D'\cup D'=u_{-\epsilon}^{-1}.(D\cup D')=x_{-\alpha}(-a^{-1}).(D\cup D')=x_{-\alpha}(-a^{-1}).\A\] is an apartment of $\I$.  \end{proof}

Let $\fQ$, $\fQ'$ be two sector-germs at infinity. We say that $\fQ$ and $\fQ'$ are adjacent if for some apartment $A$ containing $\fQ\cup \fQ'$, $\fQ$ and $\fQ'$ are adjacent in $A$. By Lemma~\ref{l_iso_sect}, this does not depend on the choice of $A$. 

\begin{Lemma}\label{l_SC}
Let $\fQ$ be a sector-germ adjacent to  $\fQ_{+\infty}$ and different from $\fQ_{+\infty}$. Then there exists $i\in I$ such that $\fQ\in U_{\alpha_i} .(r_i.\fQ_{+\infty})$. 
\end{Lemma}

\begin{proof}
Let  $A$ be an apartment containing $\fQ_{+\infty}$ and $\fQ$. By Lemma~\ref{l_iso_sect}, there exists $u\in G$ such that $A=u.\A$ and $u$ fixes $\Omega:=A\cap \A$.  By Proposition~\ref{p_wMA3}, we have $u\in U^+\cdot T_0$ and by Corollary~\ref{c_Fix_A}, we can assume $u\in U^+$. Then $u^{-1}.\fQ$ is a sector-germ adjacent to $\fQ_{+\infty}$ and thus there exists $i\in I$ such that $u^{-1}.\fQ=r_i.(\fQ_{+\infty})$. By Lemma~\ref{l_u_fix} and Corollary~\ref{c_Uma_cl}, $u\in \bigcap_{x\in \Omega}U_{x,\infty}^+=U_{\Omega,\infty}^+$. 

 By Proposition~\ref{p_Mar8.58} 3) and Proposition~\ref{p_Rou4.5}, we have $U_{\Omega,\infty}^+=U_{\alpha_i,\Omega}\ltimes U_{\Omega,\infty}(\Delta_+\setminus \{\alpha_i\})$. We write $u=u_i v$, where $u_i\in U_{\alpha_i}$ and $v\in U_{\Omega,\infty}^+(\Delta_+\setminus \{\alpha_i\})$. Using Proposition~\ref{p_Rou4.5}, we write $v=\prod_{\beta\in \Delta_+\setminus\{\alpha_i\}} X_\beta(\underline{u_\beta})$, where 
 $(\underline{u_\beta})\in \prod_{\beta\in \Delta_+\setminus \{\alpha_i\}} \ffg_{\beta,\Z}\otimes \cF_{\geq f_{\Omega}}(\beta)$.  Let $x\in \Omega$. Then by \eqref{e_Kmr_1.3.14}, we have $r_i.(\Delta_+\setminus \{\alpha_i\})=\Delta_+\setminus \alpha_i$. Therefore if $\beta\in \Delta_+\setminus \{\alpha_i\}$, we have $\beta(x+r_i.C^v_f)=[\beta(x),+\infty[$. Consequently, $\omega(\underline{u_\beta})+\beta(x+r_i.C^v_f)\geq 0$ and hence $v$ fixes $x+r_i.C^v_f$. Therefore $u.r_i.\fQ_{+\infty}=u_i.r_i.\fQ_{+\infty}=\fQ$, which proves the lemma.
\end{proof}

\begin{Proposition}\label{p_SC_gp}(\textbf{sundial configuration})
Let $\fQ_1,\fQ_2$ be two adjacent sector-germs at infinity. Let $A$ be an apartment containing $\fQ_1$. Then there exists two opposite half-apartments $D_1,D_2$ of $A$ such that for both $i\in\{1,2\}$, there exists an apartment $B_i$ containing $D_i\cup \fQ_2$. 
\end{Proposition}

\begin{proof}
By symmetry, we may assume that $\fQ_1$ and $\fQ_2$ are positive. Considering $g.\fQ_1$ for some $g\in G$, we may assume $\fQ_1=\fQ_{+\infty}$. Then by Lemma~\ref{l_SC}, we have $\fQ_2=u.(r_i.\fQ_{+\infty})$ for some $i\in I$ and $u\in U_{\alpha_i}$. Let $B_1=u.\A$. Then by Lemma~\ref{l_int_hlf_apt_gp}, $D_1:=\A\cap B_1$ is a half-apartment of $\A$. By Proposition~\ref{p_EC}, $B_2:=(B_1\setminus \A)\cup (\A\setminus B_1)\cup M$ is an apartment, where $M$ is the wall of $D_1$. Then $(B_1,B_2,\A)$ satisfies the conditions of the proposition. 
\end{proof}

\begin{Corollary}\label{c_w_mas_pf}
Assume that  $(\alpha_i^\vee)_{i\in I}$ is positively free. The masure $\I$ is a weak masure in the sense of Definition~\ref{d_w_mas}. It satisfies~\ref{a_ma2} and every filter on $\A$ satisfies \ref{a_TF}. Moreover every apartment isomorphism (resp. Weyl-isomorphism) provided by \ref{a_ma2} is induced by an element of $G$ (resp. $G^u$). 
\end{Corollary}

\begin{proof}
We proved that $\I$ satisfies~\ref{a_oco},  \ref{a_wma3} and \ref{a_sc} (see  Propositions~\ref{p_SOMA}, \ref{p_wMA3} and \ref{p_SC_gp}). Therefore by Theorem~\ref{t_MA2}, $\I$ satisfies~\ref{a_ma2}.  Using Lemma~\ref{l_TF}, we deduce that every filter on $\A$ satisfies \ref{a_TF}. The ``moreover part'' follows from Definition/Propositions~\ref{dp_MA1_tildeA} and \ref{dp_MA1_A}.
\end{proof}

\section{Description of parahoric subgroups using the minimal group only}\label{ss_para_minGp}

We used Mathieu's completions $G^{ma+}$ and $G^{ma-}$ to describe the fixators $G_x$, for $x\in \A$ and we have description  of $G_{\cV}$ for filters $\cV$ with good fixators involving $U_{\cV,\infty}^+$ and $U_{\cV,\infty}^-$. 
In this subsection, we describe $G_F$ entirely in terms of the minimal group, when $F$ is an alcove or a point of $\A$.

We follow \cite[2.4]{bardy2025twin}. We start with the case of an alcove and then deduce the case of a point.

\begin{Lemma}\label{l_fx_nar_fltr}
Let $\cV$ be an inseparable filter on $\A$. Then we have: \[G_{\cV}=U_{\cV,\infty}^+\cdot U_{\cV,\fin}^-\cdot N_\cV=U_{\cV,\infty}^-\cdot U_{\cV,\fin}^+\cdot N_{\cV}.\]
\end{Lemma}

\begin{proof}
If $\epsilon\in \{-,+\}$, we have  $U_{\cV,\fin}^\epsilon \subset U_{\cV,\infty}^\epsilon\subset G_{\cV}$ (by Proposition~\ref{p_GR08_3.4}) and $N_{\cV}\subset G_{\cV}$ and thus the sense $\supset$ is clear. Let $\epsilon\in \{-,+\}$.  Let $g\in G_{\cV}$.  Then by Proposition~\ref{p_Iwa}, we can write $g=u_\epsilon nu$, with $u_\epsilon\in U^\epsilon$, $n\in N$ and $u\in U_{\cV,\fin}$. We have $U_{\cV,\fin}=\langle U_{\alpha,\cV}\mid \alpha\in \Phi\rangle\subset G_{\cV}$.

 Let $\Omega\in \cV$ be fixed by both $g$ and $u$. Then for $x\in \Omega$, we have $g.x=x=u_\epsilon n.x$ and hence $\rho_{\epsilon\infty}(u_\epsilon n.x)=n.x=x$, since $\rho_{\epsilon\infty}(u_\epsilon n.x)$ is the unique element of $U^\epsilon.u_\epsilon n.x\cap \A$. Therefore $n\in N_\Omega$ and thus $u_\epsilon=gu^{-1}n^{-1}\in G_{\Omega}$. By \eqref{e_cnj_U_fin}, we can write $nu=u'n$, with $u'\in U_{\cV,\fin}$. Then $u_\epsilon u'\in U_{\cV,\infty\epsilon}$ and thus by Proposition~\ref{p_GR08_3.4}, $g\in U_{\cV,\infty}^\epsilon\cdot U_{\cV,\fin}^{-\epsilon}\cdot N_{\cV}$. This proves the sense $\subset$ and hence the lemma.
\end{proof}

 \begin{Proposition}(\cite[Proposition 2.3]{bardy2025twin})\label{p_bphr_2.4}
 Let $\cV$ be a filter on $\A$.  Assume that $\cV$ is inseparable and that $N_{\cV}=T_0$ (for example $\cV$ is an alcove). Then:\begin{enumerate}
 \item $G_{\cV}=U_{\cV,\fin}^ {+}\cdot U_{\cV,\fin}^ {-}\cdot T_{0}=U_{\cV,\fin}^ {+}\cdot T_{0}\cdot U_{\cV,\fin}^ {-}=U_{\cV,\fin}\cdot T_{0}=\langle T_{0}, (U_{\alpha,\cV})_{\alpha\in\Phi} \rangle$,
 
 \item We have $U_{\cV,\infty}^+=U_{\cV,\fin}^+=G_{\cV}\cap U^+$ and $U_{\cV,\infty}^-=U_{\cV,\fin}^-=G_{\cV}\cap U^-$. 
 \end{enumerate}

 \end{Proposition}

 \begin{proof}

   By  Lemma~\ref{l_fx_nar_fltr} and the fact that $T_{0}$ normalizes $U_{\cV,\fin}^ {\pm},U_{\cV,\infty}^+$ and $U_{\cV,\infty}^-$ (by \ref{a_KMT4'}, Proposition~\ref{p_cnj_N}) and \eqref{e_cnj_U_fin}), we have $G_{\cV}=U_{\cV,\infty}^+\cdot T_0\cdot  U_{\cV,\fin}^-=U_{\cV,\fin}^+\cdot T_{0}\cdot U_{\cV,\infty}^-.$
 By  uniqueness in the Birkhoff decomposition (see \eqref{e_Birkhoff}), we deduce that $U_{\cV,\infty}^\epsilon=U_{\cV,\fin}^\epsilon=G_{\cV}\cap U^\epsilon$, for $\epsilon\in \{-,+\}$. 
 \end{proof}
 
  In particular the Iwahori group $K_I=G_{C_{0}^+}$ (fixator in $G$ of the fundamental alcove $C_{0}^+=\germ_{0}(C^v_{f})$) is $\langle T_{0}, (U_{\alpha,C_{0}^+})_{\alpha\in\Phi} \rangle$.
 This is the same definition as in \cite{braverman2016iwahori} (given there in the untwisted affine case).
 This result was also proved in \cite[7.2.2]{bardy2016iwahori}, using the results of \cite{braverman2016iwahori}.

\begin{Example}
In \cite[2.2.4 (3)]{bardy2025twin} it is claimed that if $\cV$ is an almost-open filter on $\A$, then $N_{\cV}=T_0$. It is actually wrong, by the example below.

Consider the affine (non-cofree) Kac--Moody  datum $\cS=(\begin{psmallmatrix} 2 & -2\\ -2 & 2\end{psmallmatrix},X,Y,(\alpha_i)_{i\in \{0,1\}}, (\alpha_i^\vee)_{i\in \{0,1\}})$, where $X$, $Y$, $(\alpha_i)_{i\in \{0,1\}}$ and $(\alpha_i^\vee)_{i\in \{0,1\}}$ are defined as follows. Set $X=\Z \alpha_1\oplus \Z  \delta$ and $Y=\Z\alpha_1^\vee \oplus \Z d$, where $\alpha_1(\alpha_1^\vee)=2$, $\alpha_1(d)=0$, $\delta(\alpha_1^\vee)=0$ and  $\delta(d)=1$. We set $\alpha_0=\delta-\alpha_1$ and $\alpha_0^\vee=-\alpha_1^\vee$. 

Set $\cV=\Omega=\{0,\alpha_1^\vee\}$. Let  $\alpha\in \Phi$. Write $\alpha=k\alpha_0+\ell\alpha_1$, with $k,\ell\in \Z$. Then $\alpha\notin \Z\delta$ and hence $k\neq \ell$. Therefore $\alpha(\alpha_1^\vee)=2(k-\ell)\neq 0$. Consequently, $\Omega$ is almost-open. However, we have $r_1r_0.\alpha_1^\vee=r_1(2\alpha_0^\vee+\alpha_1^\vee)=3\alpha_1^\vee+2\alpha_0^\vee=\alpha_1^\vee$. Therefore if $n\in N$ induces $r_1r_2$ on $\A$, $\langle n\rangle \subset  N_\Omega\not\subset T_0$. 
\end{Example}

\begin{Lemma}\label{l_U_F_fin}
Let $x\in \A$ and $\epsilon\in \{-,+\}$. Let $F$ be a local face of $\A$ based at $x$. Assume that for some alcove of $\A$ based at $x$, we have $F\subset \overline{C}$ and  $U_{F,\infty}^\epsilon=U_{C,\infty}^\epsilon$. Then $U_{F,\fin}^\epsilon=U_{F,\infty}^\epsilon$.
\end{Lemma}

\begin{proof}
By Proposition~\ref{p_GR08_3.4}, we have $U_{F,\fin}^\epsilon\subset U_{F,\infty}^\epsilon$. Conversely, take $u\in U_{F,\infty}^\epsilon$. Then $u\in U_{C,\infty}^\epsilon=U_{C,\fin}^\epsilon$, by Proposition~\ref{p_bphr_2.4}. Therefore we can write $u=u_1\ldots u_k$, with $k\in \N$ and for all $i\in \llbracket 1,k\rrbracket$, $u_i\in U_{\beta_i,C}$, for some $\beta_i\in \Phi$. Then as $F\subset \overline{C}$, we have  $u_i\in U_{\beta_i,F}$ for all $i\in \llbracket 1,k\rrbracket$ and thus $u\in U_{F,\fin}^\epsilon$. 
\end{proof}

\begin{Corollary}\label{c_BPHR_2.4.1}\cite[2.4.1 (2)]{bardy2025twin}
Let $x\in \A$. Then $U_{x,\infty}^+=U_{x,\fin}^+$ and $U_{x,\infty}^-=U_{x,\fin}^-$. 
\end{Corollary}

\begin{proof}
We have $U_{\germ_x(x+C^v_f),\infty}^+=U_{x,\infty}^+$ and $U_{\germ_x(x-C^v_f),\infty}^-=U_{x,\infty}^-$, by Proposition~\ref{p_Rou4.5}. Thus by Lemma~\ref{l_U_F_fin}, we get the result in this case.
\end{proof}

\section{Extensions of masures and axiom~\ref{a_ma2} in the  non-cofree case}\label{s_ext_mas}

In order to prove that the masure $\I$ satisfies~\ref{a_ma2}, we want to use Theorem~\ref{t_MA2}. However, an assumption of this theorem is that the family $(\alpha_i^\vee)_{i\in I_A}$ is (positively)-free. This assumption is a bit restrictive since for example, $\mathring{\fG}(\cF[t,t^{-1}])\rtimes \cF^\times$ has a non-free set of coroots, if $\mathring{\fG}$ is a split reductive group. To avoid this assumption, we prove that if $G$ has a non-free set of coroots, then $\I$ is a quotient of the masure $\widetilde{\I}$ of $\tilde{G}$, where $\tilde{G}$ is an extension of $G$ with free simple coroots. Then we deduce~\ref{a_ma2} for $\I$ from~\ref{a_ma2} for $\tilde{\I}$.

Recall that $\cS=(A,X,Y,(\alpha_i)_{i\in I},(\alpha_i^\vee)_{i\in I})$ be a free Kac--Moody datum. Let $\tilde{\cS}=(A,\tilde{X},\tilde{Y},(\tilde{\alpha}_{i})_{i\in I},(\tilde{\alpha_i^\vee})_{i\in I})$ be an other free Kac--Moody datum. We assume that it satisfies the following properties. 

\begin{Assumption}\label{as_ext}
\begin{enumerate}

\item $\tilde{X}=X\oplus \overline{X}$ and $\tilde{Y}=Y\oplus \overline{Y}$, for some lattices $\overline{X}$ and $\overline{Y}$. 

\item $\tilde{\alpha}_i|_{X}=\alpha_i$, $\tilde{\alpha}_i(\overline{X})=\{0\}$, and $\tilde{\alpha}_i^\vee\in \alpha_i^\vee+\overline{Y}$, for all $i\in I$,

\item $\fG_{\tilde{\cS}}$ is a central extension of $\fG_{\cS}$ and the Kernel of the projection map $\fG_{\tilde{\cS}}\twoheadrightarrow \fG_{\cS}$ is $\overline{\fT}:\cR\mapsto \Hom(\overline{X},\cR^\times)$. 
\end{enumerate}
\end{Assumption}

There always exists a cofree Kac--Moody datum satisfying these properties, for example 
the cofree extension of $\cS$ defined in \cite[7.3.2]{marquis2018introduction}, by  \cite[Proposition 7.64]{marquis2018introduction} (see also the proof of Proposition~\ref{p_rt_sp_dec}, where this extension is used).

 Let $\tilde{\A}=\tilde{Y}\otimes \R=\A\oplus \overline{\A}$, where $\overline{\A}=\overline{Y}\otimes \R$.  By assumption, we have $\tilde{\alpha_i^\vee}\in \alpha_i^\vee+\overline{\A}, \forall i\in I.$ Let $\tilde{\Phi}$ (resp. $\Phi$) be the root system of $\tilde{\cS}$ (resp. $\cS$).   We have: \begin{equation}\label{e_tilde_alp_A_bar}
\tilde{\alpha}(\overline{\A})=\{0\}, \forall \tilde{\alpha}\in \tilde{\Phi}. 
\end{equation}

Recall that $\cF$ is a field equipped with a non-trivial valuation $\omega:\cF^\times \rightarrow \R$. Let  $\tilde{G}:=\fG_{\tilde{\cS}}(\cF)$. Let $\tilde{\I}=\I(\fG_{\tilde{\cS}},\cF,\omega)=\tilde{G}\times \tilde{\A}/\sim$ and  $\I=\I(\fG_{\cS},\cF,\omega)=G\times \A/\sim$ be the masures of $\tilde{G}$ and $G$ defined in Subsection~\ref{ss_d_mas}. We denote with a tilde the objects corresponding to $\tilde{\fG}$ (for example $\tilde{U}^+$ is the positive standard unipotent subgroup of $\tilde{G}$).

 Recall that $\tilde{T}=\Hom_{\mathrm{Gr}}(\tilde{X},\cF^\times)$, $\overline{T}=\Hom_{\mathrm{Gr}}(\overline{X},\cF^\times)$ and $T=\Hom_{\mathrm{Gr}}(X,\cF^\times)$. We regard $T$ and $\overline{T}$ as subgroups of $\tilde{T}$ by setting $t(\overline{\chi})=1=\overline{t}(\chi)$, for $\chi\in X$, $\overline{\chi}\in \overline{X}$, $t\in T$ and  $\overline{t}\in \overline{T}$.

\begin{Lemma}\label{l_pi_tilde_phi}
The map $\pi:\tilde{\Phi}\rightarrow \Phi$ defined by $\pi(\sum_{j\in I} n_j\tilde{\alpha}_j)=\sum_{j\in I} n_j \alpha_j$, for $\sum_{j\in I} n_i\tilde{\alpha}_i\in \tilde{\Phi}$ is well-defined and is a bijection.
\end{Lemma}

\begin{proof}
For $i\in I$, denote by $\tilde{r}_i$ the simple reflection of $\tilde{\A}$ associated with $i$. Let $k\in \Z_{\geq 2}$ and  $i_1,\ldots,i_k\in I$. Write $\tilde{r}_{i_{k-1}}\ldots \tilde{r}_{i_2}.\tilde{\alpha}_{i_1}=\sum_{j\in I} n_j\tilde{\alpha}_j$, where $(n_j)\in \Z^I$.  Then $\tilde{r}_{i_k}\ldots \tilde{r}_{i_2}.\tilde{\alpha}_{i_1}=\sum_{j\in I} n_j \tilde{r}_{i_k}.\tilde{\alpha}_j=\sum_{j\in I} n_j(\tilde{\alpha}_j-a_{j,i_k}\tilde{\alpha}_{i_k})$. Assume that $r_{i_{k-1}}\ldots r_{i_2}.\alpha_{i_1}=\sum_{j\in I} n_j\alpha_j$. Then $r_{i_k}\ldots r_{i_2}.\alpha_{i_1}=r_{i_k}. \sum_{j\in I} n_j\alpha_j=\sum_{j\in I} n_j(\alpha_j-a_{j,i_k}\alpha_{i_k})$. By induction on $k$, we deduce that for all $k\in \N$ and $i_1,\ldots, i_k\in I$, if $\tilde{\alpha}=\tilde{r}_{i_k}\ldots \tilde{r}_{i_2}.\tilde{\alpha}_{i_1}$ and $\alpha=r_{i_k}\ldots r_{i_2}.\alpha_{i_1}$, then  the coefficients of $\alpha$ and $\tilde{\alpha}$ are the same, in the free families $(\tilde{\alpha}_j)_{j\in J}$ and $(\alpha_j)_{j\in J}$ respectively. Lemma follows.
\end{proof}

From now on, if $\tilde{\alpha}\in \tilde{\Phi}$, we write $\alpha$ instead of $\pi(\tilde{\alpha})$. Let $\pi:\tilde{\A}\twoheadrightarrow \A=\tilde{\A}/\overline{Y}\otimes \R$ and $\pi:\tilde{G}\twoheadrightarrow G=\tilde{G}/\overline{T}$ be the canonical projections.

\begin{Lemma}\label{l_pi_tilde_G_x}

\begin{enumerate}
\item Let $\tilde{\cV}$ be a filter on $\tilde{\A}$ and $\epsilon\in \{-,+\}$. Then $\pi(\tilde{U}_{\tilde{\cV},\fin}^\epsilon)=U^\epsilon_{\pi(\tilde{\cV}),\fin}$ and $\pi(\tilde{U}_{\tilde{\cV},\infty}^\epsilon)=U^\epsilon_{\pi(\tilde{\cV}),\infty}$.

\item Let $\tilde{x}\in \tilde{\A}$. Then $\pi(\tilde{G}_{\tilde{x}})\subset G_{\pi(\tilde{x})}$. 
\end{enumerate}

\end{Lemma}

\begin{proof}
(1) By \cite[1.10]{rousseau2016groupes},  $\pi$ restricts to an  isomorphism $\tilde{U}^\epsilon\rightarrow U^\epsilon$. Let $\tilde{\Omega}\in \tilde{\cV}$ and $\Omega=\pi(\tilde{\Omega})$.  Let   $\tilde{\alpha}\in \tilde{\Phi}$. We have $\tilde{\alpha}=\alpha\circ\pi$ and thus $f_{\tilde{\Omega}}(\tilde{\alpha})=f_{\Omega}(\alpha)$. Consequently, $\pi(\tilde{U}_{\tilde{\Omega},\fin}^{\epsilon})=U_{\Omega,\fin}^{\epsilon}$ and thus $\pi(\tilde{U}_{\tilde{\cV},\fin}^{\epsilon})=\bigcup_{\tilde{\Omega}\in \tilde{\cV}} U_{\pi(\tilde{\Omega}),\fin}^\epsilon=U_{\pi(\cV),\fin}^{\epsilon}$. 

Let $\tilde{x}\in \tilde{\A}$ and $x=\pi(\tilde{x})$. Then by Corollary~\ref{c_BPHR_2.4.1}, $\pi(\tilde{U}_{\tilde{x},\infty}^\epsilon)=\pi(\tilde{U}_{\tilde{x},\fin}^\epsilon)=U_{x,\fin}^\epsilon=U_{x,\infty}^{\epsilon}$.

Let $\tilde{\Omega}\in \tilde{\cV}$. Then by Corollary~\ref{c_U_ma_inter} and since $\pi$ restricts to an isomorphism from $\tilde{U}^\epsilon$ to  $U^\epsilon$, we have $\pi(U^{\epsilon}_{\tilde{\Omega},\infty})=\pi(\bigcap_{\tilde{x}\in \tilde{\Omega}} \tilde{U}_{\tilde{x},\infty}^{\epsilon})=\bigcap_{\tilde{x}\in \tilde{\Omega}}\pi(\tilde{U}_{\tilde{x},\infty}^\epsilon)=\bigcap_{x\in \pi(\tilde{\Omega})}U_{x,\infty}^{\epsilon}=U_{\pi(\tilde{\Omega}),\infty}^{\epsilon}$. Now $\pi(\tilde{U}_{\tilde{\cV},\infty}^{\epsilon})=\pi(\bigcup_{\tilde{\Omega}\in \tilde{\cV}} \tilde{U}_{\tilde{\cV},\infty})=\bigcup_{\tilde{\Omega}\in \tilde{\cV}} U_{\pi(\tilde{\Omega}),\infty}^{\epsilon}=U_{\pi(\tilde{\cV}),\infty}^{\epsilon}$, which proves (1).

(2)  Let $\fN=\fN_{\cS}$ and $\tilde{\fN}=\fN_{\tilde{\cS}}$ be as defined in Subsection~\ref{ss_Tits_minKM}.   Let $N=\fN(\cF)$ and $\tilde{N}=\tilde{\fN}(\cF)$. By \cite[1.10]{rousseau2016groupes}, we have $\pi(\tilde{T})=T$, $\pi(\tilde{N})=N$ and $\tilde{N}=\pi^{-1}(N)$.

 Denote by $\tilde{\nu}:\tilde{N}\rightarrow \mathrm{Aut}(\tilde{\A})$ (resp. $\nu:N\rightarrow \mathrm{Aut}(\A)$) the actions defined in Subsection~\ref{ss_Action_N_A}. Let $\overline{T}=\mathrm{Hom}(\overline{X},\cF^\times)$. We have $\tilde{T}=T\overline{T}$ and $\overline{T}$ is contained in the center of $\tilde{G}$. Let  $\tilde{\chi}\in \tilde{X}$. Write $\tilde{\chi}=\chi \overline{\chi}$, with $\chi\in X$ and $\overline{\chi} \in \overline{X}$. Then $\tilde{\chi}(t\overline{t})=\chi(t)\overline{\chi}(\overline{t})$, for $(t,\overline{t})\in T\times \overline{T}$. Hence by definition of $\tilde{\nu}$ and $\nu$, we have $\pi(\tilde{\nu}(\tilde{t}))=\nu(\pi(\tilde{t}))$, for all $\tilde{t}\in \tilde{T}$. Let $i\in I$.  We have    $\pi(\tilde{s}_{\tilde{\alpha}_i})=\tilde{s}_{\alpha_i}$ (we write $\tilde{s}_{\alpha_i}$ and $\tilde{s}_{\tilde{\alpha_i}}$ instead of $\tilde{s}_i$ used so far to distinguish $\fN$ and $\tilde{\fN}$). Therefore $\pi\circ\tilde{\nu}(\tilde{s}_{\tilde{\alpha}_i})=r_i=\nu(\pi(\tilde{s}_{\tilde{\alpha}_i})))$. 
 
Let $\tilde{t}\in \tilde{T}$ and $\chi\in X\subset \tilde{X}$. Write $\tilde{t}=t\overline{t}$, with $t\in T$ and $\overline{t}\in \overline{T}$. Then $\chi(\tilde{t})=\tilde{t}(\chi)=t\overline{t}(\chi)=t(\chi)=\chi(t)=\chi(\pi(\tilde{t}))$.  Thus by \eqref{e_act_T},  we have $\pi(\tilde{\nu}(\tilde{t}))=\nu(\pi(\tilde{t}))$.  Using \eqref{e_dec_N}, we deduce:\begin{equation}\label{e_action_tilde_N}
  \pi\circ \tilde{\nu}(\tilde{n})=\nu(\pi(\tilde{n})),\forall \tilde{n}\in \tilde{N}.
  \end{equation}

Let $\epsilon\in \{-,+\}$.  By \eqref{e_action_tilde_N}, we have $\pi(\tilde{N}_{\tilde{x}})\subset N_{\pi(\tilde{x})}$. By Theorem~\ref{t_G_x} and (1), we have $\tilde{G}_{\tilde{x}}= \tilde{N}_{\tilde{x}}\cdot \tilde{U}_{\tilde{x},\infty}^{-}\cdot \tilde{U}_{\tilde{x},\infty}^+$. Therefore $\pi(\tilde{G}_{\tilde{x}})\subset G_{\pi(\tilde{x})}.$
\end{proof}

\begin{Proposition}\label{p_pr_ext_mas}
 Let $\pi:\tilde{\A}\twoheadrightarrow \A=\tilde{\A}/\overline{Y}\otimes \R$ and $\pi:\tilde{G}\twoheadrightarrow G=\tilde{G}/\overline{T}$ be the canonical projections. Then $\pi$ extends to a map $\pi:\tilde{\I}\rightarrow \I$ defined by $\pi(\tilde{g}.\tilde{x})=\pi(\tilde{g}).\pi(\tilde{x})$, for $\tilde{g}\in \tilde{G}$ and $\tilde{x}\in \tilde{\A}$. 

\end{Proposition}

\begin{proof}

  Let $\tilde{x},\tilde{y}\in \tilde{\A}$ and $\tilde{g}\in \tilde{G}$ be such that $\tilde{y}=\tilde{g}.\tilde{x}$. By \eqref{e_def_sim}, there exists $\tilde{n}\in \tilde{N}$ such that $\tilde{y}=\tilde{n}.\tilde{x}$. Then $\tilde{n}^{-1}\tilde{g}\in \tilde{G}_{\tilde{x}}$.  Therefore, by Lemma~\ref{l_pi_tilde_G_x}, we have $\pi(\tilde{n}^{-1}\tilde{g}).x=x$ and thus $\pi(\tilde{g}).\pi(\tilde{x})=\pi(\tilde{n}).\pi(\tilde{x})=\pi(\tilde{n}.\tilde{x})=\pi(\tilde{y})$, by \eqref{e_action_tilde_N}. Proposition follows. 
\end{proof}

\begin{Lemma}\label{l_lft_stb_apt}
Let $\tilde{z}\in \tilde{\A}$ and  $z=\pi(\tilde{z})\in \A$. Let $\tilde{g}\in \tilde{G}$  be such that $g.z=z$, where $g=\pi(\tilde{g})$. Then $\tilde{g}.\tilde{z}\in\tilde{z}+\overline{\A} \subset \tilde{\A}$. 
\end{Lemma}

\begin{proof}
By Theorem~\ref{t_G_x} and Corollary~\ref{c_BPHR_2.4.1}, we can write $g=u_+ u_- n$, with $u_+\in U_{z,\fin}^+$, $u_-\in U_{z,\fin}^-$ and $n\in N_{z}$. Let $\tilde{n}\in \pi^{-1}(\{n\})$ and let $\tilde{u}_+$ and $\tilde{u}_-$ be the elements of $\tilde{U}^+$ and $\tilde{U}^-$ corresponding to $u_+$ and $u_-$. Let $\tilde{h}=\tilde{u}_+\tilde{u}_- \tilde{n}$.  Then $\pi(\tilde{h})=\pi(\tilde{g})$ and thus $\tilde{h}\in \tilde{g}\overline{T}$. 

By \eqref{e_action_tilde_N}, we have $\tilde{n}.\tilde{z}\in \tilde{z}+\overline{\A}$. We have $\tilde{\alpha}(\tilde{z}+\overline{\A})=\alpha(z)$ for all $\tilde{\alpha}\in \tilde{\Phi}$ and thus $\tilde{u_+}$ and $\tilde{u_-}$ fix $\tilde{n}.\tilde{z}$. Therefore $\tilde{g}.\tilde{z}\in \tilde{h}.\tilde{z}+\overline{\A}=\tilde{z}+\overline{\A}$.

\end{proof}

\begin{Lemma}\label{l_pr_img_apt}
We have $\pi^{-1}(\A)=\tilde{\A}$. 
\end{Lemma}

\begin{proof}
Let $\tilde{y}\in \tilde{\I}$ be such that $y:=\pi(\tilde{y})\in \A$. Write $\tilde{y}=\tilde{g}.\tilde{z}$, where $\tilde{g}\in \tilde{G}$ and $\tilde{z}\in \tilde{\A}$. Let $\tilde{x}\in \tilde{\A}$ be such that $\pi(\tilde{x})=y$. Let $g=\pi(\tilde{g})$.  By \eqref{e_def_sim}, we can write $g.z=n.z$, where $n\in N$ and $z=\pi(\tilde{z})$. Let $\tilde{n}\in \pi^{-1}(\{n\})$. Then $\pi(\tilde{g}.\tilde{z})=\pi(\tilde{n}.\tilde{z})$ and thus by Lemma~\ref{l_lft_stb_apt}, we have $\tilde{y}=\tilde{g}.\tilde{z}\in \tilde{n}.(\tilde{z}+\overline{\A})\subset \tilde{\A}$. 
\end{proof}

\begin{Proposition}\label{p_ext_inter}
\begin{enumerate}
\item The projection $\pi$ induces a bijection $\tilde{A}\mapsto \pi(\tilde{A})$ from the set of apartments of $\tilde{\I}$ to the set of apartments of $\I$. Its reciprocal is the map $A\mapsto \pi^{-1}(A)$. 

\item If $A$ and $B$ are two apartments of $\I$, and $\tilde{A}=\pi^{-1}(A)$, $\tilde{B}=\pi^{-1}(B)$, then we have $A\cap B=\pi(\tilde{A}\cap \tilde{B})$.
\end{enumerate} 
\end{Proposition}

\begin{proof}
(1) Denote by $\tilde{\sA}$ (resp. $\sA$) the set of apartments of $\tilde{\I}$ (resp. $\I$). Let $\pi^{\sA}:\tilde{\sA}\rightarrow \sA$ be defined by $\pi^{\sA}(\tilde{A})=\pi(\tilde{A})$, for $\tilde{A}\in \tilde{\sA}$.  Let $\tilde{A}\in \tilde{\sA}$. Then $\tilde{A}=\tilde{g}.\tilde{\A}$, for some $\tilde{g}\in \tilde{G}$. By definition of $\pi$, we have  $\pi(\tilde{A})=\pi(\tilde{g}).\A$ and thus $\pi^{\sA}(\tilde{A})\in \sA$, which proves that $\pi^{\sA}$ is well-defined. 

If $A\in \sA$, then $A=g.\A$, for some $g\in G$. Let $\tilde{g}\in \pi^{-1}(\{g\})$. Then $\pi(\tilde{g}.\tilde{\A})=g.\A=A$, which proves that $\pi^{\sA}$ is surjective. 

Let $A\in \sA$. Write $A=g.\A$, with $g\in G$. Let $\tilde{g}\in \pi^{-1}(\{g\})$.  Then by Lemma~\ref{l_pr_img_apt},  $\pi^{-1}(A)=\tilde{g}.\tilde{\A}$. Now if $\tilde{h}\in \tilde{G}$ is such that $\pi(\tilde{h}.\tilde{\A})=g.\A$, then $\pi(\tilde{h}^{-1}\tilde{g}.\tilde{\A})=\A$ and thus by Lemma~\ref{l_pr_img_apt}, $\tilde{h}^{-1}\tilde{g}.\tilde{\A}\subset \tilde{\A}$. Therefore $\tilde{g}.\tilde{\A}\subset \tilde{h}.\tilde{\A}$ and hence $\tilde{g}.\tilde{\A}=\tilde{h}.\tilde{\A}$. Consequently $\pi^\sA$ is injective and $(\pi^{\sA})^{-1}(A)=\tilde{g}.\tilde{\A}$. This proves (1). 

(2) Translating by an element of $G$, we can assume that $A=\A$. We have $\pi(\tilde{\A}\cap \tilde{B})\subset \pi(\tilde{\A})\cap \pi(\tilde{B})=\A\cap B$. In particular, if $\tilde{\A}\cap \tilde{B}=\emptyset$, then $\A\cap B=\emptyset$.

Assume now that $\A\cap B\neq \emptyset$. Let $x\in \A\cap B$. Let $g\in G$ be such that $g.\A=B$ and $g$ fixes $x$. Let $\tilde{g}\in \pi^{-1}(\{g\})$. Then by (1), we have $\tilde{g}.\tilde{\A}=\tilde{B}$. By Lemma~\ref{l_lft_stb_apt}, if $\tilde{x}\in \pi^{-1}(\{x\})$, we have $\tilde{g}.\tilde{x}\in \tilde{x}+\overline{\A}\subset \tilde{\A}$ and thus $\tilde{g}.\tilde{x}\in \tilde{\A}\cap \tilde{B}$. We have $\pi(\tilde{g}.\tilde{x})=g.x=x$, which proves that $\A\cap B\subset \pi(\tilde{\A}\cap \tilde{B})$. Consequently, $\tilde{\A}\cap \tilde{B}$ is empty if and only if $\A\cap B$ is empty, and we have $\A\cap B=\pi(\tilde{\A}\cap \tilde{B})$. 
\end{proof}

\begin{Corollary}\label{c_w_mas_gen}
Let $\cS$ be a free Kac--Moody datum and $\cF$ be a field equipped with a non-trivial valuation $\omega:\cF^\times \rightarrow \R$. Then  the masure $\I$ of $(\fG_{\cS}(\cF),\omega)$ is a weak masure in the sense of Definition~\ref{d_w_mas}. It satisfies \ref{a_ma2} and every filter on $\A$ satisfies \ref{a_TF}. Moreover every apartment isomorphism (resp. Weyl-isomorphism) provided by \ref{a_ma2} is induced by an element of $G$ (resp. $G^u$). 
\end{Corollary}

\begin{proof}
As explained in the beginning of this section, there exists a cofree extension $\tilde{\cS}$ of $\cS$ satisfying Assumption~\ref{as_ext}. Let $A$ and $B$ be two apartments of $\I$. Up to translating by an element of $G$, we can assume that $A=\A$. By Corollary~\ref{c_w_mas_pf}, $\tilde{\I}$ satisfies~\ref{a_ma2}. Let $\tilde{B}=\pi^{-1}(B)$. Then $\tilde{\A}\cap \tilde{B}$ is enclosed: there exists $k\in \N$, $\tilde{\beta}_1,\ldots,\tilde{\beta}_k\in \tilde{\Phi}$ and $\lambda_1,\ldots,\lambda_k\in \Lambda$ such that $\tilde{\A}\cap \tilde{B}=\bigcap_{i=1}^k D(\tilde{\beta}_i,\lambda_i)$. By Proposition~\ref{p_ext_inter}, we have $\A\cap B=\pi(\tilde{\A}\cap \tilde{B})=\bigcap_{i=1}^k D(\beta_i,\lambda_i)$: $\A\cap B$ is enclosed. 

Let $\tilde{g}\in \tilde{G}$ be such that $\tilde{g}.\tilde{\A}=\tilde{B}$ and $\tilde{g}$ fixes $\tilde{\A}\cap \tilde{B}$. Let $g=\pi(\tilde{g})\in G$. Assume that $\A\cap B$ is non-empty and take $x\in \A\cap B$. Let $\tilde{x}\in \pi^{-1}(\{x\})$. Then $\tilde{x}\in \pi^{-1}(\A)\cap \pi^{-1}(B)=\tilde{\A}\cap \tilde{B}$, thus $\tilde{g}.\tilde{x}=\tilde{x}$ and hence $g.x=\pi(\tilde{g}.\tilde{x})=\pi(\tilde{x})=x$. Therefore $g$ fixes $\A\cap B$ and $\I$ satisfies~\ref{a_ma2}. 
\end{proof}

The following corollary extends \cite[Proposition 7.4.8]{bruhat1972groupes} to masures.

\begin{Corollary}\label{c_BT_7.4.8}
\begin{enumerate}
\item Let $g\in G$ (resp. $g\in G^u$). Then there exists $n\in N$ (resp. $n\in G^u\cap N$) and $\bw\in \widetilde{W^a}$ (resp. $\bw\in W^a$) such that $g.x=n.x$, for all $x\in \A\cap g^{-1}.\A$. 

\item Let $\cV$ be a filter on $\A$. Then  $G_\cV\cap G^u$ (and thus $G_\cV$) acts transitively on the set of apartments containing $\cV$. 
\end{enumerate} 
\end{Corollary}

\begin{proof}
(1) Let $A=g^{-1}.\A$ and $\Omega=A\cap \A$. By \ref{a_ma2}, there exists an apartment isomorphism $f:A\rightarrow \A$  fixing $\Omega$. By Definition/Proposition~\ref{dp_MA1_tildeA} (resp. Definition/Proposition~\ref{dp_MA1_A}), there exists $h\in G$ (resp. $h\in G^u$) inducing $f$ on $\A$. Let $n=gh^{-1}$.  We have $n.\A=\A$ and hence by Proposition~\ref{p_N_stab_A}, $n\in N$ (resp. $n\in G^u\cap N$). Then for all $x\in \Omega$, we  have $n.x=gh^{-1}.x=g.f^{-1}(x)=g.x$. We then set $\bw=\nu(n)$. Then $\bw\in \widetilde{W^a}$ and if $n\in G^u$, we have $\bw\in W^a$, by Lemma~\ref{l_Gu}. 

(2) Let $A$ be an apartment containing $\cV$. Write $A=g.\A$, where $g\in G^u$.  Then $\Omega:=\A\cap g.\A\supset \cV$. Applying (1) with $g^{-1}$ instead of $g$, we deduce the existence of $n\in G^u\cap N$ such that for all $x\in \Omega$, we have $g^{-1}.x=n.x$. Then $gn\in (G^u\cap G_{\Omega})\subset (G^u \cap G_{\cV})$ and $gn.\A=g.\A=\A$. 
\end{proof}

\section{Filters have good fixators}

In this section, we prove that every filter $\cV$ on $\A$ has a good fixator. For this, we first treat the case where $\cV=\Omega$ is a set with non-empty interior. In this case, we have $G_{\Omega}\subset G_{C}=U_{C,\infty}^+\cdot U_{C,\infty}^-\cdot T_0$, for every alcove $C$ contained in $\Omega$. We  then commute $U_{C,\infty}^-$ and $T_0$, which enables to use uniqueness in the Birkhoff decomposition in this case. We then generalize the result to every filter on $\A$.

\begin{Lemma}\label{l_elt_fx_cl}
Let $g\in G$. We assume that $g$ fixes a subset $\Omega$ of $\A$ such that $\Omega$ has non-empty interior. Then $g$ fixes $\cl(\Omega)$. Moreover, if $B=g.\A$, then $g$ fixes $\A\cap A$. 
\end{Lemma}

\begin{proof}
Let $B=g.\A$. Let $h\in G$ be such that $h.\A=B$ and such that $h$ fixes $\A\cap B$, which exists by~\ref{a_ma2}. Then $h^{-1}g.\A=\A$ and $h^{-1}g$ acts as an affine automorphism on $\A$. As it fixes $\Omega$, which has non-empty interior, it is the identity.   Moreover $\A\cap B$ is enclosed (by~\ref{a_ma2}) and contains $\Omega$, so $\A\cap B\in \cl(\Omega)$ and $h$ fixes $\A\cap B$. Therefore $g$ fixes $\A\cap B$, which proves the lemma. 
\end{proof}

\begin{Proposition}\label{p_eq_gd_fx_cl}
Let $\cV$ be a filter on $\A$. Then $G_{\cV}=G_{\cl(\cV)}\cdot  N_{\cV}$. 
\end{Proposition}

\begin{proof}
 We have $N_\cV\subset G_\cV$ and $G_{\cl(\cV)}\subset G_{\cV}$  and  thus $G_\cV\supset G_{\cl(\cV)}\cdot N_\cV$. Let $g\in G_\cV$. Let $B=g.\A$. Let $h\in G$ be such that $h.\A=B$ and $h$ fixes $\A\cap B$. As $\A\cap B\supset \cV$ and as $\A\cap B$ is enclosed, $\A\cap B\supset \cl(\cV)$ and thus $h\in G_{\cl(\cV)}$. Now $h^{-1}g.\A=\A$ thus $h^{-1}g\in N_\cV$, which proves that $g\in hN_{\cV}\subset G_{\cl(\cV)}\cdot N_{\cV}$. Therefore  $G_{\cV}=G_{\cl(\cV)}\cdot N_{\cV}$. 
\end{proof}

For $\Omega$ a subset of $\A$, we denote by $\overline{\conv(\Omega)}$ its closed convex hull.

\begin{Lemma}\label{l_fx_ClsdCv_hl}
Let $\Omega$ be a non-empty subset of $\A$. Then $G_{\Omega}=G_{\overline{\conv(\Omega)}}$. 
\end{Lemma}

\begin{proof}
As $\Omega\subset \overline{\conv(\Omega)}$, we have $G_{\overline{\conv(\Omega)}}\subset G_{\Omega}$. By Proposition~\ref{p_eq_gd_fx_cl}, we have $G_{\Omega}= G_{\cl(\Omega)}\cdot N_{\Omega}$. As an affine automorphism of $\A$ fixes $\Omega$ if and only if it fixes $\overline{\conv(\Omega)}$, we have $N_{\Omega}=N_{\overline{\conv(\Omega)}}$. We also have $\overline{\conv(\Omega)}\subset \cl(\Omega)$ and thus $G_{\cl(\Omega)}\subset G_{\overline{\conv(\Omega)}}$. Therefore $G_{\Omega}\subset G_{\overline{\conv(\Omega)}}\cdot N_{\overline{\conv(\Omega)}}\subset G_{\overline{\conv(\Omega)}}$, and the lemma follows.
\end{proof}

Let $\cV$ be a filter on $\A$. We say that $\cV$ has non-empty interior if for every $\Omega\in \cV$, $\Omega$ has non-empty interior.

\begin{Lemma}\label{l_op_flt_gdFx}
Let $\cV$ be a filter on $\A$ which has non-empty interior. Then $G_{\cV}=U_{\cV,\infty}^+\cdot U_{\cV,\infty}^-\cdot T_0=U_{\cV,\infty}^+\cdot T_0\cdot  U_{\cV,\infty}^-$. In particular, $\cV$ has a good fixator.
\end{Lemma}

\begin{proof}
By Proposition~\ref{p_cnj_N}, $T_0$ normalizes $U_{\cV,\infty}^-$. 
Let $g\in G_{\cV}$. Let $\Omega\in \cV$ be such that $g\in G_\Omega$. Let $B=g.\A$. Then $B$ contains $\Omega$. By~\ref{a_ma2} and Definition/Proposition~\ref{dp_MA1_tildeA}, there exists $h\in G$ such that $B=h.\A$ and $h$ fixes $\A\cap B$. Then $h^{-1}g\in N$ and fixes $\Omega$, which has non-empty interior. Therefore $h^{-1}g\in T_0$ and $g$ fixes $\A\cap B$: $g\in G_{\A\cap B}$. We now set $\Omega=\A\cap B$. 

For $x\in \mathring{\Omega}$, $g$ fixes $C_x:=\germ_x(x+C^v)$ and thus $g\in G_{C_x}$. By Lemma~\ref{l_OCO}, $G_{C_x}=U_{C_x,\infty}^+\cdot U_{C_x,\infty}^-\cdot N_{C_x}$. But as every element of $C_x$ has non-empty interior, $N_{C_x}=T_0$. Therefore we can write $g=u_x^+ u_x^- t_x$, where $u_x^+\in U_{C_x,\infty}^+$, $u_x^-\in U_{C_x,\infty}^-$ and $t_x\in T_0$. Now if $y\in \mathring{\Omega}$, we have: \[
u_x^+u_x^- t_x=u_y^+u_y^- t_y.
\] Using Lemma~\ref{l_utu}, we get $u_x^+=u_y^+$ and $u_x^-=u_y^-$ and $t_x=t_y$. Set $u^+=u_x^+$, $u^-=u_x^-$ and $t=t_x$. 

We have $u^+\in \bigcap_{z\in \mathring{\Omega}} U_{z,\infty}^+$ and $u^-\in \bigcap_{z\in \mathring{\Omega}} U_{z,\infty}^-$. Using Corollary~\ref{c_U_ma_inter}, we deduce $u^+\in U_{\mathring{\Omega},\infty}^+$ and $u^-\in U_{\mathring{\Omega},\infty}^-$. We also have $U_{\mathring{\Omega},\infty}^+=U_{\overline{\mathring{\Omega}},\infty}^+$ and $U_{\mathring{\Omega},\infty}^-=U_{\overline{\mathring{\Omega}},\infty}^{-}$ and since $\Omega$ is convex and closed we deduce $g\in U_{\overline{\mathring{\Omega}},\infty}^+\cdot U_{\overline{\mathring{\Omega}},\infty}^-\cdot T_0= U_{\Omega,\infty}^+\cdot U_{\Omega,\infty}^- \cdot T_0\subset U_{\cV,\infty}^+ \cdot U_{\cV,\infty}^-\cdot T_0$. 
\end{proof}

\begin{Lemma}\label{l_csq_char_apt}(see \cite[Lemma 4.7]{bardy2025twin})
Assume that $\cS$ is positively cofree. Let $x\in \A$, $u_+\in U^+$ and $u_-\in U^-$. Assume that $u_+ u_-.x=x$. Then $u_+.x=u_-.x=x$. 
\end{Lemma}

\begin{proof}
Set $a=u_-.x=(u_+)^{-1}.x$. By definition, $\rho_{-\infty}(a)=\rho_{-\infty}(u_-.x)=x=\rho_{+\infty}((u_+)^{-1}.x)=\rho_{+\infty}(a)$. Using Corollary~\ref{c_char_aprtmt} we deduce $a\in \A$. Therefore $a=\rho_{-\infty}(a)=x$, which proves the lemma.
\end{proof}

\begin{Theorem}\label{t_fltrs_gd_fix}
Let $\cV$ be a filter on $\A$. Then we have:
 \begin{equation}\label{e_good_fixator}G_{\cV}=U_{\cV,\infty}^+\cdot U_{\cV,\infty}^- \cdot N_{\cV}=U_{\cV,\infty}^-\cdot U_{\cV,\infty}^+\cdot  N_{\cV}.
\end{equation}

Moreover: \begin{enumerate}
\item  $U_{\cV,\infty}^+=G_\cV\cap U^+$, $U_{\cV,\infty}^-=G_{\cV}\cap U^-$ and $N_\cV=N\cap G_\cV$.

\item If $w\in W^v$, then we can replace $(U_{\cV,\infty}^+,U_{\cV,\infty}^-)$ by $(U_{\cV,\infty}(w.\Delta_+),U_{\cV,\infty}(w.\Delta_-))$ in \eqref{e_good_fixator}. 
\end{enumerate} 

\end{Theorem}

\begin{proof}
We first assume that $\cS$ is positively cofree. We start by proving \eqref{e_good_fixator} in the case where $\cV$ is a convex closed set. Let $\Omega$ be a closed and convex non-empty subset of $\A$. Let $g\in G_\Omega$.  Let $B=g.\A$. Using Proposition~\ref{p_splt_apt},  we  write $B=\bigcup_{i=1}^k P_i$, where $k\in \N$ and  for $i\in \llbracket 1,k\rrbracket$, $P_i$ is enclosed, have non-empty interior and is contained in an apartment containing $\fQ_{+\infty}$.   Let $S=\supp(\Omega)$. As $\Omega=\bigcup_{i=1}^k \Omega\cap P_i$ and as the $\Omega\cap P_i$ are convex, there exists $i\in \llbracket 1,k\rrbracket$ such that $\Omega_i:=\Omega\cap P_i$ has non-empty interior in $S$. 

Let $A_i$ be an apartment containing $P_i\cup\fQ_{+\infty}$. Let $u\in U^+$ be such that $A_i=u.\A$. Let $h\in G$ be such that $h.A_i=B$ and $h$ fixes $A_i\cap B$. Then $u^{-1}hu.\A=u^{-1}.B$ and $u^{-1}hu$  fixes $\Omega'_i:=u^{-1}.P_i$ pointwise. In other words, $u^{-1}hu \in G_{\Omega'_i}$.  As $\Omega'_i\subset u^{-1}. (A_i\cap B)$ has non-empty interior we can apply Lemma~\ref{l_op_flt_gdFx} and we can write $u^{-1}hu=v_+ v_-t$, with $v_+\in U_{\Omega'_i,\infty}^+$, $v_-\in U_{\Omega'_i,\infty}^-$ and $t\in T_0$. Then $hu=uv_+ v_-t$. Note that since $u$ fixes $\A\cap A_i$, we have \[\Omega_i\subset \Omega'_i\] and $hu$ fixes $\Omega_i$ pointwise. Now  $hu.\A=B=g.\A$ and thus $n:=(hu)^{-1}g$ belongs to $N$ and fixes $\Omega_i$: $n\in N_{\Omega_i}$. As $n$ induces an affine automorphism of $\A$, it fixes $\supp(\Omega_i)=S$ and in particular, $n\in N_\Omega$.  Therefore $g=hun=uv_+ v_-tn=u'_+u'_-n'$, where $u'_+=uv_+\in U^+$, $u'_-=v_-\in U^-$ and $n'=tn\in N_{\Omega}$. 

If $x\in \Omega$, then $g.x=x=u'_+u'_-.x$. Using Lemma~\ref{l_csq_char_apt} we deduce $u'_+.x=u'_-.x=x$. Therefore $u'_+\in \bigcap_{x\in \Omega}U_{x,\infty}^+$ and  $u'_-\in \bigcap_{x\in \Omega}U_{x,\infty}^-$.  Using Corollary~\ref{c_U_ma_inter}, we deduce $u^+\in U_{\Omega,\infty}^+$ and $u^-\in U_{\Omega,\infty}^-$. We thus proved $g\in U_{\Omega,\infty}^+\cdot U_{\Omega,\infty}^-\cdot  N_\Omega$. Therefore $G_\Omega\subset U_{\Omega,\infty}^+\cdot  U_{\Omega,\infty}^-\cdot N_\Omega$. Using Remark~\ref{r_chng_ord_fix}, we get  the left equality of \eqref{e_good_fixator}: $G_\Omega=U_{\Omega,\infty}^+\cdot  U_{\Omega,\infty}^-\cdot  N_\Omega$. Similarly, replacing $\fQ_{+\infty}$ by $\fQ_{-\infty}$,  we obtain \eqref{e_good_fixator} in this case.

Let now $\Omega$ be a non-empty subset of $\A$. Then by  Lemma~\ref{l_fx_ClsdCv_hl}, we have: \[G_{\Omega}=G_{\overline{\conv(\Omega)}}=U_{\overline{\conv(\Omega)},\infty}^+\cdot U_{\overline{\conv(\Omega)},\infty}^-\cdot N_{\overline{\conv(\Omega)}}\subset U_{\Omega,\infty}^+\cdot U_{\Omega,\infty}^-\cdot N_{\Omega}.\] Using Remark~\ref{r_chng_ord_fix}, we deduce \eqref{e_good_fixator}. 

Let now $\cV$ be a filter on $\A$ and $g\in G_{\cV}$. Let $\Omega\in \cV$ be such that $g\in G_\Omega$. Then $g\in U_{\Omega,\infty}^+\cdot  U_{\Omega,\infty}^{-} \cdot  N_\Omega\subset U_{\cV,\infty}^+\cdot  U_{\cV,\infty}^- \cdot N_\cV$. By Remark~\ref{r_chng_ord_fix}, we have $G_{\cV}=U_{\cV,\infty}^+\cdot U_{\cV,\infty}^{-} \cdot N_\cV$ and similarly, we obtain \eqref{e_good_fixator} when $\cS$ is positively cofree.  

We no longer assume $\cS$ to be positively cofree. We consider a cofree extension $\tilde{\cS}$ of $\cS$ satisfying Assumption~\ref{as_ext}, which exists by  \cite[7.3.2 and Proposition 7.64]{marquis2018introduction}.  We keep the same notations as in Section~\ref{s_ext_mas}. In particular, $\widetilde{G}$ is a central extension of $G$ and we have a projection $\pi:\widetilde{G}\twoheadrightarrow G=\widetilde{G}/\overline{T}$ (where $\overline{T}$ is central in $\widetilde{G}$), which induces a projection $\pi:\widetilde{\I}\twoheadrightarrow \I$ (see Proposition~\ref{p_pr_ext_mas}). Let $\cV$ be a filter on $\A$. We also regard $\cV$ as a filter on $\tilde{\A}$. 

Let $g\in G_{\cV}$. Let $\Omega\in \cV$ be such that $g\in G_{\Omega}$. Let $\tilde{g}\in \pi^{-1}(\{g\})$. Let $x\in \Omega$. Then $\pi(\tilde{g}.x)=g.x=x$ and thus $\tilde{g}.x\in \tilde{\A}$, by Lemma~\ref{l_pr_img_apt}. Therefore $\tilde{g}\in \tilde{G}(\Omega\curvearrowright  \A)$. Using Corollary~\ref{c_w_mas_pf} and \eqref{e_good_fixator}, we have $\tilde{G}(\Omega\curvearrowright \A)=\tilde{N}\cdot \tilde{G}_{\Omega}=\tilde{N}\cdot \tilde{U}_{\Omega,\infty}^-\cdot \tilde{U}_{\Omega,\infty}^+$. Consequently, \[
g=\pi(\tilde{g})\in \pi(\tilde{N})\cdot \pi(\tilde{U}^-_{\Omega,\infty})\cdot \pi(\tilde{U}^+_{\Omega,\infty}).\]

Using Lemma~\ref{l_pi_tilde_G_x}, we deduce $g\in N\cdot U_{\Omega,\infty}^-\cdot U_{\Omega,\infty}^+$. Thus $G_{\cV}\subset N\cdot U_{\cV,\infty}^-\cdot U_{\cV,\infty}^+$, which implies the left side of \eqref{e_good_fixator}, by Lemma~\ref{l_sfct_cnd_GF}. By symmetry, we deduce \eqref{e_good_fixator}, without cofreeness assumption on $\cS$.

By definition, we have $N_\cV=N\cap G_\cV$. We have $U_{\cV,\infty}^+\subset G_\cV\cap U^+$. Let $u_+\in G_\cV\cap U^+$. By \eqref{e_good_fixator}, we can write $u_+=v_+ v_-n$, with $v_+\in U_{\Omega,\infty}^+\subset U^+$, $v_-\in U_{\Omega,\infty}^-\subset U^-$ and $n\in N$. Then $v_+^{-1}u_+=v_-n$ and $v_-n u_+^{-1} v_+=1$. Using the Birkhoff decomposition of $G$, we deduce $u_+=v_+$ and thus $u_+\in U_{\cV,\infty}^+$.  Therefore $U_{\cV,\infty}^+= G_\cV\cap U^+$. Similarly we get $U_{\cV,\infty}^-=G_\cV\cap U^-$, which proves (1).

Let $w\in W^v$. Set    $\sU_{\cV,\infty}=U_{\cV,\infty}^+\cdot U_{\cV,\infty}^-\cdot N^u_\cV$. We have $\sU_{\cV,\infty}=U_{\cV,\infty}^+\cdot U_{\cV,\infty}^-\cdot N^u_\cV=U_{\cV,\infty}(w.\Delta_+)\cdot U_{\cV,\infty}(w.\Delta_-)\cdot N^u_\cV$ by  Lemma~\ref{l_pre_GR08}. Therefore $G_{\cV}=\sU_{\cV,\infty}\cdot N_{\cV}= U_{\cV,\infty}(w.\Delta_+)\cdot U_{\cV,\infty}(w.\Delta_-)\cdot N_{\cV}$ which completes the proof of (2).
\end{proof}

\begin{Corollary}\label{c_U_clV=UV}
Let $\cV$ be a filter on $\A$. Then $U_{\cV,\infty}^+=U_{\cl(\cV),\infty}^+$ and $U_{\cV,\infty}^-=U_{\cl(\cV),\infty}^-$.
\end{Corollary}

\begin{proof}
We have $\cV\subset \cl(\cV)$ and thus $U_{\cl(\cV),\infty}^+\subset U_{\cV,\infty}^+$. Conversely, let $u\in U_{\cV,\infty}^+$. Let $\Omega\in \cV$ be such that $u\in U_{\Omega,\infty}^+$. As $u$ fixes $\fQ_{+\infty}$ which has non-empty interior, we can apply Lemma~\ref{l_elt_fx_cl} to  $u$: $u$  fixes $\cl(\Omega)$.  Therefore $u\in G_{\cl(\cV)}\cap U^+$ and hence $u\in U_{\cl(\cV),\infty}^+$, by Theorem~\ref{t_fltrs_gd_fix}. Therefore $U_{\cV,\infty}^+=U_{\cl(\cV),\infty}^+$ and the same holds for $U_{\cV,\infty}^-$, by symmetry. 
\end{proof}

\begin{Corollary}\label{c_fltrs_gd_fix}
Every filter on $\A$ has a good fixator.  In particular, $G_\cV$ acts transitively on the set of apartments containing $\cV$. 
\end{Corollary}

\begin{proof}
This is  a combination of Corollary~\ref{c_w_mas_gen}, Theorem~\ref{t_fltrs_gd_fix} and Lemma~\ref{l_TF}. 
\end{proof}

Note that Corollary~\ref{c_fltrs_gd_fix} contradicts \cite[Exemple 4.12 3 (c)]{rousseau2016groupes}. Actually, there is a mistake in this counter-example, see \cite[7]{rousseau2025split} for more explanations.

\section{Chimneys and axiom~\ref{a_ma3}}\label{s_chimney}

In this section, we prove that $\I$ satisfies axiom~\ref{a_ma3}. We start by describing the fixator of a chimney and of its germ (see Proposition~\ref{p_fx_chimn}). We then prove a decomposition of $G$, which implies~\ref{a_ma3} (see Proposition~\ref{p_chmny_frndl}).

\subsection{Masure associated with a set of simple roots}\label{ss_mas_smpl_rt}

Let $J\subset I$. We set  $\Phi_J=(\bigoplus_{i\in J} \R \alpha_i)\cap \Phi$\index[notation]{p@$\Phi_J$}.  By \cite[7.4]{marquis2018introduction}, the group $G_J=\langle T,U_\alpha\mid \alpha\in \Phi_J\rangle$ naturally identifies with $\fG_{\cS|_{J}}(\cF)$, where $\cS|_J=((a_{i,j})_{i,j\in J},X,Y,(\alpha_i)_{i\in J},(\alpha_i^\vee)_{i\in J})$. Recall that $\A_J$\index[notation]{a@$\A_J$} is the set $\A$ equipped with the hyperplane arrangement $\{\alpha^{-1}(\{k\})\mid k\in \Lambda,\alpha\in \Phi_J\}$. The group $G_J$ acts on two masures: the masure $\I_J=\I(\fG_{\cS|_{J}},\cF,\omega)$, whose standard apartment is  $\A_J$, and the ``masure'' $G_J.\A\subset \I$. We prove in this subsection that they can be identified. In order to distinguish these two masures, we denote by $g\odot x$ the element of $\I_J$,
 which is the translate of $x$ by $g$. We add an index $J$ to denote the objects relative to $\I_J$: for example $U_{J,\fin}^+=\langle U_\alpha\mid \alpha\in\Phi_+\cap  \Phi_J\rangle$,  $U_{x,\fin,J}^+=\langle U_{\alpha,x}\mid \alpha\in \Phi_J\rangle\cap U^+$, etc. Note that for $n\in N_J=\langle
 T,\tilde{s_j}\mid j\in J\rangle$ and $x\in \A_{J}$, we have $n\odot x=n.x$ since these actions are defined similarly. Therefore for $x\in \A$, we have: \begin{equation}\label{e_N_x_J}
 N_{x,J}=N_x\cap N_J
 \end{equation}

\begin{Lemma}\label{l_GJ_cap_G_x}
Let $x\in \A$. Then $G_J\cap G_x=G_{x,J}$, where $G_{x,J}$ is the fixator of $x$ in $(G_J,\odot)$.
\end{Lemma}

\begin{proof}
By  Corollary~\ref{c_BPHR_2.4.1}, we have $U_{x,\infty,J}^\epsilon=U_{x,\fin,J}^\epsilon \subset U_{x,\fin}^\epsilon=U_{x,\infty}^\epsilon$ for both $\epsilon\in \{-,+\}$. By definition, the action of  $N_J$ on $\A$ is the restriction of the action of $N$ on $\A$. We deduce that $G_{x,J}=\langle U_{x,\infty,J}^+,U_{x,\infty,J}^-,N_{x,J}\rangle\subset\langle U_{x,\infty}^+,U_{x,\infty}^-,N_{x}\rangle=G_x$.

Let $g\in G_J\cap G_x$. By the Iwasawa decomposition (Proposition~\ref{p_Iwa}) applied to $G_J$, we can  write $g=unk$, with $u\in U_J^+$, $n\in N_J$ and $k\in U_{x,\fin,J}$. Then $g.x=un.x=x$, since $U_{x,\fin,J}\subset G_{x,J}\subset G_x$. Therefore $\rho_{+\infty}(g.x)=x=n.x$. Consequently $n\in N_{x}\cap N_J=N_{x,J}$ and hence $u\in U_{x,\infty}^+\cap U_J^+$.

 Let $\Delta_{J,+}=\bigoplus_{j\in J}\R\alpha_j\cap \Delta_+$. Then $\Delta_{J,+}$ is closed and coclosed in $\Delta_+$ and thus by Proposition~\ref{p_Mar8.58}, we have a unique decomposition $U^{ma+}=\fU_{\Delta_{J,+}}(\cF)\cdot \fU_{\Delta_+\setminus \Delta_{J,+}}(\cF)$. Then $u\in \langle U_{\alpha_j}\mid j\in J\rangle\subset \fU_{\Delta_{J,+}}(\cF)$. Write $u=\prod_{\alpha\in \Delta_{J,+}} X_\alpha(\underline{u_\alpha})$, where $(\underline{u_\alpha})\in \prod_{\alpha\in \Delta_{J,+}}\ffg_{\alpha,\Z}\otimes \cF$. Then by Proposition~\ref{p_Rou4.5}, $\omega(\underline{u_\alpha})+\alpha(x)\geq 0$ for all $\alpha\in \Delta_{J,+}$ (where $\omega(\underline{u_\alpha})$ is defined in \ref{sss_ext_val}) and thus $u\in U_{x,\infty,J}^+\subset G_{x,J}$. Therefore $g\in G_{x,J}$, and the lemma follows.
\end{proof}

\begin{Proposition}\label{p_I_J}
The map $\iota:\I_J\rightarrow G_J.\A$ defined by $g\odot x\mapsto g.x$, for $g\in G_J$ and $x\in \A_J$ is well-defined and is a bijection. 
\end{Proposition}

  \begin{proof}
  Let $(g,x),(h,y)\in G_J \times \A_J$ be such that $g\odot x=h\odot y$. Then $h^{-1}g\odot x=y$. Therefore by \eqref{e_def_sim} applied to $\I_J$, there exists $n\in N_J$ such that $y=n\odot x$ and $n^{-1}hg\in G_{x,J}$. By Lemma~\ref{l_GJ_cap_G_x}, $n^{-1}h^{-1}g.x=x$ and $g.x=hn.x=h.y$. This proves that $\iota$ is well-defined and it  is clearly surjective (since $\A_J=\A$ as a set).

Let $x,y\in \A_J$ and $g\in G_J$ be such that $y=g.x$.  Using the Iwasawa decomposition in $G_J$ (Proposition~\ref{p_Iwa}), we write $g=u_+nk$, with $u_+\in U_J^+$, $n\in N_J$ and $k\in G_{x,J}$. Then $y=g.x=\rho_{+\infty}(g.x)=n.x$. Therefore $n^{-1}g\in G_J\cap G_x=G_{x,J}$ (by Lemma~\ref{l_GJ_cap_G_x}). Thus $n^{-1}g\odot x=x$ and hence $g\odot x=n\odot x=n.x=y$. This proves that $\iota$ is injective and completes the proof of the proposition.
  \end{proof}

 \begin{Corollary}\label{c_act_G_J}
 Let $x\in \A$ and $g\in G_x\cap G_J$. Then $g$ fixes pointwise $x+V_J$, where $V_J=\bigcap_{i\in J} \ker(\alpha_i)$.
 \end{Corollary}
 
 \begin{proof}
By Proposition~\ref{p_I_J}, we  have $g\in G_{x,J}$. We have $\A_{\ines,J}=V_J$. Thus by Lemma~\ref{l_Fx_Ain}, we have $g\odot (x+y)=x+y$, for all $y\in V_J$ and hence $g.(x+y)=x+y$ for all $y\in V_J$. 
 \end{proof}

\subsection{Action of $U^+_{\infty}(\Delta^u_+(F^v))$ on $\I$}

Let $F^v$ be a vectorial face of $\pm\cT$. Recall that $\Delta^u_+(F^v)$, $\Phi^u_+(F^v)$, $\Delta_m^+(F^v)$ and $\Phi^m_+(F^v)$ are defined in Subsection~\ref{ss_Lv_dec}.

In this subsection, we study the action of the (defined below) subgroup $U^+_{\infty}(\Delta^u_+(F^v))$  of $U^+$ on $\I$. We give an interpretation of this group as the fixator of a filter on $\A$ (see Proposition~\ref{p_fix_V_Fv}). This study is a first step in the description of the fixator of a chimney.

 We set: \[\Phi(F^v)=\{\alpha\in \Phi\mid \alpha(F^v)\geq 0\}, \Phi^m(F^v)=\{\alpha\in \Phi\mid \alpha(F^v)=\{0\}\}\]\index[notation]{p@$\Phi(F^v)$}\index[notation]{p@$\Phi^v(F^v)$} and \[\Phi^u_+(F^v)=\Phi(F^v)\setminus \Phi^m(F^v).\]\index[notation]{p@$\Phi^u_+(F^v)$}

We call a subset of $\Phi$ of the form $\Phi(F^v)$ a \textbf{parabolic subset of roots}\index{parabolic subset of roots}.

By symmetry we can assume that $F^v\subset \cT$. Up to changing the choice of the fundamental chamber, we can assume that $F^v$ is contained in $\overline{C^v_f}$. The \textbf{parabolic subgroup} $P(F^v)$\index[notation]{p@$P(F^v)$} is the group $\langle T,U_\alpha\mid \alpha\in \Phi(F^v)\rangle$.   The \textbf{Levi factor} $M(F^v)$\index[notation]{m@$M(F^v)$} of $P(F^v)$ is $\langle T,U_\alpha\mid \alpha\in \Phi^m(F^v)\rangle$. We have $M(F^v)=G_J$, where $J=\{i\in I\mid \alpha_i(F^v)=\{0\}\}$, for the notation of Subsection~\ref{ss_mas_smpl_rt}.  Recall that the  \textbf{unipotent radical} $U^+(F^v)$\index[notation]{u@$U^+(F^v)$} of $P(F^v)$ is the smallest normal sub-group of $U^+$ containing the $U_\alpha$, for $\alpha\in \Phi^u_+(F^v)$. The group $M(F^v)$ can be identified with $\fG_{\cS|_J}(\cF)$ (see \cite[7.4.6]{marquis2018introduction}) and we have the following decomposition: \begin{equation}\label{e_Rem_6.2}
P(F^v)=M(F^v)\ltimes U^+(F^v),
\end{equation} by \cite[6.2.2]{remy2002groupes}.

If $J\subset I$, we set $W^v_J=\langle r_j\mid j\in J\rangle\subset W^v$. This is the fixator of $F^v(J)$ in $W^v$.

\begin{Lemma}\label{l_parb_rt_sys}
Let $J\subset I$ be non-empty. Then $(\bigoplus_{j\in J} \R\alpha_j)\cap \Phi=W^v_J\cdot \{\alpha_j\mid j\in J\}$.
\end{Lemma}

\begin{proof}
Set $\Phi_J=(\bigoplus_{j\in J} \R\alpha_j)\cap \Phi$ and $\Phi_J' =W^v_J\cdot \{\alpha_j\mid j\in J\}\subset (\bigoplus_{j\in J} \R\alpha_j)$. It  is clear that $\Phi_J'\subset \Phi_J$.  We prove the converse inclusion by induction on $|\htt(\alpha)|$. 

Let $\alpha\in \Phi_J$. If $|\htt(\alpha)|=1$, then $\alpha\in  \pm \{\alpha_j\mid j\in J\}\subset \Phi_J'$. Let $k\in \Z_{\geq 1}$. Assume that for all $\alpha\in \Phi_J$ such that $|\htt(\alpha)|\leq k$, we have $\alpha\in \Phi_J'$.  Let $\alpha\in \Phi_J$. Up to replacing $\alpha$ by $-\alpha$, we can assume that $\alpha\in \Phi_+$.  Assume that $\htt(\alpha)= k+1$. Then by \cite[Proposition 6.2 (5)]{marquis2018introduction}, there exists $i \in I$ such that $|\htt(r_i.\alpha)|<\alpha$. By \eqref{e_Kmr_1.3.14}, we have $r_i.\alpha\in \Phi_+$. Write $\alpha=\sum_{j\in J} n_j\alpha_j$, where $(n_j)\in \N^{J}$. We have $r_i.\alpha=\sum_{j\in J} n_j r_i.\alpha_j=\sum_{j\in J}(\alpha_j-\alpha_j(\alpha_i^\vee)\alpha_i)$. As the $\alpha_j(\alpha_i^\vee)$, $j\in J$ are the coefficients of the Kac--Moody matrix defining $\cS$, there are non-positive if $i\in I\setminus J$. Therefore we necessarily have $i\in J$. Thus by assumption there exist $w\in W^v_J$ and $k\in J$ such that $r_i.\alpha=w.\alpha_k$ and hence $\alpha=r_iw.\alpha_k\in \Phi_J'$. We deduce that $\Phi_J\subset \Phi_J'$. 
\end{proof}

We now fix $F^v$ and we assume that $F^v\subset \overline{C^v_f}$. Recall that $W^v_{F^v}$ is the fixator of $F^v$ in $W^v$. Let $\sC_{F^v}=\conv(\bigcup_{w\in W^v_{F^v}}w.C^v_f)$.\index[notation]{c@$\sC_{F^v}$}

\begin{Lemma}\label{l_cl_Delta_u}
 Let $\alpha\in \Phi$. Assume that there exists $x\in \A$ and $\lambda\in \R$ such that $D(\alpha,\lambda)\supset x+\sC_{F^v}$. Then $\alpha(F^v)>0$. 
\end{Lemma} 

\begin{proof}
As $\sC_{F^v}\supset F^v$, we necessarily have $\alpha(F^v)+\alpha(x)+\lambda\geq 0$ and thus $\alpha(F^v)\geq 0$. Let $J=\{i\in I\mid \alpha_i(F^v)=\{0\}\}$.  Let $\alpha\in \Phi$ be such that $\alpha(F^v)=0$. Then $\alpha\in \Phi_J=\Phi\cap \bigoplus_{i\in J}\R \alpha_i$ and thus we can write $\alpha=w.\alpha_i$, for some $w\in W^v_J$ and $i\in J$, by Lemma~\ref{l_parb_rt_sys}. Then $w.C^v_f\cup wr_i C^v_f\subset \sC_{F^v}$. As $\alpha(w.C^v_f\cup wr_i C^v_f)=\R^*$, we deduce that $D(\alpha,\lambda)$ does not contain $x+\sC_{F^v}$, for every $x\in \A$. 
\end{proof}

We set: \[\cV_{F^v}=\{\Omega\subset \A\mid  \exists k\in \N, \beta_1,\ldots,\beta_k\in \Phi^u_+(F^v),\lambda_1,\ldots,\lambda_k\in\R, \Omega\supset \bigcap_{i=1}^k D(\beta_i,\lambda_i)\}.\]\index[notation]{v@$\cV_{F^v}$} Note that every element of $\cV_{F^v}$ is non-empty since if $\xi\in F^v$, $k\in \N$, $\beta_1,\ldots,\beta_k\in \Phi^u_+(F^v)$ and $\lambda_1,\ldots,\lambda_k\in \R$, we have $t\xi \in \bigcap_{i=1}^k D(\beta_i,\lambda_i)$ for $t\in \R_{>0}$ large enough. Therefore $\cV_{F^v}$ is a filter on $\A$. It corresponds to the germ at infinity of  $\sC_{F^v}$.

\begin{Lemma}\label{l_fR_sub_V_Fv}
Let $\Omega\in \cV_{F^v}$. Then for every  subset $\Omega'$ of $\A$ such that $\Omega'/\A_{\ines}$ is bounded, there exists $\xi\in F^v$ such that $\Omega'+\xi+F^v\subset \Omega$. In particular, if $\fr$ is a chimney of direction $F^v$, then $\cV_{F^v}$ contains its germ at infinity $\fR$. 
\end{Lemma}

\begin{proof}
By assumption, there exist $k\in \N$,  $\beta_1,\ldots,\beta_k\in \Phi^u_+(F^v)$ and $\lambda_1,\ldots,\lambda_k\in \R$ such that $\Omega\supset  \bigcap_{i=1}^k D(\beta_i,\lambda_i)$. Let $\Omega'$ be a non-empty  subset of $\A$ such that $\A/\A_{\ines}$ is bounded. Then if $\xi\in F^v$ is such that $\beta_i(\xi)\geq \sup (-\beta_i(\Omega'))$ for all $i\in \llbracket 1,k\rrbracket$, we have $\Omega'+\xi+F^v\subset \Omega$, which proves the first part of the lemma.

Let $\fr=\fr(F,F^v)$ be a chimney of direction $F^v$ where $F$ is a face of $\A$. Let $\fR=\germ_\infty(\fr)$. Let $\Omega\in \cV_{F^v}$. Then there exists $\xi\in F^v$ such that $F+\xi+F^v\subset \Omega$. Therefore $\Omega\in \fR$, which proves that $\fR$ is contained in $\cV_{F^v}$ (as filters).
\end{proof}

We defined  $U^+_{\infty}(\Delta^u_+(F^v))=\fU^{ma}_{\Delta^u_+(F^v)}(\cF)\cap U^+$ and proved that it is equal to $U^+(F^v)$ in Proposition~\ref{p_Levi_dec_U}.

\begin{Proposition}\label{p_fix_V_Fv}
We have $U^+_{\infty}(\Delta^u_+(F^v))=U_{\cV_{F^v},\infty}^+$. In other words,  $U^+_{\infty}(\Delta^u_+(F^v))$ is the fixator of $\cV_{F^v}$ in $U^+$. 
\end{Proposition}

\begin{proof}
Let $u\in U^+_{\infty}(\Delta^u_+(F^v))$ and $A:=u.\A$. By Lemma~\ref{l_elt_fx_cl},  $u$ fixes $\A\cap A$. By~\ref{a_ma2}, there exist $k\in  \N$, $\beta_1,\ldots,\beta_k\in \Phi$ and $\lambda_1,\ldots,\lambda_k\in \R$ such that $u.\A\cap A=\bigcap_{i=1}^k D(\beta_i,\lambda_i)$.

 Write $u=\prod_{\alpha\in \Delta^u_+(F^v)} X_\alpha(\underline{u_\alpha})$, where $(\underline{u_\alpha})\in \prod_{\alpha\in \Delta^u_+(F^v)} \ffg_{\alpha,\Z}\otimes \cF$. Let $w\in W^v_{F^v}$ and $\alpha\in \Delta^u_+(F^v)$. Then $\overline{\alpha(w.C^v_f)}\supset \alpha(w.F^v)=\alpha(F^v)=\R_{>0}$. Therefore $\alpha(w.C^v_f)=\R_{>0}$. As $u\in U^+$, there exists $x\in \A$ such that $u$ fixes $x$. Then for all $\alpha\in \Delta^u_+(F^v)$, we have $\alpha(x)+\omega(\underline{u_\alpha})\geq 0$ and thus  if $w\in W^v_{F^v}$, we have $\alpha(x+w.C^v_f)+\omega(\underline{u_\alpha})\geq 0$. Therefore $u$ fixes $x+\sC_{F^v}$. Using Lemma~\ref{l_cl_Delta_u} we deduce that  for all $i\in \llbracket 1,k\rrbracket$, we have $\beta_i(F^v)>0$. Therefore $u.\A\cap \A\in \cV_{F^v}$. Using Theorem~\ref{t_fltrs_gd_fix}, we deduce that $u\in U_{\cV_{F^v},\infty}^+$ and hence $U^+_{\infty}(\Delta^u_+(F^v))\subset U_{\cV_{F^v},\infty}^+$.
 
 Conversely, let $u\in U_{\cV_{F^v},\infty}^+$. Let $\Omega\in \cV_{F^v}$ be fixed by $u$. Write $u=\prod_{\alpha\in \Delta_+}X_\alpha(\underline{u_\alpha})$, where $(\underline{u_\alpha})\in \prod_{\alpha\in \Delta_+}(\ffg_{\alpha,\Z}\otimes \cF)$. Let $\alpha\in \Delta^m_+(F^v)$. Let $x\in \A$. Then by Lemma~\ref{l_fR_sub_V_Fv}, there exists $\xi\in F^v$ such that $x+\xi+F^v\subset \Omega$. Therefore $\omega(\underline{u_\alpha})+\alpha(x+\xi+F^v)=\omega(\underline{u_\alpha})+\alpha(x)\geq 0$. As this is true for every $x\in \A$, we deduce that $\omega(\underline{u_\alpha})$ is infinite, which proves that $\underline{u_\alpha}=0$. As $\Delta_+\setminus \Delta^m_+(F^v)=\Delta^u_+(F^v)$, we deduce that $u\in U^+_{\infty}(\Delta^u_+(F^v))$, which proves the inclusion $U_{\cV_{F^v},\infty}^+\subset U^+_{\infty}(\Delta^u_+(F^v))$. Lemma follows.
\end{proof}

\subsection{Fixator of a chimney}

In this subsection, we describe the fixator of a chimney and of its germ.

\begin{Lemma}\label{l_U_fR}
Let $\fr=\fr(F,F^v)$ be a chimney and $\fR$ be its germ at infinity. Then $U_{\fR,\infty}^+=U^+_{\infty}(\Delta^u_+(F^v))\rtimes U_{F,\infty}^+(\Delta^m_+(F^v))$. 
\end{Lemma}

\begin{proof}
Let $u\in U^+$. Then by Proposition~\ref{p_Levi_dec_U}, we can write $u=u_1u_2$, where $u_1\in U^+_{\infty}(\Delta^u_+(F^v))$ and $u_2\in U^+_{\infty}(\Delta^m_+(F^v))$. By Proposition~\ref{p_fix_V_Fv} and Lemma~\ref{l_fR_sub_V_Fv}, $u_1$ fixes $\fR$. Therefore $u\in U_{\fR,\infty}^+$ if and only if $u_2\in U_{\fR,\infty}^+$ and $U_{\fR,\infty}^+=U^+_{\infty}(\Delta^u_+(F^v))\rtimes U_{\fR,\infty}^+(\Delta^m_+(F^v))$

Let $u\in U^+_{\infty}(\Delta^m_+(F^v))$. Write $u=\prod_{\alpha\in \Delta^m_+(F^v)} X_{\alpha}(\underline{u_\alpha})$, with $(\underline{u_\alpha})\in \prod_{\alpha\in \Delta^m_+(F^v)} \ffg_{\alpha,\Z}\otimes \cF$. Let $V_{F^v}$ be the vector space spanned by $F^v$. Let $x\in \A$. Then: \begin{align*}
u.x=x &\Leftrightarrow \forall \alpha\in \Delta^m_+(F^v), \omega(\underline{u_\alpha})+\alpha(x)\geq 0 \\
&\Leftrightarrow \forall \alpha\in \Delta^m_+(F^v),\omega(\underline{u_\alpha})+\alpha(x+V_{F^v})\geq 0\\
& \Leftrightarrow \forall y\in V_{F^v}, u.(x+y)=x+y.
\end{align*}

Therefore $U_{\fR,\infty}^+(\Delta^m_+(F^v))=U_{\fR+V_{F^v},\infty}^+(\Delta^m_+(F^v))=U_{F+F^v,\infty}^+(\Delta^m_+(F^v))=U_{F,\infty}^+(\Delta^m_+(F^v))$, which proves the lemma.
\end{proof}

For $F$ a face of $\A$ and $F^v$ a vectorial face of $\A$, we denote by $M(F,F^v)$ the fixator of $F$ in $M(F^v)$. 

The following proposition generalizes \cite[Proposition 6.8]{rousseau2011masures}.
 
\begin{Proposition}\label{p_fx_chimn}
Let $F^v$ be a vectorial face of $\overline{C^v_f}$. Let $F$ be a face of $\A$ and $\fr=\fr(F,F^v)$. Let $\fR$ be the germ at infinity of $\fr$. Then $G_{\fr}=M(F,F^v)\ltimes U_{F,\infty}^+(\Delta^u_+(F^v))$ and \begin{align*}
G_{\fR}&=M(F,F^v)\ltimes U^+_{\infty}(\Delta^u_+(F^v))\\
&=\left(N_{\fR}\cdot U_{F,\infty}^-(\Delta^m_-(F^v))\cdot U^+_{F,\infty}(\Delta^m_+(F^v))\right)\ltimes U^+_{\infty}(\Delta^u_+(F^v)).
\end{align*} 
\end{Proposition} 

\begin{proof}
By Theorem~\ref{t_fltrs_gd_fix}, we have $G_{\fR}=N_{\fR}\cdot U_{\fR,\infty}^-\cdot U_{\fR,\infty}^+$.  For $w\in W^v\setminus W^v_{F^v}$, we have $w.F^v\neq F^v$ and thus $w.\fR\neq \fR$. Therefore $N_{\fR}\subset \{n\in N\mid \nu^v(n)\in W^v_{F^v}\}$. Let  $J=\{i\in I\mid \alpha_i(F^v)=\{0\}\}$. By \eqref{e_Wv_FJ}, we have $W^v_{F^v}=\langle r_j\mid j\in J\rangle$. Therefore  we have $N_{\fR}\subset N_J$.  Consequently, $N_{\fR}=N_{F,J}$, the fixator of $F$ in $N_J$.

Let $u\in U_{\fR,\infty}^-$. Then we can write $u=\prod_{\alpha\in \Delta_-} X_\alpha(\underline{u_\alpha})$, with $(\underline{u_\alpha})\in \prod_{\alpha\in \Delta_-}\ffg_{\alpha,\Z}\otimes \cF$. Let $\Omega=\{x\in \A\mid u.x=x\}$. Then by Theorem~\ref{t_fltrs_gd_fix}, $\Omega\in \fR$ and thus there exists $\xi\in F^v$ such that $F+\xi+F^v\subset \Omega$. Let $\alpha\in \Delta_-$. Then we have $\omega(\underline{u_\alpha})+\alpha(F+\xi+F^v)\geq 0$. If $\alpha(F^v)=\R_{<0}$, this implies $\underline{u_\alpha}=0$, which proves that $u\in U^{\fR,\infty}(\Delta^m_-(F^v))$. Moreover, if $\alpha\in \Delta^m_-(F^v)$,
 we have $\omega(\underline{u_\alpha})+\alpha(F+\xi+F^v)\geq 0$ and thus $\omega(\underline{u_\alpha})+\alpha(F)\geq 0$. Therefore $u\in U_{F,\infty}^-(\Delta^m_-(F^v))$ and hence $U_{\fR,\infty}^-\subset U^{-}_{F,\infty}(\Delta^m_-(F^v))$. Conversely, $U_{F,\infty}^-(\Delta^m_-(F^v))$ fixes $F+V_{F^v}$ pointwise, 
 where $V_{F^v}$ is the vector space spanned by $F^v$. In particular, it fixes $\fR$ and $\fr$. Therefore $U_{\fR,\infty}^-=U_{\fr,\infty}^-=U^{-}_{F,\infty}(\Delta^m_-(F^v))$. 

Combined with Lemma~\ref{l_U_fR}, we deduce: \[G_\fR=\left(N_{F,J}\cdot U_{F,\infty}^-(\Delta^m_-(F^v))\cdot U_{F,\infty}^+(\Delta^m_+(F^v))\right)\ltimes  U^+_{\infty}(\Delta^u_+(F^v)).\]

By Proposition~\ref{p_Levi_dec_U} and its analogue for $U^-$, we have $ U_{F,\infty}^-(\Delta^m_-(F^v)), U_{F,\infty}^+(\Delta^m_+(F^v))\subset G_{J}$, with the notation of Subsection~\ref{ss_mas_smpl_rt}. Therefore we proved $G_{\fR}\subset G_J\ltimes U^+_{\infty}(\Delta^u_+(F^v))$. By Corollary~\ref{c_act_G_J}, $G_{\fR}\cap G_J=G_{F,J}$ which proves the second part of the proposition. 

We have $G_{\fr}\subset G_{\fR}$. We clearly have  $G_{\fr}\cap U^+_{\infty}(\Delta^u_+(F^v))\subset U_{F,\infty}^+(\Delta^u_+(F^v))$, since an element of $U^+_{\infty}(\Delta^u_+(F^v))$  fixing $\fr$ fixes $F$. Conversely, let $u\in  U_{F,\infty}^+(\Delta^u_+(F^v))$.  Then $u$ fixes $F$ and it fixes $\fR$. As $u$ fixes $\A\cap u.\A$, it fixes $\cl(F,\fR)\supset \fr$, by Lemma~\ref{l_elt_fx_cl}. This proves the proposition, since $M(F,F^v)$ fixes $F+F^v$ and $\cl(F+F^v)\supset \fr$.
\end{proof}

\subsection{Axiom~\ref{a_ma3}}\label{ss_MA3}

In this subsection, we prove axiom~\ref{a_ma3}. We mainly rewrite the proof of \cite{rousseau2011masures}.

So far, we considered groups associated with a vectorial subface of $\overline{C^v_f}$ (like $P(F^v)$, $U^+(F^v)$, etc.). However, when we consider two vectorial faces at the same times, we cannot assume that they are contained in a closed vectorial chamber in general. Let $F^v$ be a positive vectorial face. Write $F^v=w.F^v_1$, with $F^v_1\subset \overline{C^v_f}$ and $w\in W^v$. Then we set $P(F^v)=\tilde{w} P(F_1^v)\tilde{w}^{-1}$, $M(F^v)=\tilde{w} M(F_1^v) \tilde{w}^{-1}$ and $U(F^v)=\tilde{w} U^+(F_1^v)\tilde{w}^{-1}$, where $\tilde{w}$ is defined in \eqref{e_tilde_w}. This is well-defined, independently of the choice of $w$, by \ref{a_KMT7}. Symmetrically, we define $P(F^v)$, $M(F^v)$ and $U(F^v)$ when $F^v$ is a negative vectorial face.

\begin{Lemma}\label{l_parab}
Let $F_1^v,F_2^v$ be two vectorial faces of $\A$. Then $P(F_1^v)\cdot N\cdot P(F_2^v)=G$. 
\end{Lemma}

\begin{proof}
Up to changing the choice of the fundamental chamber $C^v_f$ and by symmetry if $F_1^v$ is negative, we can assume that $F_1^v$ is contained in $\overline{C^v_f}$. Then $P(F_1^v)$ contains $B$. Let $n\in N$ be such that $w.F_2^v\subset \pm \overline{C^v_f}$, where $w=\nu^v(n)$. We have $P(w.F_2^v)\supset B^\epsilon$, where $\epsilon$ is the sign of $F_2^v$.  Let $g\in G$. By the Bruhat or the Birkhoff decomposition, depending on the sign of $F_2^v$ (\eqref{e_Bruhat} and \eqref{e_Birkhoff}), we have $gn^{-1}\in P(F_1^v)N P(w.F_2^v)=P(F_1^v) NnP(F_2^v)n^{-1}$ and hence $g\in P(F_1^v)\cdot N \cdot P(F_2^v)$.
\end{proof}

\begin{Proposition}\label{p_chmny_frndl}(\cite[Proposition 6.7]{rousseau2011masures})
Let $\fr_1=\fr(F_1,F_1^v)$ and $\fr_2=\fr(F_2,F_2^v)$ be two chimneys, with $\fr_1$ splayed, i.e $F_1^v$ spherical. Let $\fR_1$ and $\fR_2$ be their respective germs at infinity. Then we have $G=G_{\fR_1}\cdot N\cdot  G_{\fR_2}$.
\end{Proposition}

\begin{proof}
Let $g\in G$. By Lemma~\ref{l_parab}, there exists $n\in N$ such that $g\in P(F_1^v) n P(F_2^v)=P(F_1^v)\cdot n M(F_2^v)\cdot U(F_2^v)$. The set $\Phi(n^{-1}F_1^v)\cap \Phi^m(F_2^v)$ is a parabolic root subset of $\Phi^m(F_2^v)$. By the Iwasawa decomposition (Proposition~\ref{p_Iwa}) of $M(F_2^v)$, we have: \[M(F_2^v)\subset (M(F_2^v)\cap P(n^{-1}.F_1^v))\cdot (N\cap M(F_2^v))\cdot 
M(F_2,F_2^v).\] Therefore $g\in P(F_1^v)\cdot N \cdot M(F_2,F_2^v)\cdot U(F_2^v)$ and there exist $m\in M(F_1^v)$ and  $n_1\in N$ such that: \[g\in U(F_1^v) mn_1 G_{\fR_2}=U(F_1^v)m G_{n_1.\fR_2} n_1. \]

We have $M(F_1^v)=\langle T,U_\alpha\mid \alpha \in \Phi^m(F^v_1)\rangle$ and thus $M(F_1^v)$ is reductive, since $F_1^v$ is spherical. The set $\Phi(n_1.F^v_2)\cap \Phi^m(F_1^v)$ is a parabolic root subset of $\Phi^m(F_1^v)$. Denote by $M(n_1.F_2^v,F_1^v)$  the Levi subgroup of $M(F_1^v)$ associated with the finite root subsystem  $\Phi^m(n_1.F_2^v)\cap \Phi^m(F_1^v)$. Note that since $\Phi^m(n_1.F_2^v)\cap \Phi^m(F_1^v)$ is finite, the group $M(n_1.F_2^v,F_1^v)$ is reductive. By \eqref{e_Rem_6.2}, we have $P_ {M(F_1^v)}(n_1.F_2^v)=M(n_1.F_2^v,F_1^v)\cdot U_{M(F_1^v)}(n_1.F_2^v)$, where $P_{M(F_1^v)}(n_1.F_2^v)$ is the parabolic subgroup   associated with the vectorial face of $(\A,\Phi^m(F_1^v))$ containing $n_1.F_2^v$ and $U_{M(F_1^v)}(n_1.F_2^v)$ is its unipotent radical. We have $U_{M(F_1^v)}(n_1.F_2^v)\subset U(n_1.F_2^v)$ and thus by the Iwasawa decomposition of $M(F_1^v)$, we have: \[M(F_1^v)\subset M(F_1,F_1^v)\cdot (N\cap M(F_1^v))\cdot M(n_1.F_2^v,F_1^v)\cdot (U(n_1.F_2^v)\cap M(F_1^v)).\] Therefore there exists $n_2\in N\cap M(F_1^v)$ such that: \begin{align*}
m&\in n_2M(n_2^{-1}.F_1,n_2^{-1}.F_1^v)\cdot M(n_1.F_2^v,F_1^v)\cdot (U(n_1.F_2^v)\cap M(F_1^v))\\ &=n_2M(n_2^{-1}.F_1,F_1^v)\cdot M(n_1.F_2^v,F_1^v)\cdot (U(n_1.F_2^v)\cap M(F_1^v)).
\end{align*} 

Let $F_1'$ be the face of $(\A,\Phi^m(n_1.F^v_2)\cap \Phi^m(F_1^v))$ containing $n_2^{-1}.F_1$. Let $M_{F_1'}$ be the fixator of $F_1'$ in $M(n_1.F_2^v,F_1^v)$. We have $M_{F_1'}\subset M(n_2^{-1}.F_1,F_1^v)$. Similarly, let $F_2'$ be the face of $(\A,\Phi^m(n_1.F_2^v)\cap \Phi^m(F_1^v))$ containing $n_1.F_2$. The fixator $M_{F_2'}$ of $F_2'$ in  $M(n_1.F_2^v,F_1^v)$ satisfies $M_{F_2'}\subset M(n_1.F_2,n_1.F_2^v)$. 

By the affine Bruhat decomposition in $M(n_1.F_2^v,F_1^v)$ (\cite[Théorème 7.3.4]{bruhat1972groupes}), we have: \begin{align*}
M(n_1.F_2^v,F_1^v)\subset M_{F_1'}\cdot (N\cap M(n_1.F_2^v,F_1^v))\cdot M_{F_2'}\subset M(n_2^{-1}.F_1,F_1^v)\cdot N\cdot M(n_1.F_2,n_1.F_2^v).
\end{align*}

Therefore: \begin{align*}
m&\in n_2M(n_2^{-1}.F_1,F_1^v)\cdot N \cdot M(n_1.F_2,n_1.F_2^v)\cdot U(n_1.F_2^v)\\
&=M(F_1,F_1^v)\cdot N\cdot G_{n_1.\fR_2} \text{ by Proposition~\ref{p_fx_chimn}}
\end{align*}

and \[g\in U(F_1^v)\cdot M(F_1,F_1^v)\cdot N\cdot G_{n_1.\fR_2}n_1=G_{\fR_1}\cdot N\cdot G_{\fR_2},\] which completes the proof of the  proposition.
\end{proof}
 
\begin{Corollary}\label{c_MA3}
The masure $\I$ of $G$ satisfies axiom~\ref{a_ma3}: if $\fR_1$ is a splayed chimney germ and if $\fR_2$ is a chimney germ or a local face, then there exists an apartment containing $\fR_1$ and $\fR_2$. 
\end{Corollary}

\begin{proof}
If $\fR_2$ is a chimney germ, then this is a consequence of Proposition~\ref{p_chmny_frndl} and of Lemma~\ref{l_frndly_prs}. Assume now that $\fR_2=F$ is  a face. Up to translating by an element of $G$, we can assume that $F$ is contained in $\A$.  Let $\fr_2=\fr(F,\A_{\ines})$. Then $\fr_2$ is a chimney. Thus there exists an apartment $A$ containing $\fR_1$ and the germ $\fR_2'$ of $\fr_2$ at infinity. Then $A$ contains $\fR_1$ and $F+\xi$, for some $\xi\in \A_{\ines}$. By~\ref{a_ma2}, $A\cap \A$ is a finite intersection of half-apartments and thus $A\cap \A$ contains $F+\xi+\A_{\ines}=F+\A_{\ines}$, which proves that $A$ contains $F$ and $\fR_1$. 
\end{proof}

We already know that $\I$ satisfies~\ref{a_ma2}, by Corollary~\ref{c_w_mas_gen}. As it satisfies~\ref{a_ma3}, we deduce:

\begin{Theorem}\label{t_I_abst_mas}
Let $\cS=(A,X,Y,(\alpha_i)_{i\in I},(\alpha_i^\vee)_{i\in I})$ be a free Kac--Moody datum in the sense of Subsection~\ref{ss_KM_data} and $\fG=\fG_{\cS}:\{\mathrm{Fields}\}\rightarrow \{\mathrm{Groups}\}$ be the constructive Kac--Moody group functor associated to $\cS$ in Subsection~\ref{ss_Tits_minKM}. Let $\cF$ be  a field equipped with a non-trivial valuation $\omega:\cF\rightarrow \R\cup \{+\infty\}$. For $x\in \A$, define $G_x=\langle U_{x,\infty}^+,U_{x,\infty}^-,N_x\rangle$, where $U_{x,\infty}^{\pm}$ is defined in Definition~\ref{d_U_infty_pm} and $N_x$ is defined as the fixator of $x$ in $N$ for the action $\nu$ of $N$ on $\A$ defined in Section~\ref{s_act_N_A}.
Then the masure $\I$ defined in Subsection~\ref{ss_d_mas} is an abstract masure in the sense of Definition~\ref{d_w_mas}.
\end{Theorem}

\section{Towards a definition of a masure associated with maximal Kac--Moody group}

It seems natural to define a masure associated $\widehat{\I}$ with $\fG^{ma+}(\cF)$ instead of $G$ for at least two reasons. Firstly because it seems to be a natural tool to study $\widehat{G}:=\fG^{ma+}(\cF)$ and secondly because the masure $\I$ of $G$ is defined via the use of $\widehat{G}$. Thus it would be natural to define $\widehat{\I}$ and then to define $\I$ as the subset $G.\A$. It would be interesting to define such a masure as $(\widehat{G}\times \A)/\sim$ for some equivalence relation $\sim$, which would be equipped with an action of the entire $\widehat{G}$. Here we simply define an object $\widehat{\I}$ equipped with an action of $\widehat{U}=\fU^{ma+}(\cF)$.

 First note that as a set equipped with an action of $U^+$, the definition of $\I$ is easier, by the lemma below.

\begin{Lemma}
Define $\sim$ on $U\times \A$ by $(u,x)=(v,x)$ if and only if $v^{-1}u\in U_{x,\infty}^+$. Then the map $\iota:(U^+\times \A)/\sim \rightarrow \I$ defined by $[u,x]\rightarrow u.x$ is well-defined and is a bijection.
\end{Lemma}

\begin{proof}
Let $u,v\in U^+$ and $x\in \A$ be such that $u.x=v.x$. Then $u^{-1}v\in G_x\cap U^+=U_{x,\infty}^+$ (by Theorem~\ref{t_fltrs_gd_fix}) and thus $(u,x)\sim (v,x)$. Therefore $\iota$ is well-defined. By the Iwasawa decomposition (Proposition~\ref{p_Iwa}), we have $\I=\bigcup_{u\in U^+} u.\A$, thus $\iota$ is surjective. Let $[u,x],[v,y]\in (U^+\times \A)/\sim$ be such that $u.x=v.y$. Then  $\rho_{+\infty}(u.x)=x=\rho_{+\infty}(v.y)=y$ and thus $x=y$. Now $u^{-1}v\in G_x\cap U^+=U_{x,\infty}^+$ (by Theorem~\ref{t_fltrs_gd_fix}) and thus $[u,x]=[v,x]$.
\end{proof}

For $x\in \A$, set $\widehat{U}_x=\fU^{ma+}_x(\cF)\subset \widehat{U}$. Let $\sim$ be the equivalence relation on $\widehat{U}\times \A$ defined by $(u,x)\sim (v,x)$ if and only if $u^{-1}v\in \widehat{U}_x$, for $u,v\in \widehat{U}$ and $x\in \A$.  We define $\widehat{\I}$ as $(\widehat{U}\times \A)/\sim$. Then $\A$ naturally embeds in $\widehat{\I}$ via the map $x\mapsto [1,x]$, for $x \in \A$, where $[u,x]$ is the class of $(u,x)$ for $u\in \widehat{U}$. We now write $u.x$ instead of $[u,x]$.  Then the natural map $\I\rightarrow \widehat{\I}$ defined by $u.x\mapsto u.x$ is well-defined and is an embedding since $U_{x,\infty}^+\subset \widehat{U}_x$, for every $x\in \A$.

If $\Omega\subset \A$ is a non-empty set, then its fixator in $\widehat{U}$ is $\widehat{U}_\Omega:=\bigcap_{x\in \Omega}\widehat{U}_x=\fU^{ma+}_{\Omega}(\cF)$ and if $\cV\subset \A$ is a filter, its fixator is $\widehat{U}_\cV:=\bigcup_{\Omega\in \cV} \widehat{U}_\Omega$. Therefore $U_{\cV,\infty}^+=\widehat{U}_\cV\cap U^+$.

Note however that we can have $u.\A\cap \A=\emptyset$, for $u\in \widehat{U}$. Indeed, take $\lambda\in C^v_f\cap Y$ and take $(\alpha_n)\in (\Phi_+)^\N$  such that $((\htt(\alpha_n))$ is a strictly increasing sequence. For $n\in \N$, choose $a_n\in \cF$ such that $\omega(a_n)=-2n\alpha_n(\lambda)$. Set $u=\prod_{n\in \N} x_{\alpha_n}(a_n)$. Then $u\notin \widehat{U}_{n\lambda}$, for every $n\in \N$ and thus $u.(n\lambda)\notin \A$, for every $n\in \N$. If $u.\A$ contained an element $x$ of $\A$, then it would fix $x+C^v_f$ and thus if would fix $n\lambda$, for $n\gg 0$, which is impossible.

\chapter{Parahoric subgroups and minimal group}\label{C_para_min_gp}

In this chapter, we give a description of $G_F$ entirely in terms of the minimal group $G$, for a spherical local $F$ (a notion that we describe below), generalizing the case of an alcove, that we already treated (see Proposition~\ref{p_bphr_2.4}). We also give a description of the fixator of $0$ and of the Iwahori subgroup involving
the “minimal group” $\fG(\cO)$.

\section{Description of the parahoric subgroup associated with a spherical face}

We call a face $F\subset \A$ \textbf{spherical} if its direction is, i.e if $F=\germ_x(x+F^v)$, for some spherical face $F^v$. We now prove that if $F$ is  a point or a spherical local face, then $U_{F,\infty}^\epsilon=U_{F,\fin}^\epsilon$, for $\epsilon\in \{-,+\}$.   To do that, we prove that $U_{F,\infty}^+=U_{C_+,\infty}^+$ and $U_{F,\infty}^-=U_{C_-,\infty}^-$, for some well-chosen alcoves $C_+$ and $C_-$ and then we use the case of alcoves.

Recall that for $J\subset I$, we set $\Phi_J=(\bigoplus_{i\in J} \R \alpha_i)\cap \Phi$ and that  $\A_J$ is the set $\A$ equipped with the hyperplane arrangement $\{\alpha^{-1}(\{k\})\mid \alpha\in \Phi_J, k\in \Lambda\}$.

\begin{Lemma}\label{l_w0_J}
Let $F^v$ be a spherical subface of $\overline{C^v_f}$. Let $V_{F^v}=\mathrm{vect}(F^v)\subset \A$. Let $J=\{i\in I\mid \alpha_i(F^v)=\{0\}\}$. Let $w_{0,J}$ be the longest element of $W^v_J$ (defined in \cite[Proposition 2.3.1]{bjorner2005combinatorics}). Then $F^v-\overline{C^v_f}=w_{0,J}.\overline{C^v_f}-\overline{C^v_f}$. 
\end{Lemma}

\begin{proof}
Consider the apartment $\A_J$ defined in Subsection~\ref{ss_mas_smpl_rt}. Let $C^v_{f,J}$ be the vectorial chamber of $\A_J$ containing $C^v_f$. Let us prove: \begin{equation}\label{e_chmbr} 
C^v_{f,J}=C^v_f+V_{F^v}.
\end{equation} 
For $x\in C^v_f+V_{F^v}$, we have $\alpha_j(x)>0$ for all $j\in J$ and thus $C^v_f+V_{F^v}\subset C^v_{f,J}$. Let $x\in C^v_{f,J}$.  Let $\xi\in F^v$. Then $\alpha_j(x)>0$ for all $j\in J$. For $i\in I\setminus J$, we have $\alpha_i(\xi)>0$ and thus $x+t\xi\in C^v_f$ for $t\in \R_{>0}$ large enough. Therefore $x\in C^v_f+V_{F^v}$ and we have \eqref{e_chmbr}.

By definition, we have $\ell(w_{0,J} r_i)<\ell(w_{0,J})$, for all $i\in J$ and thus by \cite[1.3.13 Lemma]{kumar2002kac}, we have $w_{0,J}.\alpha_i\in \Phi_{-,J}$. Consequently, $w_{0,J}.C^v_{f,J}=-C^v_{f,J}$ and thus $w_{0,J}.\overline{C^v_f}+V_{F^v}=-\overline{C^v_f}+V_{F^v}$. We have $V_{F^v}=F^v-F^v$ and thus:
 \begin{align*}
 F^v-\overline{C^v_f}&=F^v-\overline{C^v_f}+V_{F^v}-\overline{C^v_f}\\
 &= w_{0,J}.\overline{C^v_f} +V_{F^v}+F^v-\overline{C^v_f}\\
 &=w_{0,J}.\overline{C^v_f}+V_{F^v}-\overline{C^v_f}\\
 &=w_{0,J}.\overline{C^v_f}+F^v-F^v-\overline{C^v_f}\\
 &= w_{0,J}.\overline{C^v_f}-\overline{C^v_f}. 
 \end{align*}
\end{proof}

\begin{Lemma}\label{l_U_F_-}

Let $x\in \A$ and $F^v$ be a spherical subface of $\overline{C^v_f}$. Let $F=\germ_x(x+F^v)$. Then $U_{F,\infty}^-=U_{F,\fin}^-$.
\end{Lemma}

\begin{proof}
 Let $C=\germ_x(x+w_{0,J}.C^v_f)$, with the notation of Lemma~\ref{l_w0_J}. Let us prove $U_{F,\infty}^-=U_{C,\infty}^-$.

We have $F\subset \overline{C}$ and thus $U_{F,\infty}^-\supset U_{\overline{C},\infty}^-=U_{C,\infty}^-$. Let $u\in U_{F,\infty}^-$ and $A=u.\A$.   Then by~\ref{a_ma2}, there exist $k\in \N$, $\beta_1,\ldots,\beta_k\in \Phi$ and $\lambda_1,\ldots,\lambda_k\in \R$ such that $\Omega:=A\cap \A=\bigcap_{i=1}^k D(\beta_i,\lambda_i)$. As $\Omega$ contains $\fQ_{-\infty}$, we deduce that $\beta_i\in \Phi_-$, for all $i\in \llbracket 1,k\rrbracket$. 

Let $i\in \llbracket 1,k\rrbracket$. If $\beta_i(F^v)<0$, then $\beta_i(x)+\lambda_i>0$, thus $x\in \mathring{D}(\beta_i,\lambda_i)$  and hence $C\subset D(\beta_i,\lambda_i)$.

 Assume now $\beta_i(F^v)=0$. We have $\beta_i(F)=\beta_i(x+V_{F^v})$, where $V_{F^v}=\mathrm{vect}(F^v)$. As $\beta_i\in \Phi_-$, we have:  \[D(\beta_i,\lambda_i)\supset x+V_{F^v}-\overline{C^v_f}\supset x+F^v-\overline{C^v_f}=x+w_{0,J}.\overline{C^v_f}-\overline{C^v_f}\supset C-\overline{C^v_f}\supset C,\] by Lemma~\ref{l_w0_J}. Consequently, $C$ is contained in $\Omega$.  By Lemma~\ref{l_elt_fx_cl}, $u$ fixes $\Omega$, thus $u\in G_C\cap U^-=U_{C,\infty}^-$, by Theorem~\ref{t_fltrs_gd_fix}. Consequently, $U_{F,\infty}^-\subset U_{C,\infty}^-$ and thus $U_{F,\infty}^-= U_{C,\infty}^-$. We conclude with lemma~\ref{l_U_F_fin}.
\end{proof}

The following proposition generalizes partially Proposition~\ref{p_bphr_2.4} to the case of spherical local faces.

\begin{Proposition}\label{p_sphrcl_para}
Let $F$ be a spherical local face of $\A$. Then: \[G_F=U_{F,\fin}^+\cdot U_{F,\fin}^-\cdot N_F=\langle U_{\alpha,F}\mid \alpha\in \Phi\rangle \cdot N_F,\ U_{F,\fin}^-=U_{F,\infty}^-\text{ and }U_{F,\fin}^+=U_{F,\infty}^+.\] 
\end{Proposition}

\begin{proof}
By symmetry, we can assume that $F$ is positive. Let $F^v$ be a spherical  vectorial face and $x\in \A$ be such that $F=\germ_x(x+F^v)$. Let $w\in W^v$ be such that $F^v\subset w.\overline{C^v_f}$. By Theorem~\ref{t_fltrs_gd_fix}, we have $G_{F}=U_{F,\infty}(w.\Delta_+)\cdot U_{F,\infty}(w.\Delta_-)\cdot N_{F}$. Thus we can assume that $w=1$. We have $U_{F,\infty}^+=U_{x,\infty}^+=U_{\germ_x(x+C^v_f),\infty}^+$ (by Corollary~\ref{c_Uma_cl}) 	and thus by Lemma~\ref{l_U_F_fin}, we have $U_{F,\infty}^+=U_{F,\fin}^+$. By Lemma~\ref{l_U_F_-}, we have $U_{F,\infty}^-=U_{F,\fin}^-$. We conclude with Theorem~\ref{t_fltrs_gd_fix}.
\end{proof}

\section{Interpretation of parahoric subgroups in term of $\fG^{\min}(\cO)$}

In this section, we obtain a description of the fixator of $0$ and of the fixator of the fundamental alcove in terms minimal group $\fG^{\min}(\cO)$ defined in Definition~\ref{d_G_min}. This generalizes well-known results in the case of reductive groups.

\subsection{Fixator of $0$}

Let $\cO$ be the ring of integers of $(\cF,\omega)$, $\fm$ be its maximal ideal and $\kk=\cO/\fm$ be its residue field.

\begin{Proposition}\label{p_fxt_zr}(see \cite[Proposition 3.1]{hebert2025topologies})
\begin{enumerate}
\item The fixator $G_0$ of $0$ in $G$ is the  group $\fG^{\min}(\cO)$.

\item We have $N_0=N\cap \fG^{\min}(\cO)=\fN(\cO)$, $U_{0,\infty}^+=U^+\cap \fG^{\min}(\cO)=U_{0,\fin}^+$ and $U_{0,\infty}^-=U^-\cap \fG^{\min}(\cO)=U_{0,\fin}^-$. Moreover, $N_0$ is the normalizer of $T_0$ in $G_0$.

\end{enumerate}
\end{Proposition}

\begin{proof}
(1) For $i\in I$, $x_{\alpha_i}(\cO)$, $x_{-\alpha_i}(\cO)$ and $T_0=\fT(\cO)$ fix $0$. Therefore by \eqref{e_minimal_group_semilocal_ring},  $\fG^{\min}(\cO)\subset G_0$. 

By Theorem~\ref{t_G_x}, we have: \begin{equation}\label{e_G0}
G_0=U_{0,\infty}^+ \cdot U_{0,\infty}^-\cdot  N_0,
\end{equation} where $N_0=\{n\in N\mid n.0=0\}$.  By Proposition~\ref{p_sphrcl_para}, we have:   \begin{equation}\label{e_U0pm+}
U_{0,\infty}^+=U_{0,\fin}^+:=\langle x_\alpha(u)\mid \alpha\in \Phi, u\in \cO\rangle\cap U^+\subset \fG^{\min}(\cO)
\end{equation}
and 
\begin{equation}\label{e_U0nm-}
U_{0,\infty}^-=U_{0,\fin}^-:=\langle x_\alpha(u)\mid \alpha\in \Phi, u\in \cO\rangle\cap U^-\subset \fG^{\min}(\cO).
\end{equation} For $i\in I$,  $\tilde{s}_{i}=x_{\alpha_i}(1) x_{-\alpha_i}(-1)x_{\alpha_i}(1)\in \fG^{\min}(\cO)$.  We have  $N=\langle \fT(\cF), \tilde{s}_{i}\mid i \in I\rangle$. Let $n\in N_0$. Write $\nu(n)=w\in W^v$. Write $w=r_{i_1}\ldots r_{i_k}$, with $k=\ell(w)$ and $i_1,\ldots,i_k\in I$. Let $n'=\tilde{s}_{i_1}\ldots \tilde{s}_{i_k}\in N_0$. By definition of $\nu$ (see Proposition~\ref{p_d_nu}), $\nu(n')=w$ and  $t:=n'^{-1}n\in T\cap \ker(\nu)$. By \eqref{e_act_T}, $\ker(\nu)\cap T=T_0$. Therefore:
 \begin{equation}\label{e_N0}
N_0=\langle \tilde{s}_{i}\mid i\in I\rangle \cdot T_0.
\end{equation} and in particular, $N_0\subset \fG^{\min}(\cO)$. This proves 1). 

(2) The part before ``Moreover'' of (2) follows from Theorem~\ref{t_fltrs_gd_fix} and Proposition~\ref{p_sphrcl_para}.
Let $N'$ be the normalizer of $T_0$ in $G_0$. Let $n\in N_0$. Then $nT_0n^{-1}\subset T$, since $N_0$  is contained in the normalizer $N$ of $T$ in $G$. Moreover $nT_0 n^{-1}$ fixes $0$ and thus $nT_0 n^{-1}\subset T_0$, which proves that $N_0\subset N'$. 

Let $g\in N'\subset G_0$. Write $g=u_+ u_- n$, with $u_+\in U_{0,\infty}^+$, $u_-\in U_{0,\infty}^-$ and $n\in N_0$. Then replacing $g$ by $gn^{-1}$ if necessary, we can assume that $n=1$. Let $t\in T_0$. Then by assumption $t':=gtg^{-1}\in T_0$. We have $gtg^{-1}=t'=u_+ u_- t u_-^{-1} u_+^{-1}=t'$, thus by \ref{a_KMT4'}, we have:
\[u_-tu_-^{-1}=u_+^{-1}t'u_+=u_- \underbrace{(tu_-^{-1} t^{-1})}_{\in U^-}t=u_+^{-1}\underbrace{(t' u_+ t'^{-1})}_{\in U^+} t'\]

Using Lemma~\ref{l_utu}, we get that  $t=t'$  and that  \[ u_-t=tu_-\text{ and }u_+t=tu_+,\] for all $t\in T_0$. 

  Write $u_+=\prod_{\alpha\in \Delta_+} X_\alpha(\underline{u_\alpha})$, with $(\underline{u_\alpha})\in \prod_{\alpha\in \Delta_+} \ffg_{\alpha,\Z}\otimes \cO$. Then by \eqref{e_com_T_X}, $tu_+ t^{-1}=\prod_{\alpha\in \Delta_+}X_\alpha(\alpha(t)\underline{u_\alpha})$. As the ``coordinates '' $X_\alpha$ are unique (\eqref{e_normal_form_Upma}), we deduce that $\alpha(t)\underline{u_\alpha}=\underline{u_\alpha}$, for all $\alpha\in \Delta_+$. Let $a\in \cO^\times$ be such that $a^2\neq 1$ (for example we can take $a=1+x$, for some $x\in \cF$ such that $\omega(x)>0$). Let $t=\alpha^\vee(a)\in T_0$. Then $\alpha(t)=a^2\neq 1$ and hence $\underline{u_\alpha}=0$. This proves that $u_+=1$. Similarly, $u_-=1$ and hence $g\in N_0$. Thus we proved $N'\subset N_0$ and then $N'=N_0$, which proves (2). 
\end{proof}

\subsection{Congruence subgroups}

Let $\fm=\{x\in \cO\mid \omega(x)>0\}$\index[notation]{m@$\fm$} be the maximal ideal of $\cO$ and $\kk=\cO/\fm$\index[notation]{k@$\kk$}. Let $\pr_{\kk}:\fG^{\min}(\cO)\rightarrow \fG^{\min}(\kk)$\index{p@$\pr_{\kk}$} be the morphism induced by the natural projection $\cO\twoheadrightarrow \kk$. 
In this subsection we study $\ker \pr_{\kk}$, following \cite[3.2]{hebert2025topologies}: we prove that it decomposes as the product of its intersections with $U^-$, $U^+$ and $T$ (see Proposition~\ref{p_ker_pi_k}), using the masure $\I$ of $G$. We also describe  $U^-\cap \ker \pr_{\kk}$ and $U^+\cap \ker \pr_{\kk}$ through their actions on $\I$ and we deduce that $\ker \pr_{\kk}$ fixes $C_0^+\cup C_0^-$.  

Let $C,C'$ be two positive alcoves based at $0$. By Proposition~\ref{p_frndly_fac}, there exists an apartment $A$ containing $C$ and $C'$. Let $g\in G$ be such that  $g.A=\A$ and $g.C=C_0^+$. Then by Corollary~\ref{c_BT_7.4.8}, there exists $w\in \widetilde{W^a}$ such that $g.x=w.x$, for all $x\in \A\cap g^{-1}.\A$. As $w.0=0$, we have $w\in W^v$.  Then $g.C'=w.C'$ and we set $\dw(C,C')=w$\index[notation]{d@$\dw$}, which is well-defined, independently of the choices we made, by~\ref{a_ma2}. Then $\dw$ is $G$-invariant. 
 Note that this ``distance'' is a restriction of the ``$W$-distance'' $d^{W}$ of \cite[1.11]{bardy2016iwahori}.

\begin{Lemma}\label{l_fiber_retraction}
Let $C$ be a positive alcove of $\I$ based at $0$ and $w\in W^v$.  Write $w=r_{i_1}\ldots r_{i_k}$, with $k=\ell(w)$ and $i_1,\ldots,i_k\in I$. Let $\beta_{1}=\alpha_{i_1}$, $\beta_2=r_{i_1}.\alpha_{i_2}$, $\ldots$, $\beta_k=r_{i_1}\ldots r_{i_{k-1}}.\alpha_{i_k}$. Then  we have $\rho_{+\infty}(C)=w.C_0^+$ if and only if  there exists $a_1,\ldots,a_k\in \cO$ such that $C=x_{\beta_1}(a_1)\ldots x_{\beta_k}(a_k).w.C_0^+$. In other words, \[\rho_{+\infty}^{-1}(\{w.C_0^+\})=\fU_{\Inv(w^{-1})}(\cO).w.C_0^+,\] for the notation of \eqref{e_U_Inv_w}.
\end{Lemma}

\begin{proof}
As $x_{\beta_1}(\cO).\ldots.x_{\beta_k}(\cO)$ fix $0$, an element of $x_{\beta_1}(\cO)\ldots x_{\beta_k}(\cO).w.C_0^+$ is a positive alcove based at $0$. By \eqref{e_Kmr_1.3.14}, $\beta_1,\ldots,\beta_k\in \Phi_+$. Thus $x_{\beta_1}(\cO).\ldots.x_{\beta_k}(\cO).w.C_0^+\subset U^+.w.C_0^+$  and we have one implication. 

 We prove the reciprocal by induction on $\ell(w)$. Assume $w=1$. Then $\rho_{+\infty}(C)=C_0^+$. Let $A$ be an apartment containing $\fQ_{+\infty}$ and $C$. Let $g\in G$ be such that $g.A=\A$ and $g$ fixes $A\cap \A$. Then $C=g^{-1}.C_0^+$, by definition of $\rho_{+\infty}$. Moreover, $A$ contains $0$ and $\fQ_{+\infty}$ and thus it contains $\conv(0,\fQ_{+\infty})\supset C_0^+$. Therefore $C=C_0^+$ and the lemma is clear in this case. Assume now that $k:=\ell(w)\geq 1$. 

Let $C_0'=C_0^+$, $C_1'=r_{i_1}.C_0^+$, $\ldots$, $C_k'=r_{i_1}\ldots r_{i_k}.C_0^+=C'$. Let $C$ be a positive alcove based at $0$ and such that $\rho_{+\infty}(C)=C_k'$. Let $A$ be an apartment containing $C$ and $\fQ_{+\infty}$. Let $g\in G$ be such that $g.A=\A$ and $g$ fixes $A\cap \A$.  Set $C_i=g^{-1}.C_i'$, for $i\in \llbracket 0,k\rrbracket$. Then $g$ fixes $\fQ_{+\infty}$ and hence $\rho_{+\infty}(g.x)=x$ for every $x\in \A$. Therefore $\rho_{+\infty}(C_i)=C_i'$, for $i\in \llbracket 0,k\rrbracket$. In particular, $\rho_{+\infty}(C_{k-1})=C_{k-1}'$. By induction, we may assume that there exist  $a_1,\ldots,a_{k-1}\in \cO$ such that  $C_{k-1}=uv.C_0^+$, where $u=x_{\beta_{1}}(a_1)\ldots x_{\beta_k}(a_{k-1})$ and $v=r_{i_1}\ldots r_{i_{k-1}}.$  Moreover we have: \[\dw(C_{k-1},C_k)=\dw(C_{k-1}',C_k')=r_{i_k}=\dw\left(u^{-1}.C_{k-1},u^{-1}.C_k\right)=\dw\left(v.C_0^+,u^{-1}.C_k\right).\] 

Let $P$ be the panel common to $v.C_0^+$ and $u^{-1}.C_k$. Then $P\subset \beta_k^{-1}(\{0\})$. Let $D$ be the half-apartment delimited by $\beta_{k}^{-1}(\{0\})$ and containing $C_{k-1}$. Then as \[\beta_{k}(C_{k-1})=r_{i_1}\ldots r_{i_{k-1}}.\alpha_{i_k}(r_{i_1}\ldots r_{i_{k-1}}.C_0^+)>0,\] $D$ contains $\fQ_{+\infty}$. By Proposition~\ref{p_loc_SC}, there exists an apartment $B$ containing $D$ and $u^{-1}.C_k$. Let $g'\in G$ be such that $g'.B=\A$ and $g'$ fixes $\A\cap B$. We have $g'.u^{-1}.C_k=\rho_{+\infty}(C_k)=C_k'$. By Theorem~\ref{t_fltrs_gd_fix}, $g'\in  U_{\beta_k,0}\cdot T_0$ and as $T_0$ fixes $\A$, we can assume that $g'\in U_{\beta_k,0}=x_{\beta_k}(\cO)$. Write $g'=x_{\beta_k}(-a_k)$, with $a_k\in \cO$. Then $C_k=u.x_{\beta_k}(a_k).C_k=x_{\beta_1}(a_1)\ldots x_{\beta_k}(a_k).w.C_0^+$, which proves that $\rho_{+\infty}^{-1}(\{w.C_0^+\})=x_{\beta_1}(\cO)\ldots x_{\beta_k}(\cO).w.C_0^+$. We conclude with Corollary~\ref{c_dec_inv}.
\end{proof}

\begin{Proposition}\label{p_ker_pi_k}
\begin{enumerate}
\item We have $U_{C_0^+,\infty}^-=U^-\cap \ker(\pr_{\kk})$ and $U_{-C_0^+,\infty}^+=U^+\cap \ker(\pr_{\kk})$.

\item We have $\ker(\pr_{\kk})=\left(\ker(\pr_{\kk})\cap U^+\right).\left(\ker(\pr_{\kk})\cap U^-\right).\left(\ker(\pr_{\kk})\cap \fT(\cO)\right).$

\item We have $\ker(\pr_{\kk})\subset G_{C_0^+\cup C_0^-}$. 
\end{enumerate}  
\end{Proposition}

\begin{proof}
2) Let $u\in U_{C_0^+,\infty}^-$.   By definition, there exists $\Omega\in C_0^+$ such that $u\in U_{\Omega,\infty}^{-}$. Let $x\in C^v_f\cap \Omega$. Then   $u\in \prod_{\alpha\in \Delta_-} X_\alpha(\ffg_{\alpha,\Z}\otimes \cF_{\omega\geq -\alpha(x)})\cap G_0$. As $-\alpha(x)>0$ for every $\alpha\in \Delta_-$, Proposition~\ref{p_mrphsm} implies: 
\begin{equation}\label{e_relation_kerpi_UCO}
  U_{C_0^+,\infty}^-\subset \ker(\pr_{\kk}).
  \end{equation}

 Let $g\in \ker(\pr_{\kk})\subset \fG^{\min}(\cO)$. Then $g$ fixes $0$ and $g.C_0^+$ is a positive alcove based at $0$. Write $\rho_{+\infty}(g.C_0^+)=w.C_0^+$, with $w\in W^v$. Write $w=r_{i_1}\ldots r_{i_m}$, with $m=\ell(w)$ and $i_1,\ldots,i_m\in I$. Let $\tilde{n}=\tilde{r}_{i_1}\ldots \tilde{r}_{i_m}\in \fN(\cF)$. By Lemma~\ref{l_fiber_retraction}, there exists $u\in \langle x_\beta(\cO)\mid \beta\in \Phi_+\rangle\subset U_{0,\infty}^+$ such that $g.C_0^+=u\tilde{n}.C_0^+$. Then $g=u\tilde{n}i$, with $i\in G_{C_0^+}$. By Theorem~\ref{t_fltrs_gd_fix}, we have: \[G_{C_0^+}=U_{C_0^+,\infty}^+ \cdot  U_{C_0^+,\infty}^- \cdot N_{C_0^+}.\] As every element of $C_0^+$ has non empty interior, we have $N_{C_0^+}=T_0$.  Moreover, $U_{C_0^+,\infty}^+=U_{0,\infty}^+$ and $T_0$ normalizes $U_{0,\infty}^+$ and $U_{C_0^+,\infty}^-$, by Proposition~\ref{p_cnj_N}. Therefore, \[G_{C_0^+}=T_0\cdot U_{0,\infty}^+\cdot U_{C_0^+,\infty}^-.\] Write $i=tu_+ u_-$, with $t\in T_0$, $u_+\in U_{0,\infty}^+$ and $u_-\in U_{C_0^+,\infty}^-$.

According to  \eqref{e_relation_kerpi_UCO}, we have: \[\pr_{\kk}(g)=1=\pr_{\kk}(u\tilde{n}tu_+u_-)=\pr_{\kk}(u)\pr_{\kk}(\tilde{n}t)\pr_{\kk}(u_+)\pr_{\kk}(u_-)=\pr_{\kk}(u)\pr_{\kk}(\tilde{n}t)\pr_{\kk}(u_+) .\]
 
 By \cite[3.16 Proposition]{rousseau2016groupes}  or  \cite[Theorem 8.118]{marquis2018introduction}, 
\[(\fG^{ma+}(\kk),\fN(\kk),\fU^{ma+}(\kk),\fU^{-}(\kk),\fT(\kk),\{r_i\mid i\in I\})\] is a refined Tits system. By \cite[3.16 Remarque]{rousseau2016groupes}, we have the Bruhat decomposition: \[\fG^{ma+}(\kk)=\bigsqcup_{n\in \fN(\kk)} \fU^{+}(\kk) n \fU^{ma+}(\kk).\]  As $\pr_{\kk}(u)\in \fU^+(\kk)$, $\pr_{\kk}(\tilde{n}t)\in \fN(\kk)$ and $\pr_{\kk}(u_+)\in \fU^+(\kk)\subset  \fU^{ma+}(\kk)$, we deduce $\pr_{\kk}(\tilde{n}t)=1$.

By Proposition~\ref{p_d_nu_v}, there exists a group morphism $\nu^v_\kk:\fN(\kk)\rightarrow W^v$ such that $\nu^v_\kk(\tilde{s}_{i})=r_i$ for $i\in I$ and $\nu^v_\kk(\fT(\kk))=1$. Then $\nu_\kk^v(\tilde{n}t)=w=1$. Therefore $w=1$, $\tilde{n}=1$ and $g=utu_+u_-=u'tu_-$, for some $u'\in U_{0,\infty}^+$, since $t$  normalizes $U_{0,\infty}^+$. By Lemma~\ref{l_utu} and by symmetry of the roles of $U^{-}_{\cdot,\infty}$ and $U^{+}_{\cdot ,\infty}$, we have $u'\in U_{-C_0^+,\infty}^+$. By \eqref{e_relation_kerpi_UCO} applied to $U_{C_0^+,\infty}^{-}$ and to  $U_{-C_0^+,\infty}^+$, we have $\pr_{\kk}(g)=1=\pr_{\kk}(u')\pr_{\kk}(t)\pr_{\kk}(u_-)=\pr_{\kk}(t)$ and thus: \[g\in U_{-C_0^-,\infty}^+\cdot  (T_0\cap \ker (\pr_{\kk}))\cdot  U_{C_0^+,\infty}^-= U_{-C_0^-,\infty}^+\cdot U_{C_0^+,\infty}^-\cdot (T_0\cap \ker (\pr_{\kk})).\]  
By \eqref{e_relation_kerpi_UCO}, we deduce 2). 

3) We have $U_{-C_0^+,\infty}^+=U_{-C_0^+ + C^v_f,\infty}^+\subset G_{C_0^+\cup -C_0^+}$, $U_{C_0^+,\infty}^-=U_{C_0^+ -C^v_f,\infty}^{-}\subset G_{C_0^+\cup-C_0^+}$ and $T\cap \ker (\pr_{\kk}) \subset T_0\cap G_0\subset G_\A$, which proves 3). 

1) We already proved one inclusion. Let $u\in \ker (\pr_{\kk})\cap U^-$. Then by what we proved above, $u\in U_{-C_0^-,\infty}^+. (T_0\cap \ker (\pr_{\kk})). U_{C_0^+,\infty}^-$. By Lemma~\ref{l_utu}, we have $u\in U_{C_0^+,\infty}^-$, and the proposition follows. 
\end{proof}

\begin{Remark}
Let $u\in U^-\cap \ker \pr_{\kk}=U_{-C_0^+,\infty}^{-}$. Write $u=\prod_{\alpha\in \Delta_-}X_\alpha(\underline{v_\alpha})$, where $\underline{v_\alpha}\in \ffg_{\alpha,\Z}\otimes \cO$, for every $\alpha\in \Delta_-$.  Let $\lambda\in Y$ be such that $\alpha_i(\lambda)=1$ for every $i\in I$. Let $\Omega\in C_0^+$ be such that $u\in U_{\Omega,\infty}^-$ and  $\eta\in \R^*_+$ be such $\eta\lambda\in \Omega$. Then $u\in U_{\Omega,\infty}^-$ implies $\omega(v_\alpha)\geq |\alpha(\eta \lambda)|=\eta \htt(\alpha)$ for every $\alpha\in \Delta_-$. In particular, $\omega(\underline{v_\alpha})$ goes to $+\infty$ when $-\htt(\alpha)$ goes to $+\infty$.
\end{Remark}

\subsection{The Iwahori subgroup}

The \textbf{Iwahori subgroup}\index{Iwahori subgroup} is the fixator in $G$ of the positive fundamental alcove $C_0^+$. By \cite[Proposition 8.128]{marquis2018introduction}, $\fG$ is a non-degenerate Tits functor. In particular, the natural projection $\pr_\kk:\cO\rightarrow \kk$ induces a projection $\pr_\kk:\fG^{\min}(\cO)\rightarrow \fG(\kk)$. We now obtain the usual description of the Iwahori subgroup in the Kac--Moody frameworks.
 
\begin{Proposition}\label{p_Iwahori}
We have $G_{C_0^+}=\pr^{-1}_{\kk}(\fB(\kk))$. 
\end{Proposition}

\begin{proof}
By Theorem~\ref{t_fltrs_gd_fix}, we have $G_{C_0^+}=U_{C_0^+,\infty}^+\cdot U_{C_0^+,\infty}^-\cdot T_0$. Let $u_-\in U_{C_0^+,\infty}^-$. Let $\Omega\in C_0^+$. By Proposition~\ref{p_Rou4.5}, we can write $u_-=\prod_{\alpha\in \Delta_-}X_{\alpha}(\underline{u_\alpha})$, with $(\underline{u_\alpha})\in \prod_{\alpha\in \Delta_-}(\ffg_{\alpha,\Z}\otimes \cF)$ satisfying $\omega(\underline{u_\alpha})+\alpha(\Omega)\geq 0$, for all $\alpha\in \Delta_-$. For $\alpha\in \Delta_-$, we have  $\omega(\underline{u_\alpha})+\alpha(\Omega)\geq 0$ and thus $\omega(\underline{u_\alpha})>0$. By Proposition~\ref{p_mrphsm}, we have $\pr_{\kk}(U_{C_0^+,\infty}^{-})=\{1\}$. Consequently, \[\pr_{\kk}(G_{C_0^+})\subset \fU^+(\kk)\cdot \{1\} \cdot \pr_{\kk}(T_0)\subset \fB(\kk).\]

  We need to prove the reverse inclusion. Let $g\in \fG^{\min}(\cO)$ be such that $\pr(g)\in \fB(\kk)$. Then $g\in G_0$ and thus  by Theorem~\ref{t_G_x}, we can write $g=u_+u_-n$, with $u_+\in U_{0,\infty}^+$, $u_-\in U_{0,\infty}^-$ and $n\in N_0$. By Proposition~\ref{p_fxt_zr}, we have $u_+$, $u_-$ and $n\in \fG^{\min}(\cO)$.  By assumption, we can write $\pr_{\kk}(g)=v_+ t$, with $v_+\in \fU^+(\kk)$ and $t\in \fT(\kk)$. Then: \[\pr_\kk(u_+)^{-1}v_+ t\pr_\kk(n)^{-1}\pr_\kk(u_-^{-1})=1.\] By the Birkhoff decomposition in $\fG(\kk)$ (\cite[Corollary 7.70 and Proposition B33]{marquis2018introduction}), we deduce: \[\pr_\kk(u_-)=1\text{ and }t\pr_\kk(n)^{-1}=1.\]

 As $\pr_\kk(\cO^\times)=\kk^\times$, $\pr_\kk:T_0\rightarrow \fT(\kk)$ is surjective. Therefore there exists $\tilde{t}\in T_0$ such that $\pr_\kk(\tilde{t})=t$. Then $t':=\tilde{t}n^{-1}\in \ker(\pr_\kk)\cap N$. By  Proposition~\ref{p_ker_pi_k} (3), $t'$ fixes pointwise $C_0^+\cup -C_0^+$. Therefore $t'$ stabilizes $\A$ and fixes an open subset of $\A$. As $t'$ induces an affine automorphism of $\A$, we deduce that $t'$ fixes $\A$ pointwise. Therefore $t'\in T_0$ and $n\in T_0$. 

We also proved that $u_-\in \ker(\pr_\kk)\cap U^-$. By Proposition~\ref{p_ker_pi_k} (1), we deduce $u_-\in U_{C_0^+,\infty}^-$. Moreover $U_{0,\infty}^+=U_{C^v_f,\infty}^+=U_{C_0^+,\infty}^+$ and thus $g=u_+u_-n\in U_{C_0^+,\infty}^+\cdot U_{C_0^+,\infty}^-\cdot  T_0=G_{C_0^+}$. Therefore $G_{C_0^+}=\pr^{-1}(\fB(\kk))$. 
\end{proof}

\addcontentsline{toc}{section}{Index}

\printindex



\addcontentsline{toc}{section}{Index of Symbols}

\printindex[notation]

\printindex[axiom]

\bibliographystyle{amsalpha}
\bibliography{bibliographie.bib}

\providecommand{\bysame}{\leavevmode\hbox to3em{\hrulefill}\thinspace}
\providecommand{\MR}{\relax\ifhmode\unskip\space\fi MR }
\providecommand{\MRhref}[2]{%
  \href{http://www.ams.org/mathscinet-getitem?mr=#1}{#2}
}
\providecommand{\href}[2]{#2}
\begin{thebibliography}{BPHR25}

\bibitem[AFJ26]{appenzeller2026morphisms}
Raphael Appenzeller, Xenia Flamm, and Victor Jaeck, \emph{Morphisms of
  generalized affine buildings}, Preprint, {arXiv}:2601.04130 [math.{GR}]
  (2026), 2026.

\bibitem[AH19]{abdellatif2019completed}
Ramla Abdellatif and Auguste H\'ebert, \emph{Completed {I}wahori-{H}ecke
  algebras and parahoric {H}ecke algebras for {K}ac-{M}oody~groups over local
  fields}, Journal de l'\'Ecole polytechnique --- Math\'ematiques \textbf{6}
  (2019), 79--118 (en).

\bibitem[AVM10]{abramenko2010intersection}
Peter Abramenko and Hendrik Van~Maldeghem, \emph{Intersections of apartments},
  J. Comb. Theory, Ser. A \textbf{117} (2010), no.~4, 440--453 (English).

\bibitem[BB05]{bjorner2005combinatorics}
Anders Bj\"{o}rner and Francesco Brenti, \emph{Combinatorics of {C}oxeter
  groups}, Graduate Texts in Mathematics, vol. 231, Springer, New York, 2005.
  \MR{2133266}

\bibitem[BG09]{bruns2009polytopes}
Winfried Bruns and Joseph Gubeladze, \emph{Polytopes, rings, and {$K$}-theory},
  Springer Monographs in Mathematics, Springer, Dordrecht, 2009. \MR{2508056}

\bibitem[BG23]{bozec2023MV}
Tristan Bozec and St{\'e}phane Gaussent, \emph{Mirkovi{\'c}-vilonen polytopes
  and masures}, J. Lie Theory \textbf{33} (2023), no.~3, 719--746 (English).

\bibitem[BGR16]{bardy2016iwahori}
Nicole Bardy{-Panse}, St\'ephane Gaussent, and Guy Rousseau,
  \emph{Iwahori-{H}ecke algebras for {K}ac-{M}oody groups over local fields},
  Pacific J. Math. \textbf{285} (2016), no.~1, 1--61. \MR{3554242}

\bibitem[BK11]{braverman2011spherical}
Alexander Braverman and David Kazhdan, \emph{The spherical {Hecke} algebra for
  affine {Kac}-{Moody} groups. {I}}, Ann. Math. (2) \textbf{174} (2011), no.~3,
  1603--1642 (English).

\bibitem[BKP16]{braverman2016iwahori}
Alexander Braverman, David Kazhdan, and Manish~M. Patnaik,
  \emph{Iwahori-{H}ecke algebras for {$p$}-adic loop groups}, Invent. Math.
  \textbf{204} (2016), no.~2, 347--442. \MR{3489701}

\bibitem[Bou81]{bourbaki1981elements}
Nicolas Bourbaki, \emph{\'{E}l\'{e}ments de math\'{e}matique}, Masson, Paris,
  1981, Groupes et alg\`ebres de Lie. Chapitres 4, 5 et 6. [Lie groups and Lie
  algebras. Chapters 4, 5 and 6]. \MR{647314}

\bibitem[Bou06]{bourbaki2006groupes2}
\bysame, \emph{{\'E}l{\'e}ments de math{\'e}matique. {Groupes} et alg{\`e}bres
  de {Lie}. {Chapitres} 2 et 3}, reprint of the 1972 original ed., Berlin:
  Springer, 2006 (French).

\bibitem[Bou07]{bourbaki2007groupes1}
\bysame, \emph{{\'E}l{\'e}ments de math{\'e}matique. {Groupes} et alg{\`e}bres
  de {Lie}. {Chapitre} 1}, reprint of the 1972 original ed., Berlin: Springer,
  2007 (French).

\bibitem[BPHR25]{bardy2025twin}
Nicole Bardy-Panse, Auguste H{\'e}bert, and Guy Rousseau, \emph{Twin masures
  associated with {Kac}-{Moody} groups over {Laurent} polynomials}, Ann.
  Represent. Theory \textbf{2} (2025), no.~3, 281--353 (English).

\bibitem[BS14]{bennett2014axiomatic}
Curtis~D. Bennett and Petra~N. Schwer, \emph{On axiomatic definitions of
  non-discrete affine buildings}, Adv. Geom. \textbf{14} (2014), no.~3,
  381--412, With an appendix by Koen Struyve. \MR{3228889}

\bibitem[BT72]{bruhat1972groupes}
Fran{\c{c}}ois Bruhat and Jacques Tits, \emph{Groupes r{\'e}ductifs sur un
  corps local}, Publications Math{\'e}matiques de l'IH{\'E}S \textbf{41}
  (1972), no.~1, 5--251.

\bibitem[BT84]{bruhat1984groupes}
\bysame, \emph{Groupes r{\'e}ductifs sur un corps local}, Publications
  Math{\'e}matiques de l'IH{\'E}S \textbf{60} (1984), no.~1, 5--184.

\bibitem[BV25]{bouthier2025geometric}
Alexis Bouthier and Eric Vasserot, \emph{On the geometric {Satake} equivalence
  for {Kac}-{Moody} groups}, Preprint, {arXiv}:2510.11466 [math.{RT}] (2025),
  2025.

\bibitem[Cha10]{charignon2010immeubles}
Cyril Charignon, \emph{Immeubles affines et groupes de {K}ac-{M}oody}, Ph.D.
  thesis, universit{\'e} Henri Poincar{\'e} {N}ancy 1, 2010.

\bibitem[Chi01]{chiswell2001introduction}
Ian Chiswell, \emph{Introduction to {$\Lambda$}-trees}, World Scientific
  Publishing Co., Inc., River Edge, NJ, 2001. \MR{1851337}

\bibitem[EP05]{engler2005valued}
Antonio~J. Engler and Alexander Prestel, \emph{Valued fields}, Springer Monogr.
  Math., Berlin: Springer, 2005 (English).

\bibitem[Gar95]{garland1995cartan}
H.~Garland, \emph{A {C}artan decomposition for {$p$}-adic loop groups}, Math.
  Ann. \textbf{302} (1995), no.~1, 151--175. \MR{1329451}

\bibitem[GI63]{goldman1963space}
O.~Goldman and N.~Iwahori, \emph{The space of p-adic norms}, Acta Math.
  \textbf{109} (1963), 137--177 (English).

\bibitem[GR08]{gaussent2008kac}
St{\'e}phane Gaussent and Guy Rousseau, \emph{Kac-{M}oody groups, hovels and
  {L}ittelmann paths}, Annales de l'institut Fourier, vol.~58, 2008,
  pp.~2605--2657.

\bibitem[GR14]{gaussent2014spherical}
\bysame, \emph{Spherical {H}ecke algebras for {K}ac-{M}oody groups over local
  fields}, Annals of Mathematics \textbf{180} (2014), no.~3, 1051--1087.

\bibitem[H{\'e}b17]{hebert2017gindikin}
Auguste H{\'e}bert, \emph{{G}indikin-{K}arpelevich finiteness for {K}ac-{M}oody
  groups over local fields}, International Mathematics Research Notices
  \textbf{2017} (2017), no.~22, 7028--7049.

\bibitem[H{\'e}b18]{hebert2018study}
\bysame, \emph{{Study of masures and of their applications in arithmetic.
  English version}}, hal.archives ouvertes tel-01856620v1 (2018).

\bibitem[H{\'{e}}b20]{hebert2020new}
Auguste H{\'{e}}bert, \emph{A {N}ew {A}xiomatics for {M}asures}, Canad. J.
  Math. \textbf{72} (2020), no.~3, 732--773. \MR{4098599}

\bibitem[H{\'e}b21]{hebert2021distances}
Auguste H{\'e}bert, \emph{Distances on a masure}, Transform. Groups \textbf{26}
  (2021), no.~4, 1331--1363 (English).

\bibitem[H{\'e}b22]{hebert2022new}
\bysame, \emph{A new axiomatics for masures. {II}}, Adv. Geom. \textbf{22}
  (2022), no.~4, 513--522 (English).

\bibitem[H{\'e}b25]{hebert2025topologies}
\bysame, \emph{Topologies on split {Kac}-{Moody} groups over valued fields},
  Ann. Math. Blaise Pascal \textbf{32} (2025), no.~1, 77--121 (English).

\bibitem[HO89]{hahn1989classical}
Alexander~J. Hahn and O.~Timothy O'Meara, \emph{The classical groups and
  {$K$}-theory}, Grundlehren der mathematischen Wissenschaften [Fundamental
  Principles of Mathematical Sciences], vol. 291, Springer-Verlag, Berlin,
  1989, With a foreword by J. Dieudonn\'{e}. \MR{1007302}

\bibitem[HP24]{hebert2024affine}
Auguste H{\'e}bert and Paul Philippe, \emph{On affine {Kazhdan}-{Lusztig}
  {R}-polynomials for {Kac}-{Moody} groups}, Preprint, {arXiv}:2410.04872
  [math.{RT}] (2024), 2024.

\bibitem[HUL12]{hiriart2012fundamentals}
Jean-Baptiste Hiriart-Urruty and Claude Lemar{\'e}chal, \emph{Fundamentals of
  convex analysis}, Springer Science \& Business Media, 2012.

\bibitem[Kac94]{kac1994infinite}
Victor~G Kac, \emph{Infinite-dimensional {L}ie algebras}, vol.~44, Cambridge
  university press, 1994.

\bibitem[KM08]{kapovich2008path}
Michael Kapovich and John~J. Millson, \emph{A path model for geodesics in
  {E}uclidean buildings and its applications to representation theory}, Groups
  Geom. Dyn. \textbf{2} (2008), no.~3, 405--480. \MR{2415306 (2009d:20069)}

\bibitem[KP85]{kac1985defining}
V.~G. Kac and D.~H. Peterson, \emph{Defining relations of certain infinite
  dimensional groups}, {\'E}lie {Cartan} et les math{\'e}matiques
  d'aujourd'hui, {The} mathematical heritage of {Elie} {Cartan}, {Semin}.
  {Lyon} 1984, {Ast{\'e}risque}, {No}.{Hors} {S{\'e}r}. 1985, 165-208 (1985).,
  1985.

\bibitem[KP23]{kaletha2023bruhat}
Tasho Kaletha and Gopal Prasad, \emph{Bruhat-tits theory. {A} new approach},
  New Math. Monogr., vol.~44, Cambridge: Cambridge University Press, 2023
  (English).

\bibitem[Kum02]{kumar2002kac}
Shrawan Kumar, \emph{Kac-{M}oody groups, their flag varieties and
  representation theory}, Progress in Mathematics, vol. 204, Birkh\"auser
  Boston, Inc., Boston, MA, 2002. \MR{1923198}

\bibitem[Mar18]{marquis2018introduction}
Timoth{\'e}e Marquis, \emph{An introduction to {Kac}-{Moody} groups over
  fields}, EMS Textb. Math., Z{\"u}rich: European Mathematical Society (EMS),
  2018 (English).

\bibitem[Mat89]{mathieu1989construction}
Olivier Mathieu, \emph{Construction d'un groupe de {K}ac-{M}oody et
  applications}, Compositio Math. \textbf{69} (1989), no.~1, 37--60.
  \MR{986812}

\bibitem[MM25]{marquis2025presentation}
Timoth{\'e}e Marquis and Bernhard M{\"u}hlherr, \emph{Presentation and
  uniqueness of {Kac}-{Moody} groups over local rings}, Preprint,
  {arXiv}:2510.11272 [math.{GR}] (2025), 2025.

\bibitem[Mut19]{muthiah2019double}
Dinakar Muthiah, \emph{Double-affine {K}azhdan-{L}usztig polynomials via
  masures}, arXiv preprint arXiv:1910.13694 (2019).

\bibitem[Par00]{parreau2000immeubles}
Anne Parreau, \emph{Immeubles affines: construction par les normes},
  Crystallographic Groups and Their Generalizations: Workshop, Katholieke
  Universiteit Leuven Campus Kortrijk, Belgium, May 26-28, 1999, vol. 262,
  American Mathematical Soc., 2000, p.~263.

\bibitem[Par23]{parreau2023affine}
\bysame, \emph{Affine buildings: construction by norms and study of
  isometries}, Innov. Incidence Geom. \textbf{20} (2023), no.~2-3, 471--517
  (English).

\bibitem[R{\'e}m02]{remy2002groupes}
Bertrand R{\'e}my, \emph{Groupes de {K}ac-{M}oody d\'eploy\'es et presque
  d\'eploy\'es}, Ast\'erisque (2002), no.~277, viii+348. \MR{1909671}

\bibitem[Rou06]{rousseau2006groupes}
Guy Rousseau, \emph{Kac-{Moody} groups split over a local field, microaffine
  buildings}, Compos. Math. \textbf{142} (2006), no.~2, 501--528 (French).

\bibitem[Rou11]{rousseau2011masures}
\bysame, \emph{Masures affines}, Pure and Applied Mathematics Quarterly
  \textbf{7} (2011), no.~3, 859--921.

\bibitem[Rou16]{rousseau2016groupes}
\bysame, \emph{Groupes de {K}ac-{M}oody d\'eploy\'es sur un corps local {II}.
  {M}asures ordonn\'ees}, Bull. Soc. Math. France \textbf{144} (2016), no.~4,
  613--692. \MR{3562609}

\bibitem[Rou17]{rousseau2017almost}
\bysame, \emph{{A}lmost split {K}ac–{M}oody groups over ultrametric fields},
  Groups Geometry, and Dynamics \textbf{11} (2017), 891--975.

\bibitem[Rou23]{rousseau2023euclidean}
\bysame, \emph{Euclidean buildings. {Geometry} and group actions}, EMS Tracts
  Math., vol.~35, Berlin: European Mathematical Society (EMS), 2023 (English).

\bibitem[Rou25]{rousseau2025split}
Guy Rousseau, \emph{Split kac-moody groups over a local field, ii. ordered
  masures}, 2025.

\bibitem[RTW15]{remy2015BTbuildings}
Bertrand R{\'e}my, Amaury Thuillier, and Annette Werner, \emph{Bruhat-tits
  buildings and analytic geometry}, Berkovich spaces and applications. Based on
  a workshop, Santiago de Chile, Chile, January 2008 and a summer school,
  Paris, France, June 2010, Cham: Springer, 2015, pp.~141--202 (English).

\bibitem[Ser66]{serre1966algebres}
Jean-Pierre Serre, \emph{Alg\`ebres de {L}ie semi-simples complexes}, W. A.
  Benjamin, inc., New York-Amsterdam, 1966. \MR{0215886}

\bibitem[Ser77]{serre1977arbres}
\bysame, \emph{Arbres, amalgames, {{\(\mathrm{SL}_2\)}}. {R{\'e}dig{\'e}} avec
  la collaboration de {Hyman} {Bass}}, Ast{\'e}risque, vol.~46, Soci{\'e}t{\'e}
  Math{\'e}matique de France (SMF), Paris, 1977 (French).

\bibitem[Ser80]{serre1980trees}
\bysame, \emph{Trees. {Transl}. from the {French} by {John} {Stillwell}},
  Berlin-{Heidelberg}-{New} {York}: {Springer}-{Verlag}. {IX}, 142 p., 1980.

\bibitem[Ser01]{serre2001complex}
\bysame, \emph{Complex semi-simple lie algebras, translated from the french by
  ga jones, reprint of the 1987 edition}, Springer Monographs in Mahematics,
  Springer Verlag, Berlin (2001).

\bibitem[Ste16]{steinberg2016lectures}
Robert Steinberg, \emph{Lectures on {Chevalley} groups}, Univ. Lect. Ser.,
  vol.~66, Providence, RI: American Mathematical Society (AMS), 2016 (English).

\bibitem[Tit87]{tits1987uniqueness}
Jacques Tits, \emph{Uniqueness and presentation of {K}ac-{M}oody groups over
  fields}, J. Algebra \textbf{105} (1987), no.~2, 542--573. \MR{873684}

\bibitem[Tit13]{tits2013resume}
\bysame, \emph{R\'esum\'es des cours au {C}oll\`ege de {F}rance 1973--2000},
  Documents Math\'ematiques (Paris) [Mathematical Documents (Paris)], vol.~12,
  Soci\'et\'e Math\'ematique de France, Paris, 2013. \MR{3235648}

\end{thebibliography}

\end{document}